\documentclass{amsart}
\usepackage[utf8]{inputenc}

\usepackage{comment}

\usepackage{geometry}
\usepackage{stmaryrd}
\usepackage{accents}
\usepackage{upgreek}
\usepackage[normalem]{ulem}

\usepackage[dvipsnames]{xcolor}

\newcommand\myshade{85}
\colorlet{mylinkcolor}{Red}
\colorlet{mycitecolor}{Green}
\colorlet{myurlcolor}{Plum}

\usepackage{epsdice}

\usepackage{enumitem}

\usepackage[titletoc]{appendix}

\usepackage[OT2,T1]{fontenc}

\usepackage[final]{microtype}

\usepackage[all,cmtip]{xy}
\usepackage{tikz}
\usepackage{tikz-cd}

\tikzcdset{
  cells={font=\everymath\expandafter{\the\everymath\displaystyle}},
}

\tikzcdset{
  circled number/.style={
    description,
    /tikz/circle,
    /tikz/draw,
    /tikz/fill=white,
    /tikz/inner sep=1pt,
    /tikz/font=\scriptsize
  }
}

\makeatletter
\tikzcdset{
  eq node/.style={
    commutative diagrams/math mode=false, anchor=center},
  eq/.style={
    phantom, description,
    /tikz/every to/.append style={
      edge node={node[commutative diagrams/eq node]
        {\@eqnswtrue\make@display@tag\ltx@label{{#1}}}}}}}
        
\tikzcdset{scale cd/.style={every label/.append style={scale=#1},
    cells={nodes={scale=#1}}}}
\makeatother

\tikzset{
  trim node/.default=1cm,
  trim node/.style={
    overlay,
    append after command={
      ([xshift={+#1}]\tikzlastnode.north west)
      ([xshift={+-#1}]\tikzlastnode.south east)}},
  down and trim/.default=1cm,
  down and trim/.style={
    yshift=-(\pgfmatrixcurrentcolumn-1)*1.5\baselineskip,
    trim node={#1}},
  downup and trim/.default=1cm,
  downup and trim/.style={
    yshift=iseven(\pgfmatrixcurrentcolumn) ? -1.5\baselineskip : 0pt,
    trim node={#1}},
  -|/.style={to path={-|(\tikztotarget)\tikztonodes}},
  |-/.style={to path={|-(\tikztotarget)\tikztonodes}},
  -| sl/.style={-|, xslant=-1},
  |- sl/.style={|-, xslant= 1},
  center picture/.style={
    trim left=(current bounding box.center),
    trim right=(current bounding box.center)}}

\usepackage{rotating}
\usepackage{floatpag}

\usepackage[draft=false,linktocpage=true,pdfusetitle]{hyperref}

\hypersetup{
    linkcolor = mylinkcolor!\myshade!black,
    citecolor = mycitecolor!\myshade!black,
    urlcolor = myurlcolor!\myshade!black,
    colorlinks = true
}
\usepackage[capitalise]{cleveref}
\crefformat{equation}{\ensuremath{(#2#1#3)}}
\crefmultiformat{equation}{\ensuremath{(#2#1#3)}}{ and~\ensuremath{(#2#1#3)}}{, \ensuremath{(#2#1#3)}}{, and~\ensuremath{(#2#1#3)}}
\crefname{enumi}{\unskip}{\unskip}

\AtBeginEnvironment{appendices}{\crefalias{section}{appendix}}

\AtBeginDocument{
    \def\MR#1{}
}

\usepackage{colonequals}
\usepackage{mathrsfs}

\usepackage{relsize}
\usepackage[bbgreekl]{mathbbol}
\usepackage{amsfonts}
\DeclareSymbolFontAlphabet{\mathbb}{AMSb}
\DeclareSymbolFontAlphabet{\mathbbl}{bbold}

\mathchardef\mhyphen="2D

\renewcommand{\AA}{\mathbf{A}}
\newcommand{\BB}{\mathbf{B}}
\newcommand{\CC}{\mathbf{C}}
\newcommand{\DD}{\mathbf{D}}

\newcommand{\FF}{\mathbf{F}}
\newcommand{\GG}{\mathbf{G}}
\newcommand{\NN}{\mathbf{N}}

\newcommand{\PP}{\mathbf{P}}
\newcommand{\QQ}{\mathbf{Q}}
\newcommand{\Q}{\mathbf{Q}}

\newcommand{\TT}{\mathbf{T}}
\newcommand{\ZZ}{\mathbf{Z}}
\newcommand{\Z}{\mathbf{Z}}
\newcommand{\bL}{\mathbf{L}}

\newcommand{\cA}{\mathcal{A}}
\newcommand{\cB}{\mathcal{B}}
\newcommand{\cC}{\mathcal{C}}
\newcommand{\cD}{\mathcal{D}}
\newcommand{\cE}{\mathcal{E}}
\newcommand{\cF}{\mathcal{F}}
\newcommand{\F}{\mathcal{F}}
\newcommand{\cG}{\mathcal{G}}
\newcommand{\G}{\mathcal{G}}
\newcommand{\cH}{\mathcal{H}}
\newcommand{\cI}{\mathcal{I}}

\newcommand{\cL}{\mathscr{L}}
\newcommand{\cM}{\mathscr{M}}
\newcommand{\cO}{\mathcal{O}}
\let\DanishO\O
\renewcommand{\O}{\mathcal{O}}
\newcommand{\cP}{\mathcal{P}}
\newcommand{\cQ}{\mathcal{Q}}

\newcommand{\cV}{\mathscr{V}}
\newcommand{\cX}{\mathscr{X}}

\newcommand{\rL}{\mathrm{L}}

\newcommand{\FFF}{\mathrm{FF}}

\newcommand{\m}{\mathfrak{m}}

\newcommand{\n}{\mathfrak{n}}

\newcommand{\rR}{\mathrm{R}}

\newcommand{\uGamma}{\underline{\Gamma}}
\newcommand{\uRGamma}{\rR\underline{\Gamma}}
\newcommand{\uQ}{\ud{\Q}}
\newcommand{\uZ}{\ud{\Z}}

\newcommand{\Adic}{\mathrm{Adic}}
\newcommand{\aff}{\mathrm{aff}}
\newcommand{\alg}{\mathrm{alg}}
\newcommand{\Ab}{\mathrm{Ab}}
\newcommand{\an}{\mathrm{an}}
\newcommand{\AN}{\mathrm{AN}}
\newcommand{\BBC}{\mathscr{BC}}

\newcommand{\can}{\mathrm{can}}
\newcommand{\Cat}{\mathrm{Cat}}
\newcommand{\CCart}{\mathrm{CCart}}
\newcommand{\cHom}{\mathscr{H}\kern-.5pt om}
\newcommand{\cl}{\mathrm{cl}}

\newcommand{\cond}{\mathrm{cond}}
\newcommand{\Cond}{\mathrm{Cond}}

\newcommand{\cont}{\mathrm{cont}}

\newcommand{\diam}{\lozenge}
\newcommand{\Diam}{\mathrm{Diam}}

\newcommand{\et}{\mathrm{\acute{e}t}}

\newcommand{\Idem}{\mathrm{Idem}}

\newcommand{\incl}{\mathrm{incl}}

\newcommand{\op}{\mathrm{op}}

\newcommand{\Perf}{\mathrm{Perf}}
\newcommand{\Perfd}{\mathrm{Perfd}}
\newcommand{\proet}{\mathrm{pro\acute{e}t}}
\newcommand{\Pro}{\mathrm{Pro}}

\newcommand{\qc}{\mathrm{qc}}

\newcommand{\qproet}{\mathrm{qpro\acute{e}t}}
\newcommand{\rad}{\mathrm{rad}}
\newcommand{\red}{\mathrm{red}}
\renewcommand{\ring}{\mathrm{ring}}

\newcommand{\sol}{\square}
\newcommand{\Solid}{\mathrm{Solid}}
\renewcommand{\sp}{\mathrm{sp}}
\newcommand{\Vect}{\mathrm{Vect}}
\newcommand{\wphi}{\widetilde{\varphi}}
\newcommand{\wX}{\widetilde{X}}
\newcommand{\zc}{\mathrm{zc}}

\DeclareMathOperator{\Adj}{Adj}
\DeclareMathOperator{\Ani}{Ani}
\DeclareMathOperator{\Aut}{Aut}
\DeclareMathOperator{\BC}{BC}
\DeclareMathOperator{\CAlg}{CAlg}

\DeclareMathOperator{\coker}{coker}
\DeclareMathOperator*{\colim}{colim}

\DeclareMathOperator{\End}{End}
\DeclareMathOperator{\eq}{eq}

\DeclareMathOperator{\Ev}{Ev}

\DeclareMathOperator{\Ext}{Ext}

\DeclareMathOperator{\fib}{fib}
\DeclareMathOperator{\Frac}{Frac}
\DeclareMathOperator{\Fun}{Fun}

\DeclareMathOperator{\grp}{grp}
\DeclareMathOperator{\Grp}{Grp}

\DeclareMathOperator{\Hh}{H}
\DeclareMathOperator*{\hocolim}{hocolim}
\DeclareMathOperator{\Hom}{Hom}
\DeclareMathOperator{\hyp}{hyp}
\DeclareMathOperator{\id}{id}
\DeclareMathOperator{\KM}{KM}
\DeclareMathOperator{\Latt}{Latt}
\DeclareMathOperator{\Map}{Map}

\DeclareMathOperator{\Mod}{Mod}

\DeclareMathOperator{\Ob}{Ob}

\DeclareMathOperator{\PD}{PD}
\DeclareMathOperator{\PF}{PF}
\DeclareMathOperator{\Pic}{Pic}

\DeclareMathOperator{\PShv}{PShv}
\DeclareMathOperator{\res}{res}
\DeclareMathOperator{\RH}{RH}
\DeclareMathOperator{\Ring}{Ring}

\DeclareMathOperator{\Set}{Set}
\DeclareMathOperator{\Shv}{Shv}
\DeclareMathOperator{\Sol}{Sol}

\DeclareMathOperator{\Spa}{Spa}
\DeclareMathOperator{\Spd}{Spd}
\DeclareMathOperator{\Spec}{Spec}

\DeclareMathOperator{\supp}{supp}

\DeclareMathOperator{\tr}{tr}
\DeclareMathOperator{\Tr}{Tr}

\newcommand{\cal}{\mathcal}
\newcommand{\ov}{\overline}
\renewcommand{\rm}{\mathrm}
\newcommand{\ud}{\underline}

\newcommand{\wdh}{\widehat}
\newcommand{\xr}{\xrightarrow}

\newcommand{\abs}[1]{\lvert#1\rvert}
\newcommand{\suchthat}{\;\ifnum\currentgrouptype=16 \middle\fi\vert\;}
\newcommand\restr[2]{{\left.\kern-\nulldelimiterspace#1\vphantom{\big|}\right|_{#2}}}
\newcommand{\blank}{{-}}

\usepackage{mathtools}
\usepackage{amsmath}

\usepackage{amsthm,amssymb}

\newcommand{\xtwoheadrightarrow}[2][]{
  \xrightarrow[#1]{#2}\mathrel{\mkern-14mu}\rightarrow
}

\numberwithin{equation}{subsection}

\theoremstyle{plain}
\newtheorem{theorem}[equation]{Theorem}
\newtheorem{thm}[equation]{Theorem}

\newtheorem{proposition}[equation]{Proposition}
\newtheorem{lemma}[equation]{Lemma}
\newtheorem{claim}[equation]{Claim}
\newtheorem{corollary}[equation]{Corollary}
\newtheorem{cor}[equation]{Corollary}
\newtheorem{construction/proposition}[equation]{Construction/Proposition}

\theoremstyle{definition}
\newtheorem{definition}[equation]{Definition}
\newtheorem{defn}[equation]{Definition}
\newtheorem{notation}[equation]{Notation}

\newtheorem{construction}[equation]{Construction}
\newtheorem{variant}[equation]{Variant}

\newtheorem{example}[equation]{Example}
\newtheorem{examples}[equation]{Examples}

\newtheorem{remark}[equation]{Remark}
\newtheorem{rmk}[equation]{Remark}

\newtheorem{warning}[equation]{Warning}
\newtheorem{setup}[equation]{Setup}

\title{ Poincar\'e Duality for pro-\'etale \texorpdfstring{$\mathbf{Q}_p$}{Qp}-local systems}

\author{Shizhang Li}
\author{Wies\l{}awa Nizio\l{}}
\author{Emanuel Reinecke}
\author{Bogdan Zavyalov}

\address[Shizhang Li]{Morningside Center of Mathematics and State Key Laboratory of Mathematics Sciences,
Academy of Mathematics and Systems Science, Chinese Academy of Sciences, Beijing 100190, China}
\email{lishizhang@amss.ac.cn}

\address[Wies\l{}awa Nizio\l{}]{CNRS, IMJ-PRG, Sorbonne Universit\'e, 4 place Jussieu, 75005 Paris, France}
\email{wieslawa.niziol@imj-prg.fr}

\address[Emanuel Reinecke]{Institut for Matematiske Fag, K{\o}benhavns Universitet, Universitetsparken 5, 2100 K{\o}benhavn \DanishO, Denmark}
\email{reinecke@math.ku.dk}

\address[Bogdan Zavyalov]{University of Maryland, Department of Mathematics, 4176 Campus Dr, College Park, MD 20742, USA}
\email{bogd.zavyalov@gmail.com}

\begin{document}

\begin{abstract} 
We prove finiteness and Poincar\'e duality for pro-\'etale $\QQ_p$-local systems on proper $p$-adic rigid-analytic spaces in the framework of Banach--Colmez spaces and establish their optimal cohomological vanishing bounds.
As a consequence, we also obtain ordinary finiteness and duality for their arithmetic pro-\'etale cohomology.
We deduce these results from their analogs for perfect complexes on the Fargues--Fontaine curve, which in turn reduce to Poincar\'e duality for perfect complexes over period rings.
We give a simple proof of the latter duality via the same diagrammatic argument as in our previous work for finite coefficients. 
Along the way, we establish optimal $v$-descent for perfect complexes over certain period rings on affinoid perfectoid spaces. 
\end{abstract}

\maketitle

\tableofcontents
\addtocontents{toc}{\protect\setcounter{tocdepth}{1}}

\section{Introduction}

\subsection{Main results}

Let $K$ be a nonarchimedean field of mixed characteristic $(0,p)$ and let $C\coloneqq \wdh{\ov{K}}$ be its completed algebraic closure. Berkovich and Huber independently developed theories of \'etale cohomology for rigid-analytic spaces over $K$ in \cite{Berkovich} and \cite{Huber-etale}, respectively; this theory was later extended to diamonds and small $v$-sheaves in \cite{diamonds}. When the coefficients are $\ZZ/n$ with $n$ prime to $p$, the resulting theory closely parallels the algebraic theory of \'etale cohomology developed in \cite{SGA4, SGA41/2, SGA5}. For $p$-power torsion coefficients, however, many of these favorable properties break down.

Indeed, \'etale cohomology with $\FF_p$-coefficients is considerably subtler: the cohomology groups may be infinite-dimensional, may depend on the choice of algebraically closed ground field, and need not satisfy either the K\"unneth formula or a general form of Poincar\'e duality. Substantial progress in understanding this theory began with Scholze's seminal paper \cite{Scholze-Hodge}, which, among many other results, established the finiteness of \'etale $\FF_p$-cohomology for proper rigid-analytic spaces. Subsequently, Bhatt--Hansen \cite{BH} developed a robust theory of Zariski-constructible sheaves and Verdier duality; Zavyalov \cite{Z-thesis} and Mann \cite{Mann-thesis} independently established Poincar\'e duality for smooth proper spaces; and \cite{LRZ24} proved a general form of Poincar\'e duality for arbitrary proper spaces.

The goal of this paper is to extend this general form of Poincar\'e duality from torsion coefficients to $\uQ_p$-local systems on the pro-\'etale site of $X$, or equivalently on its $v$-site.\footnote{
For technical reasons, we prefer to work with the $v$-site in this paper. The category of $\uQ_p$-local systems as well as their cohomology are insensitive to the choice of the pro-\'etale or $v$-topology; see \cref{lemma:no-difference-proetale-v} and \cref{qproet-v-pushforward-etale}. However, to have the correct notion of Banach--Colmez spaces below; one has to use either the big pro-\'etale topology or the $v$-topology.
For consistency, we use the $v$-topology in this introduction, but all results also remain valid with the big pro-\'etale topology.}
The torsion-coefficient theory already yields the expected statements for constant $\uQ_p$-coefficients. In particular, if $X$ is proper over $C$, then $\Hh^i_v(X,\uQ_p)$ is finite-dimensional; if $X$ is moreover smooth, proper, and of equidimension $d$, there is a canonical Poincar\'e duality isomorphism 
\[ 
    \Hh^i_v(X,\uQ_p)^{\vee} \simeq \Hh^{2d-i}_v(X,\uQ_p(d)). 
\] 
The same passage from torsion to rational coefficients applies more generally to $\uQ_p$-local systems which admit a $\uZ_p$-lattice. In fact, under an extra mild hypothesis, these finiteness and duality statements were established in \cite{LLZ} even before the work of Mann and Zavyalov. 

However, cohomology of general $\uQ_p$-local systems, whose study was initiated by Kedlaya--Liu in \cite{Kedlaya-Liu-3} and later by Colmez--Gilles--Nizio{\l} in \cite{OB},
present genuinely new phenomena.
Unlike their algebraic counterparts, these local systems need not admit a $\uZ_p$-lattice, and their cohomology can be infinite-dimensional even when $X$ is a smooth projective curve; see \cref{intro:main-example}. Thus, neither naive finiteness nor naive Poincar\'e duality can hold for arbitrary pro-\'etale $\uQ_p$-local systems.
The appropriate replacement is to retain the variation over the base: rather than considering the individual groups $\Hh^i(X,\bL)$, we consider the associated cohomology sheaves on the $v$-site of $\Spa(C,\O_C)$. We then formulate and prove finiteness and duality statements for these sheaves, whose natural finiteness class is provided by Banach--Colmez spaces. We recall their definition below. 

\begin{defn}[{\cite{Colmez, LB-thesis}}] The category of \emph{Banach--Colmez spaces} is the smallest weak Serre subcategory $\BBC(C) \subset \mathrm{Shv}(\Spa(C, \O_C)_v; \uQ_p)$ containing $\ud{\QQ}_p$ and $\O_{\Spa(C, \O_C)}$. We denote by 
\[
\cD^{b}_{\BBC}(\Spa(C, \O_C)_v; \uQ_p) \subset \cD(\Spa(C, \O_C)_v;\uQ_p)
\]
the full subcategory consisting of bounded objects whose cohomology sheaves lie in $\BBC(C)$. 
\end{defn}

We refer to \cite[\textsection 14]{diamonds} for the definition of the $v$-topology and to \cref{defn:period-sheaves} for the definition of $\O_{\Spa(C,\O_C)}$. Throughout the introduction, we freely identify an analytic adic space over $\Spa(\QQ_p,\ZZ_p)$ with its associated diamond.

\begin{thm}[Geometric Duality; {\cref{singular-PD} and \cref{singular-tor-amplitude}}]\label{thm:intro-PD-proper-qp} Let $f \colon X \to \Spa(C,\cO_C)$ be a proper rigid-analytic space of dimension $d$, let $\omega^{\Q_p}_X$ be the dualizing complex on $X$, and let $\cE \in \Perf_v(X;\uQ_p) \cap \cD^{[r, r']}(X_v; \uQ_p)$ for some $r,r'\in \ZZ\cup \{-\infty, \infty\}$. Then the following statements hold:
\begin{enumerate}[leftmargin=*,label={\upshape{(\roman*)}}]
    \item\label{thm:intro-PD-proper-qp-1} $\rR f_{v, *} \cE$ lies in $\cD^{[r, r'+2d]}_{\BBC}\bigl(\Spa(C, \O_C)_v; \uQ_p\bigr)$;
    \item\label{thm:intro-PD-proper-qp-2} there is a canonical trace map $\Tr_f\colon \rR f_{v, *} \omega^{\Q_p}_X \to \uQ_p$ such that the induced duality morphism
    \[
    \PD_f \colon \rR f_{v, *} \rR \cHom_{\cD(X_v;\uQ_p)}(\cE, \omega^{\Q_p}_X) \longrightarrow \rR\cHom_{\cD(\Spa(C)_v;\uQ_p)}(\rR f_{v, *} \cE,\uQ_p)
    \]
    defined in \cref{BI-proper-trace}\cref{BI-proper-PD-morphism} is an equivalence. 
\end{enumerate}
\end{thm}

\begin{remark} 
The recent paper \cite{ALBM} establishes the finiteness and duality statements of \cref{thm:intro-PD-proper-qp} as well as their analogs for arbitrary smooth proper morphisms $f\colon X\to Y$. 
The paper \cite{ALBM}, however, does not address singular spaces or the cohomological bounds for $\rR f_{v,*}\cE$. 

For smooth proper morphisms $f\colon X\to Y$, our methods also provide an alternative and more elementary proof.
Rather than relying on the full strength of an appropriate six-functor formalism, 
our argument builds on the diagrammatic approach of \cite{LRZ24}, together with the trace and cycle class maps constructed in \emph{loc.~cit.} 
\end{remark}

Part~\ref{thm:intro-PD-proper-qp-1} of \cref{thm:intro-PD-proper-qp} should be viewed as the appropriate replacement for finite-dimensionality in this setting. When $\cE$ is a $\uQ_p$-local system, it implies that $\rR\Gamma_v(X,\cE)$ is concentrated in degrees $[0,2d]$ and that each $\Hh^i_v(X,\cE)$ can be constructed from $\QQ_p$ and $C$ in finitely many steps.
The appearance of $C$ is essential: such contributions occur even for (analytic) local systems on the (analytification) of a smooth proper algebraic curve, as the following example shows.

\begin{example}[{\cite[Comment of PS]{Hansen:Counterexample}}]\label{intro:main-example}
Consider the Tate elliptic curve 
\[ 
    f\colon E\coloneqq \mathbf{G}_m^{\an}/p^{\mathbf{Z}} \longrightarrow \Spa(\QQ_p,\ZZ_p). 
\] 
This is a smooth proper rigid-analytic space over $\Spa(\QQ_p,\ZZ_p)$, equipped with the natural $\ZZ$-torsor 
\[ 
    \pi\colon \GG_m^{\an}\longrightarrow E. 
\] 
Let $\cL_+$ and $\cL_-$ be the rank-one $\uQ_p$-local systems trivialized by $\pi$ whose monodromy characters are, respectively, 
\[ 
    \rho_+\colon \ZZ\longrightarrow \QQ_p^\times, \qquad 1\longmapsto p, 
\] and 
\[ 
\rho_-\colon \ZZ\longrightarrow \QQ_p^\times, \qquad 1\longmapsto p^{-1}. 
\] 
Throughout this example, the subscript $(-)_C$ denotes base change to $\Spa(C,\O_C)$.

Since neither $p$ nor $p^{-1}$ lies in $\ZZ_p^\times$, neither $\cL_+$ nor $\cL_-$ admits a $\uZ_p$-lattice. Their cohomology groups are, in fact, infinite-dimensional. To see this, the Hochschild--Serre spectral sequence associated with $\pi$ gives 
\[ 
    \mathrm{E}_2^{i,j} = \Hh^i\bigl( \ZZ, \Hh^j_v(\GG_{m,C}^{\an},\uQ_p)\otimes_{\QQ_p} \rho_+ \bigr) \Longrightarrow \Hh^{i+j}_v(E_C,\cL_{+,C}), 
\] and there is an analogous spectral sequence for $\cL_-$. On the other hand, \cite[Th.~1.8]{CDN20} computes the $v$-cohomology, or equivalently the pro-\'etale cohomology, of $\GG_{m,C}^{\an}$. More precisely, 
\[ 
    \Hh^0_v(\GG_{m,C}^{\an},\uQ_p)=\Q_p, \qquad \Hh^i_v(\GG_{m,C}^{\an},\uQ_p)=0 \quad\text{for }i\geq 2, 
\] 
while $\Hh^1_v(\GG_{m,C}^{\an},\uQ_p)$ fits into a short exact sequence 
\[ 
    0 \longrightarrow \bigl(\O(\GG_{m,C}^{\an})/C\bigr)(-1) \longrightarrow \Hh^1_v(\GG_{m,C}^{\an},\uQ_p) \longrightarrow \QQ_p(-1) \longrightarrow 0. 
\] 
Combining this computation with the two Hochschild--Serre spectral sequences yields \begin{equation}\label{eqn:computation-cohomology} \Hh^i_v(E_C,\cL_{+,C}) \simeq \Hh^i_v(E_C,\cL_{-,C}) \simeq \begin{cases} C(-1), & i=1 \text{ or } 2, \\ 0, & \text{otherwise}. \end{cases} \end{equation} In particular, the nonzero cohomology groups in \eqref{eqn:computation-cohomology} are infinite-dimensional as $\QQ_p$-vector spaces.

The preceding computation can be refined to an equivalence at the level of derived pushforwards: 
\begin{equation}\label{eqn:intro-derived-pushforward}
    \rR f_{v,C,*}\cL_{+,C} \simeq \O_{\Spa C}(-1)[-1] \oplus \O_{\Spa C}(-1)[-2] \simeq \rR f_{v,C,*}\cL_{-,C}. 
\end{equation}
This also gives a concrete illustration of the duality assertion in \cref{thm:intro-PD-proper-qp}. Indeed, $\cL_{+,C}^{\vee}\simeq \cL_{-,C}$ and the dualizing complex of $E_C$ is $\uQ_p(1)[2]$. In combination with the formula $\rR\cHom_{\cD(\Spa(C)_v;\uQ_p)}\bigl(\O_{\Spa C},\uQ_p\bigr) \simeq \O_{\Spa C}(-1)[-1]$ proved in \cite[Th.~3.8]{ALB}, this yields the following sequence of isomorphisms:
\begin{align*} 
\rR f_{v,C,*} \rR\cHom_{\cD(X_v;\uQ_p)}\bigl( \cL_{+,C}, \uQ_p(1)[2] \bigr) &\simeq \rR f_{v,C,*}\cL_{-,C}(1)[2] \\ 
&\simeq \O_{\Spa C}[1]\oplus \O_{\Spa C} \\ 
&\simeq \rR{\cHom}_{\cD(\Spa(C)_v;\uQ_p)}\bigl( \O_{\Spa C}(-1)[-1] \oplus \O_{\Spa C}(-1)[-2], \uQ_p \bigr) \\
&\simeq \rR{\cHom}_{\cD(\Spa(C)_v;\uQ_p)}\bigl( \rR f_{v,C,*}\cL_{+,C}, \uQ_p \bigr). 
\end{align*} 
\end{example}

Finally, we show that these infinite-dimensional phenomena disappear completely when one works with the absolute cohomology of proper rigid-analytic spaces over a finite extension of $\QQ_p$. More precisely, we prove the following theorem.

\begin{thm}[Arithmetic Duality; {\cref{cor:arithmetic-duality-general}}]\label{thm:intro-arithmetic-PD-proper-qp} Let $K$ be a finite extension of $\QQ_p$, let $f\colon X\to \Spa(K, \O_K)$ be a proper rigid-analytic space of dimension $d$, let $\cE \in \Perf_v(X; \uQ_p) \cap \cD^{[r, r']}(X_v; \uQ_p)$. Set $\omega_X^{\mathrm{ar}}\coloneqq \omega^{\Q_p}_X(1)[2]$. Then the following statements hold:
\begin{enumerate}[leftmargin=*,label={\upshape{(\roman*)}}]
    \item\label{thm:intro-arithmetic-PD-proper-qp-1} the cohomology complex $\rR\Gamma_v(X, \cE)$ lies in $\cD^{[r, r'+2d+2]}_{\mathrm{coh}}(\QQ_p)$;
    \item\label{thm:intro-arithmetic-PD-proper-qp-2} there is a canonical trace map $\Tr_f\colon \rR\Gamma_v(X, \omega^{\mathrm{ar}}_X) \to \Q_p$ such that the induced duality morphism
    \[
    \PD_f \colon \rR \Hom_{\cD(X_v;\uQ_p)}(\cE, \omega^{\mathrm{ar}}_X) \longrightarrow \rR \Hom_{\cD(\QQ_p)}(\rR\Gamma_v(X, \cE), \Q_p)
    \]
    defined in \cref{defn:arithmetic-singular-duality-morphism}~\cref{defn:arithmetic-singular-duality-morphism-2} is an equivalence. 
\end{enumerate}
\end{thm}

We deduce \cref{thm:intro-arithmetic-PD-proper-qp} from \cref{thm:intro-PD-proper-qp} using Fontaine's theory of almost $\mathbf{C}_p$-representations, developed in \cite{FCp, Fontaine-Almost}. When $X$ is smooth, one can also deduce the finiteness and duality parts of \cref{thm:intro-arithmetic-PD-proper-qp} from \cite{ALBM}; cf.\ \cite[Rmk.~6.2.3]{ALBM}.
However, the proper case has not been addressed before.

The observation that arithmetic pro-\'etale cohomology, or equivalently arithmetic $v$-cohomology, exhibits simpler finiteness behavior than its geometric counterpart goes back to \cite{CGW23}. There, the authors proved a version of arithmetic Poincar\'e duality with constant coefficients for smooth partially proper rigid-analytic curves. This result was extended in the work of Zhenghui Li
\cite{ZL} to arbitrary smooth partially proper rigid-analytic spaces and constant coefficients. 

In order to illustrate how the infinite-dimensionality of cohomology disappears in the arithmetic cohomology, we compute the cohomology of all Tate twists of the rank-$1$ local system $\cL_+$ from \cref{intro:main-example}.

\begin{example} Let $G$ be the absolute Galois group of $\QQ_p$. Since $E_C\to E$ is a
pro-\'etale $G$-torsor, the Hochschild--Serre spectral sequence gives
\[
\mathrm{E}_2^{i,j}
=
\Hh^i_{\cont}\bigl(
G,\Hh^j_v(E_C,\cL_{+,C}(n))
\bigr)
\Longrightarrow
\Hh^{i+j}_v(E,\cL_+(n)).
\]
Combining with spectral sequence with \cref{eqn:computation-cohomology} and Tate's computation \cite[Ths.~1 and 2, pp.~176--177]{p-divisible} of Galois cohomology groups $\Hh^i_{\cont}(G, C(m))$,
we conclude that 
\[
\Hh^i_v\bigl(E,\cL_+(n)\bigr)
\simeq
\begin{cases}
\QQ_p,            & n=1 \text{ and }i=1,3,\\
\QQ_p^{\oplus 2}, & n=1 \text{ and }i=2,\\
0,            & \text{otherwise}.
\end{cases}
\]
The same argument computes the cohomology of $\cL_-(n)$ for every $n$. Together with the isomorphism $\rR{\cHom}_{\cD(E_v; \uQ_p)}\bigl(\cL_+(n),\uQ_p(2)[4]\bigr)
\simeq \cL_-(2-n)[4]$, these computations give a direct verification of the duality assertion in \cref{thm:intro-arithmetic-PD-proper-qp} as well.
\end{example}

\subsection{Fargues--Fontaine curve perspective}\label{section:FF-perspective}

In this subsection, we explain an alternative perspective on \cref{thm:intro-PD-proper-qp} from the point of view of the Fargues--Fontaine curve. Although the language of Banach--Colmez spaces provides the most direct and intuitive formulation of the theorem, the Fargues--Fontaine perspective reveals additional structure that restores the finitess, duality, and K\"unneth formula properties for cohomology of $\uQ_p$-local systems. 

For this, we recall that \cite{FF-curve-book} and \cite[\textsection 2]{Fargues-Scholze} defined the (analytic) Fargues--Fontaine curve as
\begin{equation}\label{eqn:intro-ff-curve}
    \FFF_C \coloneqq  Y_{C}/\varphi^{\ZZ} = \Bigl(\Spa\bigl(\AA_{\inf}(C, \O_C), \AA_{\inf}(C, \O_C)\bigr) \setminus \mathrm{V}(p\cdot [p^\flat])\Bigr)/\varphi^\ZZ,
\end{equation}
which is an (infinite type) adic space over $\Spa(\QQ_p, \ZZ_p)$, see also \cref{section:fargues-fontaine-curve} for an additional discussion.

It turns out that a particularly useful framework for the $v$-cohomology
$\rR\Gamma_v(X,\cE)$ of a $\uQ_p$-local system $\cE$ is to lift this to a
cohomology theory valued in perfect complexes on the analytic
Fargues--Fontaine curve. More precisely, to every proper rigid-analytic
space $X$ over $C$ and every $\uQ_p$-local system $\cE$ on $X$, we
functorially associate a perfect complex
\[
    \cH_X(\cE)\in\Perf_\an(\FFF_C;\O)
\]
together with a canonical equivalence
\begin{equation}\label{eqn:promote-to-the-curve}
    \rR\Gamma_v(X,\cE)
    \simeq
    \rR\Gamma_\an\bigl(\FFF_C,\cH_X(\cE)\bigr).
\end{equation}
The complex $\cH_X(\cE)$ satisfies the properties expected of
a well-behaved cohomology theory: finiteness, duality, and the K\"unneth
formula. From this perspective, the possible infinite-dimensionality and failure of the K\"unneth formula of
$\rR\Gamma_v(X,\cE)$ arise only after applying the derived global-sections
functor $\rR\Gamma_\an(\FFF_C,\blank)$, which morally come from the facts that $\FFF_C$ is not of finite type over $\QQ_p$ and $\rR\Gamma_\an(\FFF_C,\blank)$ is not symmetric monoidal. 

To relate this perspective back to \cref{thm:intro-PD-proper-qp}, we recall that \cite{LB-thesis, ALB} constructed a fully faithful functor $\Sol\colon \Perf_\an(\FFF_C; \O) \to \cD\bigl(\Spa(C, \O_C)_v; \uQ_p\bigr)$ which provides an identification $\Perf_\an(\FFF_C; \O)\xrightarrow[\sim]{\Sol} \cD^b_\BBC\bigl(\Spa(C, \O_C)_v; \uQ_p\bigr)$. Using this functor, we show a more precise version of \cref{eqn:promote-to-the-curve} as well as the results announced above. In order to capture all the properties of the complex $\cH_X(\cE)$, we recall that the Fargues--Fontaine curve carries a canonical ``point at infinity''
 $\iota_\infty\colon \Spa(C, \O_C) \hookrightarrow \FFF_C$, see \cite[Def./Prop.~II.1.22]{Fargues-Scholze}. Pullback along this closed immersion defines a functor $\iota_\infty^*\colon \Perf_\an(\FFF_C; \O) \to \Perf_\an(\Spa (C, \O_C); \O) \simeq \cD^b_{\mathrm{coh}}(C)$.

\begin{thm}[{\cref{thm:main-thm-hx}}]\label{thm:intro-ff-perspective} Let $f \colon X\to \Spa(C, \O_C)$ be a proper rigid-analytic space of dimension $d$, and let $\cE\in \Perf^{[r, r']}_v(X; \uQ_p)$ for some $r, r'\in \ZZ\cup\{-\infty, \infty\}$. Then one can functorially assign a perfect complex $\cH_X(\cE)\in \Perf_\an(\FFF_C; \O)$ which satisfies the following properties:
\begin{enumerate}[leftmargin=*,label={\upshape{(\roman*)}}]
    \item there is a functorial isomorphism $\Sol\bigl(\cH_X(\cE)\bigr) \simeq \rR f_{v, *}\cE$;
    \item\label{thm:intro-ff-perspective-2} there is a functorial isomorphism $\iota_\infty^*\cH_X(\cE)\simeq \rR\Gamma_v\bigl(X, \cE \otimes^\rL_{\uQ_p} \O_X\bigr)$;
    \item the complex $\cH_X(\cE)$ lies in $\Perf^{[r, r'+2d]}_\an(\FFF_C; \O)$;
    \item if $(X', \cE')$ is another such pair, then there is a natural K\"unneth formula isomorphism 
    \[
    \cH_X(\cE) \otimes^\rL_\O \cH_{X'}(\cE') \xr{\sim} \cH_{X\times_C X'}(\cE \boxtimes \cE').
    \]
    \item\label{thm:intro-ff-perspective-5} if $X$ smooth of equidimension $d$, there is a functorial duality isomorphism $\cH_X(\cE)^{\vee} \simeq \cH_X(\cE^{\vee}(d)[2d])$;
\end{enumerate}
\end{thm}

While finiteness and duality admit formulations in terms of Banach--Colmez spaces, the corresponding K\"unneth formula fails at the level of $v$-sheaves, even in the basic example considered above.

\begin{example}\label{example:intro-no-KF} In the setting\footnote{We use $\cL_{+,C}(1)$ rather than $\cL_{+,C}$ solely to avoid carrying Tate twists throughout the formulas. A choice of a compatible system of $p$-power roots of unity in $C$ trivializes $\uQ_p(1)$ and hence yields a noncanonical isomorphism $\cL_{+,C}\simeq\cL_{+,C}(1)$.
Thus, the same discussion applies verbatim to $\cL_{+,C}$.} of \cref{intro:main-example}, the natural morphism
\[
\rR f_{v, C, *} \bigl(\cL_{+, C}(1)\bigr) \otimes_{\uQ_p}^\rL \rR f_{v, C, *} \bigl(\cL_{+, C}(1)\bigr)  \to \rR (f\times_C f)_{v, C, *} \bigl(\cL_{+, C}(1) \boxtimes \cL_{+, C}(1)\bigr) 
\]
is far from being an isomorphism. In fact, the left-hand side does not even lie in $\cD^b_{\BBC}$ as it contains $\O_{\Spa C}[-1]\otimes^\rL_{\uQ_p} \O_{\Spa C}[-1]$ as a direct summand.

 By contrast, under the equivalence induced by the solution functor,
Equation~\cref{eqn:intro-derived-pushforward} translates into an isomorphism $\cH_{E_C}\bigl(\cL_{+,C}(1)\bigr) \simeq  \iota_{\infty,*}C[-1]\oplus\iota_{\infty,*}C[-2]$. So 
\[
\cH_{E_C}\bigl(\cL_{+, C}(1)\bigr) \otimes^\rL_{\O_{\FFF_C}} \cH_{E_C}\bigl(\cL_{+, C}(1)\bigr) \simeq \iota_{\infty, *}C[-1] \oplus (\iota_{\infty, *}C[-2])^{\oplus 3} \oplus (\iota_{\infty, *}C[-3])^{\oplus 3} \oplus \iota_{\infty, *}C[-4]
\]
is manifestly a perfect complex on $\FFF_C$. Moreover, the
K\"unneth isomorphism of \cref{thm:intro-ff-perspective}~\cref{thm:intro-ff-perspective-5} identifies this tensor product with $\cH_{E_C\times_C E_C}\bigl(\cL_{+,C}(1)\boxtimes\cL_{+,C}(1)\bigr)$. After applying the solution functor, we obtain
\[
\rR (f\times_C f)_{v, C, *} (\cL_{+, C} \boxtimes \cL_{+, C}) \simeq \O_{\Spa C}[-1] \oplus (\O_{\Spa C}[-2])^{\oplus 3} \oplus (\O_{\Spa C}[-3])^{\oplus 3} \oplus \O_{\Spa C}[-4]
\]
\end{example}

\cref{thm:intro-ff-perspective}~\cref{thm:intro-ff-perspective-2}
likewise records information that is not directly visible through the
intrinsic operations on the Banach--Colmez spaces $\rR^i f_{v,*}\cE$.
Nevertheless, it yields some non-trivial information about these higher pushforwards.
For example, using \cref{thm:intro-ff-perspective}~\cref{thm:intro-ff-perspective-2}, David Hansen
\cite[Rmk.~on p.~2]{Hansen-open-problems} showed that if \(X\) is
smooth and proper and $\bL$ is a $\uQ_p$-local system
on $X$, then
\[
\sum_n(-1)^n\dim_{\QQ_p}\rR^n f_{v,*}\bL
=
\bigl(\mathrm{rk}\,\bL\bigr)
\cdot
\left(
\sum_n(-1)^n\dim_{\QQ_p}\Hh^n_v(X,\QQ_p)
\right) = \mathrm{rk}\, \bL \cdot \chi(X, \QQ_p).
\]
Here, $\dim_{\QQ_p}\rR^n f_{v,*}\bL$ denotes the
\(\QQ_p\)-dimension, in the sense of \cite[p.~46]{Colmez}, of the
Banach--Colmez space $\rR^n f_{v,*}\bL$. This is the $p$-adic analytic analogue of the usual Euler characteristic formula $\chi(X, \bL) = \mathrm{rk}\, \bL \cdot \chi(X, \QQ_p)$ in algebraic geometry (see \cite[Cor.~2.11]{Illusie-Brauer}).

\subsection{Idea of the proof and further results}
\label{section:intro-proof}

The proof of \cref{thm:intro-PD-proper-qp} consists of three main ideas of different origins: a reduction to the smooth case, a reformulation of the problem in terms of the Fargues--Fontaine curve, and a diagrammatic approach to the resulting duality statement on the curve. The passage to the Fargues--Fontaine curve is crucial to our argument because the diagrammatic approach critically relies on the existence of a K\"unneth formula. While such a formula fails at the level of $v$-sheaves (see \cref{example:intro-no-KF}), it does hold at the level of perfect complexes on the Fargues--Fontaine curve.

\subsubsection{Reduction to the smooth case} The first step in the proof of \cref{thm:intro-PD-proper-qp} is to reduce to the case in which $X$ is smooth proper over $C$. The key technical input is the construction of the proper trace map $\Tr_f\colon \rR f_{v,*}\omega_X^{\Q_p} \longrightarrow \uQ_p.$ This construction ultimately rests on the analogous trace map with finite coefficients constructed in \cite{LRZ24}, which is one of the most technically challenging parts of that paper.

Once we have the trace map from \cite{LRZ24} at hand, its extension to $\uQ_p$-coefficients follows from a limit argument, together with the comparison between \'etale and $v$-cohomology for proper rigid-analytic spaces; see \cref{zc-cohomology-et-v}. We then combine the resulting trace map with resolution of singularities and induction on the dimension of $X$ to reduce both the finiteness and duality assertions of \cref{thm:intro-PD-proper-qp} to the smooth case. In fact, in the smooth case we prove the following stronger relative result for arbitrary smooth proper morphisms of diamonds over $\Spd(\QQ_p,\ZZ_p)$.

\begin{thm}[Relative Geometric Duality; {\cref{geometric-duality-Qp} and \cref{cor:tor-amplitude-qp-smoothproper-2}}]\label{intro:relative-geometric-duality-Qp}
Let $f \colon X' \to X$ be a smooth proper of equidimension $d$ morphism of diamonds over $\Spd(\QQ_p, \ZZ_p)$ and let $\cE\in \Perf_v(X'; \uQ_p) \cap \cD^{[r, r']}(X'_v; \uQ_p)$ for some $r,r'\in \ZZ \cup \{-\infty, \infty\}$. Then the following statements hold:
\begin{enumerate}[leftmargin=*,label={\upshape{(\roman*)}}]
    \item $\rR f_{v, *} \cE$ lies in $\cD^{[r, r'+2d]}_{\BBC}\bigl(X_v; \uQ_p\bigr)$;
    \item  there is a canonical trace map $\tr_f\colon \rR f_{v, *} \uQ_p(d)[2d] \to \uQ_p$ such that the induced duality morphism
    \[
    \PD_f \colon \rR f_{v, *} \rR {\cHom}_{\cD(X'_v;\uQ_p)}(\cE, \uQ_p(d)[2d]) \longrightarrow \rR{\cHom}_{\cD(X_v;\uQ_p)}(\rR f_{v, *} \cE,\uQ_p)
    \]
defined in \cref{defn:Qp-duality-map}\cref{defn:Qp-duality-map-3} is an isomorphism. 
\end{enumerate}
\end{thm}

We refer to \cref{defn:smooth-proper} for the precise definition of a smooth proper of equidimension $d$ morphism of diamonds over $\Spd(\QQ_p,\ZZ_p)$. For the purposes of this introduction, we only note that this class includes the diamondifications of smooth proper morphisms of equidimension $d$ between locally noetherian or sousperfectoid adic spaces over $\QQ_p$. We also refer to \cref{defn:relative-BC} for the definition of the relative Banach--Colmez spaces $\cD^{(b)}_{\BBC}(X_v;\uQ_p)$. Specializing \cref{intro:relative-geometric-duality-Qp} to $X=\Spa(C,\O_C)$ gives the smooth case required in the reduction above.

\subsubsection{Reduction to duality on the Fargues--Fontaine curve}
\label{section:reduction-FF} 
The preceding subsection reduces the proof of
\cref{thm:intro-PD-proper-qp} to the relative smooth proper
case. We now carry out the second step of the argument by reformulating
this as a relative statement on the Fargues--Fontaine curve. 
The solution functor from \cref{section:FF-perspective} gives an
identification $\Perf_\an(\FFF_C; \O) \xr{\sim} \cD^b_{\BBC}\bigl(\Spa(C, \O_C)_v; \uQ_p\bigr)$. 
It therefore suffices to prove duality on the curve, where we can use the
tensor formalism required for the K\"{u}nneth formula. To do so, we introduce
``analytic perfect complexes on the relative Fargues--Fontaine curve'',
equipped with a pushforward functor that intertwines with the derived $v$-pushforward via the solution functor.

A technical hurdle arises when one attempts to carry this out:
If $X$ is an adic space that is not perfectoid, the relative Fargues--Fontaine curve 
$\FFF_X\coloneqq X^\diam \times \Spd(\QQ_p, \ZZ_p)/\varphi^\ZZ\times \id$ is only a diamond and need not carry a natural adic-space structure.
The usual definition of \emph{analytic} perfect complexes therefore does not apply directly to $\FFF_X$, 
while perfect complexes in the $v$-topology on $\FFF_X$ do not give the correct answer, see \cref{warning:analytic-vs-v}. 

To resolve this issue, motivated by the fact that $Y_{C,[1,p]}\subset Y_C$ is a fundamental domain for $\FFF_C$, we instead define 
\begin{equation}
\label{eqn:perfect-complexes-on-the-curve} 
\Perf_{\an}(\FFF_X;\O) \coloneqq \operatorname*{Eq}\left( \Perf(X_v;\BB_{[1,p]})
\rightrightarrows \Perf(X_v;\BB_{[1,1]}) \right). 
\end{equation}
Here $\BB_{[1,p]}$ and $\BB_{[1,1]}$ are the period sheaves geometrizing Berger's period rings from  \cite{Berger02}; see \cref{defn:period-sheaves}. The two arrows in 
\eqref{eqn:perfect-complexes-on-the-curve} encode the restriction and Frobenius maps.

To compare this definition with the usual category of analytic perfect complexes when $X=S$ is affinoid perfectoid, one needs, in particular, the equivalence $\Perf_v(S;\BB_I) \simeq \Perf\bigl(\BB_I(S)\bigr)$ for every closed interval $I\subset(0,\infty)$ with rational endpoints. This comparison follows from the following general descent theorem.

\begin{thm}[{see~\ref{thm:v-descent-finitary-etale},~\ref{lemma:example-z-mod-n},~\ref{lemma:example-o-plus-mod-varpi},~\ref{lemma:example-AI-mod-p},~\ref{cor:v-descent-integrally},~\ref{cor:v-descent-A_I},~\ref{thm:v-descent-rationally},~\ref{thm:v-descent-rationally-BBI},~\ref{cor:integral-descend-finite-tor-amplitude-ell-adically},~and~\ref{thm:v-descent-rationally-Qell}}]\label{thm:intro-v-descent} 
Denote by $\Perfd^{\aff}_{/\QQ_p}$ the category of affinoid perfectoid spaces over $\Spa(\QQ_p, \ZZ_p)$, let $I\subset (0, \infty)$ be a closed interval with rational endpoints, let $n$ be an integer, let $\ell$ be a prime number (possibly equal to $p$), and let $r,r'\in \ZZ\cup \{-\infty, \infty\}$. Then the assignments
\[
\begin{gathered}
\Perf^{[r, r']}_\et(\blank; \ud{\ZZ/n}),\, \Perf^{[r, r']}_\proet(\blank; \ud{\Z}_\ell),\, \Perf^{[r, r']}_\proet(\blank; \ud{\Q}_\ell),\, \Perf^{[r, r']}_\et(\blank; \O^+/p), \\
 \Perf^{[r, r']}_\et(\blank; \O^+),\, \Perf^{[r, r']}(\O(\blank)),\, \Perf^{[r, r']}_\et(\blank; \AA_I),\, \Perf^{[r, r']}(\BB_I(\blank))
\end{gathered}
\]
regarded as functors $\Perfd^{\aff,\op}_{/\QQ_p} \longrightarrow \Cat_\infty$
are hypercomplete $v$-sheaves of $\infty$-categories. 
\end{thm}

    The fact that $\Perf(\BB_I(\blank))$ is a $v$-sheaf was previously established in \cite{ALB} by different methods. Likewise, {\it loc.~cit.} showed that $\Perf(\O(\blank))$ is a $v$-sheaf when restricted to affinoid perfectoid spaces admitting a morphism to a totally disconnected perfectoid space. Restricting the class of affinoid perfectoid spaces becomes necessary because very delicate descent results from \cite{Mann-thesis} are used.

Thus, for the functor $\Perf(\O(\blank))$,
\cref{thm:intro-v-descent} strengthens the result of \cite{ALB} by removing
this restriction. In addition, our proof is self-contained and avoids the descent results of \cite{Mann-thesis}.
Its key idea is to adapt to the
derived setting the methods developed in
\cite{Z-thesis,Heuer-G-torsor} for proving $v$-descent for \'etale vector
bundles over $\O$, $\O^+$, and $\O^+/p$. To the best of our knowledge, the
corresponding statements for the remaining coefficient systems appearing in
\cref{thm:intro-v-descent} have not previously appeared in the literature.

The descent theorem shows that \eqref{eqn:perfect-complexes-on-the-curve} agrees with the usual category of analytic perfect complexes for affinoid perfectoid bases and provides a well-behaved definition for arbitrary diamonds. After developing the corresponding pushforward formalism, we reduce \cref{intro:relative-geometric-duality-Qp} to the following duality theorem on the relative Fargues--Fontaine curve.

\begin{thm}[Duality on the curve; {\cref{thm:geometric-duality-curve}}]\label{intro:relative-geometric-duality-FF}
Let $f \colon X' \to X$ be a smooth proper of equidimension $d$ morphism of diamonds over $\Spd(\QQ_p, \ZZ_p)$ and let $\cE\in \Perf_\an(\FFF_{X'}; \O)$. Then the following statements hold:
\begin{enumerate}[leftmargin=*,label={\upshape{(\roman*)}}]
    \item $\rR f_{\FFF, *} \cE$ lies in $\Perf_\an(\FFF_X; \O)$;
    \item  there is a canonical trace map $\tr^\FFF_f\colon \rR f_{\FFF, *} \O(d)[2d] \to \O$ such that the induced duality morphism
    \[
    \PD_f \colon \rR f_{\FFF, *} \rR {\cHom}_{\FFF_{X'}}(\cE, \O(d)[2d]) \longrightarrow \rR{\cHom}_{\FFF_{X}}(\rR f_{\FFF, *} \cE, \O)
    \]
defined in \Cref{defn:FF-trace-maps}\Cref{defn:FF-trace-maps-3} is an isomorphism. 
\end{enumerate}
\end{thm}

\subsubsection{Proof of the duality on the Fargues--Fontaine curve. The diagrammatic approach.}

We now explain the main argument of the paper. Unwinding the equalizer description in \cref{eqn:perfect-complexes-on-the-curve}, we reduce duality on the Fargues--Fontaine curve to a relative duality theorem for period sheaves; only the case of $\cA=\BB_{[1,p]}$ is used in the proof of \cref{intro:relative-geometric-duality-FF}.

\begin{thm}[Relative duality for period sheaves; \cref{cor:PD-for-proper-smooth-pushforward}]\label{intro:relative-geometric-duality-period-sheaves}
Let $f \colon X' \to X$ be a smooth proper of equidimension $d$ morphism of diamonds over $\Spd(\QQ_p, \ZZ_p)$, let $I\subset (0, \infty)$ be a closed interval with rational endpoints, let $\cA\in \{\O, \BB_I\}$, and let $\cE\in \Perf_v(X'; \cA) \cap \cD^{[r, r']}(X'_v; \cA)$ for some $r,r'\in \ZZ \cup \{-\infty, \infty\}$. Then the following statements hold:
\begin{enumerate}[leftmargin=*,label={\upshape{(\roman*)}}]
    \item $\rR f_{v, *} \cE$ lies in $\Perf_v(X; \cA) \cap \cD^{[r, r'+2d]}(X_v; \cA)$;
    \item  there is a canonical trace map $\tr^{\cA}_f\colon \rR f_{v, *} \cA(d)[2d] \to \cA$ such that the induced duality morphism
    \[
    \PD_f \colon \rR f_{v, *} \rR \cHom_{\cD(X'_v;\cA)}(\cE, \cA(d)[2d]) \longrightarrow
    \rR \cHom_{\cD(X_v;\cA)}(\rR f_{v, *} \cE,\cA)
    \]
defined in \Cref{construction:evaluation-coevaluation-poincare-duality-non-solid}\Cref{duality-sheafy}
is an isomorphism. 
\end{enumerate}
\end{thm}

We prove \cref{intro:relative-geometric-duality-period-sheaves} using the diagrammatic method developed in \cite[Section 6]{LRZ24}. Rather than understanding the duality morphism $\PD_f$ directly, this approach shows that $\rR f_{v, *}\cE$ is dualizable (in the categorical sense) in $\cD(X_v;\cA)$, with dual given by $\rR f_{v, *}\cE^\vee(d)[2d]$. Its implementation requires the following three ingredients, together with suitable compatibility conditions:
\begin{enumerate}
\item the trace map $\tr_f^{\cA} \colon \rR f_{v, *} \cA(d)[2d] \to \cA$ for the smooth proper $f\colon X' \to X$;
\item the cycle class map $\cl_{\Delta} \colon \rR \Delta_{v, *}\cA \to \cA(d)[2d]$ for the diagonal $\Delta\colon X' \hookrightarrow X'\times_X X'$;
\item the K\"{u}nneth formula isomorphism $\rR f_{v, *}\cE \otimes^\rL_{\cA} \rR f_{v, *} \cE \xr{\sim} \rR (f\times_X f)_{v, *}\bigl(\cE \boxtimes^\rL \cE\bigr)$.
\end{enumerate}
Of these three ingredients, the K\"unneth formula is the most delicate. We do not know how to establish it directly; instead, we show the following version using the framework of solid mathematics in the sense of Clausen--Scholze \cite{Condensed}. 
\begin{theorem}[{Solid K\"unneth formula; \cref{Kunneth formula: main thm}}]
\label{intro:kunneth}
Let $S$ be a strictly totally disconnected perfectoid space over $\Spa(\Q_p, \Z_p)$, let $f \colon X \to S$ and $f' \colon X' \to S$ be 
qcqs smooth maps, and let $\cA\in \{\BB_I, \O\}$. Let $\cE \in \Perf_v(X, \cA)$ and 
$\cE' \in \Perf_v(X', \cA)$.
Then the solid K\"{u}nneth map 
\[
\KM^\Box\colon \rR\underline{\Gamma}_v(X, \cE) \otimes^{\rL \Box}_{\underline{\Gamma}_v(S, \cA)}
\rR\underline{\Gamma}_v(X', \cE') \rightarrow
\rR\underline{\Gamma}_v(X \times_{S} X', \cE \boxtimes^{\rL} \cE').
\]
is an isomorphism. 
\end{theorem}

We refer to \cref{condensed structure definition} for the definition of condensed derived global sections.
The crucial feature of \cref{intro:kunneth} is that it requires no properness assumption, which allows us to argue locally.
Using toric charts, we reduce the solid K\"unneth theorem to the corresponding statement for the 
continuous group cohomology of $\Gamma\simeq\ZZ_p^d$, which can be seen via Koszul complexes.

Finally, we use this K\"unneth formula to establish a solid version of \cref{intro:relative-geometric-duality-period-sheaves}.  

\begin{thm}[{\cref{thm:PD-absolute-solid}}]
\label{thm:intro-PD-absolute-solid} 
Let $S$ be a strictly totally disconnected perfectoid space over $\Spa(\QQ_p, \ZZ_p)$, 
let $f\colon X \to S$ be a smooth proper morphism of equidimension $d$, 
let $I\subset (0, \infty)$ be a closed interval with rational endpoints, let $\cA\in \{\BB_I, \cO\}$, and let $\cE\in \Perf_v(X; \cA)$.
Then the condensed global sections $\rR\uGamma_v(X, \cE)$ is a dualizable object in $\cD_\sol(\ud{\cA})$,
with the dual given by $\rR\uGamma_v(X, \cE^{\vee}(d)[2d])$.
\end{thm}

\cref{thm:intro-PD-absolute-solid} is the core result of this paper; we now outline its proof, following the diagrammatic approach introduced above.
We first construct the required trace and cycle class maps by combining
the analogous constructions from \cite{LRZ24} with an appropriate form
of the primitive comparison theorem; see
\cref{Relative primitive comparison-2}. 
Together with the solid K\"{u}nneth formula, 
we can produce evaluation and coevaluation maps between $\rR\uGamma_v(X, \cE)$ and $\rR\uGamma_v(X, \cE^{\vee}(d)[2d])$;
see \Cref{construction:evaluation-coevaluation-poincare-duality}.
To check their dualizability, it therefore remains to check the induced triangle identities.
We address one of them, the other being verified similarly.
Consider the commutative diagram below, in which all products of $X$ are taken over $S$:
\[ \begin{tikzcd}
&& X \arrow[ld, "\Delta"'] \arrow[rd, "\Delta"] && \\
& X^2 \arrow[ld, "\pi_2"'] \arrow[rd, "\Delta \times \id"']
&& X^2 \arrow[ld, "\id \times \Delta"] \arrow[rd, "\pi_1"] & \\
X && X^3 && X
\end{tikzcd} \]
It induces the following diagram, in which $(\blank)^*$ denotes pullback in cohomology and $(\blank)_*$ is induced by cycle class maps:
\[
\begin{tikzcd}[column sep=large,center picture,scale cd=.85]
\rR\underline{\Gamma}_v(X, \cE) \arrow[rd, "\id"'] \arrow[r, "\pi_2^*"] & 
\rR\underline{\Gamma}_v(X^2, \cA \boxtimes \cE) \arrow[r, "\mathrm{coev} \boxtimes \id"] \arrow[d, "\Delta^*"]
& \rR\underline{\Gamma}_v(X^2, (\cE \otimes \cE^{\vee}) \boxtimes \cE) 
\arrow[r, "(\Delta \times \id)_*"] \arrow[d, "\Delta^*"] &
\rR\underline{\Gamma}_v(X^3, (\cE \boxtimes \cE^{\vee} \boxtimes \cE) \otimes \cA(d)[2d])
\arrow[d, "(\id \times \Delta)^*"] \\
& \rR\underline{\Gamma}_v(X, \cE) \arrow[rd, "\id" below] \arrow[r, "\mathrm{coev} \otimes \id"] &
\rR\underline{\Gamma}_v(X, \cE \otimes \cE^{\vee} \otimes \cE) \arrow[d, "\id \otimes \mathrm{ev}"] 
\arrow[r, "\Delta_*"]&
\rR\underline{\Gamma}_v(X^2, (\cE \boxtimes (\cE^{\vee} \otimes \cE)) \otimes \cA(d)[2d])
\arrow[d, "(\id \boxtimes \mathrm{ev}) \otimes \id"] \\
& & \rR\underline{\Gamma}_v(X, \cE) \arrow[rd, "\id"'] 
\arrow[r, "\Delta_*"]& 
\rR\underline{\Gamma}_v(X^2, (\cE \boxtimes \cA) \otimes \cA(d)[2d])
\arrow[d, "\tr_{\pi_1}"] \\
& & & \rR\underline{\Gamma}_v(X, \cE).
\end{tikzcd} \]
Up to some identifications, which we make fully explicit in \cref{gigantic-diagram}, the claim that 
\[ (\id \otimes \text{ evaluation}) \circ (\text{coevaluation} \otimes \id) = \id \]
reduces to showing the composition along top row and rightmost column is the identity.
To do so, it suffices to prove that the above diagram is commutative.
This is mostly easy, except for the top right square, which follows from the fact that cycle classes commute with transversal pullbacks (see \cite[Lem.~3.3.4]{LRZ24}) and the bottom triangle, which follows from the fact that smooth traces of cycle classes of a section are the identity (see \cite[Lem.~6.4.5]{LRZ24}).

Once \cref{thm:intro-PD-absolute-solid} is proven, a theorem of Andreychev \cite{Andreychev} identifies the dualizable objects in the relevant solid derived categories with ordinary perfect complexes. Combining this result with the $v$-descent statement of \cref{thm:intro-v-descent}, we pass from \cref{thm:intro-PD-absolute-solid} to \cref{intro:relative-geometric-duality-period-sheaves}. 

As a final note, we also mention that we have a version of
\cref{intro:relative-geometric-duality-period-sheaves}
when $X=\Spa(K, \O_K)$ 
is the adic spectrum of a nonarchimedean field and $X'$ is a general proper rigid-analytic space. 

\begin{thm}[Absolute duality for period sheaves II; {\cref{singular-PD}}]\label{intro:relative-geometric-duality-period-sheaves-point}
Let $K$ be a nonarchimedean field of mixed characteristic $(0, p)$, let $f \colon X \to \Spa(K, \O_K)$ be a proper rigid-analytic space over $K$ of dimension $d$.  Let $I\subset (0, \infty)$ be a closed interval with rational endpoints, let $\cA\in \{\O, \BB_I\}$, and let $\cE\in \Perf_v(X; \cA) \cap \cD^{[r, r']}(X'_v; \cA)$ for some $r,r'\in \ZZ \cup \{-\infty, \infty\}$. Then the following statements hold:
\begin{enumerate}[leftmargin=*,label={\upshape{(\roman*)}}]
    \item $\rR f_{v, *} \cE$ lies in $\Perf_v(\Spa(K, \O_K); \cA) \cap \cD^{[r, r'+2d]}(\Spa(K, \O_K)_v; \cA)$;
    \item  there is a canonical trace map $\tr^{\cA}_f\colon \rR f_{v, *} \omega_X^\cA \to \cA$ such that the induced duality morphism
    \[
    \PD_f \colon \rR f_{v, *} \rR {\cHom}_{\cD(X'_v;\cA)}(\cE, \omega_X^\cA) 
    \longrightarrow
    \rR {\cHom}_{\cD(X_v;\cA)}(\rR f_{v, *} \cE, \omega_X^\cA)
    \]
defined in \cref{BI-proper-trace}\cref{BI-proper-PD-morphism} is an isomorphism. 
\end{enumerate}
\end{thm}

\subsection*{Organization of the paper} The thrill-seeking reader can directly get to the key arguments of the paper by focusing their attention on \cref{dualizability}, Sections~\ref{section:pushforward}--\ref{section:cohomological-duality-geometric}, and \cref{section:absolute-duality}. The rest of the paper is mostly devoted to preparatory background results, which we break down in more detail below. 

Sections \ref{section:period-rings}, \ref{section:period-sheaves}, and \ref{trace-cycle} contain the preliminary material that we need in our proofs.
\Cref{section:period-rings} and \cref{section:period-sheaves} provide technical background on period rings and period sheaves, while \cref{trace-cycle} extends the construction of trace and cycle class maps from \cite{LRZ24} to various period sheaves on sousperfectoid spaces and diamonds. 
As we could often not find suitable references in the literature, we have tried to be very detailed, but advise even the more honest reader to skip these sections on a first reading and refer back to them as necessary.

The key argument for Poincar\'e duality for perfect complexes over period sheaves with respect to general smooth proper morphisms following the diagrammatic approach from \cite{LRZ24} is explained in \cref{dualizability}.
The resulting Poincar\'e dualities for analytic perfect complexes on relative Fargues--Fontaine curves and for pro-\'etale $\underline{\QQ}_p$-local systems are then discussed in \cref{geom-arith-duality}.
This section also contains a detailed discussion of the relevant Fargues--Fontaine curve preliminaries, again recommended to be skipped on a first reading.
Finally, we prove geometric duality (\cref{thm:intro-PD-proper-qp}) and arithmetic duality (\cref{thm:intro-arithmetic-PD-proper-qp}) for pro-\'etale $\underline{\QQ}_p$-local systems on proper rigid-analytic spaces in \cref{section:absolute-duality}.

Some technical background on limits of $\infty$-categories and sheaves of $\infty$-categories which we found difficult to locate in the literature is collected in \cref{section:appendix-limits} and \cref{section:valued-sheaves}.
Finally, \cref{section:v-sheaves} is dedicated to proving many optimal $v$-descent results for perfect complexes over various coefficient rings as in \cref{thm:intro-v-descent}, which might be interesting in their own right. 

\subsection*{Acknowledgements}
We wholeheartedly thank Pierre Colmez: his talk \cite{OB} on Banach--Colmez enhanced dualities in the context of Stein spaces and constant coefficients (based on joint work with Sally Gilles and Wies{\l}awa Nizio{\l}) was an important source of
inspiration for the correct finiteness and duality statements in our work. 
During the process of writing this paper, we have benefited from discussions with the following
mathematicians: Johannes Ansch\"utz, Bhargav Bhatt, Guido Bosco, Pierre Colmez, Gabriel Dospinescu, David Hansen, Kiran Kedlaya, Zhenghui Li, Ruochuan Liu,
Sasha Petrov, and Sebastian Wolf.

We are grateful to the following institutions for their financial support and for providing excellent
working conditions while we completed parts of this project: the School of Mathematics of the Institute for
Advanced Study (S.L., W.N., E.R., and B.Z.), the Institut des Hautes \'Etudes Scientifiques (S.L.\ and E.R.), the Institute for Advanced Study in Mathematics and Zhejiang University (S.L.), the University of Copenhagen (E.R.), Princeton University (B.Z.), and the University of Maryland (B.Z.).

S.L.\ was supported by the National Key R $\&$ D Program of China No.~2023YFA1009701 and the National Natural Science Foundation of China (No.~12288201).
W.N.\ research was partially supported by the Simons Foundation  through the Simons Collaboration on Perfection in Algebra, Geometry, and Topology.
E.R.\ was partially supported by the National Science Foundation under Grant No.\ DMS-1926686.
Part of this research was conducted during the authors' participation in the IAS special year on $p$-adic geometry.

\subsection*{AI-usage disclosure}
The authors used large language models during the preparation of this manuscript for
finding references, mathematical proofreading, checking internal consistency and cross-references.

\subsection*{Notation}\label{section:notation}

In this paper, a topological field $K$ is called \emph{nonarchimedean} if its topology
is induced by a valuation of rank $1$ on $K$ and if $K$ is complete with respect to this topology; beware that while the completeness assumption is somewhat standard, it is not imposed in \cite[Def.~1.1.3]{Huber-etale}. 

A rigid-analytic space over a nonarchimedean field $K$ is always understood to be an adic space which
is locally of finite type over $\Spa(K, \O_K)$; by \cite[Prop.~4.5]{Huber-generalization}, the resulting category of quasi-separated
rigid-analytic $K$-spaces is equivalent to the category of quasi-separated rigid-analytic $K$-spaces in the classical
sense as in, say, \cite[Def.~9.3.1/4]{BGR}. An analytic adic space is locally noetherian if every $x\in X$ is contained in an open affinoid subspace $U \subset X$ for which $\O_X(U)$ is a strongly noetherian Tate ring. For a general adic space $X$, an \emph{open affinoid subspace} $U\subset X$ is an open subspace which is isomorphic to $\Spa(A, A^+)$ for a sheafy Huber pair $(A, A^+)$. 

Following \cite{diamonds}, we denote by $\Perfd$ the category of perfectoid spaces and by $\Perf$ the category of perfectoid spaces of characteristic $p$. For a diamond $X$ (resp.~a perfectoid space $X$) and $\tau\in \{\qproet, v\}$ (resp.~$\tau\in \{\proet, v\}$), we denote by $\lambda_{\tau} \colon X_\tau \to X_\et$ and $\mu\colon X_v \to X_\qproet$ (resp.~and $\nu\colon X_\et \to X_\an$) the natural morphisms of sites. We use the diamondification functor $(\blank)^\lozenge$ from \cite[Def.~15.5]{diamonds}.

For a ring $A$ and a finite set of elements $f_1, \dots, f_n, s\in A$, we denote by $A\bigl[\frac{f_1}{s}, \dots, \frac{f_n}{s}\bigr] \subset A\bigl[\frac{1}{s}\bigr]$ the $A$-subalgebra of $A\bigl[\frac{1}{s}\bigr]$ generated by $\frac{f_1}{s},\dots, \frac{f_n}{s}$. For a ring $A$, a finitely generated ideal $I\subset A$, and a finite set of elements $f_1, \dots, f_n, s\in A$, we denote by $A\bigl\langle \frac{f_1}{s}, \dots, \frac{f_n}{s}\bigl\rangle \coloneqq A\bigl[\frac{f_1}{s},  \dots, \frac{f_n}{s}\bigr]^{\wedge}_I$ the $I$-adic completion of $A\bigl[\frac{f_1}{s}, \frac{f_2}{s}, \dots, \frac{f_n}{s}\bigr]$.

For a Grothendieck abelian category $\cA$, we denote by $\cD(\cA)$ its associated derived $\infty$-category. Its homotopy category is denoted by $D(\cA)$ and coincides with the usual triangulated derived category of $\cA$.

For a complete Huber ring $A$, we denote by $\cD_\sol(\ud{A})$ the derived $\infty$-category of solid $\ud{A}$-modules in the sense of Clausen--Scholze \cite{Condensed}. Explicitly, this category is given by $\cD\bigl((\ud{A}, \ud{\ZZ})_\sol\bigr)$  from \cite[Th.~3.28]{Andreychev}. 

Given a ring $R$, we say that an object $M\in \cD(R)$ is \emph{perfect} if there is an isomorphism $M\simeq K^\bullet$ where $K^\bullet$ is a finite complex of finite projective $R$-modules. We denote by $\Perf(R)$ the $\infty$-subcategory of $\cD(R)$ of perfect complexes. For $r,r'\in \ZZ\cup \{-\infty, \infty\}$, we say that a perfect complex $M\in \Perf(R)$ has tor-amplitude in $[r, r']$ if $M\otimes^\rL_R N$ is concentrated in cohomological degrees $[r, r']$ for any $R$-module $N$. We denote by $\Perf^{[r, r']}(R)$ the $\infty$-category of perfect complexes with tor-amplitude in $[r, r']$. 

Given a ringed site $(X, \O)$, we say that an object $\F\in \cD(X; \O)$ is \emph{strictly perfect} if there is an isomorphism $\F\simeq \cal{E}^\bullet$ where $\cal{E}^\bullet$ is a finite complex of direct summands of finite free $\O$-modules. We also denote by $\Perf(X; \O)$ the $\infty$-category of perfect complexes, i.e., the full $\infty$-subcategory $\Perf(X; \O) \subset \cD(X; \O)$ of objects which are locally isomorphic to a strictly perfect complex. 

For a ringed site $(X, \O)$, we denote by $\widetilde{(\blank)} \colon \cD(\O(X)) \to \cD(X; \O)$ the derived pullback functor along the canonical morphism $(X, \O) \to (*, \O(X))$. It restricts to the functor $\widetilde{(\blank)} \colon \Perf(\O(X)) \to \Perf(X; \O)$. 

\newpage

\addtocontents{toc}{\protect\setcounter{tocdepth}{2}}

\section{Period Rings}
\label{section:period-rings}

In this section, we define and study the period rings $\BB_I(R, R^+)$ and $\AA_I(R, R^+)$ for any perfectoid pair $(R, R^+)$ over $(\Q_p, \Z_p)$ and any closed interval $I=[a,b] \subset (0, \infty)$ with rational endpoints. We decided to include this material in the paper because these period rings are often defined under the additional assumption that $(R, R^+)$ admits a map from $(\CC_p, \O_{\CC_p})$, which is inadequate for our purposes.

\subsection{Period rings via the Fargues--Fontaine disk}

Throughout this subsection, we fix an affinoid perfectoid space $S=\Spa(R, R^+)$ over $\Spa(\QQ_p, \Z_p)$.

The main goal of this subsection is to define the period rings $\BB_I(S)$ and $\AA_I(S)$ for any closed interval $I\subset (0, \infty)$ with rational endpoints.
For this, we recall that the perfectoid pair $(R, R^+)$ admits a tilt $(R^\flat, R^{\flat, +})$ which is a perfectoid pair of characteristic $p$. It comes equipped with a continuous multiplicative map 
\[
(\blank)^\sharp \colon R^\flat \simeq \lim_{x\mapsto x^p} R \xr{(\dots, x_i, \dots, x_1, x_0) \mapsto x_0} R,
\]
where $x_i^p=x_{i-1}\in R$, see \cite[Prop.~5.17]{Scholze-perfectoid} or \cite[p.~16]{diamonds} for more details. 

\begin{notation}
\label{notation:p-flat} 
We denote by $p^\flat \in R^{\flat, +}$ any element such that $(p^\flat)^\sharp = p\cdot u$ for some unit $u\in (R^+)^\times$. 
\end{notation}

We note that \cite[Lem.~3.9]{BMS1} implies that $p^\flat$ exists for any perfectoid pair $(R, R^+)$. Any choice of $p^\flat \in R^{\flat, +}$ is a pseudo-uniformizer in $R^\flat$ (see the proof of \cite[Lem.~3.10]{diamonds}). 

\begin{lemma}\label{lemma:independence-of-flat} Let $(R, R^+)$ be an affinoid perfectoid pair over $(\QQ_p, \Z_p)$ and let $p^\flat$ and $\widetilde{p}^\flat \in R^{\flat, +}$ be two elements as in \cref{notation:p-flat}. Then there is a unit $v\in (R^{\flat, +})^\times$ such that $p^\flat = \widetilde{p}^\flat\cdot v$. 
\end{lemma}
\begin{proof}
    Using \cite[Lem.~1.6(iv)]{Huber-generalization}, we see that it suffices to show that $\abs{p^\flat(x)} = \abs{\widetilde{p}^\flat(x)}$ for any $x\in \Spa(R^\flat, R^{\flat, +})$. Consider the canonical homeomorphism $(\blank)^\flat \colon \abs{\Spa(R, R^{+})} \xr{\sim} \abs{\Spa(R^\flat, R^{\flat, +})}$ (see \cite[Th.~6.3]{Scholze-perfectoid} or \cite[Th.~3.12]{diamonds}) and set $x^\sharp\in \abs{\Spa(R, R^{+})}$ to be the unique point such that $(x^\sharp)^\flat =x$.
    Let $u,\widetilde{u} \in (R^+)^\times$ such that $(p^\flat)^\sharp = p\cdot u$ and $(\widetilde{p}^\flat)^\sharp = p\cdot \widetilde{u}$, respectively.
    Then the desired claim follows from the following sequence of equalities: 
    \[
    \abs{p^\flat(x)} = \abs{(pu)(x^\sharp)} = \abs{p(x^\sharp)} = \abs{(p\widetilde{u})(x^\sharp)} = \abs{\widetilde{p}^\flat(x)}. \qedhere
    \]
\end{proof}
We recall that $S=\Spa(R, R^+)$ denotes an affinoid perfectoid space over $\Spa(\QQ_p, \Z_p)$ and that $\AA_{\inf}(S)$ denotes the ring $W(R^{\flat, +})$.
\begin{definition}\label{defn:Y-Fargues--Fontaine}
    The \emph{relative Fargues--Fontaine punctured open disk} 
    is the adic space  
    \[ 
        Y_{S} \colonequals \Spa\big(\AA_{\inf}(S), \AA_{\inf}(S)\big) \smallsetminus \rm{V}(p\cdot[p^\flat]), 
    \]
    where $\AA_{\inf}(S)=W(R^{\flat, +})$ is equipped with the $(p,[p^\flat])$-adic topology;
    this open subset does not depend on the choice of $p^\flat$ by virtue of \cref{lemma:independence-of-flat}. 
\end{definition}
Note that \cite[Prop.~II.1.1]{Fargues-Scholze} indeed implies that $Y_S$ is an analytic \emph{adic} space instead of simply a pre-adic space. However, it is not clear whether  $\Spa\big(\AA_{\inf}(S), \AA_{\inf}(S)\big)$ is an adic space (i.e., whether it is sheafy): we will never use it in this paper. 

\begin{rmk}\label{rmk:inconsistency-with-FS} We point out that our definition of $Y_S$ is slightly different from the one used in \cite[Ch.\,II]{Fargues-Scholze}. Namely, they define $Y_S$ only when $S$ is a perfectoid space of characteristic $p$. However, $Y_S$ in the sense of \cref{defn:Y-Fargues--Fontaine} is equivalent to $Y_{S^\flat}$ in the sense of \cite[Ch.\,II.1.1]{Fargues-Scholze}. 
\end{rmk}

The space $Y_S$ is never quasi-compact, making it quite difficult to study it directly. For this purpose, we introduce the following rational open subsets of $Y_S$ whose union covers $Y_S$. 

\begin{definition}\label{defn:affinoid-Fargues--Fontaine}
    Let $I\colonequals [a,b] \subset (0,\infty)$ be an interval with rational endpoints. Set $c \colonequals \min \{ n \in \ZZ_{>0} \suchthat n \cdot a \in \ZZ \text{ and } n \cdot b \in \ZZ \}$. The  \emph{Fargues--Fontaine annulus of radius $I$} is the affinoid adic space  
    \[ 
        Y_{S,I} \colonequals \Big\{x\in \Spa\big(\AA_{\inf}(S), \AA_{\inf}(S)\big) \mid \abs{p^{cb}(x)} \le \abs{[(p^{\flat})^c](x)} \le \abs{p^{ca}(x)} \neq 0 \Big\} \subset \Spa\big(\AA_{\inf}(S), \AA_{\inf}(S)\big).
    \]
    This rational open subset does not depend on the choice of $p^\flat$ by virtue of \cref{lemma:independence-of-flat}. 
\end{definition}

\begin{rmk} By definition, we see that $Y_S = \bigcup_I Y_{S, I}$, where the union is taken over closed intervals $I\subset (0, \infty)$ with rational endpoints. Furthermore, each $Y_{S, I}$ is again an analytic \emph{adic} space.
\end{rmk}

\begin{defn}\label{defn:period-rings} We define \emph{period rings} $\BB_I(S)$ and $\AA_I(S)$ as the following topological rings:
    \[ \AA_I(S) \colonequals \Hh^0_\an(Y_{S,I},\cO^+) \quad \text{and} \quad \BB_I(S) \colonequals \Hh^0_\an(Y_{S,I},\cO). \]
\end{defn}

\begin{rmk}
\label{rmk:BI-AI-Tate-Huber-pair}
Since each $Y_{S, I}$ is an affinoid adic space over $\Spa(\QQ_p, \ZZ_p)$, 
the pair $\big(\BB_I(S), \AA_I(S)\big)$ is a complete Tate--Huber pair over $(\QQ_p, \ZZ_p)$ with pseudo-uniformizer $p\in \AA_I(S)$.
Thus, $\BB_I(S)$ admits a Banach ring structure over $\QQ_p$ (see \cite[Rmk.~2.4.4]{Kedlaya-Liu-1}). 
\end{rmk}

We finish this subsection by unraveling the definitions of $\BB_I(S)$ and $\AA_I(S)$ in more explicit terms.
Recall that, for a ring $A$ and elements $f_1, \dots, f_n, g_1, \dots, g_n\in A$, we denote by $A\big[\frac{f_1}{g_1}, \dots, \frac{f_n}{g_n}\big] \subset A[\frac{1}{g_1\dots g_n}]$ the $A$-subalgebra generated by $\frac{f_1}{g_1}, \dots, \frac{f_n}{g_n}$.

\begin{rmk}\label{rmk:explicit-construction} Let $I=[a,b]$ and $c\in \ZZ_{>0}$ be as in \cref{defn:affinoid-Fargues--Fontaine}. Then \cite[Lem.\ and def.\ on p.~516]{Huber-generalization} implies that $\AA_I(S)$ is equal to the completed integral closure of $\AA_{\inf}(S)\bigl[\tfrac{p^{cb}}{[(p^{\flat})^c]},\tfrac{[(p^\flat)^c]}{p^{ca}}\bigr]$ inside $\AA_{\inf}(S)\bigl[\tfrac{1}{[p^{\flat}]p}\bigr]$
and $\BB_I(S)=\AA_I(S)[\frac{1}{p}]$. 
We warn the reader that it is a priori not clear whether the topology on $\AA_I(S)$ is the $p$-adic one.
We return to this point in \cref{cor:almost-mathematics-AI}.
\end{rmk}
For future reference, we mention the following compatibilities between the period rings $\BB_I(S)$ for various $I$: 
\begin{construction}\label{construction:compatibility-maps} Let $S$ be an affinoid perfectoid space over $\Spa(\QQ_p, \ZZ_p)$. 
 \begin{enumerate}
   \item\label{construction:compatibility-maps-1} (\emph{restriction maps}) Let $J\subset I \subset (0, \infty)$ be an inclusion of closed intervals with rational endpoints. Then the natural rational open embedding $Y_{S, J}\subset Y_{S, I}$ induces a continuous morphism $\BB_I(S) \to \BB_J(S)$.
   \item\label{construction:compatibility-maps-2} (\emph{Frobenii maps}) Let $I=[a, b]\subset (0, \infty)$ be a closed interval with rational endpoints. Set $I/p\coloneqq [a/p, b/p] \subset (0, \infty)$. Then the Frobenius isomorphism $Y_{S} \xr{\sim} Y_S$ restricts to an isomorphism $Y_{S, I/p} \xr{\sim} Y_{S, I}$. Therefore, it induces a topological isomorphism $\varphi \colon \BB_{I}(S) \xr{\sim} \BB_{I/p}(S)$.  
\end{enumerate}
\end{construction}

We also show the following overconvergent property of $\BB_I$.

\begin{lemma}\label{lemma:BI-overconvergent} Let $I=[a,b] \subset (0, \infty)$ be a closed interval with rational endpoints and let $S=\Spa(R, R^+)$ be an affinoid perfectoid space over $\Spa(\QQ_p, \ZZ_p)$. Set $S^\circ \coloneqq \Spa(R, R^\circ)$. Then the natural map $\BB_I(S) \to \BB_I(S^\circ)$ is a topological isomorphism. 
\end{lemma}
\begin{proof}
    Both $\BB_I(S)$ and $\BB_I(S^\circ)$ are complete Tate rings, so the Banach open mapping theorem \cite[Lem.\,2.4(i)]{Huber-generalization} implies that it suffices to show that $\BB_I(S) \to \BB_I(S^\circ)$ is a ring isomorphism.

    We choose $p^\flat \in R^{\flat, +} \subset R^{\flat, \circ}$ as in \cref{notation:p-flat} and $c\in \ZZ_{\geq 0}$ as in \cref{defn:affinoid-Fargues--Fontaine}. Note that the $(p, [p^\flat])$-adic topology on $W(R^{\flat, +})\bigl[\tfrac{p^{cb}}{[(p^{\flat})^c]},\tfrac{[(p^\flat)^c]}{p^{ca}}\bigr]$ coincides with the $[p^\flat]$-adic topology because $p^{cb}=\tfrac{p^{cb}}{[(p^{\flat})^c]} \cdot [p^\flat]^c$; the same applies to $W(R^{\flat, \circ})\bigl[\tfrac{p^{cb}}{[(p^{\flat})^c]},\tfrac{[(p^\flat)^c]}{p^{ca}}\bigr]$. Therefore, \cite[Lem.\ and def.\ on p.~516]{Huber-generalization} and \cite[Lem.\,1.6]{H0} imply that the map $\BB_I(S) \to \BB_I(S^\circ)$ becomes equal to the natural morphism
    \[
    W(R^{\flat, +})\bigl[\tfrac{p^{cb}}{[(p^{\flat})^c]},\tfrac{[(p^\flat)^c]}{p^{ca}}\bigr]^{\wedge}_{([p^\flat])}\bigl[\tfrac{1}{[p^{\flat}]}\bigr] \to W(R^{\flat, \circ})\bigl[\tfrac{p^{cb}}{[(p^{\flat})^c]},\tfrac{[(p^\flat)^c]}{p^{ca}}\bigr]^{\wedge}_{([p^\flat])}\bigl[\tfrac{1}{[p^{\flat}]}\bigr].
    \]
    For it to be an isomorphism, it suffices to show that $W(R^{\flat, +})\bigl[\tfrac{p^{cb}}{[(p^{\flat})^c]},\tfrac{[(p^\flat)^c]}{p^{ca}}\bigr]^{\wedge}_{([p^\flat])} \to W(R^{\flat, \circ})\bigl[\tfrac{p^{cb}}{[(p^{\flat})^c]},\tfrac{[(p^\flat)^c]}{p^{ca}}\bigr]^{\wedge}_{([p^\flat])}$ is an injective map with cokernel annihilated by $[p^\flat]$. This reduces to proving that $W(R^{\flat, +})\bigl[\tfrac{p^{cb}}{[(p^{\flat})^c]},\tfrac{[(p^\flat)^c]}{p^{ca}}\bigr] \to W(R^{\flat, \circ})\bigl[\tfrac{p^{cb}}{[(p^{\flat})^c]},\tfrac{[(p^\flat)^c]}{p^{ca}}\bigr]$ is injective with  cokernel annihilated by $p^\flat$. For that, it is enough to show that the map $W(R^{\flat, +}) \to W(R^{\flat, \circ})$ is injective with cokernel annihilated by $[p^\flat]$. Using the explicit description of Witt vectors, we finally reduce the question to showing that the map $R^{\flat, +} \to R^{\flat, \circ}$ is injective with cokernel annihilated by $[p^\flat]$. This can be easily seen from the fact that $R^{\flat, +}$ is perfect (see, for example, \cite[Lem.\,B.13]{zav-almost}). 
\end{proof}

For later use, we will also need the following fact. 

\begin{lemma}
\label{lemma:BI-determined-by-points} 
Let $I=[a,b] \subset (0, \infty)$ be a closed interval with rational endpoints and let $S=\Spa(R, R^+)$ be an affinoid perfectoid space over $\Spa(\QQ_p, \ZZ_p)$. Then the natural morphism
\[
\BB_I(S) \to \prod_{s\in S} \BB_I\bigl(\Spa(\wdh{k(s)}, \wdh{k(s)}^+)\bigr)
\]
is injective.
\end{lemma}
\begin{proof}
    For brevity, if $(K, K^+)$ is a perfectoid field pair, we write $\BB_I(K, K^+)$ for $\BB_I(\Spa(K, K^+))$ throughout this proof. The same applies to $\AA_I(K, K^+)$. 
    
    First, we note that the proof of \cite[Prop.~II.1.1]{Fargues-Scholze} implies that $\BB_I(S)$ is sous-perfectoid. In particular, it is uniform. Therefore, \cite[Th.~5.2.1]{Berkeley} implies that the natural morphism $\BB_I(S) \to \prod_{x\in \Spa(\BB_I(S), \AA_I(S))} \wdh{k(x)}$ is injective. Hence, it suffices to show that, for any $x\in \Spa(\BB_I(S), \AA_I(S))$, the morphism $\BB_I(S) \to \wdh{k(x)}$ factors through $\BB_I(S) \to \BB_I\bigl(\wdh{k(s)}, \wdh{k(s)}^+\bigr)$ for some $s\in S$. For this, we consider the continuous morphism $\pi\colon \abs{\Spa(\BB_I(S), \AA_I(S))} \to \abs{S}$ from \cite[Prop.~II.1.3]{Fargues-Scholze}. Set $s\coloneqq \pi(x)$. Then  \cite[Prop.~II.1.3]{Fargues-Scholze} implies that $x\in \Spa(\BB_I(S), \AA_I(S))$ lies in the image of the pro-(open immersion) 
    \[
    \Spa\Bigl(\BB_I\bigl(\wdh{k(s)}, \wdh{k(s)}^+\bigr), \AA_I\bigl(\wdh{k(s)}, \wdh{k(s)}^+\bigr)\Bigr) \to \Spa\bigl(\BB_I(S), \AA_I(S)\bigr).
    \]
    Thus, the morphism $\BB_I(S) \to \wdh{k(x)}$ factors through $\BB_I(S) \to \BB_I\bigl(\wdh{k(s)}, \wdh{k(s)}^+\bigr)$. 
\end{proof}

\subsection{Period rings \`a la Berger}

Let $S$ be an affinoid perfectoid space over $\Spa(\QQ_p, \Z_p)$. The main goal of this subsection is to relate the definition of $\AA_I(S)$ given in the previous subsection to the more classical one from \cite{Berger02,Faltings-Simpson}.

To motivate the results of this subsection, we note that the integral closure involved in the explicit construction of $\AA_I(S)$ (see \cref{rmk:explicit-construction}) makes it difficult to compute $\AA_I(S)$ or $\AA_I(S)/p$. Therefore, we wish to get another formula for $\AA_I(S)$ which does not involve integral closures. 

Throughout this subsection, we fix a closed interval $I=[a, b]\subset (0, \infty)$ with rational endpoints. 

\begin{construction/proposition}[{cf.\ \cite[Def.~2.1.8]{Kedlaya-slope} and \cite[Prop.~1.4.9]{FF-curve-book}}]
\label{FF Proposition Construction}
Let $(R, R^+)$ be a perfectoid pair of characteristic $p>0$, let $\varpi\in R^+$ be a pseudo-uniformizer, and let $v \colon R \to \Gamma_v \cup \{0\}$ be a valuation lying in $\Spa(R, R^+)$. Then the map $\widetilde{v}_\varpi \colon W(R^+) \to \Gamma_v \cup \{0\}$, defined by 
\[
\widetilde{v}_\varpi\Bigl(\sum_{i\geq 0} [a_i]p^i\Bigr) = \sup_{i \suchthat a_i \neq 0}\bigl\lbrace v(a_i)\cdot v(\varpi)^i \bigr\rbrace \in \Gamma_v \cup \{0\},
\]
is a valuation.
\end{construction/proposition}
\begin{proof}
    Since $v(\varpi)$ is topologically nilpotent, $\widetilde{v}_\varpi(f)$ is well-defined for any $f\in W(R^+)$. To see that it is a valuation, we denote the completed residue field of $v$ by $\wdh{k(v)}$ and its integral subring by $\wdh{k(v)}^+$. Then $\bigl(\wdh{k(v)}, \wdh{k(v)}^+\bigr)$ is a perfectoid field. Since the valuation $\widetilde{v}_\varpi$ factors through the obvious ``reduction'' morphism $W(R^+) \to W(\wdh{k(v)}^+)$, we can replace $(R, R^+)$ with $(\wdh{k(v)}, \wdh{k(v)}^+)$ to assume that $(R, R^+)$ is a perfectoid field. In this case, the proofs of \cite[Prop.~1.4.3 and Prop.~1.4.9]{FF-curve-book} can be easily adapted to show that $\widetilde{v}_\varpi$ is a valuation. 
\end{proof}

\begin{defn}\label{defn:well-adapted} Let $(R, R^+)$ be a perfectoid pair of characteristic $p>0$, let $\varpi \in R^{+}$ be a pseudo-uniformizer, and let $I=[a,b]$ and $c\in \ZZ_{>0}$ be as in \cref{defn:affinoid-Fargues--Fontaine}.
We say a pair of elements $\alpha, \beta \in R^+$ is \emph{well-adapted} to $\varpi$ and $I$ if $v(\alpha)^{ca}=v(\varpi)^c=v(\beta)^{cb}$ for any $v\in \Spa(R, R^+)$.
\end{defn}

First, we note that well-adapted pairs exist \'{e}tale locally on $\Spa(R, R^+)$. 

\begin{rmk}
\label{rmk:often-well-adapted-pairs} 
The existence of a well-adapted pair $\alpha, \beta$ is automatic if the affinoid perfectoid
$S=\Spa(R, R^+)$ is strictly totally disconnected. 
In fact, even $\varpi^{1/a}$ and $\varpi^{1/b}$ exist in $R^+$ in this case.  More generally, there is always an \'etale covering $S'\coloneqq \Spa(R', R'^+) \to S=\Spa(R, R^+)$ such that elements $\varpi^{1/a}$ and $\varpi^{1/b}$ exist in $R'^+$.
\end{rmk}

\begin{rmk}\label{rmk:map-to-well-adapted} Furthermore, the existence of a well-adapted pair is automatic if $S$ admits a map $S \to T$ such that $T$ admits a well-adapted pair and the pseudo-uniformizer $\varpi$ comes from $T$. 
\end{rmk}

\begin{rmk}\label{rmk:well-adapted-pseudo-uniformizer} We note that \cite[Lem.~1.6(iv)]{Huber-generalization} implies that for any pair $\alpha, \beta\in R^+$ well-adapted to $\varpi$ and $I$, there are units $u_1, u_2\in (R^+)^\times$ such that $\alpha^{ca}\cdot u_1 = \varpi^c = \beta^{cb} \cdot u_2$. In particular, any such $\alpha$ and $\beta$ are pseudo-uniformizers. Similarly, we note that $\alpha/\beta$ lies in $R^+$ for any well-adapted pair $\alpha, \beta$. 
\end{rmk}

\begin{lemma}
\label{AI-intclos}
Let $(R, R^+)$ be a perfectoid pair of characteristic $p>0$, let $\varpi \in R^+$ be a pseudo-uniformizer, and let $\alpha, \beta\in R^+$ be a pair of elements well-adapted to $\varpi$ and $I$. Set $A\coloneqq W(R^+)$. Then the ring $A\bigl[\frac{[\alpha]}{p}, \frac{p}{[\beta]}\bigr]$
is integrally closed inside $A\bigl[\frac{[\alpha]}{p}, \frac{p}{[\beta]} \bigr][\frac{1}{p}]$.
\end{lemma}
\begin{proof}
   For brevity, we denote by $A_{[\alpha,\beta]}$ the ring $A\big[\frac{[\alpha]}{p}, \frac{p}{[\beta]}\big]$ and by $B_{[\alpha,\beta]}$ the ring $A\big[\frac{[\alpha]}{p}, \frac{p}{[\beta]}\big]\big[\frac{1}{p}\big]$.
   Pick an element $f\in B_{[\alpha,\beta]}$ such that $f$ is integral over $A_{[\alpha,\beta]}$.
   We wish to show that $f\in A_{[\alpha,\beta]}$.
   An easy inductive argument reduces the general case to the situation when $f\in \frac{1}{p}\cdot A_{[\alpha,\beta]}$. In this case, we use that any element $a\in A=W(R^+)$ can be written as $a=[a_0] + p\cdot a'$ with $a_0\in R^+$ and $a'\in A$. Therefore, any element $f\in \frac{1}{p}\cdot A_{[\alpha,\beta]}$ can be written as the sum
   \[
   f = \frac{1}{p} \cdot \Bigl(c + \sum_{i=1}^N a_i \big( \frac{[\alpha]}{p}\big)^i + \sum_{i=1}^N b_i \big( \frac{p}{[\beta]}\big)^i \Bigr) = \frac{1}{p} \cdot \biggl([c_0] + \sum_{i=1}^N [a_{i, 0}] \big( \frac{[\alpha]}{p}\big)^i + \sum_{i=1}^N [b_{i, 0}] \big( \frac{p}{[\beta]}\big)^i \biggr) + \biggl(c' + \sum_{i=1}^N a'_i \big( \frac{[\alpha]}{p}\big)^i + \sum_{i=1}^N b'_i \big( \frac{p}{[\beta]}\big)^i)\biggr),
   \]
   where $c, c', a_i, a'_i, b_i, b'_i\in A$ and $c_0, a_{i, 0}, b_{i, 0}\in R^+$. Since the second summand   lies in $A_{[\alpha, \beta]}$, we immediately reduce the question to the case when
    \begin{equation}\label{eqn:decomposition-f}
    f = \frac{1}{p}\cdot \Big( [c_0] + \sum_{i=1}^N[a_{i,0}]\cdot \big(\frac{[\alpha]}{p}\big)^i + \sum_{i=1}^N [b_{i,0}]\big(\frac{p}{[\beta]}\big)^i\Big),
    \end{equation}
    where $c_0$, $a_{i, 0}$, $b_{i, 0}$ lie in $R^+$ and $N$ is a positive integer.
    One checks easily that for any valuation $v\in \Spa(R, R^+)$, the valuations $\widetilde{v}_{\alpha}$ and $\widetilde{v}_{\beta}$ from \cref{FF Proposition Construction} (and \cref{rmk:well-adapted-pseudo-uniformizer}) extend uniquely to valuations $\widetilde{v}_{\alpha}, \widetilde{v}_{\beta} \colon A\big[\frac{1}{p}, \frac{1}{[\beta]}\big] \to \Gamma_v \cup \{0\}$ such that the image of $A_{[\alpha,\beta]}\subset A\big[\frac{1}{p}, \frac{1}{[\beta]}\big]$ under either of these valuations lies inside $\Gamma_{\leq 1} \cup \{0\}$.
    Since $f$ is integral over $A_{[\alpha,\beta]}$, we conclude that $\widetilde{v}_{\alpha}(f)\leq 1$ and $\widetilde{v}_\beta(f)\leq 1$.
    Using \cref{eqn:decomposition-f}, we conclude that 
    \[
    v(c_0)\leq v(\alpha), v(a_{i,0})\leq v(\alpha), \text{ and } v(b_{i,0})\leq v(\beta). 
    \]
    Since this holds for any valuation $v\in \Spa(R, R^+)$, \cite[Lem.~1.6(iv)]{Huber-generalization} implies that $c_0/\alpha$, $a_{i, 0}/\alpha$, and $b_{i, 0}/\beta$ lie in $R^+$. Therefore, we conclude that 
    \[
    f = \left[\frac{c_0}{\alpha}\right] \cdot \frac{[\alpha]}{p} + \sum_{i=1}^N \left[\frac{a_{i,0}}{\alpha}\right] \cdot \left(\frac{[\alpha]}{p}\right)^{i+1} + \sum_{i=1}^N \left[\frac{b_{i, 0}}{\beta}\right]\cdot \left(\frac{p}{[\beta]}\right)^{i-1} 
    \]
    lies in $A_{[\alpha,\beta]}$. This finishes the proof. 
\end{proof}

Finally, we are ready to give the promised description of $\AA_I(S)$.

\begin{lemma}
\label{lemma:semi-explicit-A_I} 
Let $S=\Spa(R, R^+)$ be an affinoid perfectoid space over $\Spa(\QQ_p, \Z_p)$ with a choice of $p^\flat \in R^{\flat, +}$, and let $\alpha, \beta\in R^{\flat, +}$ be a pair of elements well-adapted to $p^\flat$ and $I$. Then there is a topological isomorphism $\AA_I(S) \simeq \AA_{\inf}(S)\bigl\langle \frac{[\alpha]}{p}, \frac{p}{[\beta]} \bigr\rangle \colonequals \bigl(\AA_{\inf}(S)\bigl[ \frac{[\alpha]}{p}, \frac{p}{[\beta]} \bigr]\bigr)^\wedge_{(p)}$, which is moreover functorial in $S$.    
\end{lemma}
Recall that $\AA_I(S)$ is a topological ring due to its definition as the global sections of a sheaf of topological rings;
cf.\ \cref{rmk:BI-AI-Tate-Huber-pair}.
\begin{proof}[{Proof of \cref{lemma:semi-explicit-A_I}}]
    Pick $c$ as in \cref{defn:affinoid-Fargues--Fontaine}. Then \cref{defn:affinoid-Fargues--Fontaine} and \cref{defn:well-adapted} imply that
    \[
    Y_{S, I} = \bigl\{ x\in Y_S \mid \abs{p^{cb}(x)} \le \abs{[(p^{\flat})^c](x)} \le \abs{p^{ca}(x)} \neq 0\bigr\} = \bigl\{ x\in Y_S \mid \abs{\alpha(x)} \le \abs{p(x)} \le \abs{\beta(x)} \neq 0 \bigr\} \subset Y_S.
    \]
    Thus, \cite[Lem.\ and def.\ on p.~516]{Huber-generalization} implies that $\AA_I(S)$ is equal to the topological completion of the integral closure of $\AA_{\inf}(S)\big[\frac{[\alpha]}{p}, \frac{p}{[\beta]}\big]$ in $\AA_{\inf}(S)\big[\frac{[\alpha]}{p}, \frac{p}{[\beta]} \big][\frac{1}{p}]$, where we topologize $\AA_{\inf}(S)\big[\frac{[\alpha]}{p}, \frac{p}{[\beta]}\big]$ using the $(p, [p^\flat])$-adic topology, we topologize $\AA_{\inf}(S)\big[\frac{[\alpha]}{p}, \frac{p}{[\beta]} \big][\frac{1}{p}]$ via the unique topology which makes $\AA_{\inf}(S)\big[\frac{[\alpha]}{p}, \frac{p}{[\beta]}\big]$ open and bounded, and we topologize the integral closure via the subspace topology from $\AA_{\inf}(S)\big[\frac{[\alpha]}{p}, \frac{p}{[\beta]} \big][\frac{1}{p}]$.
    Then \cref{rmk:well-adapted-pseudo-uniformizer} implies that $[p^\flat]^c \in p\cdot \AA_{\inf}(S)\bigl[\tfrac{[\alpha]}{p},\tfrac{p}{[\beta]}\bigr]$, so the $p$-adic and the $(p, [p^\flat])$-adic topologies on $\AA_{\inf}(S)\bigl[\tfrac{[\alpha]}{p},\tfrac{p}{[\beta]}\bigr]$ agree.
    
    Furthermore, \cref{AI-intclos} implies that $\AA_{\inf}(S)\big[\frac{[\alpha]}{p}, \frac{p}{[\beta]}\big]$ is already integrally closed in $\AA_{\inf}(S)\big[\frac{[\alpha]}{p}, \frac{p}{[\beta]} \big][\frac{1}{p}]$. Thus, $\AA_{I}(S) \simeq \AA_{\inf}(S)\bigl\langle \frac{[\alpha]}{p}, \frac{p}{[\beta]} \bigr\rangle$.
\end{proof}

Even though \cref{lemma:semi-explicit-A_I} gives an explicit description of $\AA_I(S)$ under a mild assumption on $S$, it is still not quite strong enough to give an explicit formula for $\AA_I(S)/p$. In order to achieve this result, we will need the following preliminary lemma. 

\begin{lemma}[{cf.\ \cite[Lemme 2.1]{Berger02}}]
\label{Berger lemme}
Let $R$ be a perfect $\FF_p$-algebra, let $\beta, \gamma, \varpi \in R$ be three nonzerodivisors in $R$. 
Set $\alpha \coloneqq \beta \cdot \gamma$, $A\coloneqq W(R)$, and 
$J\coloneqq ([\beta]y - p, xy - [\gamma]) \subset A[x, y]$. Then the following hold:
\begin{enumerate}[label={\upshape{(\roman*)}}]
    \item\label{Berger lemme-1} $J \cap p^nA[x, y] = p^nJ$ for any integer $n\geq 0$;
    \item\label{Berger lemme-2} $J \cap [\varpi]A[x, y] = [\varpi] J$.
\end{enumerate}

\end{lemma}
\begin{proof}
We only prove \cref{Berger lemme-2}, as the proof of \cref{Berger lemme-1} is similar (but with simpler notation).\footnote{
In the proof of \cref{Berger lemme-2}, we only need to use the regularity of $\varpi$ and $\beta$.
However, the proof of \cref{Berger lemme-1} would need to use the regularity of $\beta$ and $\gamma$ instead.} 
We start by showing the following claims:

\begin{claim}
\label{claim:decomposition-in-A} 
Any polynomial $f(x, y) \in A[x, y]$ can be \emph{uniquely} written as
\begin{equation}\label{equation:decomposition-A_I}
f(x, y) = c_0 + g(x) \cdot x + h(y) \cdot y + r(x, y) \cdot (xy - [\gamma]),
\end{equation}
where $c_0\in A$, $g(x) \in A[x]$, $h(y)\in A[y]$, and $r(x,y) \in A[x,y]$. 
\end{claim}
\begin{proof}
A simple inductive argument shows that any polynomial can be written in this way. To show uniqueness, it suffices to show that $0=c_0 + g(x) \cdot x + h(y) \cdot y + r(x, y) \cdot (xy - [\gamma])$ implies that $c_0$, $g(x)$, $h(y)$, and $r(x, y)$ are all equal to $0$.
We write $g(x) = \sum_{i\geq 0} a_i x^i$, $h(y) = \sum_{i\geq 0} b_i y^i$, and $r(x,y) = \sum_{i, j\geq 0} c_{i, j} x^i y^j$. Then our assumption implies that $c_{i, j}=[\gamma] c_{i+1, j+1}$ for any $i, j\geq 0$. Since $c_{i, j}=0$ for $i,j\gg 0$, we conclude that all $c_{i, j}=0$. Similarly, our assumption also implies that $c_0=c_{0,0} [\gamma]$, $a_{i}=c_{i+1, 0}[\gamma]$, and $b_{i} = c_{0, i+1}[\gamma]$. Thus, $c_0=a_i=b_i=0$ for any $i\geq 0$. This finishes the proof of the claim. 
\end{proof}

\begin{claim}\label{claim:decomposition-I-in-A} Any polynomial $F(x, y) \in J$ can be \emph{uniquely} written as
\begin{equation}\label{equation:decomposition-I-in-A_I}
F(x, y) = G(x) \cdot (px - [\alpha]) + H(y) \cdot ([\beta]y - p) + R(x, y) \cdot (xy - [\gamma]).
\end{equation}
where $G(x)\in A[x]$, $H(y)\in A[y]$, and $R(x, y)\in A[x,y]$ . 
\end{claim}
\begin{proof}
By definition, any element $F(x, y)$ in $J$ can be written as
$F(x, y) = f_1(x, y) \cdot ([\beta]y - p) + f_2(x, y) \cdot (xy - [\gamma])$
for some $f_i(x,y) \in A[x,y]$, $i=1,2$.
By \Cref{claim:decomposition-in-A}, there are $c_0\in A$, $g(x) \in A[x]$, $h(y)\in A[y]$, and $r(x,y) \in A[x,y]$ such that
$f_1(x, y) = g(x) \cdot x + (c_0 + h(y) \cdot y) + r(x, y) \cdot (xy - [\gamma])$.
Therefore, we have
\[
F(x, y) = g(x) \cdot x \cdot ([\beta]y - p) +  (c_0 + h(y) \cdot y) \cdot ([\beta]y - p) +
\big(r(x, y) \cdot ([\beta]y - p) + f_2(x, y)\big) \cdot (xy - [\gamma]).
\]
Since $x \cdot ([\beta]y - p) = [\beta] \cdot (xy - [\gamma]) - (px - [\alpha])$,
we further write
\[
F(x, y) = -g(x) \cdot (px - [\alpha]) + (c_0 + h(y) \cdot y) \cdot ([\beta]y - p) +
\big(g(x) \cdot [\beta] + r(x, y) \cdot ([\beta]y - p) + f_2(x, y)\big) \cdot (xy - [\gamma]).
\]
This finishes the proof that any $F(x, y)\in J$ can be written in the desired form. 
Next, let us show the claim of uniqueness. 
Suppose $G(x) \cdot (px - [\alpha]) + H(y) \cdot ([\beta]y - p) + R(x, y) \cdot (xy - [\gamma]) = 0$.
Using \cref{claim:decomposition-in-A}, we see that $R(x, y) = 0$.
Then using that $p$ and $[\beta]$ are nonzerodivisors in $A$, we can run an argument similar to the uniqueness argument in the 
proof of \cref{claim:decomposition-in-A} to show that both $G(x)$ and $H(y)$ are zero.
\end{proof}

Now suppose that $f(x, y)= F(x,y)$ lies in $J\cap [\varpi]A[x,y]$. Let $c_0$, $g(x)$, $h(y)$, $r(x,y)$, $G(x)$, $H(y)$, and $R(x,y)$ be the elements from \cref{claim:decomposition-in-A} and \cref{claim:decomposition-I-in-A} respectively. Then the uniqueness claim in \cref{claim:decomposition-in-A} implies that 
\begin{equation}\label{eqn:two-decompositions}
c_0 + g(x) \cdot x + h(y) \cdot y = G(x) \cdot (px - [\alpha]) + H(y) \cdot ([\beta]y - p)
\text{\hspace{0.4 cm} and \hspace{0.4 cm}} r(x, y) = R(x, y).
\end{equation}

\begin{claim}\label{claim:divisible-by-pi} In the notation from above, for any polynomial $f(x, y) \in J\cap [\varpi]A[x,y]$, all the coefficients of $G(x)$, $H(y)$, and $R(x,y)$ are divisible by $[\varpi]$.  
\end{claim}
\begin{proof}
    We consider Equality~\cref{eqn:two-decompositions}.
    The uniqueness claim in \cref{claim:decomposition-in-A} implies that $c_0$ and all coefficients of $g(x)$, $h(y)$, and $r(x,y)$ are divisible by $[\varpi]$.
    Since $r(x, y) = R(x, y)$, we immediately deduce that all coefficients of $R(x, y)$ are divisible by $[\varpi]$.

    Now we deal with the coefficients of $G(x)$. For this, we write $g(x) = \sum_{i\geq 0} b_i x^i$ and $G(x) = \sum_{i\geq 0} B_ix^i$ with $b_i, B_i\in A$. 
    Note that $p$ and $[\varpi]$ are nonzerodivisors in $A$ and $(p, [\varpi])$ is a regular sequence in $A$, so \cite[\href{https://stacks.math.columbia.edu/tag/07DW}{Tag 07DW}]{stacks-project} implies that $([\varpi], p)$ is a regular sequence in $A$ as well.
    Now Equality~\cref{eqn:two-decompositions} implies that $b_i = B_i p - B_{i+1}[\alpha]$ for any $i\geq 0$. Using that $B_{i+1}=0$ (in particular, divisible by $[\varpi]$) for $i\gg 0$, $b_i$ is divisible by $[\varpi]$ for all $i \ge 0$, and the fact that $([\varpi], p)$ is a regular sequence, a descending induction on $i$ proves that all $B_i$ are divisible by $[\varpi]$. 

    Finally, a similar argument using ascending induction and the fact that $([\varpi], p)$ is a regular sequence shows that all coefficients of $H(y)$ are also divisible by $[\varpi]$. 
\end{proof}
\Cref{claim:divisible-by-pi} yields the inclusion $J \cap [\varpi]A[x, y] \subset [\varpi] J$ in \cref{Berger lemme-2}.
Since the other inclusion is immediate, this finishes the proof.
\end{proof}

\begin{cor}\label{cor:explicit-formula-for-naive-AI}
    In the notation of \cref{Berger lemme}, the quotient $A[x, y]/J$ is $p[\beta]$-torsionfree and is isomorphic to $A\bigl[\frac{[\alpha]}{p}, \frac{p}{[\beta]}\bigr]$. Furthermore, the $p$-adic completion $A\bigl\langle \frac{[\alpha]}{p}, \frac{p}{[\beta]}\bigr\rangle \coloneqq \bigl(A\bigl[\frac{[\alpha]}{p}, \frac{p}{[\beta]}\bigr]\bigr)^{\wedge}_{(p)}$ is isomorphic to 
    $\frac{A\langle x, y \rangle}{([\beta]y - p, xy - [\gamma])}$ and is $p$-torsionfree.
\end{cor}
\begin{proof}
    First, \cref{Berger lemme} with $\varpi = \beta$ and the fact that $p$ and $[\varpi]$ are nonzerodivisors in $A$ imply that $A[x, y]/J$ is $p$-torsionfree and $[\beta]$-torsionfree. In particular, it is flat over $\ZZ_p$ and $p[\beta]$-torsionfree. 
    
    Now we recall that $A\bigl[\frac{[\alpha]}{p}, \frac{p}{[\beta]}\bigr]$ is defined as the sub-$A$-algebra of $A\bigl[\frac{1}{p}, \frac{1}{[\beta]}\bigr]$ generated by $\frac{[\alpha]}{p}$ and $\frac{p}{[\beta]}$. Therefore, we have a natural $A$-linear surjective ring homomorphism $\pi\colon A[x, y]/J \twoheadrightarrow A\bigl[\frac{[\alpha]}{p}, \frac{p}{[\beta]}\bigr]$ which sends $x$ to $\frac{[\alpha]}{p}$ and $y$ to $\frac{p}{[\beta]}$. We wish to show that $\pi$ is also injective. For this, we note that $\pi\bigl[\frac{1}{p[\beta]}\bigr]$ is an isomorphism since 
    \[
    \Bigl(\frac{A[x,y]}{J}\Bigr)\Bigl[\frac{1}{p[\beta]}\Bigr]\simeq 
    \frac{A\bigl[\frac{1}{p}, \frac{1}{[\beta]}\bigr][x,y]}{\bigl(y-\frac{p}{[\beta]}, xy-[\gamma]\bigr)}  
    \simeq \frac{A\bigl[\frac{1}{p}, \frac{1}{[\beta]}\bigr][x,y]}{\bigl(x-\frac{[\alpha]}{p}, y-\frac{p}{[\beta]}\bigr)} \simeq A\Bigl[\frac{1}{p}, \frac{1}{[\beta]}\Bigr].
    \]
    Furthermore, we note that the natural morphism $A[x,y]/J \to (A[x,y]/J)[\frac{1}{p[\beta]}]$ is injective since $A[x,y]/J$ is $p[\beta]$-torsionfree. This formally implies that $\pi$ must be injective as well. Hence, it is bijective. 

    Finally, we note that \cite[\href{https://stacks.math.columbia.edu/tag/0315}{Tag 0315}, part~(3)]{stacks-project} implies that $A\bigl\langle \frac{[\alpha]}{p}, \frac{p}{[\beta]}\bigr\rangle\simeq A\langle x,y\rangle/J^{\wedge}_{(p)}$ because $A\bigl[\frac{[\alpha]}{p}, \frac{p}{[\beta]}\bigr]$ is flat over $\Z_p$. Then \cite[\href{https://stacks.math.columbia.edu/tag/0315}{Tag 0315}, part~(1)]{stacks-project} implies that 
    $A\langle x,y\rangle^{\oplus 2} \xr{([\beta]y-p, xy-[\gamma])} J^{\wedge}_{(p)}$ 
    is surjective. Thus, $A\bigl\langle \frac{[\alpha]}{p}, \frac{p}{[\beta]}\bigr\rangle \simeq 
    A\langle x,y\rangle/([\beta]y-p, xy-[\gamma])$. 
    Lastly, $A\bigl\langle \frac{[\alpha]}{p}, \frac{p}{[\beta]}\bigr\rangle$ is $p$-torsionfree because it is the $p$-adic completion of a $p$-torsionfree ring. 
\end{proof}
Finally, we are ready to prove the desired formula for $\AA_I(S)$. 
\begin{cor}\label{cor:AI-mod-p-rings} 
Let $S=\Spa(R, R^+)$ and $\alpha, \beta\in R^{\flat, +}$ be as in \cref{lemma:semi-explicit-A_I}.
Then 
\[ \AA_I(S) \simeq \AA_{\inf}(S)\bigl\langle \tfrac{[\alpha]}{p}, \tfrac{p}{[\beta]} \bigr\rangle \simeq 
\frac{\AA_{\inf}(S)\langle x, y \rangle}{([\beta]y - p, xy - [\alpha/\beta])}. \] 
In particular, there is an isomorphism of $R^{\flat,+}$-modules \emph{(}functorial in $S$\emph{)}
\[
\AA_I(S)/p \simeq \Bigl(\bigoplus_{i\geq 0} R^{\flat, +}/\alpha \Bigr) \bigoplus \Bigl(\bigoplus_{i\geq 1} R^{\flat, +}/\beta \Bigr).
\]
\end{cor}
\begin{proof}
    The isomorphism $\AA_I(S) \simeq 
\frac{\AA_{\inf}(S)\langle x, y \rangle}{([\beta]y - p, xy - [\alpha/\beta])}$ follows
immediately from \cref{rmk:well-adapted-pseudo-uniformizer}, \cref{lemma:semi-explicit-A_I}, and \cref{cor:explicit-formula-for-naive-AI}.
Since $\alpha = -\beta \cdot (xy-\tfrac{\alpha}{\beta}) + x \cdot \beta y$ in $R^{\flat,+}[x,y]$, this isomorphism formally implies that 
    \[
    \AA_I(S)/p\simeq 
    \frac{\AA_{\inf}(S)\langle x, y \rangle}{\big(p, [\beta]y - p, xy - [\alpha/\beta]\big)} \simeq \frac{(R^{\flat, +}/\alpha)[x,y]}{\big(\beta y, xy - \alpha/\beta\big)} 
    \simeq \Big(\bigoplus_{i\geq 0} (R^{\flat, +}/\alpha) \cdot x^i \Big) \bigoplus \Big( \bigoplus_{i\geq 1} (R^{\flat, +}/\beta) \cdot y^i \Big). \qedhere
    \]
\end{proof}
\begin{cor}\label{cor:almost-mathematics-AI}
Let $S=\Spa(R, R^+)$ and $p^\flat ,\alpha, \beta\in R^{\flat, +}$ be as in \cref{lemma:semi-explicit-A_I}.
Then the natural, $p$-adic and $[p^\flat]$-adic topology on $\AA_I(S)$ coincide (in particular, $\AA_I(S)$ is $p$-adically complete), the elements $p$ and $[p^\flat] \in \AA_I(S)$ are nonzerodivisors, $\AA_I(S)^{\circ\circ} = \cup_{n\geq 1} [(p^\flat)^{1/p^n}] \AA_I(S)$, and $\AA_I(S)^{\circ\circ} \subset \AA_I(S)$ is an ideal of almost mathematics \emph{(}see \cite[2.1.1]{Gabber-Ramero}\emph{)}.
\end{cor}
\begin{proof}
    \cref{lemma:semi-explicit-A_I} implies that the topology on $\AA_I(S)$ is the $p$-adic topology and $\AA_I(S)$ is the $p$-adic completion of $\AA_{\inf}(S)\bigl[\frac{[\alpha]}{p}, \frac{p}{[\beta]}\bigr]$. In order to see that the $p$-adic and the $[p^\flat]$-adic topologies on $\AA_I(S)$ coincide, we choose a positive integer $c$ such that $c\cdot a$ and $c\cdot b$ are positive integers.
    By \Cref{rmk:explicit-construction}, we see that
    $[p^\flat]^c\in p^{ac} \cdot \AA_{I}(S)$ and $p^{bc} \in [p^\flat]^c\AA_{I}(S)$, so the two topologies coincide. Furthermore, this argument shows the equality of vanishing loci $\abs{\rm{V}(p)} = \abs{\rm{V}([p^\flat])} \subset \abs{\Spec \AA_I(S)}$. 
    
    Now we note that $p$ is a nonzerodivisor of $\AA_I(S)$ due to \cref{cor:explicit-formula-for-naive-AI}. Furthermore, the previous paragraph implies that $[p^\flat]$ becomes invertible in $\BB_I(S)=\AA_I(S)[\frac{1}{p}]$. Since $p$ is regular in $\AA_I(S)$, this implies that $[p^\flat]$ is regular in $\AA_I(S)$ as well. Furthermore, the previous paragraph also implies that $\BB_I(S)=\AA_I(S)\bigl[\frac{1}{[p^\flat]}\bigr]$. 
    
    Now we show that $\AA_I(S)^{\circ\circ} = \cup_{n\geq 1} [(p^\flat)^{1/p^n}] \AA_I(S)$ and that $\AA_I(S)^{\circ\circ} \subset \AA_I(S)$ is an ideal of almost mathematics. Since $R^{\flat, +}$ is a perfect ring, we conclude that $(p^{\flat})^{1/p^n}$ is well-defined for any $n\geq 1$. Furthermore, since $[p^\flat]$ lies in $\AA_I(S)^{\circ\circ}$ and $\AA_I(S)^{\circ\circ}$ is a radical ideal, we have $\cup_{n\geq 1} [(p^\flat)^{1/p^n}] \AA_I(S) \subset \AA_I(S)^{\circ\circ}$. Now we show the other inclusion. Let $f\in \AA_I(S)^{\circ\circ}$. Then there is an integer $n$ such that $f^{p^n}\in [p^\flat]\AA_I(S)$ because the topology on $\AA_I(S)$ coincides with the $[p^\flat]$-adic topology. Then 
    \[
    \Big(\frac{f}{[(p^\flat)^{1/p^n}]} \Big)^{p^n} \in \AA_I(S).
    \]
    Since $\AA_I(S)\subset \AA_I(S)\bigl[\frac{1}{[p^\flat]}\bigr] = \BB_I(S)$ is integrally closed, we conclude that $\frac{f}{[(p^\flat)^{1/p^n}]} \in \AA_I(S)$. In other words, $f\in [(p^\flat)^{1/p^n}] \AA_I(S)$. Hence, $\cup_{n\geq 1} [(p^\flat)^{1/p^n}] \AA_I(S) = \AA_I(S)^{\circ\circ}$. Finally, the fact that $\AA_I(S)^{\circ\circ}$ defines an ideal of almost mathematics follows directly from the observation that 
    \[
    \big( \bigcup_{n\geq 1} [(p^\flat)^{1/p^n}] \AA_I(S) \big)^2 = \bigcup_{n\geq 1} [(p^\flat)^{1/p^n}] \AA_I(S)
    \]
    and \cite[Prop.~2.1.7(i)]{Gabber-Ramero}. 
\end{proof}

As an application, we show that $\BB_I(S)$ is a PID when $S=\Spa(C, C^+)$.

\begin{cor}[{cf.\,\cite[Prop.\,2.6.8]{Kedlaya-slope}}]\label{cor:BI-PID} Let $I=[a,b] \subset (0, \infty)$ be a closed interval with rational endpoints, let $(C, C^+)$ be a complete Tate--Huber pair such that $C$ is an algebraically closed nonarchimedean field of mixed characteristic $(0, p)$. Set $S\coloneqq \Spa(C, A^+)$. Then $\BB_I(S)$ is a principal ideal domain. 
\end{cor}

We emphasize that we do not need to assume that $C^+$ is a valuation subring of $C$. 

\begin{proof}
    \cref{lemma:BI-overconvergent} implies that it suffices to assume that $C^+=C^\circ$. For brevity, we set $B_I\coloneqq \BB_I\bigl(\Spa(C, C^\circ)\bigr)$ and $A_I\coloneqq \AA_I\bigl(\Spa(C, C^\circ)\bigr)$, and $A_{\inf} \coloneqq W(C^{\flat, \circ})$. We also set $\m_{C^\flat} \subset C^{\flat, \circ}$ to be the maximal ideal of $C^{\flat, \circ}$ and $k_{C^\flat} \coloneqq  C^{\flat, \circ}/\m_{C^\flat}$ to be its residue field. 
    Then \cite[Cor.\,II.1.12]{Fargues-Scholze} or \cite[Th.\,7.11]{Kedlaya2016} imply that it suffices to show that $\Spa(B_I, A_I)$ is connected or, equivalently, that $B_I$ has no non-trivial idempotents. Since $A_I$ is integrally closed in $B_I$, it suffices to show that $A_I$ has no non-trivial idempotents. 

    For this, we choose $\alpha, \beta, p^\flat\in C^{\flat, \circ}$ as in \cref{AI-intclos}. Then \cref{cor:AI-mod-p-rings} ensures that $A_I = A_{\inf}\bigl[\frac{[\alpha]}{p}, \frac{p}{[\beta]}\bigr]^{\wedge}_{(p)}$. Note that $[\alpha] = p\cdot \frac{[\alpha]}{p}\in pA_I$. Since $C^{\flat, \circ}$ is a rank-$1$ valuation ring, we conclude that $[a]\in \mathrm{rad}(pA_I)$ for any $a\in \m_{C^\flat}$. In particular, we have $\mathrm{rad}\bigl((p,[a])_{a\in \m_{C^\flat}}\bigr) = \mathrm{rad}(p) \subset A_I$. Since $A_I$ is $p$-adically complete, we see that \cite[\href{https://stacks.math.columbia.edu/tag/0ALJ}{Tag 0ALJ}]{stacks-project} and \cite[\href{https://stacks.math.columbia.edu/tag/09XJ}{Tag 09XJ}]{stacks-project} imply that $A_I$ is $K\coloneqq (p, [a])_{a\in \m_{C^\flat}}$-adically henselian. Therefore, \cite[\href{https://stacks.math.columbia.edu/tag/09XI}{Tag 09XI}]{stacks-project} ensures that the natural map $\Idem(A_I) \to \Idem(A_I/K)$ is a bijection. Finally, the explicit formula for $A_I$ from \cref{cor:AI-mod-p-rings} implies that 
    \[
    A_I/K \simeq k_{C^\flat}[X, Y]/(XY-\ov{\alpha/\beta}),
    \]
    where $\ov{\alpha/\beta}$ is the residue class of $\alpha/\beta$ in $k_C$. By inspection, we see that $A_I/K$ has no non-trivial idempotents. This finishes the proof. 
\end{proof}

\subsection{K\"unneth formula for period rings}

Throughout this subsection, we fix a closed interval $I=[a,b] \subset (0, \infty)$ with rational endpoints. 
The main goal of this subsection is to show a version of the K\"unneth formula for $\BB_I(\blank)$. 
This will be one of our key ingredients in proving the K\"unneth formula 
(see \Cref{Kunneth formula: main thm}) for the condensed $v$-cohomology of perfect $\BB_I$-complexes.  

\begin{notation}\label{notation:kummer-field} A {\it Kummer field} $\QQ_{p, \infty}$ is the field  $\big(\bigcup_{n\geq 1}\QQ_p(p^{1/p^n})\big)^{\wedge}$ for some choice of compatible $p$-power roots $p^{1/p^n}$ in $\overline{\QQ}_p$. We denote its ring of integers by $\Z_{p, \infty}$.
\end{notation}

We first show that the construction $Y_{\blank, I}$ commutes with fiber products of affinoid perfectoid spaces. 

\begin{lemma}
\label{lemma:kunneth-for-FF-balls}
    Let $S' \to S$ and $S''\to S$ be morphisms of affinoid perfectoid spaces over $\Spa(\QQ_p,\ZZ_p)$. Then the natural morphism $Y_{S'\times_S S'', I} \to Y_{S', I} \times_{Y_{S, I}} Y_{S'', I}$ is an isomorphism. 
\end{lemma}
\begin{proof}
    Fix a choice of Kummer field $\QQ_{p,\infty}$.
    The spaces $Y_{S, I}$, $Y_{S', I}$, and $Y_{S'', I}$ are $\QQ_{p,\infty}$-sousperfectoid (in the sense of \cite[Def.~3.4]{GIZ}) by virtue of \cite[Prop.~II.1.1]{Fargues-Scholze}.
    As fiber products exist in the category of $\QQ_{p,\infty}$-sousperfectoid spaces and are given by completed tensor products on affinoids \cite[Prop.~3.6]{GIZ}, we conclude that $Y_{S', I} \times_{Y_{S, I}} Y_{S'', I}$ is also a $\QQ_{p,\infty}$-sousperfectoid space. Since the diamondification functor is conservative on $\QQ_{p,\infty}$-sousperfectoid spaces \cite[Prop.~3.12]{GIZ} and preserves fiber products, the question gets reduced to showing that the natural morphism $Y_{S'\times_S S'', I}^\lozenge \to Y_{S', I}^\lozenge \times_{Y_{S, I}^\lozenge} Y_{S'', I}^\lozenge$ is an isomorphism.
    This follows immediately from \cite[Prop.\,II.1.17]{Fargues-Scholze}.
\end{proof}

\begin{lemma}\label{lemma:almost-integral-kunneth-formula} Let $S' \to S$ and $S'' \to S$ be maps of affinoid perfectoid spaces. Then $\O^+(S') \wdh{\otimes}^{\rL}_{\O^+(S)} \O^+(S'') \to \O^+(S' \times_S S'')$ is an almost isomorphism with respect to the ideal $\O^{\circ\circ}(S) \subset \O^+(S)$. 
\end{lemma}
Note that $\O^{\circ\circ}(S) \subset \O^+(S)$ is an ideal of almost mathematics (as in \cite[2.1.1]{Gabber-Ramero}) by \cite[Lem.~B.12]{zav-almost}.
\begin{proof}[{Proof of \cref{lemma:almost-integral-kunneth-formula}}]
    For brevity, we write $S=\Spa(R, R^+)$, $S'=\Spa(R', R^{\prime,+})$, $S''=\Spa(R'', R^{\prime\prime,+})$, and $S' \times_S S'' = \Spa(T, T^+)$. Then \cite[Lem.~B.13]{zav-almost} implies that it suffices to show that $R^{\prime,\circ} \wdh{\otimes}^{\rL}_{R^\circ} R^{\prime\prime,\circ} \to T^\circ$ is an almost isomorphism. 

    We choose a pseudo-uniformizer $\varpi\in R^+$ such that $\varpi^p | p$ and $\varpi$ admits all $p$-power roots (see \cite[Lem.~3.10]{diamonds}). A choice of such roots defines an element $\varpi^\flat \in R^{\flat, \circ}$. Then \cite[Lem.~A.5]{zav-almost} implies that it suffices to show that
    \[
    \bigl(R^{\prime,\circ} \wdh{\otimes}^{\rL}_{R^\circ} R^{\prime\prime,\circ}\bigr)/^{\rL}\varpi \to T^{\circ}/^\rL\varpi 
    \]
    is an almost isomorphism. Since $R^{\circ}/^\rL\varpi \simeq R^{\circ}/\varpi \simeq R^{\flat, \circ}/\varpi^\flat \simeq R^{\flat, \circ}/^\rL \varpi^\flat$ and the same applies to $R'$, $R''$, and $T$, it is enough to prove that the natural morphism 
    \[
    \bigl(R^{\prime,\flat, \circ} \wdh{\otimes}^{\rL}_{R^{\flat, \circ}} R^{\prime\prime,\flat, \circ}\bigr)/^\rL \varpi^\flat \to T^{\flat, \circ}/^\rL \varpi^\flat
    \]
    is an almost isomorphism. For this, it suffices to show that $R^{\prime,\flat, \circ} \wdh{\otimes}^{\rL}_{R^{\flat, \circ}} R^{\prime\prime,\flat, \circ} \to T^{\flat, \circ}$ is an almost isomorphism. 

    Then \cite[Prop.~11.10]{Bhatt-Scholze-Grassmannian} implies that $R^{\prime,\flat, \circ} \otimes_{R^{\flat, \circ}} R^{\prime\prime,\flat, \circ} \simeq R^{\prime,\flat, \circ} \otimes^\rL_{R^{\flat, \circ}} R^{\prime\prime,\flat, \circ}$ is a perfect ring. Therefore, it has bounded $\varpi^\flat$-torsion, so its derived $\varpi^\flat$-adic completion coincides with the classical $\varpi^\flat$-adic completion and we conclude that 
    \[
    R^{\prime,\flat, \circ} \wdh{\otimes}^{\rL}_{R^{\flat, \circ}} R^{\prime\prime,\flat, \circ} \simeq R^{\prime,\flat, \circ} \wdh{\otimes}_{R^{\flat, \circ}} R^{\prime\prime,\flat, \circ}.
    \]
    Finally, the proof of \cite[Prop.~6.18]{Scholze-perfectoid}\footnote{\cite[Prop.~6.18]{Scholze-perfectoid}  assumes that $S$ is defined over a perfectoid field, but this is not necessary to run the argument.} implies that the natural map $R^{\prime,\flat, \circ} \wdh{\otimes}_{R^{\flat, \circ}} R^{\prime\prime,\flat, \circ} \to T^{\flat, \circ}$ is an almost isomorphism. This finishes the proof. 
\end{proof}
Recall that for any complete Huber pair $(A, A^+)$, the condensed abelian groups $\ud{A}$ and $\ud{A}^+$ are both solid by virtue of \cite[Lem.~4.5]{Andreychev}.

\begin{lemma}
\label{lemma:rational-kunneth-for-perfectoids} 
Let $S' \to S$ and $S'' \to S$ be maps of affinoid perfectoid spaces. Then the natural map $\ud{\O(S')} \otimes^{\rL \Box}_{\ud{\O(S)}} \ud{\O(S'')} \to \ud{\O(S' \times_S S'')}$ is an isomorphism. 
\end{lemma}
\begin{proof}
    Since (derived) solid tensor products commute with all colimits, it suffices to show that the natural map $\ud{\O^+(S')} \otimes^{\rL \Box}_{\ud{\O^+(S)}} \ud{\O^+(S'')} \to \ud{\O^+(S' \times_S S'')}$ is an almost isomorphism. We choose a pseudo-uniformizer $\varpi \in \O^+(S)$. Now \cite[Prop.~A.3]{bosco-phd} implies that $\ud{\O^+(S')} \otimes^{\rL \Box}_{\ud{\O^+(S)}} \ud{\O^+(S'')}$ is derived $\varpi$-adically complete. Therefore, the natural map
    \[
    \ud{\O^+(S')} \otimes^{\rL \Box}_{\ud{\O^+(S)}} \ud{\O^+(S'')} \to \ud{\O^+(S') \wdh{\otimes}^{\rL}_{\O^+(S)} \O^+(S'')} 
    \]
    is an isomorphism because it is an isomorphism modulo $\varpi^n$ for any $n$.
    Hence, \cref{lemma:almost-integral-kunneth-formula} implies that the map $\ud{\O^+(S')} \otimes^{\rL \Box}_{\ud{\O^+(S)}} \ud{\O^+(S'')} \to \ud{\O^+(S' \times_S S'')}$ is an almost isomorphism. This finishes the proof. 
\end{proof}

\begin{proposition} 
\label{Kunneth for BI rings}
Let $S' \to S$ and $S'' \to S$ be maps of affinoid perfectoid spaces over $\Spa(\Q_p, \Z_p)$. Then the natural morphism 
\[
\ud{\BB_{I}(S')} \otimes^{\rL \Box}_{\ud{\BB_{I}(S)}} \ud{\BB_{I}(S'')} \to \ud{\BB_{I}(S'\times_S S'')}
\]
is an isomorphism. 
\end{proposition}
\begin{proof}
    First, we note that \cite[Prop.~A.68]{bosco-master} and \cref{lemma:kunneth-for-FF-balls}
    imply that the natural morphism $\ud{\BB_{I}(S')} \otimes^{\Box}_{\ud{\BB_{I}(S)}} \ud{\BB_{I}(S'')} \to \ud{\BB_{I}(S'\times_S S'')}$ is an isomorphism. 
    Therefore, it only remains to show that $\ud{\BB_{I}(S')} \otimes^{\rL \Box}_{\ud{\BB_{I}(S)}} \ud{\BB_{I}(S'')}$ is concentrated in degree $0$. 
    
    Fix a Kummer field $\QQ_{p, \infty}$. Then the proof of \cite[Prop.~II.1.1]{Fargues-Scholze} implies that $\widetilde{Y}_{S, I} \coloneqq Y_{S, I} \times_{\Spa(\Q_p, \Z_p)} \Spa(\Q_{p, \infty}, \Z_{p,\infty})$ is an affinoid perfectoid space. Furthermore, \cite[Th.~A.43(ii)]{bosco-master} implies that $\ud{\Q}_{p, \infty}$ is a flat $\ud{\Q}_p$-module (with respect to the solid tensor product). Since $\ud{\Q}_p \to \ud{\Q}_{p, \infty}$ admits a section as condensed $\ud{\Q}_p$-modules, it suffices to show that 
    \begin{align*}
    \ud{\BB_{I}(S')} \otimes^{\rL \Box}_{\ud{\BB_{I}(S)}} \ud{\BB_{I}(S'')} \otimes_{\ud{\Q}_p}^{\Box} \ud{\Q}_{p, \infty} & \simeq \big(\ud{\BB_{I}(S')} \otimes^{\Box}_{\ud{\QQ}_p} \ud{\QQ}_{p, \infty} \big) \otimes^{\rL \Box}_{\big(\ud{\BB_{I}(S)} \otimes^{\Box}_{\ud{\QQ}_p} \ud{\QQ}_{p, \infty} \big)} \big(\ud{\BB_{I}(S'')} \otimes^{\Box}_{\ud{\QQ}_p} \ud{\QQ}_{p, \infty} \big) \\
    & \simeq \ud{\O(\widetilde{Y}_{S', I})} \otimes^{\rL \Box}_{\ud{\O(\widetilde{Y}_{S, I})}} \ud{\O(\widetilde{Y}_{S'', I})}
    \end{align*}
    is concentrated in degree $0$. This follows from \cref{lemma:rational-kunneth-for-perfectoids}. 
\end{proof}

\begin{proposition} 
\label{pullback for BI rings}
Let $J \subset I \subset (0, \infty)$ be an inclusion of closed intervals
with rational endpoints,
and let $f \colon S' \to S$ be a map of affinoid perfectoid spaces over $\Spa(\Q_p, \Z_p)$.
Then the natural morphisms
\[
\ud{\BB_{I}(S')} \otimes^{\rL \Box}_{f^*, \ud{\BB_{I}(S)}, \res} \ud{\BB_{J}(S)} 
\xrightarrow{\res \otimes f^*} \ud{\BB_{J}(S')}
\quad \text{and} \quad
\ud{\BB_{I}(S')} \otimes^{\rL \Box}_{f^*, \ud{\BB_{I}(S)}, \varphi} \ud{\BB_{I/p}(S)} 
\xrightarrow{\varphi \otimes f^*} \ud{\BB_{I/p}(S')}
\]
induced by the maps in \cref{construction:compatibility-maps} are isomorphisms.
\end{proposition}

\begin{proof}
The construction of the maps in \cref{construction:compatibility-maps} implies that
\[
Y_{S', J} \xrightarrow[\sim]{(Y_{f}, \res)} Y_{S, J} \times_{\res, Y_{S, I}} Y_{S', I}
\quad \text{and} \quad
Y_{S', I/p} \xrightarrow[\sim]{(Y_f, \varphi)} Y_{S, I/p} \times_{\varphi, Y_{S, I}} Y_{S', I}.
\]
Now the proof is completely analogous to that of \Cref{Kunneth for BI rings}, using the
fiber product formulas above in place of \Cref{lemma:kunneth-for-FF-balls}.
\end{proof} 

\begin{construction}\label{construction:theta-rings} Let $S=\Spa(R, R^+)$ be an affinoid perfectoid space over $\Spa(\QQ_p, \ZZ_p)$, let $p^\flat\in R^{\flat, +}$ and $u\in (R^+)^\times$ be as in \cref{notation:p-flat} (so $(p^\flat)^\sharp = p\cdot u$), and let $I=[a,b]\subset (0, \infty)$ be a closed interval with rational endpoints such that $1\in I$. 
\begin{enumerate}[leftmargin=*]
    \item There is a natural surjective continuous ring homomorphism
\[
\theta \colon \mathbf{A}_{\inf}(S) \twoheadrightarrow \cO^+(S)
\]
sending $\sum_{i=0}^n [r_i]p^i \in \AA_{\inf}(S)$ to $\sum_{i=0}^n r_i^\sharp p^i\in R^+$. We refer to \cite[discussion after Lem.~3.3]{BMS1} for the details. We denote by $\iota_\infty \colon S \to \Spa\bigl(\AA_{\inf}(S), \AA_{\inf}(S)\bigr)$ the induced map of (pre-)adic spaces. 
    \item By construction, we see that 
    \[
        v\bigl(\theta(p)\bigr)=v(p) =v(p\cdot u) = v\bigl((p^\flat)^\sharp\bigr)=v\bigl(\theta([p^\flat])\bigr) \quad \text{and} \quad v(p)\neq 0
    \]
    for any $v\in \Spa(R, R^+)$. Since $a\leq 1\leq b$ by our assumption, the morphism $\iota_\infty$ factors uniquely through $Y_{S, I}\subset Y_S$, and thus defines the morphism 
    \[
    \iota_{\infty, I} \colon S \to Y_{S, I}. 
    \]
    We denote the induced morphism on Huber pairs by 
    \[
    \theta_I\colon \bigl(\BB_I(S), \AA_I(S)\bigr) = \bigl(\O(Y_{S, I}), \O^+(Y_{S, I}) \bigr) \to \bigl(\O(S), \O^+(S) \bigr) = (R, R^+).
    \]
\end{enumerate}
\end{construction}
By \cite[Prop.~II.1.4]{Fargues-Scholze}, the map $\iota_{\infty, I}\colon  S \to Y_{S, I}$ is a closed immersion 
which realizes $S$ as a closed Cartier divisor inside $Y_{S, I}$.
Below, we give a variant of this claim with a simpler proof. 
\begin{proposition}
\label{base change along theta}
Assume one of the following:
\begin{enumerate}[label=\upshape{(\arabic*)}]
\item $I = [a, 1]$ and there exists a pair of elements $\alpha, \beta = p^{\flat} \in R^{\flat, +}$ 
well-adapted to $p^\flat$ and $I$; or
\item $I = [1, b]$ and there exists a pair of elements $\alpha = p^{\flat}, \beta \in R^{\flat, +}$ 
well-adapted to $p^\flat$ and $I$.
\end{enumerate}
Then the map $\theta_I \colon \AA_I(S) \to \cO^+(S)$
is surjective with kernel generated by a nonzerodivisor $\xi'$.
Moreover, for any $S' \to S$, the element $\xi'$ remains a nonzerodivisor generating $\ker\bigl(\AA_I(S') \to \cO^+(S')\bigr)$.
\end{proposition}

\begin{proof}
Let us only prove the statement for the case $I = [a, 1]$, the proof for the other case being similar.
Let $\alpha, \beta = p^{\flat} \in R^{\flat, +}$ be a pair of elements well-adapted to $p^\flat$ and $I$.
By the condition on $p^{\flat}$ in \Cref{notation:p-flat}, we may choose a sequence of elements $\{a_n \in R^{\flat, +}\}_{n \geq 0}$ 
such that $a_0 \in (R^{\flat,+})^{\times}$ and
$p = (p^{\flat})^{\sharp}\cdot \sum_{n \geq 0} p^n (a_n)^{\sharp} \in R^{+}$.
Then the element $\xi' \coloneqq \frac{p}{[p^\flat]} - \sum_{n \geq 0} p^n [a_n]$ lies in $\ker(\theta_I)$.
Our statements now all follow from the claim that the following sequence is short exact:
\[
0 \to \AA_I(S) \xrightarrow{\cdot \xi'} \AA_I(S) \to \cO^+(S) \to 0
\]

Let us justify this claim. Since all terms involved are $p$-torsionfree (by definition)
and $p$-complete (for $\AA_I(S)$ use \Cref{lemma:semi-explicit-A_I}), this is equivalent to showing
that the sequence is exact after reducing modulo $p$.
By \Cref{cor:AI-mod-p-rings}, modulo $p$ the sequence becomes
\[
0 \to \frac{R^{\flat, +}[x, y]}{(p^{\flat}y, xy - \alpha/p^{\flat})} \xrightarrow{\cdot(y - a_0)}
\frac{R^{\flat, +}[x, y]}{(p^{\flat}y, xy - \alpha/p^{\flat})} 
\xrightarrow[y \mapsto a_0]{x \mapsto a_0^{-1} \cdot \alpha/p^{\flat}} R^{\flat, +}/p^{\flat} \to 0.
\]
First, let us show that the kernel of the first map is indeed $0$. 
Any element in 
\[
\frac{R^{\flat, +}[x, y]}{(p^{\flat}y, xy - \alpha/p^{\flat})} \simeq \bigl(\bigoplus_{i\geq 0} R^{\flat, +}/\alpha \bigr) \bigoplus \bigl(\bigoplus_{i\geq 1} R^{\flat, +}/p^\flat \bigr)
\]
can be written uniquely
as a finite sum $r_n y^n + r_{n-1} y^{n-1} + \cdots + r_1 y + c + s_1 x + \cdots + s_{m-1} x^{m-1} + s_m x^m$
for some $r_i \in R^{\flat, +}/p^\flat$  and $c, s_k \in R^{\flat, +}/\alpha$. The claim that $y - a_0$ is a nonzerodivisor
then follows from the fact that it is monic and $a_0 \in (R^{\flat, +})^{\times}$.
To see that the cokernel is as described, we compute
\[
\frac{R^{\flat, +}[x, y]}{(p^{\flat}y, xy - \alpha/p^{\flat}, y - a_0)} = 
\frac{R^{\flat, +}[x]}{(p^{\flat}a_0, a_0 \cdot x - \alpha/p^{\flat})} = R^{\flat, +}/p^{\flat}. \qedhere
\]
\end{proof}

\begin{cor}\label{cor:base change along theta} Let $f \colon S' \to S$ be a map of affinoid perfectoid spaces over $\Spa(\Q_p, \Z_p)$ and let $I\subset (0, \infty)$ be a closed interval with rational endpoints such that $1 \in I$. Then the natural morphism 
    \[
        \ud{\BB_{I}(S')} \otimes^{\rL \Box}_{f^*, \ud{\BB_{I}(S)}, \theta_I} \ud{\cO(S)} 
        \xrightarrow{\theta_I \otimes f^*} \ud{\cO(S')}
    \]
is an isomorphism.
\end{cor}
\begin{proof}
    First, note that the morphism $\theta_I\colon \BB_I(S) \to \cO(S)$ factors as the composition $\BB_I(S) \xr{\mathrm{res}} \BB_{[1,1]}(S) \xr{\theta_{[1,1]}} \cO(S)$. Thus, \cref{pullback for BI rings} ensures that 
    \[
    \ud{\BB_I(S')} \otimes^{\rL \Box}_{f^*, \ud{\BB_I(S)}, \theta_I} \ud{\cO(S)} \simeq \ud{\BB_I(S')} \otimes^{\rL \Box}_{f^*, \ud{\BB_I(S)}, \mathrm{res}} \ud{\BB_{[1, 1]}(S)} \otimes^{\rL \Box}_{\id, \ud{\BB_{[1,1]}(S)}, \theta_{[1,1]}} \ud{\cO(S)} \simeq \ud{\BB_{[1,1]}(S')} \otimes^{\rL \Box}_{f^*, \ud{\BB_{[1,1]}(S)}, \theta_{[1,1]}} \ud{\cO(S)},
    \]
    so we can assume that $I=[1,1]$. In this case, the result follows immediately from \cref{base change along theta} since 
    \[
     \ud{\BB_{[1,1]}(S')} \otimes^{\rL \Box}_{f^*, \ud{\BB_{[1,1]}(S)}, \theta_{[1,1]}} \ud{\cO(S)} \simeq \ud{\BB_{[1,1]}(S')} \otimes^{\rL, \Box}_{\ud{\BB_{[1,1]}(S)}} \ud{\BB_{[1,1]}(S)/(\xi')} \simeq \ud{\BB_{[1,1]}(S')/(\xi')} \simeq \ud{\cO(S')}. \qedhere
    \]
\end{proof}

\subsection{Permanence properties of period rings}\label{subsection:rings}

Throughout this subsection, fix a closed interval $I=[a,b] \subset (0, \infty)$ with rational endpoints. We show that the functor $S\mapsto Y_{S, I}$ preserves open immersions, \'etale morphisms, finite \'etale morphisms, and $v$-coverings.

By \cite[Lem.~15.1(ii)]{diamonds}, for a perfectoid space $T$ of characteristic $p$, a map $T \to \Spd(\Q_p)$ corresponds to a choice of an untilt $T^\sharp$ over $\Spa(\Q_p, \Z_p)$. 
Next, we recall the definition of structure $v$-sheaves.
\begin{definition}[{cf.\ \cite[Th.~8.7]{diamonds}}]
    A diamond $X$ over $\Spd(\QQ_p, \Z_p)$ is equipped with the following sheaves for the $v$-topology:
    \begin{enumerate}
        \item $\cO^+_X$, given by $\cO^+_X(T) \colonequals \cO^+_{T^\sharp}(T^\sharp)$ for any $T\in \Perf_{/X}$;
        \item $\cO_X$, given by $\cO_X(T) \colonequals \cO_{T^\sharp}(T^\sharp)$ for any $T\in \Perf_{/X}$.
    \end{enumerate}
\end{definition}

\begin{definition}[{cf.\ \cite[Def.\,9.6]{Hansen-Kedlaya}}]\label{defn:v-completion} The \emph{$v$-completion} of a Tate--Huber pair $(A, A^+)$ over $(\QQ_p, \ZZ_p)$ is the pair $(\check{A}, \check{A}^+)\coloneqq \bigl(\Hh^0_v(\Spd(A, A^+), \O), \Hh^0_v(\Spd(A, A^+), \O^+)\bigr)$. 
We topologize $\check{A}^+$ via the $p$-adic topology and $\check{A}$ by requiring $\check{A}^+$ to be open and bounded in $\check{A}$.
\end{definition}

\begin{lemma}\label{lemma:v-completion-complete} Let $(A, A^+)$ be a Tate--Huber pair over $(\QQ_p, \ZZ_p)$. Then $\check{A}^+$ is $p$-adically complete and $(\check{A}, \check{A}^+)$ is a complete uniform Tate--Huber pair.  
\end{lemma}
\begin{proof}
    Since affinoid perfectoid spaces form a basis of the $v$-topology, we can find a $1$-truncated $v$-hypercovering $\Spd(S, S^+) \rightrightarrows \Spd(R, R^+) \rightarrow \Spd(A, A^+)$ for some perfectoid pairs $(R, R^+)$ and $(S, S^+)$ over $(\QQ_p, \ZZ_p)$. We denote by $\pi_1^*, \pi_2^*\colon R^+ \to S^+$ the maps corresponding to the two maps $\Spd(S, S^+) \rightrightarrows \Spd(R, R^+)$.
    Then the sheaf axiom implies that  
    \[
    \check{A}^+ = \Hh^0_v(\Spd(A, A^+), \O^+) \simeq  \ker\bigl(R^+ \xr{\pi_1^* - \pi_2^*} S^+\bigr). 
    \]
    Since $R^+$ and $S^+$ are $p$-adically complete and $p$-torsionfree, $\check{A}^+$ is $p$-adically complete and $p$-torsionfree as well. Furthermore, using that $R^+$ is integrally closed in $R \simeq R^+[\frac{1}{p}]$ and $S^+$ is $p$-torsionfree, we deduce that $\check{A}^+$ is integrally closed in $\check{A}\simeq \check{A}^+[\frac{1}{p}]$. 

    Thus, in order to conclude that $(\check{A}, \check{A}^+)$ is a complete Tate--Huber pair, it suffices to observe $\check{A}^+$ is open in $\check{A}$ by definition. Finally, \cite[Cor.\,1.4(ii)]{H0} and \cite[Lem.\,C.1.3]{zav-almost} imply that $\check{A}$ is uniform. 
\end{proof}

\begin{definition}[{cf.\ \cite[Def.\,9.6]{Hansen-Kedlaya}}]\label{defn:v-complete}
A Tate--Huber pair $(A, A^+)$ over $(\QQ_p, \ZZ_p)$ is {\it $v$-complete} if the natural morphism $(A, A^+) \to (\check{A}, \check{A}^+)$ is a topological isomorphism. 
\end{definition}

\begin{rmk} The discussion in the second paragraph of \cite[Def.\,9.1]{Hansen-Kedlaya} implies that \cref{defn:v-complete} coincides with \cite[Def.\,9.6]{Hansen-Kedlaya} for any pair Tate--Huber pair $(A, A^+)$ over $(\QQ_p, \ZZ_p)$.
\end{rmk}

Our next goal is to show that the pair $\bigl(\BB_I(S), \AA_I(S)\bigr)$ is $v$-complete for any $S\in \Perfd^\aff_{\QQ_p}$. For this, we will need some preliminary lemmas. 

\begin{lemma}\label{lemma:criterion-for-iso} Let $f \colon (A, A^+)\to (B, B^+)$ be a morphism of complete Tate--Huber pairs such that $A \to B$ is an isomorphism of rings and the induced morphism $\abs{f} \colon \abs{\Spa(B, B^+)} \to \abs{\Spa(A, A^+)}$ is a bijection of adic spectra. Then $f^\sharp\colon (A, A^+) \to (B, B^+)$ is an isomorphism of Tate--Huber pairs.
\end{lemma}
\begin{proof}
    First, \cite[Lem.\,2.4]{Huber-generalization} implies that $A \to B$ is an isomorphism of topological rings.
    Therefore, we only need to show that the natural morphism $A^+ \to B^+$ is an isomorphism of rings. For this, we note that \cite[Prop.\,1.6]{Huber-generalization} and the assumption that $\abs{\Spa(B, B^+)} \to \abs{\Spa(A, A^+)}$ and $A \to B$ are bijections imply that 
    \[
    A^+ = \bigl\{f\in A  \suchthat \abs{f(x)}\leq 1 \text{ for all } x\in \abs{\Spa(A, A^+)} \bigr\} \simeq \bigl\{f\in B \suchthat \abs{f(x)}\leq 1 \text{ for all } x\in \abs{\Spa(B, B^+)} \bigr\} = B^+. 
    \]
    This finishes the proof. 
\end{proof}

For the next lemma, we recall that a morphism of complete Tate--Huber pairs $(A, A^+) \to (B, B^+)$ is finite \'etale if $B$ is a finite \'etale $A$-algebra and $B^+$ is the integral closure of $A^+$ in $B$.

\begin{lemma}\label{lemma:finite-etale-after-v-completion} Let $(A, A^+) \to (B, B^+)$ be a morphism of Tate--Huber pairs over $(\QQ_p, \ZZ_p)$. If the associated morphism of diamonds $\Spd(B, B^+) \to \Spd(A, A^+)$ is finite \'etale, then $(\check{A}, \check{A}^+) \to (\check{B}, \check{B}^+)$ is a finite \'etale morphism of complete Tate--Huber pairs.
\end{lemma}
\begin{proof}
    Without loss of generality, we can assume that $(A, A^+)$ and $(B, B^+)$ are complete Tate--Huber pairs as the passage to the completion does not change the associated diamonds (and, thus, the associated $v$-completions). Now we observe that \cite[Lem.\,15.6]{diamonds}, \cite[8.2.17]{Kedlaya-Liu-1}, and \cite[Lem.\,C.1.1]{zav-almost} imply that there is a finite \'etale morphism $(A, A^+) \to (B', B^{\prime,+})$ such that there is an isomorphism $\Spd(B, B^+)\simeq \Spd(B', B^{\prime,+})$ of diamonds over $\Spd(A, A^+)$. Therefore, we can replace $(B, B^+)$ by $(B', B^{\prime,+})$ to assume that $(A, A^+) \to (B, B^+)$ is a finite \'etale morphism of complete Tate--Huber pairs. 
    
    Then \cite[Lem.\,B.3.5]{Z-quotients} and \cite[Lem.\,C.1.1]{zav-almost} imply that $\check{A} \otimes_A B \simeq \check{A} \wdh{\otimes}_A B$. We denote by $(\check{A}^+ \otimes_{A^+} B^+)^+$ the integral closure of $\rm{Im}(\check{A}^+ \otimes_{A^+} B^+ \to \check{A} \otimes_A B)$ inside $\check{A} \otimes_A B$.
    It suffices to show that the natural morphism
    \[
    (\check{A}_B, \check{A}_B^+) \colonequals \bigl(\check{A} \otimes_A B, (\check{A}^+ \otimes_{A^+} B^+)^+\bigr) \to (\check{B}, \check{B}^+)
    \]
    is an isomorphism. Now \cite[Lem.\,9.8(a)]{Hansen-Kedlaya} guarantees that $\check{A}_B \to \check{B}$ is an isomorphism. 
    Combined with \cref{lemma:v-completion-complete} and \cref{lemma:criterion-for-iso}, it remains to prove that the induced morphism $\abs{\Spa(\check{B}, \check{B}^+)} \to \abs{\Spa(\check{A}_B, \check{A}_B^+)}$ is a bijection. Furthermore, \cite[Lem.\,15.6]{diamonds} ensures that it suffices to show that $\Spd(\check{B}, \check{B}^+) \to \Spd(\check{A}_B, \check{A}_B^+)$ is an isomorphism. This follows from the observation that 
    \[
    \Spd(\check{A}_B, \check{A}_B^+) \simeq \Spd(\check{A}, \check{A}^+) \times_{\Spd(A, A^+)} \Spd(B, B^+)
    \]
    and the fact that $\Spd(\check{A}, \check{A}^+) \to \Spd(A, A^+)$ and $\Spd(\check{B}, \check{B}^+) \to \Spd(B, B^+)$ are isomorphisms (see \cite[Lem.\,9.8(c)]{Hansen-Kedlaya}). 
\end{proof}

\begin{cor}\label{cor:Fargues--Fontaine-v-complete} Let $S$ be an affinoid perfectoid space over $\Spa(\QQ_p, \ZZ_p)$. Then the Huber pair $\big(\BB_I(S), \AA_I(S)\big)$ is $v$-complete. 
\end{cor}
\begin{proof}
    First, \cite[Lem.\,9.8(c)]{Hansen-Kedlaya} and \cref{lemma:criterion-for-iso} imply that it suffices to show that the natural morphism 
    \[
    \BB_I(S) \to \Hh^0_v(Y_{S, I}, \O)
    \]
    is an isomorphism. For this, we fix a choice of a Kummer field $\QQ_{p,\infty}$ and denote by $G$ the absolute Galois group of $\QQ_p$.

    Then $\BB_I(S) \wdh{\otimes}_{\QQ_p} \QQ_{p, \infty}$ is a perfectoid Tate algebra (see \cite[Prop.~II.1.1]{Fargues-Scholze}), so \cite[Prop.\,6.18]{Scholze-perfectoid} guarantees that 
    \[
    \BB_I(S) \wdh{\otimes}_{\QQ_p} \CC_p \simeq \BB_I(S)\wdh{\otimes}_{\QQ_p} \QQ_{p, \infty}\wdh{\otimes}_{\QQ_{p, \infty}} \CC_p
    \]
    is also a perfectoid Tate algebra. In particular, the adic space $\widetilde{Y}_{S, I}\coloneqq Y_{S, I}\times_{\Spa(\QQ_p, \ZZ_p)} \Spa(\CC_p, \O_{\CC_p})$ is affinoid perfectoid.
    Since $\widetilde{Y}_{S, I}$ is an affinoid perfectoid space and $\widetilde{Y}_{S, I}^\diam \to Y_{S, I}^\diam$ is a $\ud{G}$-torsor, we conclude that $\Hh^0_v(Y_{S, I}, \O) = \Hh^0_v(\widetilde{Y}_{S, I}, \O)^G = \big(\BB_I(S) \wdh{\otimes}_{\QQ_p} \CC_p\big)^G$. Thus, the question boils down to showing that the natural morphism
    \[
    \BB_I(S) \to \big(\BB_I(S) \wdh{\otimes}_{\QQ_p} \CC_p\big)^G
    \]
    is an isomorphism. Now we note that \cref{rmk:BI-AI-Tate-Huber-pair} ensures that $\BB_I(S)$ is a Banach module over $\QQ_p$. Therefore, \cite[Th.\,2.5.4]{PGS} implies that $\BB_I(S) \simeq \wdh{\bigoplus}_J \QQ_p$ for some set $J$. So we reduce the question to showing that the natural map $\wdh{\bigoplus}_J \QQ_p \to \bigl(\wdh{\bigoplus}_J \CC_p\bigr)^G$ is an isomorphism. Tate's calculation (see \cite[Th.\,1 on p.176]{p-divisible}) implies that $\QQ_p \to (\CC_p)^G$ is an isomorphism. Hence, the result follows from the explicit calculation that $\wdh{\bigoplus}_J \bigl( \CC_p^G \bigr) \to \bigl(\wdh{\bigoplus}_J \CC_p\bigr)^G$ is an isomorphism. 
\end{proof}

\begin{lemma}\label{lemma:permanence-properties-Y}
    Let $S' \to S$ be a morphism of affinoid perfectoid spaces over $\Spa(\QQ_p,\ZZ_p)$. If $S' \to S$ is an open immersion (resp.\ finite \'etale map, resp.\ \'etale map, resp. surjective), then so is $Y_{S', I} \to Y_{S, I}$.
\end{lemma}
\begin{proof}
    The case of open immersions follows from \cite[Prop.\,II.1.3]{Fargues-Scholze}. Now suppose that $S' \to S$ is a surjective map. Then \cite[Prop.\,11.15 and Lem.\,15.6]{diamonds} imply that, for the purpose of showing that $Y_{S', I} \to Y_{S, I}$ is a surjective map, it suffices to show that $Y^\diam_{S', I} \to Y^\diam_{S, I}$ is an epimorphism in the $v$-topology. This follows from \cite[Prop.\,II.1.17]{Fargues-Scholze} and the fact that epimorphisms in any topos are preserved under base change \cite[Claim a) in the proof of Exp.~II, Prop.~4.3.2)]{SGA4}.

    Lastly, we recall that, locally on the source and target, any \'etale map can be written as an open immersion followed by a finite \'etale map. Therefore, the only thing we are left to show is that $Y_{S', I} \to Y_{S, I}$ is a finite \'etale map if $S' \to S$ is finite \'etale. This follows immediately from \cref{cor:Fargues--Fontaine-v-complete} and \cref{lemma:finite-etale-after-v-completion}. 
\end{proof}

\cref{lemma:kunneth-for-FF-balls} and \cref{lemma:permanence-properties-Y} (applied to open immersions) ensure that we can extend the construction of $Y_S$ and $Y_{S, I}$ to general perfectoid spaces $S$ over $\Spa(\Q_p, \Z_p)$ via gluing; see \cite[p.\,49]{Fargues-Scholze} for more details.

\begin{defn}\label{defn:period-rings-2} For a perfectoid space $S$ over $\Spa(\Q_p, \Z_p)$, the \emph{period rings} $\AA_I(S)$ and $\BB_I(S)$ are
    \[ \AA_I(S) \colonequals \Hh^0_\an(Y_{S,I},\cO^+) \quad \text{and} \quad \BB_I(S) \colonequals \Hh^0_\an(Y_{S,I},\cO). \]
\end{defn}

\subsection{Period rings are \texorpdfstring{$v$}{v}-sheaves}

Throughout this subsection, we fix a closed interval $I=[a,b] \subset (0, \infty)$ with rational endpoints. The main goal of this subsection is to  show that the functors 
\[
\AA_I(\blank), \BB_I(\blank) \colon \Perfd^{\aff,\op}_{/\QQ_p} \to \Ring_{\rm{top}}
\]
are $v$-sheaves of topological rings. We begin with the following preliminary lemma:

\begin{lemma}
\label{lemma:almost-acyclic} 
Let $S=\Spa(R, R^+)$ be an affinoid perfectoid space over $\Spa(\QQ_p, \ZZ_p)$, let $p^\flat\in R^{\flat, +}$ be an element as in \cref{notation:p-flat}, and let $\alpha, \beta\in R^{\flat, +}$ 
be a pair of elements well-adapted to $p^\flat$ and $I$. Let $J$ be a finite set and let $\{S_j \to S\}_{j\in J}$ 
be a $v$-covering in $\Perfd^\aff_{/\QQ_p}$.
Then the complex of $\AA_I(S)$-modules
    \[
    \AA_I(S) \to \prod_{j\in J} \AA_I(S_j) \to \prod_{j, j' \in J} \AA_I(S_{j}\times_S S_{j'}) \to \dots
    \]
    is almost exact with respect to the ideal of almost mathematics $\AA_I^{\circ\circ}(S) = \cup_{n\geq 1} [(p^\flat)^{1/p^n}] \AA_I(S)$. 
\end{lemma}

\cref{cor:almost-mathematics-AI} ensures that $\AA_I^{\circ\circ}(S)$ is indeed an ideal of almost mathematics. 
\begin{proof}
    First, \cite[Lem.\,A.5]{zav-almost} and the fact that $\AA_I(T)$ is $p$-torsionfree for any $T\in \Perfd^{\aff}_{/\QQ_p}$ (\cref{cor:almost-mathematics-AI}) imply that it suffices to show that the complex of $\AA_I(S)/p$-modules
    \[
    \AA_I(S)/p \to \prod_{j\in J} \AA_I(S_j)/p \to \prod_{j, j' \in J} \AA_I(S_{j}\times_S S_{j'})/p \to \dotsb
    \]
    is almost exact with respect to the ideal $\cup_{n\geq 1} (p^\flat)^{1/p^n} \AA_I(S)/p$. 
    This follows directly from \cref{cor:AI-mod-p-rings} and \cite[Prop.\,8.8]{diamonds}. 
\end{proof}

\begin{cor}[{cf.\,\cite[Prop.\,II.2.1]{Fargues-Scholze}}]\label{cor:AI-etale-sheaf}
    The presheaves $\AA_I(\blank), \BB_I(\blank) \colon \Perfd^{\aff, \op}_{/\QQ_p} \to \mathrm{Ring}$ are $v$-sheaves of topological rings.
\end{cor}
\begin{proof}
    \begin{enumerate}[wide,label={\textit{Step}~\arabic*.},ref={Step~\arabic*}]
    \item\label{BI-etale-sheaf} \textit{The presheaf $\BB_I(\blank)$ is an \'etale sheaf.} 
    Let $\{S_j \to S\}_{j\in J}$ be an \'etale covering in $\Perfd^\aff_{/\QQ_p}$. Then $\{Y_{S_j, I} \to Y_{S, I}\}_{j\in J}$ is an \'etale covering as well by virtue of \cref{lemma:permanence-properties-Y}. Furthermore, \cref{lemma:kunneth-for-FF-balls} ensures that $Y_{S_j\times_{S} S_{j'}, I} \simeq Y_{S_j, I} \times_{Y_{S, I}} Y_{S_{j'}, I}$ for any $j,j'\in J$. Therefore, the question boils down to showing that $\O_{(Y_{S, I})_\et}$ is an \'etale sheaf. This follows directly from the fact that $Y_{S, I}$ is an adic space (see \cite[Prop.\,II.1.1]{Fargues-Scholze}) and  \cite[Th.\,8.2.2(c)]{Kedlaya-Liu-1}.

    \item\label{BI-v-sheaf} \textit{The presheaf $\BB_I(\blank)$ is a $v$-sheaf.}
    Let $\{S_j \to S\}_{j\in J}$ be a $v$-covering. We wish to show that the sequence
    \begin{equation}\label{eqn:sheaf-B_I-complex}
    0 \to \BB_I(S) \to \prod_{j\in J} \BB_I(S_j) \to \prod_{j, j'\in J} \BB_I(S_j \times_S S_{j'})
    \end{equation}
    is exact. Since $S$ is quasi-compact, we can assume that $J$ is a finite set. Then \cref{BI-etale-sheaf} ensures that $\BB_I(\blank)$ is an \'etale sheaf. Thus, we can check that \cref{eqn:sheaf-B_I-complex} is exact \'etale locally on $S$.
    By \cref{rmk:often-well-adapted-pairs}, we can therefore assume that there exists a pair $\alpha, \beta\in \O_S^+(S)^\flat$ which is well-adapted to $p^\flat\in \O_S^+(S)^\flat$ and $I$. In this case, \cref{lemma:almost-acyclic} immediately implies exactness of \cref{eqn:sheaf-B_I-complex} because $\BB_I(\blank)=\AA_I(\blank)[\frac{1}{p}]$ and $p\in \AA_I^{\circ\circ}(S)$. 

    \item \textit{The presheaf $\AA_I(\blank)$ is a $v$-sheaf.}
    Since $S$ is quasi-compact and $\AA_I(S_1\sqcup S_2) \simeq \AA_I(S_1)\times \AA_I(S_2)$, it suffices to show that for every $v$-surjection $S' \to S$, the complex $0 \to \AA_I(S) \to \AA_I(S') \to \AA_I(S'\times_S S')$ is exact. This follows from the fact that $\BB_I(\blank)$ is a sheaf (by \cref{BI-v-sheaf}) together with the observations that $Y_{S', I} \to Y_{S, I}$ is surjective for any $v$-covering $S' \to S$ (see \cref{lemma:permanence-properties-Y}) and that 
    \[ 
    \AA_I(T) = \Hh^0_\an(Y_{T,I},\cO^+) = \{ f \in \Hh^0_\an(Y_{T,I},\cO) \suchthat \abs{f(x)}\leq 1 \text{ for all } x \in Y_{T, I} \}
    \]
    for any $T\in \Perfd^\aff_{/\QQ_p}$. 
    
    \item \textit{The $v$-sheaves $\AA_I(\blank)$ and $\BB_I(\blank)$ are sheaves of topological rings.}
    It suffices to show that, for any $v$-covering $\{S_j \to S\}_{j\in J}$, the morphisms $\AA_I(S) \to \prod_{j\in J} \AA_I(S_j)$ and $\BB_I(S) \to \prod_{j\in J} \BB_I(S_j)$ are homeomorphisms onto their images. Since $S$ is quasi-compact, we can assume that $J$ is a finite set. In this case, the result follows immediately from \cite[Lem.\,2.4]{Huber-generalization}. \qedhere 
    \end{enumerate}
\end{proof}

\section{Period Sheaves}
\label{section:period-sheaves}

The main goal of this construction is to globalize the notion of period rings from \cref{section:period-rings}, introduce the notion of period sheaves, and study properties of perfect complexes over these period sheaves. 

Throughout this section, we will freely use the language of diamonds as developed in \cite{diamonds}. In particular, we will freely use that, for a perfectoid space $S$ of characteristic $p$, a map $S \to \Spd(\Q_p, \ZZ_p)$ is equivalent to a choice of an untilt $S^\sharp$ over $\Spa(\Q_p, \Z_p)$ (see \cite[\textsection 15]{diamonds}).

Throughout this subsection, we fix a closed interval $I=[a,b] \subset (0, \infty)$ with rational endpoints.

\subsection{Basic properties of period sheaves}

In this subsection, we introduce the period sheaves that will be important for us throughout this paper and study their basic properties.

\begin{definition}\label{defn:period-sheaves}
    A diamond $X$ over $\Spd(\QQ_p, \Z_p)$ is equipped with the following sheaves for the $v$-topology:
    \begin{enumerate}
        \item for a topological space $T$, the sheaf $\ud{T}$ is given by $\ud{T}(S) = \Map_\cont(\abs{S},T)$ for any $S\in \Perf_{/X}$;
        \item $\cO^{\flat,+}_X$ (resp.\,$\cO^{\flat}_X$), given by $\cO^{\flat,+}_X(S) = \cO^+_S(S)$ (resp.\,$\cO^{\flat}_X(S)=\O_S(S)$) for any $S\in \Perf_{/X}$;
        \item $\cO^+_X$ (resp.\,$\cO_X$), given by $\cO^+_X(S) = \cO^+_{S^\sharp}(S^\sharp)$ (resp.\,$\cO_X(S)=\cO_{S^\sharp}(S^\sharp)$) for any $S\in \Perf_{/X}$;
        \item $\AA_I$ (resp.\,$\BB_I$), given by $\AA_I(S)=\AA_I(S^\sharp)$ (resp.\,$\BB_I(S)=\BB_I(S^\sharp)$) for any $S\in \Perf_{/X}$ (see \cref{defn:period-rings-2}).
    \end{enumerate}
\end{definition}

In what follows, we will consider sheaves of the form $\ud{T}$ only for $T=\ZZ/p^n\ZZ$, $\ZZ_p$, or $\QQ_p$ (with the ``$p$-adic'' topologies on $\ZZ_p$ and $\QQ_p$). 

\begin{rmk} Using that the map $\abs{f}\colon \abs{X} \to \abs{Y}$ is a quotient map for any $v$-covering $X\to Y$, one easily concludes that $\ud{T}$ is a $v$-sheaf for any topological space $T$. Furthermore, \cite[Th.\,8.7]{diamonds} ensures that $\cO^+_X$, $\cO_X$, $\cO^{\flat, +}$, and $\cO_X^\flat$ are $v$-sheaves. Finally, \cref{cor:AI-etale-sheaf} implies that $\AA_I$ and $\BB_I$ are $v$-sheaves as well.
\end{rmk}

For the next definition, we fix a diamond $X$ over $\Spd(\QQ_p, \ZZ_p)$.

\begin{definition} An element $\varpi\in \O_X^+(X)$ (resp. $\pi\in \O_X^{\flat, +}(X)$) is a \emph{pseudo-uniformizer} if, for any affinoid perfectoid space $S$ with a map $S \to X$, the element $\varpi|_S\in \O_X^+(S)=\O^+_{S^\sharp}(S^\sharp)$ (resp.\,$\pi|_S\in \O_X^{\flat, +}(S)=\O_{S}^+(S)$) is a pseudo-uniformizer.    
\end{definition}

For future reference, we include the following lemma:

\begin{lemma}\label{lemma:uniformizer-divisible-by-p} Let $X$ be a quasi-compact diamond over $\Spd(\QQ_p, \ZZ_p)$ and let $\varpi\in \O_X^+(X)$ be a pseudo-uniformizer. Then we can find an integer $n$ and an element $a\in \O_X^+(X)$ such that $p^n= \varpi \cdot a$.
\end{lemma}
\begin{proof}
    Choose a quasi-pro-\'etale covering $S \to X$ with an affinoid perfectoid $S$. We also choose a quasi-pro-\'etale covering $\{T_i \to S\times_X S\}_{i\in I}$ with each $T_i$ affinoid perfectoid.
    Then since $p\in \O_{S^\sharp}^+(S^\sharp)$ is a pseudo-uniformizer, there is an integer $n$ and $a\in \O_{S^\sharp}^+(S^\sharp)$ such that $p^n=\varpi|_S\cdot a$ in $\O_{S^\sharp}^+(S^\sharp)$. It suffices to show that $a$ lies in the image of $\O_X^+(X) \hookrightarrow \O_{S^\sharp}^+(S^\sharp)$. To see this, it is enough to show that the two possible immages of $a$ under the map $\O_{S^\sharp}^+(S^\sharp) \to \prod_{i\in I} \O_{T_i^\sharp}^+(T_i^\sharp)$. This follows from the fact that the same holds for $\varpi|_S$ and $p^n$ and the fact that $\prod_{i\in I} \O_{T_i^\sharp}^+(T_i^\sharp)$ is $\varpi$-torsionfree (as $\varpi|_{T_i}\in \O_{T_i^\sharp}^+(T_i^\sharp)$ is regular since it is a pseudo-uniformizer). 
\end{proof}

For the next definition, we fix a diamond $X$ over $\Spd(\QQ_p, \ZZ_p)$ and pseudo-uniformizers $\varpi\in \O_X^+(X)$ and $\pi\in \O_X^{\flat, +}(X)$. 

\begin{definition}\label{defn:period-sheaves-2} The $v$-sheaves $\O_X^+/\varpi \in \Shv(X_v; \ZZ)$ and $\O_{X}^{\flat,+}/\pi\in \Shv(X_v; \ZZ)$ are defined as $v$-sheaf quotients of $\O_X^+$ and $\O_X^{\flat, +}$ respectively. 
\end{definition}

\begin{setup}\label{setup:period-sheaves} Let $X$ be a diamond over $\Spd(\Q_p, \Z_p)$, let $\varpi\in \O_X^+(X)$ and $\pi\in \O_X^{\flat, +}(X)$ be pseudo-uniformizers, let $\ell$ be a prime (possibly $\ell=p$), $n>0$ a positive integer, and let 
\[
\cA \in \{\ud{\QQ}_\ell, \ud{\ZZ}_\ell, \ud{\ZZ/n}, \BB_I, \AA_I, \O_X, \O_X^+, \O_X^+/\varpi, \O_X^\flat, \O_X^{\flat, +}, \O_X^{\flat, +}/\pi\}
\]
be one of the above $v$-sheaves. 
\end{setup}

\begin{rmk}\label{rmk:no-map-to-qp} If $\cA\in\{\ud{\QQ}_\ell, \ud{\ZZ}_\ell, \ud{\ZZ/n}\}$ in \cref{setup:period-sheaves}, then one can drop the assumption that $X$ is defined over $\Spd(\QQ_p, \ZZ_p)$. 
\end{rmk}

\begin{notation} Let $X$ and $\cA$ be as in \cref{setup:period-sheaves}.  By a slight abuse of notation, we also denote by $\cA\in \Shv(X_\qproet; \ZZ)$ the restriction of $\cA$ to the (small) quasi-pro-\'etale site of $X$. We will sometimes write $\cA_{\tau}$ for $\tau\in \{\qproet, v\}$ to emphasize the topology we are working with.
\end{notation}

\begin{rmk} We note that \cref{cor:pullback-etale-finitary} and \cref{lemma:example-o-plus-mod-varpi} implies that $\Hh^i_\qproet(S; \O_{S_\qproet}^+/\varpi) \simeq \Hh^i_v(S; \O_{S_v}^+/\varpi)$ for any perfectoid $S$ over $\Spd(\QQ_p, \ZZ_p)$. This implies that $(\O_{X_v}^+/\varpi)|_{X_\qproet} \simeq \O_{X_\qproet}^+/\varpi$ for any diamond $X$ over $\Spd(\QQ_p, \ZZ_p)$ and a pseudo-uniformizer $\varpi\in \O_X^+(X)$. Therefore, there is no ambiguity in the meaning of $\O_X^+/\varpi$. The same applies to $\O_X^{\flat, +}/\pi$ for a pseudo-uniformizer $\pi\in \O_X^{\flat, +}(X)$. 
\end{rmk}

Now we are ready to establish some basic properties of period sheaves.

\begin{proposition}\label{prop:AI-properties}
    Let $S$ be an affinoid perfectoid space with a morphism $S\to \Spd(\Q_p, \Z_p)$, let $S^\sharp$ be the corresponding untilt, and let $X$ be a diamond over $S$. Let $p^\flat \in \O_S^{\flat, +}(S)$ be an element as in \cref{notation:p-flat} corresponding to $S^\sharp$, let $\alpha, \beta \in \O_S^{\flat, +}(S)$ be two elements well-adapted to $p^\flat$ and $I$ (see \cref{defn:well-adapted}), and let $\tau \in \{\et, \qproet, v\}$. Then there is a canonical isomorphism 
    \[
    \AA_{I}/p \simeq \left(\bigoplus_{i\geq 0} \cO^{\flat,+}_{X}/\alpha \right) \oplus \left(\bigoplus_{i\geq 1} \cO^{\flat,+}_{X}/\beta \right)
    \]
    in $\Shv(X_\tau; \ZZ)$. 
\end{proposition}
\begin{proof}
    This follows immediately from \cref{cor:AI-mod-p-rings} and the fact that sheafification commutes with any colimits.     
\end{proof}

\begin{lemma}\label{lemma:A_inf-is-p-adically-complete} Let $X$ be a diamond over $\Spd(\Q_p, \Z_p)$ and let $\tau\in \{\qproet, v\}$. Then $\BB_I \simeq \AA_I[\frac{1}{p}]$ in $\Shv(X_\tau; \Z)$ and the natural morphisms
\[
\AA_I \to \lim_\NN \AA_I/p^n \text{ and } \AA_I \to \rR\lim_\NN \AA_I/^\rL p^n.
\]
are isomorphisms in $\Shv(X_\tau; \Z)$ and $\cD(X_\tau; \ZZ)$, respectively. In particular, $\AA_I$ is a derived $p$-adically complete sheaf.
\end{lemma}
\begin{proof}
    All these questions are $\tau$-local on $X$, so we can assume that $X$ is an affinoid perfectoid space. Then the isomorphism $\AA_{I}[\frac{1}{p}] \simeq \BB_{I}$ follows from the fact that $X_\tau$ is a coherent topos. 

    Now we note that $\AA_{I}$ is a $p$-torsionfree sheaf due to \cref{cor:almost-mathematics-AI}. Therefore, $\AA_{I}/p^n \simeq \AA_{I}/^\rL p^n$. Therefore, it suffices to show that $\AA_{I}$ is a derived $p$-adically complete and the natural map $\lim_\NN \AA_{I}/p^n \to \rR\lim_\NN \AA_{I}/p^n$ is an isomorphism. The first claim follows from \cref{lemma:derived-p-complete-AI}, while the second claim follows from \cite[Prop.~3.1.10]{proetale} and the fact that $X_\tau$ is a replete topos (see \cref{lemma:replete}). 
\end{proof}

\begin{cor}\label{cor:continuity-cohomology} Let $X$, $\varpi\in \O_X^+(X)$, $\pi\in \O_X^{\flat, +}(X)$,  $\ell$, $n$, and $\cA$ be as in \cref{setup:period-sheaves}, and let $T$ be a profinite set. Then there is a functorial bijection $\cA(X\times T) \simeq \Map_{\cont}\big(T, \cA(X)\big)$.
\end{cor}
\begin{proof}
    Since any such $\cA$ is a quasi-pro-\'etale sheaf of topological abelian groups, we conclude that both $\cA(-\times T)$ and $\Map_{\cont}(T, \cA(-))$ are also quasi-pro-\'etale sheaves of topological abelian groups. Therefore, it suffices to construct the desired functorial bijection only for strictly totally disconnected $X$. 
    
    In this case, \cite[Lem.\,4.3.7]{proetale} reduces the cases $\cA \in \{\ud{\QQ}_\ell, \BB_I, \O_X, \O_X^\flat\}$ to $\cA\in \{\ud{\ZZ}_\ell, \AA_I, \O_X^+, \O_X^{\flat, +}\}$. If $\cA = \O_X^+$ (resp. $\O_X^{\flat, +}$, resp. $\ud{\ZZ}_\ell$), we use that $\cA$ is $\varpi$-adically complete (resp. $\pi$-adically complete, resp. $\ell$-adically complete) to reduce this case to the case $\cA = \O_X^+/\varpi$ (resp.\,$\O_X^{\flat, +}/\pi$, resp.\,$\ud{\ZZ/\ell^n}$). If $\cA = \AA_I$, we use that $\cA$ is $p$-adically complete (see \cref{lemma:A_inf-is-p-adically-complete}), \cref{rmk:often-well-adapted-pairs}, and \cref{prop:AI-properties} to reduce this case to $\cA = \O_X^{\flat, +}/\pi$. In other words, it suffices to treat the cases when $\cA= \O_X^+/\varpi$, $\cA = \O_X^{\flat, +}/\pi$, and $\cA=\ud{\ZZ/n}$. 
    
    The sheaves $\O_X^{\flat, +}/\pi$, $\O_X^+/\varpi$, and $\ud{\ZZ/n}$ come as pullbacks from the etale site by virtue of \cref{lemma:finitary-everything-is-etale}, \cref{lemma:example-o-plus-mod-varpi}, and \cref{lemma:example-z-mod-n}. Therefore, \cite[Prop.~8.5]{diamonds}  ensures that 
    \[
    \cA(X\times T) \simeq \colim_{I^\op} \cA(X\times T_i) \simeq \colim_{I^\op} \cA(X)^{\abs{T_i}} \simeq \colim_{I^\op} \Map_{\cont}(T_i, \cA(X)) \simeq \Map_{\cont}(T, \cA(X)),
    \]
    where $T=\lim_I T_i$ is a cofiltered limit of finite sets. 
\end{proof}

\begin{cor}\label{cor:no-higher-cohomology-period-sheaves} Let $X$, $\varpi\in \O^+_X(X)$, $\pi\in \O_X^{\flat, +}(X)$,  $\ell$, $n$, and $\cA$ be as in \cref{setup:period-sheaves}. If $X$ is a strictly totally disconnected perfectoid space, then $\rR\Gamma_\tau(X, \cA) \simeq \cA(X)[0]$ for $\tau\in \{\qproet, v\}$.
\end{cor}
In particular, this result implies that $\rR\Gamma_\qproet(X, \cA) \simeq \rR\Gamma_v(X, \cA)$.
\begin{proof}
    Arguing as in \cref{cor:continuity-cohomology}, we can reduce to the case $\cA = \O_X^+/\varpi$, $\cA = \O_X^{\flat, +}/\pi$, or $\cA=\ud{\ZZ/n}$. In either of these cases, the result follows immediately from \cref{cor:pullback-etale-finitary}, \cref{lemma:example-o-plus-mod-varpi}, and \cref{lemma:example-z-mod-n}, and the fact that $\Hh^i_\et(X, \cF)=0$ for any \'etale sheaf $\F$ and $i>0$. 
\end{proof}

For future reference, we also mention the following result.

\begin{lemma}
\label{lemma:no-difference-proetale-v} 
Let $X$, $\varpi\in \O_X^+(X)$, $\pi\in \O_X^{\flat, +}(X)$,  $\ell$, $n$, and $\cA$ be as in \cref{setup:period-sheaves}. 
Then the pullback functor $\mu^*\colon \Perf_\qproet(X; \cA) \to \Perf_v(X; \cA)$ is an equivalence. 
\end{lemma}
\begin{proof}
The claim is quasi-pro-\'etale local, so we can assume that $X$ is a strictly totally disconnected perfectoid space. Then the cases when $\cA=\O^+_X/\varpi$, $\O^{\flat, +}_X/\pi$, or $\ud{\ZZ/n}$ follow from \cref{cor:pullback-etale-finitary}, \cref{lemma:example-o-plus-mod-varpi}, and \cref{lemma:example-z-mod-n}, the case of $\uZ_\ell$ follows from \Cref{cor:pullback-Z-ell}~\Cref{cor:pullback-Z-ell-2}, the case of $\uQ_\ell$ follows from \Cref{thm:v-perfect-complexes-vs-etale-perfect-complexes-rationally-Qell}~\cref{thm:v-perfect-complexes-vs-etale-perfect-complexes-rationally-Qell-2}, the case $\cA=\O^+_X$ or $\cA=\O_X^{\flat, +}$ follows from \cref{thm:v-perfect-complexes-vs-etale-perfect-complexes-integrally}, the case $\cA=\O_X$ or $\cA = \O_X^\flat$ follows from \cref{thm:v-perfect-complexes-vs-etale-perfect-complexes-rationally}, the case $\cA=\AA_I$ follows from \cref{thm:v-perfect-complexes-vs-etale-perfect-complexes-integrally-A_I}, and the case $\cA=\BB_I$ follows from \cref{thm:v-perfect-complexes-vs-etale-perfect-complexes-rationally-BBI}. 
\end{proof}

\subsection{Condensed global sections}

In this subsection, we recall the definition of condensed (derived) global sections and prove some properties of condensed global sections of perfect complexes over period rings. 

\begin{definition}
\label{condensed structure definition}
    Let $X$ be a diamond and let $\tau\in \{\qproet, v\}$.
    Let $X_\tau \to *_\proet$ be the natural morphism of sites which maps a profinite set $T = \lim_i T_i \in *_\proet$ to $X \times \underline{T} = \lim_i X \times T_i$.
    The \emph{condensed global sections} functor
    \[ \underline{\Gamma}_\tau(X,\blank) \colon \Shv(X_\tau; \underline{\ZZ}) \to \Mod_{\underline{\ZZ}}(\Cond) \]
    is the pushforward along $X_\tau \to *_\proet$.
    We denote the corresponding \emph{derived condensed global sections} functor by
    \[ \rR\underline{\Gamma}_\tau(X,\blank) \colon \cD(X_\tau;\underline{\ZZ}) \to \cD(\Mod_{\underline{\ZZ}}(\Cond)). \]
\end{definition}

\begin{variant}\label{variant:condensed-global-sections} Let $S$ be an affinoid perfectoid space over $\Spd(\Q_p, \Z_p)$, let $\varpi\in \O_{S}^+(S)$ and $\pi\in \O_S^{\flat, +}(S)$ be pseudo-uniformizers, let $X$ be a diamond over $S$, let $\ell$ be a prime number, let $n>0$ be a positive integer, let $\tau\in \{\qproet, v\}$, and let $\cA$ be as in \cref{setup:period-sheaves}. Then the functor of the derived condensed global sections factors as 
\[
\rR\uGamma_{\tau}(X, \blank) \colon \cD(X_\tau; \cA) \to \cD(\Mod_{\ud{\cA(S)}}(\Cond)).
\]
\end{variant}

We first show that quasi-pro-\'etale derived condensed global sections coincide with $v$ derived condensed global sections in our cases of interest. 

\begin{lemma}
\label{lemma:properties-condensed-global-sections} Let $X$, $\varpi\in \O_X^+(X)$, $\pi\in \O_X^{\flat, +}(X)$,  $\ell$, $n$, and $\cA$ be as in \cref{setup:period-sheaves}. Let $\tau\in \{\qproet, v\}$. 
Then 
\begin{enumerate}[label=\upshape{(\roman*)}]
    \item\label{lemma:properties-condensed-global-sections-1} there is a functorial isomorphism $\uGamma_\tau(X, \cA) \simeq \ud{\cA(X)}$, where $\ud{\cA(X)}$ is the condensed $\ud{\Z}$-module corresponding to the topological $\Z$-module $\cA(X)$;
    \item\label{lemma:properties-condensed-global-sections-2} the natural morphism $\uGamma_\tau(X, \cA)[0] \to \uRGamma_\tau(X, \cA)$ is an isomorphism if moreover $X$ 
    is strictly totally disconnected.
\end{enumerate}
\end{lemma}
\begin{proof}
    Part~\cref{lemma:properties-condensed-global-sections-1} follows immediately from \cref{cor:continuity-cohomology}. 
    For part~\cref{lemma:properties-condensed-global-sections-2}, we note that $X\times T$ is strictly totally disconnected for any profinite set $T$ (see \cite[Lem.~7.19]{diamonds}). Therefore, it suffices to show that $\cA(X)[0] \to \rR\Gamma_\tau(X, \cA)$ for any strictly totally disconnected $X$. This follows immediately from \cref{cor:no-higher-cohomology-period-sheaves}. 
\end{proof}

In \cref{setup:period-sheaves}, \cref{lemma:properties-condensed-global-sections} can be bootstrapped to the following variants. 

\begin{lemma}
\label{lemma:properties-condensed-derived-global-sections} 
Let $X$, $\varpi\in \O_X^+(X)$, $\pi\in \O_X^{\flat, +}(X)$,  $\ell$, $n$, and $\cA$ be as in \cref{setup:period-sheaves},
and let $\cP\in \Perf_\qproet(X; \cA)$.
Then the natural morphism
\[
\rR\uGamma_\qproet(X, \cal{P}) \to \rR\uGamma_v(X, \mu^*\cal{P})
\]
is an isomorphism. 
\end{lemma}
\begin{proof}
    The question is quasi-pro-\'etale local on $X$, so we can assume that $X$ is a strictly totally disconnected perfectoid space and $\cP$ is represented by a finite complex of direct summands of finite free $\cA$-modules.
    By induction, we reduce to the case when $\cP$ is a direct summand of $\cA^{\oplus n}$, which could be further reduced to the case $\cP=\cA$. In this case, the result follows immediately from \cref{lemma:properties-condensed-global-sections}. 
\end{proof}

We finish the subsection with the following result:

\begin{lemma}
\label{RGamma-solid} 
    Let $X$, $\varpi\in \O_X^+(X)$, $\pi\in \O_X^{\flat, +}(X)$, $\ell$, $n$, and $\cA$ be as in \cref{setup:period-sheaves}, and let $\cP\in \Perf_v(X; \cA)$. Then the derived condensed global sections $\rR\uGamma_v(X, \cP)\in \cD(\Cond(\Ab))$ are solid. 
\end{lemma}
\begin{proof}
    Since $\cD(\Solid) \subset \cD(\Cond(\Ab))$ is closed under limits (see \cite[Def.\,5.1]{Condensed}), we can argue
    $v$-locally on $X$, so we can assume that $X$ is a strictly totally disconnected perfectoid space and $\cP$ is a strictly perfect complex. Then an inductive argument reduces the question to the case when $\cP$ is a direct summand of $\cA^{\oplus n}$ for some integer $n\geq 0$. Since $\cD(\Solid) \subset \cD(\Cond(\Ab))$ is closed under retracts and direct sums, we can further reduce to the case $\cP=\cA$. Since $\cD(\Solid) \subset \cD(\Cond(\Ab))$ is closed under all limits and colimits, we can run a reduction as in the proof of \cref{cor:continuity-cohomology} to reduce the question to the case $\cA=\ud{\ZZ/n}$, $\cA=\O_X^+/\varpi$, or $\cA=\O_X^{\flat, +}/\pi$. In either of these cases, we note that \cref{lemma:properties-condensed-global-sections} implies that $\rR\uGamma_v(X, \cA) \simeq \ud{\cA(X)}$, where we topologize $\cA(X)$ using the discrete topology. Therefore, the result follows from the fact that discrete abelian groups are solid (which follows from the fact that $\ud{\Z}$ is solid and solid abelian groups are closed under colimits; see \cite[Th.\,5.8]{Condensed}). 
\end{proof}

\subsection{Digression: \'etale sheaves on diamonds}
In this subsection, we give a brief overview of \'etale sheaves on diamonds and their base change properties.
Recall that given a locally spatial diamond $X$, one can define the \'etale site $X_\et$ whose objects are
diamonds $X'$ together with an
\'etale map $X' \to X$ and whose covers are jointly surjective maps \cite[Def.~14.1]{diamonds}.
The resulting category $\cD^+(X_\et;\ZZ/n)$ embeds fully faithfully into $\cD^+(X_\tau;\ZZ/n)$ for $\tau \in \{ \qproet, v\}$ via the natural pullback functor \cite[Prop.~14.10]{diamonds}.
Moreover, containment in this full subcategory can be tested $v$-locally;
that is, if $f \colon Y \to X$ is a $v$-covering of locally spatial diamonds, then $F \in \cD^+(X_\tau;\ZZ/n)$ is contained in $\cD^+(X_\et,\ZZ/n)$ if and only if $f^*F \in \cD^+(Y_\tau;\ZZ/n)$ is contained in $\cD^+(Y_\et;\ZZ/n)$ \cite[Th.~14.12.(i)]{diamonds}.

Later in \cref{trace-cycle}, we will also need to consider diamonds which are not necessarily spatial for the construction of compactly supported cohomology;
see \cref{not-spatial-warning}.
In this general case, the \'etale site is not large enough because $\Pro(X_\et)$ may not yield a basis for $X_\qproet$ anymore. 
Thus, \'etale sheaves should instead be defined via descent following \cite{diamonds}.

Throughout this subsection, we fix an integer $n>0$. 

\begin{definition}[{cf.\ \cite[Def.~14.13]{diamonds}}]
    Let $X$ be a diamond.
    The \emph{\'etale derived category of left-bounded complexes} $\cD^+_\et(X;\ZZ/n) \subseteq \cD^+(X_\qproet; \ZZ/n)$ is the full subcategory consisting of those $F \in \cD^+(X_\qproet; \ZZ/n)$ such that $f^*F \in \cD^+(Y_\et;\ZZ/n)$ for all maps from strictly totally disconnected perfectoid spaces $f \colon Y \to X$.
\end{definition}

Before we begin studying the properties of this category, we want to make some remarks.

\begin{remark}\label{etale-pushforward}
Given any qcqs morphism of diamonds $f \colon X \to Y$, the derived quasi-pro-\'etale pushforward functor  $\rR f_{\qproet, *}\colon \cD(X_\qproet; \ZZ/n) \to \cD(Y_\qproet; \ZZ/n)$ preserves $\cD^+_\et$ by the same argument as in \cite[Cor.~16.8.(ii)]{diamonds}.
In this setting, this quasi-pro-\'etale pushforward therefore agrees with the abstractly defined pushforward $\rR f_*$ from \cite[Lem.~17.5]{diamonds}.
\end{remark}

\begin{rmk}\label{rmk:etale-proetale-pushforward}
If $X$ is locally spatial, then $\cD^+_\et(X;\ZZ/n) \simeq  \cD^+(X_\et;\ZZ/n)$ \cite[Rmk.~14.14]{diamonds}. Throughout this paper, we will freely use this identification. 
Moreover, if $f\colon X \to Y$ is a qcqs map of locally spatial diamonds, then \cite[Cor.\,16.7]{diamonds} and \cref{etale-pushforward} imply that $\rR f_{\et, *}=\rR f_{\qproet, *}=\rR f_* \colon \cD^+_\et(X; \ZZ/n) \to \cD^+_\et(Y; \ZZ/n)$
\end{rmk}

\begin{rmk}[{cf.\ \cite[Cor.~16.10]{diamonds}, \cite[Prop.~3.1]{AM}}]\label{rmk:quasi-proetale-base-change} Since the (derived) quasi-pro-\'etale pushforward commutes with arbitrary quasi-pro-\'etale pullback, the above discussion implies that, for any qcqs morphism of diamonds $f\colon X \to Y$, the derived pushforward functor $\rR f_{\et, *} = \rR f_{\qproet, *}\colon \cD^+_\et(X; \ZZ/n) \to \cD^+_\et(Y; \ZZ/n)$ commutes with quasi-pro-\'etale base change.     
\end{rmk}

\begin{warning}
The analog of \cref{etale-pushforward} for $\rR f_{v,*}$ does not hold because, unlike \cite[Cor.~16.8.(ii)]{diamonds}, we do not assume that $n$ is prime to $p$.
Indeed, if this were the case, then similarly to  \cref{rmk:quasi-proetale-base-change}, it would follow that for any qcqs map $f \colon X \to Y$ of locally spatial diamonds, $\rR f_{\et,*}$ commutes with arbitrary base change.
This is already false for the projection morphism $\DD^{1,\lozenge}_C \to \Spd(C,C^\circ)$ and the base change along an extension $(C,C^\circ) \subset (C',C^{\prime,\circ})$ of algebraically closed nonarchimedean fields. 
\end{warning}

\begin{definition}\label{geometric-points}
    Let $X$ be a diamond and $x \in \abs{X}$ be a point.
    A \emph{qpro\'et-geometric point over $x$} is a quasi-pro-\'etale map $\bar{x} \colon \Spa(C,C^+) \to X$ for some algebraically closed nonarchimedean field $C$ and bounded open valuation subring $C^+ \subset C$ which sends the closed point to $x$.
    A \emph{qpro\'et-geometric point of $X$} is a qpro\'et-geometric point over some $x \in X$.
\end{definition}
\begin{definition}
    Let $X$ be a diamond and $\ov{x} \colon \Spa(C, C^+) \to X$ be a qpro\'et-geometric point of $X$.
    The \emph{stalk functor} $(\blank)_{\ov{x}} \colon \cD^+_\et(X, \ZZ/n) \to \cD^+(\ZZ/n)$ is the composite functor 
    \[
    \rR\Gamma_{\qproet}\big(\Spa(C, C^+), \ov{x}^*(\blank)\big) \colon \cD^+_\et(X, \ZZ/n) \to \cD^+(\ZZ/n).
    \]
\end{definition}
\begin{lemma}[{cf.\ \cite[Prop.~14.3]{diamonds}}]\label{enough-points}
    Let $X$ be a diamond.
    \begin{enumerate}[leftmargin=*,label=\upshape{(\roman*)}]
        \item\label{enough-points-existence} Any point $x \in \abs{X}$ has a qpro\'et-geometric point over it. 
        \item\label{enough-points-conservative} Choose a qpro\'et-geometric point $\ov{x}\colon \Spa(C_x, C_x^+) \to X$ over each point $x\in \abs{X}$. Then the family of functors $\{ (\blank)_{\bar{x}}\}_{x\in \abs{X}} \colon \cD^+_\et(X, \ZZ/n) \to \cD^+(\ZZ/n)$ is conservative. 
    \end{enumerate}
\end{lemma}
\begin{proof}
    Let $f \colon Y \to X$ be a quasi-pro-\'etale covering by a strictly totally disconnected perfectoid space $Y$.
    Given any $x \in \abs{X}$, we can pick $y \in \abs{Y}$ with $f(y) = x$.
    Let $(C,C^+)$ be a completed algebraic closure of the residue field of $y$ (as in \cite[(2.5.2)]{Huber-etale}).
    Then the induced map $\bar{y} \colon \Spa(C,C^+) \to Y$ is quasi-pro-\'etale because it is the limit over the (cofiltered) category of \'etale neighborhoods of $y \in \abs{Y}$.\footnote{
    Here we use that, for a finite separable extension $K\subset L$ and a valuation subring $L^+ \subset L$ extending $K^+$, the map $\Spa(L, L^+) \to \Spa(K, K^+)$ is a composition of a pro-open immersion and a finite etale map (it does not need to be finite \'etale).}
    Hence, the composition $\bar{x} \colonequals f \circ \bar{y}$ is quasi-pro-\'etale as well, yielding \cref{enough-points-existence}. 

    Next, we claim that for any given qpro\'et-geometric point $\bar{x} \colon \Spa(C_x,C^+_x) \to X$ over $x$, there is a qpro\'et-geometric point $\ov{y}\colon \Spa(C_y, C_y^+) \to Y$ over $y$ and a morphism $g\colon \Spa(C_y, C_y^+) \to \Spa(C_x, C_x^+)$ such that the diagram
    \[
    \begin{tikzcd}
    \Spa(C_y, C_y^+) \arrow{d}{g} \arrow{r}{\ov{y}} & Y \arrow{d}{f}\\
    \Spa(C_x, C_x^+) \arrow{r}{\ov{x}} & X
    \end{tikzcd}
    \]
    commutes:
    Indeed, consider the fiber product $Y_{\ov{x}} \coloneqq Y\times_X \Spa(C_x, C_x^+)$. 
    Since $\ov{x}$ is quasi-pro-\'etale, we conclude that $Y_{\ov{x}} \to Y$ is pro-\'etale.
    On the other hand, the map $\abs{Y_{\ov{x}}} \to \abs{\Spa(C_x, C_x^+)} \times_{\abs{X}} \abs{Y}$ is surjective, so there is a point $y'\in \abs{Y_{\ov{x}}}$ lying above $y$.
    Choose a qpro\'et-geometric point $\ov{y'}\colon \Spa(C_y, C_y^+) \to Y_{\ov{x}}$ over $y'$. Then the composition $\Spa(C_y, C_y^+) \xr{\ov{y'}} Y_{\ov{x}} \to Y$ does the job.
    
    Thanks to the claim from the last paragraph and the conservativity of $\cD_\et^+(X; \ZZ/n) \to \cD_\et^+(Y; \ZZ/n)$, it suffices to prove \cref{enough-points-conservative} for $Y$ instead of $X$.
    In this case (or, more generally, for any locally spatial diamond), the assertion follows from \cite[Prop.~14.3]{diamonds}.
\end{proof}
\begin{lemma}\label{pushforward-stalks}
    Let $f \colon X \to Y$ be a qcqs morphism of diamonds and $F \in \cD^+_\et(X;\ZZ/n)$. 
    Let $\bar{y} \colon \Spa(C,C^+) \to Y$ be a qpro\'et-geometric point of $Y$ and set $X_{\bar{y}} \colonequals X \times_{Y,\bar{y}} \Spa(C,C^+)$.
    Then the natural map
    \[ (\rR f_{\qproet, *} F)_{\bar{y}} \longrightarrow \rR\Gamma_{\qproet}(X_{\bar{y}},F_{\bar{y}}) \]
    in $\cD(\ZZ/n)$ is an equivalence. 
\end{lemma}

\begin{proof}
    Denote by $f_{\bar{y}} \colon X_{\bar{y}} \to \Spa(C,C^+)$ and $\iota \colon X_{\bar{y}} \to X$ the morphisms induced by $f$ and $\bar{y}$, respectively. After unwinding the definitions, \cite[Prop.\,14.10]{diamonds} implies that it suffices to show that the natural morphism 
    \[
    \rR\Gamma_{\qproet}(\Spa(C, C^+), \ov{y}^* \rR f_{\qproet, *} F) \to \rR\Gamma_\qproet(X_{\ov{y}}, F_{\ov{y}})
    \]
    is an isomorphism. This morphism is induced by the base change map $\bar{y}^*\rR f_{\qproet, *}F \to \rR f_{\bar{y},\qproet, *}\iota^* F$, so it suffices to show that this latter map is an isomorphism. 
    This follows directly from \Cref{rmk:quasi-proetale-base-change}.
\end{proof}

For future reference, we include proofs for two lemmas for which we could not find precise references. The first of these is widely known.

\begin{lemma}\label{lemma:replete} Let $X$ be a diamond and let $\tau\in \{\qproet, v\}$. Then $X_\tau$ is a replete topos (in the sense of \cite[Def.\,3.1.1]{proetale}).
\end{lemma} 

\begin{proof}
    If $\tau=\qproet$, this was proven in \cite[Lem.\,1.2]{Mann-Werner}. Thus, we only need to deal with $\tau=v$. For a tower of surjective morphisms
    \[
    \dots \twoheadrightarrow F_{n+1} \twoheadrightarrow \dots \twoheadrightarrow F_1 \twoheadrightarrow F_0,
    \]
    we need to show that $\lim_n F_n \to F_0$ is surjective. Since affinoid perfectoid spaces form a basis for $X_v$, it suffices to show that for any affinoid perfectoid $S_0$ and a section $s_0\in F_0(S_0)$, there is a $v$-covering $S \twoheadrightarrow S_0$ such that $s_0|_S\in F_0(S)$ lies in the image of $(\lim_n F_n)(S) \to F_0(S)$.
    
    For $S_0$ and $s_0\in F_0(S_0)$ as above, there is an affinoid perfectoid space $S_1$ and a $v$-covering $S_1 \to S_0$ such that $s_0$ lifts to an element $s_1\in F_1(S_1)$. By induction, we can construct a tower of $v$-coverings $S_n\to S_{n-1}$ of affinoid perfectoid spaces and sections $s_n\in F_n(S_n)$ such that $s_n$ lifts $s_{n-1}$ for each $n$. Then $S\coloneqq \lim_n S_n \to S_0$ is still a $v$-covering and the system $\{s_n|_{S}\in F_n(S)\}$ defines a section $s\in (\lim F_n)(S)$ which lifts $s_0|_S\in F_0(S)$. This finishes the proof.
\end{proof}
\begin{lemma}\label{cartesian-cocartesian}
    Consider a fiber square of small $v$-sheaves
    \[ \begin{tikzcd}
            Z' \arrow[r,hook,"i'"] \arrow[d,"\pi'"] & X' \arrow[d,"\pi"] \\
            Z \arrow[r,hook,"i"] & X
    \end{tikzcd} \]
    in which all maps are proper in the sense of \cite[Def.~18.1]{diamonds}, $i$ is a monomorphism, and $\pi$ is an isomorphism over the open $\abs{X} \smallsetminus \abs{i}(\abs{Z})$.
    Then the fiber square is also a pushout square in $\Shv(\Perf_v;\Set)$ and in $\Shv(\Perf_v;\Ani)$.
\end{lemma}
\begin{proof}
By assumption, the map of $v$-sheaves $Z \sqcup X' \xrightarrow{i \sqcup \pi} X$ is quasi-compact and induces a surjection $\abs{Z \sqcup X'} \simeq \abs{Z} \sqcup \abs{X'} \xtwoheadrightarrow{\abs{i} \sqcup \abs{\pi}} \abs{X}$.
Therefore, it is surjective thanks to \cite[Lem.~12.11]{diamonds}. In particular, 
    it suffices to prove that the diagram becomes a pushout square after base change along $i$ and $\pi$.
    Moreover, colimits in a topos (resp.~$\infty$-topos) are universal \cite[Prop.~II.4.3.1)]{SGA4} (resp.\ \cite[Prop.~6.1.3.10]{HTT}), so this is equivalent to proving that the diagrams of $v$-sheaves
    \[ \begin{tikzcd}
            Z \times_X Z' \arrow[r,hook] \arrow[d] & Z \times_X X' \arrow[d] \\
            Z \times_X Z \arrow[r,hook] & Z
    \end{tikzcd}
    \quad \text{and} \quad
    \begin{tikzcd}
            X' \times_X Z' \arrow[r,hook] \arrow[d] & X' \times_X X' \arrow[d] \\
            X' \times_X Z \arrow[r,hook] & X'
    \end{tikzcd} \]
    are pushout squares.
    Now we proceed in two steps.
    \begin{enumerate}[wide,label={\textit{Step~\arabic*}.},ref={Step~\arabic*}]
        \item\label{cartesian-cocartesian-monomorphism} \textit{Proof for the left square.}
        Since $i$ is a monomorphism, the bottom horizontal map $Z \times_X Z \to Z$ is an isomorphism.
        As the diagram is already a fiber square, one can then verify directly the universal property that it is also a pushout square.
        
        \item \textit{Proof for the right square.} 
        This square identifies with the right square in the larger diagram
        \[ \begin{tikzcd}
            Z' \arrow[d,hook,"i'"] \arrow[r,hook,"\Delta"] & Z' \times_Z Z' \arrow[r] \arrow[d,hook,"i' \times i'"] & Z' \arrow[d,hook,"i'"] \\
            X' \arrow[r,hook,"\Delta"] & X' \times_X X' \arrow[r] & X'.
        \end{tikzcd} \]
        Since the horizontal compositions are the identity, one can check again directly from the universal properties that it suffices to show that the left square of the larger diagram is a pushout square.
        In this square, all morphisms are proper monomorphisms, and 
        $X' \sqcup (Z' \times_Z Z') \xrightarrow{\Delta \sqcup (i' \times i')} X' \times_X X'$ is again surjective,
        thanks to the assumption that $\pi$ is an isomorphism over the open $\abs{X} \smallsetminus \abs{i}(\abs{Z})$ and hence $\Delta$ is an isomorphism over the open $\abs{X' \times_X X'} \smallsetminus \abs{i' \times i'}(\abs{Z' \times_Z Z'})$.
        Therefore, the same proof as in \cref{cartesian-cocartesian-monomorphism} applied to the pullbacks along both $\Delta$ and $i' \times i'$ shows that the left square of the larger diagram is a pushout square. \qedhere
    \end{enumerate}
\end{proof}

\subsection{Digression: perfect complexes on locally spatial diamonds}

In this subsection, we study perfect $\ud{\ZZ/n}$, $\O_X^+/\varpi$, and $\O_{X}^{\flat, +}/\pi$-complexes on locally spatial diamonds. In particular, we show that these categories are independent of whether we work with the \'etale, quasi-pro-\'etale, or $v$-topologies. We also compare different versions of (derived) pushforwards. We start with the case of $\ud{\ZZ/n}$-complexes. 

\begin{lemma}\label{lemma:locally-spatial-Fp-perfect-complexes} Let $X$ be a locally spatial diamond, let $\tau\in \{\qproet, v\}$, let $n\geq 1$ be an integer, and let $\lambda_\tau\colon (X_\tau; \ud{\ZZ/n}) \to (X_\et; \ud{\ZZ/n})$ be the natural morphism of ringed sites. Then the functor
\[
\lambda_\tau^* \colon \Perf_\et(X; \ud{\ZZ/n}) \to \Perf_\tau(X; \ud{\ZZ/n})
\]
is an equivalence.
\end{lemma}
\begin{proof}
    First, \cite[Prop.\,14.10]{diamonds} implies that $\lambda_\tau^*$ is fully faithful. Therefore, we only need to show that $\lambda_\tau^*$ is essentially surjective. Let $\cal{P} \in \Perf_\tau(X; \ud{\ZZ/n})$.
    We proceed by first proving that the counit of adjunction $\epsilon_\cP \colon \lambda_\tau^*\rR\lambda_{\tau, *} \cP \to \cP$ is an equivalence, and then that $\rR\lambda_{\tau, *} \cP \in \cD(X_\et; \ud{\ZZ/n})$ is a perfect complex.
    
    For the first claim, note that there exists $\cP_\et\in \cD^+(X_\et; \ZZ/n)$ such that $\lambda^*_\tau \cP_\et \simeq \cP$: since this can be checked $\tau$-locally thanks to \cite[Th.\,14.12(i)]{diamonds}, it suffices to prove it when $X$ is a perfectoid space, where we can apply \cref{cor:pullback-etale-finitary} and \cref{lemma:example-z-mod-n}.
    Thus, we are reduced to showing that 
    \[
    \epsilon_{\lambda^*_\tau \cP_\et} \colon \lambda_\tau^*\rR\lambda_{\tau, *}\lambda^*_\tau \cP_\et \longrightarrow \lambda^*_\tau \cP_\et
    \]
    is an equivalence. The axioms of units-counits imply that the composition
    \[
    \lambda_\tau^* \cP_\et \xr{\lambda^*_\tau(\eta_{\cP_\et})} \lambda_\tau^*\rR\lambda_{\tau, *}\lambda^*_\tau \cP_\et \xr{\epsilon_{\lambda^*_\tau \cP_\et}} \lambda^*_\tau \cP_\et
    \]
    is equal to $\id_{\lambda^*_\tau \cP_\et}$, where $\eta_{\cP_\et} \colon \cP_\et \to \rR\lambda_{\tau, *}\lambda^*_\tau \cP_\et$ is the unit of the $(\lambda_\tau^*, \rR\lambda_{\tau, *})$-adjunction. Now $\eta_{\cP_\et}$ is an isomorphism because $\lambda_\tau^*$ is fully faithful. In particular, we conclude $\lambda_\tau^*(\eta_{\cP_\et})$ is an isomorphism. This yields the claim $\epsilon_{\lambda_\tau^*\cP_\et}=\epsilon_{\cP}$ is isomorphisms as well. 
    
    Finally, to see that $\rR\lambda_{\tau, *} \cP \in \cD(X_\et; \ud{\ZZ/n})$ is a perfect complex, it suffices to show that $\rR\lambda_{\tau, *}\cP$ is dualizable \cite[\href{https://stacks.math.columbia.edu/tag/0FPV}{Tag 0FPV}]{stacks-project}.
    For this, we set $\cP^{\vee} \coloneqq \rR\cHom(\cP, \ud{\ZZ/n})$ and consider the cup product map
    \[
    \cup \colon \rR\lambda_{\tau,*}\cP \otimes_{\ud{\ZZ/n}}^\rL \rR\lambda_{\tau,*}\cP^\vee \longrightarrow \rR\lambda_{\tau, *}\bigl(\cP \otimes_{\ud{\ZZ/n}}^\rL \cP^{\vee}\bigr).
    \]
    We claim that it is an equivalence. Since $\lambda_{\tau}^*$ is fully faithful, it suffices to show that $\lambda_\tau^*(\cup)$ is an equivalence. Using that $\lambda_\tau^*\rR\lambda_{\tau, *} \cal{Q} \simeq \cal{Q}$ for any $\cal{Q}\in \Perf_\tau(X; \ud{\ZZ/n\ZZ})$, we conclude that $\lambda_\tau^*(\cup) \simeq \id_{\cP \otimes^\rL \cP^{\vee}}$, so $\cup$ is an equivalence. 

    Now let $c\colon \ud{\ZZ/n} \to \cP \otimes_{\ud{\ZZ/n}}^\rL \cP^\vee$ and $e\colon \cP \otimes_{\ud{\ZZ/n}}^\rL \cP^\vee \to \ud{\ZZ/n}$ be the co-evaluation and evaluation morphisms for $\cP$. Then we define the coevaluation and evaluation maps for $\rR\lambda_{\tau, *}\cP$ as follows
    \[
    c' \colon \ud{\ZZ/n} \xrightarrow[\sim]{\eta_{\ud{\ZZ/n}}} \rR\lambda_{\tau,*}\bigl(\ud{\ZZ/n}\bigr) \xrightarrow{\rR\lambda_{\tau, *}(c)}\rR\lambda_{\tau, *}\bigl(\cP \otimes_{\ud{\ZZ/n}}^\rL \cP^\vee) \xr[\sim]{\cup^{-1}} \rR\lambda_{\tau,*}\cP \otimes_{\ud{\ZZ/n}}^\rL \rR\lambda_{\tau,*}\cP^\vee,
    \]
    \[
    e'\colon \rR\lambda_{\tau,*}\cP \otimes_{\ud{\ZZ/n}}^\rL \rR\lambda_{\tau,*}\cP^\vee \xrightarrow[\sim]{\cup} \rR\lambda_{\tau, *}\bigl(\cP \otimes_{\ud{\ZZ/n}}^\rL \cP^\vee) \xr{\rR\lambda_{\tau,*}(e)} \rR\lambda_{\tau, *}\bigl(\ud{\ZZ/n}\bigr) \xr[\sim]{\eta^{-1}_{\ud{\ZZ/n}}} \ud{\ZZ/n}.
    \]
One checks directly that this defines the structures of a dualizable object on $\rR\lambda_{\tau, *}\cP$ and, thus, it is perfect. 
\end{proof}

Now we start showing that quasi-pro-\'etale and $v$ pushforwards of perfect $\ud{\ZZ/n}$-complexes are canonically identified. For this, we would need the following definition.  

\begin{defn}\label{defn:fiberwiseproper} A map of diamonds $f\colon X \to Y$ over $\Spd(\ZZ_p, \ZZ_p)$ is \emph{fiberwise proper} if for every morphism $\bar{y} \colon \Spd(C(\bar{y}),C(\bar{y})^+) \to Y$, where $C(\bar{y})$ is an algebraically closed nonarchimedean field and $C(\bar{y})^+ \subset C(\bar{y})$ a bounded open valuation subring, the fiber product $X_{\bar{y}} \coloneqq X\times_Y \Spd(C(\bar{y}),C(\bar{y})^+)$ is the diamond attached to a proper finite type adic space $X^\sharp_{\bar{y}}$ over $\Spa(C(\bar{y})^\sharp, C(\bar{y})^{\sharp,+})$. 
\end{defn}

\begin{warning} In \cref{defn:fiberwiseproper}, the adic space $X^\sharp_{\bar{y}}$ over $\Spa(C(\bar{y})^\sharp, C(\bar{y})^{\sharp,+})$ is not unique. For example, for any such $X^\sharp_{\bar{y}}$, we have $(X^\sharp_{\bar{y}})^{\diam} \simeq (X^\sharp_{\bar{y}, \red})^{\diam}$. 
\end{warning}

\begin{lemma}
\label{lemma:proper-base-change-etale}
Let $g\colon S' \to S$ be a morphism of strictly totally disconnected perfectoid spaces over $\Spd(\ZZ_p, \ZZ_p)$, let $X$ be a spatial diamond with a morphism $f\colon X \to S$ which is fiberwise proper, let $n\geq 1$ be an integer, and let $\cP\in \Perf_\et(X; \ud{\ZZ/n})$. Set $X_{S'}\coloneqq X\times_S S'$ with the projection morphisms $g'\colon X_{S'} \to X$ and $f'\colon X_{S'} \to S'$. Then the natural morphisms
\begin{equation}
\label{eqn:normal-base-change}
    \alpha\colon g^*_\et \rR f_{\et, *} \cP \to \rR f'_{\et, *} g'^* \cP, \text{ and }
\end{equation}
\begin{equation}\label{eqn:global-section-base-change}
    \beta\colon \rR\Gamma_\et(X, \cP) \otimes^\rL_{\ud{\ZZ/n}(S)} \ud{\ZZ/n}(S') \to \rR\Gamma_\et(X_{S'}, g'^*_\et \cP)
\end{equation}
are isomorphisms.
\end{lemma}
\begin{proof}
    Since $S'$ and $S$ are strictly totally disconnected, one sees that $\rR\Gamma_\et(X_{S'}, \alpha)=\beta$, so it suffices to show that $\alpha$ is an isomorphism. Since $X$ is spatial (in particular, it is qcqs), the morphism $f\colon X \to S$ is qcqs. Hence, \cref{enough-points}, \cref{rmk:etale-proetale-pushforward}, and \cref{pushforward-stalks} imply that we can argue on stalks at qpro\'et-geometric points and so we can reduce the question to showing that the base change map 
    \[ \rR\Gamma_\et(X,\cP) \to \rR\Gamma_\et(X \times_{\Spa(C, C^+)} \Spa(C', C'^+),g'^*_\et \cP) \]
    is an isomorphism, where $C\subset C'$ is an extension of algebraically closed nonarchimedean fields with open bounded valuation subrings $C^+\subset C \cap C'^+$ and $C'^+\subset C'$, and $X$ is the diamond attached to a proper finite type adic space over $\Spa(C^\sharp, C^{\sharp, +})$. Then \cite[Lem.\,15.6]{diamonds}, \cite[Lem.\,10.3]{adic-notes} imply that it suffices to show an analogous claim when $C^+=\O_C$, $C'^+=\O_{C'}$, and $X$ is a proper rigid-analytic space over $C^\sharp$.
    In this case, the result follows immediately from \cite[Th.~1.4]{Heu24}.
\end{proof}

\begin{proposition}
\label{qproet-v-pushforward-etale}
Let $Y$ be a diamond over $\Spd(\ZZ_p, \ZZ_p)$ and let 
$f \colon X \to Y$ be a fiberwise proper morphism which is representable in spatial diamonds.
Let $n\geq 1$ be an integer, let $\cP\in \Perf_\qproet(X; \ud{\ZZ/n})$,
    and let $\mu_X\colon (X_v, \ud{\ZZ/n}) \to (X_\qproet, \ud{\ZZ/n})$ be the projection morphism. 
    Then the base change morphism
    \[ c\colon \mu^*_Y \rR f_{\qproet,*} \cP \to \rR f_{v,*} \mu^*_X\cP \]
    is an isomorphism. 
\end{proposition}
\begin{proof}
    The statement is quasi-pro-\'etale local on $Y$, so we may assume that $Y$ is a strictly totally disconnected perfectoid space.
    Then it suffices to show that for any strictly totally disconnected space $S \to Y$ and any factorization $S \to T \to Y$ such that $T$ is an affinoid perfectoid space and $T \to Y$ is pro-\'etale, the natural base change map
    \[ \rR\Gamma_\qproet(X\times_Y T,\cP) \otimes^\rL_{\ud{\ZZ/n}(T)} \ud{\ZZ/n}(S) \to \rR\Gamma_v(X\times_Y S,\mu^* \cP) \]
    is an isomorphism. Now \cite[Lem.\,7.19]{diamonds} implies that $T$ is strictly totally disconnected as well. Therefore, we can replace $T$ with $Y$ and reduce the question to showing that the natural map
    \[
    \rR\Gamma_\qproet(X,\cP) \otimes^\rL_{\ud{\ZZ/n}(Y)} \ud{\ZZ/n}(S) \to \rR\Gamma_v(X\times_Y S,\mu^* \cP)
    \]
    is an isomorphism for any morphism of strictly totally disconnected perfectoid spaces $S \to Y$. In this case, both $X$ and $X_S$ are spatial, and so \cref{lemma:locally-spatial-Fp-perfect-complexes} implies that the question boils down to showing that, for any $\cP_\et\in \Perf_\et(X; \ud{\ZZ/n})$, the natural map
    \[
    \rR\Gamma_\et(X,\cP_\et) \otimes^\rL_{\ud{\ZZ/n}(Y)} \ud{\ZZ/n}(S) \to \rR\Gamma_\et(X\times_Y S,g'^*_\et\cP_\et)
    \]
    is an isomorphism, where $g'\colon X \times_Y S \to X$ is the projection morphism. This follows immediately from \cref{lemma:proper-base-change-etale}. 
\end{proof}

Our next goal is to show analogs of \cref{lemma:locally-spatial-Fp-perfect-complexes} and \cref{qproet-v-pushforward-etale} for perfect $\O^+/\varpi$ and $\O^{\flat, +}/\pi$-complexes. 

\begin{lemma}
\label{lemma:comes-from-etale}
    Let $X$ be a locally spatial diamonds over $\Spd(\QQ_p, \ZZ_p)$, let $\varpi\in \O_X^+(X)$ \emph{(}resp.\,$\pi\in \O_X^{\flat, +}(X)$\emph{)} be a pseudo-uniformizer, let $\tau\in \{\qproet, v\}$, and let $\lambda_{\tau} \colon X_\tau \to X_\et$ be the projection morphism. Set $\cA \coloneqq \O_X^+/\varpi$ \emph{(}resp.\,$\cA\coloneqq \O_X^{\flat, +}/\pi$\emph{)} and $\cA_\et \coloneqq \lambda_{\tau, *}\cA$. Then the following hold: 
    \begin{enumerate}
        \item\label{lemma:comes-from-etale-1} if $\tau=\qproet$, then the natural morphism $\lambda^{-1}_\qproet \cA_\et \to \cA$ is an isomorphism;
        \item\label{lemma:comes-from-etale-2} the functor $\lambda_\tau^* \colon \Perf_\et(X; \cA_\et) \to \Perf_\tau(X; \cA)$ is an equivalence. 
    \end{enumerate}
\end{lemma}
\begin{proof}
    We explain the proof when $\cA = \O_X^+/\varpi$, the other case being very similar. Since both claims are \'etale local on $X$, we can assume that $X$ is a spatial diamond. In particular, it is qcqs. 
    
    We first show \cref{lemma:comes-from-etale-1}. First, we know that \cref{lemma:uniformizer-divisible-by-p} implies that $\cA$ is a sheaf of $\ZZ/p^n$-modules for some integer $n$. Thus, \cite[Prop.\,14.8]{diamonds} implies that it suffices to show that $\cA$ lies in the essential image of the fully faithful functor $\Shv_\et(X; \ud{\ZZ/p^n}) \hookrightarrow \Shv_\qproet(X; \ud{\ZZ/p^n})$. Then \cite[Th.\,14.12]{diamonds} ensures that we can check this $\qproet$-locally on $X$, so we can assume that $X$ is an affinoid perfectoid space. In this case, the result follows immediately from \cref{lemma:finitary-everything-is-etale} and \cref{lemma:example-o-plus-mod-varpi}. This finishes the proof of \cref{lemma:comes-from-etale-1}

    Now we address Part~\cref{lemma:comes-from-etale-2}. \cref{lemma:no-difference-proetale-v} ensures that it suffices to deal with $\tau=\qproet$. In this case, the proof is completely analogous to the proof of \cref{lemma:locally-spatial-Fp-perfect-complexes}. The only difference is that we use \cref{lemma:example-o-plus-mod-varpi} instead of \cref{lemma:example-z-mod-n}.     
\end{proof}

We next show the following version of \cref{qproet-v-pushforward-etale}.
 
\begin{proposition}\label{qproet-v-pushforward}
    Let $f \colon X \to Y$ be a qcqs morphism of diamonds over $\Spd(\QQ_p,\ZZ_p)$, let $\varpi\in \O_Y^+(Y)$ (resp.\,$\pi\in \O_Y^{\flat, +}(Y)$) be a pseudo-uniformizer, let $\cP\in \Perf_\qproet(X; \O^+_X/\varpi)$ (resp.\,$\cP\in \Perf_\qproet(X; \O_X^{\flat, +}/\pi)$). Then the natural base change map
    \[ c\colon \mu^*_Y \rR f_{\qproet,*} \cP \to \rR f_{v,*} \mu^*_X\cP \]
    is an $\cO^{\circ\circ}_Y$-almost (resp.\ $\cO^{\flat,\circ\circ}_Y$-almost) isomorphism.\footnote{This means that, for any affinoid perfectoid $S$ with a map $S \to Y$, we have $\rm{fib}(c)(S) \otimes^\rL \cO^{\circ\circ}_Y(S)=0$ (resp.\,$\rm{fib}(c)(S) \otimes^\rL \cO^{\flat, \circ\circ}_Y(S)=0$).}
\end{proposition}
\begin{proof}
    For brevity, we set $\cA \coloneqq \O_Y^+/\varpi$ (resp. $\cA\coloneqq \O_Y^{\flat, +}/\pi$). Then, arguing as in \cref{qproet-v-pushforward-etale}, we reduce the question to showing that the natural morphism
        \[
        \rR\Gamma_\qproet(X,\cP) \otimes^\rL_{\cA(Y)} \cA(S) \to \rR\Gamma_v(X\times_Y S,\cP)
        \]
    is an almost isomorphism when $S \to Y$ is a morphism of strictly totally disconnected perfectoid spaces over $\Spd(\QQ_p, \ZZ_p)$. 
    
    For this, we choose a simplicial quasi-pro-\'etale hypercovering $X_\bullet \to X$ by strictly totally disconnected perfectoid spaces such that $\restr{\cP}{X_n}$ is a strictly perfect complex for all $n \in \Delta$.
    The products $X_n \times_Y S$ are affinoid perfectoid, so we have $\rR\Gamma_\qproet(X_n \times_Y S,\cP) \simeq \rR\Gamma_v(X_n \times_Y S,\cP)$ by \cref{cor:pullback-etale-finitary} and \cref{lemma:example-o-plus-mod-varpi}.
    In particular, 
    \[
    \rR\Gamma_v(X \times_Y S,\cP) \simeq \rR\lim_\Delta \rR\Gamma_v(X_n \times_Y S,\cP) \simeq \rR\lim_\Delta \rR\Gamma_\qproet(X_n \times_Y S,\cP).
    \]
    Therefore, it suffices to show that the morphism $\rR\Gamma_\qproet(X,\cP) \otimes^\rL_{\cA(Y)} \cA(S) \to \rR\lim_\Delta \rR\Gamma_\qproet(X_n \times_Y S,\cP)$ is an almost isomorphism. Now we write this morphism as a composition
    \[
    \rR\Gamma_\qproet(X,\cP) \otimes^\rL_{\cA(Y)} \cA(S) \xr{\alpha} \rR\lim_\Delta \bigl(\rR\Gamma_\qproet(X_n,\cP) \otimes^\rL_{\cA(Y)} \cA(S)\bigr) \xr{\beta} \rR\lim_\Delta \rR\Gamma_\qproet(X_n \times_Y S,\cP)
    \]
    and separately show that $\alpha$ and $\beta$ are almost isomorphisms. First, we show that $\alpha$ is an (actual) isomorphism. Since quasi-pro-\'etale cohomology satisfy quasi-pro-\'etale descent, it suffices to prove that the natural map
    \[
    \bigl(\rR\lim_\Delta \rR\Gamma_\qproet(X_n,\cP)\bigr) \otimes^\rL_{\cA(Y)} \cA(S) \to \rR\lim_\Delta \bigl(\rR\Gamma_\qproet(X_n,\cP) \otimes^\rL_{\cA(Y)} \cA(S)\bigr)
    \]
    is an isomorphism.
    This is a consequence of the facts that $\cD\bigl(\cA(S)\bigr)$ is right $t$-complete (see e.g.\ \cite[Prop.~1.3.5.21.(3)]{HA}), that $\cA(Y) \to \cA(S)$ is flat (see \cite[Prop.~7.23]{diamonds}), and that totalizations of uniformly bounded below complexes commute with left $t$-exact functors as long as the target
    is right $t$-complete (see \cite[Cor.\,C.0.6]{Mitya-Artem} for a precise statement). To conclude the proof, it is then enough to show that the natural base change map
    \[ \rR\Gamma_\qproet( X_n,\cP) \otimes_{\cA(Y)} \cA(S) \to \rR\Gamma_\qproet(X_n \times_Y S,\cP) \]
    is an almost isomorphism for all $n \in \Delta$.
    Since $\restr{\cP}{X_n}$ is strictly perfect, a standard argument allows us to reduce to the case $\cP = \cA[0]$.
    In that case, \cite[Prop.\,8.5(iii)]{diamonds} implies that $\rR\Gamma_\qproet(X_n,\cP) \overset{a}{\simeq} \cA(X_n)[0]$ and $\rR\Gamma_\qproet(X_n \times_Y S,\cP) \overset{a}{\simeq} \cA(X_n\times_Y S)[0]$ because $X_n$ and $X_n\times_Y S$ are affinoid perfectoid, respectively. Thus, the assertion follows from \cref{lemma:almost-integral-kunneth-formula}.
\end{proof}

\subsection{Digression: lattices in perfect complexes of \texorpdfstring{$\uQ_\ell$}{Ql}-modules}

Throughout this subsection, we fix an arbitrary prime number $\ell$. We show that every perfect complex of $\ud{\QQ}_\ell$-sheaves on a locally spatial diamond $X$ admits an integral lattice \'etale locally. We then apply this result to deduce that $\Perf_v(X; \ud{\QQ}_\ell)$ possesses a standard $t$-structure when $X$ is the diamondification of a rigid-analytic space.

In \cref{notation:sheaf-perfect-complexes-finite-amplitude}, we defined the $\infty$-category $\Perf^{[r,r']}_\tau(X; \cA)$ for any $r,r'\in \ZZ\cup \{-\infty, \infty\}$, $\cA \in \{\uQ_\ell, \uZ_\ell\}$, and $\tau\in \{\qproet, v\}$. By \cref{lemma:sheafification-of-fully-faithful}, this is a full $\infty$-subcategory of $\Perf^{[-\infty, \infty]}_\tau(X; \cA) \simeq \Perf_\tau(X; \cA)$, and thus embeds naturally into $\cD(X_\tau; \cA)$ as a full $\infty$-subcategory.
Furthermore, when $X$ is a strictly totally disconnected perfectoid space, \cref{cor:ell-adic-integral-perfect-complex-on-strictly-totally-disconnected-tor-amplitude} and \cref{cor:v-perfect-complexes-vs-etale-perfect-complexes-rationally-Qell-tor-amplitude} imply that $\Perf^{[r, r']}_v(X; \cA) \simeq \Perf^{[r, r']}_\qproet(X; \cA) \simeq \Perf^{[r, r']}\bigl(\cA(X)\bigr)$, where the right-hand side denotes perfect complexes with tor-amplitude in $[r, r']$. Consequently, for an arbitrary diamond $X$, the $\infty$-category $\Perf^{[r, r']}_\tau(X; \cA)$ naturally forms a full $\infty$-subcategory of $\cD^{[r,r']}(X_\tau ; \cA)$, where $\cD^{[r,r']}(X_\tau ; \cA)$ denotes the $\infty$-subcategory consisting of objects concentrated in degrees $[r, r']$.

\begin{lemma}\label{lemma:lattice-etale-locally-lsd} Let $X$ be a locally spatial diamond, let $\ell$ be a prime number, let $\tau\in \{\qproet, v\}$, and let $\cP\in \Perf^{[r,r']}_\tau(X; \ud{\QQ}_\ell)$ for some $r, r'\in \ZZ\cup \{-\infty, \infty\}$. Then there is an \'etale covering $X' \to X$ and an object $\cP^+ \in  \Perf^{[r,r']}_\tau(X'; \ud{\Z}_\ell)$ such that $\cP^+[\frac{1}{\ell}] \simeq \restr{\cP}{X'}$.
\end{lemma}
\begin{proof}
    First, \cref{lemma:no-difference-proetale-v} and \cref{rmk:no-map-to-qp} imply that we can assume that $\tau=\qproet$. Furthermore, the question is \'etale local on $X$, so we can assume that $X$ is spatial (in particular, that it is qcqs). Then \cref{lemma:colimit-of-sheaves} implies that we can assume that $r, r' \in \ZZ$ are finite. Then we consider the quasi-pro-\'etale sheaf of lattices 
    \[
    \Latt^{[r, r']}_{\qproet, \cal{P}}(\blank) \coloneqq * \times_{\Perf^{[r, r']}_\qproet(-; \ud{\QQ}_\ell)} \Perf^{[r, r']}_\qproet(-; \ud{\ZZ}_\ell) \colon \Perf^{\op}_{/X} \to \Cat_\infty,
    \]
    where the morphism $* \to \Perf^{[r, r']}_\qproet(-; \ud{\QQ}_\ell)$ corresponds to the element $\cP\in \Perf_{\qproet}^{[r, r']}(X; \ud{\QQ}_\ell)$. 
    
    It suffices to show that $\Latt^{[r, r']}_{\qproet, \cal{P}}(X')\neq \varnothing$ for some \'etale covering $X' \to X$. Note that $\Latt^{[r, r']}_{\qproet, \cal{P}}(\blank)$ is a sheaf of $(r'-r+1, 1)$-categories.

    {\it Step~$1$. The natural map $\colim_{I^\op} \Latt^{[r,r']}_{\qproet, \cal{P}}(Y_i) \to \Latt^{[r, r']}_{\qproet, \cal{P}}(Y)$ is an equivalence for any $Y=\lim_I Y_i\in \Pro(X_{\et, \mathrm{qc}, \mathrm{sep}})$.} First, \cref{examples:compactly-generated-by-cotruncated-objects}~\cref{examples:compactly-generated-by-cotruncated-objects-2}, \cref{lemma:finite-limits-vs-filtered-colimits}, and quasi-pro-\'etale descent imply that we can argue quasi-pro-\'etale locally on $X$. So we can assume that $X$ is a strictly totally disconnected perfectoid space. In this case, all $Y_i$ are also strictly totally disconnected perfectoid spaces by virtue of \cite[Lem.~7.19]{diamonds}. The claim then follows directly from \cref{lemma:lattices-are-finitary-Qell}. 

    {\it Step~$2$. End of proof.} First, \cite[Prop.~11.24]{diamonds} ensures that there is a strictly totally disconnected perfectoid space $Y$ and a quasi-pro-\'etale covering $Y \to X$ such that $Y=\lim_I Y_i$ is a cofiltered limit of quasi-compact, separated, \'etale morphisms $Y_i \to X$. Each of $Y_i \to X$ is necessarily surjective since so is $Y\to X$. Then Step~$1$ implies that it suffices to show that $\Latt^{[r, r']}_{\qproet, \cal{P}}(Y) \not\simeq \varnothing$. This follows from \cref{lemma:lattice-etale-locally-Qell} and the fact that every \'etale covering of $Y$ splits. 
\end{proof}

For the next definition, we fix a diamond $X$, a prime number $\ell$, and $\tau\in \{\qproet, v\}$.

\begin{defn} A \emph{$\uQ_\ell$-local system} (resp.~a \emph{$\uZ_\ell$-local system}) on $X$ is a $\tau$-sheaf of $\uQ_\ell$-modules (resp.~$\uZ_\ell$-modules) $\F$ which is $\tau$-locally isomorphic to $M \otimes_{\QQ_\ell} \uQ_\ell$ (resp. $M\otimes_{\ZZ_\ell} \uZ_\ell$) for some finite-dimensional $\QQ_\ell$-vector space $M$ (resp.~finitely generated $\ZZ_\ell$-module $M$).
\end{defn}

\begin{lemma}\label{lemma:characterizing-perfect-complexes} Let $K$ be a nonarchimedean field of residue characteristic $p$, let $X^\sharp$ be a rigid-analytic space over $K$, and let $X\coloneqq (X^\sharp)^\diam$ be the associated diamond. Let $\ell$ be a prime number, let $\tau\in \{\qproet, v\}$, let $\cA\in \{\uZ_\ell, \uQ_\ell\}$, and let $\cP\in \cD^b(X_\tau; \cA)$. Then $\cP$ lies in $\Perf_\tau(X; \cA)$ if and only if $\cal{H}^i(\cP)$ is a $\cA$-local system for every $i\in \ZZ$.
\end{lemma}
We warn the reader that the analogous statement does not hold for a general (locally spatial) diamond $X$. In fact, it is already false for $X = \ud{S} \times \Spa(C, \O_C)$ for any infinite metrizable pro-finite set $S$ and an algebraically closed nonarchimedean field $C$ of residue characteristic $p$. 

\begin{example} Let $S$ and $C$ be as above, and choose a non-isolated point $s\in S$. Since $S$ is a metrizable profinite set, there exists a decreasing sequence of clopen neighborhoods $S=U_0 \supset U_1 \supset U_2 \supset \dots$ such that $\{s\} = \cap_{n\geq 0} U_n$. Define the continuous function $f\colon S \to \QQ_p$ via the formula 
\[
f(x) = \begin{cases} 0 \text{ if $x=s$} \\ p^n \text{ if $x\in U_n\smallsetminus U_{n+1}$}\end{cases}.
\]
Then the  multiplication morphism $\cdot f\colon \Map_\cont(S, \QQ_p) \to \Map_\cont(S, \QQ_p)$ is injective and, thus, $M\coloneqq \coker(\cdot f)$ defines a perfect complex $M[0]\in \Perf\bigl(\Map_\cont(S, \QQ_p)\bigr)$. Then the associated complex of $\uQ_p$-modules $\widetilde{M[0]}$ lies in $\Perf\bigl(\ud{S} \times \Spa(C, \O_C)\bigr)$, but its (only) cohomology sheaf is not a $\uQ_p$-local system.
\end{example}

\begin{proof}
    First, any $\cA$-local system is a perfect complex. The ``if'' direction therefore follows immediately from the $2$-out-of-$3$ property for perfect complexes. 
    Hence, we only need to show that cohomology sheaves of a perfect complex are locally constant in the $\tau$-topology. The claim is \'etale local on $X$, so \cref{lemma:lattice-etale-locally-lsd} and \cite[Lem.~15.6]{diamonds} imply that it suffices to treat the case when $\cP\in \Perf_\tau(X; \uZ_\ell)$. In this case, \cref{cor:surjective-ell-adically}~\cref{cor:surjective-ell-adically-1}, \cref{lemma:locally-spatial-Fp-perfect-complexes}, and \cite[Lem.~11.1]{adic-notes} imply that $\Perf_\tau(X; \uZ_\ell)$ is equivalent to $\cD^b_{\mathrm{lisse}}(X, \Z_\ell)$ in the sense of \cite[Def.~3.32]{BH}. In this case, the result follows from \cite[Lem.~3.35]{BH}.
\end{proof}

\begin{cor} Keep the notation of \cref{lemma:characterizing-perfect-complexes}. Then the standard $t$-structure on $\cD(X_\tau; \ud{\QQ}_\ell)$ (resp.~on $\cD(X_\tau; \uZ_\ell)$) restricts to a $t$-structure on $\Perf_\tau(X; \ud{\QQ}_\ell)$ (resp.~$\Perf_\tau(X; \uZ_\ell)$).
\end{cor}

\begin{cor}\label{cor:local-system-tor-ampl-0} Keep the notation of \cref{lemma:characterizing-perfect-complexes}. Then $\cP\in \Perf^{[0, 0]}_\tau(X; \uQ_\ell)$ (resp.~$\cP\in \Perf^{[0, 0]}_\tau(X; \uZ_\ell)$) if and only if $\cP$ is $\uQ_\ell$-local system concentrated in degree $0$ (resp. $\cP$ is an $\ell$-torsionfree $\uZ_\ell$-local system concentrated in degree $0$).
\end{cor}
\begin{proof}
    We give a proof for $\uZ_\ell$-local systems. The case of $\uQ_\ell$-local systems is analogous (and in fact simpler). 

     Clearly, any $\ell$-torsionfree $\uZ_\ell$-local system (concentrated in degree 0) lies in $\Perf^{[0, 0]}_\tau(X; \uZ_\ell)$. Thus, we only need to show the any object in $\Perf^{[0, 0]}_\tau(X; \uZ_\ell)$ is isomorphic to an $\ell$-torsionfree $\uZ_\ell$-local system concentrated in degree $0$. 
     
     Now \cref{lemma:characterizing-perfect-complexes} implies that it suffices to show that any element of $\Perf^{[0, 0]}_\tau(X; \uZ_\ell)$ is concentrated in degree $0$ and is $\ell$-torsionfree. This can be checked $\tau$-locally, so we can assume that $X$ is a strictly totally disconnected perfectoid space. In this case, \cref{cor:ell-adic-integral-perfect-complex-on-strictly-totally-disconnected-tor-amplitude} implies that 
     \[
     \Perf^{[0, 0]}_v(X; \uZ_\ell) \simeq \Perf^{[0, 0]}_\qproet(X; \uZ_\ell) \simeq \Perf^{[0, 0]}(\Map_\cont(\abs{X}, \ZZ_\ell)) \simeq \mathrm{Mod}_{\Map_\cont(\abs{X}, \ZZ_\ell)}^{\mathrm{fin.proj.}},
     \]
     where as before $\Perf^{[0, 0]}(\Map_\cont(\abs{X}, \ZZ_\ell))$ denotes the category of perfect complexes with tor-amplitude in $[0, 0]$. In particular, any such object is $\ell$-torsionfree and concentrated in degree $0$.  
\end{proof}

For future reference, we also establish a stronger version of \cref{cor:local-system-tor-ampl-0} for $\uQ_\ell$-perfect complexes. 

\begin{cor}\label{cor:local-system-tor-ampl} Keep the notation of \cref{lemma:characterizing-perfect-complexes}. Assume that $\cP\in \Perf_\tau(X; \uQ_\ell)$ and let $r,r'\in \ZZ$. Then $\cP \in  \Perf_\tau^{[r,r']}(X; \uQ_\ell)$ if and only if $\cP\in \cD^{[r,r']}(X_\tau; \uQ_\ell)$.
\end{cor}
\begin{proof}
    First, we assume that $\cP\in \Perf_\tau^{[r,r']}(X; \uQ_\ell)$. We wish to show that $\cal{H}^i(\cP)=0$ for $i\notin [r,r']$. This can be checked locally, so we can assume that $X$ is a strictly totally disconnected perfectoid space. In this case, \cref{cor:v-perfect-complexes-vs-etale-perfect-complexes-rationally-Qell-tor-amplitude} implies that  $\Perf_\tau^{[r,r']}(X; \uQ_\ell) \simeq \Perf^{[r,r']}\bigl(\Map_\cont(X, \QQ_\ell)\bigr)$. Thus, the result follows from  \cite[\href{https://stacks.math.columbia.edu/tag/0658}{Tag 0658}]{stacks-project}.

    Now assume that $\cP\in \cD^{[r,r']}(X_\tau; \uQ_\ell) \cap \Perf_\tau(X; \uQ_\ell)$. Then \cref{cor:local-system-tor-ampl-0} ensures that $\cal{H}^i(\cP) \in \Perf^{[0,0]}_\tau(X; \uQ_\ell)$ and $\cal{H}^i(\cP) =0$ for $i\notin[r, r']$. We wish to show that $\cP\in \Perf^{[r, r']}_\tau(X; \uQ_\ell)$. This now follows from a simple inductive argument. 
\end{proof}

\subsection{Relative primitive comparison}

The main goal of this subsection is to prove two relative versions of the primitive comparison theorem.
We start with the following preliminary lemma:

\begin{lemma}\label{lemma:almost-overconvergent}
Let $K$ be a nonarchimedean field of mixed characteristic $(0, p)$, let $K^+\subset K$ be an open bounded valuation subring, and let $X$ be a qcqs adic space over $\Spa(K, K^+)$.
Set $X^\circ\coloneqq X\times_{\Spa(K, K^+)} \Spa(K, K^\circ)$ and let $j\colon X^\circ \hookrightarrow X$ be the natural pro-(open immersion). Let $\F\in D^b(X_\et; \FF_p)$ be a bounded complex such that each $\cal{H}^i(\F)$ is overconvergent. Then the natural morphism
\begin{equation}\label{eqn:almost-overconvergent}
\rR\Gamma_\et\bigl(X, \F \otimes^{\rL} \O_{X_\et}^+/p\bigr) \to \rR\Gamma_\et\bigl(X^{\circ}, j^*\F \otimes^\rL \O^+_{X^\circ_\et}/p\big)
\end{equation}
is an almost isomorphism.
\end{lemma}
\begin{proof}
    First, we may assume that $\F\in\Shv(X_\et; \FF_p)$ is an overconvergent sheaf of $\FF_p$-modules. Then \cite[Lem.\,6.7.10]{zav-almost} implies that both sides of \cref{eqn:almost-overconvergent} satisfy $v$-descent, so we can assume that $X=\Spa(R, R^+)$ is a strictly totally disconnected perfectoid space. In this case, \cite[Lem.\,6.7.6]{zav-almost} implies that $\F=\colim_I \F_i$ is a filtered colimit of special sheaves on $X$ (see \cite[Def.\,6.7.5]{zav-almost} for the definition of special sheaves). Since we have reduced to the case where $X$ is qcqs and hence both sides of \cref{eqn:almost-overconvergent} commute with uniformly bounded filtered colimits, it suffices to prove the claim when $\F$ is a special sheaf. In other words, there is a finite disjoint clopen decomposition $X= \sqcup_{r=1}^n X_r$ such that $\F = \bigoplus_{r=1}^n i_{X_r \to X, *} \ud{\FF}_p^{\oplus d_r}$. Then we reduce to the case when $\F = i_* \ud{\FF}_p$, where $i\colon Y \hookrightarrow X$ is a clopen subspace. In this case, we can replace $X$ with $Y$ and $X^\circ$ with $Y^\circ \coloneqq Y\times_{\Spa(K, K^+)} \Spa(K, K^\circ)$ to assume that $\F=\ud{\FF}_p$. 

    We have reduced the situation to the case when $X=\Spa(R, R^+)$ is a strictly totally disconnected perfectoid space and $\F = \ud{\FF}_p$. Since strictly totally disconnected perfectoid spaces do not have higher \'etale cohomology, the question now boils down to showing that the natural morphism
    \[
    \Gamma(X_\et, \O_{X_\et}^+/p) \to \Gamma(X_\et^{\circ}, \O_{X^{\circ}_\et}^+/p)
    \]
    is an almost isomorphism. Since $K^\circ\subset R^{\circ}$ and $R$ is complete, we deduce that 
    \[
    X^{\circ} = \Spa(R, R^+) \times_{\Spa(K, K^+)} \Spa(K, K^{\circ}) \simeq \Spa\bigl(R\wdh{\otimes}_K K,  (R\wdh{\otimes}_K K)^+\bigr) \simeq \Spa(R, R'^+)
    \]
    for some open integrally closed subring $R^+\subset R'^+\subset R^{\circ}$. Therefore, \cite[Th.\,6.3]{diamonds} implies that we only need to show that $R^+/p \to R'^+/p$ is an almost isomorphism. This follows immediately from \cite[Lem.\,B.13]{zav-almost}.
\end{proof}

The following statement is bootstrapped from \cite{Scholze-Hodge}.

\begin{lemma}[Relative primitive comparison I]
\label{Relative primitive comparison-1}
    Let $Y$ be a diamond over $\Spd(\QQ_p, \ZZ_p)$, let $f\colon X \to Y$ be a fiberwise proper morphism which is representable in spatial diamonds, let $\varpi\in \O_Y^+(Y)$ be a pseudo-uniformizer \emph{(}resp.\ $\pi\in\O_Y^{\flat, +}(Y)$ be a pseudo-uniformizer\emph{)}, let $n\geq 1$ be an integer such that $p^n=\varpi\cdot a$ for some $a \in \O_Y^+(Y)$ \emph{(}resp.\ $n=1$\emph{)}, let $\tau\in \{\qproet, v\}$, and let $\bL \in \Perf_\tau(X, \ud{\ZZ/p^n})$. Set $\cA \coloneqq \O^+/\varpi$ \emph{(}resp. $\cA\coloneqq \O^{\flat, +}/\pi$\emph{)}. Then the natural morphism 
    \begin{equation}\label{eqn:primitive-comparison}
    \bigl(\rR f_{\tau,*} \bL\bigr) \otimes^\rL_{\ud{\ZZ/p^n}} \cA \to \rR f_{\tau,*}\bigl(\bL \otimes^\rL_{\ud{\ZZ/p^n}} \cA\bigr)
    \end{equation}
    is an almost isomorphism with respect to $\O_Y^{\circ\circ}$ \emph{(}resp.\, $\O_Y^{\flat, \circ\circ}$\emph{)}.  
\end{lemma}
\begin{proof}
    First, \cref{lemma:no-difference-proetale-v}, \cref{qproet-v-pushforward-etale}, and \cref{qproet-v-pushforward} imply that it suffices to consider the case $\tau=\qproet$.
    Since the question is quasi-pro-\'etale local on $Y$, we can immediately reduce to the case where $S=Y$ is a strictly totally disconnected perfectoid space. In this case, the diamond $X$ is spatial, so \cref{lemma:locally-spatial-Fp-perfect-complexes} implies that $\lambda^*_\qproet \colon \Perf_\et(X;\ud{\ZZ/p^n}) \xr{\sim} \Perf_\qproet(X;\ud{\ZZ/p^n})$ is an equivalence. We denote by $\bL_\et \in \Perf_\et(X; \ud{\ZZ/p^n})$ the \'etale perfect complex corresponding to $\bL$. Furthermore, \cref{lemma:comes-from-etale} implies that 
    \begin{equation}\label{eqn:A-is-etale}
    \cA \simeq \lambda^{-1}_\qproet \lambda_{\qproet, *} \cA. 
    \end{equation}
    We set $\cA_\et \coloneqq \lambda_{\qproet, *} \cA$. Then \cite[Cor.~16.7]{diamonds} implies that it suffices to show that the natural morphism 
    \[
    (\rR f_{\et,*} \bL_\et) \otimes^\rL_{\ud{\ZZ/p^n}} \cA_\et \to \rR f_{\et,*}(\bL_\et \otimes^\rL_{\ud{\ZZ/p^n}} \cA_\et)
    \]
    is an almost isomorphism.
    
    Using \cref{enough-points} and \cref{pushforward-stalks}, it suffices to prove that given $s \in \abs{S}$ and a qpro\'et-geometric point $\bar{s} \colon \Spa(C,C^+) \to S$ over $s$, the natural map
    \[ \rR\Gamma_\et(X_{\bar{s}},\bL_\et) \otimes^\rL_{\ZZ/p^n} \cA_\et(\Spa(C, C^+)) \to \rR\Gamma_\et(X_{\bar{s}},\bL_\et \otimes^\rL_{\underline{\ZZ/p^n}} \cA_\et)\]
    is an almost isomorphism. By assumption, we know that $X_{\ov{s}}\simeq (X_{\ov{s}}^\sharp)^\diam$ for some proper finite type $\Spa(C^\sharp, C^{\sharp, +})$-adic space $X_{\ov{s}}^\sharp$. Now \cite[Lem.\,15.6]{diamonds} implies that $\cD\bigl((X_{\ov{s}}^\sharp)_\et; \ZZ/p^n\bigr) \simeq \cD(X_{\ov{s}}; \ZZ/p^n)$, so the question boils down to showing that 
    \begin{equation}\label{eqn:primitive-comparison-cohomology}
        \rR\Gamma_\et(X_{\ov{s}}^\sharp,\bL_\et) \otimes^\rL_{\ZZ/p^n} \cA_\et\bigl(\Spa(C^\sharp, C^{\sharp, +})\bigr) \to \rR\Gamma_\et(X_{\ov{s}}^\sharp,\bL_\et \otimes^\rL_{\underline{\ZZ/p^n}} \cA_\et)
    \end{equation}
    is an almost isomorphism. 
    
    {\it Step~$1$. $\cA=\O^+/p\simeq \O^{\flat, +}/p^\flat$.} In this case, we verify \cref{eqn:primitive-comparison-cohomology}. First, \cite[Th.\,6.5.9]{zav-almost} implies that $\cA_\et\simeq \O_{X_{\ov{s}, \et}^\sharp}^+/p$. Since the pro-(open immersion) $X_{\ov{s}}^\sharp\times_{\Spa(C^\sharp,C^{\sharp, +})} \Spa(C^\sharp,C^{\sharp, \circ}) \to X_{\ov{s}}^\sharp$ induces almost isomorphisms on the source and on the target (see \cite[Lem.\,10.3]{adic-notes} and \cite[Lem.\,B.13]{zav-almost} for the left side and \cref{lemma:almost-overconvergent} for the right hand side), we may moreover assume that $C^{\sharp, +} = C^{\sharp, \circ}$. Furthermore, using that $\cA_\et$ is a sheaf of $\ud{\ZZ/p}$-modules and the projection formula $\rR\Gamma_\et(X_{\ov{s}}^\sharp, \bL_\et) \otimes^\rL_{\ZZ/p^n} \ZZ/p \simeq \rR\Gamma_\et(X_{\ov{s}}^\sharp, \bL_\et \otimes^\rL_{\ud{\ZZ/p^n}} \ud{\ZZ/p})$ from \cite[Prop.~5.5.1(iv)]{Huber-etale}, we can immediately reduce to the case $n=1$. So we are left to show that the natural morphism
    \[
       \rR\Gamma_\et(X_{\ov{s}}^\sharp,\bL_\et) \otimes^\rL_{\ZZ/p} C^{\sharp, \circ}/p \to \rR\Gamma_\et(X_{\ov{s}}^\sharp,\bL_\et \otimes^\rL_{\underline{\ZZ/p}} \O^+_{X^\sharp_{s,\et}}/p)
    \]
    is an almost isomorphism for a proper rigid-analytic space $X_{\ov{s}}^\sharp$ over $C^\sharp$ and $\bL_\et\in \Perf_\et(X_{\ov{s}}^\sharp; \ZZ/p)$. This follows directly from 
    \cite[Lem.~7.3.7]{zav-almost} (see also \cite[Th.~5.1]{Scholze-Hodge} and \cite[Th.~3.17]{Scholze-CDM}).
    
    {\it Step~$2$. $\cA=\O^+/p^m$ or $\cA=\O^{\flat, +}/p^{\flat, m}$ for some integer $m\geq 1$.} In this case, we verify \cref{eqn:primitive-comparison} directly. The result follows immediately from Step~$1$ via a simple inductive argument using the short exact sequences
    \[ 0 \to \cO^+/p \xr{\cdot p^{m-1}} \cO^+/p^m \to \cO^+/p^{m-1} \to 0 \quad \text{and} \quad 0 \to \cO^{\flat,+}/p^{\flat} \xr{\cdot (p^{\flat})^{m-1}} \cO^{\flat,+}/p^{\flat,m} \to \cO^{\flat,+}/p^{\flat,m-1} \to 0. \]
        
    {\it Step~$3$. General $\cA=\O^+/\varpi$ and $\cA=\O^{\flat, +}/\pi$.} In this case, we verify again \cref{eqn:primitive-comparison-cohomology}. We explain the proof for $\cA = \cO^+/\varpi$, the case $\cA = \cO^{\flat,+}/\pi$ being analogous. For brevity, we rename $\Spa(C, C^+)$ as $Y$, $X_{\ov{s}}$ as $X$, and $X_{\ov{s}}^\sharp$ as $X^\sharp$. Then we consider the resolution
    \[
    K_X \coloneqq \Bigl(\dots \xr{\cdot \varpi} \O_{X}^+/p^n \xr{\cdot a} \O_{X}^+/p^n \xr{\cdot \varpi} \O_{X}^+/p^n\Bigr) \xr{\sim} \O_{X}^+/\varpi.
    \]
    Now \cref{eqn:A-is-etale} and \cite[Prop.\,14.8]{diamonds} imply that 
    \[
    K_{X, \et}\coloneqq \Bigl(\dots \xr{\cdot \varpi} \lambda_{\qproet, *}(\O_{X}^+/p^n) \xr{\cdot a} \lambda_{\qproet, *}(\O_{X}^+/p^n) \xr{\cdot \varpi} \lambda_{\qproet, *}(\O_{X}^+/p^n)\Bigr) \xr{\sim} \lambda_{\qproet, *}(\O_{X}^+/\varpi) = \cA_\et
    \]
    is still a resolution. We set $K_{X, \et, i}\coloneqq \sigma^{\geq -i} K_{X, \et}$, where $\sigma^{\geq -i}$ is the naive filtration. Then we see that $\hocolim_i K_{X,\et, i} \simeq \cA_{X, \et}$. We define $K_{Y, \et, i}$ in a similar way. In what follows, we freely use \cite[Lem.\,15.6]{diamonds} and the identifications $\cD(X^\sharp_\et; \ZZ/p^n) \simeq \cD(X_\et; \ZZ/p^n)$, $\cD(\Spa(C^\sharp, C^{\sharp, +})_{\et}; \ZZ/p^n) \simeq \cD(Y_\et; \ZZ/p^n)$. Then an inductive argument and Step~$2$ imply that the natural morphism
    \[
        \rR\Gamma_\et(X^\sharp,\bL_\et) \otimes^\rL_{\ZZ/p^n} K_{Y, \et, i}(\Spa(C^\sharp, C^{\sharp, +})) \to \rR\Gamma_\et(X^\sharp,\bL_\et \otimes^\rL_{\underline{\ZZ/p^n}} K_{X, \et, i})
    \]
    is an almost isomorphism for any $i\geq 1$. Since $\rR\Gamma_\et(X^\sharp, \blank)$ and $\cA_\et(\blank)$ commute with colimits (see \cite[Lem.\,9.1(2)]{adic-notes}), we get the following sequence of (almost) isomorphisms:
    \begin{align*}
    \rR\Gamma_\et(X^\sharp,\bL_\et) \otimes^\rL_{\ZZ/p^n} \cA_\et\bigl(\Spa(C^\sharp, C^{\sharp, +})\bigr)& \simeq \rR\Gamma_\et(X^\sharp, \bL_\et) \otimes^\rL_{\ZZ/p^n}\Bigl( \hocolim_i K_{Y, \et, i} \bigl(\Spa(C^\sharp, C^{\sharp, +})\bigr)\Bigr) \\
    & \simeq \hocolim_i \Bigl(\rR\Gamma_\et(X^\sharp, \bL_\et) \otimes^\rL_{\ZZ/p^n} K_{Y, \et, i} \bigl(\Spa(C^\sharp, C^{\sharp, +})\bigr)\Bigr) \\
    & \overset{\mathrm{a}}{\simeq} \hocolim_i \Bigl(\rR\Gamma_\et(X^\sharp, \bL_\et \otimes^\rL_{\ud{\ZZ/p^n}} K_{X, \et, i}) \Bigr) \\
    & \simeq \rR\Gamma_\et(X^\sharp, \bL_\et \otimes^\rL_{\ud{\ZZ/p^n}} \hocolim_i K_{X, \et, i})\\
    &  \simeq \rR\Gamma_\et(X^\sharp, \bL_\et \otimes^\rL_{\ud{\ZZ/p^n}} \cA_\et) \qedhere
    \end{align*}
\end{proof}

\begin{lemma}[Relative primitive comparison II]
\label{Relative primitive comparison-2}
    Let $f\colon X \to Y$ be a fiberwise proper morphism of diamonds over $\Spd(\QQ_p, \ZZ_p)$ which is representable in spatial diamonds, 
    let $(\mathcal{A}, \mathcal{B}) \in \{(\mathcal{O}^+, \mathcal{O}), (\AA_I, \BB_I)\}$,
    let $\tau\in \{\qproet, v\}$, and let $\bL \in \Perf_\tau(X, \uZ_p)$. Then the natural morphism 
    \begin{equation}\label{eqn:primitive-comparison-map}
        \bigl(\rR f_{\tau,*} \bL\bigr) \widehat{\otimes}^\rL_{\underline{\ZZ}_p} \cA \to \rR f_{\tau,*}\bigl(\bL \otimes^\rL_{\underline{\ZZ}_p} \cA\bigr)
    \end{equation}
    is an almost isomorphism with respect to $\AA_I^{\circ\circ}$ if $\cA=\AA_I$ and with respect to $\O_Y^{\circ\circ}$ if $\cA = \cO_Y^+$. 
    In particular, after inverting $p$, the map
    \[
    \bigl((\rR f_{\tau,*} \bL) \widehat{\otimes}^\rL_{\underline{\ZZ}_p} \cA \bigr)\bigl[\tfrac{1}{p}\bigr] \to \rR f_{\tau,*}(\bL \otimes^\rL_{\underline{\ZZ}_p} \mathcal{B})
    \]
    is an isomorphism.
\end{lemma}
Here, the symbol $\wdh{\otimes}^{\rL}$ denotes the derived $p$-adic completion of the derived tensor product.
See also \cite[Prop.~5.2.2]{ALBM} for a related statement.
\begin{proof}
    The question is $\tau$-local on $Y$, so we can immediately reduce to the case when $Y$ is a strictly totally disconnected perfectoid space. In this case, we choose elements $\alpha, \beta\in \O_Y^{\flat,+}(Y)$ which are well-adapted to $p^\flat$ and $I$ (see \cref{rmk:often-well-adapted-pairs}). Then \cref{cor:almost-mathematics-AI} and \cref{prop:AI-properties} imply that 
    \[
    \AA_I^{\circ\circ}/p \simeq \bigcup_{n \ge 1} (p^\flat)^{1/p^n}(\AA_I/p) \simeq \bigoplus_\NN \Bigl(\bigcup_{n\geq 1} (p^\flat)^{1/p^n}\bigl(\O_{Y}^{\flat, +}/\alpha\bigr)\Bigr)\oplus \bigoplus_\NN \Bigl(\bigcup_{n\geq 1}(p^\flat)^{1/p^n}\bigl(\O_{Y}^{\flat, +}/\beta\bigr) \Bigr).
    \]
    Therefore, we can use \cref{lemma:A_inf-is-p-adically-complete} and \cite[Lem.~A.5]{zav-almost} to reduce the case $\cA=\AA_I$ to the cases $\cA=\O_Y^{\flat, +}/\alpha$ and $\cA=\O_Y^{\flat, +}/\beta$. Then the result follows immediately from \cref{Relative primitive comparison-1}. 

    A similar argument reduces the case $\cA=\O_Y^+$ to the case $\cA = \O_Y^{+}/p$, which also follows immediately from  \cref{Relative primitive comparison-1}. 
    For the last statement, we note that $\rR f_{\tau,*}(\bL \otimes^\rL_{\underline{\ZZ}_p} \cA)\bigl[\tfrac{1}{p}\bigr] 
    \simeq \rR f_{\tau,*}(\bL \otimes^\rL_{\underline{\ZZ}_p} \mathcal{B})$ as $f$ is qcqs.
\end{proof}

\section{Trace and cycle class maps}\label{trace-cycle}

The purpose of this section is to construct the trace and cycle morphisms for the period sheaves $\O$ and $\BB_I$. We build upon the construction of $\ZZ/n$-trace from \cite[Section 6]{LRZ24} for smooth morphisms between locally noetherian analytic adic spaces. Since this level of generality is inadequate for our purposes, we first extend this construction to smooth morphisms of sousperfectoid spaces and then to smooth proper morphisms of diamonds over $\Spd(\QQ_p, \ZZ_p)$. Then we use this $\ZZ/n$-trace to define trace morphisms for $\O$ and $\BB_I$-coefficients.

Throughout this section, we fix an integer $n>0$ and a closed interval $I\subset (0, \infty)$ with rational endpoints. 

\subsection{Proper pushforward for diamonds}\label{diamond-preliminaries}

In our discussion of trace maps in \cref{traces}, we will need a theory of proper pushforwards in \'etale cohomology for maps between suitable analytic adic spaces.
While the construction of proper pushforwards in \cite{Huber-etale} probably carries over to our desired setting, Huber imposes the blanket assumption that all his spaces be locally noetherian, which is not sufficient for our purposes.
To address this issue, we review in this subsection a slightly different perspective on proper pushforwards using diamonds, which applies in our context.

In \cite{diamonds}, proper pushforwards are discussed in the generality of $v$-stacks, but only for coefficients whose torsion order is coprime to $p$.
We will explain that many results (and their proofs) also work for $\FF_p$-coefficients as long as one only uses suitable diamonds instead of $v$-stacks.
In many places, we impose the assumption that our morphisms are quasi-compact;
although the extension of the theory from quasi-compact separated to all separated morphisms in \cite{Huber-etale,diamonds} should still go through in our setting \cref{not-quasi-compact}, we do not discuss this here as the quasi-compact case is sufficient for the purposes of this paper.

\begin{definition}[{\cite[Prop.~18.6, Prop.~18.8.(iii), Prop.~22.3.(i)]{diamonds}}]\label{canonical-compactification}
    Let $f \colon X \to Y$ be a separated morphism of diamonds.
    The \emph{canonical compactification} of $X$ over $Y$ is the diamond $\overline{X}^{/Y}$ whose values on a totally disconnected perfectoid space $\Spa(R, R^+)$ are given by 
    \[ \overline{X}^{/Y}(R,R^+) = X(R,R^\circ) \times_{Y(R,R^\circ)} Y(R,R^+). \]
    We call $f$ \emph{compactifiable} if it is separated and the natural map $X \to \overline{X}^{/Y}$ is an open immersion.
\end{definition}
The induced map $\overline{X}^{/Y} \to Y$ is always partially proper, and it is proper if $f$ is quasi-compact \cite[Cor.~18.8.(i),(vi)]{diamonds}.
Moreover, if $X' \to X$ is a surjective (resp.\ quasi-pro-\'etale) morphism of diamonds which are separated over $Y$, then so is the induced morphism $\overline{X'}^{/Y} \to \overline{X}^{/Y}$ \cite[Cor.~18.8.(v),(vii)]{diamonds}.

\begin{example}\label{canonical-compactification-affinoid}
    For any morphism $(A,A^+) \to (B,B^+)$ of analytic Huber pairs over $\ZZ_p$ and any perfectoid Huber pair $(R,R^+)$, it follows immediately from the universal properties that
    \[ \Hom\bigl((B,B^+),(R,R^\circ)\bigr) \times_{\Hom((A,A^+),(R,R^\circ))} \Hom\bigl((A,A^+),(R,R^+)\bigr) \simeq \Hom\bigl((B,B^+_{A^+}),(R,R^+)\bigr), \]
    where $B^+_{A^+} \subset B^+$ is the integral closure of $A^+ + B^{\circ\circ} \subset B$ and all the homomorphisms are taken in the category of Huber pairs over $\ZZ_p$.
    Unwinding \cref{canonical-compactification} and \cite[Lem.~15.1.(ii)]{diamonds}, it follows that the natural map
    \[ \overline{\Spd(B,B^+)}^{/\Spd(A,A^+)} \longrightarrow \Spd(B,B^+_{A^+}) \]
    is an isomorphism.
\end{example}
The next lemma shows that there are plenty of examples of compactifiable morphisms of diamonds.
We will need the following definition.
\begin{definition}[{cf.\ \cite[Def.~1.2.1.(ii)]{Huber-etale}}]\label{+-weakly-finite-type}
    A morphism $f \colon X \to Y$ of analytic adic spaces is called \emph{locally of +-weakly finite type} if there exist Tate affinoid open covers $Y = \bigcup V_\alpha$ and $f^{-1}(V_\alpha) = \bigcup U_{\alpha\beta}$ such that the induced morphisms $\cO_Y(V_\alpha) \to \cO_X(U_{\alpha\beta})$ are topologically of finite type (in the sense of \cite[Lem.~3.3.(iii)]{Huber-generalization}) and there exist finitely many elements $b_1,\dotsc,b_r \in \cO^+_X(U_{\alpha\beta})$ such that $\cO^+_X(U_{\alpha\beta}) \subset \cO_X(U_{\alpha\beta})$ is the smallest integrally closed subring containing the image of $\cO^+_X(V_\alpha)$ and $b_1,\dotsc,b_r$.
    If $f$ is in addition quasi-compact, it is called of \emph{+-weakly finite type}.
\end{definition}
\begin{lemma}\label{analytic-adic-compactifiable}
    Let $f \colon X \to Y$ be a separated morphism of analytic adic spaces over $\Spa \ZZ_p$.
    Assume that $f$ is locally of +-weakly finite type as in \cref{+-weakly-finite-type}.
    Then the morphism $f^\lozenge \colon X^\lozenge \to Y^\lozenge$ of associated diamonds (in the sense of \cite[Def.~15.5]{diamonds}) is compactifiable. 
\end{lemma}
\begin{proof}
    By \cite[Prop.~22.3.(iii),(v), Lem.~15.6]{diamonds}, we can check the statement after passing to open covers of both $X$ and $Y$, so we may assume that $Y = \Spa(A,A^+)$ and $X = \Spa(B,B^+)$ are affinoid, that $A \to B$ is topologically of finite type, and that there exist $b_1,\dotsc,b_r \in B^+$ such that $B^+ \subset B$ is the smallest integrally closed subring of $B$ containing the image of $A^+$ and $b_1,\dotsc,b_r$.
    Let $B^+_{A^+}\subset B$ be the integral closure of $A^+ + B^{\circ\circ} \subset B$.
    By \cref{canonical-compactification-affinoid}, we have $\overline{X^\lozenge}^{/Y^\lozenge} \simeq \Spd(B,B^+_{A^+})$. 
    In particular, the natural map $X^\lozenge \to \overline{X^\lozenge}^{/Y^\lozenge}$ is the diamondification of the inclusion of the rational open subspace $\{ \abs{b_1},\dotsc,\abs{b_r} \le 1 \} \subset \Spa(B,B^+_{A^+})$.
\end{proof}
In the locally noetherian setting, the canonical compactifications in the sense of \cref{canonical-compactification} agree with the canonical compactifications from \cite{Huber-etale}.
\begin{lemma}\label{canonical-compactification-Huber}
    Let $f \colon X \to Y$ be a separated +-weakly finite type morphism of locally noetherian analytic adic spaces over $\ZZ_p$.
    Then $\overline{X^\lozenge}^{/Y^\lozenge}$ is naturally isomorphic to the diamond $(X')^\lozenge$ attached to the canonical compactification $X'$ of $f$ in the sense of Huber \cite[Th.~5.1.5]{Huber-etale}.
\end{lemma}
\begin{proof}
    We first show that $(X')^\lozenge \to Y^\lozenge$ is proper. The claim is local on $Y^\lozenge$, so we can assume that $Y$ is an affinoid space. Then \cite[Lem.~5.1.3(i)]{Huber-etale} and \cite[Lem.~15.6]{diamonds}  imply that $\abs{X'}=\abs{X'^\lozenge}$ is taut. Therefore, one can use \cite[Prop.~18.10]{diamonds} and a similar valuative criterion in adic spaces (see \cite[Lem.\,1.3.10]{Huber-etale}) to conclude that $(X')^\lozenge \to Y^\lozenge$ is proper. In particular, the universal property of $\overline{X^\lozenge}^{/Y}$ from \cite[Prop.~18.6]{diamonds} produces a unique natural morphism $\overline{X^\lozenge}^{/Y^\lozenge} \to (X')^\lozenge$ such that the diagram
    \[ \begin{tikzcd}
        X^\lozenge \arrow[r] \arrow[d] & \overline{X^\lozenge}^{/Y^\lozenge} \arrow[d] \arrow[ld] \\
        (X')^\lozenge \arrow[r] & Y^\lozenge
    \end{tikzcd} \]
    commutes. By virtue of \cite[Lem~11.11]{diamonds}, it suffices to show that $\overline{X^\lozenge}^{/Y^\lozenge}(K,K^+) \to (X')^\lozenge(K, K^+)$ is bijective for any perfectoid field $K$ and any open bounded valuation subring $K^+\subset K$. Using the valuative criterion on both sides, we reduce to showing that the map $X^\lozenge(K, K^\circ) = \overline{X^\lozenge}^{/Y^\lozenge}(K, K^\circ) \to (X')^\lozenge(K, K^\circ)$ is bijective. For this, it suffices to show that the natural map $X \to X'$ is bijective on rank-$1$ points. This follows from \cite[Lem.~1.1.10.(iv), Th.~5.1.5]{Huber-etale}.
\end{proof}
\begin{warning}\label{not-spatial-warning}
    If $f \colon X \to Y$ is a spatial separated morphism of diamonds, it is not clear whether the induced morphism $\overline{X}^{/Y} \to Y$ is still spatial;
    cf.\ e.g.\ the resulting complication mentioned below \cite[Th.~22.5]{diamonds}.
    In our setting, we will be able to work around this issue.
\end{warning}

\begin{lemma}\label{open-immersion-proper-pushforward}
Let $j \colon X \hookrightarrow Y$ be an open immersion of diamonds and $j_!\colon \cD(X_\qproet; \ZZ/n) \to \cD(Y_\qproet; \ZZ/n)$ be the left adjoint to $j^*\colon \cD(Y_\qproet; \ZZ/n) \to \cD(X_\qproet; \ZZ/n)$. Then $j_!$ preserves $\cD_\et^+$.
\end{lemma}
\begin{proof}
    First, \cite[\href{https://stacks.math.columbia.edu/tag/04H1}{Tag 04H1}]{stacks-project} and \cite[Th.\,14.12(i)]{diamonds} imply that the claim is quasi-pro-\'etale locally on $Y$. So we can assume that $X$ and $Y$ are  perfectoid spaces. In this case, we have $\cD^+_\et(X; \ZZ/n)\simeq \cD^+(X_\et; \ZZ/n)$ and $\cD^+_\et(Y; \ZZ/n)\simeq \cD^+(Y_\et; \ZZ/n)$ (see \cite[Prop.\,14.15]{diamonds}). Therefore, the question reduces to showing that the natural transformation $j_{!} \lambda^*_{\qproet,X} \to \lambda^*_{\qproet,Y} j_{\et,!}$ is an isomorphism, where $\lambda_{\qproet, X} \colon X_{\qproet} \to X_\et$ (resp.\,$\lambda_{\qproet, Y}\colon Y_\qproet \to Y_\et$) is the natural projection of topoi and $j_{\et, !}$ is the left adjoint to $j^*_\et\colon \cD(Y_\et; \ZZ/n) \to \cD(X_\et; \ZZ/n)$. This can be verified for the associated transformation of right adjoints $j^*_\et \lambda_{\qproet,Y,*} \to \lambda_{\qproet,X,*}j^*_\qproet$, which is formal.
\end{proof}

\begin{definition}[{cf.\ \cite[Def.~22.4]{diamonds}}]
\label{proper-pushforward}
    Let $f \colon X \to Y$ be a compactifiable quasi-compact morphism of diamonds.
    Let $j \colon X \hookrightarrow \overline{X}^{/Y}$ and $\overline{f}^{/Y} \colon \overline{X}^{/Y} \to Y$ be the induced morphisms.
    Then using \cref{etale-pushforward} and
    \cref{open-immersion-proper-pushforward}, the \emph{proper pushforward} along $f$ is the functor
    \[ \rR f_! \colonequals \rR \overline{f}^{/Y}_{\qproet, *} \circ j_! \colon \cD^+_\et(X;\ZZ/n) \to \cD^+_\et(Y;\ZZ/n). \]

    If $Y = \Spa(C,C^+)$ for some algebraically closed nonarchimedean field $C$ and bounded open valuation subring $C^+ \subset C$, we also sometimes use the notation
    \[ \rR\Gamma_c(X,\blank) \colonequals \rR f_!(\blank) \colon \cD^+_\et(X;\ZZ/n) \to \cD^+_\et(Y;\ZZ/n) \simeq \cD^+(\ZZ/n). \]
\end{definition}

\begin{rmk}\label{Huber-diamond-comp-supp}
Let $f\colon X \to Y$ be a +-weakly finite type separated morphism of locally noetherian analytic adic spaces. Then $\rR f_!\colon \cD^+(X_\et; \ZZ/n) \to \cD^+(Y_\et; \ZZ/n)$ coincides with $(\rR f_!)^{\lozenge}$ under the equivalence of bounded below \'etale derived categories from \cite[Lem~15.6]{diamonds} and \cite[Rmk.~14.14]{diamonds}. Indeed, \cref{canonical-compactification-Huber} implies that it suffices to treat the case of an open immersion and proper morphisms separately. In the former case, both $j_!$ and $(j_!)^\lozenge$ are left adjoints to $j^*$, and in the latter case, both $\rR f_! = \rR f_*$ and $\rR f^\lozenge_! = \rR f^\lozenge_{\qproet, *} =\rR f^{\lozenge}_{\et, *}$ are right adjoints to $f^*$. 
\end{rmk}

\begin{remark}\label{not-quasi-compact}
    Assume that $X$ is a locally spatial diamond, $Y$ is a quasi-separated locally spatial diamond, and $f \colon X \to Y$ is a compactifiable morphism (not necessarily quasi-compact).
    Then the proper pushforwards in the sense of \cref{proper-pushforward} on the restrictions to opens with quasi-compact restriction $\restr{f}{V} \colon V \xhookrightarrow{j_V} X \to Y$ are all compatible by \cref{proper-pushforward-composition} below.
    Therefore, they define a functor $\rR f_! \colon \cD^+_{\et,\qc/Y}(X; \ZZ/n) \to \cD^+_\et(Y; \ZZ/n)$ on ``complexes with quasi-compact support over $Y$.''
    Since $X$ has a basis of quasi-compact opens and $Y$ is quasi-separated, any $F \in \cD^+_\et(X; \ZZ/n)$ is a filtered colimit of $j_{V,!}j^*_V F \in \cD^+_{\et,\qc/Y}(X; \ZZ/n)$, with $j_V \colon V \hookrightarrow X$ as before.\footnote{
    Note that the statement would in general be false if $Y$ is not quasi-separated:
    For example, let $K^+$ be an infinite rank microbial valuation ring such that the complement of the closed point in $\Spec(K^+)$ is not quasi-compact.
    Set $K \colonequals \Frac(K^+)$ and $X \colonequals \Spa(K,K^+)$, let $U \subset X$ be the complement of the closed point, and let $Y$ be the gluing of two copies of $X$ along $U$.
    Then for either of the open immersions $X \hookrightarrow Y$, there is no basis of opens of $X$ which are quasi-compact over $Y$.}
    Similarly to \cite[Def.~22.13]{diamonds}, one could thus define $\rR f_!\colon \cD^+(X_\et; \ZZ/n) \to \cD^+(Y_\et; \ZZ/n)$ by left Kan extension from $\rR f_! \colon \cD^+_{\et,\qc/Y}(X_\et; \ZZ/n) \to \cD^+_\et(Y; \ZZ/n)$, which would be informally given by the formula $\rR f_! F = \colim_{V} \rR (f\circ j_V)_! j_V^* F$.
    Under the diamondification functor, this again agrees with the proper pushforward functor defined by Huber in \cite[\S~5.4]{Huber-etale} for all locally of +-weakly finite type separated taut morphisms of locally noetherian analytic adic spaces (which are compactifiable by \cite[Th.~5.1.5]{Huber-etale}).
    Here, we will focus on the case of quasi-compact $f$ because it is the one relevant for our argument and allows us to avoid the $\infty$-categorical subtleties involved in the left Kan extension.
\end{remark}

We discuss some of the properties of the proper pushforward, following \cite[\S~22]{diamonds}.
\begin{lemma}[{cf.\ \cite[Prop.~22.9]{diamonds}}]
\label{proper-pushforward-composition}
    Let $f \colon X \to Y$ and $g \colon Y \to Z$ be compactifiable quasi-compact morphisms of diamonds.
    Then there is a natural equivalence
    \[ \rR g_! \circ \rR f_! \simeq \rR(g \circ f)_! \]
    of functors $\cD^+_\et(X;\ZZ/n) \to \cD^+_\et(Z;\ZZ/n)$.
\end{lemma}
\begin{proof}
    Unwinding the constructions and applying \cite[Cor.~18.8.(iii),(vi)]{diamonds}, we get a commutative diagram of diamonds
    \[ \begin{tikzcd}[column sep=large]
        X \arrow[r,hook,"j_1"] \arrow[d,"f"'] & \overline{X}^{/Y} \arrow[r,hook,"j_2"] \arrow[ld,near start,"\overline{f}^{/Y}"] & \overline{X}^{/Z} \arrow[ld,"\overline{f}^{/Z}"] \\
        Y \arrow[r,hook,"j_3"] \arrow[d,"g"'] & \overline{Y}^{/Z} \arrow[ld,"\overline{g}^{/Z}"] & \\
        Z &&
    \end{tikzcd} \]
    in which the horizontal maps are open immersions and the diagonal maps are proper.
    By definition,
    \[ \rR g_! \circ \rR f_! = \rR \overline{g}^{/Z}_{\qproet, *} \circ j_{3,!} \circ \rR \overline{f}^{/Y}_{\qproet, *} \circ j_{1,!} \quad \text{and} \quad \rR (g \circ f)_! = \rR \overline{g}^{/Z}_{\qproet, *} \circ \rR \overline{f}^{/Z}_{\qproet, *} \circ j_{2,!} \circ j_{1,!}. \]
    Therefore, it suffices to produce a natural equivalence
    \[ j_{3,!} \circ \rR\overline{f}^{/Y}_{\qproet, *} \simeq \rR\overline{f}^{/Z}_{\qproet, *} \circ j_{2,!} \]
    of functors $\cD^+_\et\bigl(\overline{X}^{/Y};\ZZ/n\bigr) \to \cD^+_\et\bigl(\overline{Y}^{/Z};\ZZ/n\bigr)$.
    This is exactly the statement of the proper base change theorem \cite[Prop.~3.3]{AM} (whose proof is similar to that of \cite[Th.~19.2]{diamonds}).
\end{proof}
We have the following analog of \cref{pushforward-stalks} for proper pushforwards.
\begin{lemma}[{cf.\ \cite[Prop.~22.8]{diamonds}}]\label{proper-pushforward-stalks}
    Let $f \colon X \to Y$ be a compactifiable quasi-compact morphism of diamonds and $F \in \cD^+_\et(X;\ZZ/n)$.
    Let $\bar{y} \colon \Spa(C,C^+) \to Y$ be a qpro\'et-geometric point of $Y$ and set $X_{\bar{y}} \colonequals X \times_{Y,\bar{y}} \Spa(C,C^+)$.
    Then there is a natural equivalence in $\cD_\et(\Spa(C,C^+);\ZZ/n)$
    \[ (\rR f_! F)_{\bar{y}} \simeq \rR\Gamma_c(X_{\bar{y}},F_{\overline{y}}). \]
\end{lemma}
\begin{proof}
    The natural map $\overline{(X_{\bar{y}})}^{/\Spa(C,C^+)} \to \overline{X}^{/Y} \times_Y \Spa(C,C^+)$ is an isomorphism, as one can check directly on $(R,R^+)$-valued points for all totally disconnected perfectoid $\Spa(R,R^+)$.
    Therefore, the statement follows from a combination of \cref{pushforward-stalks} and \cite[Prop.~19.1]{diamonds}.
\end{proof}
\begin{lemma}[{cf.\ \cite[Prop.~22.12]{diamonds}}]\label{direct-sums}
    Let $f \colon X \to Y$ be a compactifiable quasi-compact morphism of diamonds.
    Let $F_j \in \cD^{\ge 0}_\et(X;\ZZ/n)$ be a collection of objects indexed by $j\in J$.
    Then the natural map $\bigoplus_{j \in J} \rR f_! F_j \to \rR f_!(\bigoplus_{j \in J} F_j)$ in $\cD^+_\et(Y;\ZZ/n)$ is an isomorphism.
\end{lemma}
\begin{proof}
    Thanks to \cref{enough-points} and \cref{proper-pushforward-stalks}, it suffices to show that for any qpro\'et-geometric point $\bar{y} \colon \Spd(C,C^+) \to Y$, the induced morphism $\bigoplus_{j \in J} \rR\Gamma_c(X_{\bar{y}},F_j) \to \rR\Gamma_c(X_{\bar{y}},\bigoplus_{j \in J} F_j)$ is an isomorphism.
    Since extension by zero along $X_{\bar{y}} \to \overline{(X_{\bar{y}})}^{/\Spd(C,C^+)}$ preserves direct sums and is $t$-exact, it suffices to show that the functor
    \[
\rR\Gamma_\qproet\Bigl(\overline{(X_{\bar{y}})}^{/\Spd(C,C^+)},\blank\Bigr) \colon \cD^{\geq 0}\Bigl(\overline{(X_{\bar{y}})}^{/\Spd(C,C^+)}_\qproet;\ZZ/n\Bigr) \to \cD^{\geq 0}(\ZZ/n)
    \]
    commutes with infinite direct sums. Now an easy induction reduces to the case when each $F_j$ is concentrated in degree $0$. Since compactifiable morphisms are separated, \cite[Cor.~18.8]{diamonds} ensures that $\overline{(X_{\bar{y}})}^{/\Spd(C,C^+)} \to \Spd(C, C^+)$ is proper (in particular, it is qcqs). Thus the desired commutation with direct sums follows from the commutation with finite direct sums and \cite[\href{https://stacks.math.columbia.edu/tag/0739}{Tag~0739}]{stacks-project}. 
\end{proof}
\begin{proposition}[{cf.\ \cite[Th.~22.5]{diamonds}}]\label{cohomological-amplitude}
    Fix $d \in \NN$.
    Let $f \colon X \to Y$ be a compactifiable quasi-compact morphism of diamonds over $\Spd(\ZZ_p, \ZZ_p)$. Assume that for every qpro\'et-geometric point $\bar{y} \colon \Spa(C(\bar{y}),C(\bar{y})^+) \to Y$, the fiber $X_{\bar{y}}$ is the diamond associated with a separated noetherian analytic adic space of +-weakly finite type over $\Spa(C(\bar{y})^\sharp,C(\bar{y})^{\sharp, +})$ of relative dimension $\le d$ (see \cite[Def.,~1.8.1(ii)]{Huber-etale}).
    Then $\rR f_! \colon \cD^+_\et(X;\ZZ/n) \to \cD^+_\et(Y;\ZZ/n)$ has cohomological amplitude $[0,2d]$;
    that is, for any $F \in \cD^+_\et(X;\ZZ/n)$ concentrated in degree $0$, we have $\rR^if_! F = 0$ for $i \not\in [0,2d]$.
\end{proposition}

\begin{proof}
    Let $F \in \cD^+_\et(X;\ZZ/n)$ be concentrated in degree $0$.
    By \cref{enough-points} and \cref{proper-pushforward-stalks}, it suffices to show that $\Hh^i_c(X_{\bar{y}},F) = 0$ for $i \not\in [0,2d]$ and every qpro\'et-geometric point $\ov{y}$ of $Y$.
    We may therefore assume that $Y = \Spa(C,C^+)$ for some algebraically closed nonarchimedean field $C$ and bounded open valuation subring $C^+ \subset C$ and that $X$ is the spatial diamond attached to a separated noetherian analytic adic space of +-weakly finite type over $\Spa(C^\sharp, C^{\sharp, +})$ of relative dimension $\le d$.
    In this case, the definition of compactly supported cohomology from \cref{proper-pushforward} agrees with Huber's definition from \cite[Th.~5.4.3]{Huber-etale} (see \cref{Huber-diamond-comp-supp}), so the assertion reduces to the analogous one in Huber's setup, which is \cite[Prop.~5.5.8]{Huber-etale}.
\end{proof}
\begin{corollary}\label{derived-abelian}
    Let $d \in \NN$ and $f \colon X \to Y$ be as in \cref{cohomological-amplitude}.
    \begin{enumerate}[leftmargin=*,label=\upshape{(\roman*)}]
        \item\label{derived-abelian-factorization} Any morphism $\rR f_! \underline{\ZZ/n}(d)[2d] \to \underline{\ZZ/n}$ uniquely factors as
        \[ \rR f_! \underline{\ZZ/n}(d)[2d] \xlongrightarrow{\tau^{\ge 0}} \rR^{2d}f_! \underline{\ZZ/n}(d) \longrightarrow \underline{\ZZ/n}. \]
        \item\label{derived-abelian-sequence} Let $X = \bigcup_{j \in J} U_j$ be a retrocompact open covering by subdiamonds.
        Denote by $f_j \colon U_j \to Y$ and $f_{jk} \colon U_j \cap U_k \to Y$ the induced morphisms.
        Then the following sequence of natural maps is exact:
        \[ \bigoplus_{j,k \in J} \rR^{2d}f_{jk,!} \underline{\ZZ/n}(d) \to \bigoplus_{j \in J} \rR^{2d}f_{j,!} \underline{\ZZ/n}(d) \to \rR^{2d}f_!\underline{\ZZ/n}(d) \to 0 \]
    \end{enumerate}
\end{corollary}
\begin{proof}
    Part \cref{derived-abelian-factorization} follows immediately from the fact that $\rR f_! \underline{\ZZ/n}(d)[2d] \in D^{[-2d,0]}_\et(Y;\ZZ/n)$ thanks to \cref{cohomological-amplitude}.
    Part \cref{derived-abelian-sequence} is a direct consequence of \Cref{proper-pushforward-composition},
    \cref{direct-sums}, and the right-exactness of $\rR^{2d}f_!(\blank)$.
\end{proof}

\subsection{Trace maps for sousperfectoid spaces}
\label{traces}
Our next goal is to define trace maps for the period sheaves $\BB_I$ for a certain class of smooth proper morphisms between diamonds. This construction will rely heavily on the trace maps for $\ZZ/n$-coefficients from \cite[Th.~1.2.1]{LRZ24}. However, the setting of locally noetherian analytic adic spaces is not flexible enough for our applications because we will later apply it to smooth spaces over a perfectoid base. Therefore, we dedicate this subsection to extending the statement for $\ZZ/n$-coefficients from \cite[Th.~1.2.1]{LRZ24} to the generality of sousperfectoid analytic adic spaces in the sense of \cite[\S 7]{Hansen-Kedlaya} and \cite[\S IV.4.1]{Fargues-Scholze}.
Let us recall this notion.

\begin{definition}[{\cite[Def.~7.1]{Hansen-Kedlaya}}]
    An analytic adic space $X$ with $p \in \cO^{\circ\circ}(X)$ is called \emph{sousperfectoid} if there exists an affinoid open covering $X = \bigcup_{j \in J} \Spa(R_j, R^+_j)$ and for each $j \in J$ a continuous algebra homomorphism $R_j \to \widetilde{R}_j$ to a perfectoid $R_j$-algebra which splits in the category of topological $R_j$-modules.
\end{definition}
We caution the reader that it is unknown whether the functions on a sousperfectoid space $X = \Spa(R,R^+)$ which is affinoid admit a split injection $R \to \widetilde{R}$ to a perfectoid $R$-algebra (so one may have to pass to an open covering even in that case).
Furthermore, it is even unclear whether any nonarchimedean field of mixed characteristic is sousperfectoid;
however, this holds if the field is perfectoid (trivial) or topologically countably generated over $\QQ_p$ (\cite[Rmk.~7.8]{Hansen-Kedlaya}).
On the other hand, sousperfectoid spaces enjoy useful permanence properties under several operations:

\begin{proposition}[{\cite[Prop.~3.8, Lem.~7.5, Lem.~7.3]{Hansen-Kedlaya}}]\label{sousperfectoid-permanence}
    Let $X$ be a sousperfectoid space.
    Then the following spaces are sousperfectoid as well:
    \begin{enumerate}[label=\upshape{(\alph*)}]
        \item any \'etale\footnote{We recall that a morphism $Y\to X$ is \'etale if for each $y\in Y$, there is an open neighborhood $y\in U\subset Y$ such that $f|_U \colon U \to X$ factors as a composition of open immersions and finite \'etale maps, see \cite[Def.\,8.2.16]{Kedlaya-Liu-1}.} $Y\to X$,
        \item the relative unit disk $\DD^1_X \colonequals \Spa(\ZZ[T],\ZZ[T]) \times_{\Spa(\ZZ,\ZZ)} X \to X$
    \end{enumerate}
\end{proposition}
By virtue of the local structure of smooth morphisms between locally noetherian analytic adic spaces \cite[Cor.~1.6.10]{Huber-etale}, \cref{sousperfectoid-permanence}
implies that for any smooth morphism $Y \to X$ of locally noetherian analytic adic spaces,
if $X$ is sousperfectoid, then $Y$ is sousperfectoid as well.
This local structure also motivates the following definition of smooth morphisms between general sousperfectoid spaces:
\begin{definition}
    A morphism $f \colon X \to Y$ of sousperfectoid spaces is called \emph{smooth} if there exists an open covering $X = \bigcup_{i \in I} U_i$ and \'etale maps $U_i \to \DD^{d_i}_Y$ for some $d_i \in \ZZ_{\ge 0}$ such that $\restr{f}{U_i}$ factors as $U_i \to \DD^{d_i}_Y \to Y$.  
\end{definition}
We will again need some permanence properties.
\begin{lemma}[{\cite[Prop.~IV.4.10]{Fargues-Scholze}}]\label{sousperfectoid-properties}
    \begin{enumerate}[label=\upshape{(\roman*)},leftmargin=*]
        \item Let $f \colon X \to Y$ and $g \colon Y \to Z$ be smooth morphisms of sousperfectoid spaces.
        Then $g \circ f$ is smooth as well.
        \item\label{sousperfectoid-properties-product} Let $X \to Y$ be a smooth morphism and $Y' \to Y$ be an arbitrary morphism of sousperfectoid spaces.
        Then fiber product $X \times_Y Y'$ exists in the category of adic spaces, is sousperfectoid, and the natural morphism $X \times_Y Y' \to Y'$ is again smooth.
    \end{enumerate}
\end{lemma}

For later reference, we also record the following fact: 

\begin{lemma}\label{lemma:uniqueness-untilt} Let $S$ be a perfectoid space over $\Spa(\QQ_p, \ZZ_p)$, let $\rm{Sm}_{/S}$ be the category of smooth sousperfectoid $S$-spaces, and let $\rm{Diam}_{/S^\flat}$ be the category of diamonds over $S^\flat=S^\diam$. Then the functor
\[
(\blank)^\diam \colon \rm{Sm}_{/S} \to \rm{Diam}_{/S^\flat}
\]
is fully faithful. 
\end{lemma}
\begin{proof}
    Pick $X, Y\in \rm{Sm}_{/S}$. We wish to show that the natural morphism $\Hom_S(X, Y) \to \Hom_{S^\flat}(X^\diam, Y^\diam)$ is an isomorphism. Since $|X|=|X^\diam|$, $|Y|=|Y^\diam|$, $|S|=|S^\flat|$ \cite[Lem.\,15.6]{diamonds}, a standard gluing argument allows us to reduce to the case when  $X=\Spa(B, B^+)$, $Y=\Spa(A, A^+)$ and $S=\Spa(R, R^+)$ are affinoid spaces. Furthermore, we can localize even further to assume that the morphisms $(R, R^+) \to (A, A^+)$ and $(R, R^+) \to (B, B^+)$ are smooth in the na\"{\i}ve sense in the sense of \cite[Def.\,5.11]{Hansen-Kedlaya}, i.e., both morphisms can be obtained as compositions of rational localizations, finite \'etale morphisms, and Tate algebra extensions. 
    
    In this case, \cite[Prop.\,2.1]{Huber-generalization} implies that $\Hom_S(X, Y)$ is in bijection with the set of $(R, R^+)$-linear continuous homomorphisms $(A, A^+) \to (B, B^+)$. While the set $\Hom_{S^\flat}(X^\diam, Y^\diam)$ is in bijection with the set of $(R, R^+)$-linear continuous homomorphisms $(A, A^+) \to (\O_{X^\diam}(X^{\diam}), \O^+_{X^\diam}(X^\diam))$, see \cref{defn:period-sheaves}. Therefore, it suffices to show that the natural morphism
    \[
    (B, B^+) \to (\O_{X^\diam}(X^{\diam}), \O^+_{X^\diam}(X^\diam))
    \]
    is a homeomorphism, i.e., $X$ is $v$-complete in the sense of \cref{defn:v-complete}. This follows from \cite[Th.\,11.18]{Hansen-Kedlaya}. 
\end{proof}

With these preliminaries in place, we now extend the construction of traces with $\ZZ/n$-coefficients to the setting of smooth morphisms of sousperfectoid spaces.
To make sense of the corresponding \'etale sites, we will work with the associated diamonds using \cite[Lem.~15.6]{diamonds} and put the theory explained in \cref{diamond-preliminaries} to use.

To simplify notation, we will conflate $X$ with $X^\lozenge$ and drop the symbols $(\blank)^\lozenge$ throughout. In particular, for a separated $+$-weakly finite type morphism $f\colon X \to Y$ of sousperfectoid spaces, we write $\rR f_!$ instead of $\rR f_!^\diam$; we refer to \cref{analytic-adic-compactifiable} and \cref{proper-pushforward} for the existence of $\rR f_!^\diam$. If $f\colon X \to Y$ is a map between locally noetherian analytic adic spaces, then \cref{Huber-diamond-comp-supp} ensures $\rR f_!$ coincides with $\rR f_!$ from \cite{Huber-etale}. We also write $\rR f_{\qproet, *}$ instead of $\rR f^\diam_{\qproet, *}$. 

\begin{proposition}\label{sousperfectoid-trace}
    There is a unique way to assign to any quasi-compact separated smooth of equidimension $d$ morphism $f \colon X \to Y$ of sousperfectoid spaces with $n \in \cO^\times_Y$ a trace map $\tr_f \colon \rR f_!\underline{\ZZ/n}(d)[2d] \to \underline{\ZZ/n}$ in $\cD^+_\et(Y;\ZZ/n)$ such that for any $y \in \abs{Y}$ and qpro\'et-geometric point $\bar{y} \colon \Spa(C,C^+) \to Y$ over $y$, the pullback \[(\tr_f)_{\bar{y}} \colon (\rR f_!\underline{\ZZ/n}(d)[2d])_{\bar{y}} \simeq \rR \Gamma_c(X_{\bar{y}},\underline{\ZZ/n}(d))[2d] \to \ZZ/n\] agrees under the equivalence from \cref{Huber-diamond-comp-supp} with the trace map for the morphism of noetherian analytic adic spaces $f_{\bar{y}} \colon X_{\bar{y}} \to \Spa(C,C^+)$ constructed in \cite[Th.~1.2.1]{LRZ24}. 
\end{proposition}
Here, the base change isomorphism $(\rR f_!\underline{\ZZ/n}(d)[2d])_{\bar{y}} \simeq \rR \Gamma_c(X_{\bar{y}},\underline{\ZZ/n}(d))[2d]$ follows from \cref{proper-pushforward-stalks}.
Before giving the proof of \cref{sousperfectoid-trace}, we explain the example of $\PP^{1,\an}_Y \to Y$;
this special case will also play an important role in the treatment of the general statement.
\begin{construction}\label{c1-construction}
    Let $Y$ be sousperfectoid space with $n \in \cO^\times_Y$ and $f \colon \PP^{1,\an}_Y \to Y$ be the relative analytic projective line (see \cite[Prop.~IV.4.9]{Fargues-Scholze}).
    Then there is a natural element $\cO_{\PP^{1,\an}_Y}(1) \in \Pic(\PP^{1,\an}_Y) \simeq \Hh^0(Y, \rR f_{\et, *}\GG_m[1])$, given as the pullback of $\cO_{\PP^1_\ZZ}(1) \in \Pic(\PP^1_\ZZ)$ under the natural map of locally ringed spaces $\PP^{1,\an}_Y \to \PP^{1,\an}_{\Spa(\ZZ,\ZZ)} \to \PP^1_\ZZ$. 
    This defines a morphism $\ud{\ZZ} \to \rR f_{\et, *} \GG_m[1]$.
    Since $n \in \cO^\times_Y$, the Kummer sequence produces a map of \'etale sheaves $\rR f_{\et, *}\GG_m[1] \to \rR f_{\et, *}\mu_n[2]$.
    Composing, we obtain a natural map $\underline{\ZZ} \to \rR f_{\et, *}\mu_n[2]$ which uniquely determines a map $c\colon \ud{\ZZ/n} \to \rR f_{\et, *}\mu_n[2]$.
    Combined with \cite[Lem.~15.6 and Prop.~16.6]{diamonds}, we get a map
    \[
    c\colon \ud{\ZZ/n} \to \rR f_{\qproet, *}\mu_n[2].
    \]
\end{construction}
\begin{lemma}\label{P1-trace}
    Let $Y$ be sousperfectoid space with $n \in \cO^\times_Y$ and $f \colon \PP^{1,\an}_Y \to Y$ be the relative analytic projective line.
    Let $\eta \colon \ud{\ZZ/n} \to \rR f_{\qproet, *}\ud{\ZZ/n}$ be the natural morphism coming from adjunction and $c$ be the one from \cref{c1-construction}.
    Then
    \[ \eta(1)[2] \oplus c \colon \ud{\ZZ/n}(1)[2] \oplus \underline{\ZZ/n} \xlongrightarrow{\sim} \rR f_{\qproet, *} \ud{\ZZ/n}(1)[2] \]
    is an isomorphism such that for any $y \in \abs{Y}$ and qpro\'et-geometric point $\bar{y} \colon \Spa(C,C^+) \to Y$ over $y$, the induced projection $\bigl(\rR f_{\qproet, *} \underline{\ZZ/n}(1)[2]\bigr)_{\bar{y}} \simeq \rR\Gamma_\et\bigl(\PP^{1,\an}_{(C,C^+)},\underline{\ZZ/n}(1)\bigr)[2] \to \ZZ/n$ agrees under the equivalence of \'etale sites from \cite[Lem.~15.6]{diamonds} with the trace map from \cite[Th.~1.2.1]{LRZ24}.
\end{lemma}
\begin{proof}
    The Kummer sequence and the unit of adjunction $\underline{\ZZ/n} \to \rR f_{\qproet, *}\underline{\ZZ/n}$ are both compatible with pullbacks, so $\eta_{\bar{y}}$ and $c_{\bar{y}}$ are again given by the same construction on the fibers under the isomorphism from \cref{pushforward-stalks}.
    In this locally noetherian setting, the resulting map is known to be an isomorphism (cf.\ e.g.\ \cite[Prop.~6.1.6]{zav-6-functors}), so we can conclude that $\eta(1)[2] \oplus c$ is an isomorphism (thanks to \cref{enough-points}.\cref{enough-points-conservative}).
    Moreover, as the algebraic trace map in \cite[Exp.~XVIII, Th.~2.9]{SGA4} (which in this case ultimately reduces to \cite[Exp.~XVIII, Prop.~1.1.6]{SGA4}) or \cite[p.~424]{Fu-Etale} is constructed via the algebraic Kummer sequence and the algebraic Kummer sequence is compatible with the analytic one under the analytification map, we conclude that the projection given by $``c^{-1}_{\bar{y}}\text{''}\colon \rR\Gamma_\et\bigl(\PP^{1,\an}_{(C,C^+)},\underline{\ZZ/n}(1)\bigr)[2] \to \ZZ/n$ is the analytification of the algebraic trace map and therefore coincides with the one from \cite[Th.~1.2.1]{LRZ24}.
\end{proof}
\begin{proof}[{Proof of \cref{sousperfectoid-trace}}]
    Note that $f^\lozenge$ is compactifiable by \cref{analytic-adic-compactifiable}, so the functor $\rR f_!$ is well-defined. We first show uniqueness. Since $f$ is of equidimension $d$, \cref{derived-abelian}.\cref{derived-abelian-factorization} implies that maps $\rR f_!\underline{\ZZ/n}(d)[2d] \to \underline{\ZZ/n}$ amount to maps of abelian sheaves of $\ZZ/n$-modules $\rR^{2d}f_!\underline{\ZZ/n}(d) \to \underline{\ZZ/n}$. Thus, \cref{enough-points} guarantees that one can check on stalks at all $\bar{y}$ whether two maps $\rR f_!\,\underline{\ZZ/n}(d)[2d] \to \underline{\ZZ/n}$ coincide.
    This yields the uniqueness.

    To prove the existence statement, we first construct the trace when there exists an \'etale map $X \to \DD^d_Y$ such that $f$ factors as $X \to \DD^d_Y \to Y$.
    In this case, we get the following commutative diagram of analytic adic spaces in which $j$ is qcqs and \'etale:
    \[ \begin{tikzcd}
        X \arrow[r] \arrow[rd,"j"] \arrow[d,"f"] & \DD^d_Y \arrow[d,hook] \\
        Y & \bigl(\PP^{1,\an}_Y\bigr)^d \arrow[l,"\pi_Y"]
    \end{tikzcd} \]
    Applying \cref{P1-trace} iteratively to the projections $\bigl(\PP^{1,\an}_Y\bigr)^{r+1} \to \bigl(\PP^{1,\an}_Y\bigr)^r$, we obtain a natural map
    \[ \tr_{\pi_Y} \colon \rR \pi_{Y,!}\underline{\ZZ/n}(d)[2d] \simeq \rR \pi_{Y,\qproet, *}\underline{\ZZ/n}(d)[2d] \xlongrightarrow{\tau^{\ge 0}} \rR^{2d}\pi_{Y,\qproet, *}\underline{\ZZ/n}(d) \xlongrightarrow{\sim} \underline{\ZZ/n}. \]
    Moreover, for any $y \in \abs{Y}$ and qpro\'et-geometric point $\bar{y}$ over $y$, the stalk $(\tr_{\pi_Y})_{\bar{y}}$ identifies with the algebraic trace map under the base change isomorphism from \cref{proper-pushforward-stalks}, and hence with the trace map from \cite{LRZ24} by properties (1) and (4) of \cite[Th.~1.2.1]{LRZ24}.
    
    On the other hand, the natural counit for the adjunction between $j_!$ and $j^*$ gives a map
    \[ \tr_j \colon j_!\underline{\ZZ/n} \to \underline{\ZZ/n}. \]
    Since the adjunction of $j_!$ and $j^*$ is compatible with passing to stalks at any qpro\'et-geometric point $\bar{y}$ of $Y$, the stalk $(\tr_j)_{\bar{y}}$ identifies with the trace map from \cite{LRZ24} under the base change isomorphism from \cref{proper-pushforward-stalks} by property (3) of \cite[Th.~1.2.1]{LRZ24}.
    Thanks to \cref{proper-pushforward-composition}, the composition yields the desired map
    \[ \tr_f \colon \rR f_!\underline{\ZZ/n}(d)[2d]  \simeq \rR\pi_{Y,!}j_!\underline{\ZZ/n}(d)[2d] \xlongrightarrow{\rR\pi_{Y,!}(\tr_j(d)[2d])} \rR\pi_{Y,!}\underline{\ZZ/n}(d)[2d] \xlongrightarrow{\tr_{\pi_Y}} \underline{\ZZ/n}, \]
    whose stalks are again compatible with the traces from \cite{LRZ24}.
    This finishes the construction in the case when $f$ factors through $\DD^d_Y$.

    Returning to the general case, \cref{derived-abelian}\cref{derived-abelian-factorization} implies that the maps $\rR f_!\underline{\ZZ/n}(d)[2d] \to \underline{\ZZ/n}$ amount to maps of abelian sheaves of $\ZZ/n$-modules $\rR^{2d}f_!\underline{\ZZ/n}(d) \to \underline{\ZZ/n}$. Thus, the uniqueness claim established above allows us to reduce to the case when $Y$ is an affinoid space. In particular, it is qcqs. Then the definition of smoothness implies that we can find a covering $X = \bigcup U_i$ by retrocompact open subspaces such that $f_i \colonequals \restr{f}{U_i}$ factors through $\DD^d_Y$.
    The previous construction therefore supplies unique trace maps
    \[ \sum_i \tr_{f_i} \colon \bigoplus_i \rR f_{i,!}\underline{\ZZ/n}(d)[2d] \longrightarrow \underline{\ZZ/n}. \]
    By \cref{derived-abelian}\cref{derived-abelian-sequence} and the uniqueness of the traces, this map descends to a map
    \[ \tr_f \colon \rR f_!\underline{\ZZ/n}(d)[2d] \longrightarrow \underline{\ZZ/n}. \]
    At all $\bar{y}$, the stalks $(\tr_f)_{\bar{y}}$ still coincide with the traces from \cite[Th.~1.2.1]{LRZ24} because those satisfy the same gluing property (cf.\ \cite[Lem.~6.1.3]{LRZ24}).
    This finishes the construction in the general case.
\end{proof}
\begin{remark}\label{trace-loc-noeth}
    If $f \colon X \to Y$ is a quasi-compact separated smooth of equidimension $d$ morphism of \emph{locally noetherian} sousperfectoid spaces with $n \in \cO^\times_Y$, then the trace map $\tr_f$ from \cref{sousperfectoid-trace} agrees with the one constructed in \cite[Th.~1.2.1]{LRZ24}.
    Using that both traces amount to maps of abelian sheaves $\rR^{2d}f_!\underline{\ZZ/n}(d) \to \underline{\ZZ/n}$, this can again be checked on stalks at qpro\'et-geometric points $\bar{y}$ over all $y \in \abs{Y}$ thanks to \cref{enough-points} and thus be reduced to the compatibility in \cref{sousperfectoid-trace}.
\end{remark}
\begin{remark}
\label{trace-compatibilities}
    The trace maps $\tr_f$ from \cref{sousperfectoid-trace} still satisfy the same compatibilities as the ones constructed in \cite[Th.~1.2.1]{LRZ24}:
    \begin{enumerate}[label=(\roman*)]
        \item\label{tr-comp} compatibility with compositions;
        \item\label{tr-pb} compatibility with pullbacks;
        \item\label{tr-et} compatibility with \'etale traces; and
        \item\label{tr-P1} compatibility with the algebraic trace for $\PP^1 \to \Spec C$ where $C$ is any complete, algebraically closed extension of $\QQ_p$. 
    \end{enumerate}
    The compatibility \cref{tr-P1} follows immediately from the construction.
    
    For compatibility \cref{tr-pb}, \cref{proper-pushforward-stalks} and the compatibility of the trace with fibers imply that it suffices to prove compatibility with pullbacks of the form $\Spa(C', C'^+) \to \Spa(C, C^+)$, where $C$ and $C'$ are algebraically closed nonarchimedean fields and $C^+\subset C$, $C'^+ \subset C'$ are open and bounded valuation subrings. Since these spaces are locally noetherian, compatibility follows from \cite[Th.~1.2.1]{LRZ24}.  
    
    Finally, compatibilities \cref{tr-comp} and \cref{tr-et} can once more be checked on stalks at qpro\'et-geometric points over all points of the base. In this case, it follows from \cite[Th.~1.2.1]{LRZ24}.
\end{remark}

\subsection{Trace maps for diamonds}

In this subsection, we extend the construction of the $\ZZ/n$-trace to smooth proper morphisms of diamonds (see \cref{defn:smooth-proper}). Using this, we then also construct a trace morphism for $\O$ and $\BB_I$-coefficients. 

Before we make these constructions, we need to introduce some new definitions regarding smooth and proper morphisms of sousperfectoid spaces and diamonds. 

\begin{defn}\label{defn:proper-sousperfectoid} A morphism of sousperfectoid spaces $f\colon X \to Y$ is \emph{proper} (resp.\,\emph{separated}) if it is qcqs, locally of $+$-weakly finite type (resp. quasi-separated, locally of $+$-weakly finite type), and for every algebraically closed nonarchimedean field $C$ with an open and bounded valuation subring $C^+\subset C$, and any diagram
\[
\begin{tikzcd}
    \Spa(C, \O_C)\arrow{d} \arrow{r} & X \arrow{d}{f} \\
    \Spa(C, C^+) \arrow[ru, dashed] \arrow{r} & Y,
\end{tikzcd}
\]
there exists a unique (resp.\,at most one) dotted arrow making the diagram commute.
\end{defn}

\begin{lemma}\label{lemma:different-proper} Let $f\colon X \to Y$ be a proper (resp.\,separated) morphism of sousperfectoid spaces or locally noetherian analytic adic spaces. Then $f^\diam \colon X^\diam \to Y^\diam$ is a proper (resp.\,separated) morphism of diamonds in the sense of \cite[Def.\,18.1]{diamonds}.
\end{lemma}
\begin{proof}
    We give a proof for proper morphisms; the case of separated morphisms is analogous. First, \cite[Lem.\,6.1.13]{zav-almost} implies that $f^\diam\colon X^\diam \to Y^\diam$ is a qcqs morphism. Therefore, \cite[Prop.\,18.3]{diamonds} implies that it suffices to show that $f^\diam$ satisfies the valuative criterion of properness. \emph{The proof of} \cite[Prop.\,18.3]{diamonds} ensures that it suffices to test the valuative criterion only for algebraically closed nonarchimedean fields $K$ with an open bounded valuation subring $K^+\subset K$. Then the result follows immediately from \cref{defn:proper-sousperfectoid} in the case of sousperfectoid spaces or from \cite[Lem.\,1.3.10]{Huber-etale} in the case of locally noetherian analytic spaces.
\end{proof}

Now we are ready to define the notion of a smooth proper morphism between diamonds over $\ZZ_p$.

\begin{defn}\label{defn:smooth-proper}
A morphism of diamonds $f\colon X \to Y$ over $\Spd(\QQ_p, \ZZ_p)$ is \emph{smooth proper (resp.~smooth proper of equidimension $d$, resp.~smooth proper of dimension $d$)} if, for every strictly totally disconnected perfectoid space $S$ with a morphism $S \to Y$, the fiber product $f_S\colon X_{S}\coloneqq X\times_Y S \to S$ is equal to the diamondification of a smooth proper morphism $f^\sharp_{S}\colon X_S^\sharp \to S^\sharp$ (resp.~smooth proper of equidimension $d$, resp.~smooth proper of relative dimension $d$). 
\end{defn}

\begin{rmk} \cref{lemma:uniqueness-untilt} ensures that $f^\sharp_S\colon X_S^\sharp \to S^\sharp$ is unique. Furthermore, it is functorial in $f$. 
\end{rmk}

\begin{example} Let $f\colon X \to Y$ be a smooth proper morphism of equidimension $d$ between analytic noetherian (resp.\,sousperfectoid) adic $\Spa(\QQ_p, \ZZ_p)$-spaces. Then $f^\diam \colon X^\diam \to Y^\diam$ is a smooth proper of equidimension $d$ morphism in the sense of \cref{defn:smooth-proper}.
\end{example}

Let $f\colon X \to Y$ be a smooth proper of equidimension $d$ morphism between diamonds over $\Spd(\QQ_p, \ZZ_p)$. Then \cref{proper-pushforward}, \cref{lemma:different-proper}, and \cref{qproet-v-pushforward-etale}, imply that 
\begin{equation}\label{eqn:lower-shriek-push}
\mu_Y^*(\rR f_{!} \ud{\ZZ/n\ZZ}) \simeq \mu_Y^*(\rR f_{\qproet, *}\ud{\ZZ/n\ZZ}) \simeq \rR f_{v, *} \ud{\ZZ/n\ZZ},
\end{equation}
where $\mu_Y\colon Y_v \to Y_\qproet$ is the natural projection. 

\begin{proposition}\label{prop:trace-diamonds} There is a unique way to assign to any smooth proper of equidimension $d$ morphism $f \colon X \to Y$ between diamonds over $\Spd(\QQ_p, \ZZ_p)$ a trace map $\tr_f \colon \rR f_{v, *}\underline{\ZZ/n}(d)[2d] \to \underline{\ZZ/n}$ in $\cD^+(Y_v;\ZZ/n)$ such that 
\begin{enumerate}[label={\upshape{(\arabic*)}}]
    \item\label{prop:trace-diamonds-1} for any pullback diagram of diamonds over $\Spd(\QQ_p, \ZZ_p)$
    \[
    \begin{tikzcd}
        X' \arrow{r}{g'} \arrow{d}{f'} & X \arrow{d}{f} \\
        Y' \arrow{r}{g} & Y,
    \end{tikzcd}
    \]
    in which $f$ and $f'$ are smooth proper of equidimension $d$, the following diagram is commutative:
    \[
    \begin{tikzcd}
        g^* \rR f_{v, *} \ud{\ZZ/n\ZZ}(d)[2d] \arrow{r}{g^*\tr_f} \arrow{d}{\wr} &  g^* \ud{\ZZ/n\ZZ} \arrow{d}{\wr} \\
        \rR f'_{v, *}\ud{\ZZ/n\ZZ}(d)[2d] \arrow{r}{\tr_{f'}} & \ud{\ZZ/n\ZZ}
    \end{tikzcd}
    \]
    \item\label{prop:trace-diamonds-2} if $S$ is a strictly totally disconnected perfectoid space over $\Spd(\Q_p, \Z_p)$ and $f\colon X \to S$ is a smooth proper morphism of equidimension $d$, then we have under the isomorphism \cref{eqn:lower-shriek-push}
    \[
    \tr_{f}=\mu^*_S(\tr_{f^\sharp})
    \]
    in $D(S_v;\ZZ/n)$, where $f^\sharp\colon X^\sharp \to S^\sharp$ is the morphism from \cref{defn:smooth-proper} and $\tr_{f^\sharp}$ is the trace map from \cref{sousperfectoid-trace}. 
\end{enumerate}
    Furthermore, the trace morphism satisfies the following compatibilities:
    \begin{enumerate}
        \item Let $S^\sharp$ be a sousperfectoid adic space over $\Spa(\QQ_p, \ZZ_p)$ and let $f^\sharp \colon X^\sharp \to S^\sharp$ be a smooth proper morphism of equidimension $d$. Set $S\coloneqq (S^\sharp)^\diam$, $X\coloneqq (X^\sharp)^\diam$, and $f\coloneqq (f^\sharp)^\diam$. Under the isomorphism \cref{eqn:lower-shriek-push}, we then have $\tr_f = \mu^*_S(\tr_{f^\sharp})$ in $D(S_v;\ZZ/n)$, where $\tr_{f^\sharp}$ is the trace map from \cref{sousperfectoid-trace}.
        \item Let $S^\sharp$ be a locally noetherian analytic adic space over $\Spa(\QQ_p, \ZZ_p)$ and let $f^\sharp \colon X^\sharp \to S^\sharp$ be a smooth proper morphism of equidimension $d$. Set $S\coloneqq (S^\sharp)^\diam$, $X\coloneqq (X^\sharp)^\diam$, and $f\coloneqq (f^\sharp)^\diam$. Let $\lambda_{S^\sharp} \colon S_v \to S_\et \simeq S^\sharp_\et$ be the natural projection morphism. Under the isomorphisms from \cref{eqn:lower-shriek-push} and \cite[Prop.~16.6]{diamonds}, we then have $\tr_f = \lambda^*_{S^\sharp}(\tr_{f^\sharp})$ in $D(S_v;\ZZ/n)$, where $\tr_{f^\sharp}$ is the trace map from \cite[Th.\,1.2.1]{LRZ24}.
    \end{enumerate}
\end{proposition}

\begin{proof} Fix a smooth proper of equidimension $d$ morphism $f\colon X \to Y$ between diamonds over $\Spd(\QQ_p, \ZZ_p)$.
   \begin{enumerate}[wide,label={\textit{Step~\arabic*}.},ref={Step~\arabic*}]
    \item\label{prop:trace-diamonds-step-1} \textit{The complex $\rR f_{v, *} \ud{\ZZ/n\ZZ}(d)[2d]$ is concentrated in degrees $[-2d, 0]$.} The claim is $v$-local on $Y$, so we can assume that $S=Y$ is a strictly totally disconnected perfectoid space. 
    In this case, \cref{defn:smooth-proper} guarantees that $f\colon X\to S$ is equal to $(f_S^\sharp)^\diam$ for a smooth proper morphism $f_S^\sharp \colon X_S^\sharp \to S^\sharp$ of equidimension $d$. Then the claim follows immediately from \cref{cohomological-amplitude} and \cref{eqn:lower-shriek-push}.

    \item\label{prop:trace-diamonds-step-2} \textit{The presheaf ${\cHom}\bigl(\rR f_{v, *} \ud{\ZZ/n\ZZ}(d)[2d], \ud{\ZZ/n\ZZ}\bigr)$ is a $v$-sheaf.} \cref{prop:trace-diamonds-step-1} implies that 
    \[
    \rR{\cHom}\bigl(\rR f_{v, *} \ud{\ZZ/n\ZZ}(d)[2d], \ud{\ZZ/n\ZZ}\bigr) \in \cD^{\geq 0}(Y_v; \ud{\ZZ/n\ZZ}). 
    \]
    Therefore, the claim follows from \cite[Prop.\,3.2.2]{BBDG}. 

    \item\label{prop:trace-diamonds-step-3} \textit{Construction and uniqueness of $\tr_f$.} For any smooth proper of equidimension $d$ morphism $f \colon X \to Y$, \cref{prop:trace-diamonds-step-2} implies that ${\cHom}\bigl(\rR f_{v, *} \ud{\ZZ/n\ZZ}(d)[2d], \ud{\ZZ/n\ZZ}\bigr)$ forms a sheaf on the $v$-site $Y_v$. We wish to construct an element $\tr_f\in {\cHom}\bigl(\rR f_{v, *} \ud{\ZZ/n\ZZ}(d)[2d], \ud{\ZZ/n\ZZ}\bigr)(Y)$ satisfying Conditions~\cref{prop:trace-diamonds-1} and \cref{prop:trace-diamonds-2}. Since strictly totally disconnected perfectoid spaces form a basis on $Y_v$, any such element is equivalent to a system of elements 
    \[
    \{\tr_{f, S}\in \Hom\bigl(\rR f_{S, v, *} \ud{\ZZ/n\ZZ}(d)[2d], \ud{\ZZ/n\ZZ}\bigr) = {\cHom}\bigl(\rR f_{v, *} \ud{\ZZ/n\ZZ}(d)[2d], \ud{\ZZ/n\ZZ}\bigr)(S) \ | \ \alpha^*_v(\tr_{f, S}) =\tr_{f, S'} \text{ for any $S'\xrightarrow{\alpha} S$}\}_{S\to Y},
    \]
    where the set is indexed by strictly totally disconnected perfectoid spaces $S$ with a map $S \to Y$. Conditions~\cref{prop:trace-diamonds-1} and \cref{prop:trace-diamonds-2} tell us that we must have $\tr_{f, S} = \mu_{S}^*(\tr_{f_S^\sharp})$. Therefore, $\tr_f$ is unique if it exists, and it exists if and only if 
    \[
    \alpha^*_v\bigl(\mu_{S}^*(\tr_{f_S^\sharp}) \bigr) = \mu_{S'}^*\bigl(\tr_{f_{S'}^\sharp}\bigr)
    \]
    for any morphism $\alpha \colon S' \to S$ of strictly totally disconnected perfectoid spaces over $Y$. To prove this, we note that \cref{lemma:uniqueness-untilt} implies that $X_{S}^\sharp \times_{S^\sharp} S'^\sharp \simeq X_{S'}^\sharp$ as adic $(S')^\sharp$-spaces. Therefore, the claim follows from the sequence of equalities 
    \[
    \alpha^*_v\bigl(\mu_{S}^*(\tr_{f_S^\sharp}) \bigr) = \mu_{S'}^* \bigl(\alpha^*_{\qproet}(\tr_{f_S^\sharp}) \bigr) = \mu_{S'}^*(\tr_{f_{S'}^\sharp}),
    \]
    where the second equality is due to \cref{trace-compatibilities}~\cref{tr-pb}. 

    This argument shows that we get a unique element $\tr_f\in \Hom\bigl(\rR f_{v, *} \ud{\ZZ/n\ZZ}(d)[2d], \ud{\ZZ/n\ZZ}\bigr)$ such that $(\tr_f)|_S = \mu_{S}^*(\tr_{f_S^\sharp})$ for any strictly totally disconnected perfectoid $S$ with a map $S \to Y$. By construction, this trace map satisfies Condition~\ref{prop:trace-diamonds-2}. To show that it also satisfies Condition~\ref{prop:trace-diamonds-1}, we see that the uniqueness property from the previous sentence ensures that it suffices to verify Condition~\ref{prop:trace-diamonds-1} when $Y'=S$ is a strictly totally disconnected perfectoid space. Then it follows automatically from the construction of $\tr_f$. 

    \item\label{prop:trace-diamonds-step-4} \textit{Compatibility with \cref{sousperfectoid-trace} and \cite{LRZ24}} We pick a sousperfectoid (resp.\,locally noetherian) analytic adic space $S^\sharp$ over $\Spa(\QQ_p, \ZZ_p)$ and a smooth proper morphism $f^\sharp \colon X^\sharp \to S^\sharp$ of equidimension $d$. We wish to show that $\mu^*_{S}(\tr_{f^\sharp})=\tr_{f}$ (resp. $\lambda_{S^\sharp}^*(\tr_{f^\sharp})=\tr_f$). Now \cref{qproet-v-pushforward-etale} and \cref{etale-pushforward} imply that $\rR f_{v, *} \ud{\ZZ/n\ZZ}(d)[2d]$ lies in $\cD_\et(S; \ud{\ZZ/n\ZZ})$. Since all trace maps commute with pullbacks (see \cref{pushforward-stalks} and \cref{trace-compatibilities}~\cref{tr-pb}, resp.\,\cite[Th.\,6.1.1(2) and Lem.\,6.3.1]{LRZ24}), we can check that two trace maps agree on each qpro\'et-geometric point of $S$ (see \cref{enough-points}).
    In other words, we can assume that $S=\Spa(C, C^+)$ for an algebraically closed nonarchimedean field $C$ and an open bounded valuation subring $C^+\subset C$. In this case, the trace map from \cref{sousperfectoid-trace} coincides (up to a canonical identification) with the trace map from \cite{LRZ24} and, furthermore, we have $\tr_{f}=\mu^*(\tr_{f^\sharp})$ by the construction of $\tr_{f}$. \qedhere
\end{enumerate}
\end{proof}

Finally, we are almost ready to define the desired trace maps. 

\begin{definition}
\label{trace-maps}
Let $f \colon X \to Y$ be a smooth proper of equidimension $d$ morphism of diamonds over $\Spd(\QQ_p, \ZZ_p)$ and let $n\geq 1$ be an integer.
\begin{enumerate}[leftmargin=*,label=(\arabic*)]
\item The \emph{mod-$n$ trace map} of $f$ is the morphism in $\cD^+_\et(Y;\ZZ/n)$
\[ 
\tr^{\ZZ/n}_f \colon \rR f_{v, *}\underline{\ZZ/n}(d)[2d] \to \underline{\ZZ/n} 
\]
obtained from the maps $\tr_f$ in \cref{prop:trace-diamonds} under the identification $\rR f_!\, \ud{\ZZ/n} \simeq \rR f_{v, *} \ud{\ZZ/n}$.
\item\label{trace-maps-Zp} The \emph{$p$-adic trace map} of $f$ is the morphism in $\cD^+(Y_{v};\ZZ_p)$
\[ 
\tr^{\ZZ_p}_f \colon \rR f_{v, *}\underline{\ZZ}_p(d)[2d] \to \underline{\ZZ}_p
\]
which corresponds to $\rR\lim_r \tr^{\ZZ/p^r}_f$ under the natural identifications $\underline{\ZZ}_p \simeq \rR\lim_r \underline{\ZZ/p^r}$ (see e.g.\ \cite[Prop.~3.1.10]{proetale} and \cref{lemma:replete}) and $\rR\lim_r \rR f_{v, *}\bigl(\underline{\ZZ/p^r}\bigr) \simeq \rR f_{v, *}\bigl(\rR\lim_r \underline{\ZZ/p^r}\bigr)$ (see e.g.\ \cite[\href{https://stacks.math.columbia.edu/tag/0A07}{Tag~0A07}]{stacks-project}). 
\item\label{trace-maps-Qp} The \emph{$\uQ_p$-trace map} of $f$ is the morphism in $\cD(Y_v; \uQ_p)$ 
    \[
        \tr_f^{\Q_p} \colon \rR f_{v, *}\uQ_p(d)[2d] \to \uQ_p
    \]
    obtained as $\tr_f^{\ZZ_p}[\frac{1}{p}]$. Here  we are using the fact that $\rR f_{v, *}$ commutes with uniformly bounded below colimits
as $f$ is qcqs.
\item\label{trace-maps-AI} 
Let $(\mathcal{A}, \mathcal{B}) \in \{(\mathcal{O}^+, \mathcal{O}), (\AA_I, \BB_I)\}$.
The \emph{$\mathcal{A}$-trace map} of $f$ is the morphism in the almost\footnote{Here the almost mathematics
is with respect to $\AA_I^{\circ\circ}$ if $\cA=\AA_I$ and with respect to $\O_Y^{\circ\circ}$ if $\cA = \cO_Y^+$} category $\cD^+(Y_{v};\mathcal{A})^a$
\[
\tr_f^{\mathcal{A}} \colon \rR f_{v, *} \mathcal{A}(d)[2d] \xr{\sim} (\rR f_{v, *} \underline{\ZZ}_p(d)[2d]) \widehat{\otimes}^{\rL}_{\underline{\ZZ}_p}
\mathcal{A} \xrightarrow{\tr_f^{\ZZ_p} \otimes \id} \mathcal{A}.
\]
Here the first arrow is the inverse of the almost equivalence from \cref{Relative primitive comparison-2}.
The \emph{$\mathcal{B}$-trace map}
of $f$ is the morphism in $\cD^+(Y_{v};\mathcal{B})$
\[
\tr_f^{\mathcal{B}} \colon \rR f_{v, *} \mathcal{B}(d)[2d] \to \mathcal{B}.
\]
obtained as $\tr_f^{\cA}[\frac{1}{p}]$. Here we are using the fact that $\rR f_{v, *}$ commutes with uniformly bounded below colimits
as $f$ is qcqs.
\end{enumerate}
\end{definition}

All trace maps from \cref{trace-maps} are compatible with pullbacks in the following sense.
\begin{proposition}
\label{B_I trace compatible under pullback}
    Let $f \colon X \to Y$ be a smooth proper of equidimension $d$ morphism of diamonds over $\Spd(\QQ_p, \ZZ_p)$. 
    Let $g \colon Y' \to Y$ be a map of diamonds over $\Spd(\QQ_p, \ZZ_p)$, giving rise to the Cartesian diagram
    \[ \begin{tikzcd}
        X' \coloneqq Y' \times_Y X \arrow[r, "g'"] \arrow[d, "f'"] & X \arrow[d, "f"] \\
        Y' \arrow[r, "g"] & Y.
    \end{tikzcd}\]
Let $\cA \in \{ \underline{\ZZ/n},\underline{\ZZ}_p, \uQ_p, \mathcal{O}, \BB_I\}$ (resp. $\cA\in \{\O^+, \AA_I\}$).
Then the following diagram
\[ 
\begin{tikzcd}[column sep=large]
g^* \rR f_{v, *} \cA_X(d)[2d] \arrow[r, "\BC"] \arrow[d,"g^*(\tr^\cA_f)"] & \rR f'_{v, *} \cA_{X'}(d)[2d] \arrow[d,"\tr^\cA_{f'}"] \\
g^*\cA_{Y} \arrow[r, "\sim"] & \cA_{Y'}
\end{tikzcd} 
\]
commutes in (the homotopy category of) $\cD^+(Y'_v; \cA)$ (resp.\,in $\cD^+(Y'_v; \cA)^a$). 
\end{proposition}
\begin{proof}
For $\cA=\underline{\ZZ/n}$, this is the compatibility with pullbacks explained in 
\cref{prop:trace-diamonds}~\cref{prop:trace-diamonds-1}.
Thanks to the functoriality of $\rR\lim$, the statement for all $\cA = \underline{\ZZ/p^r}$ yields the statement for $\cA = \underline{\ZZ}_p$ by taking $\rR\lim_r \tr^{\ZZ/p^r}$.

Next, we treat the case $\cA = \mathcal{O}^+ \text{ or } \AA_I$. First, we note that the $v$-site $Y'_v$ is naturally identified with the slice site $(Y_v)_{/Y'}$. Therefore, we have $g^{-1}\cA_Y \simeq \cA_{Y'}$. In particular, we see that $g^*=g^{-1}$. This implies that $g^*$ commutes with limits and, hence, it commutes with completed derived tensor products. Combining it with \cref{Relative primitive comparison-2}, we see that the natural morphism $(g^* \rR f_{v, *} \ud{\ZZ}_p(d)[2d]) \widehat{\otimes}^\rL_{\ud{\ZZ}_p} g^*\cA_Y \to g^* \rR f_{v, *} \cA_X(d)[2d]$ is an almost isomorphism. Now we consider the following diagram:
\[
\begin{tikzcd}[column sep=huge]
(g^* \rR f_{v, *} \ud{\ZZ}_p(d)[2d]) \widehat{\otimes}^\rL_{\ud{\ZZ}_p} g^*\cA_Y \arrow[r, "\BC \widehat{\otimes}^\rL \id"] \arrow{d}{\wr} & (\rR f'_{v, *} \ud{\ZZ}_p(d)[2d]) \widehat{\otimes}^\rL_{\ud{\ZZ}_p} \cA_{Y'} \arrow{d}{\wr} \\
g^* \rR f_{v, *} \cA_X(d)[2d] \arrow{d}{g^*(\tr^\cA_f)} \arrow{r}{\BC} & \rR f'_{v, *} \cA_{X'}(d)[2d] \arrow{d}{\tr^\cA_{f'}}\\
g^* \cA_Y \arrow[r, "\sim"] & \cA_{Y'},
\end{tikzcd} 
\]
where both upper vertical arrows are isomorphisms. Therefore, in order to show that the bottom square commutes, it suffices to show that the top square and the outer square commute. The top square commutes by functoriality of the base change morphisms. As for the outer square, we note that \cref{trace-maps}~\cref{trace-maps-AI} implies that the composition of the left vertical arrows is equal to $g^*(\tr_f^{\ZZ_p}) \wdh{\otimes}^\rL \id$ while the composition of the right vertical arrows is equal to $\tr_{f'}^{\ZZ_p} \wdh{\otimes}^\rL \id$. Therefore, the outer square commutes due to the case $\cA=\ud{\ZZ}_p$ treated above.

Finally, since $g^*$ preserves colimits (being a left adjoint), we get the statement for $\cA = \uQ_p, \mathcal{O} \text{ or }  \BB_I$. 
\end{proof}

\subsection{Cycle class maps}\label{cycle classes}

In this subsection, we discuss cycle class maps for the period sheaves $\BB_I$.
As for trace maps, we must begin by extending the construction in \cite[\S~3]{LRZ24} for $\ZZ/n$-coefficients from the setting of locally noetherian analytic adic spaces to that of sousperfectoid spaces.

Recall that for any closed immersion (in the sense of \cite[Def.~10.7.(ii)]{diamonds}) of locally spatial diamonds $i \colon X \hookrightarrow Y$ the pushforward functor $i_{\et, *} \colon \Shv(X_\et;\ZZ/n) \to \Shv(Y_\et;\ZZ/n)$ is exact and preserves small direct sums (using \cref{enough-points}, this can be checked on qpro\'et-geometric points), and thus preserves all small colimits.
For any $G \in \Shv_\et(Y;\ZZ/n)$, the functor $\Hom(i_{\et, *}(\blank),G)$ is therefore representable \cite[Exp.~XVIII, Lem.~3.1.3]{SGA4}, producing a right adjoint to $i_{\et, *}$.
\begin{definition}\label{upper-shriek}
    Let $i \colon X \to Y$ be a closed immersion of locally spatial diamonds.
    We denote by
    \[ i^! \colon \Shv(Y_\et;\ZZ/n) \longrightarrow \Shv(X_\et;\ZZ/n) \quad \text{and} \quad \rR i^! \colon \cD^+_\et(Y;\ZZ/n) \longrightarrow \cD^+_\et(X;\ZZ/n) \]
    the right adjoint of $i_{\et, *} \colon \Shv(X_\et;\ZZ/n) \to \Shv(Y_\et;\ZZ/n)$ and its derived functor, respectively.
\end{definition}
\begin{remark}
    Concretely, we have $i_{\et, *}\rR i^! = \rR\cHom(i_{\et, *}\underline{\ZZ/n},\blank)$.
    As we will not need this description, we leave the proof to the interested reader.
\end{remark}
\begin{remark}
    It follows directly from \cref{upper-shriek} that for a composition of two closed immersions $X \xhookrightarrow{i} Y \xhookrightarrow{j} Z$, we have a natural isomorphism of functors $\rR(j \circ i)^! \simeq \rR i^! \circ \rR j^!$.
\end{remark}
\begin{lemma}\label{upper-shriek-bc}
    Consider a cartesian diagram of locally spatial diamonds
    \[ \begin{tikzcd}
        U \arrow[r,"u'"] \arrow[d,hook,"i'"] & X \arrow[d,hook,"i"] \\
        V \arrow[r,"u"] & Y,
    \end{tikzcd} \]
    in which $i$ and $i'$ are closed immersions.
    \begin{enumerate}[leftmargin=*,label=\upshape{(\roman*)}]
        \item\label{upper-shriek-bc-map} There is a natural transformation of functors
        \[ u^{\prime,*}\rR i^! \to \rR i^{\prime,!} u^* \colon \cD^+_\et(Y;\ZZ/n) \longrightarrow \cD^+_\et(U;\ZZ/n). \]
        \item\label{upper-shriek-bc-composition} The transformation from \cref{upper-shriek-bc-map} is compatible with compositions in $u$ and $i$.
        \item\label{upper-shriek-bc-etale} If $u$ and $u'$ are \'etale, the base change transformation from \cref{upper-shriek-bc-map} is an isomorphism.
    \end{enumerate}
\end{lemma}
\begin{proof}
    \cref{upper-shriek-bc-map}.
    Since $i_{\et, *}$ and $i'_{\et, *}$ are exact, \cite[Cor.~16.10, case (i)]{diamonds} shows that the natural base change transformation of functors
    \[ u^*i_{\et, *} \to i'_{\et, *}u^{\prime,*} \colon \cD^+_\et(X;\ZZ/n) \longrightarrow \cD^+_\et(V;\ZZ/n) \]
    is an isomorphism.
    Composing its inverse with the counit $\epsilon_i$ for the adjunction between $i_{\et, *}$ and $\rR i^!$, we obtain a natural transformation
    \[ i'_{\et, *}u^{\prime,*}\rR i^! \simeq u^*i_{\et, *}\rR i^! \xrightarrow{u^*(\epsilon_i)} u^* \colon \cD^+_\et(Y;\ZZ/n) \longrightarrow \cD^+_\et(V;\ZZ/n). \]
    Its adjoint under the adjunction of $i'_{\et, *}$ and $\rR i^{\prime,!}$ is the desired natural transformation.

    \cref{upper-shriek-bc-composition}.
    Denote the unit of the adjunction of $i'_{\et, *}$ and $\rR i^{\prime,!}$ by $\eta_{i'}$.
    Then the natural transformation from \cref{upper-shriek-bc-map} is given by
    \[ u^{\prime,*}\rR i^! \xrightarrow{\eta_{i'}} \rR i^{\prime,!}i'_{\et, *}u^{\prime,*}\rR i^! \simeq \rR i^{\prime,!}u^*i_{\et, *}\rR i^! \xrightarrow{\rR i^{\prime,!}u^*(\epsilon_i)} \rR i^{\prime,!}u^*. \]
    The statement therefore follows from the naturality of $\eta_{i'}$ and $\epsilon_i$ and their compatibility with compositions.

    \cref{upper-shriek-bc-etale}.
    Passing to left adjoints, we may instead show that the natural transformation of functors
    \begin{equation}\label{upper-shriek-bc-adjoint}
        u_!i'_{\et, *} \to i_{\et, *}i^*u_!i'_{\et, *} \simeq i_{\et, *}u'_!i^{\prime,*}i'_{\et, *} \to i_{\et, *}u'_! \colon \cD^+_\et(U;\ZZ/n) \longrightarrow \cD^+_\et(Y;\ZZ/n)
    \end{equation}
    is an isomorphism, where the first (resp.\ last) arrow is induced by the unit (resp.\ counit) of the adjunction between $i^*$ and $i_{\et, *}$ (resp.\ $i^{\prime,*}$ and $i'_{\et, *}$).
    By \cref{enough-points}, it suffices to check this on qpro\'et-geometric stalks. Then \cref{etale-pushforward}, \cref{rmk:etale-proetale-pushforward}, and \cite[\href{https://stacks.math.columbia.edu/tag/04H1}{Tag 04H1}]{stacks-project} imply that we can assume that $Y=\Spa(C, C^+)$ for an algebraically closed nonarchimedean field with open bounded valuation subring $C^+\subset C$. But then either $X$ is empty or $i\colon X \to Y$ is an isomorphism. In either of these cases, \cref{upper-shriek-bc-adjoint} is an equivalence. 
\end{proof}
The cycle class maps will be defined for Zariski-closed immersions of sousperfectoid spaces.
\begin{definition}[{\cite[Prop.~IV.4.19, Def.~IV.4.20]{Fargues-Scholze}}]
    Let $S$ be a sousperfectoid space and $f \colon X \to S$ and $g \colon Y \to S$ be smooth morphisms of sousperfectoid spaces.
    Then an $S$-morphism $i \colon X \to Y$ is called a \emph{Zariski-closed immersion} (over $S$) if the following equivalent conditions are satisfied:
    \begin{itemize}
        \item For some open covering $Y = \bigcup_{\alpha} V_\alpha$ with $V_\alpha = \Spa(A_\alpha,A^+_\alpha)$ affinoid, the fiber products $X \times_Y V_\alpha \equalscolon \Spa(B_\alpha,B^+_\alpha)$ are again affinoid, the induced maps $A_\alpha \to B_\alpha$ are surjective, and the $B^+_\alpha$ are the integral closures of the image of $A^+_\alpha$ in $B_\alpha$.
        \item For any affinoid open $V = \Spa(A,A^+) \subset Y$, the fiber product $X \times_Y V \equalscolon \Spa(B,B^+)$ is again affinoid, the induced map $A \to B$ is surjective, and $B^+$ is the integral closure of the image of $A^+$ in $B$.
    \end{itemize}
\end{definition}
\begin{proposition}[{\cite[Prop.~IV.4.19]{Fargues-Scholze}}]\label{local-structure-Zc}
    Let $S$ be a sousperfectoid space and $f \colon X \to S$ and $g \colon Y \to S$ be smooth morphisms of sousperfectoid spaces of equidimension $d_X$ and $d_Y$, respectively.
    Let $i \colon X \hookrightarrow Y$ be a Zariski-closed immersion over $S$.
    Then there exists an open covering $Y = \bigcup_{\alpha} V_\alpha$ and \'etale maps $h_\alpha \colon V_\alpha \to \DD^{d_Y}_S$ such that the preimages $U_\alpha \colonequals X \times_Y V_\alpha$ fit into cartesian diagrams
    \[ \begin{tikzcd}[column sep=huge]
        U_\alpha \arrow[r,hook,"i|_{U_\alpha}"] \arrow[d] & V_\alpha \arrow[d,"h_\alpha"] \\
        \DD^{d_X}_S \arrow[r,hook,"i_0"] & \DD^{d_Y}_S,
    \end{tikzcd}\]
    in which $i_0$ is the natural inclusion of the first $d_X$ coordinates given by $(z_1,\dotsc,z_{d_X}) \mapsto (z_1,\dotsc,z_{d_X},0,\dotsc,0)$.
\end{proposition}
In particular, \cref{local-structure-Zc} implies (together with \cite[Lem.~IV.4.13]{Fargues-Scholze}) that the pullback of a Zariski-closed immersion $i \colon X \hookrightarrow Y$ over $S$ along $S' \to S$ is again a Zariski-closed immersion.
\begin{proof}
    By \cite[Prop.~IV.4.19]{Fargues-Scholze}, we can find an open covering $Y = \bigcup_{\alpha} V_\alpha$ and for each $V_\alpha$ sections $f_1,\dotsc,f_{d_Y} \in \cO(V_\alpha)$ such that $df_1,\dotsc,df_{d_Y}$ is a basis of $\Omega^1_{V_\alpha/S}$ and $U_\alpha = \{ f_{d_X + 1} = \dotsb = f_{d_Y} = 0 \}$.
    After rescaling, we may assume that $f_r \in \cO^+(V_\alpha)$ for all $1 \le r \le d_Y$.
    Then the induced maps $h_\alpha \colonequals (f_1,\dotsc,f_{d_Y}) \colon V_\alpha \to \DD^{d_Y}_S$ are \'etale thanks to \cite[Prop.~IV.4.15]{Fargues-Scholze} and fit into the desired cartesian diagrams.
\end{proof}
For the construction of cycle class maps, we will need the computation of the $\ZZ/n$-\'etale cohomology of the relative affine line.
We remind the reader that we identify the \'etale sites of our sousperfectoid spaces with those of the associated diamonds using \cite[Lem~15.6]{diamonds} and leave the diamondification $(\blank)^\lozenge$ implicit to ease notation.
\begin{lemma}\label{cohomology-A1}
    Let $S$ be a sousperfectoid space with $n\in \O_S^\times$ and $h \colon \AA^{1, \an}_S \to S$ be the relative affine line.
    Then the natural map $\underline{\ZZ/n} \to \rR h_{\et, *}\underline{\ZZ/n}$ in $\cD^+_\et(S;\ZZ/n)$ is an isomorphism.
\end{lemma}
\begin{proof}
    When $S$ is a locally noetherian analytic adic space over $\ZZ_p$, this was shown in part (III) of the proof of \cite[Th.~3.9.1]{Huber-etale}.
    Let $b$ be the same constant as in that proof, which comes from \cite[Th.~2.1]{Luetkebohmert} and only depends on $p$ and $n$.
    Let
    \[ \begin{tikzcd}
        U(m)_S \arrow[rr,hook,"j_m"] \arrow[rd,"h_m"] && \AA^1_S \arrow[ld,"h"'] \\
        & S &
    \end{tikzcd} \]
    be the relative disk of radius $p^{mb}$ over $S$.
    By arguing on each $U(m)$ separately, one can check that the natural map $\underline{\ZZ/n} \to \rR\lim_m \rR j_{m, \et *}\underline{\ZZ/n}$ is an isomorphism.
    Since $\rR h_*$ preserves limits, we conclude that
    \[ \rR h_{\et, *} \underline{\ZZ/n} \simeq \rR h_{\et, *}\bigl(\rR\lim_m \rR j_{m,\et, *}\underline{\ZZ/n}\bigr) \simeq \rR\lim_m \rR h_*\rR j_{m,\et, *}\underline{\ZZ/n} \simeq \rR\lim_m \rR h_{m,\et, *}\underline{\ZZ/n}. \]
    
    We claim that for all $m \in \NN$
    \begin{itemize}
        \item the natural maps $\underline{\ZZ/n} \to \rR^0h_{m,\et, *} \underline{\ZZ/n}$ are isomorphisms and
        \item the natural maps $\rR^{>0}h_{m+1,\et, *}\underline{\ZZ/n} \to \rR^{>0}h_{m,\et, *}\underline{\ZZ/n}$ are 0.
    \end{itemize}
    Indeed, by \cref{enough-points} again, it suffices to check both statements on stalks at qpro\'et-geometric points $\bar{s} \colon \Spa(C,C^+) \to S$.
    Since the $h_m$ are quasi-compact, we can then use base change (\cref{pushforward-stalks}) to reduce to the case when $S = \Spa(C,C^+)$, which is covered by the locally noetherian statement in part (III) of the proof of \cite[Th.~3.9.1]{Huber-etale} (more precisely, statements $(\alpha)'$, $(\beta)'$ and $(\gamma)'$ therein).
    Put together, these two claims show that the natural map
    \[ \{ \underline{\ZZ/n} \}_{m \in \NN} \to \{ \rR h_{m,\et, *} \underline{\ZZ/n} \}_{m \in \NN} \]
    from the constant pro-system is an isomorphism of pro-systems.
    In particular, the natural map $\underline{\ZZ/n} \to \rR\lim_m \rR h_{m,\et, *} \underline{\ZZ/n}$ is an isomorphism in $\cD^+_\et(S;\ZZ/n)$, finishing the proof.
\end{proof}
\begin{proposition}\label{sousperfectoid-purity}
    Let $S$ be a sousperfectoid space with $n \in \cO_S^\times$ and $f \colon X \to S$ and $g \colon Y \to S$ be smooth morphisms of sousperfectoid spaces of equidimension $d_X$ and $d_Y$, respectively.
    Set $c \colonequals d_Y - d_X$.
    Let $i \colon X \hookrightarrow Y$ be a Zariski-closed immersion over $S$.
    Then:
    \begin{enumerate}[leftmargin=*,label=\upshape{(\roman*)}]
        \item\label{sousperfectoid-purity-pullback} (compatibility with pullbacks in $S$) For any morphism of sousperfectoid spaces $u \colon T \to S$ and induced pullbacks $u_X \colon X_T \to X$ and $i_T \colon X_T \hookrightarrow Y_T$, the maps $u^*_X \rR i^! \underline{\ZZ/n} \to \rR i^!_T \underline{\ZZ/n}$ from \cref{upper-shriek-bc}.\cref{upper-shriek-bc-map} are isomorphisms.
        \item\label{sousperfectoid-purity-map} (purity) There is a unique isomorphism $\gamma \colon \underline{\ZZ/n} \xrightarrow{\sim} \rR i^! \underline{\ZZ/n}(c)[2c]$ such that for any map $u \colon T \to S$ from a locally noetherian sousperfectoid space, the induced map in $\cD^+_\et(X_T;\ZZ/n)$
        \[ \underline{\ZZ/n} \xrightarrow{\sim} u^*_X \rR i^!\underline{\ZZ/n}(c)[2c] \simeq \rR i^!_T \underline{\ZZ/n}(c)[2c]  \]
        given by $\gamma_T$ and the isomorphism from \cref{sousperfectoid-purity-pullback} is adjoint to the cycle class map from \cite[Var.~3.3.3]{LRZ24} for the lci immersion $i_T$ of locally noetherian analytic adic spaces.
    \end{enumerate}
\end{proposition}
Recall that by \cref{sousperfectoid-properties}.\cref{sousperfectoid-properties-product}, the base changes $f_T \colon X_T \to T$ and $g_T \colon Y_T \to T$ are again smooth morphisms of sousperfectoid spaces.
\begin{proof}[{Proof of \cref{sousperfectoid-purity}}]
    We prove \cref{sousperfectoid-purity-pullback} and \cref{sousperfectoid-purity-map} simultaneously, proceeding in several steps.
   \begin{enumerate}[wide,label={\textit{Step~\arabic*}.},ref={Step~\arabic*}]
    \item\label{sousperfectoid-purity-uniqueness} \textit{The isomorphism $\gamma$ in \cref{sousperfectoid-purity-map} is unique.} 
    First, note that once the existence of an isomorphism $\underline{\ZZ/n} \simeq \rR i^! \underline{\ZZ/n}(c)[2c]$ in \cref{sousperfectoid-purity-map} is established, we know in particular that $\rR i^! \underline{\ZZ/n}(c)[2c]$ is concentrated in degree $0$.
    By virtue of \cref{enough-points}, the map $\gamma$ is then uniquely determined by its restrictions $\gamma_{\bar{s}}$ to the fibers over qpro\'et-geometric points $\bar{s} \colon \Spa(C,C^+) \to S$ over all points $s \in \abs{S}$. As the $\Spa(C,C^+)$ are in particular locally noetherian sousperfectoid spaces, it therefore suffices to show \cref{sousperfectoid-purity-pullback} and the existence of an isomorphism $\gamma$ in \cref{sousperfectoid-purity-map} whose pullbacks recover the locally noetherian cycle classes from \cite{LRZ24}. 
    
    \item\label{sousperfectoid-purity-etale} \textit{The conclusion is preserved under \'etale maps.}
    More precisely, given a fiber square
    \begin{equation}\label{sousperfectoid-purity-etale-diagram} \begin{tikzcd}
        V \arrow[r,"h'"] \arrow[d,hook,"i'"] & X \arrow[d,hook,"i"] \\
        W \arrow[r,"h"] & Y
    \end{tikzcd} \end{equation}
    in which $i$ and $i'$ are Zariski-closed immersions and $h$ and $h'$ are \'etale, we show that if the conclusion holds for $i$, then it also holds for $i'$.
    Moreover, if $h$ and $h'$ are \'etale covers, the converse is also true.
    
    Set $c \colonequals d_Y-d_X = d_W - d_V$.
    To deal with statement \cref{sousperfectoid-purity-pullback}, we extend \cref{sousperfectoid-purity-etale-diagram} to a diagram of fiber squares
    \[ \begin{tikzcd}
        V_T \arrow[rr,"h'_T"] \arrow[rd,"u_V"] \arrow[dd,hook,"i'_T"] && X_T \arrow[dd,hook,near end,"i_T"] \arrow[rd,"u_X"] & \\[-1em]
        & V \arrow[rr,crossing over,near start,"h'"] && X \arrow[dd,hook,"i"] \\
        W_T \arrow[rr,near start,"h_T"] \arrow[rd,"u_W"'] && Y_T \arrow[rd,"u_Y"] & \\[-1em]
        & W \arrow[rr,"h"] \arrow[from=uu,hook,crossing over,near start,"i'"] && Y.
    \end{tikzcd} \]
    By \cref{upper-shriek-bc}.\cref{upper-shriek-bc-map}, it induces a diagram of natural maps
    \[ \begin{tikzcd}
        h^{\prime,*}_Tu^*_X\rR i^!\underline{\ZZ/n} \arrow[r,phantom,"\simeq"] \arrow[d] &[-2em] u^*_Vh^{\prime,*}\rR i^!\underline{\ZZ/n} \arrow[r,"\sim"] & u^*_V \rR i^{\prime,!}\underline{\ZZ/n} \arrow[d] \\
        h^{\prime,*}_T \rR i^!_T\underline{\ZZ/n} \arrow[rr,"\sim"] && \rR i^{\prime,!}_T\underline{\ZZ/n}.
    \end{tikzcd} \]
    \Cref{upper-shriek-bc}.\cref{upper-shriek-bc-composition} shows that this diagram commutes and \cref{upper-shriek-bc}.\cref{upper-shriek-bc-etale} applied to the front and back shows that the horizontal maps are isomorphisms.
    Thus, the left vertical map is an isomorphism if and only if the right vertical map is an isomorphism.
    Moreover, if $h$ (and hence $h'_T$) is an \'etale cover, then the left vertical map is an isomorphism if and only the natural map $u^*_X\rR i^!\underline{\ZZ/n} \to \rR i^!_T\underline{\ZZ/n}$ is an isomorphism.
    This finishes the assertion for \cref{sousperfectoid-purity-pullback}.

    Next, we discuss \cref{sousperfectoid-purity-map}.
    Given an isomorphism $\gamma \colon \underline{\ZZ/n} \xrightarrow{\sim} \rR i^!\underline{\ZZ/n}(c)[2c]$ which realizes statement \cref{sousperfectoid-purity-map} for $i$, we can use the isomorphism $h^{\prime,*}\rR i^!\underline{\ZZ/n} \simeq \rR i^{\prime,!}\underline{\ZZ/n}$ to define
    \[ \gamma' \colon \underline{\ZZ/n} \simeq h^{\prime,*}\underline{\ZZ/n} \xlongrightarrow[\sim]{h^{\prime,*}\gamma} h^{\prime,*}\rR i^!\underline{\ZZ/n}(c)[2c] \simeq \rR i^{\prime,!}\underline{\ZZ/n}(c)[2c]. \]
    We now check the conditions on the pullbacks along maps $u \colon T \to S$ from locally noetherian sousperfectoid spaces for $\gamma'$.
    Since everything in sight commutes with pullback along $u$, we may assume that $T = S$ is locally noetherian and verify the adjointness to the cycle classes from \cite{LRZ24} for $\gamma'$ (instead of $\gamma'_T$).
    By chasing diagrams, one sees that the adjoint of $\gamma'$ is equal to the pullback of the adjoint of $\gamma$.
    On the other hand, the cycle classes from \cite{LRZ24} satisfy $h^{\prime,*}\cl_i = \cl_{i'}$ by transversal base change (\cite[Lem~3.3.4]{LRZ24}).
    This shows that if the part \cref{sousperfectoid-purity-map} of the conclusion holds for $i$, then it also holds for $i'$.

    Now assume that $h$ and $h'$ are \'etale covers and let $\gamma' \colon \underline{\ZZ/n} \xrightarrow{\sim} \rR i^{\prime,!}\underline{\ZZ/n}(c)[2c]$ be an isomorphism which realizes statement \cref{sousperfectoid-purity-map} for $i'$.
    Let $\pi_1,\pi_2 \colon \widetilde{W} \colonequals W \times_Y W \rightrightarrows W$ be the two \'etale projections and $\widetilde{\imath} \colon \widetilde{V} \hookrightarrow \widetilde{W}$ be the base change of $i'$ along $\pi_j$;
    the two base changes of $i'$ along $\pi_1$ and $\pi_2$ agree because $i'$ is the base change of $i$.
    By our previous discussion, both $\pi_1$ and $\pi_2$ produce isomorphisms $\widetilde{\gamma}_1,\widetilde{\gamma}_2 \colon \underline{\ZZ/n} \rightrightarrows \rR \widetilde{\imath}^!\underline{\ZZ/n}(c)[2c]$ which realize \cref{sousperfectoid-purity-map} for $\widetilde{\imath}$.
    On the other hand, such maps are unique thanks to \cref{sousperfectoid-purity-uniqueness}, so $\widetilde{\gamma}_1 = \widetilde{\gamma}_2$.
    As $\cHom(\underline{\ZZ/n}, \rR i^!\underline{\ZZ/n}(c)[2c])$ is an \'etale sheaf, this implies that we can descend $\gamma'$ to an isomorphism $\gamma \colon \underline{\ZZ/n} \xrightarrow{\sim} \rR i^!\underline{\ZZ/n}(c)[2c]$.
    As before, the pullback of the adjoint of $\gamma$ under $h'$ is the adjoint of $\gamma'$, so the identity $h^{\prime,*}\cl_i = \cl_{i'}$ for the cycle classes from \cite{LRZ24} guarantees that part \cref{sousperfectoid-purity-map} of the conclusion also holds for $i$.
    
    \item\label{sousperfectoid-purity-D1} \textit{Proof when $X = S$, $Y = \DD^1_S$, and $i \colon S \hookrightarrow \DD^1_S$ is the $0$-section.}
    Consider the commutative diagram
    \[ \begin{tikzcd}
        S \arrow[r,equals] \arrow[d,hook,"i"] & S \arrow[d,hook] \\
        \DD^1_S \arrow[r,hook] & \PP^{1,\an}_S, 
    \end{tikzcd} \]
    in which the bottom (resp.\ right) map is the natural open immersion (resp.\ the $0$-section).
    By \cref{sousperfectoid-purity-etale} applied to this fiber square, it suffices to prove instead the statement when $X = S$, $Y = \PP^{1,\an}_S$, and $i \colon S \hookrightarrow \PP^{1,\an}_S$ is the $0$-section.
   
    Denote by $j \colon \AA^{1,\an}_S \hookrightarrow \PP^{1,\an}_S$ the complementary affine line around $\infty$ and set $h \colonequals g \circ j$.
    The (twisted and shifted) pushforward of the standard exact triangle
    \[ i_{\et, *}\rR i^! \underline{\ZZ/n} \longrightarrow \underline{\ZZ/n} \longrightarrow \rR j_* \underline{\ZZ/n} \]
    in $\cD^+_\et(\PP^1_S;\ZZ/n)$ under $g \colon \PP^{1,\an}_S \to S$ together with the identity $\rR g_{\et, *} \circ i_{\et, *} = \rR f_{\et, *} = \id$ yields the exact triangle
    \[ \rR i^! \underline{\ZZ/n}(1)[2] \longrightarrow \rR g_{\et, *} \underline{\ZZ/n}(1)[2] \longrightarrow \rR h_{\et, *} \underline{\ZZ/n}(1)[2] \]
    in $\cD^+_\et(S;\ZZ/n)$.
    Under the natural isomorphisms from \cref{P1-trace}, \cref{cohomology-A1}, and \cite[Cor.\,16.7]{diamonds}, the second map identifies with the projection $\underline{\ZZ/n}(1)[2] \oplus \underline{\ZZ/n} \twoheadrightarrow \underline{\ZZ/n}(1)[2]$ onto the first factor (corresponding to $\rR^0g_{\et, *}\underline{\ZZ/n} \simeq \rR^0h_{\et, *}\underline{\ZZ/n}$). Therefore, we get a canonical isomorphism $\gamma^{-1}\colon \rR i^! \ud{\Z/n}(1)[2] \xr{\sim} \ud{\Z/n}$. We denote its inverse by $\gamma$. 
    As this construction is compatible with pullbacks along arbitrary maps of sousperfectoid spaces $T \to S$, the isomorphism $\gamma$ is also compatible with pullbacks, yielding \cref{sousperfectoid-purity-pullback} in this case.
    Moreover, since $\gamma$ and the cycle classes in \cite{LRZ24} in the divisor case are both defined via the Kummer sequence (\cref{c1-construction} and \cite[Def.~3.1.4, Rmk.~3.1.9]{LRZ24}) and $\cO_{\PP^{1,\an}_S}\bigl(i(S)\bigr) = \cO_{\PP^{1,\an}_S}(1)$, one can check that the pullbacks $\gamma_T$ along maps $T \to S$ from locally noetherian sousperfectoid spaces are adjoint to the cycle class maps from \cite{LRZ24}.\footnote{
    Alternatively, one can use that both constructions are compatible with the algebraic cycle classes by \cref{P1-trace} and \cite[Lem.~3.1.7]{LRZ24}, respectively.} 
    Thus, $\gamma$ satisfies the requirements of \cref{sousperfectoid-purity-map} in this case.

    \item\label{sousperfectoid-purity-disks} \textit{Proof when $X = \DD^e_S$, $Y = \DD^d_S$, and $i \colon \DD^e_S \hookrightarrow \DD^d_S$ is the natural inclusion of the first $e$ coordinates given by $(z_1,\dotsc,z_e) \mapsto (z_1,\dotsc,z_e,0,\dotsc,0)$.}
    We proceed by induction on $d \ge 0$.
    The base case $d = 0$ (and thus $i = \id$) is trivial.
    For the induction step, assume that the statement has been proven for $d-1 \ge 0$.

    We may assume that $d > e$ and factor $i$ as $i \colon \DD^e_S \xhookrightarrow{i_1} \DD^{d-1}_S \xhookrightarrow{i_2} \DD^d_S$.
    First, we claim that it suffices to prove the statement for $i_1$ and $i_2$ separately:
    Indeed, once the statement is established for both $i_1$ and $i_2$, part \cref{sousperfectoid-purity-pullback} for $i$ follows directly from the compatibility of the maps from \cref{upper-shriek-bc}.\cref{upper-shriek-bc-map} with compositions (\cref{upper-shriek-bc}.\cref{upper-shriek-bc-composition}).
    Likewise, the isomorphism $\gamma$ for $i$ in part \cref{sousperfectoid-purity-map} can simply be defined as
    \[ \gamma \colon \underline{\ZZ/n} \xlongrightarrow{\sim} \rR i^!_1\underline{\ZZ/n}(d-e-1)[2d-2e-2] \xlongrightarrow{\sim} \rR i^!_1\rR i^!_2\underline{\ZZ/n}(d-e)[2d-2e] \simeq \rR (i_2 \circ i_1)^!\underline{\ZZ/n}(d-e)[2d-2e]. \]
    For any locally noetherian $T \to S$, the adjoint of $\gamma_T$ is given by the composition of the adjoints, so the compatibility with the cycle class maps from \cite[Var.~3.3.3]{LRZ24} follows from the fact that these cycle class maps are compatible with compositions as well \cite[Cor.~3.3.8]{LRZ24}, establishing the claim.
   
    Next, the statement for $i_1$ follows from the induction hypothesis.
    As for $i_2$, since we show the statement for arbitrary sousperfectoid $S$ and the relative unit disk over a (locally noetherian) sousperfectoid space is again (locally noetherian) sousperfectoid, we can replace $S$ by $\DD^{d-1}_S$ and identify $\DD^{d-1}_S \hookrightarrow \DD^d_S$ with the $0$-section $\DD^{d-1}_S \hookrightarrow \DD^1_{\DD^{d-1}_S}$.
    The assertion then follows again from the case of the $1$-dimensional unit disk (relative over $\DD^{d-1}_S$) proved in \cref{sousperfectoid-purity-D1}.
     
    \item \textit{General case.}
    By the second part of \cref{sousperfectoid-purity-etale}, the statement is \'{e}tale-local on $Y$.
    Thanks to the local structure of Zariski-closed immersions between smooth sousperfectoid spaces from \cref{local-structure-Zc}, it is therefore a direct consequence of \cref{sousperfectoid-purity-etale} and \cref{sousperfectoid-purity-disks}. \qedhere
    \end{enumerate}
\end{proof}
With \cref{sousperfectoid-purity} in place, it is now natural to define cycle classes for sousperfectoid spaces as follows:
\begin{definition}
\label{cycle-class-definition}
    Let $S$ be a sousperfectoid space over $\Spa(\QQ_p, \ZZ_p)$ and $i \colon X \hookrightarrow Y$ be a Zariski-closed immersion of sousperfectoid spaces which are smooth over $S$ of equidimension $d_X$ and $d_Y$, respectively.
    Set $c \colonequals d_Y - d_X$. Note  that in this case \Cref{rmk:etale-proetale-pushforward} and \cref{qproet-v-pushforward-etale} imply that $i_{\et, *} \ud{\ZZ/n} \simeq \rR i_{\qproet, *} \ud{\ZZ/n} \simeq \rR i_{v, *} \ud{\ZZ/n}$ for any integer $n\geq 1$.\footnote{Unlike $i_{\et, *}$, we use the notation $\rR i_{\qproet, *}$ and $\rR i_{v, *}$ because we do not know whether the functors $i_{\qproet, *}$ and $i_{v, *}$ are exact.}
    \begin{enumerate}[leftmargin=*]
        \item The \emph{mod-$n$ cycle class map} of $i$ is the morphism
        \[ \cl^{\ZZ/n}_i \colon \rR i_{v, *}\underline{\ZZ/n} \simeq i_{\et, *} \ud{\ZZ/n} \longrightarrow \underline{\ZZ/n}(c)[2c] \]
        in $\cD^+_\et(Y;\ZZ/n) \subset \cD^+(Y_v; \ZZ/n)$ which is adjoint to the isomorphism $\gamma \colon \underline{\ZZ/n} \xrightarrow{\sim} \rR i^!\underline{\ZZ/n}(c)[2c]$ from \cref{sousperfectoid-purity}.\cref{sousperfectoid-purity-map}.
        \item The \emph{$p$-adic cycle class map} of $i$ is the morphism
        \[ \cl^{\ZZ_p}_i \colon \rR i_{v, *}\underline{\ZZ}_p \longrightarrow \underline{\ZZ}_p(c)[2c] \]
        in $\cD^+(Y_v; \ud{\ZZ}_p)$  under the natural identifications $\underline{\ZZ}_p \simeq \lim_r \underline{\ZZ/p^r} \simeq \rR\lim_r \underline{\ZZ/p^r}$ (see e.g.\ \cite[Prop.~3.1.10]{proetale} and \cite[Lem.~1.2]{Mann-Werner}) and $\rR \lim_r \rR i_{v, *} \underline{\ZZ/p^r} \simeq \rR i_{v, *} \rR \lim_r \underline{\ZZ/p^r}$. 
\item Let $(\mathcal{A}, \mathcal{B}) \in \{(\mathcal{O}^+, \mathcal{O}), (\AA_I, \BB_I)\}$.
The \emph{$\mathcal{A}$-cycle class map} of $i$ is the morphism 
\[ 
\cl^{\cA}_i \colon \rR i_{v, *}\cA \simeq \bigl((\rR i_{v, *} \underline{\ZZ}_p) \widehat{\otimes}_{\underline{\ZZ}_p} 
\cA\bigr) \xlongrightarrow{\cl^{\ZZ_p}_i \otimes \id} \cA(c)[2c] 
\]
in the almost\footnote{Here the almost mathematics
is with respect to $\AA_I^{\circ\circ}$ if $\cA=\AA_I$ and with respect to $\O_Y^{\circ\circ}$ if $\cA = \cO_Y^+$.} category
$\cD^+(Y_{v};\mathcal{A})^a$ in which the first almost isomorphism uses the relative primitive comparison
\cref{Relative primitive comparison-2}.
The \emph{$\mathcal{B}$-cycle class map} of $i$ is the morphism
\[ 
\cl^{\mathcal{B}}_i \colon \rR i_{v, *}\mathcal{B} \to \mathcal{B}(c)[2c] 
\]
in $\cD^+(Y_{v};\mathcal{B})$ obtained
as $\cl^{\cA}_i[1/p]$.
Here we are using the fact that $\rR f_{v, *}$ commutes with uniformly bounded below colimits
as $f$ is quasi-compact and quasi-separated.
\end{enumerate}
\end{definition}

We caution the reader that the cycle class maps from \cref{cycle-class-definition} might a priori depend on the structure morphisms to $S$ (though their fibers over each qpro\'et-geometric point $\bar{s}$ of $S$ are uniquely determined by \cref{sousperfectoid-purity}.\cref{sousperfectoid-purity-map}).
However, this will not cause any problems in this paper. 

To finish this section, we record the following compatibility between the trace maps from \cref{traces} and the cycle class maps from \cref{cycle classes}. For this, we will need the following general lemma. 

\begin{lemma}\label{lemma:section-is-closed} Let $Y$ be a sousperfectoid space, let $f\colon X \to Y$ be a separated qcqs smooth morphism, and let $s\colon Y \to X$ be a section of $f$. Then $s$ is a Zariski-closed immersion. 
\end{lemma}

We refer to \cref{defn:proper-sousperfectoid} for the definition of a separated morphism of sousperfectoid spaces. 
\begin{proof}
    {\it Step~$1$. $s$ is a locally closed immersion, i.e., $s= j\circ i$ such that $i$ is a Zariski-closed immersion and $j$ is an open immersion.} First, \cite[Prop.\,IV.4.19]{Fargues-Scholze} implies that the claim is local on $Y$, so we can assume that $Y=\Spa(R, R^+)$ is an affinoid. Since $s(Y)$ is quasi-compact, we can replace $X$ with a quasi-compact open neighborhood of $s(Y)$ to assume that $X$ is qcqs. In particular, we can assume that $s$ is a spectral map of spectral spaces.
    
    Choose a point $y\in Y$. Then \cite[Lem.\,IV.4.14]{Fargues-Scholze} implies that there are open subspaces $y\in U_y\subset Y$ and $s(y)\in V_y\subset X$ such that $V_y\simeq \DD^d_{U_y}$ as an $U_y$-space. Without loss of generality we can assume that $U_y$ is a quasi-compact open subset (and, therefore, $V_y\simeq \DD^d_{U_y}$ is also quasi-compact). Since $s$ is a spectral map of spectral spaces, \cite[Lem.\,2.3]{diamonds} ensures that $s(U_y) \subset X$ is pro-constructible. Therefore, \cite[\href{https://stacks.math.columbia.edu/tag/0A2X}{Tag 0A2X}]{stacks-project} implies that we can shrink $U_y$ to assume that $y\in U_y$ is affinoid and $s(U_y)\subset V_y$.

    Since the claim is local on $Y$, it suffices to show that $s|_{U_y} \colon U_y \to f^{-1}(U_y)$ is a locally closed immersion. For this, it is enough to show that the map $s|_{U_y} \colon U_y \to V_y$ is a Zariski-closed immersion. Note that $V_y\simeq \DD^d_{U_y}$ is affinoid since $U_y$ is affinoid. Therefore, \cite[Prop.\,IV.4.19]{Fargues-Scholze} guarantees that it suffices to show that $s^\sharp \colon \O(V_y) \to \O(U_y)$ is surjective and $\O^+(U_y)$ is the integral closure of $s^\sharp\bigl(\O^+(V_y)\bigr)$. Both results easily follow from the observation that the composition
    \[
    \bigl( \O(U_y), \O^+(U_y) \bigr) \xr{f^\sharp} \bigl( \O(V_y), \O^+(V_y) \bigr) \xr{s^\sharp} \bigl( \O(U_y), \O^+(U_y) \bigr)
    \]
    is equal to the identity.

    {\it Step~$2$. $s$ is topologically a closed immersion.} Step~$1$ already implies that $s$ is a topological immersion (i.e., $s$ defines a homeomorphism of $Y$ onto its image), so it suffices to show that the image of $s$ is closed. For this, \cite[Lem.\,15.6]{diamonds} implies that it suffices to show that $s^\diam \colon Y^\diam \to X^\diam$ has a closed image. Using the cartesian diagram
    \[
    \begin{tikzcd}[column sep = huge]
        Y^\diam \arrow{d}{s^\diam} \arrow{r}{s^\diam} & X^\diam \arrow{d}{(\id, s^\diam \circ f^\diam)} \\
        X^\diam \arrow{r}{\Delta_{X^\diam/Y^\diam}} & X^\diam \times_{Y^\diam} X^\diam,
    \end{tikzcd}
    \]
    we see that it suffices to show that $f^\diam \colon X^\diam \to Y^\diam$ is separated in the sense of \cite[Def.\,10.7]{diamonds}. This follows from \cref{lemma:different-proper}. 

    {\it Step~$3$. $s$ is a Zariski-closed immersion.} Step~$1$ implies that we can decompose $s$ as a composition $Y \xhookrightarrow{i} V \xhookrightarrow{j} X$ such that $i$ is a Zariski-closed immersion and $j$ is an open immersion. Step~$2$ implies that $s(Y)$ is a closed subset of $X$. Therefore, $X=V\cup (X\smallsetminus s(Y))$ gives a covering of $X$ by two open subsets such that $Y\cap V =Y \to V$ and $Y\cap (X\smallsetminus s(Y)) = \varnothing \to (X\smallsetminus s(Y))$ are Zariski-closed immersions. Thus, \cite[Prop.\,IV.4.19]{Fargues-Scholze} implies that $s$ is a Zariski-closed immersion as well. 
\end{proof}

Finally, we are ready to prove the desired compatibility between the cycle and trace morphisms. 

\begin{proposition}\label{prop:cycle-class-of-section}
Let $S$ be a sousperfectoid space over $\Spa(\QQ_p,\ZZ_p)$, let $f_X\colon X \to S$ and $f_Y \colon Y \to S$ be smooth morphisms of sousperfectoid spaces of equidimension $d_X$ and $d_Y$, respectively. Let $f \colon X \to Y$ be a smooth proper $S$-morphism of equidimension $d\coloneqq d_X - d_Y$, let $s\colon Y \to X$ be a section of $f$, and let $\cA \in \{ \underline{\ZZ/n},\underline{\ZZ}_p,\BB_I, \O\}$.  Then the following composition
\[
\cA \simeq \rR f_{v, *} \rR s_{v, *} \cA\xr{\rR f_{v, *}(\cl^\cA_s)} \rR f_{v, *}\cA(d)[2d] \xr{\tr_f^\cA} \cA
\]
is equal to identity. 
\end{proposition}

We note that  $\cl_s^{\cA}$ is well-defined due to \cref{lemma:section-is-closed}. 

\begin{proof}
    \begin{enumerate}[wide,label={\textit{Step~\arabic*}.},ref={Step~\arabic*}]
    \item\label{prop:cycle-class-of-section-1} \textit{The case of $\cA=\ud{\ZZ/n}$.} \cref{rmk:etale-proetale-pushforward} and \cref{qproet-v-pushforward-etale} imply that $\rR f_{v, *}\ud{\ZZ/n\ZZ}(d)[2d]$ lies in $\cD^+_\et(Y^\diam; \ZZ/n\ZZ)$. Since the composition $\tr_f^{\cA}\circ \rR f_{v, *}(\cl_s^{\ZZ/n})$ is a map between abelian sheaves, \cref{enough-points} implies that we can check that it is equal to the identity on stalks. Since $\tr$ and $\cl$ commute with arbitrary base change (see \cref{B_I trace compatible under pullback} and \cref{sousperfectoid-purity}), we can reduce to the case when $S=\Spa(C, C^+)$ for an algebraically closed nonarchimedean field $C$ and an open bounded valuation subring $C^+\subset C$. In this case, all spaces involved are locally noetherian and, furthermore, the trace and cycle class maps identify with those from \cite{LRZ24} by virtue of \cref{prop:trace-diamonds} and \cref{sousperfectoid-purity}.\cref{sousperfectoid-purity-map}, respectively. Thus, the assertion follows from \cite[Th.~7.2.19]{LRZ24}.
    
    \item\label{prop:cycle-class-of-section-2} \textit{The case of $\cA\in \{\ud{\ZZ}_p, \O^+, \AA_I\}$.} By the naturality of $\rR \lim$, the case $\cA = \underline{\ZZ}_p$ follows from \cref{prop:cycle-class-of-section-1}. Now we consider the case of $\cA \in \{\O^+, \AA_I\}$. From the construction of cycle class maps and trace maps with $\cA$-coefficients in \cref{cycle-class-definition} and \cref{trace-maps}, we have a commutative diagram
    \[ \begin{tikzcd}[column sep=8em,row sep=small]
        \bigl(\rR f_{v, *}\rR s_{v, *}(\ud{\ZZ}_p) \widehat{\otimes}^{\rL}_{\ud{\ZZ}_p} \cA\bigr) \arrow[r,"{\rR f_{v, *}(\cl^{\ZZ_p}_s) \wdh{\otimes}^\rL \id}"] \arrow[d,sloped,"\sim"] & \bigl(\rR f_{v, *}(\ud{\ZZ}_p(d)[2d]) \widehat{\otimes}^{\rL}_{\underline{\ZZ}_p} \cA\bigr) \arrow[r,"\tr^{\ZZ_p}_f \wdh{\otimes}^\rL \id"] \arrow[d,sloped,"\sim"] & \underline{\ZZ}_p \wdh{\otimes}^{\rL}_{\underline{\ZZ}_p} \cA \arrow[d,sloped,"\sim"] \\
        \rR f_{v, *}\rR s_{v, *} \cA \arrow[r,"{\rR f_{v, *}(\cl^{\cA}_s)}"] & \rR f_{v, *} \cA(d)[2d] \arrow[r,"\tr^{\cA}_f"] & \cA,
    \end{tikzcd} \]
    where the vertical (almost) isomorphisms follow from the relative primitive comparison theorem (see \cref{Relative primitive comparison-2}). Therefore, the claim for $\cA\in \{\AA_I, \cO^+\}$ follows from the case $\cA = \underline{\ZZ}_p$.

    \item\label{prop:cycle-class-of-section-3} \textit{The case of $\cA\in \{\O, \BB_I\}$.} This follows from \cref{prop:cycle-class-of-section-2} by inverting $p$.  \qedhere
\end{enumerate} 
\end{proof}

\section{Duality for period sheaves}\label{dualizability}

The purpose of this section is to show that the derived pushforward along smooth proper maps preserves perfect $\O$ and $\BB_I$-complexes. We argue by adapting the strategy from \cite[Sections 6.3-6.4]{LRZ24}. In \cref{subsection:singular-spaces}, we extend these results to proper singular rigid-analytic spaces over a nonarchimedean field $K$ using the proper trace map constructed in \cite[Section 7]{LRZ24}.

Throughout this section, we fix a closed interval $I\subset (0, \infty)$ with rational endpoints. 

\subsection{Digression: \texorpdfstring{$G$}{G}-equivariant sheaves}

The goal of this subsection is to define the notion of $G$-equivariant sheaves for a topos $(X, \O)$
and a sheaf of groups $G\in \Shv(X; \Grp)$. We also verify certain properties of this definition; these properties will be crucially used in our proof of the K\"unneth formula
\Cref{Kunneth formula: main thm}.

For the next definition, we fix a sheaf of groups $G\in \Shv(X; \Grp)$ with the multiplication morphism $m\colon G \times G \to G$ and an object $\widetilde{X}\in \Ob(X)$ with an action $a\colon G\times \wX \to \wX$ in $X$. We denote by $q\colon G\times \widetilde{X} \to \wX$ the projection onto $\wX$. 

For a morphism $f\colon Y' \to Y$ in the topos $X$, we also denote by $f\colon X_{/Y'} \to X_{/Y}$ the induced morphism of sliced topoi. 

\begin{defn}\label{defn:quotient-topos} The {\it quotient topos} $[\widetilde{X}/G]$ is defined as the
following $2$-colimit in the $2$-category of topoi:
\[
[\widetilde{X}/G] \coloneqq 
\begin{tikzcd}[column sep=large]
2\mhyphen\operatorname*{colim}\Bigl(
X_{/G \times G \times \wX} 
  \arrow[r, shift left=2.8ex, "\id \times a"]
  \arrow[r, "m\times \id"]
  \arrow[r, shift right=2.8ex, swap, "\id \times q"]
&
X_{/G\times \wX }
  \arrow[r, shift left=1.4ex, "a"]
  \arrow[r, shift right=1.4ex, swap, "q"]
&
X_{/\wX}
\Bigr).
\end{tikzcd}
\]
When $\widetilde{X}=*$ is the final object of $X$, we denote $[*/G]$ by $[X/G]$.
\end{defn}

The desired $2$-colimit in \cref{defn:quotient-topos} exists due to \cite[Th.\,1]{moerdijk}. Furthermore, \emph{loc.\,cit.\,} implies that the underlying category of $[\widetilde{X}/G]$ is equivalent to $\begin{tikzcd}
2\mhyphen\lim\bigl(X_{/\wX}
\arrow[r, shift left=0.8ex, "a^{-1}"]
  \arrow[r, shift right=0.8ex, swap, "q^{-1}"] &
 X_{/G\times \wX }
   \arrow[r, shift left=1.8ex, ""]
  \arrow[r, ""]
  \arrow[r, shift right=1.8ex, swap,]
&
X_{/G \times G \times \wX}\bigr) 
\end{tikzcd}
$. This implies the following observation: 

\begin{rmk}\label{rmk:quotient-of-torsor} If $\widetilde{X} \to *$ is a $G$-torsor, then the induced morphism $[\widetilde{X}/G] \to X$ is an equivalence. Indeed, combining descent along the effective epimorphism $\widetilde{X} \to *$ and the above description of the underlying category of $[\widetilde{X}/G]$, we conclude that
\[
X \simeq \begin{tikzcd}
2\mhyphen\lim\bigl(X_{/\wX}
\arrow[r, shift left=0.8ex, "a^{-1}"]
  \arrow[r, shift right=0.8ex, swap, "q^{-1}"] &
 X_{/G\times \wX }
   \arrow[r, shift left=1.8ex, ""]
  \arrow[r, ""]
  \arrow[r, shift right=1.8ex, swap,]
&
X_{/G \times G \times \wX}\bigr) 
\end{tikzcd} \simeq [\widetilde{X}/G]. 
\]
\end{rmk}

\begin{notation}\label{notation:group-ring} The \emph{group ring functor} $\ud{\ZZ}[\blank] \colon\Shv(X; \Grp) \to \Shv(X; \Ring)$ is the left adjoint to the functor $(\blank)^\times \colon \Shv(X; \Ring) \to \Shv(X;\Grp)$.

For a sheaf of commutative rings $R\in \Shv(X; \rm{CAlg})$ and a sheaf of groups $G\in \Shv(X; \Grp)$, we set $R[G]\coloneqq R\otimes_{\ud{\ZZ}}\ud{\ZZ}[G]$.
\end{notation}

The next lemma gives a clean way to think about sheaves of abelian groups on $[X/G]$ in terms of sheaves of abelian groups on $X$. Before we formulate this lemma, we note that there is always a canonical morphism of topoi $\pi \colon [X/G] \to X$ which corresponds to the pair of the identity morphism $\id_X\colon X \to X$ and the identity equivalence $\id \colon \id_X \circ a \simeq \id_X \circ q \colon X_{/G} \to X$ under the universal property of $[X/G]$. Similarly, for any sheaf of commutative rings $R\in \Shv(X; \CAlg)$, the pullback $\pi^{-1}R\in \Shv([X/G]; \CAlg)$ canonically defines a sheaf of commutative rings on $[X/G]$, which we also denote by $R$.

\begin{lemma}
\label{lemma:description-quotient} 
Let $X$ be a topos, let $G\in \Shv(X; \Grp)$ be a sheaf of groups on $X$, let $R\in \Shv(X; \rm{CAlg})$ be a sheaf of commutative rings, and let $e\colon X \to [X/G]$ and $\pi\colon [X/G] \to X$ be the natural morphisms of topoi. Then there is an equivalence of categories $\Shv\bigl([X/G]; R\bigr) \simeq \Shv(X; R[G])$ such that under this equivalence  
\begin{enumerate}
    \item $e^{-1}$ is identified with the forgetful functor $\Shv(X; R[G]) \to \Shv(X; R)$;
    \item $e_*$ is identified with $\underline{\Hom}_{R}(R[G], \blank) \colon \Shv(X; R) \to \Shv(X; R[G])$;
    \item $\pi^{-1}$ is identified with the functor $\Shv(X; R) \to \Shv(X; R[G])$ induced by the augmentation morphism $R[G] \to R$;
    \item $\pi_*$ is identified with the functor $\underline{\Hom}_{R[G]}(R, \blank)\colon \Shv(X; R[G]) \to \Shv(X; R)$.
\end{enumerate}
\end{lemma}
\begin{proof}
    Denote by $a\colon X_{/G} \to X$ the natural morphism and by $p_{1}, p_2, m\colon X_{/G\times G} \to X_{/G}$ the morphisms induced by the projections $p_1, p_2\colon G\times G \to G$ and the multiplication map $m\colon G \times G \to G$ respectively.
    
    For any topos $\mathcal{T}$ and $A\in\Shv(\mathcal{T}; \rm{CAlg})$, the category $\Shv(\mathcal{T}; A)$ is canonically identified with $A$-module objects in $\mathcal{T}$. Therefore, \cite[Th.\,1]{moerdijk} implies that 
    \[
    \Shv\bigl([X/G]; R\bigr)\simeq \begin{tikzcd}
2\mhyphen\lim\bigl(\Shv(X; R)
\arrow[r, shift left=0.8ex, "a^{-1}"]
  \arrow[r, shift right=0.8ex, swap, "a^{-1}"] &
 \Shv(X_{/G};R)
   \arrow[r, shift left=1.8ex, ""]
  \arrow[r, ""]
  \arrow[r, shift right=1.8ex, swap,]
&
\Shv(X_{/G\times G};R)\bigr) 
\end{tikzcd}.
\]
    In other words, we have:
    \begin{equation}\label{eqn:descent-2}
    \Shv\bigl([X/G]; R\bigr)\simeq  \bigl\{\F\in \Shv(X; R), \, \theta \colon a^{-1}\F \xr{\sim}a^{-1}\F \,\mid\, m^{-1}(\theta)=p_2^{-1}(\theta) \circ p_1^{-1}(\theta)\bigr\}
    \end{equation}
    Note that $\theta$ defines a morphism of sheaves of sets $G \to \ud{\Aut}_{R}(\F)$ and the cocycle condition is equivalent to the fact that $\theta$ is a group homomorphism. Therefore, we conclude that 
    \[
    \Shv\bigl([X/G]; R\bigr)\simeq \{\F\in \Shv(X; R), \varphi \in \Map_{\grp}\bigl(G,\ud{\Aut}_{R}(\F)\bigr)\}.
    \]
    On the other hand, we have $\Shv(X; R[G]) \simeq \{\F\in \Shv(X; R), \varphi\in \Map_{\ring}\bigl(R[G], \ud{\End}_{R}(\F)\bigr)\}$. Combining it with the $\otimes$-$\rm{forget}$-adjunction and the adjunction from \cref{notation:group-ring}, we conclude that 
    \[
    \Shv(X; R[G])\simeq \{\F\in \Shv(X; R), \varphi \in \Map_{\grp}\bigl(G,\ud{\Aut}_{\ud{R}}(\F)\bigr)\} \simeq \Shv\bigl([X/G]; R\bigr). 
    \]
    By construction, one sees that $e^{-1}$ is identified with the forgetful functor $\Shv(X; R[G]) \to \Shv(X; R)$. By passing to right adjoints, we conclude that $e_*$ is identified with $\underline{\Hom}_{R}(R[G], \blank) \colon \Shv(X; R) \to \Shv(X; R[G])$. 

    Similarly, we see that $\pi^{-1}$ is identified with the functor that sends a sheaf $\F\in \Shv(X; R)$ to the pair $(\F, \id\colon a^{-1}\F \xr{\sim}a^{-1}\F)$ under the equivalence from \cref{eqn:descent-2}. After unraveling all the identifications, one sees that this functor corresponds to the functor $\Shv(X; R) \to \Shv(X; R[G])$ induced by the augmentation morphism $R[G] \to R$. By passing to right adjoints, we see that $\pi_*$ is identified with the functor $\underline{\Hom}_{R[G]}(R, \blank)\colon \Shv(X; R[G]) \to \Shv(X; R)$. 
\end{proof}

Now we are ready to formulate the main result of this subsection.  

\begin{proposition}
\label{G-torsor and cohomology}
Let $f\colon X \to Y$ be a morphism of topoi, let $G\in \Shv(Y, \Grp)$ be a sheaf of groups on $Y$, let $R\in \Shv(Y, \rm{CAlg})$ be a sheaf of commutative rings on $Y$, let $\widetilde{X} \to *$ be a $f^{-1}G$-torsor on $X$ with the induced morphism $j\colon X_{/\widetilde{X}} \to X$, and let $\widetilde{f}\colon X_{/\widetilde{X}} \to Y$ be the composite morphism. Then there is a natural morphism of topoi $g\colon X \to [Y/G]$ such that the diagram
\begin{equation}\label{eqn:action-on-torsor}
\begin{tikzcd}[column sep=7em]
& D(X; f^{-1}R) \arrow[ld, "\rR f_*" swap] \arrow[d, "\rR g_*"] \arrow{r}{j^{-1}}& D\bigl(X_{/\widetilde{X}}; f^{-1}R\bigr) \arrow[d, "\rR\widetilde{f}_*"] \\
D(Y; R) & D([Y/G]; R)\simeq D(Y; R[G]) \arrow[l, "\rR\ud{\Hom}_{R[G]}(R\mathpunct{\raisebox{0ex}{,}} \blank)"] \arrow[r, "\mathrm{forget}", swap] & D(Y; \ZZ)
\end{tikzcd}
\end{equation}
commutes (up to a canonical equivalence). 
\end{proposition}
\begin{proof}
        We denote by $a\colon Y_{/G} \to Y$ the slice morphism, we denote by $\alpha, q\colon f^{-1}G\times \wX \to \wX$ the action and the projection morphisms, and (by slight abuse of notation) we also denote by $\alpha, q \colon X_{/f^{-1}G\times \wX} \to X_{/\wX}$ the associated morphisms of slice topoi. Finally, we denote by $\widetilde{f}_{/G} \colon X_{/f^{-1}G\times \wX} \to Y_{/G}$ and $\widetilde{f}_{/G\times G} \colon X_{/f^{-1}G\times f^{-1}G \times \wX} \to Y_{/G\times G}$ induced by $\widetilde{f}$ (see \cite[\href{https://stacks.math.columbia.edu/tag/04H1}{Tag 04H1}]{stacks-project}). Then the diagram of topoi 
    \begin{equation}\label{eqn:alotoftopoi}
    \begin{tikzcd}
    X_{/f^{-1}G \times f^{-1}G \times \wX} \arrow[d, shift left=1.8ex, ""]
  \arrow[d, ""]
  \arrow[d, shift right=1.8ex, swap]  \arrow{r}{\widetilde{f}_{/G\times G}} & Y_{/G\times G} \arrow[d, shift left=1.8ex, ""]
  \arrow[d, ""]
  \arrow[d, shift right=1.8ex, swap] \\
    X_{/f^{-1}G \times \widetilde{X}} \arrow[d, shift left=0.8ex, "q"]
  \arrow[d, shift right=0.8ex, swap, "\alpha"] \arrow{r}{\widetilde{f}_{/G}}& Y_{/G} \arrow[d, shift left=0.8ex, "a"]
  \arrow[d, shift right=0.8ex, swap, "a"] \\
    X_{/\widetilde{X}} \arrow{r}{\widetilde{f}} & Y
    \end{tikzcd}
    \end{equation}
    commutes (up to a canonical equivalence). We wish to show that two possible bottom squares in Diagram~\cref{eqn:alotoftopoi} are cartesian. First, the diagram 
    \begin{equation}\label{eqn:cartesian-1}
    \begin{tikzcd} 
    X_{/f^{-1}G \times \widetilde{X}} \arrow[d, "q"]  \arrow{r}{\widetilde{f}_{/G}}& Y_{/G} \arrow[d, "a"] \\
    X_{/\widetilde{X}} \arrow{r}{\widetilde{f}} & Y
    \end{tikzcd}
    \end{equation}
    is cartesian by virtue of \cite[Exp.~IV, Prop.~5.11]{SGA4}. Furthermore, let $\widetilde{\alpha} \colon f^{-1}G \times \wX \to f^{-1}G\times \wX$ be the morphism defined via the formula $\widetilde{\alpha}(g, x) = (g, \alpha(g, x))$ and let $\widetilde{\alpha} \colon X_{/f^{-1}G\times \wX} \to X_{/f^{-1}G\times \wX}$ be the associated morphism of slice topoi. Then the assumption that $\wX$ is a $f^{-1}G$-torsor implies that $\widetilde{\alpha}$ is an isomorphism and that the diagram
    \[
    \begin{tikzcd}
         X_{/f^{-1}G \times \widetilde{X}} \arrow[r, "\widetilde{\alpha}"] \arrow[d, "\alpha"] & X_{/f^{-1}G \times \widetilde{X}} \arrow[dl, "q"] \\
         X_{/\wX}
    \end{tikzcd}
    \]
    commutes (up to a canonical equivalence). Combining this with the fact that Diagram~\cref{eqn:cartesian-1} is cartesian, we deduce that the diagram 
    \begin{equation}\label{eqn:cartesian-2}
    \begin{tikzcd} 
    X_{/f^{-1}G \times \widetilde{X}} \arrow[d, "\alpha"]  \arrow{r}{\widetilde{f}_{/G}}& Y_{/G} \arrow[d, "a"] \\
    X_{/\widetilde{X}} \arrow{r}{\widetilde{f}} & Y
    \end{tikzcd}
    \end{equation}
    is cartesian as well. 
    
    Now the universal property of $2$-colimits, \cref{rmk:quotient-of-torsor}, Diagram~\cref{eqn:alotoftopoi} induce an essentially unique morphism of topoi $g \colon X\simeq [X_{/\widetilde{X}}/f^{-1}G] \to [Y/G]$ such that the following diagram of topoi commutes
    \[
    \begin{tikzcd}
    X_{/\widetilde{X}} \arrow{r}{\widetilde{f}} \arrow{d}{j} & Y \arrow{rd}{\id} \arrow{d}{e} & \\
    X\arrow[r, "g"] \arrow[rr, bend right, "f"]& {[}Y/G{]} \arrow[r, "\pi"] & Y,
    \end{tikzcd}   
    \]
    where the morphisms $\pi$ and $e$ are as in \cref{lemma:description-quotient}. Combining this with the identifications $\Shv([Y/G]; R)\simeq\Shv(Y; R[G])$ and $\pi_*=\ud{\Hom}_{R[G]}(R, \blank)$ from \cref{lemma:description-quotient}, we immediately see that the left triangle in \cref{eqn:action-on-torsor} commutes (up to a canonical equivalence). To see that the right square commutes (up to a canonical equivalence), it suffices to show that the natural transformation of functors
    \[
    e^{-1}\circ g_* \to \widetilde{f}_* \circ j^{-1} \colon \Shv(X; \ZZ) \to \Shv(Y; \ZZ)
    \]
    is an equivalence. For this, we note that the functor $g^{-1}$ sends a descent pair
    $\bigl\{F\in \Shv(Y; \ZZ), \varphi \colon a^{-1}F \xr{\sim} a^{-1}F\bigr\}$ to the pair 
    \[
    \bigl\{\widetilde{f}^{-1}F\in \Shv(X_{/\widetilde{X}}; \ZZ), \, \alpha^{-1}\widetilde{f}^{-1} F\simeq \widetilde{f}_{/G}^{-1}a^{-1}F \xr{\widetilde{f}_{/G}^{-1}(\varphi)} \widetilde{f}_{/G}^{-1}a^{-1}F \simeq q^{-1}\widetilde{f}^{-1}F\bigr\}.
    \]
    Since $g_*$ is the right adjoint to this functor, we see that 
    \[
    g_*\Bigl(\bigl\{G\in \Shv(X_{/\widetilde{X}}; \ZZ), \psi \colon \alpha^{-1}G \xr{\sim} q^{-1}G \bigr\}\Bigr) = \bigl\{\widetilde{f}_*G\in \Shv(Y; \ZZ), a^{-1}\widetilde{f}_*G \simeq \widetilde{f}_{/G, *}\alpha^{-1} G \xr{\widetilde{f}_{/G, *}(\psi)} \widetilde{f}_{/G, *} q^{-1}G \simeq a^{-1}\widetilde{f}_{*}G\bigr\},
    \]
    where the isomorphisms $a^{-1}\widetilde{f}_*G \simeq \widetilde{f}_{/G, *}\alpha^{-1} G$ and $\widetilde{f}_{/G, *} q^{-1}G \simeq a^{-1}\widetilde{f}_{*}G$ come from the fact that Diagrams~\cref{eqn:cartesian-1} and \cref{eqn:cartesian-2} are cartesian and the trivial base change (see \cite[\href{https://stacks.math.columbia.edu/tag/04H1}{Tag 04H1}]{stacks-project}). 

    With this description, one sees directly that the natural transformation of functors $e^{-1}\circ g_* \to \widetilde{f}_* \circ j^{-1}$ is an equivalence. 
\end{proof}

\subsection{K\"unneth formula}

The main goal of this subsection is to prove a solid version of the K\"unneth formula
for perfect $\O$ and $\BB_I$-complexes. 

\begin{construction}[{K\"{u}nneth map, cf.\ \cite[Construction 6.3.4]{LRZ24}}]
\label{Construction:Kunneth map}
Let \begin{equation}
\label{eqn:kunneth-data}
\begin{tikzcd}
W \arrow[r,"g"] \arrow[d,"g'"] \arrow[rd, "h"] & X \arrow[d,"f"] \\
X' \arrow[r,"f'"] & S
\end{tikzcd}
\end{equation} 
be a commutative diagram of diamonds over $\Spd(\QQ_p, \ZZ_p)$, let $\cA\in \{\BB_I, \O\}$, and let $\cE \in \cD(X_v; \cA)$ and $\cE'\in \cD(X'_v; \cA)$.\footnote{
We caution the reader that throughout this construction, all instances of
$(-)^*$ mean the (derived) pullback in the setting of maps of ringed sites.
Namely, $f^*(\cE) \coloneqq f^{-1}(\cE) \otimes^\rL_{f^{-1}(\cA_Y)} \cA_X$
for a map of ringed sites $f \colon (X, \cA_X) \to (Y, \cA_Y)$.}

\begin{enumerate}[leftmargin=*,label=\upshape{(\arabic*)}]
\item\label{Kunneth-map-definition} 
We define the \emph{K\"{u}nneth map}
\[
\rm{KM}_{f, f'}\colon \rR f_{v, *}\cE \otimes^\rL_{\cA} \rR f'_{v, *}\cE' \to \rR h_{v, *}(g^*\cE \otimes^\rL_{\cA} g'^*\cE')
\]
as the adjoint to the map
\[
h^*\bigl(\rR f_{v, *}(\cE) \otimes^\rL_{\cA} \rR f'_{v, *}(\cE')\bigr) \simeq (g^* f^* \rR f_{v, *}\cE) \otimes^\rL_{\cA} (g'^* f'^* \rR f'_{v, *} \cE') 
\xrightarrow{g^*(\epsilon_f) \otimes^\rL g'^*(\epsilon_{f'})} g^*\cE \otimes^\rL_{\cA} g'^*\cE',
\]
where $\epsilon_f$ (resp.\ $\epsilon_{f'}$) denotes the counit of the $(f^*, \rR f_{v, *})$-adjunction (resp.\ the $(f'^*, \rR f'_{v, *})$-adjunction).

\item Combining the above K\"{u}nneth map with \Cref{variant:condensed-global-sections}, we get the \emph{condensed K\"unneth map} as the following composition:
\[ \begin{tikzcd}
\rR\underline{\Gamma}_v(X, \cE) \otimes^\rL_{\underline{\Gamma}_v(S, \cA)} \rR\underline{\Gamma}_v(X', \cE')
\arrow[rd,swap, "\KM^\cond_{X, X'}"] \arrow[r, "\sim"] &
\rR\underline{\Gamma}_v(S, \rR f_*\cE) \otimes^\rL_{\underline{\Gamma}_v(S, \cA)} 
\rR\underline{\Gamma}_v(S, \rR f'_*\cE') \arrow[r, "\cup"] &
\rR\underline{\Gamma}_v(S,\rR f_{v, *}\cE \otimes^\rL_{\cA} \rR f'_{v, *}\cE') \arrow[d, "\rm{KM}_{f, f'}"] \\
& \rR\underline{\Gamma}_v(W, g^*\cE \otimes^\rL_{\cA} g'^*\cE') &
\rR\underline{\Gamma}_v(S, \rR h_{v, *}(g^*\cE \otimes^\rL_{\cA} g'^*\cE')), \arrow[l, "\sim"']
\end{tikzcd} \]
where $\cup$ is the cup-product map from 
\cite[\href{https://stacks.math.columbia.edu/tag/0B6C}{Tag~0B6C}]{stacks-project}.

\item\label{Kunneth-map-definition-3} If $\rR\underline{\Gamma}_v(X, \cE)$, $\rR\underline{\Gamma}_v(X', \cE')$, and $\rR\underline{\Gamma}_v(W, g^*\cE \otimes^\rL_{\cA} g'^*\cE')$ are solid, then $\KM^\cond$ factors through the solid tensor product
\[
\KM^\Box_{X, X'}\colon \rR\underline{\Gamma}_v(X, \cE) \otimes^{\rL \sol}_{\underline{\Gamma}_v(S, \cA)}
\rR\underline{\Gamma}_v(X', \cE') \to \rR\underline{\Gamma}_v(W, g^*\cE \otimes^\rL_{\cA} g'^*\cE'),
\]
which we will refer to as the \emph{solid K\"{u}nneth map}. If there is no risk of confusion, we denote it simply by $\KM^\sol$. 
\end{enumerate}
We will also often denote $g^*\cE \otimes^\rL_\cA g'^*\cE'$ simply by $\cE \boxtimes^\rL \cE'$ unless there is some risk for confusion. 
\end{construction}

\begin{rmk} The solid assumption in \cref{Construction:Kunneth map}\cref{Kunneth-map-definition-3} is automatic if $\cE\in \Perf_v(X; \cA)$ and $\cE'\in \Perf_v(X'; \cA)$ (see \cref{RGamma-solid}).
\end{rmk}

The main result of this subsection is the following K\"{u}nneth formula statement. 

\begin{theorem}
\label{Kunneth formula: main thm}
Let $S$ be a strictly totally disconnected perfectoid space over $\Spd(\Q_p, \Z_p)$, let $S^\sharp$ be the corresponding untilt.
Let $f^\sharp \colon X^\sharp \to S^\sharp$ and $f'^\sharp \colon X'^\sharp \to S^\sharp$ be 
qcqs smooth maps of sousperfectoid spaces and let $\cA\in \{\BB_I, \O\}$. Set $X\coloneqq (X^\sharp)^\diam$ and $X'\coloneqq (X'^\sharp)^\diam$. Let $\cE \in \Perf_v(X, \cA)$ and 
$\cE' \in \Perf_v(X', \cA)$.
Then the solid K\"{u}nneth map 
\[
\KM^\Box\colon \rR\underline{\Gamma}_v(X, \cE) \otimes^{\rL \Box}_{\underline{\Gamma}_v(S, \cA)}
\rR\underline{\Gamma}_v(X', \cE') \rightarrow
\rR\underline{\Gamma}_v(X \times_{S} X', \cE \boxtimes^{\rL} \cE').
\]
is an isomorphism. 
\end{theorem}

We start the proof of \cref{Kunneth formula: main thm} with the following preliminary lemmas.

\begin{lemma}\label{lemma:cohomology-of-zp} Let $\gamma\in \ZZ_p$ be a topological generator of $\ZZ_p$, let $M\in \cD(\ud{\ZZ}_p[\ud{\ZZ}_p])$, and let $\gamma^*\colon M \to M$ be the associated automorphism of $M$. If the underlying condensed complex of $M$ is solid (i.e., it lies in $\cD_\sol(\ud{\ZZ})$), then the natural morphism
\[
\xymatrix{
\rR\ud{\Hom}_{\ud{\ZZ}_p[\ud{\ZZ}_p]}(\ud{\ZZ}_p, M) \ar[r]^-{\sim} &
\lim \Bigl( M
\ar@<.4ex>[r]^-{\gamma^*} \ar@<-.4ex>[r]_-{\id} &
M \Bigr)
}
\]
is an isomorphism.
\end{lemma}

We remind the reader that $\ud{\ZZ}_p$ denotes the condensed ring associated with the topological ring $\ZZ_p$ and $\ud{\ZZ}_p[\ud{\ZZ}_p]$ denotes the associated $\ud{\ZZ}_p$-linear group algebra (see \cref{notation:group-ring}). 

\begin{proof}
    First, \cite[Lem.\,B.4(i)]{bosco-master} implies that $\ud{\ZZ}_p[\ud{\ZZ}_p]^\Box \simeq \ud{\ZZ}_p[\![\ud{\ZZ}_p]\!] \colonequals \lim_n \ud{\ZZ}_p[\ud{\ZZ/p^n}]$. Therefore, \cite[Lem.\,A.16, B.4(ii)]{bosco-master} and our assumption that $M$ is solid imply that $\rR\ud{\Hom}_{\ud{\ZZ}_p[\ud{\ZZ}_p]}(\ud{\ZZ}_p, M) \simeq \rR\ud{\Hom}_{\ud{\ZZ}_p[\![\ud{\ZZ}_p]\!]}(\ud{\ZZ}_p, M)$. In other words, it suffices to show that the natural map
 \[
\xymatrix{
\rR\ud{\Hom}_{\ud{\ZZ}_p[\![\ud{\ZZ}_p]\!]}(\ud{\ZZ}_p, M) \ar[r]^-{\sim} &
\lim \Bigl( M
\ar@<.4ex>[r]^-{\gamma^*} \ar@<-.4ex>[r]_-{\id} &
M \Bigr)
}
\]
is an isomorphism. This follows directly from the fact that $\ud{\ZZ}_p[\![\ud{\ZZ}_p]\!] \xr{\gamma-1} \ud{\ZZ}_p[\![\ud{\ZZ}_p]\!]$ is a projective resolution of $\ud{\ZZ}_p$ in the category of $\ud{\ZZ}_p[\![\ud{\ZZ}_p]\!]$-modules. 
\end{proof}

\begin{lemma}
\label{cohomology and Koszul: one variable}
Let $\wX \to X$ be a $\underline{\ZZ}_p$-torsor
of diamonds over $\Spd(\QQ_p, \ZZ_p)$, let $\gamma\in \ZZ_p$ be a topological generator, and let $\gamma\colon \wX \to \wX$ be the associated $X$-automorphism of $\wX$. Let $\cE \in \cD(X_v;\ud{\ZZ}_p)$ be a complex of sheaves of $\ud{\ZZ}_p$-modules such that $\rR\underline{\Gamma}_v(\wX, \cE)$ is solid, i.e., it lies $\cD_\Box(\ud{\ZZ})$.
Then the natural morphism
\[
\xymatrix{
\rR\underline{\Gamma}_v(X, \cE) \ar[r]^-{\sim} &
\lim \Bigl( \rR\uGamma_v(\wX, \restr{\cE}{\wX})
\ar@<.4ex>[r]^-{\gamma^*} \ar@<-.4ex>[r]_-{\id} &
\rR\underline{\Gamma}_v(\wX, \restr{\cE}{\wX}) \Bigr)
}
\]
is an isomorphism and $\rR\uGamma_v(X, \cE)$ is solid as well. 
\end{lemma}
\begin{proof}
We apply \cref{G-torsor and cohomology} to the morphism $X_v \to *_\proet$, the group object $\ud{\ZZ}_p$ on $*_\proet$, and the $\ud{\ZZ}_p$-torsor $\wX \to X$ to conclude that $\rR\uGamma_v(\wX, \cE)$ is naturally an object of $\cD(\ud{\ZZ}_p[\ud{\ZZ}_p])$ and that $\rR\uGamma_v(X, \cE) \simeq \rR\underline{\Hom}_{\ud{\Z}_p[\ud{\ZZ}_p]}\bigl(\ud{\ZZ}_p, \rR\uGamma_v(\wX, \restr{\cE}{\wX}\bigr)$. Since $\rR\uGamma_v(\wX, \cE)$ is assumed to be solid, the isomorphism
\[
\xymatrix{
\rR\underline{\Gamma}_v(X, \cE) \ar[r]^-{\sim} &
\lim \Bigl( \rR\uGamma_v(\wX, \restr{\cE}{\wX})
\ar@<.4ex>[r]^-{\gamma^*} \ar@<-.4ex>[r]_-{\id} &
\rR\underline{\Gamma}_v(\wX, \restr{\cE}{\wX}) \Bigr)
}
\]
follows immediately from \cref{lemma:cohomology-of-zp}. Now the fact that $\rR\uGamma_v(X, \cE)$ is solid follows from the fact that $\cD_\sol(\ud{\Z})\subset \cD(\Cond(\Ab))$ is closed under limits (see \cite[Def.\,5.1]{Condensed}).  
\end{proof}

\begin{definition}
A map of diamonds $\widetilde{X} \to X$ is called a \emph{successive $\ud{\ZZ}_p$-torsor}
if there exists $n \in \ZZ_{>0}$ and a factorization $\widetilde{X} = X_n \to X_{n-1}
\to \cdots \to X_0 = X$ such that each $X_i \to X_{i-1}$ is a $\ud{\ZZ}_p$-torsor.
\end{definition}

Here is the main example of a successive $\ud{\ZZ}_p$-torsors.

\begin{example}
\label{successive Zp-torsor}
Let $S=\Spa(R, R^+)$ be a strictly totally disconnected perfectoid space over $\Spa(\QQ_p, \ZZ_p)$. We fix a compatible system $(\zeta_{p^r}\in R)_{r\in \NN}$ of $p$-power roots of $1$. This defines a trivialization $\ZZ_p \xr{\sim} \ZZ_p(1)$. Let $\TT^d_S = \Spa(R\langle T_1^{\pm 1}, \dots, T_d^{\pm 1}\rangle, R^+\langle T_1^{\pm 1}, \dots, T_d^{\pm 1}\rangle)$ be the $d$-dimensional relative torus, and let  
\[
\widetilde{\TT}^d_S \coloneqq \Spa(R\langle T_1^{\pm 1/p^\infty}, \dots, T_d^{\pm 1/p^\infty}\rangle, R^+\langle T_1^{\pm 1/p^\infty}, \dots, T_d^{\pm 1/p^\infty}\rangle)
\]
be the $d$-dimensional perfectoid torus. Then the natural morphism
\[
\bigl(\widetilde{\TT}^d_S\bigr)^\diam \simeq \lim_m \Spd\bigl(R\langle T_1^{\pm 1/p^m}, \dots, T_d^{\pm 1/p^m}\rangle, R^+\langle T_1^{\pm 1/p^m}, \dots, T_d^{\pm 1/p^m}\rangle\bigr) \to (\TT^d_S)^\diam
\]
is a quasi-pro-\'etale covering and an $\ud{\ZZ}_p^d \simeq \ud{\ZZ}_p(1)^d$-torsor. Furthermore, this is a successive $\ZZ_p$-torsor since we may successively adjoin $p$-power roots of one variable at a time.  
\end{example}

\begin{lemma}
\label{Kunneth formula: one variable}
Let $S$ be a strictly totally disconnected perfectoid space over $\Spd(\QQ_p,\ZZ_p)$ and let \begin{equation}
\begin{tikzcd}
\widetilde{W} \arrow[dd, bend right = 50, "\widetilde{g'}"'] \arrow{d}{t'} \arrow{r}{\widetilde{g}} \arrow[dr, phantom, "\lrcorner", very near start] & \wX \arrow{d}{t} \\
W \arrow[r,"g"] \arrow[d,"g'"] & X \arrow[d,"f"] \\
X' \arrow[r,"f'"] & S
\end{tikzcd}
\end{equation} 
be a commutative diagram of diamonds over $\Spd(\QQ_p, \ZZ_p)$ such that the upper square is cartesian and the morphism $t$ is a successive $\ud{\ZZ}_p$-torsor. 
Let $\cA\in \{\BB_I, \O\}$, let $\cE \in \Perf_v(X; \cA)$ and $\cE'\in \Perf_v(X'; \cA)$ such that the solid K\"unneth map 
\[
\KM^\sol_{\wX, X'}\colon \rR\underline{\Gamma}_v(\wX, \restr{\cE}{\wX})
\otimes^{\rL \sol}_{\underline{\Gamma}_v(S, \cA)}
\rR\underline{\Gamma}_v(X', \cE') \rightarrow
\rR\underline{\Gamma}_v\bigl(\widetilde{W},\widetilde{g}^*(\restr{\cE}{\wX})
\otimes^\rL \widetilde{g}^{\prime,*}\cE'\bigr)
\]
is an isomorphism. Then the following solid K\"{u}nneth map 
\[
\KM^\sol_{X, X'}\colon \rR\underline{\Gamma}_v(X, \cE) \otimes^{\rL \sol}_{\underline{\Gamma}_v(S, \cA)}
\rR\uGamma_v(X', \cE') \to
\rR\uGamma_v(W, g^*\cE \otimes^\rL g^{\prime,*}\cE')
\]
is an isomorphism as well. 
\end{lemma}

Note that all $\rR\uGamma_v(\blank, \blank)$ involved in the formulation of \cref{Kunneth formula: one variable} are solid due to \cref{RGamma-solid}. Thus, the solid K\"unneth map is indeed well-defined.  

\begin{proof}
A straightforward inductive argument immediately reduces the question to the case when $\wX \to X$ is a $\ud{\ZZ}_p$-torsor. For an object $M\in \cD_\sol(\Cond(\Ab))$ with an endomorphism $f\colon M \to M$, we set \[
\xymatrix{
(M)^{f = 1} \coloneqq 
\lim \biggl( M
\ar@<.4ex>[r]^-{f} \ar@<-.4ex>[r]_-{\rm{id}} &
M \biggr) \simeq \mathrm{fib}(M \xrightarrow{f - \rm{id}} M).
}
\]
With this notation at hand, the functoriality of the solid K\"{u}nneth map from \Cref{Construction:Kunneth map} in the input diagram yields a commutative diagram
\[
\begin{tikzcd}[column sep=10 em]
\rR\uGamma_v(X, \cE) \otimes^{\rL \sol}_{\uGamma_v(S, \cA)}
\rR\uGamma_v(X', \cE') \arrow[r,"\KM^\sol_{X, X'}"] \arrow[d] & 
\rR\uGamma_v(W, g^*\cE \otimes^\rL g'^*\cE') \arrow[d] \\
\left(\rR\uGamma_v(\wX,\restr{\cE}{\wX})
\otimes^{\rL \sol}_{\uGamma_v(S, \cA)}
\rR\uGamma_v(X', \cE') \right)^{\gamma = 1}\arrow[r,"(\KM^\sol_{\wX, X'})^{\gamma = 1}"] & 
\left(\rR\uGamma_v(\widetilde{W}, \widetilde{g}^*(\restr{\cE}{\wX}) 
\otimes^\rL \widetilde{g}'^*\cE')\right)^{\gamma = 1}.
\end{tikzcd}
\]
Our assumption implies that $(\KM^\sol_{\wX, X'})^{\gamma = 1}$ is an isomorphism. Therefore, it suffices to show that both vertical maps are isomorphisms. The right vertical map is an isomorphism due to \cref{cohomology and Koszul: one variable}. To see that the left vertical map is an isomorphism, we decompose it as a composition
\[
\rR\underline{\Gamma}_v(X, \cE) \otimes^{\rL \sol}_{\underline{\Gamma}_v(S, \cA)}
\rR\underline{\Gamma}_v(X', \cE') \to \left(\rR\underline{\Gamma}_v(\widetilde{X},\restr{\cE}{\wX}) \right)^{\gamma = 1}
\otimes^{\rL \sol}_{\underline{\Gamma}_v(S, \cA)}
\rR\underline{\Gamma}_v(X', \cE') \to \left(\rR\underline{\Gamma}_v(\widetilde{X},\restr{\cE}{\wX})
\otimes^{\rL \sol}_{\underline{\Gamma}_v(S, \cA)}
\rR\underline{\Gamma}_v(X', \cE') \right)^{\gamma = 1}.
\]
Then the first map is an isomorphism due to \cref{cohomology and Koszul: one variable}, while the second map is an isomorphism since $\gamma$ does not act on $\rR\underline{\Gamma}_v(X', \cE')$.  This implies that both vertical maps are isomorphism and, therefore, finishes the proof. 
\end{proof}

For the next definition, we fix an affinoid sousperfectoid space $S$ and a smooth morphism $f\colon X \to S$ of sousperfectoid spaces.

\begin{defn} The $S$-space $X$ {\it admits good coordinates} if there is an integer $d$ and an $S$-morphism $g\colon X\to \TT^d_S$ such that $g$ is a composition of rational subdomains and finite \'etale morphisms. 
\end{defn}
Similarly to the classical situation, every smooth space over a sousperfectoid affinoid space admits good coordinates locally in the analytic topology. 
\begin{lemma}\label{lemma:local-good-coordinates} Let $S=\Spa(R, R^+)$ be an affinoid sousperfectoid space over $\Spa(\QQ_p, \ZZ_p)$, let $f\colon X \to S$ be a smooth morphism of sousperfectoid spaces, and let $x\in X$ be a point. Then there is an open affinoid neighborhood $x\in U_x\subset X$ such that $U_x$ admits good coordinates. 
\end{lemma}
\begin{proof}
    The definition of smooth morphisms implies that there is an open affinoid $x\in U_x \subset X$ such that $f|_{U_x}\colon U_x \to S$ factors as a composition of an \'etale morphism $g_x\colon U_x \to \DD^d_{S}$ followed by the projection $\DD^d_{S} \to S$. Furthermore, the definition of an \'etale map implies that (possibly after further shrinking $U_x$) we can assume that $g_x$ is a composition of rational subdomains and finite \'etale maps. Furthermore, we can identify $\DD^d_{S}$ with a rational subdomain $\TT^d_{S}\bigl(\frac{T_1-1}{p},\dots, \frac{T_d-1}{p} \bigr) \subset \TT^d_{S}$, where $T_1, \dots, T_d$ are the coordinates on $\TT^d_{S}$. Therefore, the composition
    \[
    U_x \xr{g_x} \DD^d_{S} \hookrightarrow \TT^d_{S}
    \]
    defines the desired good coordinates on $U_x$. 
\end{proof}

Before we start proving  \Cref{Kunneth formula: main thm}, we need one more preliminary lemma.

\begin{lemma}\label{Kunneth formula:preliminary}
Let $S$ be a strictly totally disconnected perfectoid space over $\Spd(\QQ_p, \ZZ_p)$ corresponding to an untilt $S^\sharp \to \Spa(\Q_p, \ZZ_p)$.
Let $f^\sharp \colon X^\sharp \to S^\sharp$ and $f'^\sharp \colon X'^\sharp \to S^\sharp$ be 
qcqs smooth maps of sousperfectoid spaces, and let $X^\sharp=U^\sharp_1\cup U^\sharp_2$ be an open covering with $U^\sharp_3\coloneqq U^\sharp_1\cap U^\sharp_2$. Set $X\coloneqq (X^\sharp)^\diam$, $X'\coloneqq (X'^\sharp)^\diam$ and $U_i\coloneqq (U_i^\sharp)^\diam$. Let
 $\cE \in \Perf(X_v, \cA)$ and 
$\cE' \in \Perf(X'_v, \cA)$ for $\cA\in \{\O, \BB_I\}$. If the solid K\"unneth maps
\[
\KM^\sol_{U_i, X'}\colon \rR\uGamma_v(U_i, \restr{\cE}{U_i}) \otimes^{\rL \Box}_{\uGamma_v(S, \cA)}
\rR\uGamma_v(X', \cE') \to
\rR\uGamma_v(U_i \times_{S} X', \restr{\cE}{U_i} \boxtimes^{\rL} \cE')
\]
are isomorphisms for $1\leq i\leq 3$, then the solid K\"unneth map
\[
\KM^\sol_{X, X'} \colon \rR\uGamma_v(X, \cE) \otimes^{\rL \Box}_{\uGamma_v(S, \cA)}
\rR\uGamma_v(X', \cE') \rightarrow
\rR\uGamma_v(X \times_{S} X', \cE \boxtimes^{\rL} \cE')
\]
is an isomorphism as well. 
\end{lemma}
\begin{proof}
We consider the Mayer--Vietoris exact triangle
\[
\rR\uGamma_v(X, \cE) \to \rR\uGamma_v(U_1, \cE) \oplus \rR\uGamma_v(U_2, \cE) \to \rR\uGamma_v(U_3, \cE).
\]
Since the solid K\"unneth map is functorial, we get the commutative diagram
\[
\begin{tikzcd}[column sep=15em]
\rR\uGamma_v(X, \cE) \otimes^{\rL\sol}_{\uGamma_v(S, \cA)} \rR\uGamma_v(X', \cE') \arrow{r}{\KM^\sol_{X, X'}} \arrow{d} & \rR\uGamma_v(X \times_{S} X', \cE \boxtimes^{\rL} \cE')\arrow{d} \\
{ \begin{gathered}
    \rR\uGamma_v(U_1, \restr{\cE}{U_1}) \otimes^{\rL\sol}_{\uGamma_v(S, \cA)}\rR\uGamma_v(X', \cE') \\ \oplus \\
    \rR\uGamma_v(U_2, \restr{\cE}{U_2}) \otimes^{\rL\sol}_{\uGamma_v(S, \cA)}\rR\uGamma_v(X', \cE')
  \end{gathered}
}
 \arrow{d}\arrow{r}{\KM^\sol_{U_1, X'} \oplus \KM^\sol_{U_2, X'}} & { \begin{gathered}
    \rR\uGamma_v(U_1 \times_{S} X', \restr{\cE}{U_1} \boxtimes^{\rL} \cE')\\ \oplus \\
    \rR\uGamma_v(U_2 \times_{S} X', \restr{\cE}{U_2} \boxtimes^{\rL}\cE')
  \end{gathered}
} \arrow{d} \\
\rR\uGamma_v(U_3, \restr{\cE}{U_3}) \otimes^{\rL\sol}_{\uGamma_v(S, \cA)} \rR\uGamma_v(X', \cE')\arrow{r}{\KM^\sol_{U_3, X'}} & \rR\uGamma_v(U_3 \times_{S} X', \restr{\cE}{U_3} \boxtimes^{\rL} \cE'),
\end{tikzcd}
\]
where both vertical sequences are exact triangles. Since each $\KM^\sol_{U_i, X'}$ is an isomorphism, we deduce that $\KM^\sol_{X, X'}$ must be an isomorphism as well. 
\end{proof}

Finally, we are ready to prove \Cref{Kunneth formula: main thm}. 

\begin{proof}[Proof of \Cref{Kunneth formula: main thm}]

\begin{enumerate}[wide,label={\textit{Case~\arabic*}.},ref={Case~\arabic*}]
\item\label{kunneth:case-1}\textit{$X^\sharp$ and $X'^\sharp$ admit good coordinates.} Choose good coordinates $g^\sharp\colon X^\sharp \to \TT^d_{S^\sharp}$ and $g'^\sharp\colon X'^\sharp \to \TT^{d'}_{S^\sharp}$ and set $\wX^\sharp \coloneqq \widetilde{\TT}^d_{\TT^d_{S^\sharp}} \times_{S^\sharp} X$ and $\wX'^\sharp \coloneqq \widetilde{\TT}^{d'}_{\TT^{d'}_{S^\sharp}} \times_{S^\sharp} X'$, where $\widetilde{\TT}^d_{S^\sharp}$ and $\widetilde{\TT}^{d'}_{S^\sharp}$ are the affinoid perfectoid spaces from \cref{successive Zp-torsor}.
Then \cite[Th.\,6.3(ii) and 7.9]{Scholze-perfectoid} imply\footnote{Strictly speaking, \cite[Th.\,6.3(ii) and 7.9]{Scholze-perfectoid} require the perfectoid spaces to be defined over a perfectoid field, but this is not needed in the proof.} that $\wX^\sharp$ and $\wX'^\sharp$ are affinoid perfectoid spaces. We also set $\wX \coloneqq (\wX^\sharp)^\diam$ and $\wX'\coloneqq (\wX'^\sharp)^\diam$.
Then \cref{successive Zp-torsor} implies that the morphisms $\wX \to X$ and $\wX' \to X'$ are successive $\ud{\ZZ}_p$-torsors. Now we consider the following commutative diagram
\[
\begin{tikzcd}[column sep = 7 em]
\wX' \times_S \wX \arrow[d] \arrow[rr] &&
\wX \arrow[d, "\text{successive } \ZZ_p\text{-torsor}"] \\
\wX' \times_S X \arrow[d] \arrow[r] &
X' \times_S X \arrow[d] \arrow[r] & X \arrow[d] \\
\wX' \arrow[r, "\text{successive } \ZZ_p\text{-torsor}"] & X' \arrow[r] & S,
\end{tikzcd}
\]
in which all squares are cartesian. Now we apply \Cref{RGamma-solid} and \cref{Kunneth formula: one variable} twice (first to $\wX' \to X'$ and $X \to S$ and then to $\wX \to X$ and $\wX' \to S$)
to conclude that it suffices to show that the solid K\"unneth map
\[
\KM^\sol_{\wX, \wX'} \colon \rR\uGamma_v(\wX, \restr{\cE}{\wX})
\otimes^{\rL \sol}_{\uGamma_v(S, \cA)}
\rR\uGamma_v(\wX', \restr{\cE'}{\wX'}) \to
\rR\uGamma_v\bigl((\wX \times_S \wX'), \restr{\cE}{\wX} \boxtimes^\rL \restr{\cE'}{\wX'}\bigr)
\]
is an isomorphism. Now \cref{cor:perfect-O-modules-on-affinoid-perfectoid} when $\cA=\O$ and \cref{cor:perfect-O-modules-on-affinoid-perfectoid-BBI} when $\cA=\BB_I$ imply that $\restr{\cE}{\wX}\simeq \widetilde{P}$ and $\restr{\cE'}{\wX'} \simeq \widetilde{P}'$ for some $P\in \Perf(\cA(\wX))$ and $P'\in \Perf(\cA(\wX'))$. In particular, $P$ (resp.\,$P'$) is quasi-isomorphic to a finite complex of finite projective $\cA(\wX)$-modules (resp.~$\cA(\wX')$-modules). Since each finite projective module is a retract of a finite free module, a simple inductive argument reduces the question to the case $\restr{\cE}{\wX}=\cA$ and $\restr{\cE'}{\wX'}=\cA$. In this case, the result follows immediately from \cref{lemma:rational-kunneth-for-perfectoids} when $\cA = \O$ and from \cref{Kunneth for BI rings} when $\cA=\BB_I$.

\item\label{kunneth:case-2}\textit{$X^\sharp$ admits an open embedding into a sousperfectoid affinoid $S^\sharp$-spaces $Y^\sharp$ and $X'^\sharp$ admits good coordinates.} 
We prove, by induction on $n$, that $\KM^\sol_{X, X'}$ is an isomorphism whenever $X^\sharp$
can be covered by $n$ rational subdomains of $Y^\sharp$, each of which admits good coordinates.
Note that by \cref{lemma:local-good-coordinates} and the assumption that $X^\sharp$ is quasi-compact,
we can always find a finite covering $X^\sharp = \cup_{i=1}^n U^\sharp_i$ of the desired kind.

If $n=1$, then the result follows from \cref{kunneth:case-1}. Now we assume that we proved the claim $n-1\geq 1$ and try to deduce it for $n$. 
We set $U^\sharp \coloneqq U^\sharp_1$, $U\coloneqq (U^\sharp)^\diam$, $V^\sharp\coloneqq U^\sharp_2\cup U^\sharp_3\cup \dots \cup U^\sharp_n$, and $V\coloneqq (V^\sharp)^\diam$. Let $g\colon U^\sharp \to \TT^d_{S^\sharp}$ be good coordinates on $U^\sharp$. 

We wish to show that $V^\sharp\cap U^\sharp$ has an open covering by $(n-1)$ affinoids admitting good coordinates. For this, we note that each $U^\sharp_i\cap U^\sharp$ is a rational subdomain in $U^\sharp$. Thus, the composition $U^\sharp_i\cap U^\sharp \hookrightarrow U^\sharp \xr{g} \TT^d_{S^\sharp}$ defines good coordinates on $U^\sharp_i\cap U^\sharp$ for $i=2, \dots, n$. Therefore, $\bigcup_{i=2}^n U^\sharp_i \cap U^\sharp = V^\sharp \cap U^\sharp$ is the desired open cover.
Hence, the inductive hypothesis implies that $\KM^\sol_{U, X'}$, $\KM^\sol_{V, X'}$, and $\KM^\sol_{U\cap V, X'}$ are isomorphisms. Thus, \cref{Kunneth formula:preliminary} ensures that $\KM^\sol_{X, X'}$ must be an isomorphism as well. 

\item\label{kunneth:case-3}\textit{$X^\sharp$ is a general smooth qcqs $S^\sharp$-space and $X'^\sharp$ admits good coordinates.}
We prove, by induction on $n$, that $\KM^\sol_{X, X'}$ is an isomorphism whenever $X^\sharp$
can be covered by $n$ affinoid open subsets.
Since $X^\sharp$ is a qcqs space, we can always find a finite open covering
$X^\sharp=\cup_{i=1}^n U^\sharp_i$ such that each $U^\sharp_i$ is affinoid.

If $n=1$, the claim follows from \cref{kunneth:case-2}. Now we assume that we proved the claim for $n-1\geq 1$ and try to deduce it for $n$. We set $U^\sharp\coloneqq U^\sharp_1$, $U\coloneqq (U^\sharp)^\diam$, $V^\sharp\coloneqq U^\sharp_2\cup U^\sharp_3 \cup \dots \cup U^\sharp_n$, and $V\coloneqq (V^\sharp)^\diam$. Then the induction hypothesis implies that $\KM^\sol_{U, X'}$ and $\KM^\sol_{V, X'}$ are isomorphisms, while \cref{kunneth:case-2} implies that $\KM^\sol_{U\cap V, X'}$ is an isomorphism. Therefore, \cref{Kunneth formula:preliminary} ensures that $\KM^\sol_{X, X'}$ must be an isomorphism as well. 

\item\label{kunneth:case-4}\textit{general $X^\sharp$ and $X'^\sharp$.} We can repeat the argument in \cref{kunneth:case-2} with $X^\sharp$ replaced by $X'^\sharp$ to deduce the result for any $X^\sharp$ and for $X'^\sharp$ which admit an open embedding into a sousperfectoid smooth $S^\sharp$-affinoid spaces. Then we can repeat the argument of \cref{kunneth:case-3} to deduce the result for a general qcqs smooth $S^\sharp$-space $X'^\sharp$. \qedhere
\end{enumerate}
\end{proof}

\begin{rmk} 
The only reason we need to assume that $S$ is strictly totally disconnected in the formulation of \cref{Kunneth formula: main thm} is to make sure that $\O(S^\sharp)$ admits all $p$-power roots of unity. In particular, the torsor $\widetilde{\TT}^d_{S^\sharp} \to \TT^d_{S^\sharp}$ is a successive $\ZZ_p$-torsor as opposed to a successive $\ZZ_p(1)$-torsor. In particular, the conclusion of \cref{Kunneth formula: main thm} remains true for any affinoid perfectoid $S$ which admits a map to $\Spd(\QQ_p(\mu_{p^\infty}), \QQ_p(\mu_{p^\infty})^\circ)$. We do not spell out the details as we will never need this fact.
\end{rmk}

In a similar way, we also get the following variants:

\begin{theorem}
\label{solid base change formula}
Let $S$ be a strictly totally disconnected perfectoid space over $\Spd(\QQ_p, \ZZ_p)$ corresponding to an untilts $S^\sharp$ over $\Spa(\QQ_p, \ZZ_p)$, let $S'$ be an affinoid perfectoid space with a map $S' \to S$ over $\Spd(\QQ_p, \ZZ_p)$, let $f^\sharp \colon X^\sharp \to S^\sharp$ be a qcqs smooth map of sousperfectoid spaces, and let $\cA\in \{\BB_I, \O\}$. Set $X\coloneqq (X^\sharp)^\diam$.
Let $\cE \in \Perf_v(X, \cA)$ and 
$\cE' \in \Perf_v(S', \cA)$.
Then the solid K\"{u}nneth map 
\[
\KM^\Box\colon \rR\uGamma_v(X, \cE) \otimes^{\rL \Box}_{\uGamma_v(S, \cA)}
\rR\uGamma_v(S', \cE') \rightarrow
\rR\uGamma_v(X\times_S S', \cE \boxtimes^{\rL} \cE').
\]
is an isomorphism. 
\end{theorem}
\begin{proof}
If $X^\sharp$ admits good coordinates $g^\sharp\colon X^\sharp \to \TT^d_{S^\sharp}$, then the proof is completely analogous to \cref{kunneth:case-1} in the proof of \cref{Kunneth formula: main thm}.  In general, we can argue as in \cref{kunneth:case-2} and \cref{kunneth:case-3} in the proof of \cref{Kunneth formula: main thm} to reduce the general case to the case when $X$ admits good coordinates.   
\end{proof}

\begin{theorem}
\label{solid pullback compatibility}
Let $S$ be a strictly totally disconnected perfectoid space over $\Spd(\QQ_p, \ZZ_p)$
corresponding to an untilts $S^\sharp$ over $\Spa(\QQ_p, \ZZ_p)$, 
let $J \subset I \subset (0, \infty)$ be an inclusion of closed intervals
with rational endpoints, and let $f^\sharp \colon X^\sharp \to S^\sharp$ be a qcqs smooth map of sousperfectoid spaces. Set $X\coloneqq (X^\sharp)^\diam$. Then for any object $\cE \in \Perf_v(X, \BB_I)$, the natural maps
\[
\rR\uGamma_v(X, \cE) \otimes^{\rL \Box}_{\uGamma_v(S, \BB_I), \res}
\uGamma_v(S, \BB_J) \rightarrow
\rR\uGamma_v(X, \cE \otimes^{\rL}_{\BB_I, \res} \BB_J)
\text{ and}
\]
\[
\rR\uGamma_v(X, \cE) \otimes^{\rL \Box}_{\uGamma_v(S, \BB_I), \varphi}
\uGamma_v(S, \BB_{I/p}) \rightarrow
\rR\uGamma_v(X, \cE \otimes^{\rL}_{\BB_I, \varphi} \BB_{I/p})
\]
are isomorphisms.
\end{theorem}

\begin{proof}
The proof is again analogous to that of \cref{Kunneth formula: main thm},
except we need to use \Cref{pullback for BI rings} in place of \Cref{Kunneth for BI rings}.
\end{proof}

Now we wish to give a sheaf analog of \cref{construction:theta-rings} and \cref{cor:base change along theta}.

\begin{construction}\label{construction:theta-sheaves} Let $X$ be a diamond over $\Spd(\QQ_p, \ZZ_p)$ and let $I\subset (0, \infty)$ be a closed interval with rational endpoints such that $1\in I$. Then, by varying the map $\theta_I\colon \BB_I(S) \to \O_X(S)$ from \cref{construction:theta-rings} over all affinoid perfectoid spaces $S \to X$, we get a homomorphism of $v$-sheaves of rings 
\[
\theta_I \colon \BB_I \to \O_X.
\]
\end{construction}

\begin{theorem}
\label{compatibility with HT realization}
Let $S$ be a strictly totally disconnected perfectoid space over $\Spd(\QQ_p, \ZZ_p)$
corresponding to an untilts $S^\sharp$ over $\Spa(\QQ_p, \ZZ_p)$, let $f \colon X^\sharp \to S^\sharp$ be a qcqs smooth map of sousperfectoid spaces with $X\coloneqq (X^\sharp)^\diam$, let $I\subset (0, \infty)$ be a closed interval with rational endpoints such that $1\in I$, and let $\cE \in \Perf_v(X; \BB_{I})$. Then the natural map
\[
\rR\uGamma_v(X, \cE) \otimes^{\rL \Box}_{\uGamma_v(S, \BB_{I}), \theta_I}
\uGamma_v(S, \O) \rightarrow
\rR\uGamma_v(X, \cE \otimes^{\rL}_{\BB_{I}, \theta_I} \cO)
\]
is an isomorphism.
\end{theorem}

\begin{proof}
The proof is again analogous to that of \cref{Kunneth formula: main thm},
except that we need to use \Cref{cor:base change along theta} in place of \Cref{Kunneth for BI rings}.
\end{proof}

\subsection{Poincar\'e duality}\label{smooth-PD-BI}

The main goal of this subsection is to prove Poincar\'e duality for perfect $\BB_I$-complexes (resp.\,$\O$-complexes). We follow the strategy used in \cite[\textsection 6.4]{LRZ24}. Namely, we explicitly construct the evaluation and coevaluation maps, then the duality statement boils down to the verification of commutativity of certain diagrams. 

The main result of this subsection is the following Poincar\'e Duality statement. 

\begin{thm}\label{thm:PD-absolute-solid} Let $S$ be a strictly totally disconnected perfectoid space with a morphism $S \to \Spd(\QQ_p, \ZZ_p)$, let $S^\sharp$ be the corresponding untilt, let $f^\sharp \colon X^\sharp \to S^\sharp$ be a smooth proper morphism of equidimension $d$ between sousperfectoid spaces, let $\cA\in \{\BB_I, \cO\}$. Set $X\coloneqq (X^\sharp)^\diam$. Let $\cE\in \Perf_v(X; \cA)$ with $\cE^\vee\coloneqq \rR{\cHom}_{\cA}(\cE, \cA)$ being its $\cA$-dual. Then the evaluation and
coevaluation maps defined in
\cref{construction:evaluation-coevaluation-poincare-duality} turn
$\rR\uGamma_v(X, \cE) \in \cD_\sol\bigl(\ud{\cA(S)}\bigr)$ into a dualizable object 
with its left dual given by 
$\rR\uGamma_v(X, \cE^{\vee}(d)[2d])$. In particular,  there is a natural isomorphism
\[
\PD^\sol_f \colon \rR\uGamma_v(X, \cE^{\vee}(d)[2d])\xr{\sim} \rR\ud{\Hom}_{\ud{\cA(S)}}\bigl(\rR\uGamma_v(X, \cE), \ud{\cA(S)}\bigr)
\]
\end{thm}

We note that even though \cref{thm:PD-absolute-solid} is formulated using the solid formalism, it almost immediately implies Poincar\'e duality for usual cohomology $\rR\Gamma_v(X, \cA)$ (see \cref{cor:duality-non-solid}). However, it is crucial to first prove \cref{thm:PD-absolute-solid} in order to obtain a classical duality. The key reason why it is important for us to work in the solid framework is the solid K\"unneth formula \cref{Kunneth formula: main thm}, which is critically used in the proof of \cref{thm:PD-absolute-solid}. While we will later show that the classical K\"unneth formula also does hold for smooth proper $X^\sharp$ (see \cref{cor:KF-for-proper-smooth-pushforward}), we do not know how to prove it without proving \cref{thm:PD-absolute-solid} first. 

As a consequence of \cref{thm:PD-absolute-solid}, we also prove a relative version of duality that applies to any smooth proper morphism $f\colon X \to Y$ of equidimension $d$ between $\Spd(\QQ_p, \ZZ_p)$-diamonds, see \cref{cor:PD-for-proper-smooth-pushforward}.

We begin with the following general version of the projection formula. 

\begin{proposition}[{\cite[\href{https://stacks.math.columbia.edu/tag/0943}{Tag 0943}]{stacks-project}}]
\label{PF}
Let  $f \colon (X, \mathcal{O}_X) \to (Y, \mathcal{O}_Y)$ be a morphism of ringed topoi, let $\cE \in \cD(X, \mathcal{O}_X)$ and $\mathcal{K} \in \cD(Y, \cO_Y)$. Then there is a canonical arrow 
$$\PF_f \colon \rR f_* \cE \otimes^\rL_{\cO_Y} \mathcal{K} \to \rR f_*(\cE \otimes^\rL_{\cO_X} \rL f^* \mathcal{K}).$$
Furthermore, $\PF_f$ is an equivalence if $\mathcal{K} \in \Perf(Y, \cO_Y)$.
\end{proposition}

Using \cref{PF}, we can give the following constructions of the evaluation and coevaluation maps. Later, we check that these maps indeed define the structure of a dualizable
object on $\rR\Gamma_v(X, \cE)$.

\begin{construction}[{cf.\ \cite[Construction 6.4.2]{LRZ24}}]
\label{construction:evaluation-coevaluation-poincare-duality}

Let $S$ be a strictly totally disconnected perfectoid space over $\Spd(\QQ_p, \ZZ_p)$ corresponding to the untilt $S^\sharp$, let $X^\sharp$ be a sousperfectoid space, and let $f^\sharp\colon X^\sharp \to S^\sharp$ be a smooth proper morphism of equidimension $d$, let $\cA\in \{\BB_I, \cO\}$. Set $X\coloneqq (X^\sharp)^\diam$ and $f\coloneqq (f^\sharp)^\diam$. Let $\cE\in \Perf_v(X; \cA)$. We also set $\cE^{\vee} \coloneqq \rR{\cHom}(\cE, \cA)$ and denote by $ev_{\cE} \colon \cE^{\vee} \otimes^\rL \cE \to \cA$  and $coev_{\cE} \colon \cA \to \cE \otimes^\rL \cE^{\vee}$
the natural evaluation and coevaluation maps. 

Consider the following commutative diagram:
\begin{equation}\label{eqn:diagonal}
\begin{tikzcd}
X \arrow[rd,"\Delta"]\arrow[rdd, swap, bend right, "\rm{id}"] \arrow[rrd, bend left, "\rm{id}"]& & \\
& X\times_S X \arrow{r}{\pi_2} \arrow{d}{\pi_1} \arrow{rd}{h} & X \arrow{d}{f} \\
& X \arrow[r, swap, "f"] & S
\end{tikzcd} 
\end{equation}

\begin{enumerate}[label=\upshape{(\roman*)},leftmargin=*]
\item\label{evaluation map} (\emph{Evaluation map}) 
We define the \emph{evaluation map} 
\[
e^\sol(f, S, \cE) \colon 
\rR\uGamma_v(X, \cE^{\vee}(d)[2d])
\otimes^{\rL \sol}_{\uGamma_v(S, \cA)} \rR\uGamma_v(X, \cE)
\to \uGamma_v(S, \cA)
\]
as the composition
\[
\begin{tikzcd}
\rR\uGamma_v(X, \cE^{\vee}(d)[2d]) 
\otimes^{\rL \sol}_{\uGamma_v(S, \cA)} \rR\uGamma_v(X, \cE) \arrow[d, dashed] \arrow{rr}{\cup} & &
\rR\uGamma_v(X, \cE^{\vee}(d)[2d] \otimes^{\rL} \cE)
\arrow{d}{\rR\uGamma_v(ev_{\cE}(d)[2d])} \\
\uGamma_v(S, \cA) & \arrow{l}{\rR\uGamma_v(\tr^{\cA}_f)} \rR\uGamma_v(S, \rR f_{v, *} \cA(d)[2d]) & \arrow{l}{\sim}\rR\uGamma_v(X, \cA(d)[2d])
\end{tikzcd}
\]
where the arrow $\cup$
is the solid K\"{u}nneth map from \Cref{Construction:Kunneth map}
(applied to the outer square in the above commutative diagram)
and $\tr^{\cA}$ is the trace map from \Cref{trace-maps}. 
\item\label{coevaluation map} (\emph{Coevaluation map}) We define the \emph{coevaluation map} 
\[
c^\sol(f, S, \cE) \colon \uGamma_v(S, \cA) \to
\rR\uGamma_v(X, \cE)
\otimes^{\rL \sol}_{\uGamma_v(S, \cA)} \rR\uGamma_v(X, \cE^{\vee}(d)[2d])
\]
as the composition
\[
\begin{tikzcd}    
\uGamma_v(S, \cA) \arrow[r,"f^*"]\arrow[ddd, dashed] & \rR\uGamma_v(X, \cA)
\arrow[r,"\rR\uGamma_v(X\text{,} coev_{\cE})"] &
\rR\uGamma_v(X, \cE\otimes^\rL \cE^{\vee})\arrow[d, "\wr"']\\
& & \rR\uGamma_v\bigl(X^{2/S}, \rR\Delta_{v, *}\Delta^*(\pi_1^* \cE \otimes^\rL \pi_2^*\cE^{\vee})\bigr) \arrow[d, "\rR\uGamma_v\bigl(\PF^{-1}_\Delta\bigr)", "\wr"']
\\
& & \rR\uGamma_v\bigl(X^{2/S}, (\pi_1^*\cE \otimes^\rL \pi_2^*\cE^{\vee}) \otimes^\rL \rR\Delta_{v, *}\cA\bigr) \arrow{d}{\rR\uGamma_v(\id\otimes \cl_\Delta^\cA)} \\
\rR\uGamma_v(X, \cE)
\otimes^{\rL \sol}_{\uGamma_v(S, \cA)} \rR\uGamma_v(X, \cE^{\vee}(d)[2d]) & &\rR\uGamma_v(X^{2/S}, \pi_1^*\cE\otimes^\rL \pi_2^*\cE^{\vee}(d)[2d]),  \arrow[ll, "\KM^{\sol, -1}"', "\sim"]
\end{tikzcd}
\]
where $X^{2/S}$ denotes $X\times_S X$, $\cl_\Delta^\cA$ is the cycle class map from \cref{cycle-class-definition}, 
and $\KM^\sol$ is the K\"unneth map from \Cref{Construction:Kunneth map} 
(applied to the bottom square in the commutative diagram above).
We crucially use \Cref{PF} and the K\"{u}nneth formula (\cref{Kunneth formula: main thm}) to invert the projection formula map and the solid K\"{u}nneth map, respectively. 
\item\label{PD map} (\emph{Duality map}) We define the \emph{duality map} $\rm{PD}^\sol_f(\cal{E}) \colon \rR \uGamma_v(X, \cE^{\vee}(d)[2d]) \to \rR \ud{\Hom}_{\ud{\cA(S)}}\bigl(\rR \uGamma_v(X, \cE), \ud{\cA(S)}\bigr)$ as the map adjoint to $e^\sol(f, S, \cal{E})$ under the tensor-hom adjunction.
In other words, $\PD^\sol_f(\cE)$ is the composition
\[
\rR \uGamma_v(X, \cE^{\vee}(d)[2d]) \longrightarrow \rR \ud{\Hom}_{\ud{\cA(S)}}\bigl(\rR \uGamma_v(X, \cE), \rR \uGamma_v(X, \cA(d)[2d])\bigr) \xrightarrow{\tr^\cA_f \circ -} \rR \ud{\Hom}_{\ud{\cA(S)}}\bigl(\rR \uGamma_v(X, \cE), \ud{\cA(S)}\bigr),
\]
where the first map comes from \cite[\href{https://stacks.math.columbia.edu/tag/0B6D}{Tag~0B6D}]{stacks-project}.
\end{enumerate}
\end{construction}

Now we verify certain compatibilities that will be crucial for the proof of \cref{thm:PD-absolute-solid}.

\begin{lemma}[{cf.\ \cite[Lem.~6.4.3]{LRZ24}}]
\label{lemma:base-change-trace} 
Let $S$ be a strictly totally disconnected perfectoid with a morphism $S \to \Spd(\QQ_p, \ZZ_p)$, let $S^\sharp$ be the corresponding untilt, let $f^\sharp\colon X^\sharp \to S^\sharp$ be a smooth proper morphism of equidimension $d$ between sousperfectoid spaces,
let $\cA\in \{\BB_I, \cO\}$. Set $X\coloneqq (X^\sharp)^\diam$, $X^{2/S}\coloneqq X\times_S X$, and $f\coloneqq (f^\sharp)^\diam$. Let $\cE\in \Perf_v(X; \cA)$ and let $\pi_1\colon X^{2/S} \to X$ be the projection onto the first factor. Then the diagram
\[
\begin{tikzcd}[column sep = huge]
\rR\underline{\Gamma}_v(X, \cE) \otimes^{\rL\Box} \rR\underline{\Gamma}_v(X,\cA(d)[2d]) \arrow{d}{\KM^\sol} 
\arrow{r}{\id \otimes^\rL \rR\uGamma_v(\tr^{\cA}_f)} & \rR\underline{\Gamma}_v(X, \cE) \\
\rR\underline{\Gamma}_v\bigl(X^{2/S}, \pi^*_1 \cE \otimes^\rL \cA(d)[2d]\bigr)
\arrow{r}{\rR\underline{\Gamma}_v(\PF^{-1}_{\pi_1})}
& \rR\underline{\Gamma}_v\bigl(X, \cE \otimes^\rL \rR \pi_{1, v, *}\cA(d)[2d]\bigr) 
\arrow{u}{\rR\underline{\Gamma}_v(\rm{id} \otimes^\rL \tr_{\pi_1}^{\cA})}
\end{tikzcd}
\]
commutes.
\end{lemma}
\begin{proof}
We set $h\colon X^{2/S} \to S$ to be the natural morphism. Then we consider the diagram
\[ \begin{tikzcd}[scale cd=.9,center picture,row sep=huge]
\rR\uGamma_v(X, \cE) \otimes^{\rL\Box} \rR\uGamma_v(X,\cA(d)[2d]) \arrow{rd}{\cup} \arrow{dd}{\KM^\sol} \arrow{rr}{\id\otimes \rR\uGamma_v(\tr^{\cA})}& & \rR\uGamma_v(X, \cE) \arrow[d, equals] \\
& \rR\uGamma_v(S, \rR f_{v, *} \cE \otimes^\rL \rR f_{v, *}\cA(d)[2d]) \arrow[r, red, "\rR\uGamma_v(\id \otimes^\rL \tr_f^\cA)"] \arrow[d, red, "\rR\uGamma_v(\KM_{f,f})"] & \rR\uGamma_v(S, \rR f_{v, *}\cE)  \\
\rR\uGamma_v\bigl(X^{2/S}, \pi_1^*\cE \otimes^\rL \cA(d)[2d]\bigr)\arrow[r, equals] \arrow[rrd, swap, "\rR\uGamma_v(\PF_{\pi_1}^{-1})"]& \rR\uGamma_v\bigl(S, \rR h_*(\pi_1^*\cE \otimes^\rL \cA(d)[2d])\bigr) \arrow[r, red, "\alpha"] &\rR\uGamma_v\bigl(S, \rR f_{v, *}(\cE \otimes^\rL \rR \pi_{1, v, *}\cA(d)[2d])\bigr) \arrow[u, red, "\rR\uGamma_v(\rR f_{v, *}(\id \otimes^\rL \tr_{\pi_1}^\cA))"]  \arrow[d, equals] \\
& & \rR\uGamma_v\bigl(X, \cE \otimes^\rL \rR \pi_{1, v, *}\cA(d)[2d]\bigr) \arrow[uuu, start anchor=east, bend right=70, "\rR\uGamma_v\bigl(\rm{id} \otimes^\rL \tr_{\pi_1}^{\cA}\bigr)"],
\end{tikzcd} \]
where $\alpha=\rR\uGamma_v\bigl(\rR f_{v,*}(\rm{PF}^{-1}_{\pi_1})\bigr)$. After unravelling all the definitions, we see that all external ``squares'' commute. Therefore, it suffices to show that the red square commutes. For this, it suffices to show that the following diagram commutes:
\[
\begin{tikzcd}[column sep = huge]
\rR f_{v, *} \cE \otimes^\rL \rR f_{v, *}\cA(d)[2d] \arrow{d}{\rm{KM}} \arrow{r}{\id \otimes^\rL \tr^{\cA}_f}
& \rR f_{v, *} \cE \\
\rR h_{v, *}(\pi_{1}^*\cE \otimes^\rL \cA (d)[2d]) \arrow{r}{\rR f_{v, *}(\rm{PF}^{-1}_{\pi_1})} &
\rR f_{v, *}(\cE \otimes^\rL \rR \pi_{1, v, *} \cA(d)[2d]) \arrow{u}{\rR f_{v, *}(\rm{id} \otimes^\rL \tr_{\pi_1}^{\cA})} \\
\end{tikzcd}
\]

Now the proof of \cite[Lem.\,6.4.3]{LRZ24} applies essentially verbatim. 
Indeed, \cite[Lem.\,6.3.8]{LRZ24} works for any commutative diagram of
ringed sites, and \cite[Th.\,6.1.1(2)]{LRZ24} can be replaced by 
\cref{B_I trace compatible under pullback}.
\end{proof}

\begin{lemma}[{cf.\ \cite[Lem.~6.4.5]{LRZ24}}]
\label{lemma:preliminary-duality}
Let $S$ be a strictly totally disconnected perfectoid with a morphism $S \to \Spd(\QQ_p, \ZZ_p)$, let $S^\sharp$ be the corresponding untilt, let $X^\sharp$ and $Y^\sharp$ be smooth sousperfectoid spaces over $S^\sharp$ of relative equidimension $d_X$ and $d_Y$, respectively. Let $f^\sharp\colon X^\sharp \to Y^\sharp$ be a smooth proper morphism of equidimension $d\coloneqq d_X - d_Y$, let $s^\sharp \colon Y^\sharp \to X^\sharp$ be a section of $f^\sharp$, and let $\cA\in \{\BB_I, \cO\}$. We set $X\coloneqq (X^\sharp)^\diam$, $Y\coloneqq (Y^\sharp)^\diam$, $f\coloneqq (f^\sharp)^\diam$, and $s\coloneqq (s^\sharp)^\diam$. Then the compositions
\begin{gather*}
\cA_{Y} \simeq \rR f_{v, *} \rR s_{v, *}\cA_{Y}  \xlongrightarrow{\rR f_{v, *}\cl_s^\cA} \rR f_{v, *} \cA_{X}(d)[2d]
\xlongrightarrow{\tr^\cA_f} \cA_{Y} \quad \text{and} \\
\begin{multlined}[.92\textwidth]
\cA_{Y}(d)[2d]  \simeq \rR f_{v, *} \rR s_{v, *}\bigl(\cA_{Y}(d)[2d]\bigr)
\xrightarrow{\rR f_{v, *}(\cl_s^{\cA} \otimes^\rL \id)} \rR f_{v, *}(\cA_{X}(d)[2d] \otimes^\rL \cA_{X}(d)[2d]) \\
\xrightarrow{\PF_{f}^{-1}} \rR f_{v, *} \bigl(\cA_{X}(d)[2d]\bigr) \otimes^\rL \cA_{Y}(d)[2d]\xrightarrow{\tr^\cA_f \otimes^\rL \id}  \cA_{Y}(d)[2d]
\end{multlined}
\end{gather*}
are equal to the identity morphisms.
\end{lemma}

\begin{proof}
The first composition is equal to the identity morphism due to \cref{prop:cycle-class-of-section}. As for the second composite, the proof of \cite[Lem.~6.4.5]{LRZ24} works verbatim
and implies that the second composite is the first composite shifted by $2d$ and twisted by $d$, from which the claim follows.
\end{proof}

\begin{proof}[Proof of \Cref{thm:PD-absolute-solid}]
Keep the notation of
\Cref{thm:PD-absolute-solid}. We shall mimic the proof of \cite[Th.~6.4.1]{LRZ24}.
With the \Cref{construction:evaluation-coevaluation-poincare-duality},
the only thing we need to check is that the compositions\footnote{In these formulas, we implicitly make the usual identifications
$\rR\uGamma_v(X, \cE) \otimes^{\sol, \rL} \underline{\cA(S)} \simeq \rR\uGamma_v(X, \cE)$, 
$\underline{\cA(S)} \otimes^{\sol, \rL} \rR\uGamma_v(X, \cE) \simeq \rR\uGamma_v(X, \cE)$, etc.
To indicate this subtlety, we use quotation marks in the maps which implicitly use these identifications.}
\[ \begin{tikzcd}[scale cd=.85,center picture,column sep=7em]
\rR\uGamma_v(X, \cE) \arrow{r}{``c^\sol(f, S, \cE) \otimes^{\rL\sol} \rm{id}\text{''}} & 
\rR\uGamma_v(X, \cE) \otimes^{\rL\sol} \rR\uGamma_v(X, \cE^{\vee}(d)[2d]) \otimes^{\rL\sol} \rR\uGamma_v(X, \cE)
\arrow{r}{``\rm{id}\otimes^{\rL\sol} e^\sol(f, S, \cE)\text{''}}  & \rR\uGamma_v(X, \cE),
\end{tikzcd} \]
\[ \begin{tikzcd}[scale cd=.85,center picture,column sep=6em]
\rR\uGamma_v(X, \cE^{\vee}(d)[2d]) \arrow{r}{``\rm{id} \otimes^{\rL\sol} c^\sol(f, S, \cE)\text{''}} & 
\rR\uGamma_v(X, \cE^{\vee}(d)[2d]) \otimes^{\rL\sol} \rR\uGamma_v(X, \cE) \otimes^{\rL\sol} 
\rR\uGamma_v(X, \cE^{\vee}(d)[2d])
\arrow{r}{``e^\sol(f, S, \cE) \otimes^{\rL\sol} \rm{id}\text{''}} & \rR\uGamma_v(X, \cE^{\vee}(d)[2d])
\end{tikzcd} \]
are equivalent to the identity morphisms. 
For brevity, we denote the first composition by $\varphi(f, S, \cE)$ and the second composition by $\psi(f, S, \cE)$. 
We give a full justification for why $\varphi(f, S, \cE) \simeq \rm{id}$, and then only describe the necessary changes to justify that
$\psi(f, S, \cE) \simeq \rm{id}$.

Recall that we denote by $\eta$ the unit of the (derived) pullback-pushforward adjunction, by $\KM^\sol$  the solid K\"unneth map from 
\Cref{Construction:Kunneth map}, by $\PF$ the projection formula map, and, for any objects $\F, \G\in \cD(X_v; \cA)$,
by $\sigma\colon \F \otimes^\rL \G \xr{\sim} \G\otimes^\rL \F$ the morphism that ``swaps factors''.
In what follows, we also freely use \Cref{Kunneth formula: main thm} and \Cref{PF} which guarantee that 
the K\"unneth map and the projection formula map are isomorphisms under some assumptions that are always satisfied in this proof. 

That being said, we resume the notation of \cref{construction:evaluation-coevaluation-poincare-duality} and consider the diagram in \cref{gigantic-diagram}. 
Using the definitions of $\PF$, $\KM$, and the basic properties of adjunctions, one can check that this diagram commutes. 
For the most part, commutativity follows from functoriality of various constructions.
The following two points need extra arguments: 
triangle~\cref{gigantic-diagram-evcoev} commute since $ev_{\cE}$ and $coev_{\cE}$ 
define a duality datum on $\cE$ (see \cite[\href{https://stacks.math.columbia.edu/tag/0FPU}{Tag 0FPU}]{stacks-project}), and trapezoid~\cref{gigantic-diagram-base-change-trace}
commutes due to \cref{lemma:base-change-trace}. 
    
\begin{sidewaysfigure}

\thisfloatpagestyle{empty}
\tiny

\vspace{60em}
\[ \begin{tikzcd}[center picture,row sep=huge, column sep = 6em]
\underline{\cA(S)} \otimes^\sol \rR\underline{\Gamma}\bigl(X, \cE\bigr) \arrow[rr, "\KM^\sol"] 
\arrow[d,"f^* \otimes \rm{id}"] & & \rR\underline{\Gamma}\bigl(X, \cA\otimes \cE\bigr)
\arrow[d,"\pi_2^*"] \arrow[rd, start anchor=east,bend left=15,"\PF^{-1}_\Delta"] && \\
\rR\underline{\Gamma}\bigl(X,\cA\bigr) \otimes^\sol \rR\underline{\Gamma}\bigl(X, \cE\bigr) \arrow[rr,"\KM^\sol"] \arrow[d,"coev_{\cE} \otimes \id"] & & \rR\underline{\Gamma}\bigl(X^2, \cA \boxtimes \cE\bigr)
\arrow[d,"coev_{\cE} \boxtimes \id"] \arrow{r}{\eta_{\Delta} \otimes \rm{id}} & 
\rR\underline{\Gamma}\bigl(X^2, \rR\Delta_*\cA \otimes \pi_2^* \cE\bigr) 
\arrow{d}{\rR\Delta_*(coev_{\cE})\otimes \rm{id}}
\arrow[dr, start anchor=east, bend left=15, "\PF_\Delta"] &  \\
\rR\underline{\Gamma}\bigl(X, \cE \otimes \cE^\vee\bigr) \otimes^\sol \rR\underline{\Gamma}\bigl(X, \cE\bigr) \arrow[rr,"\KM^\sol"] \arrow[d,sloped, equals] & & 
\rR\underline{\Gamma}\bigl(X^2, (\cE \otimes \cE^\vee) \boxtimes \cE\bigr) \arrow[d,sloped, equals] 
\arrow{r}{\eta_{\Delta}\otimes \rm{id}} & 
\rR\underline{\Gamma}\bigl(X^2, \rR\Delta_*(\cE\otimes \cE^{\vee}) \otimes \pi_2^*\cE\bigr)
\arrow{dd}{\PF_\Delta} & \rR\underline{\Gamma}\bigl(X, \cA \otimes \cE\bigr)
\arrow{ldd}{coev_{\cE}\otimes \rm{id}}  \arrow{dd}{\sigma} \arrow[dddl, eq=gigantic-diagram-evcoev] \\
\rR\underline{\Gamma}\bigl(X, \Delta^*(\cE \boxtimes \cE^\vee)\bigr) \otimes^\sol \rR\underline{\Gamma}(X, \cE) \arrow[r,"\pi_1^* \otimes \pi_2^*"] \arrow[d,equals] & 
\rR\underline{\Gamma}\bigl(X^2, \pi_1^*\Delta^*(\cE \boxtimes \cE^\vee)\bigr) \otimes^\sol 
\rR\underline{\Gamma}\bigl(X^2, \pi_2^* \cE\bigr)\arrow{r}{\KM^\sol} \arrow{d}{\eta_\Delta \otimes \rm{id}} &  \rR\underline{\Gamma}\Bigl(X^2, \pi_1^*\Delta^*(\cE \boxtimes \cE^\vee) \otimes \pi^*_2\cE \Bigr) 
\arrow[d, "\eta_\Delta \otimes \rm{id}"]  & &  \\
\rR\underline{\Gamma}\Bigl(X^2, \Delta_* \Delta^* (\cE \boxtimes \cE^\vee)\Bigr) \otimes^\sol \rR\underline{\Gamma}(X, \cE) \arrow[r,"\rm{id} \otimes \pi_2^*"] \arrow[d,"\PF^{-1}_\Delta \otimes \id"]&  \rR\underline{\Gamma}\Bigl(X^2, \Delta_* \Delta^* (\cE \boxtimes \cE^\vee)\Bigr) \otimes^\sol \rR\underline{\Gamma}\bigl(X^2, \pi_2^* \cE\bigr)\arrow{r}{\KM^\sol} \arrow{d}{\PF^{-1}_\Delta \otimes \rm{id}} & \rR\underline{\Gamma}\Bigl(X^2, \Delta_{*} \Delta^*(\cE \boxtimes \cE^\vee) \otimes \pi^*_2\cE \Bigr) \arrow{d}{\PF^{-1}_{\Delta}\otimes \rm{id}}  \arrow{r}{\PF_{\Delta}} &
\rR\underline{\Gamma}(X, \cE \otimes \cE^\vee \otimes \cE)
\arrow{r}{\rm{id}\otimes ev_{\cE}} \arrow{d}{\PF^{-1}_\Delta} &   
\rR\underline{\Gamma}(X, \cE \otimes \cA)
\arrow{d}{\PF^{-1}_\Delta} \\
\rR\underline{\Gamma}\bigl(X^2, (\cE \boxtimes \cE^\vee)
\otimes \rR\Delta_* \cA \bigr) \otimes^\sol \rR\underline{\Gamma}(X, \cE)  \arrow[d,"(\rm{id}\otimes \cl_\Delta^{\cA}) \otimes \rm{id}"] \arrow[r,"\rm{id}\otimes  \pi_2^*"] & 
\rR\underline{\Gamma}\bigl(X^2, (\cE \boxtimes \cE^\vee)
\otimes \rR\Delta_* \cA \bigr) \otimes^\sol \rR\underline{\Gamma}\bigl(X^2, \pi_2^* \cE\bigr) 
\arrow{r}{\KM^\sol} \arrow[d,"(\rm{id}\otimes \cl_\Delta^{\cA}) \otimes \rm{id}"] & 
\rR\underline{\Gamma}\bigl(X^2, (\cE \boxtimes \cE^\vee) \otimes \rR\Delta_* \cA \otimes \pi_2^* \cE \bigr)
\arrow[d,swap, "\rm{id}\otimes \cl_\Delta^{\cA} \otimes \rm{id}"] \arrow[rr, bend right=15, "\rm{id} \otimes \pi_2^*(ev_{\cE}) \otimes \rm{id}" swap] \arrow{r}{\rm{id} \otimes \PF_{\Delta}}   &
\rR\underline{\Gamma}\bigl(X^2, \pi_1^* \cE \otimes \Delta_*(\cE^\vee \otimes \cE)\bigr)
\arrow{r}{\rm{id}\otimes \Delta_*(ev_{\cE})} & 
\rR\underline{\Gamma}\bigl(X^2, \pi_1^*\cE \otimes \Delta_* \cA \bigr)
\arrow[dddl, end anchor=east, bend left, swap, "\rm{id}\otimes \cl_\Delta^{\cA}"] \arrow[ddd, "\PF^{-1}_{\pi_1}"] \\
\rR\underline{\Gamma}\bigl(X^2, (\cE \boxtimes \cE^\vee)(d)[2d]\bigr) \otimes^\sol \rR\underline{\Gamma}(X, \cE) 
\arrow[r,"\rm{id} \otimes \pi_2^*"] \arrow[d,"(\KM^\sol)^{-1} \otimes \id"] & \rR\underline{\Gamma}\bigl(X^2,
(\cE \boxtimes \cE^\vee)(d)[2d]\bigr) \otimes^\sol \rR\underline{\Gamma}\bigl(X^2, \pi_2^* \cE\bigr)
\arrow{r}{\KM^\sol} & \rR\underline{\Gamma}\bigl(X^2, (\cE \boxtimes \cE^\vee)(d)[2d] \otimes \pi_2^* \cE\bigr) \arrow[dd,equals] && \\
\rR\underline{\Gamma}(X, \cE) \otimes^\sol \rR\underline{\Gamma}(X, \cE^\vee (d)[2d]) \otimes^\sol \rR\underline{\Gamma}(X, \cE)  \arrow[d,"\id \otimes \KM^\sol"] & & \\
\rR\underline{\Gamma}(X, \cE) \otimes^\sol \rR\underline{\Gamma}(X, \cE^\vee (d)[2d] \otimes  \cE) 
\arrow{d}{\rm{id} \otimes ev_{\cE}(d)[2d]} \arrow{rr}{\KM^\sol}  &  & 
\rR\underline{\Gamma}\bigl(X^2, (\cE \boxtimes \cE^\vee)(d)[2d] \otimes \pi_2^* \cE\bigr) \arrow{r}{\rm{id}\otimes \pi_2^*(ev_{\cE}(d)[2d])} \arrow{d}{\rm{id} \otimes \pi_2^*(ev_{\cE}(d)[2d])} &
\rR\underline{\Gamma}\bigl(X^2, \pi_1^*\cE \otimes \cA(d)[2d]\bigr)
\arrow{d}{\PF^{-1}_{\pi_1}} & \rR\underline{\Gamma}(X, \cE \otimes \cA)  \arrow{dl}{\rm{id} \otimes \rR \pi_{1, *}(\cl_\Delta^{\cA})} \\
\rR\underline{\Gamma}(X, \cE) \otimes^\sol \rR\underline{\Gamma}(X,\cA(d)[2d]) \arrow[d,"\id \otimes \tr_f^{\cA}"] \arrow{rr}{\KM^\sol} \arrow[rrd, eq=gigantic-diagram-base-change-trace] 
& & \rR\underline{\Gamma}\bigl(X^2, \pi^*_1 \cE \otimes \cA(d)[2d]\bigr)  
\arrow{r}{\PF^{-1}_{\pi_1}} &  \rR\underline{\Gamma}\bigl(X, \cE\otimes \rR \pi_{1, *}\cA(d)[2d]\bigr)  
\arrow{ld}{\rm{id} \otimes \tr_{\pi_1}^{\cA}} &   \\
\rR\underline{\Gamma}(X, \cE) \otimes^\sol \underline{\cA(S)} & & \arrow[ll, "(\KM^\sol)^{-1}"] \rR\underline{\Gamma}(X, \cE \otimes \cA). &&
\end{tikzcd} 
\]
\caption{To ease the notation, we adapt the following list of conventions for this diagram. We drop the $(-)^\rL$ symbol above $\otimes$. The symbol $X^2$ denotes $X\times_S X$. All pushforwards and pullbacks are taken with respect to the $v$-topology. In particular, the notation $\rR\uGamma$ means $\rR\uGamma_v$. For any two $\F, \G\in \cD(X_v, \cA)$, the notation $\F\boxtimes \G$ denotes $\pi_1^* \F \otimes^\rL \pi_2^*\G$. Lastly, whenever the context is clear, we label an arrow induced by a map of sheaves simply by the sheaf map itself.}\label{gigantic-diagram}
\end{sidewaysfigure}   

The map $\underline{\cA(S)} \otimes^{\rL\sol} \rR\underline{\Gamma}_v\bigl(X, \cE\bigr) \to 
\rR\underline{\Gamma}_v(X, \cE) \otimes^{\rL\sol} \underline{\cA(S)}$ obtained by going down the entire left column is equal to
$\varphi(f, S, \cE)$
by its very construction. The commutativity of the diagram in \cref{gigantic-diagram} implies that this composition can be computed by going around the outer diagram from the top left corner to the bottom left corner in a clockwise direction.
Furthermore, we see that \cref{lemma:preliminary-duality} and the formula $\rR f_*(\PF_{\pi_1}^{-1}) \circ \rR h_*(\PF_{\Delta}^{-1}) = \rm{id}$ imply that 
\[
\bigl(\rm{id}\otimes^{\rL\sol} e^\sol(f, S, \cE)\bigr) \circ \bigl(c^\sol(f, S, \cE) \otimes^{\rL\sol} \rm{id}\bigr) = \KM^{-1} \circ 
\rR\underline{\Gamma}_v(X,\sigma) \circ \KM = \sigma \colon 
\underline{\cA(S)} \otimes \rR\underline{\Gamma}_v\bigl(X, \cE\bigr) \to 
\rR\underline{\Gamma}_v(X, \cE) \otimes \underline{\cA(S)}.
\]
This formally implies that $\varphi(f, S, \cE) = \rm{id} \colon \rR\uGamma_v(X, \cE) \to \rR\uGamma_v(X, \cE)$. 

To see that $\psi(f, S, \cE)=\rm{id}$, we need to use a diagram similar to that of \cref{gigantic-diagram}; we leave it to the reader to figure out the exact shape of the diagram. We only mention that every instance of $\pi_1$ should be replaced with $\pi_2$ (and vice versa) and one needs to use the second part of \cref{lemma:preliminary-duality} (as opposed to the first part used in the proof above).
\end{proof}

\begin{remark}
\label{rmk:PD-literal-analog-of-LRZ}
Let us note that our proof of \Cref{thm:PD-absolute-solid} is modeled on
the proof of \cite[Th.~6.4.1]{LRZ24}. 
In fact, we can get similar result to \cite[Th.~6.4.1]{LRZ24}, but in the current setting
of smooth proper maps $X \to Y$ of equidimension $d$ between diamonds over $\Spd(\QQ_p, \ZZ_p)$ and perfect complexes over $\cA \in \{\ZZ/n, \ZZ_p\}$, as follows.
First, for any $\cE \in \Perf_v(X; \cA)$
and $\cE^{\vee}\coloneqq \rR {\cHom}(\cE, \cA)$, we get a natural map
$e(f, \cE) \colon \rR f_{v, *} (\cE^{\vee}(d)[2d]) \otimes^{\rL}_{\cA} \rR f_{v, *} (\cE) 
\to \cA$ defined as the composition of the cup product followed by the trace map from \Cref{trace-maps}.
Then similar to \cite[Th.~6.4.1]{LRZ24}, we claim that $\rR f_{v, *} \cE$ and 
$\rR f_{v, *} (\cE^{\vee}(d)[2d])$ are in $\Perf_v(Y, \cA)$ and
the map $e(f, \cE)$ 
induces a natural isomorphism
\[
\PD_f(\cE) \colon \rR f_{v, *} \cE^{\vee}(d)[2d] \xrightarrow{\sim} \rR{\cHom}(\rR f_{v, *} \cE, \cA \bigr).
\]

It suffices to check the statement $v$-locally on $Y$, therefore we are immediately reduced to
the case when $Y$ is a perfectoid space and hence $X$ is a sousperfectoid space smooth
of equidimension $d$ and proper over $Y$.
In this case, we claim that one can define a coevaluation map which, together with the map $e(f, \cE)$, makes
$\rR f_* (\cE^{\vee}(d)[2d])$ into a left dual of $\rR f_*\cE$.

For the construction of the coevaluation map, we shall mimic \cite[Construction 6.4.2.(iii)]{LRZ24}.
The construction made in \textit{loc.\ cit.}\ has two ingredients: the cycle class map and the K\"{u}nneth formula.
For the former in our new setting, we can use the cycle class maps defined in \Cref{cycle-class-definition}.

As for the K\"{u}nneth formula, following the strategy in \cite[Section 6.3]{LRZ24}, it suffices to prove
the base change theorem for fiberwise proper morphisms. Then we may reduce to proving such a base change theorem
when $\cA = \ZZ/n$, so we can further reduce to proving it in the \'{e}tale topology. In the end, by looking at stalks at qpro\'et-geometric points, we are eventually reduced to 
the proper base change theorem \cite[Th.~3.15]{BH} (which itself crucially relies on Scholze's comparison theorem). 

Once all these ingredients are there, we may prove the desired dualizability
by repeating the proof of \Cref{thm:PD-absolute-solid}.
\end{remark}

Now we wish to deduce some consequences of \cref{thm:PD-absolute-solid} which do not involve any condensed mathematics. We refer the reader to \cref{defn:smooth-proper} for the definition of a smooth proper morphism (of equidimension $d$) between diamonds over $\Spd(\QQ_p, \ZZ_p)$. 

\begin{construction}[{cf.\ \cite[Construction 6.4.2]{LRZ24}}]
\label{construction:evaluation-coevaluation-poincare-duality-non-solid}
Let $S$ be a diamond over $\Spd(\QQ_p, \ZZ_p)$,
let $X$ be a diamond,
let $f\colon X \to S$ be a smooth proper map of equidimension $d$, let $\cA\in \{\BB_I, \cO\}$, and let $\cE\in \Perf_v(X; \cA)$. We set $\cE^{\vee} \coloneqq \rR{\cHom}(\cE, \cA)$. We denote by $ev_{\cE} \colon \cE^{\vee} \otimes^\rL \cE \to \cA$ the natural evaluation map. 

Consider the following commutative diagram:
\begin{equation}
\begin{tikzcd}
X \arrow[rd,"\Delta"]\arrow[rdd, swap, bend right, "\rm{id}"] \arrow[rrd, bend left, "\rm{id}"]& & \\
& X\times_S X \arrow{r}{\pi_2} \arrow{d}{\pi_1} \arrow{rd}{h} & X \arrow{d}{f} \\
& X \arrow[r, swap, "f"] & S
\end{tikzcd} 
\end{equation}

\begin{enumerate}[leftmargin=*,label=\upshape{(\arabic*)}]
\item\label{evaluation map-nonsolid} (\emph{Evaluation map} I) 
We define the \emph{evaluation map} 
\[
e(f, S, \cE) \colon 
\rR\Gamma_v(X, \cE^{\vee}(d)[2d])
\otimes^{\rL}_{\Gamma_v(S, \cA)} \rR\Gamma_v(X, \cE)
\to \Gamma_v(S, \cA)
\]
as the composition
\[
\begin{tikzcd}
\rR\Gamma_v(X, \cE^{\vee}(d)[2d]) 
\otimes^{\rL}_{\Gamma_v(S, \cA)} \rR\Gamma_v(X, \cE) \arrow[d, dashed] \arrow{rr}{\cup} & &
\rR\Gamma_v(X, \cE^{\vee}(d)[2d] \otimes^{\rL} \cE)
\arrow{d}{\rR\Gamma_v(ev_{\cE}(d)[2d])} \\
\Gamma_v(S, \cA) & \arrow{l}{\rR\Gamma_v(\tr^{\cA}_f)} \rR\Gamma_v(S, \rR f_{v, *} \cA(d)[2d]) & \arrow{l}{\sim}\rR\Gamma_v(X, \cA(d)[2d])
\end{tikzcd}
\]
where $\tr^{\cA}_f$ is the trace map from \Cref{trace-maps}. 
\item\label{evaluation map-sheafy} (\emph{Evaluation map} II) 
We define the \emph{sheafy evaluation map} 
\[
e(f, \cE) \colon 
\rR f_{v, *} (\cE^{\vee}(d)[2d])
\otimes^{\rL}_{\cA_S} \rR f_{v, *} (\cE) 
\to \cA_S
\]
as the composition
\[
\begin{tikzcd}
\rR f_{v, *} (\cE^{\vee}(d)[2d])
\otimes^{\rL}_{\cA_S} \rR f_{v, *} (\cE)  \arrow[d, dashed] \arrow{rr}{\cup} & &
\rR f_{v, *} (\cE^{\vee}(d)[2d]\otimes^\rL \cE) 
\arrow{d}{\rR f_{v, *}(ev_{\cE}(d)[2d])} \\
\cA &  &\arrow{ll}{\tr^{\cA}_f}   \rR f_{v, *}\cA(d)[2d]
\end{tikzcd}
\]
where $\tr^{\cA}_f$ is the trace map from \Cref{trace-maps}. 

\item\label{duality-nonsolid} (\emph{Duality map} I)
We define the \emph{duality map} 
\[
\PD^\cA_f(\cE) \colon \rR \Gamma_v(X, \cE^{\vee}(d)[2d]) \to \rR\Hom_{\cA(S)}(\rR \Gamma_v(X, \cE), \cA(S) \bigr)
\]
as the map adjoint to $e(f, S, \cal{E})$ from \cref{evaluation map-nonsolid} under the tensor-hom adjunction.

\item\label{duality-sheafy} (\emph{Duality map} II)
We define the \emph{sheafy duality map} 
\[
\PD^\cA_f(\cE) \colon \rR f_{v, *}(\cE^{\vee}(d)[2d]) \to \rR\cHom_{\cD(S_v;\cA_S)}(\rR f_{v, *} \cE, \cA_S \bigr)
\]
as the map adjoint to $e(f, \cal{E})$ from \cref{evaluation map-sheafy} under the tensor-hom adjunction.
\end{enumerate}
\end{construction}

For the next corollary, we will need to use the condensification functor from \cite[Def.\,5.8 and Th.\,5.9]{Andreychev}. Namely, let $(A, A^+)$ be a complete analytic Huber pair. Then \emph{loc.\,cit.} constructs a fully faithful ``condesification'' functor
\[
\ud{(\blank)}\colon \cD(A) \hookrightarrow \cD_\sol\big((A, A^+)\bigr)
\]
which is exact, commutes with filtered colimits, and satisfies $\ud{M}(*)\simeq M$ for any $M\in \cD(A)$. Furthermore, this functor is symmetric monoidal: since $\ud{M\otimes^{\rL}N}$ lies in $\cD_\sol\big((A, A^+)\bigr)$ by construction, there is a canonical map $\ud{M}\otimes^{\rL, \sol} \ud{N} \to \ud{M\otimes^{\rL}N}$ that we need to show to be an isomorphism. Since all functors involved are exact and commute with colimits, we reduce the question to the case when $M=A[d]$ and $N=A[d']$ where the claim is evident. 

\begin{cor}\label{cor:duality-non-solid} In the notation of \cref{thm:PD-absolute-solid}, the complex $\rR\Gamma_v(X^\diam, \cE)$ lies in $\Perf(\cA(S))$ and the duality morphism
\[
\PD_f^\cA(\cE) \colon \rR\Gamma_v(X, \cE^{\vee}(d)[2d]) \to \rR\Hom\bigl(\rR\Gamma_v(X, \cE), \cA(S)\bigr)
\]
is an isomorphism. 
\end{cor}
\begin{proof}
    First, \cite[Th.\,5.50 and Cor.\,5.51.1]{Andreychev} imply that the fully faithful ``condensification'' functor $\ud{(\blank)}\colon \cD(\cA(S)) \to \cD_\sol(\ud{\cA(S)})$
    induces an equivalence between the subcategory of perfect $\cA(S)$-complexes and the subcategory of dualizable objects in $\cD_\sol(\ud{\cA(S)})$.
    Thus, \cref{thm:PD-absolute-solid} implies that there exists $P\in \Perf(\cA(S))$ such that $\ud{P} \simeq \rR\uGamma_v(X, \cE)$. Since $P\simeq \ud{P}(*)$, we conclude that 
    \[
    P \simeq \ud{P}(*) \simeq \rR\uGamma_v(X, \cE)(*) \simeq \rR\Gamma_v(X, \cE).
    \]
    Hence, we see that $P\simeq \rR\Gamma_v(X, \cE)\in \Perf(\cA(S))$. Similarly, we conclude that $\rR\Gamma_v(X, \cE^{\vee}(d)[2d])\in \Perf(\cA(S))$. Since the condensification functor is symmetric monoidal, we see that $e^\sol(f, S, \cE)$ is equal to the condensification of $e(f, S, \cE)$. Thus, the fact that the condensification functor is fully faithful and \cref{thm:PD-absolute-solid} imply that $\PD$ must be an isomorphism.
\end{proof}

\begin{corollary}
\label{cor:proper-smooth-pushforward-is-perfect}
Let $Y$ be a diamond over $\Spd(\QQ_p, \ZZ_p)$, 
let $f\colon X \to Y$ be a smooth proper map, 
let $\cA\in \{\BB_I, \cO\}$, and let $\cE\in \Perf_v(X; \cA)$. Then $\rR f_{v, *}\cE \in \Perf_v(Y; \cA)$.
Furthermore, for any affinoid perfectoid space $S$ with a map $S \to Y$, the complex $\rR\Gamma_v(X_S, \restr{\cE}{X_S})$ is perfect and the natural morphism
\[
\widetilde{\rR\Gamma_v(X_S, \restr{\cE}{X_S})} \to \restr{(\rR f_{v, *}\cE)}{S} 
\]
is an isomorphism. 
Consequently, for any inclusion  $J \subset I \subset (0, \infty)$ of closed intervals
with rational endpoints and any object $\cE \in \Perf_v(X, \BB_I)$, the natural maps
\[
\rR f_{v, *}\cE \otimes^{\rL}_{\BB_I, \res}
\BB_J \rightarrow
\rR f_{v, *}(\cE \otimes^{\rL}_{\BB_I, \res} \BB_J)
\text{ and }
\rR f_{v, *}\cE \otimes^{\rL}_{\BB_I, \varphi}
\BB_{I/p} \rightarrow
\rR f_{v, *}(\cE \otimes^{\rL}_{\BB_I, \varphi} \BB_{I/p})
\]
are isomorphisms.
\end{corollary}
\begin{proof}
    First, we show that $\rR f_{v, *}\cE$ is a perfect complex. We can check this $v$-locally on $Y$. Therefore, we may and do assume that $Y=S$ is a strictly totally disconnected 
    perfectoid space. Then $f$ comes as a diamondification of a smooth proper morphism $f^\sharp \colon X^\sharp \to S^\sharp$. Furthermore, $S^\sharp$ and $X^\sharp$ are quasi-compact in this case. Therefore, $X^\sharp$ can be decomposed into a disjoint union of clopen subsets $X_i^\sharp$  such that each $X_i^\sharp \to S^\sharp$ is smooth of equidimension $i$. Thus, in order to prove both claims of this corollary, it suffices to further assume $X^\sharp=X_d^\sharp$ is smooth proper of equidimension $d$.
    In this case, \cref{cor:duality-non-solid} ensures that $\rR\Gamma_v(X_S, \restr{\cE}{X_S})$ is a perfect complex. Thus, it suffices to show that the natural morphism
    \[
    \widetilde{\rR\Gamma_v(X_S, \restr{\cE}{X_S})} \to \restr{(\rR f_{v, *}\cE)}{S} 
    \]
    is an isomorphism. After unravelling the definition, the question boils down to showing that the natural morphism
    \[
    \rR\Gamma_v(X_S, \restr{\cE}{X_S}) \otimes^{\rL}_{\cA(S)} \cA(S') \to \rR\Gamma_v(X_{S'}, \restr{\cE}{X_{S'}})
    \]
    is an isomorphism for any morphism $S' \to S$ of strictly totally disconnected perfectoid spaces.
    This follows immediately from \cref{lemma:properties-condensed-global-sections}, 
    \cref{solid base change formula}, and the fact that $\rR\Gamma_v(X_S, \restr{\cE}{X_S})$ 
    is perfect (so its solid tensor product coincides with the usual one). 

    Now we wish to show that $\rR\Gamma_v(X_S, \restr{\cE}{X_S})$ is perfect and the natural morphism $\widetilde{\rR\Gamma_v(X_S, \restr{\cE}{X_S})} \to \restr{(\rR f_{v, *}\cE)}{S}$ for any affinoid perfectoid $S$. This follows immediately from \cref{cor:perfect-O-modules-on-affinoid-perfectoid} (resp.~\cref{cor:perfect-O-modules-on-affinoid-perfectoid-BBI}) and the proven above fact that $(\rR f_{v, *}\cE)$ is a perfect complex.  
    
    The last statement now follows from the statements already established so far, together with
    \Cref{solid pullback compatibility}.
\end{proof}

\begin{corollary}
\label{cor:HT realization for proper smooth pushforward}
Let $f\colon X \to Y$ be a smooth proper morphism of diamonds over $\Spd(\QQ_p, \ZZ_p)$, let $I\subset (0, \infty)$ be a closed interval with rational endpoints such that $1\in I$, and let $\cE \in \Perf_v(X; \BB_{I})$. Then the natural map
\[
\rR f_{v, *}\cE \otimes^{\rL}_{\BB_{I}, \theta_I}
\cO \rightarrow
\rR f_{v, *}(\cE \otimes^{\rL}_{\BB_{I}, \theta_I} \cO)
\]
is an isomorphism. 
Furthermore, if $Y$ is an affinoid perfectoid space, then the natural map 
\[
\rR\Gamma_v(X, \cE) \otimes^{\rL}_{\Gamma_v(S, \BB_{I}), \theta_I}
\Gamma_v(S, \O) \rightarrow
\rR\Gamma_v(X, \cE \otimes^{\rL}_{\BB_{I}, \theta_I} \cO)
\]
is an isomorphism. 
\end{corollary}

We refer to \cref{construction:theta-sheaves} for the definition of the homomorphism $\theta_I$. 

\begin{proof}
The first statement is $v$-local on $Y$. Thus, we may assume that $Y=S$ is an affinoid perfectoid space in the proof of both claims.
By \Cref{cor:proper-smooth-pushforward-is-perfect}, it suffices to show that the natural
map 
\[
\rR\Gamma_v(X, \cE) \otimes^{\rL}_{\Gamma_v(S, \BB_{I}), \theta_I}
\Gamma_v(S, \O) \rightarrow
\rR\Gamma_v(X, \cE \otimes^{\rL}_{\BB_{I}, \theta_I} \cO)
\]
is an isomorphism. This in turn follows from \Cref{compatibility with HT realization} and
the fact that $\rR\Gamma_v$'s involved in the formula above are
perfect complexes (again by \Cref{cor:proper-smooth-pushforward-is-perfect}).
\end{proof}

\begin{corollary}
\label{cor:PD-for-proper-smooth-pushforward}
Let $Y$ be a diamond over $\Spd(\QQ_p, \ZZ_p)$, let $f\colon X \to Y$ be a smooth proper map of equidimension $d$, 
let $\cA\in \{\BB_I, \cO\}$, let $\cE\in \Perf_v(X; \cA)$, and set $\cE^{\vee}\coloneqq \rR {\cHom}(\cE, \cA)$.
Then the following hold:
\begin{enumerate}
    \item\label{cor:PD-for-proper-smooth-pushforward-1} the sheafy duality morphism 
\[
\PD^\cA_f(\cE) \colon \rR f_{v, *} \cE^{\vee}(d)[2d] \to \rR{\cHom}\bigl(\rR f_{v, *} \cE, \cA \bigr)
\]
is an isomorphism;
    \item\label{cor:PD-for-proper-smooth-pushforward-2} if $Y=S$ is an affinoid perfectoid space, then the duality morphism
\[
\PD^\cA_f(\cE) \colon \rR \Gamma_v(X, \cE^{\vee}(d)[2d]) \to \rR\Hom_{\cA(S)}\bigl(\rR \Gamma_v(X, \cE), \cA(S) \bigr)
\]
is an isomorphism as well.
\end{enumerate}
\end{corollary}
\begin{proof}
    Part~\cref{cor:PD-for-proper-smooth-pushforward-1} is $v$-local on $Y$, so we can assume that 
    $Y = S$ is a strictly totally disconnected perfectoid space. 
    Then the result follows immediately from \cref{cor:proper-smooth-pushforward-is-perfect} 
    and \cref{cor:duality-non-solid}. Part~\cref{cor:PD-for-proper-smooth-pushforward-2} follows from Part~\cref{cor:PD-for-proper-smooth-pushforward-1}, \cref{cor:proper-smooth-pushforward-is-perfect}, and \cref{cor:perfect-O-modules-on-affinoid-perfectoid-BBI}. 
\end{proof}

\begin{corollary}
\label{cor:KF-for-proper-smooth-pushforward}
Let $Y$ be a diamond over $\Spd(\QQ_p, \ZZ_p)$, let $f\colon X \to Y$ and $f'\colon X'\to Y$ be smooth proper maps between diamonds, let $h\colon X\times_Y X' \to Y$ be the fiber product morphism, 
let $\cA\in \{\BB_I, \cO\}$, let $\cE\in \Perf_v(X; \cA)$, and let $\cE'\in \Perf_v(X'; \cA)$.
Then the following hold:
\begin{enumerate}
    \item\label{cor:KF-for-proper-smooth-pushforward-1} the K\"unneth morphism
    \[
    \rR f_{v, *} \cE  \otimes^\rL_{\cA} \rR f'_{v, *} \cE' \to \rR h_{v, *} (\cE\boxtimes^\rL \cE') 
    \]
    is an isomorphism;
    \item\label{cor:KF-for-proper-smooth-pushforward-2} if $Y=S$ is an affinoid perfectoid space, then the K\"unneth morphism
    \[
    \rR \Gamma_v(X, \cE) \otimes^\rL_{\cA(S)} \rR \Gamma_v(X', \cE') \to \rR \Gamma_v(X\times_S X', \cE \boxtimes^\rL \cE') 
    \]
    is an isomorphism as well. 
\end{enumerate}
\end{corollary}
\begin{proof}
    Part~\cref{cor:KF-for-proper-smooth-pushforward-1} is $v$-local on $Y$, so we can assume that 
    $Y = S$ is a strictly totally disconnected perfectoid space. Then the claim follows immediately from \cref{cor:proper-smooth-pushforward-is-perfect}, \cref{Kunneth formula: main thm}, and the fact that (derived) solid tensor product of perfect complexes coincides with the usual (derived) tensor product. Part~\cref{cor:KF-for-proper-smooth-pushforward-2} follows from Part~\cref{cor:KF-for-proper-smooth-pushforward-1}, \cref{cor:proper-smooth-pushforward-is-perfect}, and \cref{cor:perfect-O-modules-on-affinoid-perfectoid-BBI}. 
\end{proof}

\subsection{Cohomological dimension}
The main goal of this subsection is to show that a smooth proper morphism of equidimension $d$ has cohomological dimension $2d$, at least when restricted to perfect $\BB_I$- and $\cO$-complexes. 

For this, we will need some preliminary lemmas. 

\begin{lemma}\label{lemma:tor-amplitude-over-points} Let $S$ be an affinoid perfectoid space with a morphism $S \to \Spd(\QQ_p, \ZZ_p)$ corresponding to an untilt $S^\sharp \to \Spa(\QQ_p, \ZZ_p)$, let $\cA\in \{\BB_I, \O\}$, and let $P\in \Perf\bigl(\cA(S)\bigr)$ be a perfect complex. Fix $r,r'\in \ZZ \cup \{-\infty, \infty\}$.
Then $P$ lies in $\Perf^{[r,r']}\bigl(\cA(S)\bigr)$ if and only if $P\otimes^\rL_{\cA(S)} \cA\bigl(\Spa(\wdh{k(s)}, \wdh{k(s)}^+)\bigr)$ lies in $\Perf^{[r,r']}\bigl(\cA(\Spa(\wdh{k(s)}, \wdh{k(s)}^+))\bigr)$ for every $s\in S$.    
\end{lemma}
    We recall that, for a ring $R$, we use $\Perf^{[r,r']}(R)$ to denote the full $\infty$-subcategory of $\cD(R)$ consisting of perfect complexes with tor-amplitude in $[r,r']$. 
\begin{proof}
    To simplify the exposition, we will use the notation $Y_{S, I}$ for $Y_{S^\sharp, I}$ as defined in \cref{defn:affinoid-Fargues--Fontaine}. This notation is consistent with \cite{Fargues-Scholze}; see \cref{rmk:inconsistency-with-FS}.
    
    We give the proof for $\cA=\BB_I$ as the other case is similar (in fact, easier). For brevity, we write $\BB_I\bigl(\wdh{k(s)}, \wdh{k(s)}^+\bigr)$ instead of $\BB_I\bigl(\Spa(\wdh{k(s)}, \wdh{k(s)}^+)\bigr)$. The forward direction of the lemma is obvious. So we only need to show that $P$ has tor-amplitude in $[r, r']$ if $P\otimes^\rL_{\BB_I(S)} \BB_I\bigl(\wdh{k(s)}, \wdh{k(s)}^+\bigr)$ has tor-amplitude in $[r, r']$ for every $s\in S$. 
    
    Thanks to \cite[\href{https://stacks.math.columbia.e$du/tag/068V}{Tag 068V}]{stacks-project}, it suffices to prove that, for any maximal ideal $\m\subset \BB_I(S)$, the tensor product $P\otimes^{\rL}_{\BB_I(S)} \bigl(\BB_I(S)/\m\bigr)$ lies in $D^{[r, r']}\bigl(\BB_I(S)/\m\bigr)$. Then \cite[Lem.~1.4]{Huber-generalization} implies that we can find a point $x_\m \in Y_{S, I}=\Spa\bigl(\BB_I(S), \AA_I(S)\bigr)$ such that $\supp(x_\m)=\m$. Consider the continuous map $\abs{Y_{S, I}} \to \abs{S}$ from (the sentence after) \cite[Prop.~II.1.2]{Fargues-Scholze} and let $s_\m$ be the image of $x_\m$ under this map. Then \cite[Prop.~II.1.3]{Fargues-Scholze} implies that $x_\m$ lies in the image of the pro-(open immersion) 
    \[
    i_\m\colon \Spa\bigl(\BB_I\bigl(\wdh{k(s_\m)}, \wdh{k(s_\m)}^+\bigr), \AA_I\bigl(\wdh{k(s_\m)}, \wdh{k(s_\m)}^+\bigr)\bigr) = Y_{s_{\m}, I} \hookrightarrow Y_{S, I} = \Spa\bigl(\BB_I(S), \AA_I(S)\bigr).
    \]
    Pick $y_{\m}\in Y_{s_\m, I}$ such that $i_\m(y_\m)=x_\m$ and set $\n \coloneqq \supp(y_\m) \in \Spec \BB_I\bigl(\widehat{k(s_\m)}, \widehat{k(s_\m)}^+\bigr)$. Then we have the following commutative diagram:
    \[
    \begin{tikzcd}
        \BB_I(S) \arrow[r] \arrow[d] & \BB_I(S)/\m \arrow[d] \\
        \BB_I\bigl(\wdh{k(s_\m)}, \wdh{k(s_\m)}^+\bigr) \arrow[r] & \BB_I\bigl(\wdh{k(s_\m)}, \wdh{k(s_\m)}^+\bigr)/\n 
    \end{tikzcd}
    \]
    Now our assumption guarantees that $P \otimes^\rL_{\BB_I(S)} \BB_I\bigl(\wdh{k(s_\m)}, \wdh{k(s_\m)}^+\bigr)$ lies in $\Perf^{[r, r']}\Bigl(\BB_I\bigl(\wdh{k(s_\m)}, \wdh{k(s_\m)}^+\bigr)\Bigr)$. This formally implies that 
    \[
    P\otimes^\rL_{\BB_I(S)} \BB_I(S)/\m \otimes^\rL_{\BB_I(S)/\m} \BB_I\bigl(\wdh{k(s_\m)}, \wdh{k(s_\m)}^+\bigr)/\n \simeq P \otimes^\rL_{\BB_I(S)} \BB_I\bigl(\wdh{k(s_\m)}, \wdh{k(s_\m)}^+\bigr) \otimes^\rL_{\BB_I\bigl(\wdh{k(s_\m)}, \wdh{k(s_\m)}^+\bigr) }\BB_I\bigl(\wdh{k(s_\m)}, \wdh{k(s_\m)}^+\bigr)/\n
    \]
    lies in $\cD^{[r, r']}\Bigl( \BB_I\bigl(\wdh{k(s_\m)}, \wdh{k(s_\m)}^+\bigr)/\n\Bigr)$. Finally, since $\BB_I(S)/\m$ is a field, the morphism $\BB_I(S)/\m \to \BB_I\bigl(\wdh{k(s_\m)}, \wdh{k(s_\m)}^+\bigr)/\n$ is faithfully flat. Thus, we conclude that $P\otimes^\rL_{\BB_I(S)} \BB_I(S)/\m$ lies in $\cD^{[r, r']}\bigl( \BB_I(S)/\m\bigr)$, as desired. 
\end{proof}

We also need the following general lemma. 

\begin{lemma}\label{lemma:tor-dimension-tensor-product} Let $(C, C^+)$ be a pair consisting of an algebraically closed nonarchimedean field $C$ and an open bounded valuation subring $C^+\subset C$, let $S\coloneqq \Spd(C, C^+) \to \Spd(\QQ_p, \ZZ_p)$ be a morphism, and let $S' \to S$ be a morphism of affinoid perfectoid spaces. Let $\cA\in \{\O, \BB_I\}$ and let $a\in \cA(S)$ be a non-zero element. Then the image of $a$ in $\cA(S')$ is a nonzerodivisor. 
\end{lemma}
\begin{proof}
    We give a proof for $\cA=\BB_I$ as the other case is similar (in fact, easier). First, \cref{lemma:BI-overconvergent} ensures that we can replace $S$ with $S^\circ\coloneqq \Spd(C, C^\circ)$ and $S'$ with $(S')^\circ \coloneqq S'\times_{S} S^\circ$ for the purposes of proving this lemma. Therefore, we can assume that $S=\Spd(C, C^\circ)$. Moreover, we can assume that $S' \neq \varnothing$, as $0$ is a nonzerodivisor in the $0$-ring. Then \cref{lemma:BI-determined-by-points} and \cref{lemma:BI-overconvergent} ensure that we can assume that $S'=\Spa(K', K'^\circ)$ for some perfectoid field $K'$ with a map $C \to K'$. Let $(C', C'^\circ)$ be the completed algebraic closure of $(K', K'^\circ)$. Then \cref{cor:AI-etale-sheaf} ensures that the natural map $\BB_I\bigl(\Spa(K', K'^\circ)\bigr) \to \BB_I\bigl(\Spa(C', C'^\circ)\bigr)$ is injective since $\Spa(C', C'^\circ) \to \Spa(K', K'^\circ)$ is a $v$-covering. So we can assume that $S'=\Spa(C', C'^\circ)$. Since $\Spa(C', C'^\circ) \to \Spa(C, C^\circ)$ is a $v$-covering, \cref{cor:BI-PID} and \cref{cor:AI-etale-sheaf} imply that $\BB_I(S) \hookrightarrow \BB_I(S')$ is an injective map of domains. Therefore, the element $a$ remains a nonzerodivisor in $\BB_I(S')$.
\end{proof}
Before proving the main result of this section, we recall that we defined in \cref{notation:sheaf-perfect-complexes-finite-amplitude} the $\infty$-category $\Perf^{[r,r']}_v(X; \cA)$ for $r,r'\in \ZZ\cup \{-\infty, \infty\}$ and $\cA\in \{\O, \BB_I\}$.
By \cref{lemma:sheafification-of-fully-faithful}, this is a full $\infty$-subcategory of $\Perf^{[-\infty, \infty]}_v(X; \cA) \simeq \Perf_v(X; \cA)$, and thus embeds naturally as a full $\infty$-subcategory of $\cD(X_v; \cA)$. Furthermore, \cref{thm:v-descent-rationally-BBI} (resp.~\cref{thm:v-descent-rationally}) ensures that when $X$ is an affinoid perfectoid space, we have an equivalence $\Perf^{[r, r']}_v(X; \cA) \simeq \Perf^{[r, r']}\bigl(\cA(X)\bigr)$, where the latter denotes perfect complexes with tor-amplitude in $[r, r']$. Consequently, for an arbitrary diamond $X$, the $\infty$-category $\Perf^{[r, r']}_v(X; \cA)$ naturally forms a full $\infty$-subcategory of $\cD^{[r,r']}(X_v ; \cA)$, where $\cD^{[r,r']}(X_v ; \cA)$ denotes the $\infty$-subcategory consisting of objects concentrated in degrees $[r, r']$. 

Now we are ready to prove the main result of this subsection. 
\begin{corollary}
\label{cor:PD-tor-amplitude}
Let $Y$ be a diamond over $\Spd(\QQ_p, \ZZ_p)$, let $f\colon X \to Y$ be a smooth proper map of dimension $d$, let $\cA\in \{\BB_I, \cO\}$, and let $\cE\in \Perf^{[r,r']}_v(X; \cA)$ for some $r, r'\in \ZZ\cup \{-\infty, \infty\}$. 
Then the following hold:
\begin{enumerate}
    \item\label{cor:PD-tor-amplitude-1} the derived pushforward $\rR f_{v, *}\cE$ lies in $\Perf^{[r, r'+2d]}_v(Y; \cA)$; 
    \item\label{cor:PD-tor-amplitude-2} if $Y=S$ is an affinoid perfectoid space, then $\rR\Gamma_v(X, \cE)$ lies in $\Perf^{[r, r'+2d]}\bigl(\cA(Y)\bigr)$. 
\end{enumerate}
\end{corollary}
\begin{proof}
    \emph{Step~$1$. Reduce to the case of Part~\cref{cor:PD-tor-amplitude-2} with  $Y=\Spd(C, C^+)$ for an algebraically closed nonarchimedean field $C$ and an open bounded valuation subring $C^+\subset C$ and $r=0$, $r'\in \ZZ$.} First, we note that \cref{notation:sheaf-perfect-complexes-finite-amplitude}, \cref{cor:proper-smooth-pushforward-is-perfect}, and \cref{thm:v-descent-rationally-BBI} (resp.~\cref{thm:v-descent-rationally}) imply that it suffices to to show that $\rR\Gamma_v(X, \cE) \in \Perf^{[r, r'+2d]}(\cA(Y))$ under the additional assumption that $Y$ is a strictly totally disconnected perfectoid space.

    In this case, we can use \cref{cor:proper-smooth-pushforward-is-perfect} once again together with \cref{lemma:tor-amplitude-over-points} to reduce to the case when $Y=\Spa(C, C^+)$ for an algebraically closed nonarchimedean field $C$ and an open bounded valuation subring $C^+\subset C$. Now the observation that $X$ is qcqs and \cref{lemma:colimit-of-sheaves} imply that $\colim_{\NN} \Perf^{[-n, n]}_v(X; \cA) \to \Perf_v^{[-\infty, \infty]}(X; \cA) \simeq \Perf_v(X; \cA)$ is an isomorphism. Thus, we can assume that $\cE\in \Perf^{[r, r']}_v(X; \cA)$ for some finite $r, r'\in \ZZ$. After an appropriate cohomological shift, we can moreover take $r=0$. 

    \emph{Step~$2$. $\rR\Gamma_v(X, \cE) \in \Perf^{[0, k]}(\cA(Y))$ for some $k\in \ZZ$.} We already know that $\rR\Gamma_v(X, \cE)\in \Perf(\cA(Y))$ by virtue of \cref{cor:proper-smooth-pushforward-is-perfect}. If $\cA=\O$, we note that $\cA(Y)=C$ is a field. So since $\cE\in \Perf^{[0,r']}_v(X; \cA)\subset \cD^{[0, r']}(X; \cA)$, it is therefore clear that $\rR\Gamma_v(X, \cE) \in \Perf^{[0, k]}(\cA(Y)) = \cD^{[0, k]}_{\mathrm{coh}}(\cA(Y))$ when $\cA=\O$. 
    
    If $\cA=\BB_I$, then \cref{cor:BI-PID} guarantees that $\cA(Y)$ is a principal ideal domain. Thus, in this case, it suffices to show that $\rR\Gamma_v(X, \cE) \otimes^\rL_{\BB_I(Y)} \BB_I(Y)/(a)\in \cD^{\geq 0}\bigl(\BB_I(Y)\bigr)$ for any non-zero $a\in \BB_I(Y)$. Now \cref{PF} guarantees that 
    \[
    \rR\Gamma_v(X; \cE) \otimes^\rL_{\BB_I(Y)} \BB_I(Y)/(a) \simeq \rR\Gamma_v(X; \cE \otimes^\rL_{\BB_I} \BB_I/^\rL a\BB_I). 
    \]
    Hence, it suffices to prove that $\cE \otimes^\rL_{\BB_I} \BB_I/^\rL a\BB_I$ is in $\cD^{\geq 0}(X_v; \cA)$. This can be checked $v$-locally on $X$. So, until the end of this step, we drop the assumption that $X$ is smooth proper over $Y$ and instead assume that $X$ is an affinoid perfectoid space and $\cE = \widetilde{P}$ for $P\in \Perf^{[0, r']}\bigl(\BB_I(X)\bigr)$ (see \cref{lemma:locally-surjective} and \cref{cor:perfect-O-modules-on-affinoid-perfectoid-BBI}). Then the question boils down to showing that 
    \[
    P\otimes^\rL_{\BB_I(X)} \BB_I(X)/^\rL a\BB_I(X) \text{ lies in }\cD^{\geq 0}\bigl(\BB_I(X)\bigr).
    \]
    This follows immediately from \cref{lemma:tor-dimension-tensor-product} and the definition of $\Perf^{[0, r']}\bigl(\BB_I(X)\bigr)$. 

    \emph{Step~$3$. $\rR\Gamma_v(X, \cE) \in \Perf^{[0, r'+2d]}(\cA(Y))$.} Assuming again that $X$ is smooth proper over $Y=\Spd(C, C^+)$ of dimension $d$ (so it is the diamondification of a smooth proper adic space over $\Spa(C^\sharp, C^{\sharp, +})$ of relative dimension $d$), there is a finite clopen decomposition $X= \sqcup_{0\leq n\leq d} X_n$ such that $X_n \to Y=\Spd(C, C^+)$ is smooth proper of equidimension $n$. Therefore, we can replace $X$ with $X_n$ to assume that $X \to Y$ is of equidimension $n$ for some integer $n\leq d$.
    
    In this case, \cref{cor:duality-non-solid} ensures that $\rR\Gamma_v(X, \cE) \simeq \rR\Gamma_v(X, \cE^{\vee}(n)[2n])^{\vee}$. Since $\cE\in \Perf^{[0, r']}_v(X; \cA)$, we conclude that $\cE^{\vee}(n)[2n]\in \Perf^{[-r'-2n, -2n]}_v(X; \cA)$. Hence, Step~$2$ applied to $\cE^{\vee}(n)[2n]$ (shifted by $r'+2n$) implies that $\rR\Gamma_v(X; \cE^{\vee}(n)[2n]) \in \Perf^{[-r'-2n, t]}(\cA(Y))$ for some $t\in \ZZ$. From this, we deduce 
    \[
    \rR\Gamma_v(X, \cE) \simeq \rR\Gamma_v(X, \cE^{\vee}(n)[2n])^{\vee} \in \Perf^{[-t, r'+2n]}(\cA(Y)).
    \]
    Invoking Step~$2$ once more, we conclude that $\rR\Gamma_v(X, \cE) \in \Perf^{[0, r'+2n]}\bigl(\cA(Y)\bigr)$. Finally, the result follows from the fact that $n\leq d$.
\end{proof}

As a formal consequence of \cref{cor:PD-tor-amplitude}, we get a similar result for cohomological amplitude in place of the tor-amplitude. 

\begin{corollary}
\label{cor:tor-amplitude-2}
Let $Y$ be a diamond over $\Spd(\QQ_p, \ZZ_p)$, let $f\colon X \to Y$ be a smooth proper map of dimension $d$, 
let $\cA\in \{\BB_I, \cO\}$, and let $\cE\in \Perf_v(X; \cA) \cap \cD^{[r, r']}(X_v; \cA)$ for some $r, r'\in \ZZ\cup \{-\infty, \infty\}$. Then $\rR f_{v, *} \cE$ lies in $\Perf_v(Y; \cA) \cap \cD^{[r, r'+2d]}(Y_v; \cA)$. 
\end{corollary}
\begin{proof}
    We only need to show that $\rR f_{v, *} \cE$ lies in $\cD^{[r, r'+2d]}(Y_v; \cA)$ as we already know that it is a perfect complex (\cref{cor:proper-smooth-pushforward-is-perfect}).
    Clearly, this complex lies in $\cD^{\geq r}(Y_v; \cA)$, so we only need to show that it lies in $\cD^{\leq r'+2d}(Y_v; \cA)$. Now we note that $\Perf_v(X; \cA) \cap \cD^{\leq r'}(X_v; \cA)\subset \Perf_v^{[-\infty, r']}(X; \cA)$. Thus, \cref{cor:PD-tor-amplitude} implies that $\rR f_{v, *} \cE$ lies in $\Perf^{[-\infty, r'+2d]}_v(Y; \cA) \subset \cD^{\leq r'+2d}(Y_v; \cA)$. This finishes the proof. 
\end{proof}

\section{Relative duality}\label{geom-arith-duality}

In this section, we establish the first main result of the paper: relative geometric duality for $\uQ_p$-coefficients. More precisely, we prove a version of Poincar\'e duality for smooth proper morphisms and perfect complexes of $\uQ_p$-modules. For this, we first develop a theory of (analytic) perfect complexes on the relative Fargues--Fontaine curve and establish a duality theorem for such complexes with respect to smooth proper morphisms of diamonds. A crucial input in this step is the duality for $\BB_I$-coefficients established in \cref{cor:PD-for-proper-smooth-pushforward}. We then use the Riemann--Hilbert and solution functors to deduce the corresponding duality statement for $\uQ_p$-coefficients. Finally, we show that a smooth proper morphism of dimension $d$ has cohomological dimension $2d$ on the subcategory of perfect complexes of $\uQ_p$-modules.

\subsection{The Fargues--Fontaine curve}\label{section:fargues-fontaine-curve}

In this subsection, we recall the definition of the relative (analytic) Fargues--Fontaine curve $\FFF_{S}$ associated to a perfectoid space $S$ over $\Spa(\QQ_p, \ZZ_p)$, following \cite[\textsection II.1]{Fargues-Scholze}. We then collect several fundamental properties of these curves. In particular, we express the $\infty$-category of analytic perfect complexes on $\FFF_S$ in terms of the $\infty$-categories of perfect complexes of $\BB_I$-modules on the $v$-site of $S$ for $I=[1,1]$ and $I=[1,p]$. This will be crucial in later subsections for defining and studying the category $\Perf_\an(\FFF_X; \O)$ for an arbitrary diamond $X$ over $\Spd(\QQ_p, \ZZ_p)$.

Throughout this subsection, we let $S$ be a fixed affinoid perfectoid space over $\Spa(\QQ_p, \ZZ_p)$. We also let $(\Adic)_{\AN}$ denote the category of analytic adic spaces endowed with the (big) analytic topology. While this forms a site in the sense of \cite[\href{https://stacks.math.columbia.edu/tag/00VH}{Tag 00VH}]{stacks-project}, we caution that $(\Adic)_{\AN}$ does not admit arbitrary finite limits.

For $S$ as above, \cref{defn:Y-Fargues--Fontaine} yields the relative Fargues--Fontaine punctured open disk $Y_S$, equipped with the Frobenius automorphism $\varphi \colon Y_S \to Y_S$. By \cite[Prop.~II.1.16]{Fargues-Scholze}, the action of $\varphi^\ZZ$ on $Y_S$ is free and totally discontinuous. Consequently, the naive quotient $\rm{FF}_S\coloneqq Y_S/\varphi^{\ZZ}$, formed in the category of topologically locally $v$-ringed spaces (see \cite[Def.~4.1.2]{Z-quotients}\footnote{Strictly speaking, \cite[Def.~4.1.2, Rmk.~4.1.3, and Lem.~2.1.2]{Z-quotients} ensure that $X/G$ is a categorical quotient in the category of topologically locally $v$-ringed spaces only when the group $G$ is finite. However, the argument in \emph{loc.~cit.} never used that $G$ is finite and, thus, it applies to any group $G$.}), naturally inherits the structure of an analytic adic space.

\begin{defn}\label{defn:FF-curve} The \emph{relative Fargues--Fontaine curve $\rm{FF}_S \coloneqq Y_S/\varphi^\ZZ$} is the quotient of $Y_S$ by the action of Frobenius. We denote by $\pi_S \colon Y_S \to \rm{FF}_S$ the projection morphism. 
\end{defn}

Since the action of $\varphi^{\ZZ}$ is free and totally discontinuous, every point $x \in \rm{FF}_S$ admits an open neighborhood $U_x$ such that the map $\pi_S^{-1}(U_x) \to U_x$ is isomorphic over $U_x$ to the projection $U_x \times \ZZ \to U_x$. This leads to the following conclusion: 

\begin{rmk}\label{rmk:FF-quotient} The projection morphism $\pi_S \colon Y_S \to \rm{FF}_S$ realizes $\rm{FF}_S$ as the $\varphi^\ZZ$-quotient of $Y_S$ formed in the category 
$\Shv\bigl((\Adic)_\AN; \Set\bigr)$. Indeed, this claim is local on $\rm{FF}_S$ in the analytic topology, so it suffices to prove an analogous claim for the projection $U\times \ZZ \to U$, which is evident. 
\end{rmk}

The primary downside of \cref{rmk:FF-quotient} is that $Y_S$ is not affinoid or even quasi-compact: this complicates the use of descent along $\pi_S$ to compute the category of perfect complexes (or vector bundles) on $\FFF_S$. To remedy this, we present $\FFF_S$ as a gluing of two affinoid subsets by identifying a suitable ``fundamental domain'' for the action of $\varphi^\ZZ$ on $Y_S$. But before we do this, we need the following general lemma. 

\begin{lemma}\label{lemma:cartesian-cocartesian-big-analytic} Let $X$ be an adic space over $\Spa(\ZZ_p, \ZZ_p)$ which is either locally noetherian or sousperfectoid. Let $j\colon U \hookrightarrow X$ be an open immersion and let $\pi \colon X' \to X$ be an \'etale map. Set $U'\coloneqq U\times_X X'$, $Z\coloneqq \abs{X}\smallsetminus \abs{U}$, and $Z'\coloneqq \abs{X'}\smallsetminus \abs{U'}$. Suppose that $\restr{\abs{\pi}}{Z'} \colon Z' \to Z$ is bijective and the induced map $\bigl(\wdh{k(\pi(x'))}, \wdh{k(\pi(x'))}^+\bigr) \to \bigl(\wdh{k(x')}, \wdh{k(x')}^+\bigr)$ is an isomorphism for any $x'\in Z'$. Then the fiber square of analytic adic spaces
    \[ \begin{tikzcd}
            U' \arrow[r,hook,"j'"] \arrow[d,"\pi'"] & X' \arrow[d,"\pi"] \\
            U \arrow[r,hook,"j"] & X
    \end{tikzcd} \]
    is a pushout square in $\Shv\bigl((\Adic)_{\AN}; \Set\bigr)$ and in $\Shv\bigl((\Adic)_{\AN}; \Ani\bigr)$.
\end{lemma}

While the assumptions on $X$ should not be necessary, we included them to streamline the proof.

\begin{proof}
    The proof of \cref{cartesian-cocartesian} applies almost verbatim in this situation. 
    We only need to address two things: $\{U \xhookrightarrow{j} X, X' \xr{\pi} X\}$ and $\{X' \xhookrightarrow{\Delta} X'\times_{X} X', U'\times_U U' \xhookrightarrow{j'\times_X j'} X'\times_X X'\}$  are coverings in the big analytic topology. 

    For the first claim, it suffices to show that, for each $x\in Z$, there is an open subspace $x\in V_x\subset X$ such that $\restr{\pi}{\pi^{-1}(V_x)} \colon \pi^{-1}(V_x)\to V_x$ admits a section.
    By assumption,     there is $x'\in X'$ such that $\pi(x')=x$. Since $\pi$ is \'etale, we know that $\pi$ is an open map. Thus it suffices to find an open subspace $x'\in V'_x\subset X'$ such that $\restr{\pi}{V'_x}\colon V'_x \to \pi(V'_x)$ is an isomorphism. This follows immediately from our assumption that $\bigl(\wdh{k(\pi(x'))}, \wdh{k(\pi(x'))}^+\bigr) \to \bigl(\wdh{k(x')}, \wdh{k(x')}^+\bigr)$ is an isomorphism and  \cite[Lem.~IV.4.14]{Fargues-Scholze} (note that $d=0$ due to \'etaleness of $\pi$).\footnote{The statement of \cite[Lem.~IV.4.14]{Fargues-Scholze} assumes that $X$ is sousperfectoid, but the proof applies to any locally noetherian adic space over $\Spa(\ZZ_p, \ZZ_p)$.} 

    For the second claim, we first observe that $X' \xhookrightarrow{\Delta} X'\times_{X} X'$ and $U' \times_U U' \xhookrightarrow{j'\times_X j'} X'\times_X X'$ are both open immersions due to the assumption that $\pi$ is \'etale and $j$ is an open immersion. Therefore, it suffices to show that $\Delta$ and $j'\times_X j'$ are jointly surjective, i.e., $\abs{X'\times_X X'\smallsetminus U'\times_U U'}$ lies inside $\abs{\Delta}(\abs{X'})$. The assumption on residue fields and the assumption that $\restr{\abs{\pi}}{Z'} \colon Z' \to Z$ is bijective imply that $\abs{X'\times_X X'\smallsetminus U'\times_U U'} \simeq Z'\times_{Z} Z' \simeq \abs{\Delta}(Z')\subset \abs{\Delta}(\abs{X'})$. This finishes the proof. 
\end{proof}

To apply \cref{lemma:cartesian-cocartesian-big-analytic} to our setting, recall from \cref{defn:affinoid-Fargues--Fontaine} that every closed interval $I\subset (0, \infty)$ with rational endpoints determines an open affinoid subspace $Y_{S, I} \subset Y_S$, which we call the Fargues--Fontaine annulus of radius $I$.

\begin{notation}\label{notation:restriction-frobenius-space} Let $S$ be an affinoid perfectoid space over $\Spa(\QQ_p, \ZZ_p)$. 
 \begin{enumerate}
   \item\label{notation:restriction-frobenius-space-1} (\emph{Inclusion map}) We denote by $\incl \colon Y_{S, [1, 1]} \hookrightarrow Y_{S, [1, p]}$ the natural inclusion of relative Fargues--Fontaine annuli. 
   \item\label{notation:restriction-frobenius-space-2} (\emph{Frobenius map}) We denote by $\widetilde{\varphi} \colon Y_{S, [1, 1]} \hookrightarrow Y_{S, [1, p]}$ the composition $Y_{S, [1, 1]} \xr[\sim]{\varphi} Y_{S, [p, p]} \hookrightarrow Y_{S, [1, p]}$ of the Frobenius isomorphism followed by the natural inclusion.
\end{enumerate}
\end{notation}

With this notation, we are finally ready to prove the following lemma: 

\begin{lemma}\label{lemma:FF-quotient-of-affinoid} Let $S$ be an affinoid perfectoid space over $\Spa(\QQ_p, \ZZ_p)$ and let $\rm{FF}_S$ be the associated Fargues--Fontaine curve. Then the following square 
    \begin{equation}\label{eqn:fundamental-domain} \begin{tikzcd}[column sep = huge]
            Y_{S, [1, 1]} \sqcup Y_{S, [1, 1]} \arrow[r,hook,"\incl\, \sqcup\, \wphi"] \arrow[d,"\id \sqcup \id"] & Y_{S, [1, p]} \arrow[d,"\restr{\pi_S}{Y_{S, [1,p]}}"] \\
            Y_{S, [1, 1]} \arrow[r,hook,"\restr{\pi_S}{Y_{S, [1,1]}}"] & \rm{FF}_S
    \end{tikzcd} 
    \end{equation}
    is a pushout square in $\Shv\bigl((\Adic)_{\AN}; \Set\bigr)$ and in $\Shv\bigl((\Adic)_{\AN}; \Ani\bigr)$.
\end{lemma}
\begin{proof}
    We only need to verify that Diagram~\cref{eqn:fundamental-domain} satisfies the assumptions of \cref{lemma:cartesian-cocartesian-big-analytic}. First, we note that $\rm{FF}_S$ is sousperfectoid because $Y_S$ is sousperfectoid (see the proof of \cite[Prop.~II.1.1]{Fargues-Scholze}) and the action of $\varphi^{\ZZ}$ is free and totally discontinuous. Furthermore, this freeness and total discontinuity ensure that the restriction $\pi_S|_{Y_{S, [1,p]}} \colon Y_{S, [1,p]}\to \rm{FF}_S$ is \'etale and induces an isomorphism on residue fields.

    To verify the remaining properties, namely, that $\pi_S|_{Y_{S, [1,1]}}$ is an open immersion, Diagram~\ref{eqn:fundamental-domain} is cartesian, and $\pi_S|_{Y_{S, [1, p]}}$ is bijective away from $Y_{S, [1,1]} \sqcup Y_{S, [p,p]} = \rm{Im}(\incl \,\sqcup\, \wphi)$, it suffices to show that if $y\in Y_{S, [1, p]}$ and $n\in \ZZ\smallsetminus \{0\}$ satisfy $\varphi^n(y)\in Y_{S, [1, p]}$, then either $y\in Y_{S, [1,1]}$ and $n=1$, or $y\in Y_{S, [p, p]}$ and $n=-1$.

    To see this, we choose $p^\flat\in \O^{+, \flat}_S(S)$ as in \cref{notation:p-flat}. Suppose $y\in Y_{S, [1,p]}$ is a point such that $\varphi^n(y)\in Y_{S, [1,p]}$ for some $n>0$. Because $\varphi([p^\flat])=[p^\flat]^p$ and $\varphi(p)=p$, the point $y$ must satisfy the following inequalities: 
    \[
    \abs{p(y)}^p \leq \abs{[p^\flat](y)} \leq \abs{p(y)} \text{ and } \abs{p(y)}^p \leq \abs{[p^\flat](y)}^{p^n} \leq \abs{p(y)}.
    \]
    Given $n>0$, these inequalities hold if and only if $n=1$ and $\abs{[p^\flat](y)}=\abs{p(y)}$, which precisely means $y\in Y_{S, [1, 1]}$. A symmetric argument shows that if $n<0$ and $\varphi^n(y)\in Y_{S, [1, p]}$, then $y\in Y_{S, [p, p]}$ and $n=-1$. This completes the proof.
\end{proof}

\begin{cor}\label{cor:FF-coeq} Let $S$ be an affinoid perfectoid space over $\Spa(\QQ_p, \ZZ_p)$ and let $\rm{FF}_S$ be the associated relative Fargues--Fontaine curve. Then the diagram $\begin{tikzcd}
\bigl(Y_{S, [1,1]} 
\arrow[r, shift left=0.8ex, "\incl"]
  \arrow[r, shift right=0.8ex, swap, "\wphi"] &
Y_{S, [1, p]}\bigr)
\end{tikzcd} \to \rm{FF}_S$ is a coequalizer diagram in $\Shv(\Adic_\AN; \Ani)$.
\end{cor}
\begin{proof}
    This follows immediately from \cref{lemma:FF-quotient-of-affinoid} and (the dual version of) \cite[\href{https://kerodon.net/tag/03HC}{Tag 03HC}]{kerodon}. 
\end{proof}

Finally, we wish to use \cref{cor:FF-coeq} to describe $\Perf_\an(\FFF_S; \O)$ in terms of perfect complexes of $\BB_I(S)$-modules for $I=[1,1]$ and $I=[1,p]$, which we can further relate to $\Perf_v(S; \BB_I)$ due to \cref{thm:v-descent-rationally-BBI}. Before we do this, we introduce the following notation: 

\begin{notation}\label{notation:restriction-frobenius} Let $S$ be an affinoid perfectoid space over $\Spa(\QQ_p, \ZZ_p)$. 
 \begin{enumerate}
   \item\label{notation:restriction-frobenius-1} (\emph{Restriction map})  We denote by 
   \[
   \res\colon \BB_{[1, p]}(S) \to \BB_{[1, 1]}(S)
   \]
   the map induced by the map $ Y_{S, [1,1]} \xhookrightarrow{\incl} Y_{S, [1,p]}$ (see \cref{notation:restriction-frobenius-space}~\cref{notation:restriction-frobenius-space-1}) and by $\res^*\colon \Perf\bigl(\BB_{[1, p]}(S)\bigr) \to \Perf\bigl(\BB_{[1, 1]}(S)\bigr)$ the induced pullback functor. 
   \item\label{notation:restriction-frobenius-2} (\emph{Frobenius map}) We denote by 
   \[
   \widetilde{\varphi}\colon \BB_{[1, p]}(S) \to \BB_{[1, 1]}(S)
   \]
   the map induced by the Frobenius morphism $Y_{S, [1,1]} \xhookrightarrow{\wphi} Y_{S, [1,p]}$, see \cref{notation:restriction-frobenius-space}~\cref{notation:restriction-frobenius-space-2}. We denote by $\wphi^*\colon \Perf\bigl(\BB_{[1, p]}(S)\bigr) \to \Perf\bigl(\BB_{[1, 1]}(S)\bigr)$ the induced pullback functor. 
\end{enumerate}
\end{notation}

\begin{notation}
\label{notation:restriction-frobenius-diamond} 
Let $X$ be a diamond over $\Spd(\QQ_p, \ZZ_p)$. 
 \begin{enumerate}
   \item (\emph{Restriction map}) Since affinoid perfectoid spaces form a basis for the $v$-topology on $X$, \cref{notation:restriction-frobenius}~\cref{notation:restriction-frobenius-1} defines a morphism of period sheaves 
   \[
   \res\colon \BB_{[1, p]} \to \BB_{[1, 1]}
   \]
   on the $v$-site of $X$. We denote by $\blank \otimes^\rL_{\BB_{[1,p]}, \res} \BB_{[1, 1]}\colon \Perf_v\bigl(X; \BB_{[1, p]}\bigr) \to \Perf_v\bigl(X; \BB_{[1, 1]}\bigr)$ the induced tensor product functor.
   \item (\emph{Frobenius map}) Similarly, \cref{notation:restriction-frobenius}~\cref{notation:restriction-frobenius-2} defines a morphism of period sheaves 
   \[
   \wphi\colon \BB_{[1, p]} \to \BB_{[1, 1]}
   \]
   on the $v$-site of $X$. We denote by $\blank \otimes^\rL_{\BB_{[1,p]}, \wphi} \BB_{[1,1]}\colon \Perf_v\bigl(X; \BB_{[1, p]}\bigr) \to \Perf_v\bigl(X; \BB_{[1, 1]}\bigr)$ the induced tensor product functor.
\end{enumerate}
\end{notation}

\begin{cor}\label{cor:Perf-on-FF} There are equivalences 
\[ \begin{tikzcd}[scale cd=.95,center picture]
\Perf_\an(\rm{FF}_\blank; \O) \arrow[r,phantom,"\simeq"] &[-1.9em] 
\eq\bigl(\Perf(\BB_{[1,p]}(\blank)) 
\arrow[r, shift left=0.8ex, "\res^*"]
  \arrow[r, shift right=0.8ex, swap, "\wphi^*"] &
\Perf(\BB_{[1, 1]}(\blank))\bigr) \simeq  \eq\bigl(\Perf_v(\blank; \BB_{[1,p]}) 
\arrow[r, shift left=0.8ex, "\otimes^\rL_{\BB_{[1, p]}, \res} \BB_{[1, 1]}"]
  \arrow[r, shift right=0.8ex, swap, "\otimes^\rL_{\BB_{[1, p]},\wphi} \BB_{[1, 1]}"] &[3.5em]
\Perf_v(\blank; \BB_{[1,1]})\bigr)
\end{tikzcd}
\] 
of functors $\Perfd_{/\QQ_p}^{\aff, \op} \to \Cat^\otimes_\infty$. 
\end{cor}
We recall that $\Cat_\infty^\otimes\simeq \CAlg(\Cat_\infty)$ denotes the $\infty$-category of symmetric monoidal $\infty$-categories. Note that \cite[Lem.~B.2.4(i)]{Heyer-Mann} ensures that the forgetful functor $\Cat_\infty^\otimes \to \Cat_\infty$ is conservative and commutes with all limits. In particular, we note that a $\Cat_\infty^\otimes$-valued presheaf is a sheaf if and only if the associated $\Cat_\infty$-valued presheaf is a sheaf. We will freely use this observation in the proof. 
\begin{proof}
    The functor $\Perf_{\an}(\blank; \O) \colon \Adic^{\op} \to \Cat_\infty^\otimes$ is an analytic sheaf of $\infty$-categories. Thus, \cite[Lem.~A.3.4]{Mann-thesis} implies that $\Perf_{\an}(\blank; \O)$ essentially uniquely extends to a limit-preserving functor 
    \[
    \Perf_{\an}(\blank; \O) \colon \Shv(\Adic_\AN; \Ani)^{\op} \to \Cat^\otimes_\infty.
    \]
    By construction, the functor $\Perf_\an(\rm{FF}_\blank; \O)$ is equivalent to the composition of functors
    \[
    \Perfd_{\QQ_p}^{\aff, \op} \xr{\rm{FF}_{\blank}} \Shv(\Adic_\AN, \Ani)^\op \xr{\Perf_\an(\blank; \O)} \Cat^\otimes_\infty.
    \]
    Since $\Perf_{\an}(\blank; \O)\colon \Shv(\Adic_\AN, \Ani)^\op \to \Cat^\otimes_\infty$ is a sheaf, it is a limit-preserving functor. Thus, \cref{cor:FF-coeq} implies that we have an equivalence of functors
    \[
    \Perf_\an(\rm{FF}_{\blank}; \O) \simeq  \begin{tikzcd}
\eq\bigl(\Perf_\an(Y_{\blank, [1, p]}) 
\arrow[r, shift left=0.8ex, "\incl^*"]
  \arrow[r, shift right=0.8ex, swap, "\wphi^*"] &
\Perf_\an(Y_{\blank, [1, 1]})\bigr).
\end{tikzcd}
    \]
    Finally, since $Y_{S, [1, p]}$ and $Y_{S, [1,1]}$ are affinoid adic spaces, \cite[Th.~1.4]{Andreychev} and \cref{thm:v-descent-rationally-BBI} imply that 
\[ \begin{tikzcd}[/tikz/baseline=(\tikzcdmatrixname-2-2.base)]
\eq\bigl(\Perf_\an(Y_{\blank, [1, p]}) 
\arrow[r, shift left=0.8ex, "\incl^*"]
  \arrow[r, shift right=0.8ex, swap, "\wphi^*"] &
\Perf_\an(Y_{\blank, [1, 1]})\bigr) \simeq \eq\bigl(\Perf(\BB_{[1,p]}(\blank)) 
\arrow[r, shift left=0.8ex, "\res^*"]
  \arrow[r, shift right=0.8ex, swap, "\wphi^*"] &
\Perf(\BB_{[1, 1]}(\blank))\bigr) \\ 
& \simeq \eq\bigl(\Perf_v(\blank; \BB_{[1,p]}) 
\arrow[r, shift left=0.8ex, "\otimes^\rL_{\BB_{[1, p]}, \res} \BB_{[1, 1]}"]
   \arrow[r, shift right=0.8ex, swap, "\otimes^\rL_{\BB_{[1, p]}, \wphi} \BB_{[1, 1]}"]  &
\Perf_v(\blank; \BB_{[1,1]})\bigr).
\end{tikzcd} \qedhere \]
\end{proof}

\subsection{Analytic perfect complexes on the Fargues--Fontaine curve}

In this subsection, we define and study the correct category $\Perf_\an(\FFF_X; \O)$ of ``analytic perfect complexes on the relative Fargues--Fontaine curve $\FFF_X$'' for any diamond $X$ over $\Spd(\QQ_p, \ZZ_p)$. Unlike the case of an affinoid perfectoid space $X$, we do not have any actual adic space $\FFF_X$, so we cannot directly define $\Perf_\an(\FFF_X; \O)$ as perfect complexes on a certain analytic site of 
$\FFF_X$. 
In fact, the relative curve $\FFF_X$ exists only as a diamond, which is not good enough to capture the category of analytic perfect complexes. Instead, we define $\Perf_\an(\FFF_X; \O)$ via the formula from \cref{cor:Perf-on-FF}. Alternatively, one could have defined $\Perf_\an(\FFF_X; \O)$ by resolving $X$ via perfectoids in the $v$-topology; we show that these two approaches give the same answer. However, we warn the reader that $\Perf_\an(\FFF_S; \O)$ is not equivalent to $\Perf_v(\FFF_S^\diam; \O)$ even when $S$ is an affinoid perfectoid space over $\Spa(\QQ_p, \ZZ_p)$ (see \cref{warning:analytic-vs-v}).

\begin{defn}
\label{defn:analytic-perfectoid-complexes-FF} 
We define the functor $\Perf_\an(\rm{FF}_\blank; \O) \colon \Diam_{/\QQ_p}^{\op} \to \CAlg(\Cat_\infty)$ as the following equalizer 
\[
\Perf_\an(\rm{FF}_\blank; \O) \coloneqq \begin{tikzcd}[column sep = 8em]\eq\Bigl(\Perf_v(\blank; \BB_{[1,p]}) 
\arrow[r, shift left=0.8ex, "\otimes^\rL_{\BB_{[1, p]}, \res} \BB_{[1, 1]}"]
   \arrow[r, shift right=0.8ex, swap, "\otimes^\rL_{\BB_{[1, p]}, \wphi} \BB_{[1, 1]}"]  &
\Perf_v(\blank; \BB_{[1,1]})\Bigr)\end{tikzcd}.
\]
By construction, $\Perf_\an(\FFF_\blank;\O)$ comes equipped with a conservative symmetric monoidal functor 
\[
    \pi^*_{\blank}\colon \Perf_\an(\FFF_\blank; \O) \to \Perf_v(\blank; \BB_{[1, p]}).
\]
\end{defn}

\cref{defn:analytic-perfectoid-complexes-FF} might initially seem counterintuitive, as we define analytic sheaves on the relative Fargues--Fontaine curve in terms of $v$-sheaves on $X$. However, this notation is justified by \cref{rmk:complexes-FF-affinoid} below. Conceptually, one can view $\Perf_\an(\FFF_X; \O)$ as complexes that behave as $v$-sheaves in the ``$X$-direction'' and as analytic sheaves in the ``Fargues--Fontaine direction'' (see \cref{rmk:extend-from-perfectoids} and \cref{rmk:unique-extension} for precise statements).

\begin{rmk}\label{rmk:complexes-FF-affinoid} 
\cref{cor:Perf-on-FF} guarantees that, for an affinoid perfectoid space $S$ over $\Spa(\QQ_p, \ZZ_p)$, there is a natural equivalence $\Perf_\an(\rm{FF}_{S^\diam}; \O)\simeq \Perf_{\an}(\rm{FF}_S; \O)$. Here, $\Perf_\an(\rm{FF}_{S^\diam}; \O)$ is the category from \cref{defn:analytic-perfectoid-complexes-FF} and $\Perf_{\an}(\rm{FF}_S; \O)$ is the usual category of analytic perfect complexes on the adic space $\rm{FF}_S$ (see \cref{defn:FF-curve}). This equivalence justifies the notation $\Perf_\an(\FFF_X; \O)$ used in \cref{defn:analytic-perfectoid-complexes-FF}.
\end{rmk}

For the purposes of \cref{defn:analytic-perfectoid-complexes-FF}, the notation ``$\Perf_\an(\mathrm{FF}_X; \O)$'' can simply be treated as an indivisible formal symbol. While it is true that any diamond $X$ over $\Spd(\QQ_p, \ZZ_p)$ yields a diamond $\mathrm{FF}^\diam_X\coloneqq \Spd \QQ_p \times X/(\varphi^\ZZ \times \id)$
equipped with a structure $v$-sheaf $\O_{\mathrm{FF}^\diam_X}$ (via its projection to $\Spd(\QQ_p, \ZZ_p)$ as in \cref{defn:period-sheaves}), evaluating $v$-complexes on this space will not recover the analytic category we just defined.

\begin{warning}\label{warning:analytic-vs-v} 
There is a natural functor $\Perf_\an(\mathrm{FF}_X; \O) \to \Perf_v(\mathrm{FF}_X^\diam; \O)$ for any diamond $X$ over $\Spd(\QQ_p, \ZZ_p)$. However, even when $X$ is an affinoid perfectoid space over $\Spd(\QQ_p, \ZZ_p)$, this functor fails to be fully faithful or essentially surjective. This failure is visible on derived global sections: for $X=\Spd(\CC_p, \O_{\CC_p})$, we have $\rR\Gamma_\an(\FFF_{X}; \O) \simeq \QQ_p[0]$, while $\rR\Gamma_v(\FFF^\diam_{X}; \O) \simeq \QQ_p[0] \oplus \QQ_p[-1]$.    
\end{warning}

By construction, the $\infty$-category $\Perf_\an(\FFF_X; \O)$ has a natural pullback functoriality.

\begin{rmk}\label{rmk:pullback-perf} By definition, the assignment $X \mapsto \Perf_\an(\FFF_X; \O)\in \Cat_\infty^\otimes$ is functorial in $\Diam_{/\QQ_p}^\op$. In other words, for any morphism of diamonds $f\colon X' \to X$ over $\Spd(\QQ_p, \ZZ_p)$, there is a specified symmetric monoidal pullback functor $f_{\FFF}^*\colon \Perf_\an(\FFF_X; \O) \to \Perf_\an(\FFF_{X'}; \O)$ that fits into the following commutative diagram:
\begin{equation}\label{eqn:pullback-diagram-perf}
    \begin{tikzcd}[column sep = 6em, row sep = 4em]\Perf_\an(\FFF_{X}; \O) \arrow[d, "f_{\FFF}^*"] \arrow[r, "\pi_X^*"] & \Perf_v(X; \BB_{[1, p]})\arrow[d, "f_{v, [1, p]}^*"] \arrow[r, shift left=0.8ex, "\otimes^\rL_{\BB_{[1, p]}, \res} \BB_{[1, 1]}"]\arrow[r, shift right=0.8ex, swap, "\otimes^\rL_{\BB_{[1, p]}, \wphi} \BB_{[1, 1]}"] &[4em] \Perf_v(X; \BB_{[1, 1]})\arrow[d, "f_{v, [1, 1]}^*"] \\
    \Perf_\an(\FFF_{X'}; \O) \arrow[r, "\pi^*_{X'}"] & \Perf_v(X'; \BB_{[1, p]}) \arrow[r, shift left=0.8ex, "\otimes^\rL_{\BB_{[1, p]}, \res} \BB_{[1, 1]}"]\arrow[r, shift right=0.8ex, swap, "\otimes^\rL_{\BB_{[1, p]}, \wphi} \BB_{[1, 1]}"] & \Perf_v(X'; \BB_{[1, 1]}) \end{tikzcd}
\end{equation}
\end{rmk}

\begin{lemma}\label{lemma:an-perf-v-sheaf} The functor $\Perf_{\an}(\FFF_{\blank}; \O) \colon \Diam^{\op}_{/\QQ_p} \to \Cat^\otimes_\infty$ is a hypercomplete $v$-sheaf of $\infty$-categories.
\end{lemma}
\begin{proof}
    Since hypercomplete $v$-sheaves are closed under limits, it suffices to show that $\Perf_v(\blank; \BB_{[1,p]})$ and $\Perf_v(\blank; \BB_{[1,1]})$ are hypercomplete $v$-sheaves of symmetric monoidal $\infty$-categories. This follows immediately from \cref{lemma:perfect-complexes-are-hypercomplete} and the observation that $\Cat_\infty^\otimes \to \Cat_\infty$ is conservative and commutes with all limits (see \cite[Lem.~B.2.4(i)]{Heyer-Mann}). 
\end{proof}

\cref{lemma:an-perf-v-sheaf} provides a systematic way to compute $\Perf_\an(\FFF_X; \O)$ for any diamond $X$ over $\Spd(\QQ_p, \ZZ_p)$ in terms of affinoid perfectoid spaces mapping to $X$.

\begin{rmk}\label{rmk:extend-from-perfectoids} As a direct consequence of \cref{rmk:complexes-FF-affinoid} and \cref{lemma:an-perf-v-sheaf}, for any diamond $X$ over $\Spd(\QQ_p, \ZZ_p)$, there is a canonical identification $\Perf_\an(\mathrm{FF}_X; \O) \simeq \lim_{S \to X} \Perf_\an(\mathrm{FF}_{S^\sharp}; \O)$, where the limit is taken over the category of affinoid perfectoid spaces equipped with a morphism $S \to X$.
\end{rmk}

In fact $\Perf_\an(\FFF_\blank; \O)$ is the essentially unique way to extend 
$\Perf_\an(\FFF_\blank; \O) \colon \Perfd_{\QQ_p}^{\aff, \op} \to \Cat_\infty^\otimes$
to a hypercomplete $v$-sheaf on $\Diam_{/\QQ_p}^{\op}$, by the following remark:

\begin{rmk}
\label{rmk:unique-extension} 
Recall that affinoid perfectoid spaces over $\Spd(\QQ_p, \ZZ_p)$
form a basis of the $v$-site of $\Spd(\QQ_p, \ZZ_p)$. Therefore,  
\cite[Cor.\,A.7]{Ko20} ensures that $\Perf_{\an}(\rm{FF}_{\blank}; \O) \colon \Diam^{\op}_{/\QQ_p} \to \Cat_\infty$ is the (essentially) unique extension of $\Perf_{\an}(\rm{FF}_{\blank}; \O) \colon \Perfd^{\aff, \op}_{/\QQ_p} \to \Cat_\infty$ to a hypercomplete sheaf on $\Diam_{/\QQ_p}$.
\end{rmk}

While we could have used \cite[Cor.\,A.7]{Ko20} in order to \emph{define} $\Perf_\an(\FFF_X; \O)$, we find \cref{defn:analytic-perfectoid-complexes-FF} easier to work with for our purposes. 

We end the section by recording the following construction: 
\begin{construction}\label{construction:infinity-pullback} Let $X$ be a diamond over $\Spd(\QQ_p, \ZZ_p)$. We define the functor $\iota_\infty^*\colon \Perf_\an(\FFF_X; \O) \to \Perf_v(X; \O)$ as the composition of functors
\[
\Perf_\an(\FFF_X; \O) \xr{\pi_X^*} \Perf_v(X; \BB_{[1, p]}) \xr{\blank \otimes^{\rL}_{\BB_{[1,p], \theta_{[1, p]}}} \O} \Perf_v(X; \O),
\]
where $\pi_X^*$ is from \cref{defn:analytic-perfectoid-complexes-FF} and $\theta_{[1,p]}$ is from \cref{construction:theta-sheaves}. 

When $X=S$ is an affinoid perfectoid space with a map $S \to \Spd(\QQ_p, \ZZ_p)$ corresponding to an untilt $S^\sharp$ of $S$, \cref{cor:perfect-O-modules-on-affinoid-perfectoid} and \cref{rmk:complexes-FF-affinoid} allow us to identify $\iota_\infty^*$ with a functor 
\[
\Perf_\an(\FFF_{S^\sharp}; \O) \to \Perf\bigl(\O(S^\sharp)\bigr).
\]
\end{construction}
In \cref{section:pushforward}, we will also discuss the pushforward functoriality of $\Perf_\an(\FFF_X; \O)$.

\subsection{Riemann--Hilbert and solution functors}

In this subsection, we discuss the Riemann--Hilbert and solution functors relating $\Perf_v(X; \uQ_p)$ to $\Perf_\an(\FFF_X; \O)$ 
for any diamond $X$ over $\Spd(\QQ_p, \ZZ_p)$. To rigorously define these functors, 
we introduce an auxiliary $\infty$-category $\cD'(\FFF_X; \O)$ that contains $\Perf_\an(\FFF_X; \O)$ as a full subcategory. 
While our definition of $\cD'$ is somewhat ad hoc, its primary advantage is presentability, 
which allows us to construct the desired adjoints via the Adjoint Functor Theorem.

Throughout this subsection, we fix a diamond $X$ over $\Spd(\QQ_p, \ZZ_p)$. In what follows, we will freely use the fact that the $\infty$-category $\cD(X_v; \BB_I)$ is presentable for any small $v$-stack $X$ and closed interval $I\subset (0, \infty)$ with rational endpoints (see \cite[Prop.~1.3.5.21]{HA}). 

\begin{defn}
\label{defn:auxiliary-analytic-perfectoid-complexes-FF} 
We define the functor $\cD'(\FFF_\blank; \O) \colon \Diam_{/\QQ_p}^{\op} \to \CAlg(\Pr^\rL)$ as the equalizer 
\[
\cD'(\FFF_\blank; \O) \coloneqq \begin{tikzcd}[column sep = 10em]\eq\bigl(\cD(\blank_v; \BB_{[1,p]}) 
\arrow[r, shift left=0.8ex, "\otimes^\rL_{\BB_{[1, p]}, \res} \BB_{[1, 1]}"]
   \arrow[r, shift right=0.8ex, swap, "\otimes^\rL_{\BB_{[1, p]}, \wphi} \BB_{[1, 1]}"]  &
\cD(\blank_v; \BB_{[1,1]})\bigr)\end{tikzcd}.
\]
By construction, the $\infty$-category $\cD'(\FFF_\blank;\O)$ comes equipped with a conservative symmetric monoidal functor $\pi^*_{\blank}\colon \cD'(\FFF_\blank; \O) \to \cD((\blank)_v; \BB_{[1, p]})$. 
\end{defn}

Let us first unpack this definition to see what kind of structure  $\cD'(\FFF_X; \O)$ admits. 

\begin{rmk}\label{rmk:pullback-cD'} Recall that $\CAlg(\Pr^\rL)$ is the $\infty$-category of presentably symmetric monoidal $\infty$-categories. Therefore, for each diamond $X$ over $\Spd(\QQ_p, \ZZ_p)$, the $\infty$-category $\cD'(\FFF_X; \O)$ is presentable and is naturally a symmetric monoidal category such that the tensor product functor, which we denote by  
\[
    \blank \otimes^\rL_\O \blank \colon \cD'(\FFF_X; \O)\times \cD'(\FFF_X; \O) \to \cD'(\FFF_X; \O),
\]
commutes with colimits in each variable. Furthermore, the $\CAlg(\Pr^\rL)$-functoriality of this construction ensures that, for every morphism $f\colon X' \to X$ in $\Diam_{/\QQ_p}$, there is a specified symmetric monoidal colimit-preserving pullback functor $f_{\FFF}^*\colon \cD'(\FFF_X; \O) \to \cD'(\FFF_{X'}; \O)$. This pullback functor fits into a diagram analogous to \cref{eqn:pullback-diagram-perf}. 
\end{rmk}

\begin{rmk}\label{rmk:inner-hom} Since $\cD'(\FFF_X; \O)$ is a presentably symmetric monoidal category, it is a closed symmetric monoidal category in the sense of \cite[Def.~4.1.1.15]{HA}.\footnote{To justify that $\cD'(\FFF_X; \O)$ is closed, it suffices to show that $F\otimes^\rL \blank \colon \cD'(\FFF_X; \O) \to \cD'(\FFF_X; \O)$ admits a right adjoint. This follows from the fact that it commutes with colimits and \cite[Cor.~5.5.2.9]{HTT}} Thus, \cite[Rmk.~4.2.1.31 and Ex.~4.2.1.32]{HA} imply that we get an internal  Hom functor on $\cD'(\FFF_X; \O)$, which we denote by 
\[
\rR{\cHom}_{\FFF_X}(\blank, \blank) \colon \cD'(\FFF_X; \O)^{\mathrm{op}}
\times \cD'(\FFF_X; \O)
\to \cD'(\FFF_X; \O).
\]
\end{rmk}

We also provide a more explicit description of objects of $\cD'(\FFF_X; \O)$. 

\begin{rmk}\label{rmk:explicit} Since the forgetful functor $\CAlg(\Pr^\rL) \to \Cat_\infty$ commutes with limits (\cite[Prop.~5.5.3.13]{HTT}, \cite[Cor.~3.2.2.5]{HA}), we can apply \cite[\href{https://kerodon.net/tag/03H9}{Tag 03H9}]{kerodon} to explicitly describe the objects of $\cD'(\FFF_X; \O)$. Specifically, they can be identified with pairs $(F, \alpha)$, where $F \in \cD(X_v; \BB_{[1, p]})$ and $\alpha$ is an isomorphism $\alpha\colon F \otimes^\rL_{\BB_{[1, p]}, \wphi} \BB_{[1,1]}\xr{\sim} F \otimes^\rL_{\BB_{[1, p]}, \res} \BB_{[1,1]}$. 
\end{rmk}

Furthermore, the symmetric monoidal $\infty$-category $\cD'(\FFF_X; \O)$ recovers $\Perf_\an(\FFF_X; \O)$ by passing to the subcategory of dualizable objects. 

\begin{construction}\label{construction:dualizable-in-D'} \cite[Prop.~4.3.1.11]{HA} and \cite[\href{https://stacks.math.columbia.edu/tag/0FPV}{Tag 0FPV}]{stacks-project} ensure that the category of dualizable objects in $\cD'(\FFF_X; \O)$ is canonically identified with $\Perf_\an(\FFF_X; \O)$. In particular, we have a canonical fully faithful functor 
\[
j_{\FFF_X} \colon \Perf_\an(\FFF_X; \O) \hookrightarrow \cD'(\FFF_X; \O).
\]
\end{construction}

We now turn to the construction of the Riemann--Hilbert functor $\RH_X\colon \cD(X_v; \uQ_p) \to \cD'(\FFF_X; \O)$ and its right adjoint $\Sol_X \colon \cD'(\FFF_X; \O) \to \cD(X_v; \uQ_p)$. The definition of $\RH$ relies on the following universal property of our auxiliary category:

\begin{rmk}\label{rmk:maps-into-D'} For any $\cC\in \CAlg(\Pr^\rL)$, the mapping space into the equalizer $\cD'(\FFF_X; \O)$ is naturally identified as:
\[
\Map_{\CAlg(\Pr^\rL)}\bigl(\cC, \cD'(\FFF_X; \O)\bigr) \simeq \eq\left( \begin{tikzcd}[column sep = 7em] \Map_{\CAlg(\Pr^\rL)}\left(\cC, \cD(X_v; \BB_{[1, p]})\right) \arrow[r, shift left=0.8ex, "\otimes^\rL_{\BB_{[1, p]}, \res} \BB_{[1, 1]}"] \arrow[r, shift right=0.8ex, swap, "\otimes^\rL_{\BB_{[1, p]}, \wphi} \BB_{[1, 1]}"]  & \Map_{\CAlg(\Pr^\rL)}\left(\cC, \cD(X_v; \BB_{[1, 1]})\right) \end{tikzcd} \right).
\]
By \cite[\href{https://kerodon.net/tag/03HA}{Tag 03HA}]{kerodon} and \cite[\href{https://kerodon.net/tag/03H9}{Tag 03H9}]{kerodon}, this implies that constructing a colimit-preserving symmetric monoidal functor $\cD(X_v; \uQ_p) \to \cD'(\FFF_X; \O)$ is equivalent to specifying a symmetric monoidal, colimit-preserving functor $F \colon \cD(X_v; \uQ_p) \to \cD(X_v; \BB_{[1,p]})$ together with a natural equivalence 
\[
\alpha\colon F(\blank)\otimes^\rL_{\BB_{[1, p]}, \wphi} \BB_{[1,1]} \xr{\sim} F(\blank) \otimes^\rL_{\BB_{[1, p]}, \res} \BB_{[1,1]}.
\]
\end{rmk}

To use this observation, we note that, for each closed interval $I\subset (0, \infty)$ with rational endpoints, there is a canonical morphism $\iota_I\colon \ud{\ZZ}_p \to \AA_I$ coming from the fact that $\AA_I$ is a $p$-complete sheaf of $\ZZ_p$-algebras
(see \cref{lemma:A_inf-is-p-adically-complete}). In particular, it induces a canonical morphism $\iota_I\colon \ud{\QQ}_p \to \BB_I$. This leads to the following construction: 

\begin{construction}\label{constr:RH} To define the functor $\RH_X \colon \cD(X_v; \ud{\QQ}_p) \to \cD'(\FFF_X; \O)$ in $\CAlg(\Pr^\rL)$, we use the equivalence from \cref{rmk:maps-into-D'} and specify the pair $(F, \alpha)$. First, we set the functor $F$ to be:
\[ F \coloneqq \blank \otimes^\rL_{\ud{\QQ}_p, \iota_{[1, p]}} \BB_{[1,p]} \colon \cD(X_v; \ud{\QQ}_p) \to \cD(X_v; \BB_{[1,p]}). \]
Next, we define $\alpha$ as the natural isomorphism arising from the fact that the two distinct compositions mapping $\ud{\QQ}_p$ to $\BB_{[1,1]}$ are identical; namely, both
\[
\left(\ud{\QQ}_p \xr{\iota_{[1,p]}} \BB_{[1, p]} \xr{\res} \BB_{[1,1]}\right) \quad \text{and} \quad \left(\ud{\QQ}_p \xr{\iota_{[1,p]}} \BB_{[1, p]} \xr{\wphi} \BB_{[1,1]}\right)
\]
are equal to the canonical map $\iota_{[1,1]} \colon \ud{\QQ}_p \to \BB_{[1,1]}$.\end{construction}

\begin{variant}\label{RH-perfect}
Because $\RH_X \colon \cD(X_v; \ud{\QQ}_p) \to \cD'(\FFF_X; \O)$ is a symmetric monoidal functor, it preserves dualizable objects. Consequently, by \cref{construction:dualizable-in-D'} and \cite[\href{https://stacks.math.columbia.edu/tag/0FPV}{Tag 0FPV}]{stacks-project} applied to $\cD(X_v; \uQ_p)$, it restricts to a symmetric monoidal functor $\RH_X \colon \Perf_v(X; \ud{\QQ}_p) \to \Perf_\an(\FFF_X; \O)$ that fits into the following commutative diagram:
\[\begin{tikzcd}[column sep = 5em]
\Perf_v(X; \uQ_p)\arrow[d, hook, "j_X"] \arrow[r, "\RH_X"] & \Perf_\an(\FFF_X; \O) \arrow[d, hook, "j_{\FFF_X}"] \\
\cD(X_v; \uQ_p) \arrow[r, "\RH_X"] & \cD'(\FFF_X; \O),
\end{tikzcd}
\]
where $j_X$ is the natural embedding.
\end{variant}

\begin{variant} We also denote by $\RH^{\ZZ_p}_X\colon \Perf_v(X; \uZ_p) \to \Perf_\an(\FFF_X; \O)$ the composite functor 
\[
\Perf_v(X; \uZ_p) \to \Perf_v(X; \uQ_p)\xr{\RH_X} \Perf_\an(\FFF_X; \O).
\]
\end{variant}

\begin{rmk}\label{rmk:pullback-of-RH} By construction, the composition $\Perf_v(X; \uZ_p) \xr{\RH^{\ZZ_p}_X} \Perf_\an(\FFF_X; \O) \xr{\pi^*_X} \Perf_v(X; \BB_{[1, p]})$ is canonically identified with the functor $\blank \otimes^\rL_{\uZ_p} \BB_{[1, p]}\colon \Perf_v(X; \uZ_p) \to \Perf_v(X; \BB_{[1, p]})$.
\end{rmk}

We now define the solution functor. 

\begin{construction}\label{construction:R-tau} 
Because $\RH_X \colon \cD(X_v; \ud{\QQ}_p) \to \cD'(\FFF_X; \O)$ is a colimit-preserving symmetric monoidal functor between presentable categories, the Adjoint Functor Theorem \cite[Cor.~5.5.2.9]{HTT} and \cite[Lem.~B.1.12(i)]{Heyer-Mann} guarantee the existence of a lax symmetric monoidal right adjoint:
\[ 
\Sol_X\colon \cD'(\FFF_X; \O) \to \cD(X_v; \ud{\QQ}_p). 
\]
By a slight abuse of notation, we also use $\Sol_X\colon \Perf_\an(\FFF_X;\O) \to \cD(X_v; \ud{\QQ}_p)$ to denote the restriction of this adjoint to the subcategory $\Perf_\an(\FFF_X; \O)$. 
\end{construction}

This functor admits an explicit description.

\begin{rmk}\label{rmk:explicit-adjoint} 
Using the identification from \cref{rmk:explicit} and \cite[Lem.~D.4.7]{Heyer-Mann}, $\Sol_X \colon \cD'(\FFF_X; \O) \to \cD(X_v; \uQ_p)$ can be computed by
$$ (F, \alpha) \mapsto \mathrm{eq}\left(\begin{tikzcd} F \arrow[r, shift left=0.8ex, "\res"]
   \arrow[r, shift right=0.8ex, swap, "\alpha\circ \wphi"] & F\otimes^\rL_{\BB_{[1,p]}, \res} \BB_{[1,1]} \end{tikzcd}\right) \simeq \fib\left(F \xr{\res-\alpha\circ \wphi} F\otimes^\rL_{\BB_{[1,p]}, \res} \BB_{[1,1]}\right), $$
where $\res\colon F \to F\otimes^\rL_{\BB_{[1,p]}, \res} \BB_{[1,1]}$ and $\wphi \colon F \to F\otimes^\rL_{\BB_{[1,p]}, \wphi} \BB_{[1,1]}$ are the natural maps. Using \cref{rmk:extend-from-perfectoids}, \cref{cor:Perf-on-FF}, and the fact that affinoid perfectoid spaces form a basis on $X_v$, the restricted functor $\Sol_X\colon \Perf_\an(\FFF_X; \O) \to \cD(X_v; \ud{\QQ}_p)$ can be informally described by the assignment
\[
\Perf_\an(\FFF_X; \O) \ni \cP \mapsto \left(S\mapsto \rR\Gamma_\an(\FFF_{S^\sharp}, \cP|_{\FFF_{S^\sharp}})\right), 
\]
where $S$ varies over all affinoid perfectoid spaces mapping to the diamond $X$.
\end{rmk}

For future reference, we record the following result. 

\begin{lemma}\label{lemma:coh-dimension-sol} Let $S$ be an affinoid perfectoid space over $\Spd(\QQ_p, \ZZ_p)$ and let $\cE \twoheadrightarrow \cE'$ be a surjection of vector bundles on $\FFF_S$. Then the following statements hold:
\begin{enumerate}[leftmargin=*,label={\upshape{(\roman*)}}]
    \item\label{lemma:coh-dimension-sol-1} the complex $\Sol_S(\cE)$ is concentrated in degrees $0$ and $1$;
    \item\label{lemma:coh-dimension-sol-2} the natural morphism $\cal{H}^1\bigl(\Sol_S(\cE)\bigr) \to \cal{H}^1\bigl(\Sol_S(\cE')\bigr)$ is a surjection in $\Shv(S_v; \uQ_p)$.
\end{enumerate}
\end{lemma}
\begin{proof}
    Under the equivalence from \cref{cor:Perf-on-FF}, a vector bundle on $\FFF_S$ corresponds to a finite projective $\BB_{[1,p]}(S)$-module $M_{[1,p]}$ together with an isomorphism $\alpha\colon \res^*M_{[1,p]} \xr{\sim}\widetilde{\varphi}^*M_{[1,p]}$. Equivalently, it is given by a (finite locally free) $\BB_{[1,p]}$-module $\widetilde{M_{[1,p]}}$ together with an isomorphism $\alpha\colon \widetilde{M_{[1,p]}} \otimes_{\BB_{[1, p]},\res} \BB_{[1,1]} \simeq \widetilde{M_{[1,p]}} \otimes_{\BB_{[1, p]},\wphi} \BB_{[1,1]}$. Then \cref{rmk:explicit-adjoint} ensures that
    \[
    \Sol_S(\cE) \simeq \fib\left(\widetilde{M_{[1,p]}} \xr{\res-\alpha\circ \wphi} \widetilde{M_{[1,p]}} \otimes_{\BB_{[1,p]}, \res} \BB_{[1,1]}\right).
    \]
    So it is evidently concentrated in degrees $0$ and $1$. This finishes the proof of \cref{lemma:coh-dimension-sol-1}.

    Lastly, \cref{lemma:coh-dimension-sol-2} follows immediately from \cref{lemma:coh-dimension-sol-1} by observing that $\ker(\cE \to \cE')$ is a vector bundle since it is a kernel of a surjective morphism of vector bundles. 
\end{proof}

\begin{construction}
\label{construction:cup-product-dR} 
The lax symmetric monoidal structure on $\Sol_X$ provides a functorial morphism $\cup \colon \Sol_X(\F_1)\otimes^\rL_{\uQ_p} \Sol_X(\F_2) \to \Sol_X(\F_1\otimes^\rL_\O\F_2)$ for any $\F_1, \F_2 \in \cD'(\FFF_X; \O)$. Explicitly, this morphism is adjoint to the morphism
\[
\RH_X\bigl( \Sol_X(\F_1)\otimes^\rL_{\uQ_p} \Sol_X(\F_2) \bigr) \simeq \RH_X\bigl( \Sol_X(\F_1)\bigr) \otimes^\rL_\O \RH_X\bigl( \Sol_X(\F_2)\bigr) \xr{\epsilon_{\F_1} \otimes^\rL_{\O} \epsilon_{\F_2}} \F_1 \otimes^\rL_\O \F_2,
\]
where $\epsilon_{\F_i}$ are the counits of the $(\RH_X, \Sol_X)$-adjunction. 
\end{construction}

Now we study the basic properties of the $\RH$ and $\Sol$ functors. We start with the following preliminary lemma: 

\begin{lemma}\label{lemma:global-sections-FF-curve} Let $X$ be a diamond over $\Spd(\QQ_p, \ZZ_p)$. Then the natural morphism $\uQ_p \to \fib(\BB_{[1, p]} \xr{\res-\wphi} \BB_{[1,1]})$ is an isomorphism.
\end{lemma}
\begin{proof}
    It suffices to show that 
    \begin{equation}\label{eqn:exactB}
    0 \to \ud{\QQ}_p(S) \to \BB_{[1, p]}(S) \xr{\res-\wphi} \BB_{[1,1]}(S) \to 0
    \end{equation}
    is a short exact sequence for any strictly totally disconnected perfectoid space $S$ over $X$. Let $S^\sharp$ be the untilt of $S$ corresponding to the composite map $S\to \Spd(\QQ_p, \ZZ_p)$. Then \cref{cor:Perf-on-FF} and \cref{cor:perfect-O-modules-on-affinoid-perfectoid-BBI} imply that $\rR\Gamma_\an(\FFF_{S^\sharp}, \O)$ is isomorphic to $\fib\bigl(\BB_{[1, p]}(S) \xr{\res-\wphi} \BB_{[1,1]}(S)\bigr)$. Thus, \cite[Prop.~II.2.5]{Fargues-Scholze} 
    and \cref{thm:v-perfect-complexes-vs-etale-perfect-complexes-rationally-Qell}~\cref{thm:v-perfect-complexes-vs-etale-perfect-complexes-rationally-Qell-3} imply that the natural map
    \[
    \rR\Gamma_v(S, \uQ_p) \simeq \uQ_p(S) \to \fib\bigl(\BB_{[1, p]}(S) \xr{\res-\wphi} \BB_{[1,1]}(S)\bigr) \simeq \rR\Gamma_\an(\FFF_{S^\sharp}, \O)
    \]
    is an isomorphism. Thus, \cref{eqn:exactB} is indeed exact. 
\end{proof}

\begin{proposition}
\label{prop:dr-rh} 
Let $X$ be a diamond over $\Spd(\QQ_p, \ZZ_p)$ and let $\cP\in \Perf_v(X; \ud{\QQ}_p)$. Then the unit of the adjunction $\can_\cP \colon \cP \to \Sol_X(\RH_X(\cP))$ is an isomorphism.
\end{proposition}
\begin{proof}
    \cref{rmk:explicit-adjoint} implies that it suffices to show that the natural morphism
    \[
    \cP \to \fib\bigl(\cP \otimes^\rL_{\ud{\QQ}_p} \BB_{[1,p]} \xr{\id\otimes \res-\id \otimes \wphi} \cP\otimes^\rL_{\ud{\QQ}_p} \BB_{[1,1]}\bigr)
    \]
    is an isomorphism. Since $\cP$ is dualizable, we conclude that $\fib(\cP \otimes^\rL_{\ud{\QQ}_p} \BB_{[1,p]} \xr{\id\otimes \res-\id \otimes \wphi} \cP\otimes^\rL_{\ud{\QQ}_p} \BB_{[1,1]}) \simeq \cP \otimes^\rL_{\ud{\QQ}_p} \fib(\BB_{[1,p]} \xr{\res-\wphi} \BB_{[1,1]})$. Then the result follows immediately from \cref{lemma:global-sections-FF-curve}.
\end{proof}

Finally, we finish the subsection with the theorem of Ansch\"utz--Le Bras, which is indispensable for our proof of the geometric duality. Before we formulate it, we note that for any $\F_1, \F_2\in \cD'(\FFF_X; \O)$, there is a natural morphism 
\[
\Sol_X\bigl(\rR{\cHom}_{\FFF_X}({\mathcal F}_1,{\mathcal F}_2) \bigr)\to \rR\cHom_{\cD(X_v;\uQ_p)}\bigl(\Sol_X({\mathcal F}_1),\Sol_X({\mathcal F}_2)\bigr)
\]
adjoint to the composition
\[
\Sol_X\bigl(\rR{\cHom}_{\FFF_X}({\mathcal F}_1,{\mathcal F}_2) \bigr) \otimes^\rL_{\uQ_p} \Sol_X(\F_1) \xr{\cup} \Sol_X\bigl( \rR{\cHom}_{\FFF_X}({\mathcal F}_1,{\mathcal F}_2)\otimes^\rL_\O \F_1\bigr) \xr{\Sol(ev)} \Sol_X(\F_2),
\]
where $\cup$ is the cup-product map from \cref{construction:cup-product-dR}. 

\begin{theorem}(Ansch\"utz--Le Bras \cite[Cor. 3.11]{ALB})
\label{ff1}
Let $X$ be a diamond over $\Spd(\mathbf{Q}_p, \mathbf{Z}_p)$. Then the functor
\begin{equation}\label{kicia0}
\Sol_X \colon  \mathrm{Perf}_{\mathrm{an}} (\mathrm{FF}_{X};\mathcal{O} ) \to {\mathcal{D}}(X_v; \uQ_p)
\end{equation}
is fully faithful. In particular, for any $\F_1, \F_2\in \Perf_{\an}(\mathrm{FF}_{X};\mathcal{O})$, the natural morphism 
\[
\Sol_X\bigl(\rR{\cHom}_{\FFF_X}({\mathcal F}_1,{\mathcal F}_2) \bigr)\to \rR\cHom_{\cD(X_v;\uQ_p)}\bigl(\Sol_X({\mathcal F}_1),\Sol_X({\mathcal F}_2)\bigr)
\]
is an isomorphism. 

\end{theorem}

\begin{warning} 
We note a discrepancy in the formulation of \cite[Cor. 3.11]{ALB}. There, Ansch\"utz--Le Bras define a functor $\rR \tau_*$ that sends a complex $\cP\in \Perf_\an(\FFF_X; \O)$ to the $v$-sheaf (informally) given by the assignment
\[
 \Perf^{\aff, \op}_{/X} \ni S \mapsto \rR\Gamma_v(\FFF_{S^\sharp}; \cP|_{\FFF_{S^\sharp}}). 
\]
In contrast, our functor $\Sol_X$ sends $\cP$ to the $v$-sheaf:
\[
 \Perf^{\aff, \op}_{/X} \ni S \mapsto \rR\Gamma_\an(\FFF_{S^\sharp}; \cP|_{\FFF_{S^\sharp}}). 
\]
As demonstrated by \cref{warning:analytic-vs-v}, these two evaluations are different.  However,  their argument  uses  analytic Fargues--Fontaine cohomology and, in fact,   works with  $\Sol_X$ not with $\rR \tau_*$.
\end{warning}

\subsection{Pushforward functor}
\label{section:pushforward}

In this subsection, we define the (derived) pushforward functor for $\cD'(\FFF_X; \O)$ and show that it preserves $\Perf_\an$ for smooth proper morphisms. We also prove that the solution functor commutes with pushforward. 

For the next definition, we fix a morphism $f\colon X' \to X$ of diamonds over $\Spd(\QQ_p, \ZZ_p)$. This induces a colimit-preserving pullback functor $f_{\FFF}^*\colon \cD'(\FFF_X; \O) \to \cD'(\FFF_{X'}; \O)$, see \cref{rmk:pullback-cD'}. 
Since both $\cD'(\FFF_X;\O)$ and $\cD'(\FFF_{X'};\O)$ are presentable,
the Adjoint Functor Theorem \cite[Cor.~5.5.2.9]{HTT} implies that 
$f_{\FFF}^*$ admits a right adjoint functor. 

\begin{defn}\label{defn:pushforward-cD'} The \emph{derived pushforward functor} $\rR f_{\FFF,*}\colon \cD'(\FFF_{X'}; \O) \to \cD'(\FFF_{X}; \O)$ is the right adjoint to the pullback functor $f_{\FFF}^*$. 
\end{defn}

\begin{warning} We prefer to use the $\rR f_{\FFF, *}$ notation for this functor to be compatible with other pushforwards used in this paper. However, we emphasize that this functor is not a right derived functor.    
\end{warning}

First, we verify that the formation of $\rR f_{\FFF, *}$ commutes with the solution functor. 

\begin{lemma}\label{lemma:dr-pushforward} Let $f\colon X' \to X$ be a morphism of diamonds over $\Spd(\QQ_p, \ZZ_p)$. Then the diagrams
\[
\begin{tikzcd}
    \cD(X_v; \ud{\QQ}_p) \arrow[d, "f_v^*"] \arrow[r, "\RH_X"] & \cD'(\FFF_X; \O) \arrow[d, "f_{\FFF}^*"] \\
    \cD(X'_v; \ud{\QQ}_p) \arrow[r, "\RH_{X'}"] & \cD'(\FFF_{X'}; \O)
\end{tikzcd}
\quad \text{and} \quad
\begin{tikzcd}
    \cD'(\FFF_{X'}; \O)\arrow[d, "\rR f_{\FFF, *}"] \arrow[r, "\Sol_{X'}"] & \cD(X'_v; \ud{\QQ}_p) \arrow[d, "\rR f_{v, *}"] \\
    \cD'(\FFF_{X}; \O) \arrow[r, "\Sol_{X}"] & \cD(X_v; \ud{\QQ}_p)
\end{tikzcd}
\]
commute. In other words, there are canonical isomorphisms of functors 
\[
\alpha_f \colon \RH_{X'} \circ f_v^* \xr{\sim} f_{\FFF}^* \circ \RH_X \quad \text{and} \quad \alpha'_f \colon \Sol_{X}\circ \rR f_{\FFF, *}(\blank) \xr{\sim} \rR f_{v, *} \circ \Sol_{X'}(\blank).
\] 
\end{lemma}
\begin{proof}
    First, we note that the commutativity of the right diagram follows from the commutativity of the left diagram by passing to right adjoints. Therefore, it suffices to show that the left diagram commutes. Using the last sentence of \cref{rmk:pullback-cD'} and the definition of the $\RH$-functor (see \cref{constr:RH}), the question follows from the commutativity of the following diagram: 
    \[
    \begin{tikzcd}[column sep = 10em, row sep = 4em]
    \cD(X_v; \ud{\QQ}_p) \arrow[d, "f_v^*"] \arrow[r, "\otimes^\rL_{\ud{\QQ}_p, \iota_{[1,p]}} \BB_{[1, p]}"] & \cD(X_v; \BB_{[1, p]})\arrow[d, "f_v^*"] \arrow[r, shift left=0.8ex, "\otimes^\rL_{\BB_{[1, p]}, \res} \BB_{[1, 1]}"]
   \arrow[r, shift right=0.8ex, swap, "\otimes^\rL_{\BB_{[1, p]}, \wphi} \BB_{[1, 1]}"] & \cD(X_v; \BB_{[1, 1]})\arrow[d, "f_v^*"] \\
   \cD(X'_v; \ud{\QQ}_p) \arrow[r, "\otimes^\rL_{\ud{\QQ}_p, \iota_{[1,p]}} \BB_{[1, p]}"] & \cD(X'_v; \BB_{[1, p]}) \arrow[r, shift left=0.8ex, "\otimes^\rL_{\BB_{[1, p]}, \res} \BB_{[1, 1]}"]
   \arrow[r, shift right=0.8ex, swap, "\otimes^\rL_{\BB_{[1, p]}, \wphi} \BB_{[1, 1]}"] & \cD(X'_v; \BB_{[1, 1]}),
   \end{tikzcd}
   \]
    which follows immediately from the fact that the derived $\infty$-category $\cD$ is functorial in ringed topoi and the observation that $\res \circ \iota_{[1, p]} = \iota_{[1,1]} = \wphi \circ \iota_{[1, p]}$.
\end{proof}

Now we wish to show that $\rR f_{\FFF, *}$ preserves dualizable objects when $f$ is smooth and proper.
We do this in two steps: we first construct a right adjoint on the level of $\Perf_\an$ and then show that it must be compatible with $\rR f_{\FFF, *}$. 

\begin{lemma}\label{lemma:proper-pushforward-FF}
Let $f\colon X' \to X$ be a smooth proper morphism of diamonds over $\Spd(\QQ_p, \ZZ_p)$.
Then the functor $f_{\FFF}^*\colon \Perf_\an(\FFF_X; \O) \to \Perf_\an(\FFF_{X'}; \O)$ admits a right adjoint functor $\rR f^{\Perf}_{\FFF, *} \colon \Perf_\an(\FFF_{X'}; \O) \to \Perf_\an(\FFF_{X}; \O)$. Furthermore, the following diagram commutes:
\begin{equation}\label{eqn:pushforward}
\begin{tikzcd}[column sep = 8em, row sep = 4em]
    \Perf_\an(\FFF_{X'}; \O) \arrow[d, "\rR f^\Perf_{\FFF, *}"] \arrow[r, "\pi^*_{X'}"] & \Perf_v(X'; \BB_{[1, p]})\arrow[d, "\rR f_{v, *}"] \arrow[r, shift left=0.8ex, "\otimes^\rL_{\BB_{[1, p]}, \res} \BB_{[1, 1]}"]
   \arrow[r, shift right=0.8ex, swap, "\otimes^\rL_{\BB_{[1, p]}, \wphi} \BB_{[1, 1]}"] & \Perf_v(X'; \BB_{[1, 1]})\arrow[d, "\rR f_{v, *}"] \\
   \Perf_\an(\FFF_{X}; \O) \arrow[r, "\pi^*_{X}"] & \Perf_v(X; \BB_{[1, p]}) \arrow[r, shift left=0.8ex, "\otimes^\rL_{\BB_{[1, p]}, \res} \BB_{[1, 1]}"]
   \arrow[r, shift right=0.8ex, swap, "\otimes^\rL_{\BB_{[1, p]}, \wphi} \BB_{[1, 1]}"] & \Perf_v(X; \BB_{[1, 1]})
\end{tikzcd}
\end{equation}
\end{lemma}
\begin{proof}
    \cref{cor:proper-smooth-pushforward-is-perfect} implies that the functors $f_{v}^*\colon \Perf_v(X; \BB_{I}) \to \Perf_v(X'; \BB_{I})$ admit right adjoints $\rR f_{v, *} \colon \Perf_v(X'; \BB_{I}) \to \Perf_v(X; \BB_{I})$ for any closed interval $I\subset (0,\infty)$ with rational endpoints. Furthermore,
    \cref{cor:proper-smooth-pushforward-is-perfect} ensures the natural base change morphisms (i.e., Beck--Chevalley morphisms)
\[
\rR f_{v, *} (\blank) \otimes^\rL_{\BB_{[1, p]}, \res} \BB_{[1, 1]} \to \rR f_{v, *} (\blank \otimes^\rL_{\BB_{[1, p]}, \res} \BB_{[1, 1]}),
\]
\[
\rR f_{v, *} (\blank) \otimes^\rL_{\BB_{[1, p]}, \wphi} \BB_{[1, 1]} \to \rR f_{v, *} (\blank \otimes^\rL_{\BB_{[1, p]}, \wphi} \BB_{[1, 1]}),
\]
are isomorphisms. Therefore, \cite[Lem.~D.2.14]{Heyer-Mann} implies that the functor $f_{\FFF}^*\colon \Perf_\an(\FFF_X; \O) \to \Perf_\an(\FFF_{X'}; \O)$ admits a right adjoint making Diagram~\cref{eqn:pushforward} commute. 
\end{proof}

Finally, we wish to compare $\rR f^\Perf_{\FFF, *}$ to the other pushforward functor $\rR f_{\FFF, *}$ from \cref{defn:pushforward-cD'}. 

\begin{lemma}
\label{lemma:proper-pushforward-FF-coincide} 
Let $f\colon X' \to X$ be a smooth proper morphism of diamonds over $\Spd(\QQ_p, \ZZ_p)$.
Then the diagram
\begin{equation}\label{eqn:pushforward-preserves-perf}
\begin{tikzcd}[column sep = 6em]
    \Perf_\an(\FFF_{X'}; \O) \arrow[r, hook, "j_{\FFF_{X'}}"] \arrow[d, "\rR f^\Perf_{\FFF, *}"] & \cD'(\FFF_{X'}; \O) \arrow[d, "\rR f_{\FFF, *}"]\\
    \Perf_\an(\FFF_{X}; \O) \arrow[r, hook, "j_{\FFF_{X}}"] & \cD'(\FFF_{X}; \O)
\end{tikzcd}
\end{equation}
commutes, i.e., there is a natural equivalence of functors $\beta_f\colon j_{\FFF_X} \circ \rR f^\Perf_{\FFF, *} \to \rR f_{\FFF, *} \circ j_{\FFF_{X'}}$. In particular, $\rR f_{\FFF, *} \colon \cD'(\FFF_{X'}; \O) \to \cD'(\FFF_{X}; \O)$ sends dualizable objects to dualizable objects.  
\end{lemma}
\begin{proof}
    First, we note that the last sentence of the lemma follows immediately from the commutativity of Diagram~\cref{eqn:pushforward-preserves-perf} and the observation that $\Perf_\an(\FFF_{\blank}; \O)$ is exactly the subcategory of dualizable objects in $\cD'(\FFF_{\blank}; \O)$ (see the first sentence of \cref{construction:dualizable-in-D'}). 

    Now we show that Diagram~\cref{eqn:pushforward-preserves-perf} commutes. Since $j_{\FFF_X}$ and $j_{\FFF_{X'}}$ are fully faithful and both $\rR f^\Perf_{\FFF, *}$ and $\rR f_{\FFF, *}$ are defined as right adjoints to the corresponding pullback functors, it suffices to show that $\rR f_{\FFF, *}$ sends $\Perf_\an(\FFF_{X'}; \O)$ to $\Perf_\an(\FFF_{X}; \O)$. Informally, this obviously follows from \cref{cor:proper-smooth-pushforward-is-perfect} since the pullbacks $f^*_v \colon \Perf_v(X; \BB_I) \to \Perf_v(X'; \BB_I)$ for each $I \subset (0,\infty)$ admit right adjoints which commute with the pullbacks in the diagram defining $\Perf_\an(\FFF_X; \O)$. The rest of the argument is devoted to spelling out this argument in a rigorous way.  
    
    In order to do this, let $\cI$ be the $\infty$-category associated with the category $\begin{tikzcd} \bullet \arrow[r, shift left=0.8ex]
   \arrow[r, shift right=0.8ex] & \bullet \end{tikzcd}$. Furthermore, let $F_X \colon \cI \to \Cat_\infty$ be the functor associated with the diagram
    
    \[
    \begin{tikzcd}[column sep = 8em]\cD(X_v; \BB_{[1,p]}) 
    \arrow[r, shift left=0.8ex, "\otimes^\rL_{\BB_{[1, p]}, \res} \BB_{[1, 1]}"]
   \arrow[r, shift right=0.8ex, swap, "\otimes^\rL_{\BB_{[1, p]}, \wphi} \BB_{[1, 1]}"]  &
    \cD(X_v; \BB_{[1,1]})\end{tikzcd}
    \]
    and $\cE_X\coloneqq \int_\cI (F_X) \to \cI$ be the associated cocartesian fibration (see \cite[\href{https://kerodon.net/tag/026J}{Tag 026J}]{kerodon} and \cite[\href{https://kerodon.net/tag/026K}{Tag 026K}]{kerodon}). We define $F_{X'} \colon \cI \to \Cat_\infty$ and $\cE_{X'} \to \cI$ in the analogous way. The pullback functoriality gives us a natural transformation of functors $f^*\colon F_{X} \Rightarrow F_{X'}$ which translates into an $\cI$-morphism $\cE_{f^*}\colon \cE_X \to \cE_{X'}$ sending cocartesian edges to cocartesian edges
    (see \cite[Th.~8.16]{Cisinski--Nguyen}). 

    Analogously, we define $\cE^\Perf_X \subset \cE_X$ be a cocartesian subfibration over $\cI$ corresponding to $\Perf_\an(\FFF_X; \O)$ and $\cE^\Perf_{X'}\subset \cE_{X'}$ corresponding to $\Perf_\an(\FFF_X'; \O)$ respectively. Since the $v$-pullback preserves $\BB_I$-perfect complexes, the functor $\cE_{f^*}$ restricts to a functor $\cE_{f^*} \colon \cE^\Perf_X \to \cE^\Perf_{X'}$.
    
    By virtue of \cite[\href{https://kerodon.net/tag/02TK}{Tag 02TK}]{kerodon}, we can identify $\cD'(\FFF_X; \O) \simeq \lim F_X$ with the  $\infty$-category $\Fun^{\CCart}_{/\cI}(\cI, \cE_X)$ of cocartesian section of $\cE_X \to \cI$.
    Similarly, we identify $\Perf_\an(\FFF_X; \O)\simeq \Fun^{\CCart}_{/\cI}(\cI, \cE^\Perf_{X})$, $\cD'(\FFF_{X'}; \O) \simeq \Fun^{\CCart}_{/\cI}(\cI, \cE_{X'})$, and $\Perf_\an(\FFF_{X'}; \O)\simeq \Fun^{\CCart}_{/\cI}(\cI, \cE^\Perf_{X'})$. 
    
    Since $\cE_{f^*}$ preserves cocartesian section, it induces functors $\cE_{f^*} \colon \Fun^{\CCart}_{/\cI}(\cI, \cE^{(\Perf)}_X) \to \Fun^{\CCart}_{/\cI}(\cI, \cE^{(\Perf)}_{X'})$ which coincide with $f^*_{\FFF}\colon \cD'(\FFF_X; \O) \to \cD'(\FFF_{X'}; \O)$ and $f^*_{\FFF} \colon \Perf_\an(\FFF_X; \O) \to \Perf_\an(\FFF_{X'}; \O)$ under the identifications from above. Similarly, it also induces a functor $\widetilde{\cE}_{f^*} \colon \Fun_{/\cI}(\cI, \cE^{(\Perf)}_X) \to \Fun_{/\cI}(\cI, \cE^{(\Perf)}_{X'})$ on the $\infty$-categories of \emph{all} sections. 

    The functor $\rR f_{\FFF, *}$ corresponds to a functor $G\colon \Fun^{\CCart}_{/\cI}(\cI, \cE_{X'}) \to \Fun^{\CCart}_{/\cI}(\cI, \cE_{X})$ which is right adjoint to $\cE_{f^*}$. 
    Our goal is to show that $G$ sends elements of $\Fun^{\CCart}_{/\cI}(\cI, \cE^{\Perf}_{X'})$ to $\Fun^{\CCart}_{/\cI}(\cI, \cE^{\Perf}_{X})$. 
    The problem is that it is quite complicated to compute $G$ explicitly because it is not clear whether the termwise right adjoint to $\cE_{f^*}$ preserves cocartesian sections. 

    To avoid this issue, we use \cite[Def.~4.1 and Prop.~5.1]{HY17} to ensure that the functor $\widetilde{\cE}_{f^*} \colon \Fun_{/\cI}(\cI, \cE_X) \to \Fun_{/\cI}(\cI, \cE_{X'})$ admits a right adjoint $\widetilde{G} \colon \Fun_{/\cI}(\cI, \cE_{X'}) \to \Fun_{/\cI}(\cI, \cE_X)$ which is term-wise given by the right adjoint to $f^*$. 
    Therefore, it suffices to show that $\widetilde{G}$ sends elements of $\Fun^{\CCart}_{/\cI}(\cI, \cE^{\Perf}_{X'})$ to $\Fun^{\CCart}_{/\cI}(\cI, \cE^{\Perf}_{X})$. Indeed, using the fact that $\restr{\widetilde{\cE}_{f^*}}{\Fun^{\CCart}_{/\cI}(\cI, \cE_{X})} = \cE_{f^*}$, this would imply that $G(c)=\widetilde{G}(c)$ for any $c\in \Fun^{\CCart}_{/\cI}(\cI, \cE^{\Perf}_{X'})$ and, thus, we would have $G(c) \in \Fun^{\CCart}_{/\cI}(\cI, \cE^{\Perf}_{X})$ for any such $c$. 
    
    In other words, we need to show that for any $\cP\in \Fun^{\CCart}_{/\cI}(\cI, \cE^{\Perf}_{X'}) \simeq \Perf_\an(\FFF_{X'}; \O)$, the corresponding section $\widetilde{G}(\cP)$ is cocartesian (i.e., sends edges of $\cI$ to cocartesian edges of $\cE_{X}$). If we unravel what this means and use \cite[\href{https://kerodon.net/tag/01TB}{Tag 01TB}]{kerodon}, we see that this is equivalent to showing that, for any $\cP_I\in \Perf_v(X'; \BB_I)$ for $I=[1,1]$ or $I=[1, p]$, the object $\rR f_{v, *} \cP_I$ lies in $\Perf_v(X; \BB_I)$ and the natural transformations of functors 
    \[
    \rR f_{v, *}(\blank) \otimes^{\rL}_{\BB_{[1,p]}, \res}
    \BB_{[1,1]} \rightarrow
    \rR f_{v, *}(\blank \otimes^{\rL}_{\BB_{[1,p]}, \res} \BB_{[1,1]})
        \text{ and }
    \rR f_{v, *}(\blank) \otimes^{\rL}_{\BB_{[1,p]}, \wphi}
    \BB_{[1,1]} \rightarrow
    \rR f_{v, *}(\blank \otimes^{\rL}_{\BB_{[1,p]}, \wphi} \BB_{[1,1]})
    \]
    are equivalences of functors $\Perf_v(X'; \BB_{[1,p]}) \to \cD(X_v; \BB_{[1,1]})$. Both results follow immediately from \cref{cor:proper-smooth-pushforward-is-perfect}. This finally finishes the proof. 
\end{proof}

As a first application of these results, we show that the pushforward on the Fargues--Fontaine curve is compatible with the pushforward on the $v$-site with respect to the $\iota_\infty^*$-functor from \cref{construction:infinity-pullback}.
\begin{cor}
\label{cor:proper-pushforward-FF-infity-pullback} 
Let $f\colon X' \to X$ be a smooth proper morphism of diamonds over $\Spd(\QQ_p, \ZZ_p)$.
Then the diagram
\begin{equation}
\begin{tikzcd}[column sep = 6em]
    \Perf_\an(\FFF_{X'}; \O) \arrow[r, "\iota_\infty^*"] \arrow[d, "\rR f_{\FFF, *}"] & \Perf_v(X'; \O) \arrow[d, "\rR f_{v, *}"]\\
    \Perf_\an(\FFF_{X}; \O) \arrow[r, "\iota_\infty^*"] &  \Perf_v(X; \O)
\end{tikzcd}
\end{equation}
commutes, i.e., there is a canonical isomorphism of functors $\iota_\infty^*\circ \rR f_{\FFF, *} \xr{\sim}\rR f_{v, *}\circ \iota_\infty^*$.
\end{cor}
\begin{proof}
    The result follows immediately from \cref{lemma:proper-pushforward-FF-coincide},
    \cref{lemma:proper-pushforward-FF}, and \cref{cor:HT realization for proper smooth pushforward} applied to $I=[1, p]$. 
\end{proof}
As another immediate application of \cref{lemma:proper-pushforward-FF-coincide}, we show that $v$-cohomology of $\uQ_p$-local systems are always Banach--Colmez spaces. For this, we recall the precise definition of relative and absolute Banach--Colmez spaces. 

\begin{defn}\label{defn:relative-BC} The $\infty$-category of \emph{relative Banach--Colmez spaces} on $X$ is the full $\infty$-subcategory $\cD_{\BBC}^{(b)}(X_v; \uQ_p) \subset \cD(X_v;\uQ_p)$ of objects which
$v$-locally on $X$ lie in the smallest idempotent complete stable $\infty$-subcategory
generated by $\ud{\QQ}_p$ and $\O_X$. 
\end{defn}

\begin{rmk}\label{rmk:essential-image-dr} \cite[Rmk.~3.12]{ALB} ensures that $\cD_{\BBC}^{(b)}(X_v; \uQ_p)$ is equal to the essential image of the solution functor $\Sol_X\colon \Perf_\an(\FFF_X; \O) \to \cD(X_v; \uQ_p)$.
\end{rmk}

Now assume that $X=\Spd(C, \O_C)$ for an algebraically closed nonarchimedean field $C$ of mixed characteristic $(0, p)$. 

\begin{defn}\label{defn:absolute-bc} The category of \emph{Banach--Colmez spaces} on $\Spd(C, \O_C)$ is the minimal weak Serre subcategory $\BBC(C) \subset \mathrm{Shv}(\Spd(C, \O_C)_v; \uQ_p)$ containing $\ud{\QQ}_p$ and $\O_X$. We denote by $\cD^{b}_{\BBC}(\Spd(C, \O_C)_v; \uQ_p) \subset \cD(\Spd(C, \O_C)_v;\uQ_p)$ the full subcategory consisting of bounded objects whose cohomology sheaves lie in $\BBC(C)$. 
\end{defn}

\begin{rmk}
\label{rmk:relative-absolute-bc} 
A priori, the notation $\cD^{(b)}_{\BBC}(\Spd(C, \O_C)_v; \uQ_p)$ from  \cref{defn:relative-BC} 
conflicts with the notation $\cD^{b}_{\BBC}(X_v; \uQ_p)$ from \cref{defn:absolute-bc} when $X=\Spd(C, \O_C)$. 
However, there is no real conflict since these two categories coincide. Namely, \cref{rmk:essential-image-dr} ensures that it suffices to show the essential image of $\Sol_{\Spd(C, \O_C)} \colon \Perf_\an(\FFF_C; \O) \to \cD(\Spd(C, \O_C)_v; \uQ_p)$ coincides with $\cD^{b}_{\BBC}(\Spd(C, \O_C)_v; \uQ_p)$ in the sense of \cref{defn:absolute-bc}. 
Recall that \cite[\textsection 5.2, Cor.~6.10, Th.~7.1]{LB-thesis} construct a ``perverse'' $t$-structure on $\Perf_\an(\FFF_C; \O)$ 
which is bounded and
such that $\Sol$ is $t$-exact with respect to the perverse $t$-structure on $\Perf_\an(\FFF_C; \O)$ and the usual $t$-structure on $\cD(\Spd(C, \O_C)_v; \uQ_p)$ and, furthermore,
$\Sol$ identifies the heart of the perverse $t$-structure with $\BBC(C)$.
Combining this with the fact that $\Sol$ is fully faithful (see \cref{ff1}),
we conclude that the essential image of $\Sol \colon \Perf_\an(\FFF_C; \O) \to \cD(\Spd(C, \O_C)_v; \uQ_p)$
is equal to $\cD^{b}_{\BBC}(\Spd(C, \O_C)_v; \uQ_p)$. 
\end{rmk}

\begin{thm}\label{thm:cohomology-banach-colmez-spaces} Let $f\colon X \to Y$ be a smooth proper morphism of diamonds over $\Spd(\QQ_p, \ZZ_p)$ and let $\cE \in \Perf_v(X; \uQ_p)$. Then $\rR f_{v, *} \cE $ lies in $\cD^{(b)}_{\BBC}(Y_v; \uQ_p)$. In particular, if $\Q_p\subset C$ is an algebraically closed nonarchimedean field and $Y=\Spd (C, \O_C)$, then $\rR^i f_{v, *}\cE \in \BBC(C)$ for any $i\in \ZZ$.
\end{thm}
\begin{proof}
    We start by showing the first claim of the theorem. For this, \cref{prop:dr-rh} ensures that, for $\cP\coloneqq \RH_X(\cE)\in \Perf_\an(\FFF_X; \O)$, we have $\Sol_X(\cP) \simeq \cE$. Then \cref{lemma:dr-pushforward} implies that 
    \[
    \rR f_{v, *}\cE \simeq \rR f_{v, *} \Sol_X(\cP) \simeq \Sol_{Y}\bigl(\rR f_{\FFF, *} \cP\bigr).
    \]
    Therefore, \cref{rmk:essential-image-dr} ensures that it suffices to show that $\rR f_{\FFF, *} \cP$ lies in $\Perf_\an(\FFF_Y; \O)$. This follows automatically from \cref{lemma:proper-pushforward-FF-coincide}. 

    If $Y=\Spd(C, \O_C)$, the fact that $\rR^i f_{v, *}\cE$ lie in  $\BBC(C)$  for any $i\in \ZZ$ follows immediately from the previous discussion and \cref{rmk:relative-absolute-bc}.  
\end{proof}

As another application of \cref{lemma:proper-pushforward-FF-coincide}, we prove the Fargues--Fontaine version of the primitive comparison theorem. 

\begin{cor}[Primitive Comparison Theorem]\label{cor:primitive-comparison-FF} Let $f\colon X' \to X$ be a smooth proper morphism of diamonds over $\Spd(\QQ_p, \ZZ_p)$. Then the following diagram 
\[
\begin{tikzcd}
    \Perf_v(X'; \ud{\Z}_p) \arrow[d, "\rR f_{v, *}(\blank)"] \arrow[r,"\RH^{\ZZ_p}_{X'}"] & \Perf_\an(\FFF_{X'}; \O) \arrow[d, "\rR f_{\FFF, *}(\blank)"] \\
    \Perf_v(X; \ud{\Z}_p) \arrow[r, "\RH^{\ZZ_p}_X"] & \Perf_\an(\FFF_{X}; \O)
\end{tikzcd}
\]
commutes, i.e.~there is a canonical isomorphism of functors $\gamma_f\colon \RH^{\ZZ_p}_X\bigl(\rR f_{v, *}(\blank)\bigr) \xr{\sim} \rR f_{\FFF, *}\bigl(\RH^{\ZZ_p}_{X'}(\blank)\bigr)$. 
\end{cor}
\cref{rmk:PD-literal-analog-of-LRZ} guarantees that $\rR f_{v, *}(\blank)$ preserves perfect complexes of $\uZ_p$-modules.   

\begin{proof}
    We first define the transformation $\gamma'_f\colon \RH_X\bigl(\rR f_{v, *}(\blank)\bigr) \rightarrow
    \rR f_{\FFF, *}\bigl(\RH_{X'}(\blank)\bigr)$ as functors $\cD(X'_v; \uQ_p) \to \cD'(\FFF_X; \O)$. Namely, this is the transformation adjoint to the composition
    \[
    f^*_{\FFF} \circ \RH_X \circ \rR f_{v, *} \xr[\sim]{\alpha_f^{-1}\circ \rR f_{v,*}} \RH_{X'} \circ f_v^* \circ \rR f_{v, *} \xr{\RH_{X'}(\epsilon_f)} \RH_{X'},
    \]
    where $\alpha_f$ is the isomorphism from \cref{lemma:dr-pushforward} and $\epsilon_f$ is the counit of the $(f_v^*, \rR f_{v, *})$-adjunction. 

    Combining this transformation with the natural functor $\Perf_v(X', \uZ_p) \to \cD(X'_v; \uQ_p)$, we get the desired transformation $\gamma_f\colon \RH^{\ZZ_p}_X\bigl(\rR f_{v, *}(\blank)\bigr) \rightarrow \rR f_{\FFF, *}\bigl(\RH^{\ZZ_p}_{X'}(\blank)\bigr)$. 
    
    Now we show that $\gamma_f$ is an isomorphism. After unraveling all the definitions 
    (see \cref{defn:analytic-perfectoid-complexes-FF}, \cref{lemma:proper-pushforward-FF}, and \cref{constr:RH}), the question boils down to showing that the natural base change morphism
    \[
    (\rR f_{v, *} \cP)\otimes_{\uZ_p}^\rL \BB_I \to \rR f_{v, *}( \cP\otimes_{\uZ_p}^\rL \BB_I)
    \]
    is an isomorphism for $\cP\in \Perf_v(X'; \uZ_p)$ and $I=[1,1]$ or $I=[1, p]$. This follows immediately from \cref{Relative primitive comparison-2}. 
\end{proof}

\subsection{Geometric Duality}

The main goal of this subsection is to prove a geometric duality statement for a smooth proper morphism of diamonds over $\Spd(\QQ_p, \ZZ_p)$. We start with the definition of the trace morphism. 

\begin{definition}
\label{defn:FF-trace-maps}
Let $f \colon X' \to X$ be a smooth proper of equidimension $d$ morphism of diamonds over $\Spd(\QQ_p, \ZZ_p)$, and let $\cE\in \Perf_\an(\FFF_{X'}; \O)$ with the evaluation map $\mathrm{ev}_\cE \colon \cE^\vee \otimes^\rL \cE \to \O_{\FFF_{X'}}$.   
\begin{enumerate}[leftmargin=*,label=\upshape{(\arabic*)}]
\item\label{defn:FF-trace-maps-1} The \emph{trace map} of $f$ is the morphism 
    \[
        \tr_f^{\FFF} \colon \rR f_{\FFF, *} \O_{\FFF_{X'}}(d)[2d] \xr{\sim} \RH^{\ZZ_p}_X\bigl(\rR f_{v, *} \underline{\ZZ}_p(d)[2d]\bigr) \xrightarrow{\RH^{\ZZ_p}_X(\tr_f^{\ZZ_p})} \RH^{\Z_p}_X(\uZ_p) \simeq \O_{\FFF_{X}},
    \]
    where the first isomorphism comes from \cref{cor:primitive-comparison-FF} and $\tr_f$ comes from \cref{trace-maps}.
\item\label{defn:FF-trace-maps-2} The \emph{evaluation map} $e(f, \cE)  \colon \rR f_{\FFF, *} \bigl(\cE^{\vee}(d)[2d] \bigr) \otimes^\rL_{\O} \rR f_{\FFF, *}\bigl( \cE \bigr) \to \O_{\FFF_X}$ is the composition
    \[
    e(f, \cE) \colon \rR f_{\FFF, *} \bigl(\cE^{\vee}(d)[2d] \bigr) \otimes^\rL_{\O} \rR f_{\FFF, *}\bigl( \cE \bigr) \xr{\cup} \rR f_{\FFF, *}\bigl(\cE^{\vee}(d)[2d] \otimes^\rL_{\O} \cE \bigr) \xr{\rR f_{\FFF, *}(\mathrm{ev}_{\cE}(d)[2d])} \rR f_{\FFF, *}\O_{\FFF_{X'}}(d)[2d] \xr{\tr_f^\FFF} \O_{\FFF_X}.
    \]
\item\label{defn:FF-trace-maps-3} The \emph{duality morphism} 
    \[
    \PD^\FFF_f(\cE) \colon \rR f_{\FFF, *} \bigl(\cE^{\vee}(d)[2d] \bigr) \to \rR {\cHom}_{\FFF_X}\bigl(\rR f_{\FFF, *}\cE, \O_{\FFF_X}\bigr)
    \]
    is the map adjoint to $e(f, \cE)$ under the tensor-Hom adjunction. 
\end{enumerate}
\end{definition}

\begin{lemma}\label{lemma:trace-compatible} Let $f \colon X' \to X$ be a smooth proper of equidimension $d$ morphism of diamonds over $\Spd(\QQ_p, \ZZ_p)$. Then the following diagram commutes: 
\[
\begin{tikzcd}[column sep = 6em]
    \pi_X^*\bigl( \rR f_{\FFF, *} \O_{\FFF_{X'}}(d)[2d]\bigr)\arrow[d, "\wr"] \arrow[r, "\pi^*_X(\tr_f^\FFF)"] &  \pi^*_X(\O_{\FFF_X}) \arrow[d, "\wr"] \\
    \rR f_{v, *}\BB_{[1, p]}(d)[2d]\arrow[r, "\tr_f^{\BB_{[1, p]}}"] & \BB_{[1, p]},
\end{tikzcd}
\]
where $\pi^*_X$ comes from \cref{defn:analytic-perfectoid-complexes-FF}, $\tr_f^\FFF$ comes from \cref{defn:FF-trace-maps}, $\tr_f^{\BB_{[1,p]}}$ comes from \cref{trace-maps}, and the vertical isomorphisms come from \cref{lemma:proper-pushforward-FF}.
\end{lemma}
\begin{proof}
    First, \cref{rmk:pullback-of-RH} implies that $\pi_X^* \circ \RH_X^{\ZZ_p} \simeq \blank \otimes^\rL_{\uZ_p} \BB_{[1, p]} \colon \Perf_v(X; \uZ_p) \to \Perf_v(X; \BB_{[1, p]})$. Therefore, after unravelling \cref{defn:FF-trace-maps} and \cref{trace-maps}, the result follows from the following sequence of equalities
    \[
    \pi_X^*(\tr_f^{\FFF}) = \pi_X^*\bigl(\RH^{\Z_p}_X(\tr_f^{\ZZ_p})\bigr) \simeq \tr_f^{\Z_p} \otimes^\rL_{\uZ_p} \BB_{[1, p]} = \tr_f^{\BB_{[1, p]}}. \qedhere
    \]
\end{proof}

Before we prove the geometric duality on the Fargues--Fontaine curve, we need to record the following general result. 

\begin{lemma}
\label{lemma:dualizable-iso} 
Let $F\colon \cC \to \cD$ be a symmetric monoidal functor between closed symmetric monoidal 
$\infty$-categories with internal Hom functors $\ud{\Hom}_\cC$ and $\ud{\Hom}_\cD$.  
Let $C\in \cC$ be a \emph{dualizable} object and let $C'\in \cC$ be an arbitrary object. 
Then the natural map $F\ud{\Hom}_\cC(C, C') \to \ud{\Hom}_\cD(F(C), F(C'))$ is an isomorphism.
\end{lemma}
\begin{proof}
It suffices to show that the morphism is an isomorphism on the level of underlying homotopy categories. 
In this case, the result follows immediately from \cite[Prop.~3.2]{FHM03}. 
\end{proof}

Now we are ready to prove the duality on the Fargues--Fontaine curve. 

\begin{thm}\label{thm:geometric-duality-curve} Let $f \colon X' \to X$ be a smooth proper of equidimension $d$ morphism of diamonds over $\Spd(\QQ_p, \ZZ_p)$ and let $\cE\in \Perf_\an(\FFF_{X'}; \O)$. Then the duality morphism 
\[
    \PD^\FFF_f(\cE) \colon \rR f_{\FFF, *} \bigl(\cE^{\vee}(d)[2d]\bigr) \to \rR {\cHom}_{\FFF_X}\bigl(\rR f_{\FFF, *}\cE, \O_{\FFF_X}\bigr)
\]
is an isomorphism.
\end{thm}
\begin{proof}
    By construction, the functor $\pi^*_X\colon \Perf_\an(\FFF_X; \O) \to \Perf_v(X; \BB_{[1, p]})$ is symmetric monoidal and conservative. Therefore, it suffices to show that $\pi^*_X(\PD^\FFF_f(\cE))$ is an isomorphism. Then \cref{lemma:trace-compatible}, \cref{lemma:dualizable-iso}, and \cref{lemma:proper-pushforward-FF} imply that it suffices to show that the morphism
    \[
    \PD^{\BB_{[1,p]}}_f\bigl(\pi^*_{X'}(\cE)\bigr) \colon \rR f_{v, *}\bigl(\pi^*_{X'} \cE^\vee(d)[2d] \bigr) \to \rR{\cHom}\bigl(\rR f_{v, *}(\pi^*_{X'}\cE), \BB_{[1,p]} \bigr)
    \]
    is an isomorphism, where $\PD^{\BB_{[1,p]}}_f\bigl(\pi^*_{X'}(\cE)\bigr)$ is the (sheafy) duality morphism from \cref{construction:evaluation-coevaluation-poincare-duality-non-solid}. Therefore, the result follows immediately from \cref{cor:PD-for-proper-smooth-pushforward}. 
\end{proof}

Now we wish to also get a version of the geometric duality for $\uQ_p$-local systems.

\begin{definition}
\label{defn:Qp-duality-map}
Let $f \colon X' \to X$ be a smooth proper of equidimension $d$ morphism of diamonds over $\Spd(\QQ_p, \ZZ_p)$, and let $\cE\in \Perf_v(X'; \uQ_p)$ with the evaluation map $\mathrm{ev}_\cE \colon \cE^\vee \otimes^\rL \cE \to \uQ_p$.   
\begin{enumerate}[leftmargin=*,label=\upshape{(\arabic*)}]
    \item\label{defn:Qp-duality-map-1} The \emph{trace map} of $f$ is the morphism $\tr_f^{\QQ_p} \colon \rR f_{v, *} \uQ_p(d)[2d] \to \uQ_p$ defined in \cref{trace-maps}~\cref{trace-maps-Qp}.
    \item\label{defn:Qp-duality-map-2} The \emph{evaluation map} $e(f, \cE)  \colon \rR f_{v, *} \bigl(\cE^{\vee}(d)[2d] \bigr) \otimes^\rL_{\uQ_p} \rR f_{v, *} \cE \to \uQ_p$ is the composition
    \[
    e(f, \cE) \colon \rR f_{v, *} \bigl(\cE^{\vee}(d)[2d] \bigr) \otimes^\rL_{\uQ_p} \rR f_{v, *} \cE \xr{\cup} \rR f_{v, *}\bigl(\cE^{\vee}(d)[2d] \otimes^\rL_{\uQ_p} \cE \bigr) \xr{\rR f_{v, *}(\mathrm{ev}_{\cE}(d)[2d])} \rR f_{v, *}\uQ_p(d)[2d] \xr{\tr_f^{\Q_p}} \uQ_p.
    \]
    \item\label{defn:Qp-duality-map-3} The \emph{duality morphism} 
    \[
    \PD^{\Q_p}_f(\cE) \colon \rR f_{v, *} \bigl(\cE^{\vee}(d)[2d] \bigr) \to \rR {\cHom}_{\cD(X_v; \uQ_p)}\bigl(\rR f_{v, *}\cE, \uQ_p\bigr)
    \]
    is the map adjoint to $e(f, \cE)$ under the tensor-Hom adjunction. 
\end{enumerate}
\end{definition}

Our first goal is to show that compatibility of this trace map with the trace on the Fargues--Fontaine curve. 

\begin{lemma}\label{lemma:trace-qp-ff} Let $f \colon X' \to X$ be a smooth proper of equidimension $d$ morphism of diamonds over $\Spd(\QQ_p, \ZZ_p)$. Then the following diagram 
\[
\begin{tikzcd}[column sep = 6em]
\rR f_{v, *}\uQ_p(d)[2d] \arrow[r, "\tr_f^{\Q_p}"] \arrow[d, "\rR f_{v, *}(\can_{\uQ_p(d)[2d]})", "\wr"']
    & \uQ_p \arrow[dd, "\can_{\uQ_p}", "\wr"'] \\ 
    \rR f_{v, *}\bigl(\Sol_{X'}(\O_{\FFF_{X'}}(d)[2d])\bigr) \arrow[d, "(\alpha'_f)^{-1}", "\wr"'] & \\
     \Sol_X\bigl(\rR f_{\FFF, *} \O_{\FFF_{X'}}(d)[2d]\bigr) \arrow[r, "\Sol_X(\tr_f^\FFF)"] & \Sol_X(\O_{\FFF_X})
\end{tikzcd}
\]
is commutative. Here, $\alpha'_f$ is the isomorphism from \cref{lemma:dr-pushforward}, 
while $\can_{\uQ_p}$ and $\can_{\uQ_p(d)[2d]}$ are the isomorphisms from \cref{prop:dr-rh}. 
\end{lemma}
\begin{proof}
    Recall that $\tr_f^\FFF = \RH^{\Z_p}_X(\tr_f^{\Z_p})$ while $\tr_f^{\QQ_p} = \tr_f^{\ZZ_p}[\frac{1}{p}]$. Furthermore, \cref{prop:dr-rh} implies that we have a canonical identification of functors 
    \[
    \blank \otimes^\rL_{\Z_p} \uQ_p \xr[\can]{\sim}\Sol_X\circ \RH_X^{\Z_p} \colon \Perf_v(X; \uZ_p) \to \cD(X_v; \uQ_p).
    \]
    Hence, the result follows from the following sequence of identifications:
    \[
    \Sol_X(\tr_f^\FFF) = \Sol_X \circ \RH^{\Z_p}_X(\tr_f^{\Z_p}) = \tr_f^{\Z_p}\bigl[\frac{1}{p}\bigr] = \tr_f^{\Q_p} \qedhere
    \] 
\end{proof}

Let $f \colon X' \to X$ be a smooth proper of equidimension $d$ morphism of diamonds over $\Spd(\QQ_p, \ZZ_p)$ and let $\cE\in \Perf_v(X'; \uQ_p)$ be a perfect complex. 

\begin{thm}[Geometric duality]
\label{geometric-duality-Qp}
Let $f \colon X' \to X$ be a smooth proper of equidimension $d$ morphism of diamonds over $\Spd(\QQ_p, \ZZ_p)$ and let $\cE\in \Perf_v(X'; \uQ_p)$ be a perfect complex. Then the duality morphism 
\[
    \PD^{\Q_p}_f(\cE) \colon \rR f_{v, *}\bigl(\cE^{\vee}(d)[2d]\bigr) \to \rR{\cHom}_{\cD(X_v; \uQ_p)}\bigl(\rR f_{v, *}\cE, \ud{\QQ}_p \bigr)
\]
is an isomorphism.
\end{thm}
\begin{proof}
    First, \cref{prop:dr-rh} implies that there exists $\cP\in \Perf_\an(\FFF_{X'}; \O)$ such that $\Sol_{X'}(\cP)\simeq \cE$. Then the combination of \cref{ff1}, \cref{lemma:dr-pushforward}, and \cref{lemma:trace-qp-ff} implies that $\PD_f^{\Q_p}(\cE) = \Sol_X(\PD_f^{\FFF}(\cP))$. Therefore, it suffices to show that $\PD_f^{\FFF}(\cP)$ is an isomorphism. This follows automatically from \cref{thm:geometric-duality-curve}. 
\end{proof}

\subsection{Cohomological dimension}\label{section:cohomological-duality-geometric}

The main goal of this section is to show that for a smooth proper of dimension $d$ morphism $f\colon X \to Y$ of diamonds over $\Spd(\QQ_p, \ZZ_p)$, the derived pushforward has cohomological dimension $2d$ when restricted to perfect $\uQ_p$-complexes.

To that end, we first study the case where $Y=\Spa(C, \O_C)$
for an algebraically closed nonarchimedean field of mixed characteristic $(0, p)$. 
We begin by collecting some basic properties of Banach--Colmez spaces over $\Spa(C, \O_C)$.

\begin{lemma}\label{lemma:properties-BC} Let $C$ be an algebraically closed nonarchimedean field of mixed characteristic $(0, p)$ and let $V\in \BBC(C)$. Then the following statements hold:
\begin{enumerate}[leftmargin=*,label={\upshape{(\roman*)}}]
    \item\label{lemma:properties-BC-1} the complex $\rR{\cHom}_{\cD\bigl(\Spd(C, \O_C)_v; \uQ_p\bigr)}(V, \uQ_p)$ is concentrated in degrees $0$ and $1$;
    \item\label{lemma:properties-BC-2} $V\simeq \ud{\QQ}_p^n$ for some $n\in\ZZ_{\geq 0}$ if and only if $\dim_{\QQ_p} V\bigl(\Spa(C, \O_C)\bigr)<\infty$.
\end{enumerate}
\end{lemma}
\begin{proof}
    Recall that, by \cite[Prop.~7.11]{LB-thesis} and \cite[p.~45]{Colmez}, $V$
    fits into a diagram of exact sequences
    \begin{equation}\label{eqn:representation-BC}
    \begin{tikzcd}
    & & 0\arrow{d} & & \\
    & & \uQ_p^d \arrow{d} & &\\
    0 \arrow{r} & \uQ_p^{d'} \arrow{r} & V' \arrow{d} \arrow{r} & \GG_a^{d''} \arrow{r} & 0. \\
    & & V \arrow{d} & & \\
    & & 0 & &
    \end{tikzcd}
    \end{equation}

    First, we address \cref{lemma:properties-BC-1}. Note that $\rR{\cHom}_{\cD\bigl(\Spa(C, \O_C)_v; \uQ_p\bigr)}(\uQ_p, \uQ_p) \simeq \uQ_p[0]$. Therefore, Diagram~\cref{eqn:representation-BC} reduces the question to the case when $V=V'$ is an extension of $\GG_a^{d''}$ by $\uQ_p^{d'}$.
    For this, it suffices to separately treat the cases of $V=\uQ_p$ and $V=\GG_a$. The former case was already discussed above, so we only need to deal with $V=\GG_a$. This case follows immediately from \cite[Th.~3.8]{ALB} which implies that $\rR{\cHom}_{\cD\bigl(\Spa(C, \O_C)_v; \uQ_p\bigr)}(\GG_a, \uQ_p)\cong \GG_a[-1]$.

    Now we deal with \cref{lemma:properties-BC-2}. The ``only if'' direction of the claim is evident. Therefore, we only need to show that $V\simeq \uQ_p^n$ if $\dim_{\QQ_p} V\bigl(\Spa(C, \O_C)\bigr)<\infty$. Using Diagram~\cref{eqn:representation-BC} and the fact that $\rR\Gamma_v\bigl(\Spa(C, \O_C)_v, \uQ_p\bigr) \simeq \Q_p[0]$, we see that $\dim_{\QQ_p} V'\bigl(\Spa(C, \O_C)\bigr)<\infty$ as well. Using that $\rR\Gamma_v\bigl(\Spa(C, \O_C), \uQ_p\bigr) \simeq \QQ_p[0]$ once again, we deduce that $\Hh^0_v\bigl(\Spa(C, \O_C), \GG_a^{d''}\bigr) \simeq C^{d''}$ is a finite-dimensional $\QQ_p$-vector space. This is possible only if $d''=0$. Hence, $V'\simeq \uQ_p^{d'}$, so we conclude that $V\simeq \mathrm{coker}(\uQ_p^{d} \hookrightarrow \uQ_p^{d'})\simeq \uQ_p^{d'-d}$ as desired.
\end{proof}

\begin{rmk}\label{rmk:trivial-BC} \cref{lemma:properties-BC}~\cref{lemma:properties-BC-2} implies that a Banach--Colmez space $V$ is isomorphic to $0$ if and only if $V\bigl(\Spa(C, \O_C)\bigr)=0$; the forward direction is trivial, while the converse holds because the lemma first forces $V$ to 
be of the form $\uQ_p^n$, then the condition forces $n=0$.
\end{rmk}

\begin{lemma}\label{lemma:finite-coh-in-deg-0} Let $C$ be an algebraically closed nonarchimedean field of mixed characteristic $(0, p)$, let $X^\sharp$ be a qcqs rigid-analytic space over $C$, let $X\coloneqq (X^\sharp)^\diam$ be the associated diamond, and let $\bL$ be a $\uQ_p$-local system on $X$. Then $\dim_{\QQ_p} \Hh^0_v(X, \bL)<\infty$. 
\end{lemma}
\begin{proof}
    A combination of \cref{lemma:lattice-etale-locally-lsd}, \cref{cor:local-system-tor-ampl-0}, and \cite[Lem.~15.6]{diamonds} ensures that we can find an \'etale covering $X'^\sharp \to X^\sharp$ such that $\bL$ admits a $p$-torsionfree $\uZ_p$-lattice $\bL^+$ on $X'\coloneqq (X'^\sharp)^\diam$. Without loss of generality, we can assume that $X'$ is qcqs as well. Since $X' \to X$ is a $v$-covering, we conclude that $\Hh^0_v(X, \bL) \to \Hh^0_v(X', \bL|_{X'})$ is injective. Thus, it suffices to prove the claim for $X'$. In other words, we can assume that $\bL$ admits a $p$-torsionfree $\uZ_p$-lattice $\bL^+$. 
    
    In this case, it suffices to show that $\Hh^0_v(X, \bL^+)$ is a finitely generated $\ZZ_p$-module. Now the fact that we can check derived completeness $v$-locally and \cref{cor:perfect-ell-complete} implies that $\bL^+$ is derived $p$-adically complete. Then \cite[\href{https://stacks.math.columbia.edu/tag/091U}{Tag~091U}]{stacks-project} and \cite[\href{https://stacks.math.columbia.edu/tag/0BLX}{Tag 0BLX}]{stacks-project} imply that $\Hh^0_v(X, \bL^+)$ is derived $p$-adically complete. Furthermore, it is $p$-torsionfree since $\bL^+$ is $p$-torsionfree. Hence, \cite[\href{https://stacks.math.columbia.edu/tag/09BA}{Tag 09BA}]{stacks-project} ensures that it suffices to show that $\Hh^0_v(X, \bL^+)/p$ is finite. Since the natural map $\Hh^0_v(X, \bL^+)/p \to \Hh^0_v(X, \bL^+/p)$ is injective, we reduce the question to showing that $\Hh^0_v(X, \bL^+/p)$ is finite. 

    Now note that $\bL^+/p$ is an $\mathbf{F}_p$-local system. Thus, we can find a finite \'etale covering $X'^\sharp \to X^\sharp$ such that $\bL^+/p \simeq \ud{\mathbf{F}}_p^n$ over $X'\coloneqq (X'^\sharp)^\diam$. Once again, the natural morphism $\Hh^0_v(X, \bL^+/p) \to \Hh^0_v(X', \bL^+/p|_{X'})$ is injective, so we can assume that $\bL^+/p\simeq \ud{\FF}_p^n$. In this case, we have $\Hh^0_v(X, \bL^+/p) \simeq (\FF_p)^{\# \pi_0(X)\cdot n}$, which is finite due to \cite[Cor.~2.3]{adic-notes}.
\end{proof}

\begin{lemma}\label{lemma:no-BC-in-degree-0} Let $C$ be an algebraically closed nonarchimedean field of mixed characteristic $(0, p)$, let $f\colon X \to \Spd(C, \O_C)$ be a smooth proper morphism of diamonds, and let $\mathbf{V}$ be a $\uQ_p$-local system on $X$. Then $f_{v, *}\mathbf{V}\simeq \uQ_p^n$ for some $n$. 
\end{lemma}
\begin{proof}
    This follows immediately from the combination of \cref{lemma:finite-coh-in-deg-0}, \cref{thm:cohomology-banach-colmez-spaces}, and \cref{lemma:properties-BC}~\cref{lemma:properties-BC-2}.
\end{proof}

\begin{thm}\label{thm:cohomological-dimension-Qp-point} Let $C$ be an algebraically closed nonarchimedean field of mixed characteristic $(0, p)$, let $f\colon X \to \Spd(C, \O_C)$ be a smooth proper of dimension $d$ morphism of diamonds, and let $\cE \in \Perf_v(X; \uQ_p) \cap \cD^{[r, r']}(X_v; \uQ_p)$. Then $\rR f_{v, *}\cE \in \cD^{[r, r'+2d]}_{\BBC}\bigl(\Spa(C, \O_C)_v; \uQ_p\bigr)$.
\end{thm}
\begin{proof}
    We already know that $\rR f_{v, *} \cE$ lies in $\cD^b_{\BBC}\bigl(\Spa(C, \O_C)_v; \uQ_p\bigr)$ due to \cref{thm:cohomology-banach-colmez-spaces}, so we only need to show that this complex is concentrated in degrees $[r, r'+2d]$. \Cref{cor:local-system-tor-ampl}, \cref{cor:local-system-tor-ampl-0}, and \cref{lemma:characterizing-perfect-complexes} imply that it suffices to consider the case when $\cE =\bL[0]$ is a $\uQ_p$-local system concentrated in degree $0$. In this case, we need to show that $\rR f_{v, *} \bL$ is concentrated in degrees $[0, 2d]$. Since $\rR f_{v, *}$ sends $\cD^{\geq 0}$ to $\cD^{\geq 0}$, the question is equivalent to showing that $\rR f_{v, *} \bL(d)[2d]$ is concentrated in non-positive degrees. 

    Furthermore, since $X$ is a finite disjoint union of smooth proper clopen subdiamonds of equidimension $d'$ for $d' \leq d$, we can assume that $X$ is of equidimension $d$. In this case, \cref{geometric-duality-Qp} implies that $\rR f_{v, *}\bigl(\bL(d)[2d]\bigr) \simeq \rR{\cHom}_{\cD\bigl(\Spd(C, \O_C)_v; \uQ_p\bigr)}\bigl(\rR f_{v, *}(\bL^{\vee}), \ud{\QQ}_p \bigr)$, where 
    \[
    \bL^{\vee}=\rR{\cHom}_{\cD(X_v; \uQ_p)}(\bL, \uQ_p)\simeq {\cHom}_{\uQ_p}(\bL, \uQ_p)
    \]
    is concentrated in degree $0$. Now  \cref{thm:cohomology-banach-colmez-spaces} implies that $\rR f_{v, *}(\bL^{\vee})$ lies in $\cD^{[0,c]}_{\BBC}\bigl(\Spa(C, \O_C)_v,\uQ_p\bigr)$ for some integer $c\geq 0$ and, furthermore, \cref{lemma:no-BC-in-degree-0} ensures that $f_{v, *}(\bL^{\vee}) \simeq \uQ_p^n$ for some integer $n$. 

    Then \cref{lemma:properties-BC}~\cref{lemma:properties-BC-1}, the observation that $\rR{\cHom}_{\cD\bigl(\Spa(C, \O_C)_v; \uQ_p\bigr)}(\uQ_p, \uQ_p) \simeq  \ud{\QQ}_p[0]$, and the 
    spectral sequence associated with the Postnikov filtration on $\rR f_{v, *}(\bL^{\vee})$ imply that 
    \[
    \rR f_{v, *}\bigl(\bL(d)[2d]\bigr) \simeq \rR{\cHom}_{\cD\bigl(\Spd(C, \O_C)_v; \uQ_p\bigr)}\bigl(\rR f_{v, *}(\bL^{\vee}), \ud{\QQ}_p \bigr) \in \cD^{[-c, 0]}\bigl(\Spa(C, \O_C)_v; \uQ_p\bigr).
    \]
    This finishes the proof.
\end{proof}

Now we deal with the relative situation. We start with the following lemma.

\begin{lemma}\label{lemma:inject-on} Let $S$ be an affinoid perfectoid space over $\Spd(\QQ_p, \ZZ_p)$ and let $\cE$ be a vector bundle on $\FFF_S$. Then there are positive integers $n,m>0$ and a short exact sequence
\[
0 \to \cE \to \O(n)^m \to \cQ \to 0
\]
such that $\cQ$ is a vector bundle on $\FFF_S$.
\end{lemma}
\begin{proof}
    First, \cite[Th.~II.2.6]{Fargues-Scholze} implies that there exist integer $n,m >0$ and a surjection $\O(-n)^m \xtwoheadrightarrow{\varphi} \cE^{\vee}$. Note that $\cP\coloneqq \ker \varphi$ is also a vector bundle as a kernel of a surjective map of vector bundles. Then after dualizing $\varphi$, we get a short exact sequence $0 \to \cE \to \O(n)^m \to \cQ\to 0$, where $\cQ\simeq \cP^{\vee}$ is a vector bundle.
\end{proof}

\begin{lemma}\label{lemma:fiberwise-surjectivity} Let $S$ be an affinoid perfectoid space over $\Spd(\QQ_p, \ZZ_p)$ and let $u\colon \cE \to \cE'$ be a morphism of vector bundles on $\FFF_S$. Suppose that, for any algebraically closed nonarchimedean field $C$ of mixed characteristic $(0, p)$ and any morphism $s\colon \Spd(C, \O_C) \to S$ over $\Spd(\QQ_p, \ZZ_p)$, the natural morphism
\begin{equation}\label{eqn:fiberwise-surj}
    s^*\cal{H}^1\bigl(\Sol_S(\cE)\bigr) \to s^*\cal{H}^1\bigl(\Sol_S(\cE')\bigr)
\end{equation}
is surjective. Then $\cal{H}^1\bigl(\Sol_S(\cE)\bigr) \to \cal{H}^1\bigl(\Sol_S(\cE')\bigr)$ is surjective as well.
\end{lemma}
\begin{proof}
Throughout this proof, we will freely use the fact that $s^*\cal{H}^1(\Sol_S(\cE)) \simeq \cal{H}^1(\Sol_{\Spd(C, \O_C)}(\cE|_{\FFF_{C}}))$. This follows from the explicit description of the solution functor in \cref{rmk:explicit-adjoint}. 

First, we choose a short exact sequence $0 \to \cE \xr{i} \O(n)^m \to \cQ \to 0$ as in \cref{lemma:inject-on}. Let $\F \coloneqq \O(n)^m \oplus_{\cE} \cE'$ be the pushout of $\O(n)^m$ along the map $u\colon \cE \to \cE'$. In other words, $\F \simeq \coker(\cE \xr{i\oplus u} \O(n)^m\oplus \cE')$. As a pushout, it fits into the following exact sequence:
\[
0 \to \cE' \to \F \to \cQ \to 0
\]
Since both $\cE'$ and $\cQ$ are vector bundles, $\F$ must also be a vector bundle. On the other hand, since $i\colon \cE \to \O(n)^m$ is injective, the sheaf $\F$ fits into another exact sequence
\[
0 \to \cE \xr{i\oplus u} \O(n)^m \oplus \cE' \to \F \to 0.
\]
Now note that \cref{rmk:explicit-adjoint} and \cite[Prop.~II.2.5(iii)]{Fargues-Scholze} imply that $\cal{H}^1\bigl(\Sol_S(\O(n)^m)\bigr)=0$ since $n>0$. Furthermore, \cref{lemma:coh-dimension-sol}~\cref{lemma:coh-dimension-sol-1} implies that $\cal{H}^2\bigl(\Sol_S(\cE)\bigr)=0$. Thus, we conclude that the map
\[
\cal{H}^1\bigl(\Sol_S(\cE)\bigr) \to \cal{H}^1\bigl(\Sol_S(\cE')\bigr)
\]
is surjective if and only if $\cal{H}^1\bigl(\Sol_S(\F)\bigr)=0$. In particular, the surjectivity of \cref{eqn:fiberwise-surj} implies that 
\[
s^*\cal{H}^1\bigl(\Sol_S(\F)\bigr) \simeq \cal{H}^1\bigl(\Sol_{\Spd(C, \O_C)}(\F|_{\FFF_C})\bigr)=0.
\]
The classification of vector bundles on $\FFF_C$ shows that $\F|_{\FFF_C}\simeq \oplus_{\lambda\in \QQ} \O(\lambda)^{n_\lambda}$.
For $\lambda<0$, \cite[Prop.~II.2.5(i)]{Fargues-Scholze} implies that  $\cal{H}^1\bigl(\Sol_{\Spd(C, \O_C)}(\O(\lambda))\bigr)[-1]$ is isomorphic to $\Sol_{\Spd(C, \O_C)}\bigl(\O(\lambda)\bigr)$. In particular, it is never equal to zero (for the lack of a better reference, this follows from \cref{ff1}), hence $\F|_{\FFF_C}\simeq \oplus_{\lambda\in \QQ_{\geq 0}} \O(\lambda)^{n_\lambda}$ for any morphism $s\colon \Spd(C, \O_C) \to S$, where $C$ is an algebraically closed nonarchimedean field of mixed characteristic $(0, p)$. 

This implies that all fiberwise Harder--Narasimhan slopes of $\F$ are non-negative. Now we wish to show that this already implies that $\cal{H}^1\bigl(\Sol_S(\F)\bigr)=0$. Using \cref{rmk:explicit-adjoint}, we see that this sheaf is equal to the $v$-sheafification of the presheaf
\[
(S' \to S) \mapsto \Hh^1_\an\bigl(\FFF_{S'}, \F|_{\FFF_{S'}}\bigr)
\]
Therefore, \cite[Prop.~II.3.4(ii)]{Fargues-Scholze} implies that  $\cal{H}^1\bigl(\Sol_S(\F)\bigr)=0$.
\end{proof}

For the next definition, we fix a diamond $X$ over $\Spd(\QQ_p, \ZZ_p)$ and $r,r'\in \ZZ \cup \{-\infty, \infty\}$. 

\begin{defn}\label{defn:finite-tor-amplitude-ff-curve} The \emph{$\infty$-category $\Perf^{[r,r']}_\an(\FFF_X; \O)$ of perfect complexes of tor-amplitude $[r,r']$} is the full $\infty$-subcategory of $\Perf_\an(\FFF_X; \O)$ spanned by those objects $\cP$ for which $\pi_X^*\cP\in \Perf^{[r,r']}_v(X; \BB_{[1,p]})$.
\end{defn}

First, we show that this recovers the usual definition of tor-amplitude on a ringed site (see \cite[\href{https://stacks.math.columbia.edu/tag/08CG}{Tag 08CG}]{stacks-project}) when $X=S$ is an affinoid perfectoid space. We start with the following very general lemma. 

\begin{lemma}\label{lemma:perfect-ring-adic-space} Let $S=\Spa(R, R^+)$ be an affinoid adic space and let $r,r'\in \ZZ \cup \{-\infty, \infty\}$. Then $\widetilde{(-)} \colon \Perf(R) \to \Perf(S; \O)$ induces an equivalence $\widetilde{(-)}\colon \Perf^{[r, r']}(R) \to \Perf^{[r, r']}(S; \O)$. 
\end{lemma}

Here, $\Perf^{[r, r']}(S; \O)$ denotes the $\infty$-category of perfect complexes with tor-amplitude in $[r, r']$ in the sense of \cite[\href{https://stacks.math.columbia.edu/tag/08CG}{Tag 08CG}]{stacks-project}. 
\begin{proof}
    First, \cite[Th.~1.4]{Andreychev} and \cref{prop:perfect-complexes-sheafification} imply that $\widetilde{(-)}\colon \Perf(R) \to \Perf(S; \O)$ is an equivalence. Hence, we only need to show that $M$ has tor-amplitude in $[r,r']$ if and only if $\widetilde{M}$ has tor-amplitude in $[r, r']$. 

    If $M$ has tor-amplitude in $[r, r']$, then \cite[\href{https://stacks.math.columbia.edu/tag/0658}{Tag 0658}]{stacks-project} ensures that $M$ is quasi-isomorphic to a finite complex of finite projective modules $K^\bullet$ such that $K^i=0$ for $i\notin [r, r']$. Then $\widetilde{M}\simeq \widetilde{K}^\bullet$ is a finite complex of vector bundles such that $\widetilde{K}^i=0$ for $i\notin[r,r']$. Thus, $\widetilde{M}$ has tor-amplitude in $[r, r']$. 

    Now suppose that $\widetilde{M}$ has tor-amplitude in $[r, r']$. Then \cite[\href{https://stacks.math.columbia.edu/tag/068V}{Tag 068V}]{stacks-project} ensures that it suffices to show that $M\otimes^\rL_R k(\m)$ is concentrated in degrees $[r, r']$ for any maximal ideal $\m\subset R$. Choose a point $s\in S$ such that $\supp(s)=\m$; cf.\ \cite[Lem.~1.4]{Huber-generalization}. Then \cite[\href{https://stacks.math.columbia.edu/tag/09U9}{Tag 09U9}]{stacks-project} implies that $\widetilde{M}_s\simeq M\otimes^\rL_R \O_{S, s}\in \Perf^{[r, r']}(\O_{S, s})$. This implies that $M\otimes^\rL_{R} k(s)$ lies in $\cD^{[r, r']}\bigl(k(s)\bigr)$ for any $s\in S$. Then \cref{lemma:tor-amplitude-over-points} implies that $M\in \Perf^{[r, r']}(R)$ and finishes the proof. 
\end{proof}

\begin{lemma}\label{lemma:tor-amplitude-on-FFS} Let $S$ be an affinoid perfectoid space over $\Spd(\QQ_p, \ZZ_p)$, let $\cP\in \Perf_\an(\FFF_S; \O)$, and let $r,r'\in \ZZ\cup \{-\infty, \infty\}$. Then $\cP$ lies in $\Perf^{[r, r']}_\an(\FFF_S; \O)$ if and only if $\cP$ has tor-amplitude in $[r,r']$ in the sense of \cite[\href{https://stacks.math.columbia.edu/tag/08CG}{Tag 08CG}]{stacks-project}.
\end{lemma}

In other words, \cref{lemma:tor-amplitude-on-FFS} ensures that the two potentially different meanings of $\Perf^{[r, r']}(\FFF_S; \O)$ coincide.

\begin{proof}
    First, \cref{notation:sheaf-perfect-complexes-finite-amplitude} and \cref{thm:v-descent-rationally-BBI} imply that $\Perf^{[r,r']}_v(S; \BB_{[1,p]})\simeq \Perf^{[r, r']}\bigl(\BB_{[1,p]}(S)\bigr)$.
    Thus, $\Perf^{[r,r']}_v(S; \BB_{[1,p]})\simeq \Perf^{[r, r']}(Y_{S, [1,p]}; \O)$ by virtue of \cref{lemma:perfect-ring-adic-space}. 

    Note that the object $\pi_S^*\cP\in \Perf(Y_{S, [1,p]}; \O) \simeq \Perf_v(S; \BB_{[1,p]})$ can be identified with the pullback of $\cP$ along the map $Y_{S, [1,p]}\to \FFF_S$. Thus, the surjectivity of this map together with \cite[\href{https://stacks.math.columbia.edu/tag/09U9}{Tag~09U9}, \href{https://stacks.math.columbia.edu/tag/068V}{Tag~068V}]{stacks-project} imply that $\pi_S^*\cP \in \Perf^{[r, r']}(Y_{S, [1,p]}; \O) \simeq \Perf^{[r,r']}_v(S; \BB_{[1,p]})$ if and only if $\cP$ has tor-amplitude in $[r,r']$ in the sense of \cite[\href{https://stacks.math.columbia.edu/tag/08CG}{Tag 08CG}]{stacks-project}.
\end{proof}

As a corollary of these lemmas, we deduce the following result.

\begin{cor}\label{cor:representative-of-a-perfect-complex} Let $S$ be an affinoid perfectoid space over $\Spd(\QQ_p, \ZZ_p)$, let $r,r'\in \ZZ$, and let $\cP\in \Perf^{[r,r']}_\an(\FFF_S; \O)$. Then $\cP$ is isomorphic to a finite complex of vector bundles $K^\bullet$ such that $K^i\simeq 0$ when $i\notin[r, r']$.
\end{cor}
\begin{proof}
    First, \cite[Prop.~2.6]{ALB} implies that $\cP\simeq L^\bullet$, where $L^\bullet$ is a finite complex of vector bundles. Since $\cP\in \Perf^{[r, r']}_\an(\FFF_S; \O)$, we conclude that $\cP\simeq \tau^{\geq r}\tau^{\leq r'} L^\bullet\eqqcolon K^\bullet$. By construction, we have $K^i\simeq 0$ when $i\notin [r, r']$. Thus, we are only left to show that each $K^i$ is a vector bundle. For this, the proof of \cite[\href{https://stacks.math.columbia.edu/tag/08CI}{Tag 08CI}]{stacks-project} applies almost verbatim, having established \cref{lemma:tor-amplitude-on-FFS}.
\end{proof}

Finally, we are ready to prove the main theorem of this subsection. 

\begin{thm}
\label{thm:tor-amplitude-qp-smoothproper}
Let $f\colon X \to Y$ be a smooth proper of dimension $d$ morphism of diamonds over $\Spd(\QQ_p, \ZZ_p)$ and let $\cE\in \Perf^{[r,r']}_v(X; \uQ_p)$ for some $r, r'\in \ZZ\cup \{-\infty, \infty\}$. Then $\rR f_{v, *}\cE$ lies in $\cD^{[r, r'+2d]}_{\BBC}(Y_v; \uQ_p)$.
\end{thm}
\begin{proof}
First, \cref{thm:cohomology-banach-colmez-spaces} already implies that $\rR f_{v,*} \cE \in \cD^{(b)}_{\BBC}(Y_v; \uQ_p)$.
Thus, we only need to show that $\rR f_{v,*} \cE$ is concentrated in degrees $[r, r'+2d]$. 

    This question is $v$-local on $Y$, so we can assume that $Y$ is a strictly totally disconnected perfectoid space. In this case, the diamond $X$ is qcqs. Thus, we can assume that $r, r'\in \ZZ$. Now \cref{prop:dr-rh} and \cref{lemma:dr-pushforward} imply that 
    \[
    \rR f_{v, *} \cE \simeq \rR f_{v, *}\bigl( \Sol_X(\RH_X\cE) \bigr)\simeq \Sol_Y\bigl(\rR f_{\FFF, *} (\RH_X \cE)\bigr).
    \]
    By \cref{defn:finite-tor-amplitude-ff-curve} and \cref{constr:RH}, we have $\RH_X\cE\in \Perf_\an^{[r, r']}(\FFF_X; \O)$. Hence, \cref{defn:finite-tor-amplitude-ff-curve}, \cref{lemma:proper-pushforward-FF}, \cref{lemma:proper-pushforward-FF-coincide}, and \cref{cor:PD-tor-amplitude} ensure that 
    \[
    \rR f_{\FFF, *} (\RH_X \cE) \in \Perf^{[r, r'+2d]}_\an(\FFF_Y; \O).
    \]
    By virtue of \cref{cor:representative-of-a-perfect-complex}, we can then find a representative $K^\bullet \simeq \rR f_{\FFF, *} (\RH_X \cE)$ such that each $K^i$ is a vector bundle on $\FFF_Y$ and $K^i\simeq 0$ when $i\notin[r, r'+2d]$. 

    \cref{lemma:coh-dimension-sol} already implies that $\rR f_{v, *} \cE \simeq \Sol_Y(\rR f_{\FFF, *} \RH_X \cE) \simeq \Sol_Y(K^\bullet)$ is concentrated in degrees $[r, r'+2d+1]$. Thus, we only need to show that $\cal{H}^{r'+2d+1}\bigl(\Sol_Y(K^\bullet)\bigr)=0$. Now using \cref{lemma:coh-dimension-sol}~\cref{lemma:coh-dimension-sol-1} and the spectral sequence associated with the na\"{i}ve filtration on $K^\bullet$ (see \cite[\href{https://stacks.math.columbia.edu/tag/012N}{Tag 012N}]{stacks-project}), we conclude that 
    \[
    \cal{H}^{r'+2d+1}\bigl(\Sol_Y(K^\bullet)\bigr) \simeq \coker\Bigl( \cal{H}^1\bigl(\Sol_Y K^{r'+2d-1}\bigr) \to \cal{H}^1\bigl(\Sol_Y K^{r'+2d}\bigr)\Bigr).
    \]
    Therefore, the question reduces to showing that the natural map $\cal{H}^1\bigl(\Sol_Y K^{r'+2d-1}\bigr) \to \cal{H}^1\bigl(\Sol_Y K^{r'+2d}\bigr)$ is surjective. \cref{thm:cohomological-dimension-Qp-point} implies this map is surjective after a pullback along any morphism $s\colon \Spd(C, \O_C)\to Y$ for an algebraically closed field $C$ of mixed characteristic $(0, p)$. Thus, \cref{lemma:fiberwise-surjectivity} ensures that it is surjective on the nose. This finishes the proof. 
\end{proof}

\begin{corollary}
\label{cor:tor-amplitude-qp-smoothproper-2}
Let $f\colon X \to Y$ be a smooth proper of dimension $d$ morphism of diamonds over $\Spd(\QQ_p, \ZZ_p)$ and let $\cE\in \Perf_v(X; \uQ_p) \cap \cD^{[r, r']}(X_v; \uQ_p)$ for some $r, r'\in \ZZ\cup \{-\infty, \infty\}$. Then $\rR f_{v, *}\cE$ lies in $\cD^{[r, r'+2d]}_{\BBC}(Y_v; \uQ_p)$.
\end{corollary}
\begin{proof}
    This follows formally from \cref{thm:tor-amplitude-qp-smoothproper} similar to how \cref{cor:tor-amplitude-2} follows from \cref{cor:PD-tor-amplitude}.
\end{proof}

\section{Absolute duality}\label{section:absolute-duality}

In this section, we analyze the results of \cref{dualizability} and \cref{geom-arith-duality} in the case when the base is a point.
In this situation, one can make two improvements:
the finiteness and duality isomorphism, appropriately formulated, still hold for possibly singular proper rigid-analytic spaces, and one gets actual finiteness and an honest duality isomorphism in the arithmetic setting when the base $\Spa(K,\cO_K)$ is a $p$-adic local field $K$.
Our proofs of these improvements adapt the strategies of \cite[\S~7]{LRZ24} and \cite{ZL}, respectively.

Throughout the section, we fix a nonarchimedean extension $K$ of $\QQ_p$ (which will be a $p$-adic local field starting in \cref{arithmetic-duality}).
The \'etale, quasi-pro-\'etale or $v$-cohomology of a proper rigid-analytic space $X$ over $K$ is understood to be that of the associated diamond $X^\lozenge$ (keeping in mind \cite[Lem.~15.6]{diamonds}), though we will again drop the $(\blank)^\lozenge$ from the notation when talking about these cohomologies.

\subsection{Dualizing complexes: \texorpdfstring{$\ZZ_p$}{Zp}-coefficients}
We begin by recalling the dualizing complex in this setting.
\begin{construction}[{\cite[Pf.\ of Th.\ 3.36, Rmk.\ 3.37]{BH}}]\label{dualizing-construction}
Let $X$ be a rigid-analytic space over $K$ and let $n>0$ be an integer. Denote by $\lambda_v \colon \bigl(X^\lozenge_v, \underline{\ZZ/n}\bigr) \to \bigl(X^\lozenge_\et, \underline{\ZZ/n}\bigr) \simeq \bigl(X_\et, \underline{\ZZ/n}\bigr)$ the natural morphism of ringed sites. 
\begin{enumerate}[leftmargin=*]
    \item For each affinoid open $U=\Spa(A, A^\circ)$, there is an essentially unique potential dualizing complex $\omega^{\ZZ/n}_A \in \cD^+\bigl((\Spec A)_\et;\ZZ/n\bigr)$, together with pinnings $\rR\Gamma_{\bar{x}}(\Spec A,\omega^{\ZZ/n}_A) \xrightarrow{\sim} \ZZ/n(\delta(x))[2\delta(x)]$ for each geometric point $\bar{x} \to \Spec A$ lying over a point $x \in \Spec A$, which are compatible with immediate specializations, where $\delta(x)$ denotes the canonical dimension function $\delta(x) = \dim \overline{\{x\}}$ (see \cite[Th.~XVII.5.1.1]{deGabber} and \cite[Prop.~3.18]{BH});
    for every $m \mid n$, there is an essentially unique isomorphism $\omega^{\ZZ/n}_A \otimes_{\ud{\ZZ/n}}^\rL \underline{\ZZ/m} \simeq \omega^{\ZZ/m}_A$ thanks to the essential uniqueness of $\omega^{\ZZ/m}_A$ and $\rR\Gamma_x(\Spec A,\omega^{\ZZ/n}_A \otimes^\rL \underline{\ZZ/m}) \simeq \rR\Gamma_x(\Spec A,\omega^{\ZZ/n}_A) \otimes^\rL \underline{\ZZ/m}$.
    
    \item Using BBDG gluing and the fact that $\rR{\cHom}(\omega^{\ZZ/n}_A,\omega^{\ZZ/n}_A) \simeq \underline{\ZZ/n}$ is concentrated in nonnegative degrees, \cite[Th.~3.21]{BH} glues the analytifications $(\omega^{\ZZ/n}_A)^\an \in \cD^+\bigl(\Spa(A,A^\circ)_\et;\ZZ/n\bigr)$ to an essentially unique dualizing complex $\omega^{\ZZ/n}_X \in \cD^+(X_\et;\ZZ/n)$;
    using the essential uniqueness of the gluing and the compatibility of $\blank \otimes^\rL \underline{\ZZ/m}$ with pullback along affinoid opens $V \subseteq U$, we can again find essentially unique isomorphisms $\omega^{\ZZ/n}_X \otimes^\rL_{\ud{\ZZ/n}} \underline{\ZZ/m} \simeq \omega^{\ZZ/m}_X$ for any $m \mid n$.
    
    \item Finally, since $\rR\Hom^{<0}(\omega^{\ZZ/p^r}_X,\omega^{\ZZ/p^r}_X) \simeq 0$ and $\lambda^*_v \colon \cD^+(X_\et;\ZZ/p^r) \to \cD^+(X_v;\ZZ/p^r)$ is fully faithful \cite[Prop.~14.10]{diamonds}, the collection of objects $\lambda^*_v\omega^{\ZZ/p^r}_X \in \cD^+(X_v;\ZZ/p^r)$ and isomorphisms $\lambda^*_v\omega^{\ZZ/p^r}_X \otimes^\rL_{\ud{\ZZ/p^r}} \underline{\ZZ/p^s} \simeq \lambda^*_v\omega^{\ZZ/p^s}_X$
    for $r \geq s$ determines an essentially unique object 
    \[
    \omega^{\ZZ_p}_X \in \lim_r \cD^+(X_v;\ZZ/p^r) \simeq \cD^+_{p-\text{comp}}(X_v; \ud{\ZZ}_p)\subset \cD^+(X_v;\ud{\ZZ}_p),
    \]
    called the ($\ZZ_p$-)\emph{dualizing complex} of $X$;
    the equivalence above comes from \cite[Lem.\,3.5.7]{proetale} and \cref{lemma:replete}. 
    By construction, the object $\omega^{\ZZ_p}_X$ is derived $p$-adically complete and $\omega^{\ZZ_p}_X \otimes^\rL \underline{\ZZ/p^r} \simeq \lambda^*_v\omega^{\ZZ/p^r}_X$, so $\omega^{\ZZ_p}_X \simeq \rR\lim_r \lambda^*_v\omega^{\ZZ/p^r}_X$.
\end{enumerate}
\end{construction}
\begin{example}\label{dualizing-smooth-Zp}
    Assume that $X$ is a smooth rigid-analytic space over $\Spa(K,\cO_K)$ of equidimension $d$.
    Then there is a canonical isomorphism $\alpha^{\ZZ_p}_X \colon \underline{\ZZ}_p(d)[2d] \xrightarrow{\sim} \omega^{\ZZ_p}_X$:
    this follows from the existence of analogous isomorphisms $\alpha^{\ZZ/p^r}_X \colon \underline{\ZZ/p^r}(d)[2d] \xrightarrow{\sim} \omega^{\ZZ/p^r}_X$ \cite[Th.~3.21.(1)]{BH};
    as the $\alpha^{\ZZ/p^r}_X$ are constructed from essentially unique algebraic isomorphisms that are compatible with pinnings, they are compatible for varying $r$ and hence give rise to $\alpha^{\ZZ_p}_X$.

    More generally, if $f \colon X \to Y$ is a smooth morphism of equidimension $d$ between rigid-analytic spaces over $\Spa(K,\cO_K)$, there is a canonical isomorphism
    \[ \alpha^{\ZZ_p}_f \colon f^*_v \omega^{\ZZ_p}_Y(d)[2d] \xlongrightarrow{\sim} \omega^{\ZZ_p}_X, \] 
    which is induced in the same way from isomorphisms $\alpha^{\ZZ/p^r}_f \colon f^*_\et \omega^{\ZZ/p^r}_Y(d)[2d] \xrightarrow{\sim} \omega^{\ZZ/p^r}_X$ as in \cite[Cor.~7.1.2]{LRZ24}.
\end{example}
To every proper morphism $f \colon X \to Y$ of rigid-analytic spaces over $K$ and every integer $n > 0$, \cite[Th.~7.4.1, Th.~7.3.14]{LRZ24} attaches a proper trace map $\overline{\Tr}^{\ZZ/n}_f \colon \rR f_{\et,*}\omega^{\ZZ/n}_X \to \omega^{\ZZ/n}_Y$ which is essentially unique for the properties listed in \cite[Th.~7.4.1]{LRZ24}.
Our next goal is to construct similar morphisms for $\ZZ_p$-dualizing complexes.
To do so, we first need to compare the \'etale and $v$-pushforwards of the $\ZZ/n$-dualizing complex $\omega^{\ZZ/n}_X$ of a proper rigid-analytic space $X$.
While $\omega^{\ZZ/n}_X$ is in general not perfect, it follows directly from the construction in \cite[Th.~3.21]{BH} that it is still contained in the full subcategory $\cD^{(b)}_\zc(X;\ZZ/n) \subset \cD^+(X_\et;\ZZ/n)$ of locally bounded complexes with Zariski-constructible cohomology sheaves in the sense of \cite[Def.~3.1]{BH}.
Thus, its pushforwards can be compared using the following lemma.
\begin{lemma}\label{zc-cohomology-et-v}
    Let $f \colon X \to Y$ be a proper morphism of rigid-analytic spaces over $\Spa(K,\cO_K)$.
    Let $n>0$ be an integer and $\cF \in \cD^{(b)}_\zc(X;\ZZ/n)$.
    Then the natural base change map $c_{f, \cF} \colon \lambda^*_v\rR f_{\et,*} \cF \to \rR f_{v,*}\lambda^*_v \cF$ is an isomorphism.
\end{lemma}
\begin{proof}
    Since the statement is local on $Y$, we may assume that $Y$ (and hence $X$) are quasi-compact.
    In this case, $\cF \in \cD^b_\zc(X;\ZZ/n)$ is globally bounded.
    Since the full subcategory of those objects $\cF \in \cD^b_\zc(X;\ZZ/n)$ for which $c_{f,\cF}$ is an isomorphism is thick, we may assume that $\cF = g_{\et,*} \underline{\ZZ/m}$ for a finite morphism $g \colon W \to X$ and an integer $m\mid n$;
    cf.\ \cite[Prop.~3.6]{BH}.
    Moreover, by virtue of \cite[Prop.~2.6.3]{Huber-etale} and \cite[Lem.~15.6]{diamonds}, we have $\cF \simeq \rR g_{\et,*} \underline{\ZZ/m}$.
    On the other hand, the base change map $c$ is compatible with compositions in $f$, i.e., the diagram
    \[ \begin{tikzcd}
        \lambda^*_v \rR f_{\et,*} \cF \arrow[d, "\wr"] \arrow[r,"{c_{f, \cF}}"] & \rR f_{v,*} \lambda^*_v \cF \simeq 
        \rR f_{v,*} \lambda^*_v \rR g_{\et,*} \underline{\ZZ/m} \arrow[d,"{\rR f_{v,*}c_{g,\underline{\ZZ/m}}}"] \\
        \lambda^*_v \rR (f \circ g)_{\et,*}\underline{\ZZ/m}  \arrow[r,"{c_{f\circ g, \underline{\ZZ/m}}}"']& \rR (f \circ g)_{v,*} \lambda^*_v \underline{\ZZ/m}
    \end{tikzcd} \]
    commutes.
    Therefore, it suffices to prove that the right vertical and bottom horizontal maps are isomorphisms, reducing us to the case $\cF = \underline{\ZZ/m}$.
    Now the assertion follows from \cite[Cor.~16.7]{diamonds} and \cref{qproet-v-pushforward-etale}.
\end{proof}
\begin{definition}
    Let $f \colon X \to Y$ be a proper morphism of rigid-analytic spaces over $\Spa(K,\cO_K)$.
    The \emph{mod-$n$ trace} of $f$ is the morphism in $\cD^+_\et(Y;\ZZ/n) \subset \cD^+(Y_v;\ZZ/n)$
    \[ \Tr^{\ZZ/n}_f \colon \rR f_{v,*}\lambda^*_v\omega^{\ZZ/n}_X \xleftarrow[\sim]{c_{f,\omega^{\ZZ/n}_X}} \lambda^*_v\rR f_{\et,*}\omega^{\ZZ/n}_X \xlongrightarrow{\lambda^*_v\overline{\Tr}^{\ZZ/n}_f} \lambda^*_v\omega^{\ZZ/n}_Y \]
    induced by the trace $\overline{\Tr}^{\ZZ/n}_f$ from \cite[Th.~7.4.1, Th.~7.3.14]{LRZ24} under the identification from \cref{zc-cohomology-et-v}.
\end{definition}
Now we get to the construction of trace morphisms for $\ZZ_p$-dualizing complexes.
\begin{lemma}\label{proper-trace}
    Let $f \colon X \to Y$ be a proper morphism of rigid-analytic spaces over $\Spa(K,\cO_K)$.
    Then there is an essentially unique map $\Tr^{\ZZ_p}_f \colon \rR f_{v,*} \omega^{\ZZ_p}_X \to \omega^{\ZZ_p}_Y$ such that the diagram
    \[ \begin{tikzcd}
      \rR f_{v,*}\omega^{\ZZ_p}_X \otimes^\rL_{\ud{\ZZ}_p} \underline{\ZZ/p^r} \arrow[r,phantom,"\simeq"] \arrow[d,"\Tr^{\ZZ_p}_f \otimes^\rL \underline{\ZZ/p^r}"'] &[-2em] \rR f_{v,*}\lambda^*_v\omega^{\ZZ/p^r}_X \arrow[d,"\Tr^{\ZZ/p^r}_f"] \\
      \omega^{\ZZ_p}_Y \otimes^\rL_{\ud{\ZZ}_p} \underline{\ZZ/p^r} \arrow[r,phantom,"\simeq"] & \lambda^*_v \omega^{\ZZ/p^r}_Y 
    \end{tikzcd} \]
    in $\cD^+(Y_v;\ZZ/p^r)$ commutes for all $r$.
\end{lemma}
Note that since $\rR\Hom^{<0}(\rR f_{v,*}\omega^{\ZZ_p}_X \otimes^\rL_{\ud{\ZZ}_p} \underline{\ZZ/p^r},\lambda^*_v \omega^{\ZZ/p^r}_Y) \simeq \rR\Hom^{<0}(\lambda^*_v\rR f_{\et,*} \omega^{\ZZ/p^r}_X,\lambda^*_v \omega^{\ZZ/p^r}_Y) \simeq 0$ (see \cite[Th.~7.3.14]{LRZ24}, \cite[Prop.~14.10]{diamonds}, \cref{dualizing-construction}, and \cref{zc-cohomology-et-v}), the commutativity in the derived $\infty$-category $\cD^+(Y_v;\ZZ/p^r)$ of the diagram in \cref{proper-trace} is a property, not extra data. 
\begin{proof}
    For any $m \mid n$, the traces $\Tr^{\ZZ/n}_f$ and $\Tr^{\ZZ/m}_f$ are compatible via the commutative diagram of \cite[Th.~7.4.7]{LRZ24}.
    Applying $\lambda^*_v$, we see that the diagram
    \[ \begin{tikzcd}[column sep=7em,row sep=tiny]
        \rR f_{v,*}\lambda^*_v\omega^{\ZZ/n}_X \otimes^\rL_{\ud{\ZZ/n}} \underline{\ZZ/m} \arrow[d,phantom,sloped,"\simeq"] & \lambda^*_v\rR f_{\et,*}\omega^{\ZZ/n}_X \otimes^\rL_{\ud{\ZZ/n}} \underline{\ZZ/m} \arrow[l,"\sim","c_{f,\omega^{\ZZ/n}_X} \otimes_{\ud{\ZZ/n}}^\rL \ZZ/m"'] \arrow[r,"\lambda^*_v\overline{\Tr}^{\ZZ/n}_f \otimes^\rL \ZZ/m"] \arrow[d,phantom,sloped,"\simeq"] & \lambda^*_v\omega^{\ZZ/n}_Y \otimes^\rL \underline{\ZZ/m} \arrow[d,phantom,sloped,"\simeq"] \\
       \rR f_{v,*}\lambda^*_v\omega^{\ZZ/m}_X & \arrow[l,"\sim","c_{f,\omega^{\ZZ/m}_X}"'] \lambda^*_v\rR f_{\et,*}\omega^{\ZZ/m}_X \arrow[r,"\lambda^*_v\overline{\Tr}^{\ZZ/m}_f"] & \lambda^*_v\omega^{\ZZ/m}_Y
    \end{tikzcd} \]
    commutes, where the leftward horizontal arrows come from \cref{zc-cohomology-et-v} and the vertical isomorphisms are the ones from \cref{dualizing-construction}.
    Moreover, the commutativity is witnessed in an essentially unique way because $\rR\cHom^{<0}\bigl(\rR f_{v,*}\lambda^*_v\omega^{\ZZ/n}_X,\lambda^*_v\omega^{\ZZ/n}_Y\bigr) \simeq 0$ for all $n > 0$ by \cite[Th.~7.3.14]{LRZ24} (and again \cref{zc-cohomology-et-v} and the full faithfulness of $\lambda^*_v$ \cite[Prop.~14.10]{diamonds}).
    Hence, the collection of $\Tr^{\ZZ/p^r}_f$ determines an essentially unique map
    \[ \Tr^{\ZZ_p}_f \colon \rR f_{v,*}\omega^{\ZZ_p}_X \simeq \rR\lim_r\rR f_{v,*}\lambda^*_v\omega^{\ZZ/p^r}_X \xlongrightarrow{\rR\lim \Tr^{\ZZ/p^r}_f} \rR\lim_r \lambda^*_v \omega^{\ZZ/p^r}_Y = \omega^{\ZZ_p}_Y \]
    in $\cD^+(Y_v;\ZZ_p)$ which fits into the desired diagram.
\end{proof}
\begin{lemma}\label{Trace-properties-Zp}
    The trace morphisms $\Tr^{\ZZ_p}_{\blank}$ defined in \cref{proper-trace} have the following properties:
    \begin{enumerate}[leftmargin=*,label={\upshape{(\roman*)}}]
        \item\label{Trace-properties-Zp-compatibility} (Compatibility with smooth proper traces) For any smooth proper morphism $f \colon X \to Y$ of equidimension $d$ between rigid-analytic spaces over $K$, the $p$-adic trace maps $\tr^{\ZZ_p}_f$ from \cref{trace-maps}\cref{trace-maps-Zp}, the isomorphism $\alpha^{\ZZ_p}_f$ from \cref{dualizing-smooth-Zp}, and the completed projection formula map fit into a commutative diagram
        \[ \begin{tikzcd}[column sep=huge]
            \rR f_{v,*}\bigl(f^*_v\omega^{\ZZ_p}_Y(d)[2d]\bigr) \arrow[r,"\rR f_{v,*}\alpha^{\ZZ_p}_f","\sim"'] & \rR f_{v,*}\omega^{\ZZ_p}_X \arrow[d,"\Tr^{\ZZ_p}_f"] \\
            \omega^{\ZZ_p}_Y \widehat{\otimes}^\rL_{\underline{\ZZ}_p} \rR f_{v, *}\underline{\ZZ}_p(d)[2d] \arrow[u,"\PF_f"] \arrow[r,"\id \otimes^\rL \tr^{\ZZ_p}_f"] & \omega^{\ZZ_p}_Y.
        \end{tikzcd} \]
        \item\label{Trace-properties-Zp-composition} (Compatibility with compositions) For any two proper morphisms $f \colon X \to Y$ and $g \colon Y \to Z$ of rigid-analytic spaces over $K$, we have $\Tr^{\ZZ_p}_{g \circ f} \simeq \Tr^{\ZZ_p}_g \circ \rR g_{v,*}\Tr^{\ZZ_p}_f$.
    \end{enumerate}
\end{lemma} 
\begin{proof} The idea of the proof is to reduce to the case of $\ZZ/p^r$-sheaves, where these claims follow from \cite[Th.~7.4.1 (2), Prop.~7.3.2 (i)]{LRZ24} and \cite[Th.~7.4.1 (3)]{LRZ24}, respectively.  We spell out the details for \cref{Trace-properties-Zp-compatibility}, with the proof of \cref{Trace-properties-Zp-composition} being analogous. 

\begin{enumerate}[wide,label={\textit{Step~\arabic*}.},ref={Step~\arabic*}]
\item\label{Trace-properties-Zp-prelims} \textit{The natural map $\Hom\bigl(\omega^{\ZZ_p}_Y \widehat{\otimes}^\rL_{\underline{\ZZ}_p} \rR f_{v, *}\underline{\ZZ}_p(d)[2d], \omega^{\ZZ_p}_Y\bigr) \to \lim_r \Hom\bigl(\omega^{\ZZ/p^r}_Y \otimes^\rL_{\underline{\ZZ/p^r}} \rR f_{\et, *}\underline{\ZZ/p^r}(d)[2d], \omega^{\ZZ/p^r}_Y\bigr)$ is an isomorphism.}
Since $\omega^{\ZZ_p}_Y$ is derived $p$-complete, \cref{zc-cohomology-et-v} and \cite[Prop.~14.10]{diamonds} imply that 
\begin{multline*}
\rR \Hom_{\uZ_p}\bigl(\omega^{\ZZ_p}_Y \widehat{\otimes}^\rL_{\underline{\ZZ}_p} \rR f_{v, *}\underline{\ZZ}_p(d)[2d], \omega^{\ZZ_p}_Y\bigr) \simeq \rR\lim_r \rR\Hom_{\uZ_p}\bigl(\omega^{\ZZ_p}_Y \widehat{\otimes}^\rL_{\underline{\ZZ}_p} \rR f_{v, *}\underline{\ZZ}_p(d)[2d], \omega^{\ZZ/p^r}_Y\bigr) \\ \simeq \rR\lim_r \rR\Hom_{\ud{\Z/p^r}}\bigl(\omega^{\ZZ/p^r}_Y \otimes^\rL_{\underline{\ZZ/p^r}} \rR f_{\et, *}\underline{\ZZ/p^r}(d)[2d], \omega^{\ZZ/p^r}_Y\bigr).
\end{multline*}
Now note that \cite[Cor.~7.1.2]{LRZ24} and the projection formula (see \cite[Th.~5.5.9 (ii)]{Huber-etale}) ensure that $\omega^{\ZZ/p^r}_Y \otimes^\rL_{\underline{\ZZ/p^r}} \rR f_{\et, *}\underline{\ZZ/p^r}(d)[2d] \simeq \rR f_{\et, *} \omega^{\ZZ/p^r}_X$. Thus, we get the following short exact sequence:
\[
0 \to \lim_r {}^{1} \Ext^{-1}\bigl(\rR f_{\et, *} \omega^{\ZZ/p^r}_X, \omega^{\ZZ/p^r}_Y\bigr) \to \Hom\bigl(\omega^{\ZZ_p}_Y \widehat{\otimes}^\rL_{\underline{\ZZ}_p} \rR f_{v, *}\underline{\ZZ}_p(d)[2d], \omega^{\ZZ_p}_Y\bigr) \to \lim_r \Hom\bigl(\rR f_{\et, *} \omega^{\ZZ/p^r}_X, \omega^{\ZZ/p^r}_Y\bigr) \to 0
\]
In particular, it suffices to show that $\rR\Hom\bigl(\rR f_{\et, *} \omega^{\ZZ/p^r}_X, \omega^{\ZZ/p^r}_Y\bigr)$ 
lies in $\cD^{\geq 0}(\ZZ/p^r)$ for any integer $r\geq 1$. This follows directly from \cite[Th.~7.3.14]{LRZ24}.

\item \textit{End of the proof.}
\Cref{Trace-properties-Zp-prelims} implies that we can check that the diagram commutes after applying $\blank \otimes^\rL \underline{\ZZ/p^r}$ for all $r \in \ZZ_{>0}$. In this case, the claim follows from \cite[Th.~7.4.1 (2), Prop.~7.3.2 (i)]{LRZ24}. \qedhere
\end{enumerate}
\end{proof}
The next lemma will allow us later to reduce certain statements for $\omega^{\ZZ_p}$ to smooth proper rigid-analytic spaces.
Recall that a morphism $\pi \colon X' \to X$ of rigid-analytic spaces over $\Spa(K,\cO_K)$ is called a regular admissible modification if $\pi$ is proper, $X'$ is smooth over $\Spa(K,\cO_K)$, and there exists a Zariski-dense open subspace $U \subset X$ such that $\pi^{-1}(U) \subset X'$ is dense and the restriction $\restr{\pi}{\pi^{-1}(U)} \colon \pi^{-1}(U) \to U$ is an isomorphism.
\begin{lemma}
\label{dualizing-cohomology}
    Let $X$ be a quasi-compact reduced rigid-analytic space over $\Spa(K,\cO_K)$.
    \begin{enumerate}[leftmargin=*,label=\upshape{(\roman*)}]
        \item\label{dualizing-cohomology-modification-square} There exists a cartesian diagram of rigid-analytic spaces over $\Spa(K,\cO_K)$
        \begin{equation}\label{modification-square} \begin{tikzcd}
            Z' \arrow[r,"i'"] \arrow[d,"\pi'"] \arrow[rd,"h"] & X' \arrow[d,"\pi"] \\
            Z \arrow[r,"i"] & X
        \end{tikzcd} \end{equation}
        with the following properties:
        \begin{enumerate}[label={\upshape{(\alph*)}}]
            \item\label{dualizing-cohomology-modification} $\pi$ is a regular admissible modification; 
            \item\label{dualizing-cohomology-closed} $\dim Z,\dim Z' < \dim X$; and
            \item\label{dualizing-cohomology-pushout} applying the diamondification functor $(\blank)^\lozenge$ to \cref{modification-square} yields a pushout diagram in $v$-sheaves on $\Perf_{/\Spd(K,\cO_K)}$.
        \end{enumerate}
        \item\label{dualizing-cohomology-fiber-sequence} The following square in $\cD^+(X_v,\ud{\ZZ}_p)$ is bicartesian:
        \[ \begin{tikzcd}[column sep=huge]
            \rR h_{v,*}\omega^{\ZZ_p}_{Z'} \arrow[r,"{\rR\pi_{v,*}\bigl(\Tr^{\ZZ_p}_{i'}\bigr)}"] \arrow[d,"{\rR i_{v,*}\bigl(\Tr^{\ZZ_p}_{\pi'}\bigr)}"'] & \rR \pi_{v,*}\omega^{\ZZ_p}_{X'} \arrow[d,"\Tr^{\ZZ_p}_\pi"]  \\
            \rR i_{v,*} \omega^{\ZZ_p}_Z \arrow[r,"{\Tr^{\ZZ_p}_i}"] & \omega^{\ZZ_p}_X
        \end{tikzcd} \]
    \end{enumerate}
\end{lemma}
\begin{proof}
    \cref{dualizing-cohomology-modification-square}.
    Since $X$ is reduced, the smooth locus $U \subseteq X$ is Zariski-dense and open, and by Temkin's resolution of singularities, we can find a regular admissible modification $\pi \colon X' \to X$ which is an isomorphism over $U$;
    see e.g.\ \cite[Prop.~7.3.9]{LRZ24}.
    Let $i \colon Z \colonequals X \smallsetminus U \hookrightarrow X$ be the complement, endowed with the canonical reduced adic subspace structure.
    Set $Z' \colonequals Z \times_X X'$.
    This defines a diagram of the form \cref{modification-square} satisfying properties \cref{dualizing-cohomology-modification} and \cref{dualizing-cohomology-closed}.
    Under the diamondification functor $(\blank)^\lozenge$,
    \cref{modification-square} remains a fiber square of small $v$-sheaves (see e.g.\ \cite[Prop.~6.1.6]{zav-almost}), $\abs{\pi^\lozenge}$ is an isomorphism outside $\abs{U^\lozenge}$ (\cite[Lem.~15.6]{diamonds}), $i^\lozenge$ is a monomorphism (cf.\ \cite[\href{https://stacks.math.columbia.edu/tag/08LR}{Tag~08LR}]{stacks-project}), and $\pi^\lozenge$ and $i^\lozenge$ are proper due to \cref{lemma:different-proper}.
    Therefore, \cref{cartesian-cocartesian} yields property \cref{dualizing-cohomology-pushout}.

    \cref{dualizing-cohomology-fiber-sequence}.
    Choose a modification square \cref{modification-square} as in the proof of \cref{dualizing-cohomology-modification-square}.
    By the compatibility of traces with compositions (\cref{Trace-properties-Zp}\cref{Trace-properties-Zp-composition}), we have
    \[ \Tr^{\ZZ_p}_\pi \circ \rR\pi_{v,*}\bigl(\Tr^{\ZZ_p}_{i'}\bigr) \simeq \Tr^{\ZZ_p}_h \simeq \Tr^{\ZZ_p}_i \circ \rR i_{v,*}\bigl(\Tr^{\ZZ_p}_{\pi'}\bigr). \]
    This supplies the homotopy witnessing the commutativity of the displayed square, which is again necessarily unique because $\rR\Hom^{<0}\bigl(\rR h_{v,*}\omega^{\ZZ_p}_{Z'},\omega^{\ZZ_p}_X\bigr) \simeq 0$.

    In what follows, we will use repeatedly that a square in a stable $\infty$-category is bicartesian if and only if it is cartesian if and only if it is cocartesian \cite[Prop.~1.1.3.4]{HA}.
    Since all terms in the square are derived $p$-complete and the fully faithful inclusion $\cD^+_{p-\text{comp}}(X_v; \underline{\ZZ}_p)\subset \cD^+(X_v;\underline{\ZZ}_p)$ is a right adjoint (its left adjoint is given by derived $p$-completion) and hence preserves limits, it therefore suffices to prove that the square is cartesian in $\cD^+_{p-\text{comp}}(X_v; \underline{\ZZ}_p)$.
    On the other hand, $\cD^+_{p-\text{comp}}(X_v; \underline{\ZZ}_p) \simeq \lim_r \cD^+(X_v;\ZZ/p^r)$ is again stable \cite[Th.~1.1.4.4]{HA};
    as small colimits in a limit of $\infty$-categories can be detected pointwise \cite[\href{https://kerodon.net/tag/06B4}{Tag~06B4}]{kerodon}, it then suffices to show that the squares
    \[ \begin{tikzcd}[column sep=huge]
        \rR h_{v,*}\omega^{\ZZ/n}_{Z'} \arrow[r,"{\rR\pi_{v,*}\bigl(\Tr^{\ZZ/n}_{i'}\bigr)}"] \arrow[d,"{\rR i_{v,*}\bigl(\Tr^{\ZZ/n}_{\pi'}\bigr)}"'] & \rR \pi_{v,*}\omega^{\ZZ/n}_{X'} \arrow[d,"\Tr^{\ZZ/n}_\pi"]  \\
        \rR i_{v,*} \omega^{\ZZ/n}_Z \arrow[r,"{\Tr^{\ZZ/n}_i}"] & \omega^{\ZZ/n}_X
    \end{tikzcd} \]
    are cocartesian in $\cD^+(X_v;\ZZ/n)$ for all $n \in \ZZ$.
    
    Recall from \cite{LRZ24} that in the (triangulated) homotopy category of $\cD^+(X_\et,\ZZ/n)$, \cite[Lem.~7.3.12]{LRZ24} and the identity $\rR i^! \overline{\Tr}^{\ZZ/n}_\pi = \overline{\Tr}^{\ZZ/n}_{\pi'}$ (this corresponds to the identification 
    $\overline{\Tr}^{\ZZ/n}_\pi \circ \rR\pi_{\et,*}\bigl(\overline{\Tr}^{\ZZ/n}_{i'}\bigr) = \overline{\Tr}^{\ZZ/n}_i \circ \rR i_{\et,*}\bigl(\overline{\Tr}^{\ZZ/n}_{\pi'}\bigr)$
    under the adjunction between $i_*$ and $\rR i^!$)
    gives a distinguished triangle 
    \[ \rR h_{\et,*}\omega^{\ZZ/n}_{Z'} \xlongrightarrow{\bigl(-\rR\pi_{\et,*}\bigl(\overline{\Tr}^{\ZZ/n}_{i'}\bigr),\rR i_{\et,*}\bigl(\overline{\Tr}^{\ZZ/n}_{\pi'}\bigr)\bigr)} \rR \pi_{\et,*}\omega^{\ZZ/n}_{X'} \oplus \rR i_{\et,*} \omega^{\ZZ/n}_Z \xlongrightarrow{\overline{\Tr}^{\ZZ/n}_\pi \oplus \overline{\Tr}^{\ZZ/n}_i} \omega^{\ZZ/n}_X \longrightarrow \rR h_{\et,*}\omega^{\ZZ/n}_{Z'}[1]. \]
    The pullback under $\lambda^*_v$ combined with \cref{zc-cohomology-et-v} applied to the proper morphisms $h$, $\pi$ and $i$ shows that the induced sequence of maps
    \[ \rR h_{v,*}\omega^{\ZZ/n}_{Z'} \xlongrightarrow{\bigl(\rR\pi_{v,*}\bigl(\Tr^{\ZZ/n}_{i'}\bigr),-\rR i_{v,*}\bigl(\Tr^{\ZZ/n}_{\pi'}\bigr)\bigr)} \rR \pi_{v,*}\omega^{\ZZ/n}_{X'} \oplus \rR i_{v,*} \omega^{\ZZ/n}_Z \xlongrightarrow{\Tr^{\ZZ/n}_\pi \oplus \Tr^{\ZZ/n}_i} \omega^{\ZZ/n}_X \longrightarrow \rR h_{v,*}\omega^{\ZZ/n}_{Z'}[1] \]
    again forms a distinguished triangle in the homotopy category of $\cD^+(X_v;\ZZ/n)$.
    But now \cite[Th.~1.1.2.14]{HA} and Verdier's axiom (TR3) yield the desired isomorphism
    \[ \rR \pi_{v,*}\omega^{\ZZ/n}_{X'} \bigsqcup_{\rR h_{v,*}\omega^{\ZZ/n}_{Z'}} \rR i_{v,*} \omega^{\ZZ/n}_Z \simeq \bigl(\rR \pi_{v,*}\omega^{\ZZ/n}_{X'} \oplus \rR i_{v,*}\omega^{\ZZ/n}_Z\bigr) \bigsqcup_{\rR h_{v,*}\omega^{\ZZ/n}_{Z'}} 0 \xlongrightarrow{\sim} \omega^{\ZZ/n}_X \]
    in $\cD^+(X_v;\ZZ/n)$, finishing the proof.
\end{proof}

\subsection{Dualizing complexes: \texorpdfstring{$\cB$}{B}-coefficients}
Next, we give $\cO$- and $\BB_I$-versions of the dualizing complex and the trace morphism.
To do so, we need a version of the relative primitive comparison theorem that applies to dualizing complexes.
\Cref{Relative primitive comparison-2} proves a comparison for all perfect $\ZZ_p$-complexes.
Similarly to what we mentioned above for $\omega^{\ZZ/n}_X$, the dualizing complexes $\omega^{\ZZ_p}_X$ are in general not perfect, but are still contained in the full subcategory $\cD^{(b)}_\zc(X;\ZZ_p) \subset \cD^+(X_v;\ZZ_p)$ of locally bounded complexes with Zariski-constructible cohomology sheaves, this time in the sense of \cite[Def.~3.32]{BH}.
A derived $p$-complete object $\cF \in \cD^+(X_v;\ZZ_p)$ lies in $\cD^{(b)}_\zc(X;\ZZ_p)$ if its reduction $\cF \otimes^\rL_{\underline{\ZZ}_p} \underline{\FF}_p$ lies in $\cD^{(b)}_\zc(X_\et;\FF_p) \subset \cD^+(X_v;\FF_p)$ (in the previous sense).
We refer to \cite{BH} for details.
\begin{lemma}\label{dualizing-comparison-BI}
    Let $f \colon X \to Y$ be a proper morphism of rigid-analytic spaces over $\Spa(K,\cO_K)$ and $\cF \in \cD^{(b)}_\zc(X;\ZZ_p)$.
    Let $\cA \in \{ \cO^+,\AA_I \}$ and $\cB \colonequals \cA\bigl[\tfrac{1}{p}\bigr]$.
    Then the natural morphism in $\cD^+(Y_v;\cA)$
    \[ \bigl(\rR f_{v,*} \cF \bigr) \widehat{\otimes}^\rL_{\underline{\ZZ}_p} \cA \longrightarrow \rR f_{v,*}\bigl(\cF \widehat{\otimes}^\rL_{\underline{\ZZ}_p} \cA\bigr) \]
    is an almost isomorphism with respect to $\AA^{\circ\circ}_I$ if $\cA = \AA_I$ and with respect to $\cO^{\circ\circ}_Y$ if $\cA = \cO^+$.
    In particular,
    \[ \bigl((\rR f_{v,*} \cF) \widehat{\otimes}^\rL_{\underline{\ZZ}_p} \cA \bigr)\bigl[\tfrac{1}{p}\bigr] \longrightarrow \rR f_{v,*}\bigl((\cF \widehat{\otimes}^\rL_{\underline{\ZZ}_p} \cA)\bigl[\tfrac{1}{p}\bigr]\bigr) \]
    is an isomorphism in $\cD^+(Y_v;\cB)$.
\end{lemma}
\begin{proof}
    Since the statement is local on $Y$, we may assume that $Y$ (and hence $X$) are quasi-compact.
    In this case, the complex $\cF \in \cD^b_\zc(X;\ud{\ZZ}_p)$ is globally bounded. Similarly, the complex $\cF \widehat{\otimes}^\rL_{\underline{\ZZ}_p} \cA$ is globally bounded below; since it is $p$-complete, it suffices to check boundedness after reducing modulo $p$, where it becomes isomorphic to $\F/^\rL p \otimes^\rL_{\ud{\FF}_p}\cA/p\in \cD^b(X_v; \FF_p)$. 
    Since $f$ is qcqs and $\cF \widehat{\otimes}^\rL_{\underline{\ZZ}_p} \cA$ is bounded below, we have $\rR f_{v,*}\bigl((\cF \widehat{\otimes}^\rL_{\underline{\ZZ}_p} \cA)\bigl[\tfrac{1}{p}\bigr]\bigr) \simeq \rR f_{v,*}\bigl((\cF \widehat{\otimes}^\rL_{\underline{\ZZ}_p} \cA)\bigr)\bigl[\tfrac{1}{p}\bigr]$, so the second assertion follows immediately from the first one by inverting $p$.
    Moreover, it suffices to prove the first assertion after derived reduction mod $p$ because both source and target are derived $p$-complete.
    In other words, we need to show that the natural morphism
    \begin{equation}\label{Fp-zc-primitive-comparison}
        \beta_f(\cG) \colon \bigl(\rR f_{v,*} \cG\bigr) \widehat{\otimes}^\rL_{\underline{\ZZ}_p} \cA \simeq \bigl(\rR f_{v,*} \cG\bigr) \otimes^\rL_{\underline{\FF}_p} \cA/p \longrightarrow \rR f_{v,*}\bigl(\cG \otimes^\rL_{\underline{\FF}_p} \cA/p\bigr) \simeq \rR f_{v,*}\bigl(\cG \widehat{\otimes}^\rL_{\underline{\ZZ}_p} \cA\bigr)
    \end{equation}
    is an almost isomorphism for any $\cG \in \cD^b_\zc(X;\FF_p)$.

    Using that the full subcategory of $\cG \in \cD^b_\zc(X;\FF_p)$ for which $\beta_f(\cG)$ is an almost isomorphism is a thick subcategory, we can now make the same reduction as in the proof of \cref{zc-cohomology-et-v} to assume that $\cG \in \cD^b_\zc(X;\FF_p)$ is of the form $\cG = \lambda^*_v \rR g_{\et,*} \underline{\FF}_p$ for a finite morphism $g \colon W \to X$.
    The perfectness of $\underline{\FF}_p \in \cD^+(W_v;\ZZ_p)$ combined again with \cite[Cor.~16.7]{diamonds} and \cref{qproet-v-pushforward-etale} guarantees that $\cG = \lambda^*_v \rR g_{\et,*} \underline{\FF}_p \simeq \rR g_{v,*} \underline{\FF}_p.$
    
    On the other hand, the comparison morphisms from \cref{Fp-zc-primitive-comparison} are compatible with compositions in $f$, i.e.,
    \[ \bigl(\rR (f \circ g)_{v,*} \underline{\FF}_p\bigr) \widehat{\otimes}^\rL_{\underline{\ZZ}_p} \cA \xrightarrow{\beta_f(\rR g_{v,*} \underline{\FF}_p)} \rR f_{v,*} \bigl( (\rR g_{v,*} \underline{\FF}_p) \widehat{\otimes}^\rL_{\underline{\ZZ}_p} \cA \bigr) \xrightarrow{\rR f_{v,*}(\beta_g(\underline{\FF}_p))} \rR (f \circ g)_{v,*}\bigl(\underline{\FF}_p \widehat{\otimes}^\rL_{\underline{\ZZ}_p} \cA\bigr) \]
    is given by $\beta_{f \circ g}(\underline{\FF}_p)$.
    However, \cref{Relative primitive comparison-2} shows that both $\beta_g(\underline{\FF}_p)$ and $\beta_{f \circ g}(\underline{\FF}_p)$ are almost isomorphisms.
    As a consequence, $\beta_f(\rR g_{v,*} \underline{\FF}_p)$ must be an almost isomorphism as well, finishing the proof.
\end{proof}

\begin{definition}\label{needed1}
    Let $X$ be a rigid-analytic space over $\Spa(K,\cO_K)$.
    Let $\cA \in \{ \underline{\ZZ}_p,\cO^+,\AA_I \}$ and $\cB \colonequals \cA\bigl[\tfrac{1}{p}\bigr]$.
    \begin{enumerate}[leftmargin=*,label={\upshape{(\arabic*)}}]
        \item The ($\cB$-)\emph{dualizing complex} of $X$ is $\omega^{\cB}_X \colonequals \bigl(\omega^{\ZZ_p}_X \widehat{\otimes}^\rL_{\underline{\ZZ}_p} \cA \bigr)\bigl[\tfrac{1}{p}\bigr] \in \cD(X_v;\cB)$.
        \item The ($\cB$-)\emph{Verdier duality} functor for $X$ is $\DD_X(\blank) \colonequals \rR\cHom\bigl(\blank,\omega^{\cB}_X\bigr) \colon \cD(X_v;\cB)^\op \to \cD(X_v;\cB)$.
    \end{enumerate}
\end{definition}
\begin{example}\label{dualizing-smooth-B}
    Assume that $X$ is a smooth rigid-analytic space over $\Spa(K,\cO_K)$ of equidimension $d$.
    Then the canonical isomorphism $\alpha^{\ZZ_p}_X$ from \cref{dualizing-smooth-Zp} induces a canonical isomorphism
    \[ \alpha^\cB_X \colon \cB(d)[2d] \simeq \bigl(\underline{\ZZ}_p(d)[2d] \widehat{\otimes}^\rL_{\underline{\ZZ}_p} \cA\bigr)\bigl[\tfrac{1}{p}\bigr] \xlongrightarrow[\sim]{\alpha^{\ZZ_p}_X \otimes^\rL \id} \bigl(\omega^{\ZZ_p}_X \widehat{\otimes}^\rL_{\underline{\ZZ}_p} \cA\bigr)\bigl[\tfrac{1}{p}\bigr] = \omega^\cB_X. \]
    More generally, for any smooth morphism $f \colon X \to Y$ of equidimension $d$ between rigid-analytic spaces over $\Spa(K,\cO_K)$, the isomorphism $\alpha^{\ZZ_p}_f$ from \cref{dualizing-smooth-Zp} induces a canonical isomorphism in $\cD(X_v;\cB)$
    \[ \alpha^\cB_f \colon \bigl(f^*_v\omega^{\ZZ_p}_Y \widehat{\otimes}^\rL_{\underline{\ZZ}_p} \cA\bigr)\bigl[\tfrac{1}{p}\bigr](d)[2d] \simeq \bigl(f^*_v\omega^{\ZZ_p}_Y(d)[2d] \widehat{\otimes}^\rL_{\underline{\ZZ}_p} \cA\bigr)\bigl[\tfrac{1}{p}\bigr] \xlongrightarrow[\sim]{\alpha^{\ZZ_p}_f \otimes^\rL \id} \bigl(\omega^{\ZZ_p}_X \widehat{\otimes}^\rL_{\underline{\ZZ}_p} \cA\bigr)\bigl[\tfrac{1}{p}\bigr] = \omega^\cB_X. \] 
\end{example}
In order to give the $\cO$- and $\cB$-versions of the trace map, we need to recall the following notation.
\begin{notation}[{cf.\ \cite[Not.~7.5.1]{LRZ24}}]\label{Ev-Adj}
    Let $f \colon X \to Y$ be a morphism of rigid-analytic spaces over $K$ and let $\cB \in \{ \underline{\QQ}_p,\cO,\BB_I \}$.
    We denote by
    \[ \Ev_f \colon \rR f_{v,*}\rR\cHom(\blank,\blank) \longrightarrow \rR\cHom\bigl(\rR f_{v,*}(\blank),\rR f_{v,*}(\blank)\bigr) \]
    the natural \emph{evaluation transformation} of functors of $\infty$-categories $\cD(X_v;\cB)^\op \times \cD(X_v;\cB) \to \cD(Y_v;\cB)$ which corresponds under the tensor-hom adjunction to the functor of $\infty$-categories given on objects $\cE,\cE' \in \cD(X_v;\cB)$ by the composition
    \[ \rR f_{v,*}\rR\cHom(\cE,\cE') \otimes^L \rR f_{v,*}(\cE) \xlongrightarrow{\cup} \rR f_{v,*}(\rR\cHom(\cE,\cE') \otimes^L \cE) \longrightarrow \rR f_{v,*}(\cE') \]
    of the cup product map (coming from the lax symmetric monoidal structure of $\rR f_*$) and the evaluation map for internal homs;
    cf.\ e.g.\ \cite[\href{https://stacks.math.columbia.edu/tag/0B6D}{Tag~0B6D}]{stacks-project}.

    Further, we denote by
    \[ \Adj_f \colon \rR f_{v,*}\rR\cHom\bigl(f^*_v(\blank),\blank\bigr) \xlongrightarrow{\Ev_f} \rR\cHom\bigl(\rR f_{v,*}f^*_v(\blank),\rR f_{v,*}(\blank)\bigr) \xlongrightarrow{\blank \circ \eta_f} \rR\cHom\bigl(\blank,\rR f_{v,*}(\blank)\bigr) \]
    the natural \emph{adjunction transformation} of functors of $\infty$-categories $\cD(Y_v;\cB)^\op \times \cD(X_v;\cB) \to \cD(Y_v;\cB)$ given by $\Ev_f$ and composition with the unit of the adjunction between $f^*_v$ and $\rR f_{v,*}$.
\end{notation}
\begin{remark}\label{Adj-isomorphism}
    The natural adjunction transformation $\Adj_f$ from \cref{Ev-Adj} is in fact a natural isomorphism.
    Using the full faithfulness of the Yoneda embedding (\cite[\href{https://kerodon.net/tag/03NJ}{Tag~03NJ}]{kerodon}) and the fact that it can be tested on objects whether a natural transformation is a natural isomorphism (\cite[\href{https://kerodon.net/tag/01DK}{Tag~01DK}]{kerodon}), this reduces to the calculation that for any $\cE,\cF \in \cD(Y_v;\cB)$ and $\cE' \in \cD(X_v;\cB)$
    \begin{align*}
        \Map\bigl(\cF,\rR f_{v,*}\rR\cHom\bigl(f^*_v(\cE),\cE'\bigr)\bigr) &\simeq \Map\bigl(f^*_v(\cF),\rR\cHom\bigl(f^*_v(\cE),\cE'\bigr)\bigr) \simeq \Map\bigl(f^*_v(\cF) \otimes^L f^*_v(\cE),\cE'\bigr) \\
        &\simeq \Map\bigl(f^*_v(\cF \otimes^L \cE),\cE'\bigr) \simeq \Map\bigl(\cF \otimes^L \cE,\rR f_{v,*}(\cE')\bigr) \\
        &\simeq \Map\bigl(\cF,\rR\cHom\bigl(\cE,\rR f_{v,*}(\cE')\bigr)
    \end{align*}
    thanks to tensor-hom adjunction and the strict symmetric monoidality of $f^*_v$.
\end{remark}
\begin{remark}\label{Ev-Adj-composition}
    One can see exactly as in \cite[Rmk.~7.5.2]{LRZ24} that given any two morphisms $f \colon X \to Y$ and $g \colon Y \to Z$ of rigid-analytic spaces over $K$, we have $\Ev_{g \circ f} = \Ev_g \circ \rR g_{v,*}(\Ev_f)$ and $\Adj_{g \circ f} = \Adj_g \circ \rR g_{v,*}(\Adj_f)$.
\end{remark}
\begin{definition}
\label{BI-proper-trace}
    Let $f \colon X \to Y$ be a proper morphism of rigid-analytic spaces over $\Spa(K,\cO_K)$, let $\cA \in \{ \underline{\ZZ}_p,\cO^+,\AA_I \}$ and $\cB \colonequals \cA\bigl[\tfrac{1}{p}\bigr]$, and let $\cE \in \cD(X_v;\cB)$.
    \begin{enumerate}[leftmargin=*,label=\upshape{(\arabic*)}]
        \item The \emph{$\cB$-trace} of $f$ is the morphism in $\cD^+(Y_v;\cB)$
        \[ \Tr^{\cB}_f \colon \rR f_{v,*}\omega^{\cB}_X \xleftarrow{\sim} \bigl((\rR f_{v,*}\omega^{\ZZ_p}_X) \widehat{\otimes}^\rL_{\underline{\ZZ}_p} \cA\bigr)\bigl[\tfrac{1}{p}\bigr] \xlongrightarrow{\bigl(\Tr^{\ZZ_p}_f \otimes^\rL \id\bigr)\bigl[\tfrac{1}{p}\bigr]} \bigl(\omega^{\ZZ_p}_Y \widehat{\otimes}^\rL_{\underline{\ZZ}_p} \cA\bigr)\bigl[\tfrac{1}{p}\bigr] = \omega^{\cB}_Y. \]
        Here, we invert the first arrow, which is an equivalence due to \cref{dualizing-comparison-BI},
        and use the $p$-adic trace map from \cref{proper-trace} for the second arrow.
        \item\label{BI-proper-PD-morphism} The \emph{duality morphism} for $\cE$ is the composition in $\cD(Y_v;\cB)$
        \begin{multline*}
            \PD_f(\cE) \colon \rR f_{v,*}\bigl(\DD_X(\cE)\bigr) = \rR f_{v,*}\rR\cHom\bigl(\cE,\omega^\cB_X\bigr) \xlongrightarrow{\Ev_f} \rR\cHom\bigl(\rR f_{v,*}\cE,\rR f_{v,*}\omega^\cB_X\bigr) \\ \xlongrightarrow{\Tr^{\cB}_f \circ \blank} \rR\cHom\bigl(\rR f_{v,*}\cE,\omega^{\cB}_Y\bigr) = \DD_Y(\rR f_{v,*}\cE)
        \end{multline*}
        of the evaluation map from \cref{Ev-Adj} and composition with the $\cB$-trace of $f$.
    \end{enumerate}
\end{definition}
Next, we explain the analog of \cref{Trace-properties-Zp} for $\Tr^\cB_f$.
\begin{lemma}\label{Trace-properties}
    The trace morphisms $\Tr^\cB_{\blank}$ from \cref{BI-proper-trace} have the following properties:
    \begin{enumerate}[leftmargin=*,label={\upshape{(\roman*)}}]
        \item\label{Trace-properties-compatibility} (Compatibility with smooth proper traces) For any smooth proper morphism $f \colon X \to Y$ of equidimension $d$ between rigid-analytic spaces over $K$, the $\cA$-trace maps $\tr^\cA_f$ from \cref{trace-maps}\cref{trace-maps-AI}, the isomorphism $\alpha^\cB_f$ from \cref{dualizing-smooth-B}, and the completed projection formula map fit into a commutative diagram
        \[ \begin{tikzcd}[column sep=huge]
            \rR f_{v,*}\bigl((f^*_v\omega^{\ZZ_p}_Y \widehat{\otimes}^\rL_{\underline{\ZZ}_p} \cA)\bigl[\tfrac{1}{p}\bigr](d)[2d]\bigr) \arrow[r,"\rR f_{v,*}\alpha^\cB_f","\sim"'] & \rR f_{v,*}\omega^\cB_X \arrow[d,"\Tr^\cB_f"] \\
            \bigl(\omega^{\ZZ_p}_Y \widehat{\otimes}^\rL_{\underline{\ZZ}_p} \rR f_{v, *}\cA\bigr)\bigl[\tfrac{1}{p}\bigr](d)[2d] \arrow[u,"\PF_f"] \arrow[r,"\id \otimes^\rL \tr^\cA_f"] & \omega^\cB_Y.
        \end{tikzcd} \]
        \item\label{Trace-properties-composition} (Compatibility with compositions) For any two proper morphisms $f \colon X \to Y$ and $g \colon Y \to Z$ of rigid-analytic spaces over $K$, we have $\Tr^\cB_{g \circ f} \simeq \Tr^\cB_g \circ \rR g_{v,*}\Tr^\cB_f$.
    \end{enumerate}
\end{lemma}
\begin{proof}
    \cref{Trace-properties-compatibility}.
    The primitive comparison maps give rise to the following diagram:
    \[ \begin{tikzcd}[row sep=tiny,column sep=tiny]
            & \rR f_{v,*}\bigl((f^*_v\omega^{\ZZ_p}_Y \widehat{\otimes}^\rL_{\underline{\ZZ}_p} \cA)\bigl[\tfrac{1}{p}\bigr](d)[2d]\bigr) \arrow[rr,"\rR f_{v,*}\alpha^\cB_f","\sim"'] \arrow[from=dd,near end,"\PF_f"] && \rR f_{v,*}\omega^\cB_X \arrow[dd,"\Tr^\cB_f"] \\
            \bigl(\rR f_{v,*}\bigl(f^*_v\omega^{\ZZ_p}_Y(d)[2d]\bigr) \widehat{\otimes}^\rL_{\underline{\ZZ}_p} \cA \bigr)\bigl[\tfrac{1}{p}\bigr] \arrow[rr,crossing over,near end,"(\rR f_{v,*}\alpha^{\ZZ_p}_f) \otimes \id","\sim"'] \arrow[ru] && \bigl((\rR f_{v,*}\omega^{\ZZ_p}_X) \widehat{\otimes}^\rL_{\underline{\ZZ}_p} \cA \bigr)\bigl[\tfrac{1}{p}\bigr] \arrow[ru] & \\
            & \bigl(\omega^{\ZZ_p}_Y \widehat{\otimes}^\rL_{\underline{\ZZ}_p} \rR f_{v, *}\cA\bigr)\bigl[\tfrac{1}{p}\bigr](d)[2d] \arrow[rr,near start,"\id \otimes \tr^\cA_f"] && \omega^\cB_Y \\
            \bigl(\omega^{\ZZ_p}_Y \widehat{\otimes}^\rL_{\underline{\ZZ}_p} (\rR f_{v, *}\underline{\ZZ}_p(d)[2d]) \widehat{\otimes}^\rL_{\underline{\ZZ}_p} \cA \bigr)\bigl[\tfrac{1}{p}\bigr] \arrow[uu,"\PF_f"] \arrow[rr,"\id \otimes \tr^{\ZZ_p}_f \otimes \id"] \arrow[ru,red,sloped,"\sim"] && \bigl(\omega^{\ZZ_p}_Y \widehat{\otimes}^\rL_{\underline{\ZZ}_p} \cA \bigr)\bigl[\tfrac{1}{p}\bigr] \arrow[from=uu,crossing over,near start,"\Tr^{\ZZ_p}_f \otimes \id"] \arrow[ru,equals] & 
    \end{tikzcd} \]
    We want to show that the back side of the diagram commutes.
    The front side commutes by \cref{Trace-properties-Zp}\cref{Trace-properties-Zp-compatibility}.
    The right side commutes by the definition of $\Tr^\cB_f$ in \cref{BI-proper-trace} and the bottom side by the definition of $\tr^\cA_f$ in \cref{trace-maps}.
    The left side and top side commute by the naturality of the primitive comparison map.
    Lastly, the red map is an isomorphism by the relative primitive comparison theorem \cref{Relative primitive comparison-2}.
    Taken together, this shows that the back side of the diagram must commute as well.

    \cref{Trace-properties-composition}.
    The natural primitive comparison maps also give rise to the following commutative diagram:
    \[ \begin{tikzcd}[column sep=7em]
        \bigl((\rR(g \circ f)_{v,*}\omega^{\ZZ_p}_X) \widehat{\otimes}^\rL_{\underline{\ZZ}_p} \cA\bigr)\bigl[\tfrac{1}{p}\bigr] \arrow[r,"(\rR g_{v,*}\Tr^{\ZZ_p}_f) \otimes \id"] \arrow[d,sloped,"\sim"] & \bigl((\rR g_{v,*}\omega^{\ZZ_p}_Y) \widehat{\otimes}^\rL_{\underline{\ZZ}_p} \cA\bigr)\bigl[\tfrac{1}{p}\bigr] \arrow[r,"\Tr^{\ZZ_p}_g \otimes \id"] \arrow[d,sloped,"\sim"] & (\omega^{\ZZ_p}_Z \widehat{\otimes}^\rL_{\underline{\ZZ}_p} \cA)\bigl[\tfrac{1}{p}\bigr] \arrow[dd,equals] \\
        \rR g_{v,*}\bigl((\rR f_{v,*}\omega^{\ZZ_p}_X \widehat{\otimes}^\rL_{\underline{\ZZ}_p} \cA)\bigl[\tfrac{1}{p}\bigr]\bigr) \arrow[r,"\rR g_{v,*}(\Tr^{\ZZ_p}_f \otimes \id)"] \arrow[d,sloped,"\sim"] & \rR g_{v,*}\bigl((\omega^{\ZZ_p}_Y \widehat{\otimes}^\rL_{\underline{\ZZ}_p} \cA)\bigl[\tfrac{1}{p}\bigr]\bigr) \arrow[d,equals] & \\
        \rR(g \circ f)_{v,*}\omega^\cB_X \arrow[r,"\rR g_{v,*}\Tr^{\cB}_f"] & \rR g_{v,*}\omega^\cB_Y \arrow[r,"\Tr^\cB_g"] & \omega^\cB_Z
    \end{tikzcd}\]
    Since $\rR f_{v, *}\omega_X^{\ZZ_p}$ is Zariski-constructible \cite[Th.~3.36 (3)]{BH}, all vertical arrows are isomorphisms by \cref{dualizing-comparison-BI}.
    The composition of top horizontal arrows is given by $\Tr^{\ZZ_p}_{g \circ f} \otimes \id$ thanks to \cref{Trace-properties-Zp}\cref{Trace-properties-Zp-composition}.
    Keeping in mind \cref{BI-proper-trace}, the composition of bottom horizontal arrows must therefore be $\Tr^\cB_{g \circ f}$.
\end{proof}

\subsection{Geometric duality for proper rigid-analytic spaces}\label{subsection:singular-spaces}
With all definitions and compatibilities in place, we get to the first main result of this section: 
finiteness and duality statements for (possibly singular) proper rigid-analytic spaces over $K$,
at least in the absolute case.
Our approach is based on the strategy in \cite[\S~7]{LRZ24}.
\begin{theorem}
\label{singular-PD}
    Let $f \colon X \to \Spa(K,\cO_K)$ be a proper rigid-analytic space, let $\cA \in \{ \underline{\ZZ}_p,\cO^+,\AA_I \}$ and $\cB \colonequals \cA\bigl[\tfrac{1}{p}\bigr]$, and let $\cE \in \Perf_v(X;\cB)$.
    Then $\rR f_{v,*}\cE$ lies in $\cD^{(b)}_{\BBC}(\Spd(K,\cO_K)_v;\underline{\QQ}_p) \subset \cD^+(\Spd(K,\cO_K)_v;\underline{\QQ}_p)$ if $\cA = \underline{\ZZ}_p$ and in $\Perf_v(\Spd(K,\cO_K);\cB) \subset \cD^+(\Spd(K,\cO_K)_v;\cB)$ if $\cA \in \{ \cO^+,\AA_I \}$.
    Moreover, the duality morphism
    \[ \PD_f(\cE) \colon \rR f_{v,*}\bigl(\DD_X(\cE)\bigr) \longrightarrow \DD_{\Spa(K,\cO_K)}(\rR f_{v,*}\cE) \]
    is an isomorphism in $\cD^+(\Spd(K,\cO_K)_v;\cB)$.
\end{theorem}
\begin{proof}
    We proceed by induction on $d \colonequals \dim X \in \ZZ_{\ge 0} \cup \{-\infty\}$.
    If $d = -\infty$, we have $X = \varnothing$, so there is nothing to show.
    For the induction step, we assume that $d \ge 0$ and that the statement has been proven in dimensions $< d$.
    Let $X_\red \subseteq X$ be the maximal reduced closed subspace.
    Under the isomorphism $X^\lozenge \simeq X^\lozenge_\red$, the dualizing complexes $\omega^{\cB}_X$ and $\omega^{\cB}_{X_\red}$ are identified because the induced equivalence of \'etale sites identifies $\omega^{\ZZ/n}_X$ and $\omega^{\ZZ/n}_{X_\red}$ (cf.\ \cite[Lem.~7.2.4]{LRZ24}).
    Therefore, we may assume that $X$ is reduced.
    
    \Cref{dualizing-cohomology}\cref{dualizing-cohomology-modification-square} now provides the pushout square of diamonds
    \[ \begin{tikzcd}
       (Z')^\lozenge \arrow[r,"i'"] \arrow[d,"\pi'"'] \arrow[rd,"h"] & (X')^\lozenge \arrow[d,"\pi"] \\
       Z^\lozenge \arrow[r,"i"] & X^\lozenge
    \end{tikzcd} \]
    in which $X'$ is smooth proper and $Z, Z'$ are proper with $\dim Z, \dim Z'<d$. Applying the exact functor $\rR\cHom(\ZZ[\blank],\cB)$ to this diagram of $v$-sheaves on $\Perf_{/X}$, we obtain the following pullback square in $\cD^+(X_v;\cB)$:
    \[ \begin{tikzcd}
       \cB \arrow[r] \arrow[d] & \rR \pi_{v,*}\cB \arrow[d] \\
       \rR i_{v,*}\cB \arrow[r] & \rR h_{v,*}\cB
    \end{tikzcd} \]   

    Tensoring with $\cE$, using the projection formula \cref{PF} and applying $\rR f_{v,*}$, this yields the following pullback squares in $\cD^+(\Spd(K,\cO_K)_v;\cB)$, in which the $\eta$-maps are given by the units of the adjunctions between $\pi^*_v$ and $\rR \pi_{v,*}$, $i^*_v$ and $\rR i_{v,*}$, and $h^*_v$ and $\rR h_{v,*}$:
    \begin{equation}\label{v-MV} \begin{tikzcd}[column sep=huge]
       \cE \arrow[r,"\eta_\pi"] \arrow[d,"\eta_i"'] & \rR\pi_{v,*}\pi^*_v\cE \arrow[d,"\rR \pi_{v,*}\eta_{i'}"] \\
       \rR i_{v,*}i^*_v\cE \arrow[r,"\rR i_{v,*}\eta_{\pi'}"] & \rR h_{v,*}h^*_v\cE
       \end{tikzcd} \qquad \text{and} \qquad \begin{tikzcd}[column sep=huge]
       \rR f_{v,*}\cE \arrow[r,"\rR f_{v,*}\eta_\pi"] \arrow[d,"\rR f_{v,*} \eta_i"'] & \rR (f \circ \pi)_{v,*}\pi^*_v\cE \arrow[d,"\rR (f \circ \pi)_{v,*}\eta_{i'}"] \\
       \rR (f \circ i)_{v,*}i^*_v\cE \arrow[r,"\rR (f \circ i)_{v,*}\eta_{\pi'}"] & \rR (f \circ h)_{v,*}h^*_v\cE
    \end{tikzcd} \end{equation}
    By the induction hypothesis, both $\rR(f\circ i)_{v,*}i_v^*\cE$ and $\rR(f\circ h)_{v,*}h_v^*\cE$ lie in $\cD^{(b)}_{\BBC}(\Spd(K,\cO_K)_v;\underline{\QQ}_p) \subset \cD^+(\Spd(K,\cO_K)_v;\underline{\QQ}_p)$ (if $\cA = \underline{\ZZ}_p$) resp.\ in $\Perf_v(\Spd(K,\cO_K);\cB) \subset \cD^+(\Spd(K,\cO_K)_v;\cB)$ (if $\cA \in \{ \cO^+,\AA_I \}$).
    Moreover, $\rR(f\circ\pi)_{v,*}\pi_v^*\cE$ lies in $\cD^{(b)}_{\BBC}(\Spd(K,\cO_K)_v;\underline{\QQ}_p)$ resp.\ $\Perf_v(\Spd(K,\cO_K);\cB)$ by \cref{thm:cohomology-banach-colmez-spaces} resp.\ \cref{cor:proper-smooth-pushforward-is-perfect}.
    Since both categories are stable and hence closed under finite limits, the first assertion follows.

    For the proof that $\PD_f(\cE)$ is an isomorphism, consider the following diagram in $\cD^+(\Spd(K,\cO_K)_v;\cB)$, in which we have set squares of different depths in different colors to improve readability:
    \begin{equation}\label{singular-PD-modification} \begin{tikzcd}[baseline=(current  bounding  box.north),row sep=small,column sep=-7em,scale cd=.85,center picture]
        &[2em]&& \rR\cHom\bigl(\rR (f \circ h)_{v,*}h^*_v\cE,\cB \bigr) \arrow[rrrr,orange,"\blank \circ \rR (f \circ \pi)_{v,*}\eta_{i'}"] \arrow[dddd,orange,pos=.35,"\blank \circ \rR (f \circ i)_{v,*}\eta_{\pi'}"] &[3em]&[2em]&& \rR\cHom\bigl(\rR (f \circ \pi)_{v,*}\pi^*_v\cE,\cB \bigr) \arrow[dddd,orange,"\blank \circ \rR f_{v,*}\eta_\pi"] \\
        && \rR\cHom\bigl(\rR (f \circ h)_{v,*}h^*_v\cE,\rR (f \circ h)_{v,*}\omega^\cB_{Z'}\bigr) \arrow[ru,purple,pos=.1,"\Tr^\cB_{f \circ h} \circ \blank"] \arrow[dddd,olive,"\rR (f \circ i)_{v,*}\Tr^\cB_{\pi'} \circ \blank \circ \rR (f \circ i)_{v,*}\eta_{\pi'}"'{description,xshift=-1em,yshift=-1.5em}] &&&& \rR\cHom\bigl(\rR (f \circ \pi)_{v,*}\pi^*_v\cE,\rR (f \circ \pi)_{v,*}\omega^\cB_{X'}\bigr) \arrow[ru,purple,pos=.1,"\Tr^\cB_{f \circ \pi} \circ \blank"] \arrow[from=llll,crossing over,olive,"\rR (f \circ \pi)_{v,*}\Tr^\cB_{i'} \circ \blank \circ \rR (f \circ \pi)_{v,*}\eta_{i'}"] & \\
        & \rR (f \circ h)_{v,*}\rR\cHom\bigl(h^*_v\cE,\omega^\cB_{Z'}\bigr) \arrow[ru,purple,pos=.2,"\Ev_{f \circ h}"] \arrow[ld,"\rR f_{v,*}\Adj_h"'{pos=.8},"\sim"'{sloped}] \arrow[dddd,blue,"(\Tr^\cB_{\pi'} \circ \blank) \circ \rR (f \circ i)_{v,*}\Adj_{\pi'}"'{description,xshift=-1.5em,yshift=-1em}] &&&& \rR (f \circ \pi)_{v,*}\rR\cHom\bigl(\pi^*_v\cE,\omega^\cB_{X'}\bigr) \arrow[ru,purple,pos=.2,"\Ev_{f \circ \pi}"] \arrow[ld,"\rR f_{v,*}\Adj_\pi"'{pos=.8},"\sim"'{sloped}] \arrow[from=llll,crossing over,blue,"(\Tr^\cB_{i'} \circ \blank) \circ \rR (f \circ \pi)_{v,*}\Adj_{i'}"'] && \\
        \rR f_{v,*}\rR\cHom\bigl(\cE,\rR h_{v,*}\omega^\cB_{Z'}\bigr) \arrow[dddd,red,"\rR i_{v,*}\Tr^\cB_{\pi'} \circ \blank"'] &&&& \rR f_{v,*}\rR\cHom\bigl(\cE,\rR \pi_{v,*}\omega^\cB_{X'}\bigr) \arrow[from=llll,crossing over,red,"\rR \pi_{v,*}\Tr^\cB_{i'} \circ \blank",pos=.35] &&& \\
        &&& \rR\cHom\bigl(\rR (f \circ i)_{v,*}i^*_v\cE,\cB \bigr) \arrow[rrrr,orange,pos=.65,"\blank \circ \rR f_{v,*}\eta_i"] &&&& \rR\cHom\bigl(\rR f_{v,*}\cE,\cB \bigr) \\
        && \rR\cHom\bigl(\rR (f \circ i)_{v,*}i^*_v\cE,\rR (f \circ i)_{v,*}\omega^\cB_Z \bigr) \arrow[rrrr,olive,"\rR f_{v,*}\Tr^\cB_i \circ \blank \circ \rR f_{v,*}\eta_i"{xshift=1.2em}] \arrow[ru,purple,pos=.8,"\Tr^\cB_{f \circ i} \circ \blank"'] &&&& \rR\cHom\bigl(\rR f_{v,*}\cE,\rR f_{v,*}\omega^\cB_X \bigr) \arrow[ru,purple,pos=.8,"\Tr^\cB_f \circ \blank"'] \arrow[from=uuuu,crossing over,olive,"\rR f_{v,*}\Tr^\cB_\pi \circ \blank \circ \rR f_{v,*}\eta_\pi"{description,xshift=2em,yshift=1em}] & \\
        & \rR (f \circ i)_{v,*}\rR\cHom\bigl(i^*_v\cE,\omega^\cB_Z \bigr) \arrow[rrrr,blue,"(\Tr^\cB_i \circ \blank) \circ \rR f_{v,*}\Adj_i"'] \arrow[ru,purple,pos=.8,"\Ev_{f \circ i}"'] \arrow[ld,"\rR f_{v,*}\Adj_i"{pos=.2},"\sim"{sloped}] &&&& \rR f_{v,*}\rR\cHom\bigl(\cE,\omega^\cB_X \bigr) \arrow[ru,purple,pos=.8,"\Ev_f"'] \arrow[ld,equals] \arrow[from=uuuu,crossing over,blue,near start,"(\Tr^\cB_\pi \circ \blank) \circ \rR f_{v,*}\Adj_\pi"] && \\
        \rR f_{v,*}\rR\cHom\bigl(\cE,\rR i_{v,*}\omega^\cB_Z \bigr) \arrow[rrrr,red,"\Tr^\cB_i \circ \blank"'] &&&& \rR f_{v,*}\rR\cHom\bigl(\cE,\omega^\cB_X \bigr) \arrow[from=uuuu,crossing over,red,"\Tr^\cB_\pi \circ \blank"',pos=.65] &&&
    \end{tikzcd} \end{equation}
    We now explain how to construct homotopies witnessing the commutativity of \cref{singular-PD-modification}.
    
    \textit{Blue and red square.}
    In order to construct the cube connecting the blue and red square, we first use the compatibility of $\Adj$ and $\Tr^\cB$ under compositions (\cref{Ev-Adj-composition} and \cref{Trace-properties}\cref{Trace-properties-composition}) and the naturality of $\Adj$ to find $2$-simplices witnessing the commutativity of the following diagram:
    \begin{equation}\label{singular-PD-modification-1} \hspace{-2.5em}\begin{tikzcd}[column sep=large]
        \rR (f \circ h)_{v,*} \rR\cHom\bigl(h^*_v\cE,\omega^\cB_{Z'}\bigr) \arrow[r,"\rR (f \circ \pi)_{v,*}\Adj_{i'}"] \arrow[d,"\rR (f \circ i)_{v,*}\Adj_{\pi'}"] &[1em] \rR (f \circ \pi)_{v,*}\rR\cHom\bigl(\pi^*_v\cE,\rR i'_{v,*}\omega^\cB_{Z'}\bigr) \arrow[r,"\Tr^\cB_{i'} \circ \blank"] \arrow[d,"\rR f_{v,*}\Adj_\pi"] & \rR (f \circ \pi)_{v,*}\rR\cHom\bigl(\pi^*_v\cE,\omega^\cB_{X'}\bigr) \arrow[d,"\rR f_{v,*}\Adj_\pi"] \\
        \rR (f \circ i)_{v,*} \rR\cHom\bigl(i^*_v\cE,\rR \pi'_{v,*}\omega^\cB_{Z'}\bigr) \arrow[r,"\rR f_{v,*}\Adj_i"] \arrow[d,"\Tr^\cB_{\pi'} \circ \blank"] & \rR f_{v,*}\rR\cHom\bigl(\cE,\rR h_{v,*}\omega^\cB_{Z'}\bigr) \arrow[r,"\rR \pi_{v,*}\Tr^\cB_{i'} \circ \blank"] \arrow[d,"\rR i_{v,*}\Tr^\cB_{\pi'} \circ \blank"] & \rR f_{v,*}\rR\cHom\bigl(\cE,\rR \pi_{v,*}\omega^\cB_{X'}\bigr) \arrow[d,"\Tr^\cB_\pi \circ \blank"] \\
        \rR (f \circ i)_{v,*} \rR\cHom\bigl(i^*_v\cE,\omega^\cB_Z\bigr) \arrow[r,"\rR f_{v,*}\Adj_i"] & \rR f_{v,*}\rR\cHom\bigl(\cE,\rR i_{v,*}\omega^\cB_Z\bigr) \arrow[r,"\Tr^\cB_i \circ \blank"] & \rR f_{v,*}\rR\cHom\bigl(\cE,\omega^\cB_X\bigr)
    \end{tikzcd} \end{equation}
    The outer rectangle of \cref{singular-PD-modification-1} is the blue square. 
    Now precompose the functors $\Delta^1 \times \Delta^1 \to \cD^+\bigl(\Spa(K,\cO_K)_v;\cB\bigr)$ corresponding to
    \begin{itemize}
        \item the upper left square of \cref{singular-PD-modification-1} with the map of simplicial sets $\Delta^1 \times \Delta^1 \times \Delta^1 \to \Delta^1 \times \Delta^1$ given on vertices by $\{0,1\}^3 \ni (i,j,k) \mapsto \bigl(\max(i,k),\max(j,k)\bigr)$;
        \item the upper right square of \cref{singular-PD-modification-1} with the map of simplicial sets $\Delta^1 \times \Delta^1 \times \Delta^1 \to \Delta^1 \times \Delta^1$ given on vertices by $\{0,1\}^3 \ni (i,j,k) \mapsto \bigl(i,\max(j,k)\bigr)$;
        \item the lower left square of \cref{singular-PD-modification-1} with the map of simplicial sets $\Delta^1 \times \Delta^1 \times \Delta^1 \to \Delta^1 \times \Delta^1$ given on vertices by $\{0,1\}^3 \ni (i,j,k) \mapsto \bigl(\max(i,k),j\bigr)$;
        \item the lower right square of \cref{singular-PD-modification-1} with the map of simplicial sets $\Delta^1 \times \Delta^1 \times \Delta^1 \to \Delta^1 \times \Delta^1$ given on vertices by $\{0,1\}^3 \ni (i,j,k) \mapsto (i,j)$.
    \end{itemize}
    The composition of these cubes yields the desired cube connecting the blue and the red square.

    \textit{Blue and green square.}
    In order to construct the cube connecting the blue and green square, note that by \cref{Ev-Adj}, we have $\Adj = \Ev \circ (\blank \circ \eta)$.
    Therefore, we can also expand the blue square as follows:
     \begin{equation}\label{singular-PD-modification-2} \begin{tikzcd}[baseline=(current  bounding  box.north),column sep=large,scale cd=.85,center picture]
        \rR (f \circ h)_{v,*} \rR\cHom\bigl(h^*_v\cE,\omega^\cB_{Z'}\bigr) \arrow[r,"\rR (f \circ \pi)_{v,*}\Ev_{i'}"] \arrow[d,"\rR (f \circ i)_{v,*}\Ev_{\pi'}"] & \rR (f \circ \pi)_{v,*}\rR\cHom\bigl(\rR i'_{v,*} h^*_v\cE,\rR i'_{v,*}\omega^\cB_{Z'}\bigr) \arrow[r,"\Tr^\cB_{i'} \circ \blank \circ \eta_{i'}"] \arrow[d,"\rR f_{v,*}\Ev_\pi"] &[3em] \rR (f \circ \pi)_{v,*}\rR\cHom\bigl(\pi^*_v\cE,\omega^\cB_{X'}\bigr) \arrow[d,"\rR f_{v,*}\Ev_\pi"] \\
        \rR (f \circ i)_{v,*} \rR\cHom\bigl(\rR \pi'_{v,*}h^*_v\cE,\rR \pi'_{v,*}\omega^\cB_{Z'}\bigr) \arrow[r,"\rR f_{v,*}\Ev_i"] \arrow[d,"\Tr^\cB_{\pi'} \circ \blank \circ \eta_{\pi'}"] & \rR f_{v,*}\rR\cHom\bigl(\rR h_{v,*}h^*_v\cE,\rR h_{v,*}\omega^\cB_{Z'}\bigr) \arrow[r,"\rR \pi_{v,*}\Tr^\cB_{i'} \circ \blank \circ \rR \pi_{v,*}\eta_{i'}"] \arrow[d,"\rR i_{v,*}\Tr^\cB_{\pi'} \circ \blank \circ \rR i_{v,*}\eta_{\pi'}"] & \rR f_{v,*}\rR\cHom\bigl(\rR \pi_{v,*}\pi^*_v\cE,\rR \pi_{v,*}\omega^\cB_{X'}\bigr) \arrow[d,"\Tr^\cB_\pi \circ \blank \circ \eta_\pi"] \\
        \rR (f \circ i)_{v,*} \rR\cHom\bigl(i^*_v\cE,\omega^\cB_Z\bigr) \arrow[r,"\rR f_{v,*}\Ev_i"] & \rR f_{v,*}\rR\cHom\bigl(\rR i_{v,*}i^*_v\cE,\rR i_{v,*}\omega^\cB_Z\bigr) \arrow[r,"\Tr^\cB_i \circ \blank \circ \eta_i"] & \rR f_{v,*}\rR\cHom\bigl(\cE,\omega^\cB_X\bigr)
    \end{tikzcd} \end{equation}   
    Thanks to \cref{Ev-Adj-composition} and \cref{v-MV}, we can find homotopies which, together with the already chosen homotopy from \cref{Trace-properties}\cref{Trace-properties-composition}, witness the commutativity of \cref{singular-PD-modification-2}.
    In fact, by virtue of the identity $\Adj = \Ev \circ (\blank \circ \eta)$, we can construct the homotopy for the top left square in \cref{singular-PD-modification-1} from the ones chosen for the top left and bottom right square in \cref{singular-PD-modification-2}, making the two diagrams compatible.
    Exactly the same maps of simplicial sets $\Delta^1 \times \Delta^1 \times \Delta^1 \to \Delta^1 \times \Delta^1$ as before produce a cube connecting the blue square with the commutative square
    \begin{equation}\label{singular-PD-square-in-between} \begin{tikzcd}[column sep=8em]
        \rR f_{v,*}\rR\cHom\bigl(\rR h_{v,*}h^*_v\cE,\rR h_{v,*}\omega^\cB_{Z'}\bigr) \arrow[r,"\rR \pi_{v,*}\Tr^\cB_{i'} \circ \blank \circ \rR \pi_{v,*}\eta_{i'}"] \arrow[d,"\rR i_{v,*}\Tr^\cB_{\pi'} \circ \blank \circ \rR i_{v,*}\eta_{\pi'}"] & \rR f_{v,*}\rR\cHom\bigl(\rR \pi_{v,*}\pi^*_v\cE,\rR \pi_{v,*}\omega^\cB_{X'}\bigr) \arrow[d,"\Tr^\cB_\pi \circ \blank \circ \eta_\pi"] \\
        \rR f_{v,*}\rR\cHom\bigl(\rR i_{v,*}i^*_v\cE,\rR i_{v,*}\omega^\cB_Z \bigr) \arrow[r,"\Tr^\cB_i \circ \blank \circ \eta_i"] & \rR f_{v,*}\rR\cHom\bigl(\cE,\omega^\cB_X \bigr),
    \end{tikzcd} \end{equation}
    whose $2$-simplices can be constructed from the same homotopies as before.

    Now note that the evaluation transformation $\Ev_f$ from \cref{Ev-Adj} amounts to a morphism in the $\infty$-category $\Fun\bigl(\cD(X_v;\cB)^\op \times \cD(X_v;\cB),\cD(Y_v;\cB)\bigr)$.
    Via the map of simplicial sets $\Delta^1 \times \Delta^1 \to \cD(X_v;\cB)^\op \times \cD(X_v;\cB)$ classifying the (opposite of the) left square in \cref{v-MV} and the square
    \[ \begin{tikzcd}[column sep=huge]
        \rR h_{v,*}\omega^\cB_{Z'} \arrow[r,"{\rR\pi_{v,*}\Tr^\cB_{i'}}"] \arrow[d,"{\rR i_{v,*}\Tr^\cB_{\pi'}}"'] & \rR \pi_{v,*}\omega^\cB_{X'} \arrow[d,"\Tr^\cB_\pi"]  \\
        \rR i_{v,*} \omega^\cB_Z \arrow[r,"{\Tr^\cB_i}"] & \omega^\cB_X
    \end{tikzcd} \]
    from before, we obtain a morphism in the $\infty$-category $\Fun\bigl(\Delta^1 \times \Delta^1,\cD(Y_v;\cB)\bigr)$.
    Unwinding definitions, this yields a cube $\Delta^1 \times \Delta^1 \times \Delta^1 \to \cD(Y_v;\cB)$ whose front side is the square \cref{singular-PD-square-in-between} and whose back side is the green square in \cref{singular-PD-modification};
    cf.\ \cite[\href{https://kerodon.net/tag/01DG}{Tag 01DG}, \href{https://kerodon.net/tag/01DH}{Tag 01DH}]{kerodon}.
    Composing this cube with the one from the previous paragraph, we obtain the desired cube connecting the blue and green square.

    \textit{Green and orange square.}
    Using the same homotopy from \cref{v-MV} as before to witness the commutativity of the orange square as well as the equivalences $\omega^{\cB}_{\Spa(K,\cO_K)} \simeq \cB$ from \cref{dualizing-smooth-B}, we can obtain the cube connecting the green and the orange square similarly to the previous two cases.
    We leave the details to the careful reader.

    \textit{Red and orange square.}
    Note that the maps $\Adj_{(\blank)}$ from the blue to the red square are equivalences by \cref{Adj-isomorphism}.
    Their inverses still fit into a commutative cube, with the necessary homotopies guaranteed by \cite[\href{https://kerodon.net/tag/01DK}{Tag 01DK}]{kerodon}.
    Composing these homotopies with the ones connecting the blue to the green and the green to the orange square produces homotopies witnessing the commutativity of the outer box.

    \textit{End of argument.}
    The red square (with the chosen homotopy coming from \cref{Trace-properties}\cref{Trace-properties-composition}) is a pushout square by \cref{dualizing-cohomology}\cref{dualizing-cohomology-fiber-sequence}, \cref{dualizing-comparison-BI}, and the exactness of the functor $\rR f_{v,*}\rR\cHom\bigl(\cE,\blank\widehat{\otimes}^\rL_{\underline{\ZZ}_p} \cA \bigr)\bigl[\tfrac{1}{p}\bigr]\bigr)$.
    Likewise, we saw above that the orange square (with the chosen homotopy coming from \cref{v-MV}) is a pushout square.
    By virtue of \cite[\href{https://kerodon.net/tag/02XP}{Tag 02XP}]{kerodon}, the map $\PD_f(\cE)$ from the bottom right front to the bottom right back corner is then a pushout of the other maps from the front to the back corners.
    In particular, in order to finish the proof, it suffices to see that those three maps are equivalences.
    
    The induction hypothesis guarantees that $\PD_{f \circ h}(h^*_v\cE) = (\Tr^\cB_{f \circ h} \circ \blank) \circ \Ev_{f \circ h}$ and $\PD_{f \circ i}(i^*_v\cE) = (\Tr^\cB_{f \circ i} \circ \blank) \circ \Ev_{f \circ i}$ are equivalences, and thus so are the maps from the top and bottom left front to the top and bottom left back corner.
    Lastly, the map from the top right front to the top right back corner is the composition of $\rR f_{v,*}\Adj^{-1}_\pi$ with the purple map in the diagram\footnote{In this diagram, we denote by $\ud{d}_{X'}\colon \abs{X'} \to \ZZ_{\geq 0}$ the locally constant function which sends each point $x\in X'$ to the dimension of the connected component of $X'$ containing $x$.}
    \[ \begin{tikzcd}[row sep=.05em,column sep=huge,scale cd=.9,center picture]
        \rR (f \circ \pi)_{v,*}\rR\cHom\bigl(\pi^*_v\cE,\cB(\ud{d}_{X'})[2\ud{d}_{X'}]\bigr) \arrow[r,"\Ev_{f \circ \pi}"] \arrow[dd,"\rR (f \circ \pi)_{v,*}(\alpha^\cB_{X'} \circ \blank)","\sim"'{sloped}] & \rR\cHom\bigl(\rR (f \circ \pi)_{v,*}\pi^*_v\cE,\rR (f \circ \pi)_{v,*}\cB(\ud{d}_{X'})[2 \ud{d}_{X'}]\bigr) \arrow[rd,"\tr^\cB_{f \circ \pi} \circ \blank"] \arrow[dd,"\rR (f \circ \pi)_{v,*}(\alpha^\cB_{X'}) \circ \blank","\sim"'{sloped}] & \\
        && \rR\cHom(\rR (f \circ \pi)_{v,*}\pi^*_v \cE, \cB \bigr) \\
        \rR (f \circ \pi)_{v,*}\rR\cHom\bigl(\pi^*_v\cE,\omega^\cB_{X'}\bigr) \arrow[r,purple,"\Ev_{f \circ \pi}"] & \rR\cHom\bigl(\rR (f \circ \pi)_{v,*}\pi^*_v\cE,\rR (f \circ \pi)_{v,*}\omega^\cB_{X'}\bigr), \arrow[ru,purple,"\Tr^{\cB}_{f \circ \pi} \circ \blank"'] & 
    \end{tikzcd}\]
    which commutes thanks to \cref{Trace-properties}\cref{Trace-properties-compatibility} and the isomorphism $\omega^{\ZZ_p}_{\Spa(K,\cO_K)} \simeq \underline{\ZZ}_p$.
    The construction of the evaluation transformation in \cite[\href{https://stacks.math.columbia.edu/tag/0B6D}{Tag~0B6D}]{stacks-project} via the tensor-hom adjunction shows that the composition $(\tr^\cB_{f \circ \pi} \circ \blank) \circ \Ev_{f \circ \pi}$ in the top is precisely the duality morphism from \cref{defn:Qp-duality-map}\cref{defn:Qp-duality-map-3} (if $\cA = \underline{\ZZ}_p$) and in \cref{construction:evaluation-coevaluation-poincare-duality-non-solid}\cref{defn:Qp-duality-map-3} (if $\cA \in \{ \cO^+,\AA_I \}$).
    Thus, is is an equivalence thanks to \cref{geometric-duality-Qp} resp.\ \cref{cor:PD-for-proper-smooth-pushforward}.
    This finishes the proof.
\end{proof}
By a similar argument, we get the following compatibility in the interval $I$:
\begin{theorem}
\label{singular pullback compatibility}
Let $X$ be a proper rigid-analytic space over $\Spa(K,\cO_K)$, with structure
morphism $f$.
Let $J \subset I \subset (0, \infty)$ be an inclusion of closed intervals
with rational endpoints and let $\cE \in \Perf_v(X;\BB_I)$.
Then the natural maps
\[
\rR f_{v, *}\cE \otimes^{\rL}_{\BB_I, \res}
\BB_J \rightarrow
\rR f_{v, *}(\cE \otimes^{\rL}_{\BB_I, \res} \BB_J)
\quad \text{and} \quad
\rR f_{v, *}\cE \otimes^{\rL}_{\BB_I, \varphi}
\BB_{I/p} \rightarrow
\rR f_{v, *}(\cE \otimes^{\rL}_{\BB_I, \varphi} \BB_{I/p})
\]
are isomorphisms.
\end{theorem}

\begin{proof}
We only give a proof for the first isomorphism, as the proof for the second isomorphism
follows from the same argument.
As in the proof of \cref{singular-PD}, we may assume that $X$ is reduced and, using induction, that the statement has been shown in dimensions $< \dim X$.
Following the notation there, we have a map from the red to the orange pushout square in the next diagram:
\[ \begin{tikzcd}[row sep=small,column sep=-3em]
& \rR f_{v, *}(\cE \otimes^{\rL}_{\BB_I, \res} \BB_J) \arrow[rr,orange] \arrow[dd,orange] && \rR (f \circ \pi)_{v,*}(\pi^*_v\cE \otimes^{\rL}_{\BB_I, \res} \BB_J) \arrow[dd,orange] \\
\rR f_{v, *}(\cE) \otimes^{\rL}_{\BB_I, \res} \BB_J \arrow[ru,purple] \arrow[rr,crossing over,red] \arrow[dd,red] && \rR (f \circ \pi)_{v,*}(\pi^*_v\cE) \otimes^{\rL}_{\BB_I, \res} \BB_J \arrow[ru,purple] & \\ 
& \rR (f \circ i)_{v,*}(i^*_v\cE \otimes^{\rL}_{\BB_I, \res} \BB_J) \arrow[rr,orange] && \rR (f \circ h)_{v,*}(h^*_v\cE \otimes^{\rL}_{\BB_I, \res} \BB_J) \\
\rR (f \circ i)_{v,*}(i^*_v\cE) \otimes^{\rL}_{\BB_I, \res} \BB_J \arrow[ru,purple] \arrow[rr,red] && \rR (f \circ h)_{v,*}(h^*_v\cE) \otimes^{\rL}_{\BB_I, \res} \BB_J \arrow[ru,purple] \arrow[from=uu,crossing over,red] &
\end{tikzcd} \]
Since $f \circ \pi$ is smooth proper, the purple map in the top right corner
is an isomorphism by \cref{cor:proper-smooth-pushforward-is-perfect}.
The purple maps in the two bottom corners are both isomorphisms by the induction hypothesis.
Therefore, the purple map in the top left corner must also be an isomorphism.
\end{proof}
Over a perfectoid field, \cref{singular-PD} and \cref{singular pullback compatibility} now admit the following absolute variants:
\begin{corollary}
\label{absolute-singular-PD}
    Let $K$ be a perfectoid field extension of $\QQ_p$ and $X$ be a proper rigid-analytic space over $\Spa(K,\cO_K)$.
    Let $\cA \in \{\cO^+,\AA_I\}$ and $\cB \colonequals \cA\bigl[\tfrac{1}{p}\bigr]$, and let $\cE \in \Perf_v(X;\cB)$.
    Set $B \colonequals \cB(K)$ and $B_I \colonequals \BB_I(K)$.
    \begin{enumerate}[leftmargin=*,label=\upshape{(\roman*)}]
        \item\label{absolute-singular-PD-dualizable} The condensed derived global sections object $\rR\underline{\Gamma}_v(X,\cE) \in \cD_\sol(\underline{B})$ is dualizable with dual $\rR\underline{\Gamma}_v\bigl(X,\DD_X(\cE)\bigr)$.
        In particular, $\rR\underline{\Gamma}_v(X,\cE) \simeq \underline{\rR\Gamma_v(X,\cE)}$ lies in $\Perf(B) \subset \cD_\sol(\underline{B})$ and the natural duality map
        \[ \rR\Gamma_v\bigl(X,\DD_X(\cE)\bigr) \xlongrightarrow{\sim} \rR\Hom_B\bigl(\rR\Gamma_v(X,\cE),B\bigr) \]
        in $\cD(B)$ is an isomorphism.
        \item\label{absolute-singular-PD-interval} Suppose that $\cA = \AA_I$ and $\cB = \BB_I$.
        Then for any inclusion $J \subset I \subset (0, \infty)$ of closed intervals with rational endpoints, the natural maps
        \[ \rR\Gamma_v(X,\cE) \otimes^\rL_{B_I, \res} B_J \longrightarrow \rR\Gamma_v(X,\cE \otimes^\rL_{\BB_I,\res} \BB_J) \quad \text{and} \quad \rR\Gamma_v(X,\cE) \otimes^\rL_{B_I,\varphi} B_{I/p} \longrightarrow \rR\Gamma_v(X,\cE \otimes^\rL_{\BB_I,\varphi} \BB_{I/p}) \]
        are isomorphisms.
    \end{enumerate}
\end{corollary}
\begin{proof}
    We begin with preliminary comments on the correspondence of various notions of perfect complexes for perfectoid fields:
    Consider the factorization
    \[ \rR\Gamma_v(\Spd(K,\cO_K),\blank) \colon \cD\bigl(\Spd(K,\cO_K)_v;\cB\bigr) \xlongrightarrow{\rR\underline{\Gamma}_v(\Spd(K,\cO_K),\blank)} \cD\bigl(*_\proet;\underline{B}\bigr) \xlongrightarrow{\rR\Gamma(*_\proet,\blank)} \cD(B) \]
    through the derived condensed global sections functor $\rR\underline{\Gamma}_v(\Spd(K,\cO_K),\blank)$, which was defined as the derived pushforward along the morphism of site $\Spd(K,\cO_K)_v \to *_\proet$ from \cref{condensed structure definition}. Since the full $\infty$-subcategory $\Perf\bigl(*_\proet;\underline{B}\bigr) \subset \cD\bigl(*_\proet;\underline{B}\bigr)$ consists of those objects which are pro\'etale locally isomorphic to a strictly perfect complex, every pro\'etale covering $S \to *$ splits, and $\rR\Gamma(*_\proet, \ud{B})=B$, the second functor is essentially surjective and fully faithful, hence induces an equivalence $\rR\Gamma(*_\proet,\blank) \colon \Perf\bigl(*_\proet;\underline{B}\bigr) \xlongrightarrow{\sim} \Perf(B)$ with inverse $\underline{(\blank)}$.
    On the other hand, since $K$ is perfectoid, \cref{cor:perfect-O-modules-on-affinoid-perfectoid} and \cref{cor:perfect-O-modules-on-affinoid-perfectoid-BBI} show that $\rR\Gamma_v(\Spd(K,\cO_K),\blank)$ induces an equivalence $\Perf\bigl(\Spd(K,\cO_K)_v;\cB\bigr) \xrightarrow{\sim} \Perf(B)$ with inverse $\widetilde{(\blank)}$.
    As a consequence, $\rR\underline{\Gamma}_v(\Spd(K,\cO_K),\blank)$ also induces an equivalence $\Perf\bigl(\Spd(K,\cO_K)_v;\cB\bigr) \xrightarrow{\sim} \Perf\bigl(*_\proet;\underline{B}\bigr)$.
    By \cite[Cor.~5.51.1]{Andreychev}, all three categories of perfect complexes further identify with the dualizable objects in $\cD_\sol(\underline{B})$.

    Now denote by $f \colon X \to \Spa(K,\cO_K)$ the structure map.
    We know from \cref{singular-PD} that $\rR f_{v,*}\cE$ lies in $\Perf_v(\Spd(K,\cO_K);\cB)$ and has dual $\rR f_{v,*}\bigl(\DD_X(\cE)\bigr)$.
    Applying $\rR\Gamma_v(\Spd(K,\cO_K),\blank)$, the preliminary comments from the previous paragraph (combined with the identification $\omega^{\BB_I}_{\Spa(K,\cO_K)} \simeq \BB_I$ from \cref{dualizing-smooth-B}) guarantee that $\rR\Gamma_v(X,\cE) = \rR\Gamma_v(\Spa(K,\cO_K),\rR f_{v,*}\cE)$ is still perfect with natural duality isomorphism
    \begin{multline*}
        \rR\Gamma_v\bigl(X,\DD_X(\cE)\bigr) \simeq \rR\Gamma_v\bigl(\Spa(K,\cO_K),\rR f_{v,*}\DD_X(\cE)\bigr) \xlongrightarrow[\sim]{\rR\Gamma_v(\Spa(K,\cO_K),\PD_f(\cE))} \\ \rR\Gamma_v\bigl(\Spa(K,\cO_K),\DD_{\Spa(K,\cO_K)}(\rR f_{v,*}\cE)\bigr) \simeq \rR\Hom_{\cD(\Spd(K,\cO_K)_v;\BB)}\bigl(\rR f_{v,*}\cE,\cB\bigr) \\ \simeq \rR\Hom_B\bigl(\rR\Gamma_v(X,\cE),\rR\Gamma_v(\Spa(K,\cO_K),\cB)\bigr) \simeq \rR\Hom_B\bigl(\rR\Gamma_v(X,\cE),B\bigr)
    \end{multline*}
    and that moreover $\rR\underline{\Gamma}_v(X,\cE) \simeq \underline{\rR\Gamma_v(X,\cE)}$ is discrete.
    This makes $\rR\underline{\Gamma}_v(X,\cE)$ into a dualizable object of $\cD_\sol(\underline{B})$ with dual $\rR\underline{\Gamma}_v\bigl(X,\DD_X(\cE)\bigr)$, yielding \cref{absolute-singular-PD-dualizable}.

    For assertion \cref{absolute-singular-PD-interval}, we show that the first natural map is an isomorphism, the proof for the second map being entirely analogous.
    The natural transformation $\widetilde{(\blank)} \circ \res \Rightarrow \res \circ \widetilde{(\blank)}$ 
    between the two compositions of functors in the diagram
    \[ \begin{tikzcd}[row sep=small]
        \Perf(B_I) \arrow[r,"\widetilde{(\blank)}"] \arrow[d,"\res"] & \Perf_v(\Spd(K,\cO_K);\BB_I) \arrow[d,"\res"] \\
        \Perf(B_J) \arrow[r,"\widetilde{(\blank)}"] & \Perf_v(\Spd(K,\cO_K);\BB_J)
    \end{tikzcd}\]
    is an isomorphism because $\Perf(B_I)$ is generated by $B_I$ and the transformation is an isomorphism on $B_I$.
    In particular, applying $\widetilde{(\blank)}$ to $\rR\Gamma_v(X,\cE) \otimes^\rL_{B_I, \res} B_J \to \rR\Gamma_v(X,\cE \otimes^\rL_{\BB_I,\res} \BB_J)$ yields the natural map $\rR f_{v, *}\cE \otimes^{\rL}_{\BB_I, \res} \BB_J \to \rR f_{v, *}(\cE \otimes^{\rL}_{\BB_I, \res} \BB_J)$, which is an isomorphism by \cref{singular pullback compatibility}.
    Since $\widetilde{(\blank)}$ is an equivalence by the first paragraph, this shows that the first natural map is an isomorphism as well, as desired.
\end{proof}

Similarly, we also obtain the compatibility with $\theta \colon \BB_{I} \to \cO$ as well as the K\"unneth formula.

\begin{theorem}
\label{thm:HT realization for proper pushforward} 
Let $f \colon X \to \Spa(K,\cO_K)$ be a proper rigid-analytic space, let $I\subset (0, \infty)$ be a closed interval with rational endpoints such that $1\in I$, and let $\cE \in \Perf_v(X;\BB_{I})$. Then the natural map
\[
\rR f_{v, *}\cE \otimes^{\rL}_{\BB_{I}, \theta_I}
\cO \rightarrow
\rR f_{v, *}(\cE \otimes^{\rL}_{\BB_{I}, \theta_I} \cO)
\]
is an isomorphism.

Assume furthermore that $K$ is a perfectoid field extension of $\QQ_p$, then the natural map
\[
\rR\Gamma_v(X,\cE) \otimes^{\rL}_{B_{I}, \theta_I} K \longrightarrow \rR\Gamma_v(X,\cE \otimes^{\rL}_{\BB_{I}, \theta_I} \cO)
\]
is an isomorphism.
\end{theorem}
We refer to \cref{construction:theta-sheaves} for the construction of the $\theta_I$-homomoprhism
\begin{proof}
For the first statement, we argue exactly as the proof of \Cref{singular pullback compatibility}:
We just need to replace \Cref{cor:proper-smooth-pushforward-is-perfect} with 
\Cref{cor:HT realization for proper smooth pushforward}.
The second statement follows from the first, by the exact same argument as in the second paragraph
of the proof of \Cref{absolute-singular-PD}.
\end{proof}

\begin{thm}
\label{thm:KF-for-proper-pushforward}
Let $f\colon X \to \Spa(K, \O_K)$ and $f'\colon X'\to \Spa(K, \O_K)$ be proper rigid-analytic spaces, let $h\colon X\times_K X' \to \Spa(K, \O_K)$ be the fiber product morphism, let $\cB\in \{\BB_I, \cO\}$, let $\cE\in \Perf_v(X; \cB)$, and let $\cE'\in \Perf_v(X'; \cB)$.
Then the K\"unneth morphism
    \[
    \rR f_{v, *} \cE  \otimes^\rL_{\cB} \rR f'_{v, *} \cE' \to \rR h_{v, *} (\cE\boxtimes^\rL \cE') 
    \]
is an isomorphism.

Assume furthermore that $K$ is a perfectoid field extension of $\QQ_p$. Then the K\"unneth morphism
    \[
    \rR \Gamma_v(X, \cE) \otimes^\rL_{\cB(S)} \rR \Gamma_v(X', \cE') \to \rR \Gamma_v(X\times_K X', \cE \boxtimes^\rL \cE') 
    \]
    is an isomorphism as well. 
\end{thm}
\begin{proof}
    For the first statement, we argue similarly to the proof of \cref{singular pullback compatibility}: we just need to replace \Cref{cor:proper-smooth-pushforward-is-perfect} with \Cref{cor:KF-for-proper-smooth-pushforward}. The second statement follows from the first, by the same argument as in the second paragraph
of the proof of \Cref{absolute-singular-PD} combined with the fact that $\rR\Gamma_v(\blank)$ is a (strictly) symmetric monoidal equivalence in this case.
\end{proof}

Now we give analogs of \cref{lemma:proper-pushforward-FF} and \cref{lemma:proper-pushforward-FF-coincide} for (possibly singular) proper rigid-analytic spaces.
\begin{lemma}\label{proper-pushforward-FF}
Let $f\colon X \to \Spa(K,\cO_K)$ be a proper rigid-analytic space.
Then the pullback functor $f^*_{\FFF}\colon \Perf_\an(\FFF_{\Spd(K,\cO_K)}; \cO) \to \Perf_\an(\FFF_X; \cO)$ from \cref{rmk:pullback-perf} admits a right adjoint functor $\rR f^\Perf_{\FFF, *} \colon \Perf_\an(\FFF_X; \cO) \to \Perf_\an(\FFF_{\Spd(K,\cO_K)}; \cO)$ fitting into commutative diagrams
\[ \begin{tikzcd}[center picture,column sep=5em,row sep=huge]
    \Perf_\an(\FFF_X; \cO) \arrow[d, "\rR f^\Perf_{\FFF, *}"] \arrow[r, "\pi^*_X"] & \Perf_v(X; \BB_{[1, p]})\arrow[d, "\rR f_{v, *}"] \arrow[r, shift left=0.8ex, "\otimes^\rL_{\BB_{[1, p]}, \res} \BB_{[1, 1]}"]
   \arrow[r, shift right=0.8ex, swap, "\otimes^\rL_{\BB_{[1, p]}, \wphi} \BB_{[1, 1]}"] &[2em] \Perf_v(X; \BB_{[1, 1]})\arrow[d, "\rR f_{v, *}"] \\
   \Perf_\an(\FFF_{\Spd(K,\cO_K)}; \cO) \arrow[r, "\pi^*_{\Spd(K,\cO_K)}"] & \Perf_v(\Spd(K,\cO_K); \BB_{[1, p]}) \arrow[r, shift left=0.8ex, "\otimes^\rL_{\BB_{[1, p]}, \res} \BB_{[1, 1]}"]
   \arrow[r, shift right=0.8ex, swap, "\otimes^\rL_{\BB_{[1, p]}, \wphi} \BB_{[1, 1]}"] & \Perf_v(\Spd(K,\cO_K); \BB_{[1, 1]})
\end{tikzcd} \]
and
\[ \begin{tikzcd}[column sep=7em]
    \Perf_\an(\FFF_X; \cO) \arrow[r, hook, "j_{\FFF_X}"] \arrow[d, "\rR f^\Perf_{\FFF,*}"] & \cD'(\FFF_X;\cO) \arrow[d, "\rR f_{\FFF, *}"]\\
    \Perf_\an(\FFF_{\Spd(K,\cO_K)}; \cO) \arrow[r, hook, "j_{\FFF_{\Spd(K,\cO_K)}}"] & \cD'(\FFF_{\Spd(K,\cO_K)}; \cO).
\end{tikzcd} \]
In particular, $\rR f_{\FFF,*} \colon \cD'(\FFF_X;\cO) \to \cD'(\FFF_{\Spd(K,\cO_K)};\cO)$ sends dualizable objects to dualizable objects.  
\end{lemma}
\begin{proof}
    For the existence of $\rR f^\Perf_{\FFF, *} \colon \Perf_\an(\FFF_X; \cO) \to \Perf_\an(\FFF_{\Spd(K,\cO_K)}; \cO)$ and the commutativity of the first diagram, we can use the same proof as for \cref{lemma:proper-pushforward-FF} because the $f^*_v\colon \Perf_v(\Spd(K,\cO_K);\BB_I) \to \Perf_v(X;\BB_I)$ still admit right adjoints $\rR f_{v,*} \colon \Perf_v(X; \BB_I) \to \Perf_v(\Spd(K,\cO_K);\BB_I)$ by \cref{singular-PD} and the relevant base change maps are still isomorphisms by \cref{singular pullback compatibility}.
    For the commutativity of the second diagram and the last statement, the proof of \cref{lemma:proper-pushforward-FF-coincide} applies verbatim, after once more replacing the reference to \cref{cor:proper-smooth-pushforward-is-perfect} in the last paragraph by \cref{singular pullback compatibility}.
\end{proof}

As an immediate application of \cref{proper-pushforward-FF}, we show that the proper pushforward on the Fargues--Fontaine curve is compatible with the proper pushforward on the $v$-site with respect to the $\iota_\infty^*$-functor from \cref{construction:infinity-pullback}.

\begin{cor}
\label{cor:proper-pushforward-FF-infity-pullback-singular} 
Let $f\colon X \to \Spa(K, \O_K)$ be a proper rigid-analytic space. 
Then the diagram
\begin{equation}
\begin{tikzcd}[column sep = 6em]
    \Perf_\an(\FFF_{X}; \O) \arrow[r, "\iota_\infty^*"] \arrow[d, "\rR f_{\FFF, *}"] & \Perf_v(X; \O) \arrow[d, "\rR f_{v, *}"]\\
    \Perf_\an(\FFF_{\Spd(K, \O_K)}; \O) \arrow[r, "\iota_\infty^*"] &  \Perf_v\bigl(\Spd(K, \O_K); \O\bigr)
\end{tikzcd}
\end{equation}
commutes, i.e., there is a canonical isomorphism of functors $\iota_\infty^*\circ \rR f_{\FFF, *} \xr{\sim}\rR f_{v, *}\circ \iota_\infty^*$.
\end{cor}
\begin{proof}
    The result follows immediately from \cref{proper-pushforward-FF} and \cref{thm:HT realization for proper pushforward} applied to $I=[1, p]$. 
\end{proof}

With the help of \cref{proper-pushforward-FF}, we can describe the cohomology of any local system $\bL$ on a proper rigid-analytic space $X$ over $K$ via a perfect complex:
\begin{construction}\label{construction:HXL}
    Let $f \colon X \to \Spa(K,\cO_K)$ be a proper rigid-analytic space and $\cE \in \Perf_v(X,\underline{\QQ}_p)$.
    Then $\RH_X(\cE) \in \Perf_\an(\FFF_X;\cO)$ (\cref{RH-perfect}) and we define the \emph{cohomology complex} of $\cE$ as
    \[ \cH_X(\cE) \colonequals \rR f_{\FFF,*}\bigl(\RH_X(\cE)\bigr) \in \Perf_\an(\FFF_{\Spd(K,\cO_K)};\cO); \]
    here, the containment in $\Perf_\an(\FFF_{\Spd(K,\cO_K)};\cO)$ follows from \cref{proper-pushforward-FF}.
    This perfect complex recovers the cohomology of $\cE$ because
    \[ \Sol_{\Spd(K,\cO_K)}\bigl(\cH_X(\cE)\bigr) \simeq \Sol_{\Spd(K,\cO_K)}\bigl(\RH_{\Spd(K,\cO_K)}\bigl(\rR f_{v,*}(\cE)\bigr)\bigr) \simeq \rR f_{v,*}(\cE) \]
    by virtue of \cref{lemma:dr-pushforward} and \cref{prop:dr-rh}.
    When $\bL$ is a $\underline{\QQ}_p$-local system on $X$, we also set $\cH_X(\bL) \colonequals \cH_X\bigl(\bL[0]\bigr)$.
\end{construction}
Lastly, we give the analog of \cref{cor:tor-amplitude-2} and \cref{cor:tor-amplitude-qp-smoothproper-2} for singular proper rigid-analytic spaces.
\begin{corollary}
\label{singular-tor-amplitude}
Let $f \colon X \to \Spa(K,\cO_K)$ be a proper rigid-analytic space of dimension $d$, let $\cA \in \{\underline{\ZZ}_p,\cO^+,\AA_I\}$ and $\cB \colonequals \cA\bigl[\tfrac{1}{p}\bigr]$, and let $\cE \in \Perf_v(X; \cB) \cap \cD^{[r,r']}(X_v;\cB)$ (resp. $\cE\in \Perf_v^{[r, r']}(X; \cB)$) for some $r, r'\in \ZZ\cup \{-\infty, \infty\}$.
Then $\rR f_{v,*} \cE$ lies in $\cD^{[r,r'+2d]}_{\BBC}(\Spd(K,\cO_K)_v;\underline{\QQ}_p)$ if $\cA = \underline{\ZZ}_p$ and in $\Perf_v(\Spd(K,\cO_K);\cB) \cap \cD^{[r,r'+2d]}(\Spd(K,\cO_K)_v;\cB)$  (resp. $\cE\in \Perf_v^{[r, r'+2d]}(\Spd(K, \O_K)_v; \cB)$) if $\cA \in \{\cO^+,\AA_I\}$.
\end{corollary}
\begin{proof}
    We give a proof for the case of $\cE\in \Perf_v(X; \cB) \cap \cD^{[r,r']}(X_v;\cB)$, the case of $\cE\in \Perf_v^{[r, r']}(X; \cB)$ being analogous.

    Thanks to \cref{singular-PD}, we already know that $\rR f_{v,*}\cE$ lies in $\cD_{\BBC}(\Spd(K,\cO_K)_v;\underline{\QQ}_p)$ if $\cA = \underline{\ZZ}_p$ and in $\Perf_v(\Spd(K,\cO_K);\cB)$ if $\cA \in \{ \cO^+,\AA_I \}$.
    Thus, we only need to show that $\rR f_{v,*}\cE \in \cD^{[r,r'+2d]}(\Spd(K,\cO_K)_v;\cB)$.
    For this, we may assume as in the proof of \cref{singular-PD} that $X$ is reduced and that the statement has been proven in dimensions $<d$.
    
    By \cref{dualizing-cohomology}, we can find the modification square from \cref{modification-square}
    \[ \begin{tikzcd}
        Z' \arrow[r,"i'"] \arrow[d,"\pi'"] \arrow[rd,"h"] & X' \arrow[d,"\pi"] \\
        Z \arrow[r,"i"] & X
    \end{tikzcd} \]
    in which $\pi$ is a regular admissible modification and $\dim Z,\dim Z' < \dim X$.
    As in \cref{v-MV} in the proof of \cref{singular-PD}, this yields a fiber sequence
    \[ \rR f_{v,*}\cE \xrightarrow{\bigl((-1)\cdot \rR f_{v,*}\eta_\pi , \rR f_{v,*} \eta_i \bigr)} \rR (f \circ \pi)_{v,*}\pi^*_v\cE \oplus \rR (f \circ i)_{v,*}i^*_v\cE \xrightarrow{\rR (f \circ \pi)_{v,*}\eta_{i'} \oplus \rR (f \circ i)_{v,*}\eta_{\pi'}} \rR (f \circ h)_{v,*}h^*_v\cE. \]
    Since $\dim Z,\,\dim Z' < d$, we have $\rR (f \circ i)_{v,*}i^*_v\cE,\, \rR (f \circ h)_{v,*}h^*_v\cE \in \cD^{[r,r'+2d-2]}(\Spd(K,\cO_K)_v;\cB)$.
    On the other hand, \cref{cor:tor-amplitude-2} and \cref{cor:tor-amplitude-qp-smoothproper-2} show that $\rR (f \circ \pi)_{v,*}\pi^*_v\cE \in \cD^{[r,r'+2d]}(\Spd(K,\cO_K)_v;\cB)$.
    This yields the desired cohomological bounds.
\end{proof} 
\begin{remark}
    When $K$ is perfectoid and $\cB \in \{\cO^+,\AA_I\}$, \cref{singular-tor-amplitude} also implies that $\rR\Gamma_v(X,\cE) \in \Perf(B) \cap \cD^{[r,r'+2d]}(B)$.
    The proof of this absolute version is again analogous to the proof of \cref{absolute-singular-PD}.
\end{remark}

We can now summarize the main properties of the $\cH_X(\cE)$ construction. 
\begin{thm}\label{thm:main-thm-hx} Let $f \colon X\to \Spa(K, \O_K)$ be a proper rigid-analytic space of dimension $d$ and let $\cE\in \Perf^{[r, r']}_v(X; \uQ_p)$ for some $r,r'\in \ZZ \cup \{-\infty, \infty\}$. Then the following statements hold: 
\begin{enumerate}[leftmargin=*,label={\upshape{(\roman*)}}]
    \item\label{thm:main-thm-hx-1} there is a functorial isomorphism $\Sol\bigl(\cH_X(\cE)\bigr) \simeq \rR f_{v, *}\cE$;
    \item\label{thm:main-thm-hx-1.5} there is a functorial isomorphism $\iota_\infty^* \cH_X(\cE) \simeq \rR f_{v, *}\bigl( \cE \otimes^{\rL}_{\uQ_p} \cO\bigr)$;
    \item\label{thm:main-thm-hx-2} if $X$ is smooth of equidimension $d$, there is a functorial duality isomorphism $\cH_X(\cE)^{\vee} \simeq \cH_X(\cE^{\vee}(d)[2d])$;
    \item\label{thm:main-thm-hx-3} the complex $\cH_X(\cE)$ lies in $\Perf_\an^{[r, r'+2d]}(\FFF_{\Spd(K, \O_K)}; \cO)$;
    \item\label{thm:main-thm-hx-4} if $(X', \cE')$ is another such pair, then there is a natural K\"unneth formula isomorphism 
    \[
    \cH_X(\cE) \otimes^\rL_\O \cH_{X'}(\cE') \xr{\sim} \cH_{X\times_K X'}(\cE \boxtimes \cE').
    \]
\end{enumerate}
\end{thm}

When $K$ is a perfectoid field, the second part of \cref{construction:infinity-pullback} allows us to interpret \cref{thm:main-thm-hx-1.5} as an equivalence $\iota_\infty^* \cH_X(\cE) \simeq \rR\Gamma_v(X, \cE\otimes^{\rL}_{\uQ_p} \cO)$. 
\begin{proof}
    Part~\cref{thm:main-thm-hx-1} follows from \cref{lemma:dr-pushforward} and \cref{prop:dr-rh}. Part~\cref{thm:main-thm-hx-1.5} follows from \cref{cor:proper-pushforward-FF-infity-pullback-singular}. Part~\cref{thm:main-thm-hx-2} follows from \cref{thm:geometric-duality-curve}. Part~\cref{thm:main-thm-hx-3} follows from \cref{proper-pushforward-FF} and \cref{singular-tor-amplitude}. Finally, Part~\cref{thm:main-thm-hx-4} follows from \cref{thm:KF-for-proper-pushforward}. 
\end{proof}

\subsection{Arithmetic duality}\label{arithmetic-duality}
From now on, we assume that $K$ is a finite extension of $\mathbf{Q}_p$.
In this subsection, we prove honest finiteness and duality
for the absolute cohomology of $\uQ_p$-local systems on proper rigid spaces over $K$ in this arithmetic
setting.
We proceed again in two steps. 
First, we establish duality on the Fargues--Fontaine curve for the  
Galois equivariant perfect complexes representing $v$-cohomology of the given 
$\Q_p$-local system. This is done 
by Galois descent from the already established geometric duality on the Fargues--Fontaine curve;
the descent  itself uses  Fontaine's theory of almost $\mathbf{C}_p$-representations 
(see \cite{FCp,Fontaine-Almost}). Then, 
to deduce the final form of arithmetic duality we again evoke \cite[Cor.~3.11]{ALB}. 
The geometric duality from \cref{geometric-duality-Qp} paired with the Tate duality for almost 
${\mathbf C}_p$-representations proved by Fontaine in \cite{FCp} 
allows us  to prove arithmetic duality in the absolute case. 
Our arguments follow closely those used by Zhenghui Li in \cite{ZL}.

Let $C\coloneqq{\mathbf C}_p,$ and let $G_K\coloneqq\mathrm{Gal}(\overline{K}/K)$ be the Galois group of $K$.
We denote by 
\[
\rR\ud{\Gamma}(G_K, \blank) \coloneqq \rR\ud{\Hom}_{\ud{\ZZ}[\ud{G}_K]}(\ud{\ZZ}, \blank) \colon \cD(\ud{\ZZ}[\ud{G}_K]) \to \cD(\ud{\ZZ})
\]
the condensed group cohomology of $G_K$ (see \cite[Appendix B]{bosco-master}). We denote $\rR\Gamma(G_K, \blank) \coloneqq \rR\Hom_{\ud{\ZZ}[\ud{G}_K]}(\ud{\ZZ}, \blank)$ its underlying complex of abelian groups. 

We begin by defining the arithmetic trace map and the associated duality
morphisms.

\begin{definition}
Let $f\colon X \to \Spa(K, \O_K)$ be a smooth proper rigid-analytic space of
equidimension $d$ and let $\cE\in\Perf_v(X;\uQ_p)$. Write
$f_C\colon X_C\to\Spa(C,\O_C)$ for the base change of $f$.
\begin{enumerate}[leftmargin=*]
    \item The trace map $\tr_{f}^{\QQ_p}$ from
    \cref{trace-maps}~\cref{trace-maps-Qp} induces a $G_K$-equivariant
    geometric trace map
    \begin{align*}
    \tr_{X_C}^{\QQ_p}\colon
    \rR\uGamma_v\bigl(X_C,\uQ_p(d)[2d]\bigr)
    &\simeq
    \rR\uGamma_v\Bigl(
        \Spd(C,\O_C),
        \rR f_{C,v,*}\uQ_p(d)[2d]
    \Bigr)\\
    &\xr{\rR\uGamma_v\bigl(\Spd(C,\O_C),\tr_{f}^{\QQ_p}\bigr)}
    \rR\uGamma_v\bigl(\Spd(C,\O_C),\uQ_p\bigr)
    \simeq \uQ_p.
    \end{align*}
    We define the \emph{arithmetic trace map}
    \[
    \tr_X\colon
    \rR\Gamma_v\bigl(X,\uQ_p(d+1)[2d+2]\bigr)
    \longrightarrow \QQ_p
    \]
    as the composite
    \begin{align*}
    \rR\Gamma_v\bigl(X,\uQ_p(d+1)[2d+2]\bigr)
    &\simeq
    \rR\Gamma\Bigl(
        G_K,
        \rR\underline{\Gamma}_v
        \bigl(X_C,\uQ_p(d)[2d]\bigr)(1)[2]
    \Bigr)\\
    &\xr{\rR\Gamma\bigl(G_K,\tr_{X_C}^{\QQ_p}(1)[2]\bigr)}
    \rR\Gamma\bigl(G_K,\uQ_p(1)[2]\bigr)
    \xr{\tr_K}
    \QQ_p.
    \end{align*}
    Here, the first isomorphism comes
    from \cref{G-torsor and cohomology}, applied to the pro-\'etale
    $G_K$-torsor $X_C\to X$, and $\tr_K$ is the usual Galois trace map
    from Tate local duality.

    \item We define the \emph{duality map}
    \[
    \PD_X(\cE)\colon
    \rR\Gamma_v\bigl(X,\cE^\vee(d+1)[2d+2]\bigr)
    \longrightarrow
    \rR\Gamma_v(X,\cE)^\vee
    \]
    as the composite
    \begin{align*}
    \rR\Gamma_v\bigl(X,\cE^\vee(d+1)[2d+2]\bigr)
    &\longrightarrow
    \rR\Hom_{\QQ_p}\Bigl(
        \rR\Gamma_v(X,\cE),
        \rR\Gamma_v\bigl(X,\uQ_p(d+1)[2d+2]\bigr)
    \Bigr)\\
    &\xr{\tr_X \circ \blank}
    \rR\Hom_{\QQ_p}\bigl(\rR\Gamma_v(X,\cE),\QQ_p\bigr) = \rR\Gamma_v(X, \cE)^{\vee}.
    \end{align*}
    The first arrow is the evaluation map from 
    \cref{Ev-Adj} combined with the identification $\cE^{\vee}(d+1)[2d+2]\simeq \rR\ud{\Hom}_{\uQ_p}\bigl(\cE, \uQ_p(d+1)[2d+2]\bigr)$. 
\end{enumerate}
\end{definition}

Now we are ready to prove arithmetic duality for smooth proper rigid-analytic spaces. We will then generalize this result to all proper rigid-analytic spaces.

\begin{theorem}[Arithmetic duality]\label{thm:arithmetic-duality-smooth}
 Let $f\colon X \to \Spa(K, \O_K)$ be a smooth proper rigid-analytic space of equidimension $d$ and let $\cE\in \Perf_v(X;\underline{\Q}_p)$. Then the duality morphism 
 \[
  \PD_X(\cE) \colon\rR\Gamma_v\bigl(X, \cE^{\vee}(d+1)[2d+2]\bigr) \to \rR\Gamma_v(X, \cE)^\vee 
 \]
 is an isomorphism. 
\end{theorem}
\begin{proof}
{\it Digression.} Before turning to the proof, we recall some facts about the evaluation morphisms that will be used below. The solution functor $\Sol_C\colon
\Perf_\an(\FFF_C;\O)
\longrightarrow
\cD\bigl((\Spd C)_v;\uQ_p\bigr)$ is lax symmetric monoidal. Hence, for any
$\cP,\cP'\in\Perf_\an(\FFF_C;\O)$, there is a functorial evaluation morphism
\[
\Ev_{\Sol}(\cP,\cP')\colon
\Sol_C\bigl(\rR{\cHom}_{\FFF_C}(\cP,\cP')\bigr)
\longrightarrow
\rR{\cHom}_{\cD\bigl((\Spd C)_v;\uQ_p\bigr)}
\bigl(\Sol_C(\cP),\Sol_C(\cP')\bigr).
\]
By \cref{ff1}, this morphism is an isomorphism. 

Similarly, derived condensed global sections carry a natural evaluation morphism
\[
\Ev_{\rR\underline{\Gamma}_v}(\cP,\cP')\colon
\rR\underline{\Gamma}_v\Bigl(
\Spd C,
\rR{\cHom}_{\cD((\Spd C)_v;\uQ_p)}(\cP,\cP')
\Bigr)
\longrightarrow
\rR{\cHom}_{\cD(\uQ_p)}
\Bigl(
\rR\underline{\Gamma}_v(\Spd C,\cP),
\rR\underline{\Gamma}_v(\Spd C,\cP')
\Bigr)
\]
for any $\cP,\cP'\in\cD((\Spd C)_v;\uQ_p)$. Unlike
$\Ev_{\Sol}$, this morphism need not be an isomorphism.

Using the solution functor, we define the condensed global-sections functor on the Fargues--Fontaine curve
\[
\rR\underline{\Gamma}_\an(\FFF_C,\blank)\colon
\Perf_\an(\FFF_C;\O)
\longrightarrow
\cD(\uQ_p)
\]
to be the composite
\[
\Perf_\an(\FFF_C;\O)
\xrightarrow{\Sol_C}
\cD\bigl((\Spd C)_v;\uQ_p\bigr)
\xrightarrow{\rR\underline{\Gamma}_v(\Spd C,\blank)}
\cD(\uQ_p).
\]
Since both functors in this composite are lax symmetric monoidal, so is
$\rR\underline{\Gamma}_\an(\FFF_C,\blank)$. It therefore carries a canonical evaluation morphism
\[
\Ev_{\rR\underline{\Gamma}_\an}(\cP,\cP')\colon
\rR\underline{\Gamma}_\an\bigl(
\FFF_C,
\rR{\cHom}_{\FFF_C}(\cP,\cP')
\bigr)
\longrightarrow
\rR{\cHom}_{\cD(\uQ_p)}
\Bigl(
\rR\underline{\Gamma}_\an(\FFF_C,\cP),
\rR\underline{\Gamma}_\an(\FFF_C,\cP')
\Bigr).
\]

Finally, \cref{prop:dr-rh,lemma:dr-pushforward,construction:HXL} provide, for every
$\cP\in\Perf_v(X_C;\uQ_p)$, a canonical equivalence
\[
c_f(\cP)\colon
\rR f_{C,v,*}\cP
\xr{\sim}
\Sol_C\bigl(\cH_{X_C}(\cP)\bigr).
\]

\emph{Step~$1$. Reduction to equivariant duality on the Fargues--Fontaine curve.}
Consider the following diagram:
\begin{equation}\label{paris3}
\hspace{1.5em}\begin{tikzcd}[scale cd=.9,center picture]
\arrow[
  phantom,
  from=1-1,
  to=4-2,
  "{1}"{circled number}
]
\arrow[
  phantom,
  from=1-2,
  to=2-3,
  "{2}"{circled number}
]
\arrow[
  phantom,
  from=2-2,
  to=3-3,
  "{3}"{circled number}
]
\arrow[
  phantom,
  from=3-2,
  to=4-3,
  "{4}"{circled number}
]
\rR\Gamma\bigl(G_K,\rR\underline{\Gamma}_v\bigl(X_C,\cE_C^\vee(d)[2d]\bigr)(1)[2]\bigr) \arrow[r,equals] \arrow[d,equals]  &[-2em]
\rR\Gamma\bigl(G_K,\rR\underline{\Gamma}_v\bigl(\Spd C,\rR f_{C,v,*}\cE_C^\vee(d)[2d]\bigr)(1)[2]\bigr) \arrow[d,"{\PD_{f_C}^{\QQ_p}(\cE_C)}","\sim"'{sloped}] \arrow[r,"{c_f\bigl(\cE_C^\vee(d)[2d]\bigr)}","\sim"'] &[3em]
\rR\Gamma\bigl(G_K,\rR\underline{\Gamma}_\an\bigl(\FFF_C,\cH_{X_C}(\cE_C^\vee(d)[2d])\bigr)(1)[2]\bigr) \arrow[d,"{\PD_{f_C}^{\FFF}(\cE_C)}","\sim"'{sloped}] \\
\rR\Gamma_v\bigl(X,\cE^\vee(d+1)[2d+2]\bigr) \arrow[d,"{\PD_X(\cE)}"] &
\rR\Gamma\bigl(G_K,\rR\underline{\Gamma}_v\bigl(\Spd C,(\rR f_{C,v,*}\cE_C)^\vee\bigr)(1)[2]\bigr) \arrow[d,"{\Ev_{\rR\underline{\Gamma}_v}(\rR f_{C,v,*}\cE_C,\uQ_p)}"] &
\rR\Gamma\bigl(
G_K,
\rR\underline{\Gamma}_\an
\bigl(\FFF_C,\cH_{X_C}(\cE_C)^\vee\bigr)(1)[2]
\bigr)
  \arrow[l,"{\Ev_{\Sol}(\cH_{X_C}(\cE_C),\O)}","\sim"']
  \arrow[d,"{\Ev_{\rR\underline{\Gamma}_\an}
  (\cH_{X_C}(\cE_C),\O)}"]
\\
\rR\Gamma_v(X,\cE)^\vee
  \arrow[d,equals]
&
\rR\Gamma\bigl(
G_K,
\rR\underline{\Gamma}_v
\bigl(\Spd C,\rR f_{C,v,*}\cE_C\bigr)^\vee(1)[2]
\bigr)
  \arrow[d,"{\PD_K}"]
&
\rR\Gamma\bigl(
G_K,
\rR\underline{\Gamma}_\an
\bigl(\FFF_C,\cH_{X_C}(\cE_C)\bigr)^\vee(1)[2]
\bigr)
  \arrow[d,"{\PD_K}"]
  \arrow[l,"{c_f(\cE_C)^\vee}","\sim"']
\\
\rR\Gamma\bigl(
G_K,\rR\underline{\Gamma}_v(X_C,\cE_C)
\bigr)^\vee
  \arrow[r,equals]
&
\rR\Gamma\bigl(
G_K,
\rR\underline{\Gamma}_v
\bigl(\Spd C,\rR f_{C,v,*}\cE_C\bigr)
\bigr)^\vee
&
\rR\Gamma\bigl(
G_K,
\rR\underline{\Gamma}_\an
\bigl(\FFF_C,\cH_{X_C}(\cE_C)\bigr)
\bigr)^\vee
  \arrow[l,"{c_f(\cE_C)^\vee}","\sim"']
\end{tikzcd}\end{equation}
The labels in the diagram only indicate the morphisms from which the corresponding arrows are induced. We refer to
\cref{defn:Qp-duality-map}~\cref{defn:Qp-duality-map-3} for the definition of
$\PD^{\QQ_p}_{f_C}$ and to
\cref{defn:FF-trace-maps}~\cref{defn:FF-trace-maps-3} for that of
$\PD^{\FFF}_{f_C}$. The morphism $\PD_K$ is the Tate duality morphism. More precisely, for
$V\in\cD\bigl(\ud{\QQ}_p[\ud{G}_K]\bigr)$, we define $\PD_K(V)$ to be the composite
\[
\rR\Gamma\bigl(G_K,V^\vee(1)[2]\bigr)
\xr{\Ev}
\rR\Hom_{\cD(\QQ_p)}
\Bigl(
\rR\Gamma(G_K,V),
\rR\Gamma(G_K,\uQ_p(1)[2])
\Bigr)
\xr{\tr_K\circ\blank}
\rR\Hom_{\cD(\QQ_p)}
\bigl(\rR\Gamma(G_K,V),\QQ_p\bigr),
\]
where $\Ev$ is the evaluation morphism and $\tr_K$ is the Galois trace from Tate's duality. 

We claim that Diagram~\cref{paris3} commutes. Region~{\upshape(1)}
commutes by the definition of $\PD_X(\cE)$. The commutativity of
Square~{\upshape(2)} follows from the compatibility
\[
\Sol_C\bigl(\PD_{f_C}^{\FFF}(\cE_C)\bigr)
=
\PD_{f_C}^{\QQ_p}(\cE_C),
\]
established in the proof of \cref{geometric-duality-Qp}. Square~{\upshape(3)}
commutes by \cref{Ev-Adj-composition} and the identity $\rR\underline{\Gamma}_\an(\FFF_C,\blank) = \rR\underline{\Gamma}_v\bigl(\Spd C,\Sol_C(\blank) \bigr)$. Finally, Square~{\upshape(4)} commutes by the naturality of the Tate duality morphism $\PD_K$.

Since all the horizontal arrows are isomorphisms, it suffices to prove that the composite down the rightmost column of Diagram~\cref{paris3},
\[
\PD_K
\circ
\Ev_{\rR\underline{\Gamma}_\an}
\bigl(\cH_{X_C}(\cE_C),\O\bigr)
\circ
\PD_{f_C}^{\FFF}(\cE_C),
\]
is an isomorphism. By \cref{thm:geometric-duality-curve}, the morphism
$\PD_{f_C}^{\FFF}(\cE_C)$ is an isomorphism. It therefore remains to prove that
\[
\PD_K
\circ
\Ev_{\rR\underline{\Gamma}_\an}
\bigl(\cH_{X_C}(\cE_C),\O\bigr)
\]
is an isomorphism.

Now we consider the $\infty$-category $\mathcal{P}(G_K)$ of $G_K$-equivariant perfect complexes on $\FFF_C$ (see \cite[Sec. 4.1.4]{ZL} for the definition and basic properties). For brevity, set
\[
\rR\Gamma_{\cP(G_K)}(\cP)
\coloneqq
\rR\Hom_{\cP(G_K)}(\cO,\cP)
\]
for $\cP\in\cP(G_K)$. By \cite[Rmk.~4.19]{ZL}, there is a canonical isomorphism
\[
\rR\Gamma_{\cP(G_K)}(\cP) = \rR\Hom_{\cP(G_K)}(\cO,\cP)
\xr{\sim}
\rR\Gamma\Bigl(
G_K,
\rR\underline{\Gamma}_\an(\FFF_C,\cP)
\Bigr)
\]
for every $\cP\in\cP(G_K)$. The category $\cP(G_K)$ also carries a natural duality morphism.  We have a trace map
\[
\tr_{\cP(G_K)}\colon
\rR\Gamma_{\cP(G_K)}\bigl(\cO(1)[2]\bigr)
\simeq
\rR\Gamma\Bigl(
G_K,
\rR\underline{\Gamma}_\an(\FFF_C,\cO(1)[2])
\Bigr)
\xr{\sim}
\rR\Gamma\bigl(G_K,\uQ_p(1)[2]\bigr)
\xr{\tr_K}
\QQ_p.
\]
Together with the lax symmetric monoidal structure on
$\rR\Gamma_{\cP(G_K)}$, this trace induces a duality morphism
\[
\PD_{\cP(G_K)}(\cP)\colon
\rR\Gamma_{\cP(G_K)}\bigl(\cP^\vee(1)[2]\bigr)
\xr{\Ev}
\rR\Hom_{\cD(\QQ_p)}
\Bigl(
\rR\Gamma_{\cP(G_K)}(\cP),
\rR\Gamma_{\cP(G_K)}\bigl(\cO(1)[2]\bigr)
\Bigr)
\xr{\tr_{\cP(G_K)}\circ\blank}
\rR\Gamma_{\cP(G_K)}(\cP)^\vee.
\]

By construction, these morphisms fit into the following commutative diagram:
\begin{equation}\label{hk1}
\xymatrix{
\rR\Gamma\Bigl(
G_K,
\rR\underline{\Gamma}_\an
\bigl(\FFF_C,\cH_{X_C}(\cE_C)^\vee\bigr)(1)[2]
\Bigr)
\ar[d]^{
\Ev_{\rR\underline{\Gamma}_\an}
(\cH_{X_C}(\cE_C),\O)
}
&
\rR\Gamma_{\cP(G_K)}
\bigl(\cH_{X_C}(\cE_C)^\vee(1)[2]\bigr)
\ar[l]_-{\sim}
\ar[dd]^-{\PD_{\cP(G_K)}}
\\
\rR\Gamma\Bigl(
G_K,
\rR\underline{\Gamma}_\an
\bigl(\FFF_C,\cH_{X_C}(\cE_C)\bigr)^\vee(1)[2]
\Bigr)
\ar[d]^{\PD_K}
\\
\rR\Gamma\Bigl(
G_K,
\rR\underline{\Gamma}_\an
\bigl(\FFF_C,\cH_{X_C}(\cE_C)\bigr)
\Bigr)^\vee
\ar[r]_-{\sim}
&
\rR\Gamma_{\cP(G_K)}
\bigl(\cH_{X_C}(\cE_C)\bigr)^\vee.
}
\end{equation}
Thus, we reduce the question to showing that $\PD_{\cP(G_K)}(\cP)$ is an isomorphism for any $\cP\in \cP(G_K)$.

{\it Step~$2$. Finish the proof.} We show that $\PD_{\cP(G_K)}(\cP)$ is an isomorphism by reducing this to the duality theorem in almost $\CC_p$-representations developed in \cite{FCp}.

For this, we recall that Fontaine defined the derived category $\cD(\cM(G_K))$ of $G_K$-equivariant coherent sheaves on the algebraic Fargues--Fontaine curve $\FFF_C^\alg$. There is an evident analytification functor $\cD^b(\cM(G_K)) \to \cP(G_K)$, which is an equivalence. Indeed, the calculations of Ext groups in the two categories
\cite[Rmk.~4.19 and Eqn.~(4.12)]{ZL} show that it is fully faithful, while the GAGA theorem
\cite[Th.~II.2.6 and Prop.~II.2.7]{Fargues-Scholze} shows that it is essentially surjective.

We may therefore work instead with $\cD^b(\cM(G_K))$. Let
$\cD^b(\cC(G_K))$ denote Fontaine's bounded derived category of almost
$\CC_p$-representations, defined in \cite{FCp}. By \cite[Th.~4.23]{ZL}, derived global sections, together with the analytification equivalence above, induce an equivalence
\[
\rR\Gamma(\FFF_C,\blank)\colon
\cP(G_K) \simeq \cD^b(\cM(G_K))
\xr{\sim}
\cD^b(\cC(G_K)).
\]

Now, similarly to the case of $\cP(G_K)$, we set 
\[
\rR\Gamma_{\cC(G_K)}(V) \coloneqq \rR\Hom_{\cD^b(\cC(G_K))}(\Q_p, V)
\]
for any $V\in \cD^b(\cC(G_K))$. By \cite[Prop.~6.7(iii)]{FCp}, there is an isomorphism
\[
\rR\Gamma_{\cC(G_K)}\bigl(\QQ_p(1)[2]\bigr)
\simeq
\rR\Gamma_{\cont}\bigl(G_K,\QQ_p(1)[2]\bigr).
\]
Tate's trace therefore defines a trace map
\[
\tr_{\cC(G_K)}\colon
\rR\Gamma_{\cC(G_K)}\bigl(\QQ_p(1)[2]\bigr)
\simeq
\rR\Gamma_{\cont}\bigl(G_K,\QQ_p(1)[2]\bigr)
\xr{\tr_K}
\QQ_p.
\]
This trace, in turn, defines a duality morphism
\[
\PD_{\cC(G_K)}(V)\colon
\rR\Gamma_{\cC(G_K)}\bigl(V^\vee(1)[2]\bigr)
\longrightarrow
\rR\Gamma_{\cC(G_K)}(V)^\vee
\]
for every $V\in\cD^b(\cC(G_K))$.

Under the equivalence
$\rR\Gamma(\FFF_C,\blank)$, the trace transported from
$\cP(G_K)$ and the trace $\tr_{\cC(G_K)}$ are two morphisms
\[
\rR\Gamma_{\cC(G_K)}\bigl(\QQ_p(1)[2]\bigr)
\longrightarrow
\QQ_p.
\]
Both factor through $\Hh^2_{\cC(G_K)}(\QQ_p(1))[0]$, and each induces an isomorphism
\[
\Hh^2_{\cC(G_K)}(\QQ_p(1))
\simeq
\Hh^2_{\cont}(G_K,\QQ_p(1))
\xr{\sim}
\QQ_p.
\]
So they must differ by a nonzero scalar $c\in\QQ_p^\times$.

Now the equivalence
\[
\rR\Gamma(\FFF_C,\blank)\colon
\cP(G_K)
\longrightarrow
\cD^b(\cC(G_K))
\]
is compatible with the duals and it intertwines their trace maps up to the nonzero scalar $c$. Therefore, it suffices to prove that
\[
\PD_{\cC(G_K)}(V)\colon
\rR\Gamma_{\cC(G_K)}\bigl(V^\vee(1)[2]\bigr)
\longrightarrow
\rR\Gamma_{\cC(G_K)}(V)^\vee
\]
is an isomorphism for every $V\in\cD^b(\cC(G_K))$. For this, we can first assume that $V\in \cC(G_K)$ by a simple induction argument. 
Then we can check that $\PD_{\cC(G_K)}(V)$ is an isomorphism on each individual Ext group,
where it has been proven in \cite[Prop.~6.11]{FCp}. This finishes the proof.
\end{proof}

\begin{cor}\label{cor:coh-dimension-arithmetic-case}  Let $f\colon X \to \Spa(K, \O_K)$ be a smooth proper rigid-analytic space of dimension $d$ and let $\cE\in \Perf_v(X;\underline{\Q}_p) \cap \cD^{[r, r']}(X_v; \uQ_p)$. Then $\rR\Gamma_v(X, \cE)$ lies in $\cD^{[r, r'+2d+2]}_{\mathrm{coh}}(\QQ_p)$.
\end{cor}
\begin{proof}
    First, we note that $X$ admits a clopen decomposition $X = \sqcup_{n\leq d} X_n$ such that $X_n$ is smooth proper of equidimension $n\leq d$. Hence, without loss of generality, we can assume that $X=X_d$ is of equidimension $d$. Then \cref{thm:arithmetic-duality-smooth} implies that $\rR\Gamma_v(X, \cE) \simeq \rR\Gamma_v(X, \cE)^{\vee\vee}$. This formally implies that $\rR\Gamma_v(X, \cE) \in \cD_{\mathrm{coh}}(\QQ_p)$. We also clearly have that $\rR\Gamma_v(X, \cE) \in \cD^{\geq r}(\QQ_p)$, so it suffices to show that $\rR\Gamma_v(X, \cE) \in \cD^{\leq r'+2d+2}(\QQ_p)$. \cref{cor:local-system-tor-ampl} ensures that $\cE\in \Perf^{[r, r']}_v(X; \uQ_p)$, so $\cE^{\vee}\in \Perf^{[-r', -r]}_v(X; \uQ_p)$. In particular, we have $\cE^{\vee} \in \cD^{[-r', -r]}(X_v; \uQ_p)$. Hence, the desired cohomological bound follows formally from the observation that 
    \[
    \rR\Gamma_v(X, \cE)^{\vee} \simeq \rR\Gamma_v(X, \cE^{\vee}(d+1)[2d+2])
    \]
    lies in $\cD^{\geq -r'-2d-2}(\QQ_p)$. 
\end{proof}

\begin{remark}
\begin{enumerate}[leftmargin=*]
\item The functional analytic arguments in the proof of \cref{thm:arithmetic-duality-smooth} are subtle and are the reason why this proof is longer than one would (naively) expect. This is due to the fact that it is not clear how to make the category $\mathcal{C}(G_K)$ of Fontaine compatible with the condensed formalism (see also \cite[Rmk.~4.22]{ZL}). In particular, it seems unlikely to us that the arrows in the second and third rows of diagram \ref{paris3} are isomorphisms. 
\item One can also deduce \Cref{thm:arithmetic-duality-smooth} from \cite{ALBM}; cf.\ \cite[Rmk.~6.2.3]{ALBM}.
\end{enumerate}
\end{remark}

 The arithmetic duality from \cref{thm:arithmetic-duality-smooth} can be generalized to singular varieties. For this, we first need to make some preliminary definitions. 

 \begin{definition}\label{defn:arithmetic-dualizing-complex}
    Let $X$ be a rigid-analytic space over $\Spa(K,\cO_K)$.  
    \begin{enumerate}[leftmargin=*]
        \item The \emph{arithmetic dualizing complex} of $X$ is $\omega^{\mathrm{ar}}_X \colonequals \omega^{\QQ_p}_X(1)[2]$, where $\omega^{\QQ_p}_X$ is the $\QQ_p$-dualizing complex from \cref{needed1}.
        \item The \emph{arithmetic Verdier duality} functor for $X$ is $\DD^{\mathrm{ar}}_X(\blank) \colonequals \rR\cHom\bigl(\blank,\omega^{\mathrm{ar}}_X\bigr) \colon \cD(X_v;\uQ_p)^\op \to \cD(X_v;\uQ_p)$.
    \end{enumerate}
\end{definition}

\begin{definition}\label{defn:arithmetic-singular-duality-morphism}
Let $f\colon X \to \Spa(K, \O_K)$ be a proper rigid-analytic space and let $\cE\in\Perf_v(X;\uQ_p)$. Write
$f_C\colon X_C\to\Spa(C,\O_C)$ for the base change of $f$.
\begin{enumerate}[leftmargin=*,label=\upshape{(\arabic*)}]
    \item\label{defn:arithmetic-singular-duality-morphism-1} The trace morphism $\Tr_{f}^{\QQ_p}$ from
    \cref{BI-proper-trace} induces a $G_K$-equivariant
    geometric trace map
    \begin{align*}
    \Tr_{X_C}^{\QQ_p}\colon
    \rR\uGamma_v\bigl(X_C,\omega^{\QQ_p}_{X_C}\bigr)
    &\simeq
    \rR\uGamma_v\Bigl(
        \Spd(C,\O_C),
        \rR f_{C,v,*}\omega^{\QQ_p}_{X_C}
    \Bigr)\\
    &\xr{\rR\uGamma_v\bigl(\Spd(C,\O_C),\Tr_{f}^{\QQ_p}\bigr)}
    \rR\uGamma_v\bigl(\Spd(C,\O_C),\uQ_p\bigr)
    \simeq \uQ_p.
    \end{align*}
    We define the \emph{arithmetic trace map}
    \[
    \Tr_X\colon
    \rR\Gamma_v\bigl(X,\omega_X^{\mathrm{ar}}\bigr)
    \longrightarrow \QQ_p
    \]
    as the composite
    \begin{align*}
    \rR\Gamma_v\bigl(X, \omega_X^{\mathrm{ar}}\bigr) \simeq \rR\Gamma_v\bigl(X,\omega^{\QQ_p}_X(1)[2]\bigr)
    &\simeq
    \rR\Gamma\Bigl(
        G_K,
        \rR\underline{\Gamma}_v
        \bigl(X_C,\omega^{\QQ_p}_{X_C}\bigr)(1)[2]
    \Bigr)\\
    &\xr{\rR\Gamma\bigl(G_K,\Tr_{X_C}^{\QQ_p}(1)[2]\bigr)}
    \rR\Gamma\bigl(G_K,\uQ_p(1)[2]\bigr)
    \xr{\tr_K}
    \QQ_p.
    \end{align*}
    Here the first isomorphism comes
    from \cref{G-torsor and cohomology}, applied to the pro-\'etale
    $G_K$-torsor $X_C\to X$, and $\tr_K$ is the usual Galois trace map
    from Tate local duality.

    \item\label{defn:arithmetic-singular-duality-morphism-2} We define the \emph{duality map}
    \[
    \PD_X(\cE)\colon
    \rR\Gamma_v\bigl(X, \DD^{\mathrm{ar}}_X(\cE)\bigr)
    \longrightarrow
    \rR\Gamma_v(X,\cE)^\vee
    \]
    as the composite
        \begin{multline*}
            \PD_X(\cE) \colon \rR \Gamma_v\bigl(X, \DD^{\mathrm{ar}}_X(\cE)\bigr) = \rR\Gamma_v\bigl(X, \rR\cHom(\cE,\omega^{\mathrm{ar}}_X)\bigr) \xlongrightarrow{\Ev_f} \rR\Hom\bigl(\rR\Gamma_v(X, \cE),\rR\Gamma_v(X, \omega^{\mathrm{ar}}_X)\bigr) \\ \xlongrightarrow{\Tr_X \circ \blank} \rR\Hom\bigl(\rR\Gamma_v(X, \cE),\QQ_p\bigr) = \rR\Gamma_v(X, \cE)^\vee.
        \end{multline*}
        The first arrow above is the evaluation map from 
    \cref{Ev-Adj}. 
\end{enumerate}
\end{definition}

Finally, we are ready to prove the general version of the arithmetic duality theorem. 
 
\begin{corollary}[Arithmetic duality: singular varieties]\label{cor:arithmetic-duality-general} Let $X$ be a  proper rigid analytic  space over $K$ of dimension $d$ and let $\cE\in \Perf_v(X;\underline{\Q}_p) \cap \cD^{[r, r']}(X_v; \uQ_p)$ for some $r,r' \in \ZZ\cup \{-\infty, \infty\}$. Then the following statements hold:
\begin{enumerate}
    \item\label{cor:arithmetic-duality-general-1} the complex $\rR\Gamma_v(X, \cE)$ lies in $\cD^{[r, r'+2d+2]}_{\mathrm{coh}}(\QQ_p)$;
    \item\label{cor:arithmetic-duality-general-2} the duality mirphism $\PD_X(\cE) \colon \rR\Gamma_v\bigl(X, \DD^{\mathrm{ar}}_X(\cE)\bigr)
    \longrightarrow
    \rR\Gamma_v(X,\cE)^\vee$
    is an isomorphism. 
\end{enumerate}
\end{corollary}
\begin{proof} 
We proceed by induction on $d \colonequals \dim X \in \ZZ_{\ge 0} \cup \{-\infty\}$.
If $d = -\infty$, we have $X = \varnothing$, so there is nothing to show.
For the induction step, we assume that $d \ge 0$ and that the statement has been proven in dimensions $< d$.
Let $X_\red \subseteq X$ be the maximal reduced closed subspace.
Under the isomorphism $X^\lozenge \simeq X^\lozenge_\red$, the dualizing complexes $\omega_X^{\mathrm{ar}}=\omega^{\QQ_p}_X$ and $\omega_{X_\red}^{\mathrm{ar}}=\omega^{\QQ_p}_{X_\red}$ are identified because the induced equivalence of \'etale sites identifies $\omega^{\ZZ/n}_X$ and $\omega^{\ZZ/n}_{X_\red}$ (cf.\ \cite[Lem.~7.2.4]{LRZ24}).
Therefore, we may assume that $X$ is reduced.

\Cref{dualizing-cohomology}\cref{dualizing-cohomology-modification-square} now provides the following pushout square of diamonds:
    \[ \begin{tikzcd}
       (Z')^\lozenge \arrow[r,"i'"] \arrow[d,"\pi'"'] \arrow[rd,"h"] & (X')^\lozenge \arrow[d,"\pi"] \\
       Z^\lozenge \arrow[r,"i"] & X^\lozenge
    \end{tikzcd} \]
such that $X'$ is smooth proper with $\dim X'\leq d$ and $Z, Z'$ are proper with  $\dim Z, \dim Z'<d$. 

Now we first address Part~\cref{cor:arithmetic-duality-general-1}. Arguing as in the proof of \cref{singular-PD}, we get a cocartesian square
\[
\begin{tikzcd}
    \rR\Gamma_v(X, \cE) \arrow{r} \arrow{d} & \rR\Gamma_v(X', \pi_v^* \cE) \arrow{d}\\
    \rR\Gamma_v(Z, i_v^*\cE) \arrow{r}& \rR\Gamma_v(Z', h^*_v\cE). 
\end{tikzcd}
\]
This leads to an exact triangle 
\[
\rR\Gamma_v(X, \cE) \to \rR\Gamma_v(X', \pi_v^* \cE) \oplus \rR\Gamma_v(Z, i_v^*\cE) \to \rR\Gamma_v(Z', h^*_v\cE).
\]
Using the induction hypothesis applies to $Z$ and $Z'$ and \cref{cor:coh-dimension-arithmetic-case} applied to $X'$, we immediately conclude that $\rR\Gamma_v(X, \cE)$ lies in $\cD^{[r, r'+2d+2]}_{\mathrm{coh}}(\QQ_p)$.

Now we establish the duality claim by a similar reduction. Again, arguing as in the proof of \cref{singular-PD}, we get a commutative diagram (i.e., a functor $\Delta^1\times \Delta^1\times \Delta^1 \to \cD(\QQ_p)$)
  \[ \begin{tikzcd}[row sep=small,column sep=-3em]
        & \rR{\Gamma}_v\bigl(Z^{\prime},h^*_v\cE\bigr)^{\vee} \arrow[rr,orange] \arrow[dd,orange] && \rR{\Gamma}_v\bigl(X^{\prime},\pi^*_v\cE\bigr)^{\vee} \arrow[dd,orange] \\
        \rR{\Gamma}_{v}\bigl(Z^{\prime},\mathbf{D}^{\mathrm{ar}}_{Z^{\prime}}(h^*_v\cE)\bigr) \arrow[ru,purple,"{\mathrm{PD}_{Z'}(h^*_v\cE)}"{near start},"\sim"'{sloped}] \arrow[rr,crossing over,red] \arrow[dd,red] && \rR{\Gamma}_{v}\bigl(X^{\prime},\mathbf{D}^{\mathrm{ar}}_{X^{\prime}}(\pi^*_v\cE)\bigr) \arrow[ru,purple,"{\mathrm{PD}_{X'}(\pi^*_v\cE)}"{near start},"\sim"'{sloped}] & \\ 
        & \rR{\Gamma}_v\bigl(Z,i_v^*\cE\bigr)^{\vee} \arrow[rr,orange] && \rR{\Gamma}_v\bigl(X,\cE\bigr)^{\vee} \\
        \rR{\Gamma}_{v}\bigl(Z,\mathbf{D}^{\mathrm{ar}}_Z(i_v^*\cE)\bigr) \arrow[ru,purple,"{\mathrm{PD}_Z(i_v^*\cE)}"'{near end},"\sim"{sloped}] \arrow[rr,red] && \rR{\Gamma}_{v}\bigl(X,\mathbf{D}^{\mathrm{ar}}_{X}(\cE)\bigr) \arrow[ru,purple,near end,"{\mathrm{PD}_X(\cE)}"'] \arrow[from=uu,crossing over,red] &
    \end{tikzcd} \]
such that the red and orange squares are cocartesian.
Thus, in order to show that $\PD_X(\cE)$ is an isomorphism, it suffices to show that $\PD_{Z}(i^*_v\cE)$, $\PD_{Z'}(h^*_v\cE)$, and $\PD_{X'}(\pi^*_v\cE)$ are. The cases of $\PD_{Z}(i^*_v\cE)$ and $\PD_{Z'}(h^*_v\cE)$ follow from the induction hypothesis, while the case of $\PD_{X'}(\pi^*_v\cE)$ follows from \cref{thm:arithmetic-duality-smooth}. This finishes the proof. 
\end{proof}

\newpage

\addtocontents{toc}{\protect\setcounter{tocdepth}{1}}

\begin{appendices}

\section{Limits in compactly generated \texorpdfstring{$\infty$}{infinity}-categories}\label{section:appendix-limits}
        
In this appendix, we collect certain facts about limits in compactly generated $\infty$-categories. All facts in this appendix are certainly well-known to the experts, but they seem difficult to extract from the existing literature. For this reason, we decided to include proofs of these facts. 
    
\subsection{General facts}

In this subsection, we collect some facts about commutation of certain limits and colimits in compactly generated $\infty$-categories.  

\begin{defn}\label{defn:generated-by-cotruncated-objects} An $\infty$-category $\cC$ is {\it compactly generated by cotruncated objects} if $\cC$ is compactly generated and every compact object of $\cC$ is cotruncated (i.e., truncated in $\cC^\op$ in the sense of \cite[\href{https://kerodon.net/tag/05EC}{Tag 05EC}]{kerodon}). 
\end{defn}

\begin{examples}\label{examples:compactly-generated-by-cotruncated-objects} Fix an integer $n\geq 0$. Then 
\begin{enumerate}
    \item\label{examples:compactly-generated-by-cotruncated-objects-1} the $\infty$-category $\Ani_{\leq n}$ of $n$-truncated anima is compactly generated by the cotruncated object $*$;
    \item\label{examples:compactly-generated-by-cotruncated-objects-2} the $\infty$-category $\Cat_{(n, 1)}$ of $(n,1)$-categories is compactly generated by the cotruncated object $\Delta^1$;
    \item\label{examples:compactly-generated-by-cotruncated-objects-3} for any ring $A$, the $\infty$-category $\cal{D}^{\geq -n}(A)$ is compactly generated by the cotruncated objects $A[-i]$ for $i\geq -n$. 
\end{enumerate}
\end{examples}

\begin{lemma}\label{lemma:finite-limits-vs-filtered-colimits} 
Let $\cC$ be a compactly generated $\infty$-category and let $I$ be a finite simplicial set. Then $\lim_I \colon \Fun(I, \cC) \to \cC$ commutes with filtered colimits. 
If, in addition, $\cC$ is compactly generated by cotruncated objects, the same conclusion holds for $I=\Delta$. 
\end{lemma}
\begin{proof}
    For the first statement: Since $\cC$ is compactly generated, it suffices to show that $\lim_I$ commutes with filtered colimits after postcomposing $\lim_I \colon \Fun(I, \cC) \to \cC$ with $\Map_\cC(c, \blank)\colon \cC \to \Ani$ for every compact object $c\in \cC$. Therefore, it suffices to prove the result when $\cC=\Ani$. In this case, the result follows from \cite[Prop.~5.3.3.3]{HTT}.

    Now suppose that $\cC$ is compactly generated by cotruncated objects and $I=\Delta$. In such situation, it suffices to prove the claim after postcomposing with $\Map_\cC(c, \blank)$ for every compact cotruncated $c\in \cC$. In this case, $\Map_{\cC}(c, \blank)$ factors through $\Ani_{\leq n}$ for some $n$, so it suffices to prove the claim when $\cC=\Ani_{\leq n}$ for some integer $n\geq 0$. Then \cite[Prop.~A.1]{Dirac-1} implies that $\lim_{\Delta} = \lim_{\Delta_{\leq n+1}}$, so the limit becomes finite and follows from the previous case. 
\end{proof}

\subsection{Limits of categories}

In this subsection, we include some results which are specific to limits in $\Cat_\infty$. 

\begin{lemma}\label{lemma:isomorphism-of-cosimplicial-categories} Let $f^\bullet \colon \cC^\bullet \to \cD^\bullet$ be a morphism of cosimplicial objects in $\Cat_\infty$. Suppose that $f^0 \colon \cC^0 \to \cD^0$ is an equivalence and that $f^n\colon \cC^n \to \cD^n$ is fully faithful for any integer $n>0$. Then $\lim_\Delta f^\bullet \colon \lim_\Delta \cC^\bullet \to \lim_\Delta \cD^\bullet$ is an equivalence.
\end{lemma}
\begin{proof}
    For brevity, we denote by $f\colon \cC \to \cD$ the morphism $\lim_\Delta f^\bullet \colon \lim_\Delta \cC^\bullet \to \lim_\Delta \cD^\bullet$. Since limits preserve fully faithful functors, we conclude that $f$ is fully faithful. So we only need to show that it is essentially surjective. In other words, it suffices to prove that the fiber product $\cC \times_{\cD} \Delta^0$ is non-empty for any morphism $\Delta^0 \to \cD$. The morphism $\Delta^0 \to \cD$ defines an essentially unique morphism $\ud{\Delta}^0 \to \cD^\bullet$ where $\ud{\Delta}^0$ is the constant cosimplicial object. Since limits commute with limits, we deduce the question to showing that $\lim_{\Delta} \big(\cC^n \times_{\cD^n} \Delta^0\big)$ is non-empty for any $\ud{\Delta}^0 \to \cD^\bullet$. Our assumption on $f^\bullet$ implies that $\cC^0 \times_{\cD^0} \Delta^0 \simeq \Delta^0$ and that $\cC^n \times_{\cD^n} \Delta^0 \to \Delta^0$ is fully faithful for any $n>0$. In particular, $\cC^n \times_{\cD^n} \Delta^0$ is either empty or isomorphic to $\Delta^0$ for any $n\geq 0$. Now we note that the coface maps define morphisms $\Delta^0\simeq \cC^0 \times_{\cD^0} \Delta^0 \to \cC^n \times_{\cD^n} \Delta^0$ for any $n>0$. In particular, we conclude that $\cC^n \times_{\cD^n} \Delta^0$ cannot be empty, so all these categories are isomorphic to $\Delta^0$. Since $\Delta^0$ is the final object in $\Cat_\infty$, we conclude that all coface and codegeneracy maps are isomorphisms. Therefore, the combination of \cite[\href{https://kerodon.net/tag/02XX}{Tag 02XX}]{kerodon}, \cite[\href{https://kerodon.net/tag/02QP}{Tag 02QP}]{kerodon}, and \cite[\href{https://kerodon.net/tag/02QL}{Tag 02QL}]{kerodon} implies that $\lim_\Delta \big( \cC^n \times_{\cD^n} \Delta^0 \big)$ is isomorphic to $\Delta^0$ as well. In particular, it is non-empty. This finishes the proof. 
\end{proof}

\begin{cor}\label{cor:descent-data} Let $f^\bullet \colon \cC^\bullet \to \cD^\bullet$ be a morphism of cosimplicial objects in $\Cat_\infty$. Suppose that $f^n\colon \cC^n \to \cD^n$ is fully faithful for any integer $n\geq 0$. Then $\lim_\Delta f^\bullet \colon \lim_\Delta \cC^\bullet \to \lim_\Delta \cD^\bullet$ is fully faithful and $d\in \lim_\Delta \cD^\bullet$ lies in the essential image of $\lim_\Delta f^\bullet$ if and only if its projection $d^0\in \cD^0$ lies in the essential image of $f^0 \colon \cC^0 \to \cD^0$.
\end{cor}
\begin{proof}
    It suffices to show that the natural morphism $\lim_\Delta \cC^\bullet  \to \cC^0 \times_{\cD^0} \lim_{\Delta} \cD^\bullet \simeq \lim_{\Delta} \big(\ud{\cC}^0 \times_{\ud{\cD}^0} \cD^\bullet\big)$ is an equivalence. This follows immediately from \cref{lemma:isomorphism-of-cosimplicial-categories}.
\end{proof}

\begin{lemma}\label{lemma:commutation-limits-colimits} Let $F\colon I\to \Fun(\Delta, \Cat_\infty)$ be a filtered diagram of cosimplicial objects in $\Cat_\infty$. Suppose that the morphism $F(i)^n \to F(i')^n$ is fully faithful for any $i\to i'$ in $I$ and every $n\in \Delta$. Then the natural morphism
\[
\colim_{i\in I} \lim_{\Delta} F(i)^\bullet \xr{\alpha} \lim_{\Delta} \colim_{i\in I} F(i)^\bullet
\]
is an equivalence.
\end{lemma}
\begin{proof}
    Since limits and filtered colimits preserve fully faithful morphisms, we conclude that $\alpha$ is fully faithful. Therefore, it suffices to show that it is essentially surjective. According to \cref{cor:descent-data}, an object
    in $\lim_{\Delta} \colim_{i\in I} F(i)^\bullet$ comes from $\lim_{\Delta} F(i_0)^\bullet$ for some $i_0\in I$
    if and only if its image in $\colim_{i \in I} F(i)^0$ comes from $F(i_0)^0$.
    By definition, every object in $\colim_{i \in I} F(i)^0$ must come from $F(i_0)^0$ for some $i_0 \in I$,
    hence we see that $\alpha$ is essentially surjective. 
\end{proof}

\section{\texorpdfstring{$\cV$}{V}-valued sheaves}\label{section:valued-sheaves}

Let $\cC$ be an $\infty$-site in the sense of \cite[Def.~A.3.2]{Mann-thesis} and let $\cV$ be a compactly generated $\infty$-category. Then \cite[\textsection 1.3.1, App.~A.3]{SAG} and \cite[App.~A.3]{Mann-thesis} define the $\infty$-categories of $\cV$-valued sheaves $\Shv(\cC; \cV)$ and $\cV$-valued hypercomplete sheaves $\Shv^{\hyp}(\cC; \cV)$. In this appendix, we collect some standard facts about $\Shv(\cC; \cV)$ which seem to be somewhat difficult to extract from the existing literature. 

For the rest of this appendix, we fix an $\infty$-site $\cC$ and a compactly generated $\infty$-category $\cV$. 

\subsection{General facts}

First, we provide some basic results about sheafification of $\cV$-valued presheaves. Then we show that any $\cV$-valued sheaf is automatically hypercomplete if $\cV$ is compactly generated by cotruncated objects. 

Recall that \cite[Lem.~A.3.7]{Mann-thesis} implies that the inclusion $\Shv(\cC; \cV) \to \PShv(\cC; \cV)$ admits a left adjoint $\rL_{\cV} \colon \PShv(\cC; \cV) \to \Shv(\cC; \cV)$. 

\begin{lemma}\label{lemma:commutes-with-sheafification} Let $f\colon \cV \to \cV'$ be a functor between compactly generated $\infty$-categories, let $f_*\colon \PShv(\cC; \cV) \to \PShv(\cC; \cV')$ be the induced functor, and let $F\in \PShv(\cC; \cV)$. If $f$ commutes with filtered colimits, then $\rL_{\cV'}\big( f_*F\big) \to \rL_{\cV'}\big(f_*(\rL_{\cV} F)\big)$ is an isomorphism. If, in addition, $f$ also commutes with all limits, then $\rL_{\cV'}( f_*F) \to f_*(\rL_{\cV} F)$ is an isomorphism.
\end{lemma}
\begin{proof}
    The first claim follows automatically from \cite[Lem.~1.3.4.4]{SAG}. Now we prove the second claim. Since $\Shv(\cC; \cV) \subset \PShv(\cC;\cV)$ is a reflective subcategory (see the proof of \cite[Lem.~A.3.7]{Mann-thesis}), it suffices to show that $f_*(\rL_{\cV} F)$ is already a sheaf. This follows directly from \cite[Def.~A.3.3]{Mann-thesis} and the fact that $f$ preserves limits.  
\end{proof}

\begin{example}\label{examples:commutes-with-limits-and-filtered-colimits} Fix an integer $n\geq 0$. Then 
\begin{enumerate}
    \item\label{examples:commutes-with-limits-and-filtered-colimits-1} the inclusion $\Ani_{\leq n} \hookrightarrow \Ani$ commutes with all limits and filtered colimits (see \cite[\href{https://kerodon.net/tag/05E2}{Tag 05E2}]{kerodon} and \cite[\href{https://kerodon.net/tag/0631}{Tag 0631}]{kerodon});
    \item\label{examples:commutes-with-limits-and-filtered-colimits-2} the inclusion $\Cat_{(n, 1)} \hookrightarrow \Cat_\infty$ commutes with all limits and filtered colimits (see \cite[\href{https://kerodon.net/tag/05RZ}{Tag 05RZ}]{kerodon} and \cite[\href{https://kerodon.net/tag/062Z}{Tag 062Z}]{kerodon});
    \item\label{examples:commutes-with-limits-and-filtered-colimits-4} the functor $(\blank)^\simeq \colon \Cat_\infty \to \Ani$ commutes with limits and filtered colimits (see \cite[\href{https://kerodon.net/tag/02ST}{Tag 02ST}]{kerodon} and the proof of \cite[\href{https://kerodon.net/tag/062Z}{Tag 062Z}]{kerodon});
    \item\label{examples:commutes-with-limits-and-filtered-colimits-5} the functor $\rm{Fun}(\Delta^1, \blank) \colon \Cat_\infty \to \Cat_\infty$ commutes with limits and filtered colimits. Indeed, it commutes with limits because it is right adjoint to $\Delta^1\times(\blank)$ and it commutes with filtered colimits due to the proof of \cite[\href{https://kerodon.net/tag/062Z}{Tag 062Z}]{kerodon});
    \item\label{examples:commutes-with-limits-and-filtered-colimits-3} the inclusion $\cal{D}^{\geq -n}(A) \to \cD(A)$ commutes with all limits and filtered colimits.
\end{enumerate}
\end{example}

\begin{cor}\label{cor:finite-limits-commute-with-sheafifications} Let $\cV$ be a compactly generated $\infty$-category. Then the sheafification functor $\rL \colon \PShv(\cC; \cV) \to \Shv(\cC; \cV)$ commutes with finite limits.
\end{cor}
\begin{proof}
    Since $\Map_{\cV}(v, \blank)$ commutes with both limits and filtered colimits for a compact object $v\in \cV$, \cref{lemma:commutes-with-sheafification} implies that $\Map_{\cV}(v, \blank)$ commutes with $\rL$ for any such $v\in \cV$. Since $\cV$ is compactly generated by compact objects, we can reduce the question to the case $\cV=\Ani$. In this case, the result follows from \cite[Cor.~6.2.1.6 and Prop.~6.2.2.7]{HTT}. 
\end{proof}

Now we show the promised hypercompletness statement. 

\begin{lemma}\label{lemma:sheaves-are-hypersheaves} Let $\cV$ be an $\infty$-category compactly generated by cotruncated objects and let $F\in \Shv(\cC; \cV)$. Then $F$ is a hypercomplete $\cV$-valued sheaf.
\end{lemma}
\begin{proof}
    Since $\cV$ is compactly generated by co-truncated objects, \cite[Def.~A.3.5]{Mann-thesis} implies that it suffices to show that $F_v\coloneqq \Map_{\cV}(v, \blank) \circ F \colon \cC^\op \to \Ani$ is a hypercomplete sheaf for any cotruncated compact object $v\in \cV$. Since $\Map_{\cV}(v, \blank)$ commutes with limits, \cite[Def.~A.3.3]{Mann-thesis} implies that $F_v$ is a sheaf. Furthermore, we note that $F_v$ takes values in $\Ani_{\leq n}$ for some $n$ since $v$ is cotruncated. Therefore, we can assume that $\cV=\Ani_{\leq n}$. Then the result follows from \cite[Lem.~6.5.2.9]{HTT}.
\end{proof}

Finally, we consider the following example:

\begin{lemma}\label{lemma:derived-sheafification} Let $\cC$ be an explicit covering site such that $\cC$ is an ordinary category and let $F\colon \cC^\op \to \Mod_\Z$ be an (ordinary) sheaf of abelian groups and let $\alpha\colon F(\blank) \to \rR\Gamma(\blank, F)$ be the natural functor of $\cD(\Z)$-valued presheaves on $\cC$. Then $\alpha$ realizes the target as the sheafification of the source. 
\end{lemma}
\begin{proof}
    Denote by $Q$ the cofiber of $\alpha$. Then it suffices to show that $\rL Q \simeq 0$. Furthermore, \cref{examples:commutes-with-limits-and-filtered-colimits}~\cref{examples:commutes-with-limits-and-filtered-colimits-3} and \cref{lemma:commutes-with-sheafification} implies that we can consider all sheaves involved as $\cD^{\geq 0}(\Z)$-valued sheaves. Therefore, \cref{lemma:sheaves-are-hypersheaves} implies that $\rL Q$ is automatically hypercomplete. Using \cite[Prop.~1.3.3.3(3)]{SAG}, one sees that it suffices to show that the (usual) sheafification of the presheaf of abelian groups
    \[
    (U\in \cC) \mapsto \Hh^i\big(Q(U)\big)
    \]
    is zero for any $i\in \ZZ$. By construction, we have $\Hh^i\big(Q(U)\big)=0$ for $i \leq 0$ and $U\in \cC$. To see that the sheafification of $U\mapsto \Hh^i\big(Q(U)\big)$ is equal to zero for $i>0$, it suffices to show that each class $x\in \Hh^i\big(\rR\Gamma(U, F)\big) = \Hh^i(U, F)$ is annihilated after some covering $\{U'_i \to U\}_{i\in I}$. This follows directly from the construction of cohomology groups. 
\end{proof}

Now we recall that there is a functor $\tau^{\leq 0}_{\Ani} \colon \cD(\Z) \to \Ani$ defined as the composition $\cD(\Z) \xr{\tau^{\leq 0}} \cD(\Z)^{\leq 0} \simeq \Ani(\Ab) \to \Ani$, where the second equivalence is \cite[Ex.~5.1.6(2)]{purity}. This functor commutes with limits and filtered colimits. Therefore, \cref{lemma:derived-sheafification} and \cref{lemma:commutes-with-sheafification} immediately imply the following result. 

\begin{variant}\label{variant:sheafification-anima} Let $\cC$ and $F$ be as in \cref{lemma:derived-sheafification}. Then $\alpha \colon \tau^{\leq 0}_{\Ani}\big(F(\blank)[n]\big) \to \tau^{\leq 0}_{\Ani}\big(\rR\Gamma(\blank, F)[n]\big)$ realizes the target as the sheafification of the source for any integer $n$.
\end{variant}

\subsection{Sheaves of $\infty$-categories}

In this subsection, we collect some facts specific to sheaves of $\infty$-categories. Throughout this section, we fix an $\infty$-site $\cC$ such that $\cC$ is an \emph{ordinary} category. 

First, we discuss the objects in the sheafification of a presheaf of $\infty$-categories. We start with presheaves of anima. 

\begin{lemma}\label{lemma:locally-surjective} Let $F\in \PShv(\cC; \Ani)$ and let $\theta\colon F\to \rL F$ be the sheafification morphism. Then $\theta$ is locally surjective, i.e., for every $U\in \cC$ and every section $s\in \pi_0(\rL F(U))$, there is a covering $\{U_i \to U\}_{i\in I}$ and sections $t_i\in \pi_0(F(U_i))$ such that $\theta(t_i) = s|_{U_i}$.
\end{lemma}
\begin{proof}
    It suffices to show that $\pi_0 \circ F \to \pi_0 \circ \rL F$ becomes an isomorphism after sheafification. Now note that $\pi_0\colon \Ani \to \Set$ commutes with colimits since it is a left adjoint. Therefore, the result follows immediately from \cref{lemma:commutes-with-sheafification}. 
\end{proof}

\begin{lemma}\label{lemma:locally-essentially-surjective} Let $F\in \PShv(\cC; \Cat_\infty)$ and let $\theta\colon F\to \rL F$ be the sheafification morphism. Then $\theta$ is locally essentially surjective, i.e., for every $U\in \cC$ and every section $s\in \rm{h}\big(\rL F(U)\big)$, there is a covering $\{U_i \to U\}_{i\in I}$ and sections $t_i\in \rm{h}(F(U_i))$ such that $\theta(t_i) \simeq s|_{U_i}$.
\end{lemma}
\begin{proof}
    \cref{examples:commutes-with-limits-and-filtered-colimits}~\cref{examples:commutes-with-limits-and-filtered-colimits-4} and \cref{lemma:commutes-with-sheafification} imply that $(\rL F)^{\simeq} \simeq L(F^\simeq)$. This immediately reduces the question to showing that $F^\simeq \to \rL(F^\simeq)$ is locally surjective. This follows from \cref{lemma:locally-surjective}. 
\end{proof}

Now we wish to get some control over the mapping spaces in $\rL F$. For this, we introduce the following notation.

\begin{defn} Let $F \in \PShv(\cC; \Cat_\infty)$. The {\it pre-sheaf of maps} $\ud{\Fun}(\Delta^1, F) \colon \cC^{\op} \to \Cat_\infty$ is the composition of $F$ with $\Fun(\Delta^1, \blank) \colon \Cat_\infty \to \Cat_\infty$.

Let $U\in \cC$ and let $s,t\in F(U)$ corresponding to the morphisms $* \xr{s} F|_{\cC_{/U}}$ and $* \xr{t} F|_{\cC_{/U}}$ in $\PShv(\cC_{/U}; \Cat_\infty)$. The {\it mapping space presheaf} $\ud{\Map}_F(s,t)\colon \cC^{\op}_{/U} \to \Ani$ is the fiber product
\[
\ud{\Map}_F(s,t) \simeq * \times_{s, F|_{\cC_{/U}}} \ud{\Fun}(\Delta^1, F|_{\cC_{/U}}) \times_{F|_{\cC_{/U}}, t} *.
\]
We note that \cite[\href{https://kerodon.net/tag/01JC}{Tag 01JC}]{kerodon} 
ensures that $\ud{\Map}_F(s,t)$ indeed takes value in $\Ani$.
\end{defn}

More informally, the presheaf $\ud{\Map}_F(s,t)$ sends $V\in \cC_{/U}$ to $\Map_{F(V)}(s|_V, t|_V)$.

\begin{rmk} Since $\Fun(\Delta^1, \blank)$ commutes with limits, we see that $\ud{\Fun}(\Delta^1, F)$ is a sheaf of $\infty$-categories if $F$ is so. Likewise, we see that $\ud{\Map}_F(s,t)$ is a sheaf of anima on $\cC_{/U}$
if $F|_{\cC_{/U}}$ is a sheaf of $\infty$-categories.
\end{rmk}

\begin{cor}\label{cor:maps-in-sheafification} Let $F\in \PShv(\cC; \Cat_\infty)$, let $U\in \cC$, and let $s, t\in \F(U)$. Then the natural morphism 
\[
\ud{\Map}_F(s, t) \to \ud{\Map}_{\rL F}(s, t)
\]
realizes $\ud{\Map}_{\rL F}(s, t)$ as the sheafification of $\ud{\Map}_F(s, t)$.
\end{cor}
\begin{proof}
    This follows directly from \cref{lemma:commutes-with-sheafification}, \cref{examples:commutes-with-limits-and-filtered-colimits}\cref{examples:commutes-with-limits-and-filtered-colimits-5}, and \cref{cor:finite-limits-commute-with-sheafifications}.
\end{proof}

\begin{defn} A morphism $\alpha\colon F\to G$ in $\PShv(\cC; \Cat_\infty)$ is {\it fully faithful} if $\alpha(U) \colon F(U) \to G(U)$ is fully faithful for any $U\in \cC$.
\end{defn}

\begin{rmk}\label{rmk:iso-locally} \cref{cor:descent-data} implies that if $\alpha\colon F \to G$ is fully faithful and locally essentially surjective morphism of sheaves of $\infty$-categories, then $\alpha$ is an isomorphism.
\end{rmk}

We show that sheafification preserves fully faithful functors. 

\begin{lemma}\label{lemma:sheafification-of-fully-faithful} Let $\alpha\colon F \to G$ be a fully faithful morphism in $\PShv(\cC; \Cat_\infty)$. Then $\rL \alpha \colon \rL F \to \rL G$ is fully faithful as well. 
\end{lemma}
\begin{proof}
    It suffices to show that, for every $U\in \cC$ and any two sections $s,t\in (\rL F)(U)$, the natural morphism
    \[
    \ud{\Map}_{\rL F}(s, t) \to \ud{\Map}_{\rL G}(\rL \alpha(s), \rL \alpha(t))
    \]
    is an isomorphism. Since the source and the target are sheaves on $\cC_{/U}$, it suffices to show locally on $\cC$. Therefore, \cref{lemma:locally-essentially-surjective} implies that we can assume that $s, t$ lie in the essential image of $\theta_F \colon F \to \rL F$. Then \cref{cor:maps-in-sheafification} implies that it suffices to show that
    \[
    c\colon \rL \Big( \ud{\Map}_F\big(s, t\big) \Big) \to \rL \Big( \ud{\Map}_G\big(\alpha(s), \alpha(t)\big)\Big)
    \]
    is an isomorphism. The assumption that $\alpha$ is fully faithful implies that the map $\ud{\Map}_F(s, t) \to \ud{\Map}_G(\alpha(s), \alpha(t))$ is already an isomorphism. In particular, it remains an isomorphism after applying the sheafification functor $\rL$. 
\end{proof}

\begin{cor}\label{cor:criteria-sheafification} Let $F\in \PShv(\cC; \Cat_\infty)$, let $G\in \Shv(\cC; \Cat_\infty)$, and let $\alpha \colon F \to G$ be a morphism such that $\alpha$ is locally essentially surjective and, for each $U\in \cC$ and $s,t\in \pi_0(F(U)^\simeq)$, the natural morphism 
\[
\ud{\Map}_{F}(s, t) \to \ud{\Map}_G(\alpha(s), \alpha(t))
\]
realizes the target as the sheafification of the source. Then $\alpha$ realizes $G$ as the sheafification of $F$.
\end{cor}
\begin{proof}
    Let $\theta \colon F \to \rL F$ be the sheafification. Then its universal property implies that there is essentially unique $\rL \alpha\colon \rL F \to G$ such that $\rL\alpha\circ \theta  \simeq \alpha$. We wish to show that $\rL \alpha$ is an equivalence. 

    \cref{rmk:iso-locally} implies that it suffices to show that $\rL \alpha$ is locally essentially surjective and fully faithful. Clearly, it is locally essentially surjective since $\alpha$ is so. Therefore, it suffices to show that $\rL \alpha$ is fully faithful:
    Namely for any $U\in \cC$ and any two sections $s',t'\in (\rL F)(U)$, the natural morphism
    \begin{equation}\label{eqn:fully-faithful}
    \ud{\Map}_{\rL F}(s', t') \to \ud{\Map}_{G}(\rL \alpha(s'), \rL \alpha(t'))
    \end{equation}
    is an equivalence. This can be checked locally, so \cref{lemma:locally-essentially-surjective} implies that we can assume that $s$ and $t$ come from $s, t\in F(U)$. Then the result follows from \cref{cor:maps-in-sheafification} as both sides of \cref{eqn:fully-faithful} are identified with the sheafification of $\ud{\Map}_F(s,t)$. 
\end{proof}

\begin{cor}\label{cor:fully-faithful-and-locally-essentially-surjective} Let $F\in \PShv(\cC; \Cat_\infty)$, let $G\in \Shv(\cC; \Cat_\infty)$, and let $\alpha \colon F \to G$ be a fully faithful and locally essentially surjective morphism. Then $\alpha$ realizes $G$ as the sheafification of $F$.
\end{cor}

For future reference, we also include the following result. 

\begin{lemma}\label{lemma:colimit-of-sheaves} Let $(\cC, S)$ be an explicit covering site in the sense of \cite[Def.~A.3.12]{Mann-thesis}, let $F\colon I \to \Shv(\cC; \Cat_\infty)$ be a filtered diagram of sheaves of $\infty$-categories. Suppose that $F_i \to F_{i'}$ is fully faithful for any $i\to i'$ in $I$. Then $G\coloneqq \colim_I^p F_i \to \colim_I F_i$ is an isomorphism, where $\colim^p_I$ is the presheaf colimit. Furthermore, if each $F_i$ is a hypercomplete sheaf, then $G$ is a hypercomplete sheaf.
\end{lemma}
\begin{proof}
    It suffices to show that $G$ is a sheaf (resp.~hypersheaf). We note that \cite[Prop.~A.3.3.1]{SAG} (resp. \cite[Prop.~A.3.16]{Mann-thesis}) implies that we only need to show that $G \colon \cC^{\op} \to \Cat_\infty$ preserves finite products and satisfies descent with respect to \v{C}ech covers (resp. hypercovers) $U_\bullet \to X$. The first condition follows from \cref{lemma:finite-limits-vs-filtered-colimits}, while the second condition follows easily from \cref{lemma:commutation-limits-colimits}. 
\end{proof}

\subsection{Perfect complexes}

Throughout this subsection, we fix a classical ($1$-categorical) site $\cC$ with a sheaf of (ordinary) rings $\cA$. The main goal is to study the presheaf 
\[
\Perf_\cC(\blank; \cA) \colon \cC^\op \to \Cat_\infty
\]
which sends $U\in \cC$ to $\Perf(\cC_{/U}; \cA|_U)$ (see \nameref{section:notation}).

\begin{lemma}\label{lemma:perfect-complexes-are-hypercomplete} The presheaf $\Perf_\cC(\blank; \cA)$ is a hypercomplete sheaf of $\infty$-categories.
\end{lemma}
\begin{proof}
    By construction, we have $\Perf_\cC(\blank; \cA)\subset \cD(\blank; \cA)$. First, we note that \cite[Cor.~2.1.2.3]{SAG}, \cite[Rmk.~2.1.0.5]{SAG} applied to $\cX = \Shv^{\rm{hyp}}(\cC; \Ani)$, and \cite[Lem.~A.3.6]{Mann-thesis} imply that $\cD(\blank; \cA)$ is a hypercomplete sheaf. In order to conclude that $\Perf_\cC(\blank; \cA)$ is a hypercomplete sheaf, it suffices to note that we can test that an object $F\in \cD(\cC_{/U}; \cA|_U)$ lies in $\Perf(\cC_{/U}; \cA|_U)$ after any covering $U'\to U$. 
\end{proof}

Now we consider another $\Cat_\infty$-valued presheaf
\[
\Perf\big(\cA(\blank)\big) \colon \cC^\op \to \Cat_\infty
\]
which sends $U\in \cC$ to $\Perf\big(\cA(U)\big)$. It comes with the natural morphism
\[
\widetilde{(\blank)}\colon \Perf\big(\cA(\blank)\big) \to \Perf_\cC(\blank; \cA)
\]
which sends $P\in \Perf(\cA(U))$ to $\widetilde{P}\in \Perf(\cC_{/U}; \cA|_U)$. 

\begin{proposition}\label{prop:perfect-complexes-sheafification} The functor $\widetilde{(\blank)}\colon \Perf\big(\cA(\blank)\big) \to \Perf_\cC(\blank; \cA)$ realizes the target as the sheafification of the source.
\end{proposition}
\begin{proof}
    By definition, we see that $\widetilde{(\blank)}$ is locally essentially surjective. 
    Thus, by \cref{cor:criteria-sheafification}, it suffices to show that, for any $U\in \cC$ and $P, Q \in \Perf(\cA(U))$, the natural morphism
    \[
    \ud{\Map}_{\Perf(\cA(\blank))}(P, Q) \to \ud{\Map}_{\Perf_\cC(\blank; \cA)}(\widetilde{P}, \widetilde{Q}) 
    \]
    realizes the target as the sheafification of the source. A simple inductive argument reduces the question to the situation $P=\cA(U)$ and $Q=\cA(U)[n]$ for some integer $n$. In this case, the result follows immediately from \cref{variant:sheafification-anima}.
\end{proof}

\begin{notation}\label{notation:sheaf-perfect-complexes-finite-amplitude} Fix a ringed site $(\cC, \cA)$, elements $r,s \in \ZZ \cup \{-\infty, \infty\}$ with $r \leq s$. The sheaf 
\[
\Perf^{[r, s]}_\cC(\blank; \cA) \colon \cC^{\op} \to \Cat_\infty
\]
is defined as the sheafification of the presheaf $\Perf^{[r, s]}\big(\cA(\blank)\big)$.
\end{notation}

We note that \cref{prop:perfect-complexes-sheafification} implies that $\Perf_\cC(\blank; \cA) \simeq \Perf_\cC^{[-\infty, \infty]}(\blank; \cA)$.

\section{\texorpdfstring{$v$}{v}-descent for perfect complexes}\label{section:v-sheaves}

The main goal of this section is to show that the presheaves $\Perf\big(\O(\blank)\big) \colon \Perfd^{\aff, \op} \to \Cat_\infty$ and $\Perf\big(\BB_I(\blank)\big) \colon \Perfd^{\aff, \op}_{\QQ_p} \to \Cat_\infty$ are sheaves with respect to the $v$-topology. If one only considers affinoid perfectoid spaces that admit a map to a totally disconnected perfectoid space, this result was shown in \cite{ALB} using the descent results from \cite{Mann-thesis}. Our approach does not use any descent results from  \cite{Mann-thesis} and works for all affinoid perfectoid spaces.

Along the way, we also show that the functors $\Perf_\et(\blank; \O^+/\varpi)$, $\Perf_\et(\blank; \O^+)$, and $\Perf_\et(\blank; \AA_I)$ satisfy $v$-descent as well. These results do not seem to appear in the existing literature.

\subsection{Preliminary results}

We collect some basic facts about $v$-sheaves on affinoid perfectoid spaces. In particular, we show that a finitary \'etale sheaf is automatically a $v$-sheaf (see \cref{defn:finitary-presheaf}). As \'etale descent is usually automatic for most functors we are interested in, this result provides us with a very convenient way to verify $v$-descent for certain \'etale sheaves. For example, we will use it to show that the functor  $\Perf_\et(\blank; \O^+/\varpi)$ is a $v$-sheaf. 

However, we should warn the readers that many interesting functors, e.g.~$\Perf\big(\O(\blank)\big)$, are \emph{not} finitary. Therefore, this result would not be directly applicable. However all $v$-descent results established in this paper will eventually be reduced to the fact that a finitary \'etale sheaf is automatically a $v$-sheaf. 

Throughout this subsection, we fix an adic space $S$ and an $\infty$-category $\cV$ with all small colimits. We note that we do not require $S$ to be an analytic adic space. In particular, we allow $S=\Spa(\ZZ_p, \ZZ_p)$, in which case we have $\Perfd_{/S}^\aff=\Perfd^\aff$. 

\begin{defn}\label{defn:finitary-presheaf} A $\cV$-valued presheaf $F\colon \Perfd^{\aff, \op}_{/S} \to \cV$ is {\it finitary} if the natural map $\colim_{I^\op} F(X_i) \to F(X)$ is an isomorphism for any cofiltered limit $X = \lim_I X_i$ in $\Perfd^\aff_{/S}$.
\end{defn}

We recall that, for an affinoid perfectoid space $X$, the affinoid pro-\'etale site $X^\aff_\proet$ comes equipped with the projection morphism $\lambda_\proet \colon X^\aff_\proet \to X_\et^\aff$. 

\begin{lemma}\label{lemma:pullback-etale-proetale} Let $X$ be an affinoid perfectoid space, let $\cV$ be an $\infty$-category compactly generated by cotruncated objects and let $F \in \Shv(X^\aff_\et; \cV)$ be a $\cV$-valued \'etale sheaf on $X^\aff_\et$. Then for any $U = \lim_I U_i \to X$ affinoid pro-\'etale morphism, the natural morphism
\[
\colim_{I^\op} F(U_i) \to (\lambda_\proet^* F)(U)
\]
is an isomorphism. 
\end{lemma}
\begin{proof}
    The proof is similar to the proof of \cite[Prop.~8.5]{diamonds}; the only difference is that we need to use \cref{lemma:finite-limits-vs-filtered-colimits} instead of the fact that filtered colimits commute with finite limits in $\rm{Sets}$.
\end{proof}

\begin{cor}\label{cor:finitary-etale-implies-proetale} Let $\cV$ be an $\infty$-category compactly generated by cotruncated objects and let $F\colon \Perfd^{\aff, \op}_{/S} \to \cV$ be a finitary $\cV$-valued \'etale sheaf. Then $F$ is a $\cV$-valued pro-\'etale sheaf.
\end{cor}

\begin{proof}
    It suffices to show that, for every $X\in \Perfd^{\aff}$, the presheaf $F|_{X^\aff_\proet}$ is a pro-\'etale sheaf. \cref{lemma:pullback-etale-proetale} implies that $\lambda_{\proet}^* \big(F|_{X^\aff_\et}\big)\simeq F|_{X^\aff_\proet}$. In particular, $F|_{X^\aff_\proet}$ is indeed a pro-\'etale sheaf. 
\end{proof}

\begin{lemma}\label{lemma:finitariness-is-v-local} Let $\cV$ be an $\infty$-category compactly generated by cotruncated objects.
For $\tau\in \{\an, \et, \proet, v\}$, let $F \in \Shv_{\tau}(\Perfd^{\aff}_{/S}; \cV)$ be a $\cV$-valued $\tau$-sheaf. Suppose that $S$ is a perfectoid space and that there is a $\tau$-covering $S' \to S$ such that $F|_{\Perfd_{/S'}^\aff}$ is finitary. Then $F$ is finitary as well.
\end{lemma}
\begin{proof}
    We fix a cofiltered limit $X = \lim_I X_i$ in $\Perfd^{\aff}_{/S}$ and denote by $Y \simeq \lim_I Y_i$ its pullback along the map $S' \to S$. Then the claim follows from the following sequence of isomorphisms:
    \[
        \colim_{I^\op} F(X_i) \simeq \colim_{I^\op} \Big( \lim_{\Delta} F\big(Y_i^{\bullet+1/X_i}\big) \Big) \simeq \lim_{\Delta}\Big( \colim_{I^\op} F\big(Y_i^{\bullet+1/X_i}\big)\Big) \simeq \lim_\Delta\Big( F\big(Y^{\bullet+1/X}\big)\Big) \simeq F(X),
    \]
    where the first and the fourth isomorphisms use that $F$ is $\tau$-sheaf, the second isomorphism uses \cref{lemma:finite-limits-vs-filtered-colimits}, and the third isomorphism uses the assumption that $F|_{\Perfd^{\aff}_{/S'}}$ is finitary. 
\end{proof}

The next three results are the key results of this subsection.

\begin{thm}\label{thm:finitary-etale-sheaves-are-v-sheaves} Let $\cV$ be an $\infty$-category compactly generated by cotruncated objects and let $F\colon  \Perfd^{\aff, \op}_{/S} \to \cV$ be a finitary $\cV$-valued \'etale sheaf. Then $F$ is a $\cV$-valued hypercomplete $v$-sheaf.
\end{thm}

\begin{proof}
    \cref{lemma:sheaves-are-hypersheaves} implies that it suffices to show that $F$ is a $v$-sheaf. Now we note that since $F$ is an \'etale sheaf, we have $F(\sqcup_{i\in I}  X_i)  \simeq \prod F(X_i)$ 
    for a finite set $I$ and $X_i\in \Perfd^{\aff, \op}_{/S}$ indexed by $i\in I$.
    Therefore, \cite[Prop.~A.3.3.1]{SAG} implies that it suffices to show that, for any $v$-covering $X \to Y$, the natural morphism
    \[
    F(Y) \to \lim_{\Delta}\Big(F(X^{\bullet+1/Y}) \Big)
    \]
    is an isomorphism. In other words, we have to show that any $v$-covering $f\colon X \to Y$ is of universal $F$-descent in the sense of \cite[Def.~3.1.1]{Liu-Zheng}. Then \cref{cor:finitary-etale-implies-proetale} implies that $F$ is a pro-\'etale sheaf, i.e., any pro-\'etale surjection $X \to Y$ is of universal $F$-descent. Therefore, \cite[Lem.~3.1.2]{Liu-Zheng} guarantees that we can check that a $v$-covering $X\to Y$ is of universal $F$-descent pro-\'etale locally on $Y$. So, by \cite[Lem.~7.18]{diamonds}, we are reduced to the case when $Y$ is a strictly totally disconnected perfectoid space. In this situation, \cite[Lem.~6.4.4]{zav-almost} implies that we can write $f$ as a cofiltered limit $X=\lim_I X_i \xr{f_i} Y$ such that each $f_i$ admits a section. Therefore, \cite[Lem.~3.1.2(1)]{Liu-Zheng} implies that the natural morphism
    \[
        F(Y) \to \lim_{\Delta}\Big(F(X_i^{\bullet+1/Y}) \Big)
    \]
    is an isomorphism for every $i\in I$. Using that $F$ is finitary and that totalizations commute with filtered colimits in $\cV$ (see \cref{lemma:finite-limits-vs-filtered-colimits}), we conclude that the natural morphism $F(Y) \to \lim_{\Delta}\Big(F(X^{\bullet+1/Y}) \Big)$ is an isomorphism as well.
\end{proof}

A similar argument allows us to show the following result:

\begin{lemma}\label{lemma:sections-etale-locally} Let $S$ be an affinoid perfectoid space and let $F \colon \Perfd^{\aff, \op}_{/S} \to \Cat_\infty$ be a finitary presheaf of $\infty$-categories. Suppose that there is a $v$-covering $T \to S$ such that $F(T)\neq \varnothing$. Then there is an \'etale covering $\widetilde{T} \to S$ such that $F(\widetilde{T})\neq \varnothing$.
\end{lemma}
\begin{proof}
    First, \cite[Lem.~7.18]{diamonds} implies that there is an affinoid pro-\'etale covering $S' \to S$ such that $S'$ is strictly totally disconnected. Using that $F$ is finitary, it suffices to show that $F(S')\neq \varnothing$. We set $T'\coloneqq T\times_S S'$. 

    Then \cite[Lem.~6.4.4]{zav-almost} implies that we can write $T' \to S'$ as a cofiltered limit $\lim_I S'_i \to S'$ such that each $S'_i \to S'$
    admits a section. Now we use that $F$ is finitary and that $F(T')\neq \varnothing$ to conclude that $F(S'_i)\neq \varnothing$ for some $i\in I$. Using that $S'_i \to S'$ admits a section, we conclude that $F(S')\neq \varnothing$ as well.
\end{proof}

\begin{thm}\label{thm:sheafification-is-finitary} Let $\cV$ be an $\infty$-category compactly generated by cotruncated objects and let $F\colon \Perfd^{\aff, \op}_{/S} \to \cV$ be a finitary $\cV$-valued presheaf. Then its \'etale sheafification $\rL_\et F$ is a finitary $\cV$-valued \'etale sheaf. 
\end{thm}

\begin{proof}
    We let $X_\infty = \lim_I X_i$ be a cofiltered limit in $\Perfd^\aff_{/S}$. For any $\cV$-valued presheaf $G$, we note by  $G_i$ the restriction of $G$ to $X_{i, \et}^\aff$ and by $G_\infty$ the restriction to $X_{\infty, \et}^\aff$. We note \cite[Cor.~7.2]{etale-hyperdescent} (see also \cite[Prop.~7.6]{etale-hyperdescent}) implies that $(\rL_\et F)_i = \rL_\et ( F_i)$ for any $i\in I$ or $i=\infty$. So we wish to show that the natural morphism
    \[
    \colim_{I^\op} \big( \rL_\et F_i(X_i)\big) \to \rL_\et F_\infty(X_\infty)
    \]
    is an isomorphism. For this, we denote by $u_i\colon X_\infty \to X_i$ the natural projection, and by $u_i^{*, p} \colon \PShv(X^\aff_{i, \et}; \cV) \to \PShv(X^\aff_{\infty, \et}; \cV)$ and $u_i^{*} \colon \Shv(X^\aff_{i, \et}; \cV) \to \Shv(X^\aff_{\infty, \et}; \cV)$ the corresponding pullback morphisms. Then \cite[Prop.~6.4]{diamonds} implies that $\colim^p_{I^\op} u_i^{*, p} F_i \to F_\infty$ is an isomorphism of presheaves. By passing to sheafifications, we conclude that the natural morphism
    \[
    \rL_\et\big(\colim_{I^\op}{}^p u_i^{*, p} \rL_\et F_i\big) \simeq \rL_\et\big(\colim_{I^\op}{}^p u_i^{*, p} F_i\big) \to 
    \rL_\et F_\infty 
    \]
    is an isomorphism, where $\colim^p$ stands for the presheaf colimit. 

    Now set $Y_\infty \to X_\infty$ to be an \'etale morphism of affinoid perfectoid spaces. Then \cite[Prop.~6.4]{diamonds} implies that there is $i_0\in I$ and an \'etale morphism $Y_{i_0} \to X_{i_0}$ such that $Y_{i_0}\times_{X_{i_0}} X_\infty \simeq Y_\infty$. For each $i\in I_{/i_0}$, we also set $Y_i \coloneqq Y_{i_0} \times_{X_{i_0}} X_{i}$. Without loss of generality, we can replace $I$ with $I_{/i_0}$. Then another application of \cite[Prop.~6.4]{diamonds} implies that
    \begin{equation}\label{eqn:colimit}
    (\colim_{I^\op}{}^p u_i^{*, p}\rL_\et F_i\big)(Y_\infty) \simeq \colim_{I^\op} \big(\rL_\et F_i(Y_i) \big).
    \end{equation}
    Therefore, in order to finish the proof, it suffices to show that the natural morphism 
    \[
        \colim_{I^\op}{}^p u_i^{*, p} \rL_\et F_i \to \rL_\et\big(\colim_{I^\op}{}^p u_i^{*, p} \rL_\et F_i\big)
    \]
    is an isomorphism, i.e., $G\coloneqq \colim^p_{I^\op} u_i^{*, p} \rL_\et F_i$ is already an \'etale sheaf. Now \cite[Prop.~A.3.3.1]{SAG} implies that it suffices to show that $G(\sqcup_{i=1}^n Y_i) \simeq \prod_{i=1}^n G(Y_i)$ and $G$ satisfies descent with respect to \'etale covers $Z_\infty \to Y_\infty$ in $X_{\infty, \et}^\aff$.

    The fact that $G$ sends finite disjoint unions into finite products follows from formula~\cref{eqn:colimit} and \cref{lemma:finite-limits-vs-filtered-colimits}. Therefore, it suffices to check that $G$ satisfies descent with respect to an \'etale covering $Z_\infty \to Y_\infty$ in $X_{\infty, \et}^\aff$, i.e., the natural morphism
    \[
        G(Y_\infty) \to \rm{Tot}_{\Delta}\Big(G(Z_\infty^{\bullet+1/Y_{\infty}}) \Big)
    \]
    is an isomorphism. As above, \cite[Prop.~6.4]{diamonds} implies that we can spread out $Z_\infty \to Y_\infty$ to an \'etale covering $Z_{i_0} \to Y_{i_0}$ in $X_{i_0, \et}^\aff$ and we denote by $Z_i \to Y_i$ its pullback to $X_i$ for any $i\in I_{/i_0}$. Without loss of generality, we can replace $I$ with $I_{/i_0}$. Then Formula~\cref{eqn:colimit} implies that it suffices to show that the morphism
    \[
    \colim_{i\in I^{\op}} \rL_\et F_i (Y_i) \to \lim_{\Delta}\Big( \colim_{i\in I^{\op}} \rL_\et F_i (Z_i^{\bullet+1/Y_{i}}) \Big)
    \]
    is an isomorphism.
    This follows from the fact that each $\rL_\et F_i$ is a sheaf and \cref{lemma:finite-limits-vs-filtered-colimits}. 
\end{proof}

We finish the subsection by recording the following very useful fact.

\begin{lemma}
\label{lemma:colimit-of-hypersheaves} 
Let $F\colon I \to \PShv(\Perfd^{\aff}_{/S}; \Cat_\infty)$ be a filtered diagram of finitary presheaves of
$\infty$-categories and let $G$ be the presheaf colimit $\colim^p_I F_i$. 
Suppose that $F_i \to F_{i'}$ is fully faithful for any $i\to i'$ in $I$ and that each $F_i$ factors through
$\Cat_{(d(i), 1)}$ for some integer $d(i)\geq 0$. 
Then $\rL_\et G$ is a finitary hypercomplete $v$-sheaf of $\infty$-categories.
\end{lemma}
We warn the reader that the category $\Cat_\infty$ is \emph{not} compactly generated by cotruncated objects. 
\begin{proof}
     First, \cref{thm:finitary-etale-sheaves-are-v-sheaves}, \cref{thm:sheafification-is-finitary}, and \cref{examples:compactly-generated-by-cotruncated-objects}~\cref{examples:compactly-generated-by-cotruncated-objects-2} imply that $\rL_\et F_i$ is a finitary hypercomplete $v$-sheaf for each $i\in I$. Furthermore, \cref{lemma:sheafification-of-fully-faithful} implies that $\rL_\et F_i \to \rL_\et F_{i'}$ is fully faithful for any $i\to i'$ in $I$. Thus, the result follows from \cref{lemma:colimit-of-sheaves} and the fact that sheafification commutes with colimits. 
\end{proof}

\subsection{$v$-descent for perfect \'etale $\cal{A}$-complexes}

Throughout this section, we fix an adic space $S$ and a \emph{finitary \'etale sheaf} of rings $\cA \colon \Perfd_{/S}^{\aff, \op} \to \Ring$.  Let
\[
\rm{Perf}_\et(\blank; \cA) \colon \Perfd^{\aff, \op}_{/S} \to \Cat_\infty
\]
be the presheaf of $\infty$-categories that sends an affinoid perfectoid $S$-space $X$ to the $\infty$-category $\Perf(X_\et; \cA)$ of perfect complexes of $\cA$-modules on the \'etale site. The main goal of this section is to show that this is a hypercomplete $v$-sheaf of $\infty$-categories.

We start with the following preliminary lemma. 

\begin{lemma}\label{lemma:perfect-complexes-are-finitary} Let $A = \colim_I A_i$ be a filtered colimit of rings $A_i$. Then the natural functor
\[
\colim_{I} \rm{Perf}^{[r, s]}(A_i) \to \rm{Perf}^{[r, s]}(A)
\]
is an equivalence for any $r, s\in \Z \cup \{-\infty, \infty\}$.
\end{lemma}

We remind the reader that $\Perf^{[r, s]}(A)$ denotes the $\infty$-category of perfect $A$-complexes that have tor-amplitude in $[r, s]$. 
\begin{proof}
    This follows directly from the proof of \cite[\href{https://stacks.math.columbia.edu/tag/0BC7}{Tag 0BC7}]{stacks-project}. 
\end{proof}

\begin{cor}\label{cor:mod-perf-etale-finitary} Let $S$ be an adic space, let $\cA \colon \Perfd_{/S}^{\aff, \op} \to \Ring$ be a finitary \'etale sheaf, let $X = \lim_I X_i$ be a cofiltered limit in $\Perfd^\aff_{/S}$. Then the natural functor
\[
\colim_{I^\op} \rm{Perf}^{[r, s]}\Big( \cA(X_i) \Big) \to \rm{Perf}^{[r, s]}\Big( \cA(X) \Big)
\]
is an equivalence for any $r, s\in \Z \cup \{-\infty, \infty\}$.
\end{cor}
\begin{proof}
    This follows directly from the combination of \cref{lemma:perfect-complexes-are-finitary} and the assumption that $\cA$ is finitary.
\end{proof}

Now we introduce the notation that will be useful throughout the whole appendix. 

\begin{notation}\label{notation:sheaf-perfect-complexes-finite-amplitude-cA} Fix an adic space $S$, elements $r,s \in \ZZ \cup \{-\infty, \infty\}$ with $r \leq s$, a topology $\tau\in \{\an, \et, \proet, v\}$, and a $\tau$-sheaf of rings $\cA\colon \Perfd^{\aff, \op}_{/S} \to \Ring$. The $\tau$-sheaf 
\[
\Perf^{[r, s]}_\tau(\blank; \cA) \colon \Perfd^{\aff, \op}_{/S} \to \Cat_\infty
\]
is defined as the $\tau$-sheafification of the presheaf $\Perf^{[r, s]}\big(\cA(\blank)\big)$.
\end{notation}

We note that \cref{prop:perfect-complexes-sheafification} implies that $\Perf_\tau(\blank; \cA) \simeq \Perf_\tau^{[-\infty, \infty]}(\blank; \cA)$. 

\begin{thm}\label{thm:v-descent-finitary-etale} Let $S$ be an adic space, let $\cA \colon \Perfd_{/S}^{\aff, \op} \to \Ring$ be a finitary \'etale sheaf. Then the functor $\rm{Perf}^{[r, s]}_\et(\blank; \cA) \colon \Perfd^{\aff, \op}_{/S} \to \Cat_\infty$ is a finitary hypercomplete $v$-sheaf for any $r, s\in \Z \cup \{-\infty, \infty\}$ with $r \leq s$. In particular, the functor $\Perf_\et(\blank; \cA) \colon \Perfd^{\aff, \op}_{/S} \to \Cat_\infty$ is a finitary hypercomplete $v$-sheaf. 
\end{thm}
\begin{proof}
    We first consider the case when $r$ and $s$ are finite integers. Then we note that $\Perf^{[r, s]}\big(\cA(\blank)\big)$ takes values in $\Cat_{(s-r+1,1)}$. Thus, \cref{cor:mod-perf-etale-finitary} implies that $\Perf^{[r, s]}\big(\cA(\blank)\big)$ is a finitary presheaf of $(s-r+1, 1)$-categories. Therefore, \cref{thm:finitary-etale-sheaves-are-v-sheaves}, \cref{thm:sheafification-is-finitary}, and \cref{examples:compactly-generated-by-cotruncated-objects}~\cref{examples:compactly-generated-by-cotruncated-objects-2}
    imply that $\Perf^{[r, s]}_\et(\blank; \cA) = \rL_\et \Perf^{[r, s]}\big(\cA(\blank)\big)$ is a finitary hypercomplete $v$-sheaf. Now we deal with the case $r\in \Z$ and $s=\infty$. In this case, we note that $\Perf^{[r, \infty]}\big(\cA(X)\big) \simeq \colim_{s\geq r} \Perf^{[r, s]}\big(\cA(X)\big)$. Therefore, the result follows immediately from \cref{lemma:colimit-of-hypersheaves}. The cases $r=-\infty$, $s\in \Z$ and $r=-\infty$, $s=\infty$ are done analogously. 
\end{proof}

Now we note that \cref{thm:finitary-etale-sheaves-are-v-sheaves} ensures that the sheaf $\cA$ is automatically a $v$-sheaf. Thus, it makes sense to consider the category $\Perf_\tau(X; \cA)$ for $\tau\in \{\proet, v\}$.

\begin{cor}\label{cor:pullback-etale-finitary} Let $X$ be a perfectoid space, let $\cA \colon \Perfd_{/X}^{\aff, \op} \to \Ring$ be a finitary \'etale sheaf, let $\tau\in \{\proet, v\}$, and let $\lambda_\tau\colon (X_\tau, \cA) \to (X_\et, \cA)$ be the natural morphism of ringed sites. Then the following statements hold: 
\begin{enumerate}[leftmargin=*,label={\upshape{(\roman*)}}]
    \item\label{cor:pullback-etale-finitary-1} the functor $\lambda_\tau^*\colon \Perf_\et(X; \cA) \to \Perf_\tau(X; \cA)$ is an equivalence of $\infty$-categories;
    \item\label{cor:pullback-etale-finitary-2} the functor $\rR\lambda_{\tau, *}\colon \cal{D}(X_\tau; \cA) \to \cal{D}(X_\et; \cA)$ restricts to a functor $\rR\lambda_{\tau, *} \colon \Perf_\tau(X; \cA) \to \Perf_\et(X; \cA)$;
    \item\label{cor:pullback-etale-finitary-3} $\rR \lambda_{\tau, *}$ is the inverse functor to $\lambda^*_\tau$.
\end{enumerate}
\end{cor}
\begin{proof}
    All claims are \'etale local on $X$, so we can assume that $X$ is an affinoid perfectoid. Then we note that $\lambda^*_\tau \colon \Perf_\et(-; \cA) \to \Perf_\tau(-; \cA)$ realizes the target as the $\tau$-sheafification of the source\footnote{This can be either seen from \cref{cor:criteria-sheafification} or \cref{prop:perfect-complexes-sheafification}.}. Therefore, \cref{thm:v-descent-finitary-etale} implies that it must be an equivalence, so this proves the first claim. 

    To see the second and third claims, it suffices to show that $\cal{P}\to \rR \lambda_{\tau, *}\lambda_\tau^* \cal{P}$ is an equivalence for any $\cal{P}\in \Perf_\et(X; \cA)$. After unravelling the definitions, it boils down to showing that the natural morphism $c_{\cP} \colon \rR\Gamma_\et(X, \cal{P}) \to \rR\Gamma_\tau(X, \lambda_\tau^*\cal{P})$ is an equivalence for any affinoid perfectoid $X$. Since $\rR\Gamma_\et(X, \cP) \simeq \rR\Hom_{\cal{D}(X_\et; \cA)}(\cA, \cP)$ and $\rR\Gamma_\tau(X, \lambda^*_\tau\cP) \simeq \rR\Hom_{\cal{D}(X_\tau; \cA)}(\cA, \lambda^*_\tau\cP)$, we conclude that $c_{\cP}$ is an isomorphism for any $\cP\in \Perf_\et(X; \cA)$.
\end{proof}

\begin{cor}\label{cor:etale-perfect-complex-on-strictly-totally-disconnected} Let $X=\Spa(R, R^+)$ be a strictly totally disconnected perfectoid space, let $\cA \colon \Perfd_{/X}^{\aff, \op} \to \Ring$ be a finitary \'etale sheaf of rings, and let $\tau\in \{\et, \proet, v\}$. Then the functor $\widetilde{(\blank)} \colon \Perf\big(\cA(X)\big) \to \Perf_\tau(X; \cA)$ is an equivalence. Its inverse is given by 
\[
\rR\Gamma_\tau(X, \blank) \colon \Perf_\tau(X; \cA) \to \Perf\big(\cA(X)\big).
\]
\end{cor}
\begin{proof}
    \cref{cor:pullback-etale-finitary}\cref{cor:pullback-etale-finitary-1} implies that $\Perf_\et(X; \cA) \xr{\lambda_\tau^*} \Perf_\tau(X; \cA)$ is an equivalence. Therefore, it suffices to prove the claim when $\tau=\et$. In this case, we note that $\Perf_\et(X; \cA)$ is the \'etale sheafification of the presheaf $\Perf\big(\cA(\blank)\big)$. Since $X$ is strictly totally disconnected, every \'etale covering of $X$ splits. Thus, we conclude that the \'etale sheafification does not change the value on $X$. In particular, we see that $\widetilde{(\blank)} \colon \Perf\big(\cA(X)\big) \to \Perf_\et(X; \cA)$ is an equivalence. To see that its inverse is given by $\rR\Gamma_\et(X, \blank)$, it suffices to show that $P \to \rR\Gamma_\et(X, \widetilde{P})$ is an equivalence for any $P\in \Perf\big(\cA(X)\big)$. By induction, we reduce to the case when $P$ is a direct summand of $\cA(X)^{\oplus n}$ for some integer $n$. This can be further reduced to the case $P=\cA(X)$. In this situation, the question boils down to showing that $\cA(X) \to \rR\Gamma_\et(X, \cA)$ is an equivalence. This follows from the fact that $X$ is strictly totally disconnected, so all higher \'etale cohomology vanish. 
\end{proof}

We also include the following result that is used in the main body of this paper (see \cref{lemma:comes-from-etale}). 

\begin{lemma}\label{lemma:finitary-everything-is-etale} Let $X$ be a perfectoid space over $S$, let $\cA \colon \Perfd_{/X}^{\aff, \op} \to \Ring$ be a finitary \'etale sheaf, let $\lambda_\proet\colon X_\proet \to X_\et$ be the natural projection, and let $\cA_\et$ and $\cA_\proet$ be the restrictions of $\cA$ to $X_\et$ and $X_\proet$ respectively. Then the natural morphism $\lambda_\proet^{-1}(\cA_\et)\to \cA_\proet$ is an isomorphism.
\end{lemma}
We note that \cref{thm:finitary-etale-sheaves-are-v-sheaves} ensures that $\cA$ is automatically a pro-\'etale sheaf (even a $v$-sheaf), so $\cA_\proet$ is well-defined. We also emphasize that the functor $\lambda_\proet^{-1}$ is the pullback functor for sheaves of \emph{abelian groups}, which is usually different from the module-theoretic pullback functor $\lambda^*_\proet$ associated to the morphism of ringed sites $\lambda_\proet \colon (X_\proet; \cA_\proet) \to (X_\et; \cA_\et)$. However, \cref{lemma:finitary-everything-is-etale} \emph{implies} that these two functors are equivalent. 
\begin{proof}
    The question is local on $X$, so we can assume that $X$ is an affinoid. 
    Then \cref{lemma:pullback-etale-proetale} implies that it suffices to show that, for every affinoid pro-\'etale morphism $U=\lim_I U_i \to X$ in $X_\proet$, the natural morphism
    \[
    \colim_{I^\op} \cA(U_i) \to \cA(U)
    \]
    is an isomorphism. This follows immediately from the assumption that $\cA$ is finitary. 
\end{proof}

Now we provide some examples of finitary \'etale sheaves of rings. 

\begin{lemma}\label{lemma:example-z-mod-n} Let $n>0$ be a positive integer. Then the constant \'etale sheaf $\ud{\ZZ/n} \colon \Perfd^{\aff, \op} \to \Ring$ is finitary.
\end{lemma}
\begin{proof}
    \cref{thm:sheafification-is-finitary} implies that it suffices to show that the constant \emph{presheaf} $\ud{\Z/n}_{\rm{presheaf}}$ is finitary. This is evident since $\ud{\Z/n}_{\rm{presheaf}}(X) =\ZZ/n$ for any $X\in \Perfd^\aff$.
\end{proof}

\begin{rmk} \cref{lemma:example-z-mod-n} and \cref{cor:pullback-etale-finitary} together imply that the natural functor 
\[
\lambda_\tau^*\colon \Perf_\et(X; \Z/n) \to \Perf_\tau(X; \Z/n)
\]
is an equivalence for any perfectoid space $X$, $\tau\in \{\proet, v\}$, and a positive integer $n>0$.
\end{rmk}

\begin{lemma}\label{lemma:example-o-plus-mod-varpi} Let $S$ be an affinoid perfectoid space and let $\varpi\in \O_S^+(S)$ be a pseudo-uniformizer. Then the \'etale sheaf $\O^+_\et/\varpi \colon \Perfd^{\aff,\op}_{/S}\to \Ring$ is finitary. 
\end{lemma}
\begin{proof}
    First, we note that $\O^+_\et/\varpi$ is the \'etale sheafification of the presheaf of rings $G$ defined by the rule $X\mapsto \O_{X_\et}^+(X)/\varpi$. Therefore, \cref{thm:sheafification-is-finitary} ensures that it suffices to show that $G$ is finitary. For this, we fix $\Spa(R_\infty, R_\infty^+) = \lim_I \Spa(R_i, R_i^+)$, then we have $R_\infty^+ \simeq \big(\colim_{I^\op} R_i^+\big)^{\wedge}_{(\varpi)}$ (see \cite[Prop.~6.4]{diamonds}). So the question boils down to showing that $\colim_{I^\op} R_i^+/\varpi \to \big(\colim_{I^\op} R_i^+\big)^{\wedge}_{(\varpi)}/\varpi$ is an isomorphism; this follows directly from \cite[\href{https://stacks.math.columbia.edu/tag/05GG}{Tag 05GG}]{stacks-project}. 
\end{proof}

\begin{rmk} \cref{lemma:example-o-plus-mod-varpi} and \cref{thm:finitary-etale-sheaves-are-v-sheaves} imply that, in the notation of \cref{lemma:example-o-plus-mod-varpi}, the \'etale sheaf $\cO^+_\et/\varpi$ is already a $v$-sheaf. In particular, it implies that the quotient $\cO^+_\et/\varpi \colon \Perfd^{\aff, \op}_{/S} \to \Ring$ in the (big) \'etale topology coincides with the quotient $\cO^+_v/\varpi \colon \Perfd^{\aff, \op}_{/S} \to \Ring$ in the $v$-topology.
\end{rmk}

\begin{rmk} \cref{lemma:example-o-plus-mod-varpi} and \cref{cor:pullback-etale-finitary} together imply that the natural functor 
\[
\lambda_\tau^*\colon \Perf_\et(X; \O^+/\varpi) \to \Perf_\tau(X; \O^+/\varpi)
\]
is an equivalence for any perfectoid space $S$, a pseudo-uniformizer $\varpi\in \O_S^+(S)$, an $S$-perfectoid space $X$, and $\tau\in \{\proet, v\}$.
\end{rmk}

\begin{rmk} \cref{lemma:example-o-plus-mod-varpi} and \cref{thm:v-descent-finitary-etale} with $r=s=0$ recover the fact that $\Vect_\et(\blank; \O^+/\varpi)$ is a finitary $v$-stack of (ordinary) categories (see \cite[Prop.~2.6]{Heuer-G-torsor} and \cite[Th.~1.3.1]{zav-almost}). 
\end{rmk}

\begin{lemma}\label{lemma:example-AI-mod-p} Let $I=[a, b] \subset (0, \infty)$ be a closed interval with rational endpoints. Then the \'etale sheaf $\AA_{I, \et}/p \colon \Perfd^{\aff,\op}_{/\Q_p}\to \Ring$ is finitary. 
\end{lemma}
We refer to \cref{cor:AI-etale-sheaf} for the definition of the \'etale sheaf $\AA_I$. 
\begin{proof}
    First, it suffices to show that the restriction $\AA_{I, \et}/p\colon \Perfd^{\aff, \op}_{/S} \to \Ring$ is finitary for any affinoid perfectoid space $S$ over $\Spa(\QQ_p, \ZZ_p)$. Fix $p^\flat\in \O_{S}^{\flat, +}(S)$ as in \cref{notation:p-flat}. Now \cref{lemma:finitariness-is-v-local} implies that we can prove the claim \'etale locally on $S$. Therefore, \cref{rmk:often-well-adapted-pairs} ensures that we can assume that there is a pair $\alpha, \beta\in \O_S^{\flat, +}(S)$ which is well-adapted to $p^\flat$ and $I$. Then \cref{prop:AI-properties} implies that 
    \[
    \AA_{I, \et}/p \simeq \left(\bigoplus_{\mathbf{N}} \cO^{\flat,+}_{\et}/\alpha \right) \oplus \left(\bigoplus_{\mathbf{N}} \cO^{\flat,+}_{\et}/\beta \right).
    \]
    Now we note that (by passing to tilts) \cref{lemma:example-o-plus-mod-varpi} implies that each $\cO^{\flat,+}_{\et}/\alpha$ and $\cO^{\flat,+}_{\et}/\beta$ is finitary. Therefore, their direct sum is still finitary since the topos $\Perfd^{\aff, \op}_{/\Q_p}$ is coherent (see \cite[\href{https://stacks.math.columbia.edu/tag/0739}{Tag 0739}]{stacks-project})
\end{proof}

\begin{rmk} \cref{lemma:example-AI-mod-p} and \cref{thm:finitary-etale-sheaves-are-v-sheaves} imply that, in the notation of \cref{lemma:example-AI-mod-p}, the \'etale sheaf $\AA_{I, \et}/p$ is already a $v$-sheaf. In particular, it implies that the quotient $\AA_{I, \et}/p \colon \Perfd^{\aff, \op}_{/\QQ_p} \to \Ring$ in the (big) \'etale topology coincides with the quotient $\AA_{I, v}/p \colon \Perfd^{\aff, \op}_{/\QQ_p} \to \Ring$ in the $v$-topology.
\end{rmk}

\begin{rmk} \cref{lemma:example-AI-mod-p} and \cref{cor:pullback-etale-finitary} together imply that the natural functor 
\[
\lambda_\tau^*\colon \Perf_\et(X; \AA_I/p) \to \Perf_\tau(X; \AA_I/p)
\]
is an equivalence for any perfectoid space $X$ over $\Spa(\QQ_p, \ZZ_p)$ and $\tau\in \{\proet, v\}$.
\end{rmk}

\subsection{$v$-descent for perfect $\mathcal{O}^+$-complexes}\label{section:descend-O-plus}

Let
\[
\rm{Perf}_\et(\blank; \O^+) \colon \Perfd^{\aff, \op} \to \Cat_\infty
\]
be the presheaf of $\infty$-categories that sends an affinoid perfectoid space $X$ to the $\infty$-category of perfect complexes of $\O_{X_\et}^+$-modules. The main goal of this subsection is to show that this is a hypercomplete sheaf of $\infty$-categories. 

Our strategy for proving this fact will rely on showing that the category of perfect complexes of $\O_{X_\et}^+$-modules
is equivalent to the category of perfect complexes of $\O_{X_v}^+$-modules. 

We start the subsection by showing that the pullback functor $\lambda_v^*\colon \Perf_\et(X; \O^+) \to \Perf_v(X; \O^+)$ is fully faithful for any perfectoid space $X$. For this, we need the following preliminary lemma: 

\begin{lemma}\label{lemma:derived-completeness} Let $S=\Spa(R, R^+)$ be an affinoid perfectoid space, let $\varpi\in R^+$ be a pseudo-uniformizer, let $X$ be a perfectoid space over $S$, let $\tau\in \{\an, \et, \proet, v\}$, and let $\cal{P} \in \Perf_\tau(X; \O^+)$. Then $\cal{P}$ is derived $\varpi$-adically complete. 
\end{lemma}
\begin{proof}
    The claim is $\tau$-local on $X$, so we can assume that $\cal{P}$ can be represented by a finite complex of finite free $\O_{X_\tau}^+$-modules. Then a standard inductive argument reduces the general question to the case $\cal{P}=\O_{X_\tau}^+$. 
    In this case, according to \cite[Lem.~A.8]{zav-almost} and \cite[Prop.~3.4.4]{proetale},
    it suffices to show that $\Hh^i_\tau(X, \O^+)$ is $\varpi$-adically complete for all $i\geq 0$. This follows from \cite[Th.~3.24]{diamonds} for $\tau=\an$, \cite[Th.~6.3]{diamonds} for $\tau=\et$, \cite[Prop.~8.5]{diamonds} for $\tau=\proet$, and \cite[Prop.~8.8]{diamonds} for $\tau=v$.
\end{proof}

\begin{proposition}\label{prop:fully-faithful-integrally} Let $X$ be a perfectoid space, let $\cal{P}, \cal{Q}\in \Perf_\et(X; \O^+)$, let $\tau\in \{\proet, v\}$, and let $\lambda_\tau\colon (X_\tau, \cO^+) \to (X_\et, \cO^+)$ be the natural morphism of ringed sites. Then the following statements hold:
\begin{enumerate}[leftmargin=*,label={\upshape{(\roman*)}}]
    \item\label{prop:fully-faithful-integrally-1} the natural morphism $\cal{P} \to \rR\lambda_{\tau, *} \lambda^*_\tau\cal{P}$ is an isomorphism;
    \item\label{prop:fully-faithful-integrally-2} the natural morphism $\rR\Hom(\cal{P}, \cal{Q}) \to \rR\Hom(\lambda_\tau^*\cal{P}, \lambda_\tau^*\cal{Q})$ is an isomorphism;
    \item\label{prop:fully-faithful-integrally-3} the functor $\lambda_\tau^*\colon \Perf_\et(X; \O^+) \to \Perf_\tau(X; \O^+)$ is fully faithful.
\end{enumerate}
\end{proposition}
\begin{proof}
    Part~\cref{prop:fully-faithful-integrally-1} is \'etale local on $X$, so we can assume that $X=\Spa(R, R^+)$ is an affinoid perfectoid and choose a pseudo-uniformizer $\varpi\in R^+$. Then \cref{lemma:derived-completeness} and \cite[\href{https://stacks.math.columbia.edu/tag/0A0G}{Tag 0A0G}]{stacks-project} imply that it suffices to show that $\cal{P}/^{\rL} \varpi^n \to \rR\lambda_{\tau, *}\lambda_{\tau}^* \big(\cal{P}/^{\rL} \varpi^n\big)$ is an isomorphism for any integer $n\geq 1$. This follows directly from \cref{cor:pullback-etale-finitary} and \cref{lemma:example-o-plus-mod-varpi}. Now Part~\cref{prop:fully-faithful-integrally-2} follows directly from Part~\cref{prop:fully-faithful-integrally-1} and the $(\lambda_\tau^*, \rR\lambda_{\tau, *})$-adjunction. Lastly, Part~\cref{prop:fully-faithful-integrally-3} follows 
    from Part~\cref{prop:fully-faithful-integrally-2} by passing to $\tau^{\leq 0}$.
\end{proof}

Our next goal is to show that $\lambda_\tau^*\colon \Perf_\et(X; \O_X^+) \to \Perf_\tau(X; \O_X^+)$ is essentially surjective and, thus, is an equivalence. For this, we need to introduce the following notation:

\begin{defn} The {\it presheaf of topologically nilpotent functions}  $\O^{\circ\circ} \colon \Perfd^{\aff, \op} \to \Mod_\Z$ is the presheaf which sends $X=\Spa(R, R^+)$ to the ideal of topologically nilpotent elements $R^{\circ \circ}$. For each affinoid perfectoid space $X$, we denote by $\O_{X}^{\circ\circ}$ the restriction $\O^{\circ\circ}|_{\Perfd^{\aff, \op}_{/X}}$. For $\tau\in \{\an, \et, \proet, v\}$, we write $\O_{X_\tau}^{\circ\circ}$ when we want to emphasize that we consider $\O_{X}^{\circ\circ}$ as a presheaf in the $\tau$-topology. 
\end{defn}

The next lemma shows, among other things, that $\O^{\circ\circ}$ is a $\tau$-sheaf. 

\begin{lemma}\label{lemma:projection-formula-topologically-nilpotent-elements} Let $X=\Spa(R, R^+)$ be an affinoid perfectoid space and let $\tau\in \{\an, \et, \proet, v\}$. Then there is a pseudo-uniformizer $\varpi\in R^+$ admitting compatible $p$-power roots $\varpi^{1/p^n}\in R^+$. Furthermore, the following hold:
\begin{enumerate}[leftmargin=*,label={\upshape{(\roman*)}}]
    \item\label{lemma:projection-formula-topologically-nilpotent-elements-0} the presheaf $\O_{X_\tau}^{\circ\circ}$ is a $\tau$-sheaf of abelian groups and the natural morphism $\colim^p_\NN \varpi^{1/p^n} \O_{X_\tau}^+ \to \O_{X_\tau}^{\circ\circ}$ is an isomorphism. In particular, the morphism $\colim_\NN \varpi^{1/p^n} \O_{X_\tau}^+ \to \O_{X_\tau}^{\circ\circ}$ is an isomorphism as well;\footnote{We recall that $\colim^p$ stands for the presheaf colimit.} 
    \item\label{lemma:projection-formula-topologically-nilpotent-elements-1} the natural morphism 
    \[
    \rR\Gamma_\tau(X, \F)\otimes^\rL_{R^+} R^{\circ\circ} \to \rR\Gamma_\tau(X, \F\otimes_{\O_{X_\tau}^+}^\rL \O_{X_\tau}^{\circ\circ})
    \]
    is an isomorphism for any $\F\in \cal{D}^+(X_\tau; \O_{X_\tau}^+)$;
    \item\label{lemma:projection-formula-topologically-nilpotent-elements-2} the natural morphism $R^{\circ\circ} \to \rR\Gamma_\tau(X, \O_{X_\tau}^{\circ\circ})$ is an isomorphism;
\end{enumerate}
\end{lemma}
\begin{proof}
    The claim about the existence of a pseudo-uniformizer $\varpi\in R^+$ with compatible $p$-power roots of $\varpi$ follows from \cite[Lem.~3.10]{diamonds}. Then \cite[Lem.~B.5]{zav-almost} ensures that $S^{\circ\circ} = \bigcup_{n\geq 1} \varpi^{1/p^n} S^+$ for any morphism of affinoid perfectoid spaces $\Spa(S, S^+) \to \Spa(R, R^+)$. This already implies that the natural morphism $\colim^p_\NN \varpi^{1/p^n} \O_{X_\tau}^+ \to \O_{X_\tau}^{\circ\circ}$ is an isomorphism. Using that filtered colimits commute with finite limits and that $\O_{X_\tau}^+$ is a $\tau$-sheaf, we conclude that $\O_{X_\tau}^{\circ\circ}$ is a $\tau$-sheaf as well. This finishes the proof of part~\cref{lemma:projection-formula-topologically-nilpotent-elements-0}. Combining this with the facts that $\rR\Gamma_\tau(X; \blank)$ commutes with filtered colimits of uniformly bounded below complexes and that $\rR\Gamma_\tau(X, \F) \otimes^\rL_{R^+} \blank$ commutes with all colimits, we easily deduce part~\cref{lemma:projection-formula-topologically-nilpotent-elements-1}. Part~\cref{lemma:projection-formula-topologically-nilpotent-elements-2} follows directly from part~\cref{lemma:projection-formula-topologically-nilpotent-elements-1} applied to $\F=\O_{X_\tau}^+$ and \cite[Th.~3.24]{diamonds} for $\tau=\an$, \cite[Th.~6.3]{diamonds} for $\tau=\et$, \cite[Prop.~8.5]{diamonds} for $\tau=\proet$, and \cite[Prop.~8.8]{diamonds} for $\tau=v$.  
\end{proof}

For the next definition, we fix an affinoid perfectoid space $X$ and $\tau\in \{\an, \et, \proet, v\}$. 

\begin{defn} The {\it topologically reduced structure sheaf} $\ov{\O}_{X_\tau} \colon \Perfd^{\aff, \op}_{/X} \to \Mod_\Z$ is the quotient sheaf (in the $\tau$-topology) $\ov{\O}_{X_\tau} \coloneqq \O^+_{X_\tau}/\O_{X_\tau}^{\circ\circ}$. 
\end{defn}

For the purpose of the next definition, we fix an affinoid perfectoid space $X=\Spa(R, R^+)$, a pseudo-uniformizer $\varpi\in R^+$, and $\tau\in \{\an, \et, \proet, v\}$.

\begin{defn} Fix a perfect complex $\cal{P}\in \Perf_\tau(X; \O^+)$ (resp. $\cal{P}\in \Perf_\tau(X; \O^+/\varpi)$). Its {\it topological reduction} is the perfect complex $\ov{\cal{P}}\coloneqq \cal{P} \otimes^{\rL}_{\O_{X_\tau}^+} \ov{\O}_{X_\tau}\in \Perf_\tau(X; \ov{\O})$ (resp. $\ov{\cal{P}}\coloneqq \cal{P} \otimes^{\rL}_{\O_{X_\tau}^+/\varpi} \ov{\O}_{X_\tau}\in \Perf_\tau(X; \ov{\O})$). 
\end{defn}

The next lemma will be an important technical ingredient in our proof of $v$-descent for \'etale perfect $\O^+$-complexes. 

\begin{lemma}\label{lemma:approximation-mod-top-nilpotent} Let $X=\Spa(R, R^+)$ be an affinoid perfectoid space, let $\varpi\in R^+$ be a pseudo-uniformizer as in \cref{lemma:projection-formula-topologically-nilpotent-elements}, let $\tau\in \{\an, \et, \proet, v\}$, and let $\cal{P}, \cal{Q} \in \Perf_\tau(X; \O^+/\varpi^{1/p^n})$ for some integer $n\geq 0$. Then the following hold:
\begin{enumerate}[leftmargin=*,label={\upshape{(\roman*)}}]
    \item\label{lemma:approximation-mod-top-nilpotent-1} the natural morphism $\colim_{m\geq n} \rR\Gamma_\tau(X, \cal{P}/^{\rL} \varpi^{1/p^m}) \to \rR\Gamma_\tau(X, \ov{\cal{P}})$ is an isomorphism;
    \item\label{lemma:approximation-mod-top-nilpotent-2} the natural morphism $\colim_{m\geq n} \rR\Hom(\cal{P}/^{\rL} \varpi^{1/p^m}, \cal{Q}/^{\rL} \varpi^{1/p^m}) \to \rR\Hom(\ov{\cal{P}}, \ov{\cal{Q}})$ is an isomorphism.
\end{enumerate}
\end{lemma}
\begin{proof}
    \cref{lemma:projection-formula-topologically-nilpotent-elements}~\cref{lemma:projection-formula-topologically-nilpotent-elements-0} implies that $\colim_{m\geq n} \cal{P}/^{\rL} \varpi^{1/p^m} \to \cal{P} \otimes^{\rL}_{\O_{X_\tau}^+} \ov{\O}_{X_\tau}$ is an isomorphism. Thus, part~\cref{lemma:approximation-mod-top-nilpotent-1} follows from the observation that $\rR\Gamma_\tau(X, \blank)$ commutes with filtered colimits of uniformly bounded complexes. Part~\cref{lemma:approximation-mod-top-nilpotent-2} follows from part~\cref{lemma:approximation-mod-top-nilpotent-1} since $\rR\Hom(\cal{P}/^{\rL} \varpi^{1/p^m}, \cal{Q}/^{\rL} \varpi^{1/p^m}) \simeq \rR\Gamma_\tau(X, \cal{P}^\vee/^{\rL} \varpi^{1/p^m} \otimes^{\rL}\cal{Q}/^{\rL} \varpi^{1/p^m})$ (and similarly for $\ov{\cal{P}}$ and $\ov{\cal{Q}}$). 
\end{proof}

\cref{lemma:approximation-mod-top-nilpotent} allows us to check various properties of perfect $\O^+$-complexes modulo $\O^{\circ\circ}$. We collect a number of such applications below. These results will be crucial for our proof of the essential surjectivity of $\lambda_\tau^*$. 

\begin{cor}\label{cor:lifting-isomorphisms} Let $X=\Spa(R, R^+)$ be an affinoid perfectoid space, let $\tau\in \{\an, \et, \proet, v\}$, let $\cal{P}, \cal{Q} \in \Perf_\tau(X; \O^+)$, and let $\varphi \colon \cal{P }\to \cal{Q}$ be a morphism such that 
\begin{equation*}
\ov{\varphi} \colon \ov{\cal{P}}  \to \ov{\cal{Q}}
\end{equation*}
is an isomorphism. Then $\varphi$ is an isomorphism as well.
\end{cor}
\begin{proof}
    The $\Hh^0$-version of \cref{lemma:approximation-mod-top-nilpotent}~\cref{lemma:approximation-mod-top-nilpotent-2} and a standard approximation argument implies that there is a pseudo-uniformizer $\varpi\in R^+$ such that $\varphi/^\rL\varpi \colon \cal{P}/^\rL \varpi \to \cal{Q}/^\rL \varpi$ is an isomorphism. By induction, this implies that $\varphi/^\rL\varpi^n \colon \cal{P}/^\rL \varpi^n \to \cal{Q}/^\rL \varpi^n$ is an isomorphism for any integer $n\geq 1$. Then \cref{lemma:derived-completeness} implies that $\varphi$ is an isomorphism as well. 
\end{proof}

As another application of \cref{lemma:approximation-mod-top-nilpotent}, we show that we can lift certain presentations of perfect $\ov{\O}$-complexes to presentations of $\O^+/\varpi$-complexes for a suitable choice of a pseudo-uniformizer $\varpi$. 

\begin{cor}\label{cor:lifting-direct-summands} Let $X=\Spa(R, R^+)$ be an affinoid perfectoid space, let $\tau\in \{\an, \et, \proet, v\}$, and let $\ov{\cE}$ be an $\ov{\O}_{X_\tau}$-module which is a direct summand of a finite free $\ov{\O}_{X_\tau}$-module. Let $\varpi\in R^+$ be a pseudo-uniformizer as in \cref{lemma:projection-formula-topologically-nilpotent-elements}. Then there is an integer $n\geq 0$ and an $\O_{X_\tau}^+/\varpi^{1/p^n}$-module $\cE$ such that $\cE$ is a direct summand of a finite free $\O_{X_\tau}^+/\varpi^{1/p^n}$-module and $\cE \otimes_{\O_{X_\tau}^+/\varpi^{1/p^n}} \ov{\O}_{X_\tau} \simeq \ov{\cE}$. 
\end{cor}
\begin{proof}
    Let $\ov{\cE}$ be a direct summand of $\ov{\O}_{X_\tau}^{\oplus d}$ and let $\ov{\pi} \colon \ov{\O}_{X_\tau}^{\oplus d} \to \ov{\O}_{X_\tau}^{\oplus d}$ be the projector corresponding to $\ov{\cE}$, i.e., $\ov{\pi}^2=\ov{\pi}$ and $\mathrm{Im}(\ov{\pi})=\ov{\cE}$. Then the $\Hh^0$-version of \cref{lemma:approximation-mod-top-nilpotent}~\cref{lemma:approximation-mod-top-nilpotent-1} implies that we can find an integer $n\in \NN$ and a morphism $\pi \colon (\O_{X_\tau}^+/\varpi^{1/p^n})^{\oplus d} \to (\O_{X_\tau}^+/\varpi^{1/p^n})^{\oplus d}$ such that $\pi^2=\pi$ and $\pi \otimes_{\O_{X_\tau}^+/\varpi^{1/p^n}} \ov{\O}_{X_\tau} = \ov{\pi}$. Then $\cE \coloneqq \mathrm{Im}(\pi)$ is a direct summand of $(\O_{X_\tau}^+/\varpi^{1/p^n})^{\oplus d}$ and $\cE \otimes_{\O_{X_\tau}^+/\varpi^{1/p^n}} \ov{\O}_{X_\tau} \simeq \ov{\cE}$. This finishes the proof. 
\end{proof}

\begin{cor}\label{cor:lifting-complexes} Let $X=\Spa(R, R^+)$ be an affinoid perfectoid space, let $\tau\in \{\an, \et, \proet, v\}$, and let $\cal{P}\in \Perf_\tau(X; \O^+)$ such that $\ov{\cal{P}} \simeq \ov{\cal{E}}^\bullet$ where $\ov{\cal{E}}^\bullet$ is a complex of $\ov{\O}_{X_\tau}$-modules such that each $\ov{\cal{E}}^i$ is a direct summand of a finite free $\ov{\O}_{X_\tau}$-module and $\ov{\cal{E}}^i\simeq 0$ for $i\notin [-d, 0]$ for some integer $d\geq 0$. Then there is a pseudo-uniformizer $\varpi \in R^+$ such that $\cal{P}/^{\rL} \varpi$ is isomorphic to a finite complex of direct summands of finite free $\O_{X_\tau}^+/\varpi$-modules $\cal{E}^\bullet$ such that $\cal{E}^\bullet \otimes_{\O_{X_\tau}^+/\varpi} \ov{\O}_{X_\tau} \simeq \ov{\cal{E}}^\bullet$ and $\cal{E}^i\simeq 0$ for $i\notin [-d, 0]$.
\end{cor}
\begin{proof}
    Choose a pseudo-uniformizer $\varpi\in R^+$ as in \cref{lemma:projection-formula-topologically-nilpotent-elements}.
    Then \cref{cor:lifting-direct-summands} and the $\Hh^0$-version of 
    \cref{lemma:approximation-mod-top-nilpotent}~\cref{lemma:approximation-mod-top-nilpotent-2}
    implies that we can find an integer $n\in \NN$ and a finite complex $\cal{E}$ of direct summands of finite free $\O_{X_\tau}^+/\varpi^{1/p^n}$-modules such that $\cal{E}^\bullet \otimes_{\O_{X_\tau}^+/\varpi^{1/p^n}} \ov{\O}_{X_\tau} \simeq \ov{\cal{E}}^\bullet$ and $\cal{E}^i\simeq 0$ for $i\notin [-d, 0]$. Then the $\Hh^0$-version of \cref{lemma:approximation-mod-top-nilpotent}~\cref{lemma:approximation-mod-top-nilpotent-2} and a standard approximation argument imply that there exists $m\geq n$ such that $\cal{E}^\bullet \otimes_{\O_{X_\tau}^+/\varpi^{1/p^n}} \O_{X_\tau}^+/\varpi^{1/p^m} \simeq \cal{P}\otimes^{\rL}_{\O_{X_\tau}^+} \O_{X_\tau}^+/\varpi^{1/p^m}$. Now we redefine $\varpi$ to be $\varpi^{1/p^m}$ and $\cal{E}^\bullet$ to be $\cal{E}^\bullet \otimes_{\O_{X_\tau}^+/\varpi^{1/p^n}} \O_{X_\tau}^+/\varpi^{1/p^m}$. This finishes the proof.
\end{proof}

Now we are ready to get the following crucial almost vanishing result: 

\begin{cor}\label{cor:lifting-of-almost-vanishing} Let $X=\Spa(R, R^+)$ be an affinoid perfectoid space, let $\tau\in \{\an, \et, \proet, v\}$, and let $\cal{P}\in \Perf_\tau(X; \O_X^+)$ such that $\ov{\cal{P}} \simeq \ov{\cal{E}}^\bullet$ for a finite complex of direct summands of finite free $\ov{\O}_{X_\tau}$-modules such that $\ov{\cal{E}}^i\simeq 0$ for $i>0$. Then $\rR\Gamma_\tau(X, \cal{P})\in D^{\leq 0}(R^+)^a$.
\end{cor}
\begin{proof}
    We choose a pseudo-uniformizer $\varpi\in R^+$ as in \cref{cor:lifting-complexes}. Then \cref{lemma:derived-completeness}, \cite[\href{https://stacks.math.columbia.edu/tag/0BLX}{Tag 0BLX}]{stacks-project}, and \cite[Lem.~A.5]{zav-almost} imply that it suffices to show that $\rR\Gamma_\tau(X, \cal{P}/^{\rL}\varpi) \in D^{\leq 0}(R^+)^a$. Using the finite complex $\cal{E}^\bullet$ of direct summands of finite free $\O_{X_\tau}^+/\varpi$-modules constructed in \cref{cor:lifting-complexes}, we reduce the question to showing that $\rR\Gamma_\tau(X, \O_{X_\tau}^+/\varpi)\in D^{\leq 0}(R^+)^a$. This follows directly from \cite[Th.~3.24]{diamonds} for $\tau=\an$, \cite[Th.~6.3]{diamonds} for $\tau=\et$, \cite[Prop.~8.5]{diamonds} for $\tau=\proet$, and \cite[Prop.~8.8]{diamonds} for $\tau=v$.
\end{proof}

Our next goal is to prove a stronger version of \cref{cor:lifting-direct-summands}. For this, we will need the following preliminary lemma. 

\begin{lemma}\label{lemma:lifting-idempotents-algebra} Let $R$ be a commutative ring, let $x\in R$ be a nonzerodivisor that admits a compatible sequence of $p$-power roots $x^{1/p^n}\in R$, and let $J\coloneqq \cup x^{1/p^n}R$. Let $A$ be a (not necessarily commutative) $R$-algebra\footnote{By an $R$-algebra, we mean a ring $A$ with a ring homomorphism $R\to A$ such that the image of $R$ lies in the center of $A$.} which is $x$-adically complete and let $\ov{e}\in A/JA$ be an idempotent element. Then there is an idempotent element $e\in A$ such that its image in $A/JA$ is equal to $\ov{e}$.
\end{lemma}
\begin{proof}
    First, let $N$ be the kernel of the morphism $A/xA \to A/JA$. Then $N$ is a locally nilpotent ideal (i.e., for any $n\in N$ there is an integer $d$ such that $n^d=0$). Therefore, \cite[Prop.~3.6/1]{Lambek} implies that there is an idempotent $e_1\in A/xA$ such that its reduction in $A/JA$ is equal to $\ov{e}$. Similarly, we can apply \cite[Prop.~3.6/1]{Lambek} to the reduction morphisms $A/x^{n}A \to A/x^{n-1}A$ to get a sequence of compatible idempotent elements $e_n\in A/x^n A$. Since $A$ is $x$-adically complete, the sequence $\{e_n\}_{n\geq 1}$ uniquely determines an idempotent $e\in A$ such that the image of $e$ in $A/x^nA$ is equal to $e_n$. In particular, its image in $A/JA$ is equal to $\ov{e}$.  
\end{proof}

Now we are ready to prove a stronger version of \cref{cor:lifting-direct-summands}. 

\begin{cor}\label{cor:lifting-direct-summands-strong} Let $X=\Spa(R, R^+)$ be an affinoid perfectoid space, let $\tau\in \{\an, \et, \proet, v\}$, and let $\ov{\cE}$ be an $\ov{\O}_{X_\tau}$-module which is a direct summand of a finite free $\ov{\O}_{X_\tau}$-module. Then there is an $\O_{X_\tau}^+$-module $\cE$ such that $\cE$ is a direct summand of a finite free $\O_{X_\tau}^+$-module and $\cE \otimes_{\O_{X_\tau}^+} \ov{\O}_{X_\tau} \simeq \ov{\cE}$. 
\end{cor}
\begin{proof}
    Let $\ov{\cE}$ be a direct summand of $\ov{\O}_{X_\tau}^{\oplus d}$ and let $\ov{\pi} \colon \ov{\O}_{X_\tau}^{\oplus d} \to \ov{\O}_{X_\tau}^{\oplus d}$ be the projector corresponding to $\ov{\cE}$, i.e., $\ov{\pi}^2=\ov{\pi}$ and $\mathrm{Im}(\ov{\pi})=\ov{\cE}$. Choose a pseudo-uniformizer $\varpi\in R^+$ as in \cref{lemma:projection-formula-topologically-nilpotent-elements}. It suffices to show that we can lift $\ov{\pi}$ to a projector $\pi \colon (\O_{X_\tau}^+)^{\oplus d} \to (\O_{X_\tau}^+)^{\oplus d}$. In this case, the $\O_{X_\tau}^+$-module $\cE \coloneqq \mathrm{Im}(\pi)$ does the job. 
    
    In order to construct $\pi$, we set $\cV \coloneqq \ud{\End}\bigl((\O_{X_\tau}^+)^{\oplus d}\bigr)$, so we have $\ov{\cV} \coloneqq \cV \otimes_{\O_{X_\tau}^+} \ov{\O}_{X_\tau} \simeq \ud{\End}\bigl((\ov{\O}_{X_\tau})^{\oplus d}\bigr)$. The projector $\ov{\pi}$ corresponds to an idempotent element $\ov{e}_{\ov{\pi}} \in \End\bigl(\ov{\O}_{X_\tau}^{\oplus d}\bigr) = \Hh^0_\tau(X, \ov{\cV})$. Therefore, in order to construct the desired projector $\pi$, it suffices to show that the idempotent $\ov{e}_{\ov{\pi}} \in \Hh^0_\tau(X, \ov{\cV})$ can be lifted to an idempotent $e\in \Hh^0_\tau(X, \cV)$. 

    For this, we consider the short exact sequence $0\to \cV \otimes_{\O_{X_\tau}^+} \O^{\circ\circ}_{X_\tau} \to \cV \to \ov{\cV} \to 0$ and $J\coloneqq \cup_{n} \varpi^{1/p^n}R^+$. Then the combination of \cref{cor:lifting-of-almost-vanishing} and \cref{lemma:projection-formula-topologically-nilpotent-elements}~\cref{lemma:projection-formula-topologically-nilpotent-elements-0},\cref{lemma:projection-formula-topologically-nilpotent-elements-1} imply that 
    \[
    \Hh^0_\tau(X, \ov{\cV}) \simeq \Hh^0_\tau(X, \cV)/\bigl(J\Hh^0_\tau(X, \cV)\bigr).
    \]
    Therefore, \cref{lemma:lifting-idempotents-algebra} guarantees that, for the purpose of constructing the desired idempotent $e\in \Hh^0_\tau(X, \cV)$, it suffices to show that $\Hh^0_\tau(X, \cV)$ is $\varpi$-adically complete. Now the observation that $\cV$ is $\varpi$-torsionfree and derived $\varpi$-adically complete (see \cref{lemma:derived-completeness}) implies that it is classically $\varpi$-adically complete, i.e., $\cV\simeq \lim_n \cV/\varpi^n$. Since $\Hh^0_\tau(X, \blank)$ commutes with limits, we see that $\Hh^0_\tau(X, \cV) \simeq \lim_n \Hh^0_\tau(X, \cV/\varpi^n)$. Thus, \cite[\href{https://stacks.math.columbia.edu/tag/0G1Q}{Tag 0G1Q}]{stacks-project} ensures that $\Hh^0_\tau(X, \cV)$ is a $\varpi$-adically complete $R^+$-algebra. This finishes the proof. 
\end{proof}

\begin{cor}\label{cor:lifting-morphisms} Let $X=\Spa(R, R^+)$ be an affinoid perfectoid space, let $\tau\in \{\an, \et, \proet, v\}$, let $\cal{P}\in \Perf_\tau(X; \O^+)$ such that $\ov{\cal{P}} \simeq \ov{\cal{E}}^\bullet$ for a finite complex of direct summands of finite free $\ov{\O}_{X_\tau}$-modules such that $\ov{\cal{E}}^i\simeq 0$ for $i>0$, let $\cQ$ be a direct summand of a finite free $\O_{X_\tau}^+$-module, and let $\ov{\varphi} \colon \ov{\cQ} \to \ov{\cal{P}}$ be a morphism in $\cD(X_\tau; \ov{\O}_{X_\tau})$. Then there is a morphism $\varphi \colon \cQ \to \cal{P}$ lifting $\ov{\varphi}$.
\end{cor}
\begin{proof}
    First, note that $\rR{\cHom}(\cQ, \cP) \simeq \cal{P} \otimes^\rL_{\O_{X_\tau}^+} \cQ^{\vee}$ and $\rR{\cHom}(\ov{\cQ}, \ov{\cP}) \simeq \ov{\cP} \otimes^\rL_{\ov{\O}_{X_\tau}} \ov{\cQ}^{\vee} \simeq \rR{\cHom}(\cQ, \cP) \otimes^\rL_{\O_{X_\tau}^+} \ov{\O}_{X_\tau}$. Therefore, it suffices to show that the natural morphism $\Hh^0_\tau(X, \cal{P} \otimes^\rL_{\O_{X_\tau}^+} \cQ^{\vee}) \to \Hh^0_\tau(X, \ov{\cal{P}}\otimes^\rL_{\ov{\O}_{X_\tau}} \ov{\cQ}^{\vee})$ is surjective.  Using the fiber sequence
    \[
        \rR\Gamma_\tau(X, \cal{P} \otimes^\rL_{\O_{X_\tau}^+} \cQ^{\vee}\otimes^{\rL}_{\O_{X_\tau}^+} \O_{X_\tau}^{\circ\circ}) \to \rR\Gamma_\tau(X, \cal{P} \otimes^\rL_{\O_{X_\tau}^+} \cQ^{\vee}) \to \rR\Gamma_\tau(X, \ov{\cP} \otimes^\rL_{\ov{\O}_{X_\tau}} \ov{\cQ}^{\vee}),
    \]
    we see that it suffices to show that $\rR\Gamma_\tau(X, \cal{P} \otimes^\rL_{\O_{X_\tau}^+} \cQ^{\vee}\otimes^{\rL}_{\O_{X_\tau}^+} \O_{X_\tau}^{\circ\circ}) \in D^{\leq 0}(R^+)$. Now we observe that the complex $\ov{\cP} \otimes^\rL_{\ov{\O}_{X_\tau}} \ov{\cQ}^{\vee}$ has a representative $\ov{\cal{E}}^\bullet \otimes_{\ov{\O}_{X_\tau}} \ov{\cQ}^{\vee}$ such that each $\ov{\cal{E}}^i \otimes_{\ov{\O}_{X_\tau}} \ov{\cQ}^{\vee}$ is a direct summand of a finite free $\ov{\O}_{X_\tau}$-module and $\ov{\cal{E}}^i \otimes_{\ov{\O}_{X_\tau}} \ov{\cQ}^{\vee} \simeq 0$ for $i>0$. Therefore, the desired vanishing follows from the combination of \cref{cor:lifting-of-almost-vanishing} and \cref{lemma:projection-formula-topologically-nilpotent-elements}~\cref{lemma:projection-formula-topologically-nilpotent-elements-1}.
\end{proof}

\begin{proposition}\label{prop:strictly-perfect-implies-etale} Let $X=\Spa(R, R^+)$ be an affinoid perfectoid space, let $\tau\in \{\proet, v\}$, and let $\cal{P}\in \Perf_\tau(X; \O^+)$ be a perfect complex such that $\ov{\cal{P}}$ can be represented by a finite complex of direct summands of finite free $\ov{\O}_{X_\tau}$-modules. Then $\cal{P}$ lies in the essential image of the functor $\lambda_\tau^* \colon \Perf_\et(X; \O^+) \to \Perf_\tau(X; \O^+)$.
\end{proposition}

\begin{proof}
    We choose a finite complex $\ov{\cal{E}}^\bullet$ of direct summands of finite free $\ov{\O}_{X_\tau}$-modules such that $\ov{\cal{E}}^\bullet \simeq \ov{\cal{P}}$. We show that $\cal{P}$ lies in the essential image of $\lambda_\tau^*$ by induction on the length of $\ov{\cal{E}}^\bullet$.
    
    The base case is $\ov{\cal{E}}^\bullet$ of length $0$, i.e., $\ov{\cal{E}}^\bullet\simeq 0$. In this case, \cref{cor:lifting-isomorphisms} implies that the natural map $0 \to \cal{P}$ is an isomorphism as well. In particular, $\cal{P}$ lies in the essential image of $\lambda_{\tau}^*$.

    Now we perform the induction argument. We assume that we have proven the claim for any $\cal{P}$ with $\ov{\cal{E}}^\bullet$ of length less than $n\geq 1$ and wish to prove it for $\cal{P}$ with $\ov{\cal{E}}^\bullet$ of length $n$. Without loss of generality, we can assume that $\ov{\cal{E}}^i\simeq 0$ for $i\notin [-n+1, 0]$. We use \cref{cor:lifting-direct-summands-strong} to find an $\O_{X_\tau}^+$-module $\cE^0$ such that $\cE^0$ is a direct summand of a finite free $\O_{X_\tau}^+$-module and such that $\cE^0 \otimes_{\O_{X_\tau}^+} \ov{\O}_{X_\tau} \simeq \ov{\cE}^0$. Let $\ov{\varphi} \colon \ov{\cE}^{0} \to \ov{\cal{E}}^\bullet$ be the natural inclusion morphism, then \cref{cor:lifting-morphisms} implies that we can lift this morphism to a morphism $\varphi\colon \cE^0 \to \cal{P}$. We set $\cal{Q} \coloneqq \rm{cone}(\varphi)$. Then the $\Ext^1$-version of \cref{prop:fully-faithful-integrally}~\cref{prop:fully-faithful-integrally-2} implies that it suffices to show that $\cE^0$ and $\cal{Q}$ lie in the essential image of $\lambda_\tau^*$. 
    
    First, we show that $\cE^0$ lies in the essential image of $\lambda_\tau^*$. For this, we note that the facts that $\Perf_\et(X; \cO^+)$ is idempotent-complete and $\lambda_\tau^* \colon \Perf_\et(X; \cO^+) \to \Perf_\tau(X; \cO^+)$ is fully faithful (see \cref{prop:fully-faithful-integrally}~\cref{prop:fully-faithful-integrally-3}) imply that the essential image of $\lambda_\tau^*$ is closed under retracts. Combining this with the observations that $\cE^0$ is a direct summand of a finite free $\cO^+_{X_\tau}$-module and that $\cO^+_{X_\tau}$ lies in the essential image of $\lambda_\tau^*$, we conclude that $\cE^0$ lies in the essential image of $\lambda_\tau^*$. 
    
    Now we show that $\cQ$ lies in the image of $\lambda_\tau^*$. By the construction, we know that $\ov{\cal{Q}} \simeq \sigma^{<0} \ov{\cal{E}}^\bullet$, where $\sigma^{<0}$ is the na\"{i}ve truncation. In particular, the length of $\ov{\cal{E}}^\bullet$ representing $\ov{\cQ}$ is $n-1$, so $\cal{Q}$ lies in the essential image of $\lambda_\tau^*$ by the induction hypothesis. This finishes the proof.
\end{proof}

Finally, we are ready to show that $\lambda_\tau^* \colon \Perf_\et(X; \O_X^+) \to \Perf_\tau(X; \O_X^+)$ is essentially surjective for $\tau\in \{\proet, v\}$.

\begin{thm}\label{thm:v-perfect-complexes-vs-etale-perfect-complexes-integrally} Let $X$ be a perfectoid space, let $\tau\in \{\proet, v\}$, and let $\lambda_\tau \colon (X_\tau, \O^+) \to (X_\et,  \O^+)$ be the natural morphism of ringed sites. Then the following hold:
\begin{enumerate}[leftmargin=*,label={\upshape{(\roman*)}}]
    \item\label{thm:v-perfect-complexes-vs-etale-perfect-complexes-integrally-1} the functor $\rR\lambda_{\tau, *} \colon \cal{D}(X_\tau; \O^+) \to \cal{D}(X_\et; \O^+)$ restricts to a functor $\rR\lambda_{\tau, *}\colon \Perf_\tau(X; \O^+) \to \Perf_\et(X; \O^+)$;
    \item\label{thm:v-perfect-complexes-vs-etale-perfect-complexes-integrally-2} the functor $\lambda_\tau^* \colon \Perf_\et(X; \O^+) \to \Perf_\tau(X; \O^+)$ is an equivalence with the inverse functor $\rR\lambda_{\tau, *}$.
\end{enumerate}
\end{thm}
\begin{proof}
    We start by proving Part~\cref{thm:v-perfect-complexes-vs-etale-perfect-complexes-integrally-1}. We fix a $\cal{P}\in \Perf_\tau(X; \O^+)$ and need to show that $\rR\lambda_{\tau, *}\cal{P}$ lies in $\Perf_\et(X; \O^+)$. The question is \'etale local on $X$, so we can assume that $X=\Spa(R, R^+)$ is an affinoid perfectoid. We choose a pseudo-uniformizer $\varpi\in R^+$. Then \cref{cor:pullback-etale-finitary}\cref{cor:pullback-etale-finitary-1} and \cref{lemma:example-o-plus-mod-varpi} ensure that we can replace $X$ with its \'etale covering to assume that $\cal{P}/^{\rL} \varpi$ is strictly perfect, i.e., represented by a finite complex of direct summands of a finite free $\O_{X_\tau}^+/\varpi$-modules. 
    In particular, $\ov{\cal{P}}$ is strictly perfect as well. Then \cref{prop:strictly-perfect-implies-etale} ensures that $\cal{P} \simeq \lambda_{\tau}^*\cal{P}'$ for some $\cal{P}'\in \Perf_\et(X; \O^+)$. Therefore, \cref{prop:fully-faithful-integrally} guarantees that the natural map $\cal{P}' \to \rR\lambda_{\tau, *}\cal{P}$ is an isomorphism. In particular, $\rR\lambda_{\tau, *}\cal{P}$ indeed lies in $\Perf_\et(X; \O^+)$. 

    Now we show Part~\cref{thm:v-perfect-complexes-vs-etale-perfect-complexes-integrally-2}. \cref{prop:fully-faithful-integrally} implies that it suffices to show that, for any $\cal{P}\in \Perf_\tau(X; \O^+)$, the natural map $\lambda_\tau^* \rR\lambda_{\tau, *}\cal{P} \to \cal{P}$ is an isomorphism. The question is \'etale local on $X$, so we can assume that $X=\Spa(R, R^+)$ is an affinoid perfectoid space. We choose a pseudo-uniformizer $\varpi\in R^+$ and note that Part~\cref{thm:v-perfect-complexes-vs-etale-perfect-complexes-integrally-1} and \cref{lemma:derived-completeness} imply that it suffices to show that $\lambda_\tau^* \rR\lambda_{\tau, *}\big(\cal{P}/^{\rL}\varpi\big) \to \cal{P}/^{\rL}\varpi$ is an isomorphism. This follows from \cref{cor:pullback-etale-finitary}\cref{cor:pullback-etale-finitary-3} and \cref{lemma:example-o-plus-mod-varpi}. 
\end{proof}

\begin{cor}\label{cor:v-descent-integrally} The functor $\Perf_\et(\blank, \O^+) \colon \Perfd^{\aff,\op} \to \Cat_\infty$ is a hypercomplete $v$-sheaf of $\infty$-categories.
\end{cor}
\begin{proof}
    It follows immediately from \cref{thm:v-perfect-complexes-vs-etale-perfect-complexes-integrally}~\cref{thm:v-perfect-complexes-vs-etale-perfect-complexes-integrally-2} and \cref{lemma:perfect-complexes-are-hypercomplete}. 
\end{proof}

\begin{cor}\label{cor:integral-perfect-complex-on-strictly-totally-disconnected} Let $X=\Spa(R, R^+)$ be a strictly totally disconnected perfectoid space and let $\tau\in \{\et, \proet, v\}$. Then the functor $\widetilde{(\blank)} \colon \Perf(R^+) \to \Perf_\tau(X; \O^+)$ is an equivalence. Its inverse is given by 
\[
\rR\Gamma_\tau(X, \blank) \colon \Perf_\tau(X; \O^+) \to \Perf(R^+).
\]
\end{cor}
\begin{proof}
    The proof of this claim is identical to the proof of \cref{cor:etale-perfect-complex-on-strictly-totally-disconnected} using \cref{thm:v-perfect-complexes-vs-etale-perfect-complexes-integrally} in place of \cref{cor:pullback-etale-finitary}.
\end{proof}

For future reference, we also need to establish a version of \cref{cor:v-descent-integrally} for perfect complexes of bounded tor-amplitude. For this, we will use the sheaf $\Perf^{[r,s]}_\tau(\blank; \O^+) \colon \Perfd^{\aff, \op} \to \Cat_\infty$ from \Cref{notation:sheaf-perfect-complexes-finite-amplitude-cA}. 

Our goal is to show that the sheafification morphism $\Perf_\et^{[r, s]}(\blank; \O^+) \to \Perf_v^{[r,s]}(\blank; \O^+)$ is an equivalence.

\begin{rmk}\label{rmk:fully-faithful-finite-tor-amplitude-integral} The morphism $\Perf_\et^{[r, s]}(\blank; \O^+) \to \Perf_v^{[r, s]}(\blank; \O^+)$ is fully faithful for any $r, s\in \Z \cup \{-\infty, \infty\}$. Indeed, \cref{lemma:sheafification-of-fully-faithful} and \cref{cor:v-descent-integrally} imply that both $\Perf_\et^{[r, s]}(\blank; \O^+)$ and $\Perf_v^{[r, s]}(\blank; \O^+)$ are full subfunctors of $\Perf_\et(\blank; \O^+) \simeq \Perf_v(\blank; \O^+)$. Therefore,  $\Perf_\et^{[r, s]}(\blank; \O^+) \to \Perf_v^{[r, s]}(\blank; \O^+)$ is fully faithful as well.
\end{rmk}

\begin{lemma}\label{lemma:surjective-integrally} Let $(R, R^+) \to (S, S^+)$ be a continuous map of complete uniform Tate--Huber pairs such that the induced morphism $\Spa(S, S^+) \to \Spa(R, R^+)$ is surjective. Then the morphism $\Spec S^+ \to \Spec R^+$ is surjective on closed points.
\end{lemma}
\begin{proof}
    Fix a pseudo-uniformizer $\varpi\in R^+$. Since $R^+$ and $S^+$ are $\varpi$-adically complete, \cite[\href{https://stacks.math.columbia.edu/tag/05GI}{Tag 05GI}]{stacks-project} implies that $\varpi\in \rad(R^+)$ and $\varpi\in \rad(S^+)$. Therefore, it suffices to show that $\Spec S^+/\varpi \to \Spec R^+/\varpi$ is surjective. Now we note that the specialization morphism $\sp\colon \abs{\Spa(R, R^+)} \to \abs{\Spec R^+/\varpi}$ is surjective (see \cite[Prop.~II.3.1.5]{FujKato}). Therefore, surjectivity of $\Spec S^+/\varpi \to \Spec R^+/\varpi$ follows from the surjectivity of $\Spa(S, S^+)\to \Spa(R, R^+)$ and the functoriality of the specialization morphism. 
\end{proof}

\begin{cor}\label{cor:test-tor-amplitude-v-locally} Let $f\colon X \to Y$ be a $v$-covering in $\Perfd^\aff$, let $r, s\in \Z \cup \{-\infty, \infty\}$ with $r \leq s$, and let $\cal{P}\in \Perf^{[r, s]}_v(Y; \O^+)$ such that $f_v^*\cal{P} \in \Perf^{[r, s]}_\et(X; \O^+)$. Then $\cal{P}$ lies in $\Perf^{[r, s]}_\et(Y; \O^+)$ (see \cref{rmk:fully-faithful-finite-tor-amplitude-integral}).
\end{cor}
\begin{proof}
    Since $\Perf^{[r, s]}_\et(\blank; \O^+)$ is an \'etale sheaf, we can check that $\cal{P}$ lies in $\Perf^{[r, s]}_\et(Y; \O^+)$ \'etale locally on $Y$ (see \cref{cor:descent-data}).  \cref{prop:perfect-complexes-sheafification} and \cref{lemma:sheafification-of-fully-faithful} imply that $\Perf^{[r, s]}_\tau(\blank; \O^+) \to \Perf_\tau(\blank; \O^+)$ is fully faithful for $\tau \in \{\et, v\}$, while \cref{thm:v-perfect-complexes-vs-etale-perfect-complexes-integrally}\cref{thm:v-perfect-complexes-vs-etale-perfect-complexes-integrally-2} implies that $\Perf_\et(\blank; \O^+) \simeq \Perf_v(\blank; \O^+)$. Combining these results with \cref{lemma:locally-essentially-surjective} and \cref{prop:perfect-complexes-sheafification}, we see that there is an \'etale covering $Y' \to Y$
    such that $Y' = \Spa(R, R^+)$ is an affinoid perfectoid and $\cal{P}|_{Y'} \simeq \widetilde{P}$ for some $P\in \Perf(\O_{Y'}^+(Y'))$. We wish to show that $P$ lies in $\Perf^{[r, s]}(\O_{Y'}^+(Y'))$.
    
    Now we use \cref{lemma:locally-essentially-surjective} once again to find a $v$-covering $X' \to X\times_{Y}Y'$ such that $\cal{P}|_{X'} \simeq \widetilde{Q}$ for some $Q\in \Perf^{[r, s]}\big(\O_{X'}^+(X')\big)$. Without loss of generality, we can assume that $X' = \Spa(S, S^+)$ is strictly totally disconnected. Then \cref{cor:integral-perfect-complex-on-strictly-totally-disconnected} implies that $P\otimes^\rL_{\O_{Y'}^+(Y')} \O_{X'}^+(X')\simeq Q \in \Perf^{[r, s]}\big(\O_{X'}^+(X')\big)$. Therefore, \cref{lemma:surjective-integrally} and \cite[\href{https://stacks.math.columbia.edu/tag/068V}{Tag 068V}]{stacks-project} imply that $P\in \Perf^{[r, s]}\big(\O_{Y'}^+(Y')\big)$. This finishes the proof. 
\end{proof}

\begin{cor}\label{cor:integral-descend-finite-tor-amplitude} The natural morphism $\Perf_\et^{[r, s]}(\blank; \O^+) \to \Perf_v^{[r, s]}(\blank; \O^+)$ is an equivalence for any $r, s\in \Z \cup \{-\infty, \infty\}$ with $r \leq s$. In particular, $\Perf_\et^{[r, s]}(\blank; \O^+)$ is a hypercomplete $v$-sheaf. 
\end{cor} 
\begin{proof}
    \cref{rmk:fully-faithful-finite-tor-amplitude-integral} says that 
    $\Perf_\et^{[r, s]}(\blank; \O^+) \to \Perf_v^{[r, s]}(\blank; \O^+)$ is fully faithful. Now \cref{cor:test-tor-amplitude-v-locally} and \cref{lemma:locally-essentially-surjective} imply that this morphism is essentially surjective. Thus, it is an isomorphism. Finally, we note that $\Perf_v^{[r, s]}(\blank; \O^+)$ is a hypercomplete sheaf by \cref{lemma:sheaves-are-hypersheaves} when $r$ and $s$ are finite. The case of infinite $r$ or $s$ follows from \cref{lemma:colimit-of-sheaves} by passing to a colimit.
\end{proof}

\begin{rmk} For $r=s=0$, \cref{cor:integral-descend-finite-tor-amplitude} recovers the fact that $\Vect_\et(\blank; \O^+)$ is a $v$-stack of (ordinary) categories (see \cite[Th.~2.21]{Heuer-G-torsor} and \cite[Th.~1.3.2]{zav-almost}). 
\end{rmk}

\subsection{$v$-descent for perfect $\AA_I$-complexes}

In this subsection, we fix a closed interval $I=[a,b] \subset (0, \infty)$ with rational endpoints. We refer to \cref{defn:period-rings} for the definition of $\AA_I(X)$ for an affinoid perfectoid $X$ over $\Spa(\Q_p, \Z_p)$, we also refer to \cref{cor:AI-etale-sheaf} for the proof that $\AA_I$ is a $v$-sheaf on $\Perfd^{\aff}_{/\Q_p}$. We denote by $\AA_{X_\et}$ the restriction of $\AA_I$ to the (affinoid) \'etale site of $X$. Let 
\[
\rm{Perf}_\et(\blank; \AA_{I}) \colon \Perfd^{\aff, \op}_{/\Q_p} \to \Cat_\infty
\]
be the presheaf of $\infty$-categories that sends an affinoid perfectoid space $X$ to the $\infty$-category of perfect complexes of $\AA_{I, X_\et}$-modules. The main goal of this subsection is to show that this is a hypercomplete sheaf of $\infty$-categories. 

We will closely follow the strategy taken in \cref{section:descend-O-plus}. Since the arguments in this section are very similar to those in \cref{section:descend-O-plus}, we will sometimes not present the proofs in full detail and instead only discuss the main changes one needs to do compared to the analogous proofs in \cref{section:descend-O-plus}. 

\begin{lemma}\label{lemma:A_I-cohomology-of-nice-affinoids} Let $S=\Spa(R, R^+)$ be an affinoid perfectoid space over $\Spa(\Q_p, \Z_p)$ with a choice of $p^\flat\in R^{\flat, +}$. Suppose that there are elements $\alpha, \beta\in R^{\flat, +}$ well-adapted to $p^\flat$ and $I$. Let $\tau\in \{\an, \et, \proet, v\}$. Then $\rR\Gamma_\tau(S, \AA_I) \simeq^a \AA_I(S)[0]$, where we consider almost mathematics with respect to the ideal $\AA_I(S)^{\circ\circ} \subset \AA_I(S)$ (see \cref{cor:almost-mathematics-AI}).
\end{lemma}
\begin{proof}
    By definition, we have $\Hh^0_\tau(S, \AA_I)=\AA_I(S)$. The rest follows immediately from \cite[\href{https://stacks.math.columbia.edu/tag/03F9}{Tag 03F9}]{stacks-project} and \cref{lemma:almost-acyclic}.  
\end{proof}

\begin{lemma}\label{lemma:derived-p-complete-AI} Let $S$ be a perfectoid space over $\Spa(\Q_p, \Z_p)$, let $\tau\in \{\et, \proet, v\}$, and let $\cal{P}\in \Perf_\tau(S; \AA_I)$. Then $\cal{P}$ is derived $p$-adically complete. 
\end{lemma}
\begin{proof}
    The claim is $\tau$-local on $S$, so we can assume that $S=\Spa(R, R^+)$ is an affinoid perfectoid space and that $\cal{P}$ is a strictly perfect complex. Then a simple reduction reduces the claim to the case $\cP = \AA_I$.

    Then \cite[Lem.~A.8]{zav-almost} implies that it suffices to show that $\rR\Gamma_\tau(X, \AA_I)$ is derived $p$-complete for any $X\in \Perfd^{\aff}_{/S}$. Furthermore, \cite[Prop.\,3.4.2 and 3.4.4]{proetale} imply that it suffices to show that $\Hh^i_\tau(X, \AA_I)$ is $p$-adically complete for any $X\in \Perfd^{\aff}_{/S}$ and $i\geq 0$.  
    
    For this, we choose an element $p^\flat\in R^{+, \flat}$ as in \cref{notation:p-flat}. Then \cref{rmk:often-well-adapted-pairs} implies that we can also assume that there are elements $\alpha, \beta\in R^{\flat, +}$ well-adapted to $p^\flat$ and $I$. So \cref{cor:almost-mathematics-AI} imply that $\Hh^0_\tau(X, \AA_I) = \AA_I(X)$ is $p$-adically complete for any $X\in \Perfd^{\aff}_{/S}$. Now \cref{lemma:A_I-cohomology-of-nice-affinoids} ensures that $\Hh^i_\tau(X, \AA_I)\simeq^a 0$ for any $X\in \Perfd^{\aff}_{/S}$ and $i>0$. In particular, it is $p$-torsion, so it is $p$-adically complete. This finishes the proof. 
\end{proof}

Now we are ready to prove the main result of this subsection: 

\begin{thm}\label{thm:v-perfect-complexes-vs-etale-perfect-complexes-integrally-A_I} Let $S$ be a perfectoid space over $\Spa(\Q_p, \Z_p)$, let $\tau\in \{\proet, v\}$, and let $\lambda_\tau \colon (S_\tau, \AA_I) \to (S_\et,  \AA_I)$ be the natural morphism of ringed sites. Then the following hold:
\begin{enumerate}[leftmargin=*,label={\upshape{(\roman*)}}]
\item\label{thm:v-perfect-complexes-vs-etale-perfect-complexes-integrally-A_I-1} the functor $\lambda_\tau^* \colon \Perf_\et(S; \AA_I) \to \Perf_\tau(S; \AA_I)$ is an equivalence with the inverse given by 
\[
\rR\lambda_{\tau, *}\colon \Perf_\tau(S; \AA_I) \to \Perf_\et(S; \AA_I);
\]
\item\label{thm:v-perfect-complexes-vs-etale-perfect-complexes-integrally-A_I-2} if $S=\Spa(R, R^+)$ is strictly totally disconnected and $\tau'\in \{\et, \proet, v\}$, the functor $\widetilde{(\blank)} \colon \Perf\big(\AA_I(S)\big) \to \Perf_{\tau'}(S; \AA_I)$ is an equivalence with inverse given by $\rR\Gamma_{\tau'}(S, \blank) \colon \Perf_{\tau'}(S; \AA_I) \to \Perf\big(\AA_I(S)\big)$.
\end{enumerate}
\end{thm}

\begin{proof}
    First, we note that Part~\cref{thm:v-perfect-complexes-vs-etale-perfect-complexes-integrally-A_I-2} follows formally from Part~\cref{thm:v-perfect-complexes-vs-etale-perfect-complexes-integrally-A_I-1} (see \cref{cor:integral-perfect-complex-on-strictly-totally-disconnected} for a similar deduction). Furthermore, the fact that $\lambda_\tau^* \colon \Perf_\et(S; \AA_I) \to \Perf_\tau(S; \AA_I)$ is fully faithful follows formally from \cref{cor:pullback-etale-finitary}, \cref{lemma:example-AI-mod-p}, and \cref{lemma:derived-p-complete-AI} (see \cref{prop:fully-faithful-integrally} for a similar deduction). Therefore, we only need to show that $\lambda_\tau^*$ is essentially surjective and its inverse is given by $\rR\lambda_{\tau, *}$. These claims are \'etale local on $S$, so we can assume that $S=\Spa(R, R^+)$ is affinoid perfectoid with a choice $p^\flat\in R^{\flat, +}$ and we can also assume that there are elements $\alpha, \beta\in R^{\flat, +}$ well-adapted to $p^\flat$ and $I$ (see \cref{rmk:often-well-adapted-pairs}). In this case, the proof of \cref{thm:v-perfect-complexes-vs-etale-perfect-complexes-integrally} adapts almost verbatim. 

    For the reader's convenience, we sketch the main points of the proof. We first define the sheaf of topologically nilpotent elements $\AA_I^{\circ\circ}\colon \Perfd^{\aff, \op}_{/S} \to \Mod_\ZZ$ which sends $X$ to $\AA_I(X)^{\circ\circ}$ and the quotient $\tau$-sheaf $\ov{\AA}_I \coloneqq \AA_I/\AA_I^{\circ\circ} \colon  \Perfd^{\aff, \op}_{/S} \to \Mod_\ZZ$. Then one proves analogues of \cref{lemma:projection-formula-topologically-nilpotent-elements} and \cref{lemma:approximation-mod-top-nilpotent} with $\varpi \in R^+$ replaced with $[p^\flat]\in \AA_I(X)$. Then one formally deduces \cref{cor:lifting-isomorphisms}, \cref{cor:lifting-direct-summands}, \cref{cor:lifting-complexes}\footnote{The correct conclusion should be that there exists an $n$ such that $\cal{P}/^\rL [p^\flat]^{1/p^n}$ is isomorphic to a finite complex of direct summands of finite free $\AA_{I, X_\tau}/[p^\flat]^{1/p^n}$-modules $\cal{E}^\bullet$ such that $\cal{E}^\bullet \otimes_{\AA_{I, X_\tau}} \ov{\AA}_{I, X_\tau} \simeq \ov{\cal{E}}^\bullet$ and $\cal{E}^i=0$ for $i\notin [-d, 0]$.}, \cref{cor:lifting-of-almost-vanishing}\footnote{Instead of using that $\rR\Gamma_\tau(S, \O_{S_\tau}^+/\varpi)\in D^{\leq 0}(R^+)^a$, we need to use that $\rR\Gamma_\tau(S, \AA_I/[p^\flat]^{1/p^n}) \in D^{\leq 0}(\AA_I(S))^a$. This follows directly from \cref{lemma:A_I-cohomology-of-nice-affinoids}.}, \cref{cor:lifting-direct-summands-strong}, and \cref{cor:lifting-morphisms}. With these corollaries at hand, one can deduce an analogue of \cref{prop:strictly-perfect-implies-etale}. Then one can repeat the proof of \cref{thm:v-perfect-complexes-vs-etale-perfect-complexes-integrally} using \cref{cor:pullback-etale-finitary} applied to $\cA=\AA_I/p$ (see \cref{lemma:example-AI-mod-p}) to ensure that, for any perfect complex $\cal{P} \in \Perf_\tau(X; \AA_I)$, we can find an \'etale covering $X' \to X$ such that $\ov{\cal{P}}|_{X'}$ is strictly perfect. 
\end{proof}

Now we wish to get a version of \cref{cor:integral-descend-finite-tor-amplitude} for $\Perf_\et(\blank; \AA_I)$. For this, we will use the sheaf $\Perf^{[r,s]}_\tau(\blank; \AA_I) \colon \Perfd^{\aff, \op}_{/\QQ_p} \to \Cat_\infty$ from \Cref{notation:sheaf-perfect-complexes-finite-amplitude-cA}. 

\begin{lemma}\label{lemma:surjective-Y} Let $X \to Y$ be a surjective morphism of affinoid perfectoid spaces over $\Spa(\Q_p, \Z_p)$. Then the induced morphism $\Spec \AA_I(X) \to \Spec \AA_I(Y)$ is surjective on closed points.
\end{lemma}
\begin{proof}
    First, \cref{lemma:permanence-properties-Y} ensures that $\Spa(\BB_I(X), \AA_I(X)) \to \Spa(\BB_I(Y), \AA_I(Y))$ is surjective. Therefore, \cref{lemma:surjective-integrally} implies that $\Spec \AA_I(X) \to \Spec \AA_I(Y)$ is surjective on closed points. 
\end{proof}

Finally, we are ready to get the main result of this subsection. 

\begin{cor}\label{cor:v-descent-A_I} The natural morphism $\Perf_\et^{[r, s]}(\blank; \AA_I) \to \Perf_v^{[r, s]}(\blank; \AA_I)$ is an equivalence for any $r, s\in \Z \cup \{-\infty, \infty\}$ with $r \leq s$. In particular, $\Perf_\et^{[r, s]}(\blank; \AA_I)$ is a hypercomplete $v$-sheaf. 
\end{cor}
\begin{proof}
    When $r=-\infty$ and $s=\infty$, this follows immediately from \cref{thm:v-perfect-complexes-vs-etale-perfect-complexes-integrally-A_I}~\cref{thm:v-perfect-complexes-vs-etale-perfect-complexes-integrally-A_I-1} and \cref{lemma:perfect-complexes-are-hypercomplete}. For other $r$ and $s$, one can adapt the proof 
    of \cref{cor:integral-descend-finite-tor-amplitude} using \cref{lemma:surjective-Y} in place of \cref{lemma:surjective-integrally} (which is used in the proof of \cref{cor:test-tor-amplitude-v-locally}). 
\end{proof}

\subsection{$v$-descent for perfect $\mathcal{O}$-complexes}\label{section:descent-O}

Let
\[
\Perf\big(\O(\blank)\big) \colon \Perfd^{\aff, \op} \to \Cat_\infty
\]
be the presheaf of $\infty$-categories that sends an affinoid perfectoid space $X$ to the $\infty$-category of perfect complexes $\Perf\big(\O_X(X)\big)$. The main goal of this subsection is to show that this is a hypercomplete 
$v$-sheaf of $\infty$-categories.

\begin{lemma}\label{lemma:etale-descent-rationally} The presheaves  $\Perf^{[r, s]}\big(\O(\blank)\big) \colon \Perfd^{\aff, \op} \to \Cat_\infty$ are \'etale sheaves of $\infty$-categories for any $r, s\in \ZZ \cup \{-\infty, \infty\}$
with $r \leq s$. In particular, $\Perf\big(\O(\blank)\big)$ is an \'etale sheaf. 
\end{lemma}
\begin{proof}
    Let $F$ be any of the presheaves in the formulation of the lemma. Clearly, $F(\sqcup_{i=1}^n(X_i)) \simeq \prod^n_{i=1} F(X_i)$ for any $X_1, \dots, X_n\in \Perfd^\aff$. Therefore,  \cite[Prop.~A.3.3.1]{SAG} implies that it suffices to show that $F$ satisfies descent with respect to any \'etale covering $X \to Y$. In other words, we have to show that any \'etale covering $f\colon X \to Y$ is of universal $F$-descent in the sense of \cite[Def.~3.3.1]{Liu-Zheng}. Then \cite[Lem.~3.3.2]{Liu-Zheng} and \cite[Prop.~8.2.20]{Kedlaya-Liu-1} implies that it suffices to show that $F$ satisfies analytic and finite \'etale descent. Analytic descent follows from \cite[Th.~5.3]{Andreychev}. Now let $\Spa(S, S^+) \to \Spa(R, R^+)$ be a finite \'etale surjection, then \cite[Def.~6.2]{diamonds} and \cite[Lem.~C.1.2]{zav-almost} imply that $\Spec S \to \Spec R$ is a finite \'etale surjection. Furthermore, the combination of \cite[Lem.~C.1.1]{zav-almost} and \cite[Lem.~B.3.5]{Z-quotients} imply that $S^{\wdh{\otimes}^n_R} \simeq S^{\otimes^n_R}$. Therefore, the result follows from the usual fppf descent for perfect complexes of fixed tor-amplitude by \cite[Cor.~D.6.3.3 \& Prop.~2.8.4.2]{SAG}. 
\end{proof}

Now we wish to deduce that $\Perf\big(\O(\blank)\big)$ is a (hypercomplete) $v$-sheaf. We do this by comparing this functor to $\Perf_v(\blank; \O)$. We start with the following proposition:

\begin{proposition}\label{prop:fully-faithful-rationally} Let $X$ be a perfectoid space, let $\cal{P}, \cal{Q}\in \Perf_\et(X; \O)$, let $\tau\in \{\proet, v\}$, and let $\lambda_\tau \colon (X_\tau, \O) \to (X_\et,  \O)$ be the natural morphism of ringed sites. Then the following statements hold:
\begin{enumerate}[leftmargin=*,label={\upshape{(\roman*)}}]
    \item\label{prop:fully-faithful-rationally-1} the natural morphism $\cal{P} \to \rR\lambda_{\tau, *} \lambda^*_\tau\cal{P}$ is an isomorphism;
    \item\label{prop:fully-faithful-rationally-2} the natural morphism $\rR\Hom(\cal{P}, \cal{Q}) \to \rR\Hom(\lambda_\tau^*\cal{P}, \lambda_\tau^*\cal{Q})$ is an isomorphism;
    \item\label{prop:fully-faithful-rationally-3} the functor $\lambda_\tau^*\colon \Perf_\et(X; \O) \to \Perf_\tau(X; \O)$ is fully faithful.
\end{enumerate}
\end{proposition}
\begin{proof}
    As in the proof of \cref{prop:fully-faithful-integrally}, it suffices to prove Part~\cref{prop:fully-faithful-rationally-1}. The other parts follow formally from it. 

    Now the claim of \cref{prop:fully-faithful-rationally-1} is clearly \'etale local on $X$, so we can assume that $X=\Spa(R, R^+)$ is an affinoid perfectoid space with a pseudo-uniformizer $\varpi\in R^+$ and that $\cal{P}$ is a strictly perfect complex. Then a standard reduction argument reduces us to the case $\cal{P}=\O_{X_\et}$. In this case, we need to show that the natural morphism $\O_{X_\et} \to \rR\lambda_{\tau, *}\O_{X_\tau}$ is an isomorphism. Since $\lambda$ is a coherent morphism of coherent topoi, we conclude that $\O_{X_\et}\simeq \O_{X_\et}^+[\frac{1}{\varpi}]$ and $\rR\lambda_{\tau, *} \O_{X_\tau} \simeq \big(\rR\lambda_{\tau, *} \O_{X_\tau}^+\big)[\frac{1}{\varpi}]$. Then the result follows directly from \cref{prop:fully-faithful-integrally}~\cref{prop:fully-faithful-integrally-1} applied to $\cal{P}=\O_{X_\et}^+$. 
\end{proof}

Now we wish to show that $\lambda_\tau^*$ is also essentially surjective. We do this by studying the sheaf of lattices. For this, we will use the sheaf $\Perf^{[r,s]}_\tau(\blank; \O) \colon \Perfd^{\aff, \op} \to \Cat_\infty$ from \Cref{notation:sheaf-perfect-complexes-finite-amplitude-cA}. 

First, we have the following series of easy remarks.

\begin{rmk}\label{rmk:fully-faithful-finite-tor-amplitude} 
\begin{enumerate}
    \item\label{rmk:fully-faithful-finite-tor-amplitude-1}  \cref{lemma:sheafification-of-fully-faithful} implies that the natural functor $\Perf^{[r, s]}_\tau(\blank; \O) \to \Perf_\tau^{[r', s']}(\blank; \O)$ is fully faithful for $\tau\in \{\an, \et, \proet, v\}$ and $r, r', s, s'\in \Z \cup \{-\infty, \infty\}$ with $r'\leq r \leq s \leq s'$. Therefore, \cref{lemma:colimit-of-sheaves} implies that $\colim_\NN \Perf^{[-n,n]}_\tau(X; \O) \to \Perf_\tau(X; \O)$ is an equivalence for any $X\in \Perfd^\aff$;
    \item\label{rmk:fully-faithful-finite-tor-amplitude-2}  \cref{prop:fully-faithful-rationally}~\cref{prop:fully-faithful-rationally-3} implies that $\Perf_\et(\blank; \O) \to \Perf_\tau(\blank; \O)$ is fully faithful for $\tau\in \{\proet, v\}$. Combining it with \ref{rmk:fully-faithful-finite-tor-amplitude-1}, we conclude that $\Perf^{[r, s]}_\et(\blank; \O) \to \Perf^{[r, s]}_\tau(\blank; \O)$ is fully faithful for any $r, s\in \Z \cup \{-\infty, \infty\}$ and $\tau$ as above;
    \item\label{rmk:fully-faithful-finite-tor-amplitude-3}  \cref{lemma:etale-descent-rationally} implies that $\Perf^{[r,s]}\big(\O(\blank)\big) \simeq \Perf^{[r, s]}_\et(\blank; \O)$ for any $r, s \in \Z \cup \{-\infty, \infty\}$
    with $r \leq s$.
\end{enumerate}
\end{rmk}

Now we fix an affinoid perfectoid space $S$ and a perfect complex $\cal{P}\in \Perf^{[r, s]}_\tau(S; \O)$ for $\tau\in \{\proet, v\}$ and $r, s\in \Z \cup \{-\infty, \infty\}$
with $r \leq s$. We wish to show that $\cal{P}$ lies in the essential image of $\lambda_\tau^*$, we do this by studying the sheaf of lattices of $\cP$. To make this into a precise definition, we first note that a perfect complex $\cal{P} \in \Perf^{[r, s]}_\tau(S; \O)$ defines a morphism of sheaves 
\[
[\cal{P}] \colon * \to \Perf_\tau^{[r, s]}(\blank; \O) \colon \Perfd^{\aff, \op}_{/S} \to \Cat_\infty.
\]

Using this observation, we define the sheaf of lattices as follows: 

\begin{defn} The {\it $\tau$-sheaf of lattices} $\Latt_{\tau, \cal{P}}^{[r, s]}(\blank) \colon \Perfd^{\aff, \op}_{/S} \to \Cat_\infty$ is the fiber product
\[
\Latt^{[r, s]}_{\tau, \cal{P}}(\blank) \coloneqq * \times_{\Perf^{[r, s]}_\tau(-; \O)} \Perf^{[r, s]}_\tau(-; \O^+).
\]
\end{defn}

More informally, for every affinoid perfectoid $S$-space $X$, we have
\[
\Latt^{[r, s]}_{\tau, \cal{P}}(X)= \Big\{ \cal{P}^+ \in \Perf^{[r, s]}_\tau(X; \O^+) \text{ and } \varphi \colon \cal{P}^+ \otimes^{\rL}_{\O_{X_\tau}^+} \O_{X_\tau} \xr{\sim} \restr{\cP}{X} \Big\}.
\]

\begin{lemma}\label{lemma:lattices-are-finitary} Let $S=\Spa(A, A^+)$ be an affinoid perfectoid and let $\cal{P}$ be an element of $\Perf^{[r, s]}_\tau(S; \O)$ for some $r, s\in \Z \cup\{-\infty, \infty\}$
with $r \leq s$ and $\tau\in \{\proet, v\}$. Then the $\tau$-sheaf $\Latt^{[r, s]}_{\tau, \cal{P}}(\blank) \colon \Perfd^{\aff, \op}_{/S} \to \Cat_\infty$ is finitary.
\end{lemma}
\begin{proof}
    Since sheafification and filtered colimits in the presheaf category commute with finite limits (see \cref{cor:finite-limits-commute-with-sheafifications} and \cref{lemma:finite-limits-vs-filtered-colimits}), \cref{lemma:colimit-of-sheaves} and \cref{rmk:fully-faithful-finite-tor-amplitude}~\ref{rmk:fully-faithful-finite-tor-amplitude-1} imply that it suffices to consider the case when $r$ and $s$ are (finite) integers. Furthermore, \cref{lemma:finitariness-is-v-local} implies that we can check that $\Latt^{[r, s]}_{\tau, \cal{P}}(\blank)$ is finitary 
    $\tau$-locally on $S$. Therefore, we can assume that $\cal{P} \simeq \widetilde{P}$ for some $P\in \Perf^{[r, s]}\big(\O_S(S)\big)$. For brevity, we denote by $[P]\colon * \to \Perf^{[r, s]}\big(\O(\blank)\big)$ and $[\cal{P}]\colon * \to \Perf_\tau^{[r, s]}\big(\blank; \O\big)$ the corresponding morphisms in $\PShv(\Perfd^\aff_{/S}, \Cat_{(s-r+1, 1)})$. 
    
    Since the morphism $[\cal{P}]$ factors through the fully faithful morphism $\Perf^{[r, s]}\big(\O_S(S)\big) \simeq \Perf^{[r, s]}_\et(S; \O) \hookrightarrow \Perf^{[r, s]}_\tau(S; \O)$ (see \cref{rmk:fully-faithful-finite-tor-amplitude}~\cref{rmk:fully-faithful-finite-tor-amplitude-2}),  \cref{cor:integral-descend-finite-tor-amplitude} implies that 
    \[
    \Latt^{[r, s]}_{\tau, \cal{P}}(\blank) \simeq * \times_{\Perf^{[r, s]}_\tau(\blank; \O)} \Perf^{[r, s]}_\tau(\blank; \O^+) \simeq * \times_{\Perf^{[r, s]}_\et(\blank; \O)} \Perf^{[r, s]}_\et(\blank; \O^+).
    \]
    Therefore, \cref{cor:finite-limits-commute-with-sheafifications} and \cref{thm:sheafification-is-finitary} imply that it suffices to show that the presheaf
    \[
    * \times_{\Perf^{[r, s]}\big(\O(\blank) \big)} \Perf^{[r, s]}\big(\O^+(\blank) \big) = \rm{fib}_{[P]}\Big(\Perf^{[r, s]}\big(\O^+(\blank)\big) \to  \Perf^{[r, s]}\big(\O(\blank)\big) \Big)
    \]
    is finitary. For this, we choose a pseudo-uniformizer $\varpi \in \O_S^+(S)$ and fix $\Spa(R, R^+) = \lim_I \Spa(R_i, R_i^+)$ in $\Perfd^\aff_{/S}$. Then \cite[Prop.~6.4]{diamonds} implies that $R^+ \simeq (\colim_{I^\op} R_i^+\big)^{\wedge}_{(\varpi)}$ and $R \simeq R^+[\frac{1}{\varpi}]$. We also put $R_\infty^+ \coloneqq \colim_{I^\op} R_i^+$ and $R_\infty \coloneqq R_\infty^+[\frac{1}{\varpi}]$. In this situation, we have to show that 
    \[
    \colim_{I^\op} \rm{fib}_{[P\otimes^\rL_A R_i]} \big(\Perf^{[r, s]}(R_i^+) \to  \Perf^{[r, s]}(R_i) \big) \to \rm{fib}_{[P\otimes^\rL_A R]} \big(\Perf^{[r, s]}(R^+) \to  \Perf^{[r, s]}(R) \big)
    \]
    is an equivalence. To see this, we recall that the Beauville--Laszlo gluing (see \cite[Prop.~5.6(2)]{Bhatt-Tannaka} or \cite[Prop.~7.4.2.1]{SAG}) implies that 
    \[
    \rm{fib}_{[P\otimes^\rL_A R_\infty]} \big(\Perf^{[r, s]}(R_\infty^+) \to  \Perf^{[r, s]}(R_\infty) \big) \simeq \rm{fib}_{[P\otimes^\rL_A R]} \big(\Perf^{[r, s]}(R^+) \to  \Perf^{[r, s]}(R) \big).
    \]
    Therefore, we reduce the question to showing that 
    \[
    \colim_{I^\op} \rm{fib}_{[P\otimes^\rL_A R_i]} \big(\Perf^{[r, s]}(R_i^+) \to  \Perf^{[r, s]}(R_i) \big) \to \rm{fib}_{[P\otimes^\rL_A R_\infty]} \big(\Perf^{[r, s]}(R^+_\infty) \to  \Perf^{[r, s]}(R_\infty) \big)
    \]
    is an equivalence. This follows from \cref{lemma:perfect-complexes-are-finitary} and \cref{lemma:finite-limits-vs-filtered-colimits}. 
\end{proof}

Now we are ready to show that any perfect complex of $\O_v$-modules admits a lattice \'etale locally.

\begin{lemma}\label{lemma:lattice-etale-locally} Let $S=\Spa(A, A^+)$ be an affinoid perfectoid, let $\tau\in \{\proet, v\}$, and let $\cal{P}$ be an element of $\Perf^{[r, s]}_\tau(S; \O)$ for some $r,s\in \Z \cup\{-\infty, \infty\}$. Then there is an \'etale covering $S' \to S$ and an $\cal{P}^+\in \Perf^{[r,s]}_\tau(S'; \O^+)$ such that $\cal{P}^+\otimes^\rL_{\O_{S'_\tau}^+} \O_{S'_\tau} \simeq \cal{P}|_{S'}$.
\end{lemma} 

\begin{proof}
    It suffices to show that $\Latt^{[r, s]}_{\tau, \cal{P}}(S')\neq \varnothing$ for some \'etale covering $S' \to S$. Now \cref{lemma:lattices-are-finitary} and \cref{lemma:sections-etale-locally} imply that it suffices to show that there is a $\tau$-covering $S' \to S$ such that $\Latt^{[r, s]}_{\cal{P}}(S')\neq \varnothing$. So we can argue $\tau$-locally on $S$. Therefore, we can assume that $S$ is strictly totally disconnected and that we can find an object $P\in \Perf^{[r, s]}(\O_S(S))$ such that $\cal{P} \simeq \widetilde{P}$. Now it suffices to show that, for every point $x\in S$, there is an \'etale map $f_x \colon S_x \to S$ such that $x\in f_x(S_x)$ and $\Latt_{\tau, \cP}^{[r,s]}(S_x)\not\simeq \varnothing$. 
    
    For this, we set $C_x\coloneqq \wdh{k(x)}$ and $C^+_x\coloneqq \wdh{k(x)^+}$. Since $S$ is strictly totally disconnected, we conclude that $C_x$ is an algebraically closed nonarchimedean field. Furthermore, $\Spa(C_x, C_x^+) \simeq \lim_{x\in U\subset S} U$, where the limit is taken over
    (the opposite category of) rational subdomains of $S$ containing $x$. Therefore, \cref{lemma:lattices-are-finitary} ensures that it suffices to show that $\Latt_{\cP, \tau}^{[r,s]}(\Spa(C_x, C_x^+))\not\simeq \varnothing$. In other words, we reduce the question to the case $S=\Spa(C, C^+)$ and $\cP  =\widetilde{P}$ for some $P\in \Perf^{[r,s]}(\O_S(S))$.

    Since $C=\O_S(S)$ is a field, we note that $P$ has a representative $K^\bullet$ such that $K$ is a finite complex, $K^i = 0$ if $i\notin [r, s]$, and each $K^i$ is a finite free $C$-module. Then we note that after multiplying each differential by a large power of a pseudo-uniformizer, we can assume that there is a finite complex $K^{+, \bullet}$ of finite free $C^+=\O^+_S(S)$-modules such that $K^{+, \bullet} \otimes^\rL_{C^+} C \simeq K^\bullet$ and such that $K^{+, i}=0$ for $i\notin [r, s]$. Then $\widetilde{K^{+, \bullet}}$ is the desired lattice of $\cal{P}$.
\end{proof}

\begin{thm}\label{thm:v-perfect-complexes-vs-etale-perfect-complexes-rationally} Let $X$ be a perfectoid space, let $\tau\in \{\proet, v\}$, and let $\lambda_\tau \colon (X_\tau, \O) \to (X_\et,  \O)$ be the natural morphism of ringed sites. Then the following hold:
\begin{enumerate}[leftmargin=*,label={\upshape{(\roman*)}}]
    \item\label{thm:v-perfect-complexes-vs-etale-perfect-complexes-rationally-1} the functor $\rR\lambda_{\tau, *} \colon \cal{D}(X_\tau; \O) \to \cal{D}(X_\et; \O)$ restricts to a functor $\rR\lambda_{\tau, *}\colon \Perf_\tau(X; \O) \to \Perf_\et(X; \O)$;
    \item\label{thm:v-perfect-complexes-vs-etale-perfect-complexes-rationally-2} the functor $\lambda_\tau^* \colon \Perf_\et(X; \O) \to \Perf_\tau(X; \O)$ is an equivalence with the inverse functor $\rR\lambda_{\tau, *}$.
\end{enumerate}
\end{thm}
\begin{proof}
    We start by proving part~\cref{thm:v-perfect-complexes-vs-etale-perfect-complexes-rationally-1}. We fix a $\cal{P}\in \Perf_\tau(X; \O)$ and wish to show that $\rR\lambda_{\tau, *}\cal{P}$ lies in  $\Perf_\et(X; \O)$. The question is \'etale local on $X$, so we can assume that $X=\Spa(R, R^+)$ is an affinoid perfectoid with a pseudo-uniformizer $\varpi\in R^+$. Then \cref{lemma:lattice-etale-locally} ensures that we can replace $X$ with an \'etale covering to assume that $\cal{P}$ admits a $\O_{X_\tau}^+$-lattice $\cal{P}^+\in \Perf_\tau(X; \O^+)$. Then \cref{thm:v-perfect-complexes-vs-etale-perfect-complexes-integrally}~\cref{thm:v-perfect-complexes-vs-etale-perfect-complexes-integrally-1} implies that $\rR\lambda_{\tau, *} \cal{P}^+\in \Perf_\et(X; \O^+)$. Since $\lambda_{\tau}$ is a coherent morphism of coherent topoi, we conclude that $\rR\lambda_{\tau, *} \cal{P} \simeq \big(\rR\lambda_{\tau, *}\cal{P}^+\big)\big[\frac{1}{\varpi}\big] \in \Perf_\et(X; \O)$. 

    Now we show Part~\cref{thm:v-perfect-complexes-vs-etale-perfect-complexes-rationally-2}. \cref{prop:fully-faithful-rationally} and Part~\cref{thm:v-perfect-complexes-vs-etale-perfect-complexes-rationally-1} imply that it suffices to show that, for any $\cal{P}\in \Perf_\tau(X; \O)$, the natural map $\lambda_\tau^* \rR\lambda_{\tau, *}\cal{P} \to \cal{P}$ is an isomorphism. Using that $\cal{P}$ admits a lattice \'etale locally, we also reduce this claim to \cref{thm:v-perfect-complexes-vs-etale-perfect-complexes-integrally}~\cref{thm:v-perfect-complexes-vs-etale-perfect-complexes-integrally-2}. 
\end{proof}

\begin{rmk}\label{rmk:comparison-of-topologies-O} If we combine \cref{lemma:etale-descent-rationally} with \cref{thm:v-perfect-complexes-vs-etale-perfect-complexes-rationally}, we conclude that the natural pullback functor $\Perf_\an(X; \O) \to \Perf_\tau(X; \O)$ is an equivalence for any perfectoid space $X$ and any $\tau\in \{\et, \proet, v\}$.
\end{rmk}

Now we wish to get a version of \cref{cor:integral-descend-finite-tor-amplitude} for $\Perf\bigl(\O(\blank)\bigr)$. For this, we will need one preliminary lemma. 

\begin{lemma}\label{lemma:surjectivity-on-closed-points-rationally} Let $f\colon (R, R^+) \to (S, S^+)$ be a continuous morphism of complete Huber pairs. Suppose that the induced morphism $\Spa(S, S^+) \to \Spa(R, R^+)$ is surjective. Then the induced morphism $\Spec S \to \Spec R$ is surjective on closed points. 
\end{lemma}
\begin{proof}
    Using functoriality of the support map $\rm{supp} \colon \abs{\Spa(S, S^+)} \to \abs{\Spec S}$, we see that it suffices to show that, for each maximal ideal $\m\subset R$, there is a valuation $v\in \Spa(R, R^+)$ such that $v^{-1}(\{0\}) = \m$. This follows immediately from \cite[Lem.~1.4]{Huber-generalization}. 
\end{proof}

\begin{thm}\label{thm:v-descent-rationally} Let $r, s\in \Z \cup \{-\infty, \infty\}$. Then the presheaf $\Perf^{[r, s]}\big(\O(\blank)\big) \colon \Perfd^{\aff, \op} \to \Cat_\infty$ is a hypercomplete $v$-sheaf of $\infty$-categories. In particular, the presheaf $\Perf\big(\O(\blank)\big)$ is a hypercomplete $v$-sheaf. 
\end{thm}
\begin{proof}
    First, we treat the case $r=-\infty$ and $s=\infty$. In this case, \cref{lemma:etale-descent-rationally} implies that $\Perf\big(\O(\blank)\big)$ is an \'etale sheaf. Thus, \cref{prop:perfect-complexes-sheafification} ensures that $\Perf_\et(\blank; \O) \simeq \Perf\bigl(\O(\blank)\bigr)$. Then \cref{thm:v-perfect-complexes-vs-etale-perfect-complexes-rationally} and \cref{lemma:perfect-complexes-are-hypercomplete} ensure that 
    \[
    \Perf\bigl(\O(\blank)\bigr) \simeq \Perf_\et(\blank; \O) \simeq \Perf_v(\blank; \O)
    \] 
    is a hypercomplete $v$-sheaf. In order to deduce the result for any $r, s\in \Z\cup\{-\infty, \infty\}$, one argues as in the proof of \cref{cor:integral-descend-finite-tor-amplitude} using \cref{lemma:surjectivity-on-closed-points-rationally} in place of \cref{lemma:surjective-integrally} (which is used in the proof of \cref{cor:test-tor-amplitude-v-locally}). 
\end{proof}

\begin{rmk} For $r= s=0$, \cref{thm:v-descent-rationally} recovers the fact that $\Vect\big(\O(\blank)\big)$ is a $v$-stack of (ordinary) categories (see \cite[Th.~3.5.8]{Kedlaya-Liu-2}, \cite[Lem.~17.1.8]{Berkeley}, \cite[Th.~2.20]{Heuer-G-torsor}, and \cite[Th.~1.3.3]{zav-almost}).  For $r=-\infty$ and $s=\infty$, this provides a stronger version of \cite[Th.\,2.1]{ALB} which applies to all affinoid perfectoid spaces and not only those that admit a map to a totally disconnected perfectoid space. 
\end{rmk}

\begin{cor}
\label{cor:perfect-O-modules-on-affinoid-perfectoid} 
Let $X=\Spa(R, R^+)$ be an affinoid perfectoid space and let $\tau\in \{\an, \et, \proet, v\}$. Then the functor $\widetilde{(\blank)}\colon \Perf(R) \to \Perf_\tau(X; \O)$ is an equivalence with the inverse given by $\rR\Gamma_\tau(X, \blank)\colon \Perf_\tau(X, \O) \to \Perf(R)$.
\end{cor}
\begin{proof}
    \cref{rmk:comparison-of-topologies-O} implies that it suffices to prove the claim for $\tau=\an$. \cref{thm:v-descent-rationally} together with \cref{prop:perfect-complexes-sheafification} imply that the natural functor $\widetilde{(\blank)} \colon \Perf(R) \to \Perf_\an(X; \O)$ is an equivalence for any affinoid perfectoid $X=\Spa(R, R^+)$. To finish the proof, we only need to show that the natural morphism
    \[
    P \to    \rR\Gamma_\an(X, \widetilde{P})
    \]
    is an isomorphism for any $P\in \Perf(R)$. An easy induction step reduces the question to the case when $P$ is a finite projective $R$-module in degree $0$. Since any such $P$ is a direct summand of a finite free module, we further reduce to showing that the natural morphism $R \to \rR\Gamma_\an(X, \O)$ is an isomorphism. This follows immediately from \cite[Th.\,3.24]{diamonds}.
\end{proof}

\subsection{$v$-descent for perfect $\BB_I$-complexes}\label{section:descent-BI}

For the rest of this subsection, we fix a closed interval $I=[a,b] \subset (0, \infty)$ with rational endpoints. We refer to \cref{defn:period-rings} for the definition of $\BB_I(X)$ for an affinoid perfectoid $X$ over $\Spa(\Q_p, \Z_p)$, we also refer to \cref{cor:AI-etale-sheaf} for the proof that $\BB_I$ is a $v$-sheaf on $\Perfd^{\aff}_{/\Q_p}$. Let 
\[
\Perf\bigl(\BB_{I}(\blank)\bigr) \colon \Perfd^{\aff, \op}_{/\Q_p} \to \Cat_\infty
\]
be the presheaf of $\infty$-categories that sends an affinoid perfectoid space $X$ to the $\infty$-category of perfect complexes of $\BB_{I}(X)$-modules. The main goal of this subsection is to show that this is a hypercomplete 
$v$-sheaf of $\infty$-categories. 

We will closely follow the strategy taken in \cref{section:descent-O}. We will separately show \'etale descent and then 
use the sheaf of lattices to reduce $v$-descent to the case of perfect $\AA_I$-complexes. 

\begin{lemma}\label{lemma:etale-descent-BBI} The presheaves  $\Perf^{[r, s]}\big(\BB_I(\blank)\big) \colon \Perfd^{\aff, \op}_{/\QQ_p} \to \Cat_\infty$ are \'etale sheaves of $\infty$-categories for any $r, s\in \ZZ \cup \{-\infty, \infty\}$
with $r \leq s$. In particular, $\Perf\big(\BB_I(\blank)\big)$ is an \'etale sheaf. 
\end{lemma}
\begin{proof}
    Let $F$ be any of the presheaves in the formulation of the lemma. Since the functor $Y_{\blank, I}$ commutes with finite disjoint unions, $F(\sqcup_{i=1}^n(X_i)) \simeq \prod_{i=1}^n F(X_i)$ for any $X_1, \dots, X_n\in \Perfd^\aff_{\QQ_p}$. Therefore,  \cite[Prop.~A.3.3.1]{SAG} implies that it suffices to show that $F$ satisfies descent with respect to any \'etale covering $X \to Y$. In other words, we have to show that any \'etale covering $f\colon X \to Y$ is of universal $F$-descent in the sense of \cite[Def.~3.3.1]{Liu-Zheng}. Then \cite[Lem.~3.3.2]{Liu-Zheng} and \cite[Prop.~8.2.20]{Kedlaya-Liu-1} implies that it suffices to show that $F$ satisfies analytic and finite \'etale descent. Analytic descent follows from \cite[Th.~5.3]{Andreychev} together with the fact that functor $Y_{\blank, I}$ preserves open covers (see \cref{lemma:permanence-properties-Y}). Now let $X \to Y$ be a finite \'etale surjection in $\Perfd^\aff_{/\QQ_p}$, then \cref{lemma:finite-etale-after-v-completion}, \cref{cor:Fargues--Fontaine-v-complete}, and \cite[Lem.~C.1.2]{zav-almost} imply that $\Spec \BB_I(X) \to \Spec \BB_I(Y)$ is a finite \'etale surjection. Furthermore, the combination of \cite[Lem.~C.1.1]{zav-almost} and \cite[Lem.~B.3.5]{Z-quotients} imply that $\BB_I(X)^{\wdh{\otimes}^n_{\BB_I(Y)}} \simeq \BB_I(X)^{\otimes^n_{\BB_I(Y)}}$. Therefore, the result follows from the usual fppf descent for perfect complexes of fixed tor-amplitude by \cite[Cor.~D.6.3.3 \& Prop.~2.8.4.2]{SAG}. 
\end{proof}

Now we wish to deduce that $\Perf\big(\BB_I(\blank)\big)$ is a (hypercomplete) $v$-sheaf. We do this by comparing this functor to $\Perf_v(\blank; \BB_I)$. We start with the following proposition:

\begin{proposition}\label{prop:fully-faithful-rationally-BBI} Let $X$ be a perfectoid space over $\Spa(\QQ_p, \ZZ_p)$, let $\cal{P}, \cal{Q}\in \Perf_\et(X; \BB_I)$, let $\tau\in \{\proet, v\}$, and let $\lambda_\tau \colon (X_\tau, \BB_I) \to (X_\et,  \BB_I)$ be the natural morphism of ringed sites. Then the following statements hold:
\begin{enumerate}[leftmargin=*,label={\upshape{(\roman*)}}]
    \item\label{prop:fully-faithful-rationally-BBI-1} the natural morphism $\cal{P} \to \rR\lambda_{\tau, *} \lambda^*_\tau\cal{P}$ is an isomorphism;
    \item\label{prop:fully-faithful-rationally-BBI-2} the natural morphism $\rR\Hom(\cal{P}, \cal{Q}) \to \rR\Hom(\lambda_\tau^*\cal{P}, \lambda_\tau^*\cal{Q})$ is an isomorphism;
    \item\label{prop:fully-faithful-rationally-BBI-3} the functor $\lambda_\tau^*\colon \Perf_\et(X; \BB_I) \to \Perf_\tau(X; \BB_I)$ is fully faithful.
\end{enumerate}
\end{proposition}
\begin{proof}
    The proof is identical to the proof of \cref{prop:fully-faithful-rationally}. The only difference is that we need to use \cref{thm:v-perfect-complexes-vs-etale-perfect-complexes-integrally-A_I} in place of \cref{prop:fully-faithful-integrally}. 
\end{proof}

Now we wish to show that $\lambda_\tau^*$ is also essentially surjective. We do this by studying the sheaf of lattices. For this, we will use the sheaf $\Perf^{[r,s]}_\tau(\blank; \BB_I) \colon \Perfd^{\aff, \op}_{/\QQ_p} \to \Cat_\infty$ from \Cref{notation:sheaf-perfect-complexes-finite-amplitude-cA}. 

First, we have the following series of easy remarks.

\begin{rmk}
\label{rmk:fully-faithful-finite-tor-amplitude-BBI} 
Similar to \cref{rmk:fully-faithful-finite-tor-amplitude}, the functors $\Perf^{[r, s]}_\tau(\blank; \BB_I)$ satisfy the following properties: 
\begin{enumerate}
    \item\label{rmk:fully-faithful-finite-tor-amplitude-BBI-1}  The functor $\Perf^{[r, s]}_\tau(\blank; \BB_I) \to \Perf_\tau^{[r', s']}(\blank; \BB_I)$ is fully faithful for $\tau\in \{\an, \et, \proet, v\}$ and $r,r', s, s'\in \Z \cup \{-\infty, \infty\}$ with $r'\leq r \leq s \leq s'$. The natural functor $\colim_\NN \Perf^{[-n,n]}_\tau(X; \BB_I) \to \Perf_\tau(X; \BB_I)$ is an equivalence for any $X\in \Perfd^\aff_{/\QQ_p}$;
    \item\label{rmk:fully-faithful-finite-tor-amplitude-BBI-2}  The functors $\Perf_\et(\blank; \BB_I) \to \Perf_\tau(\blank; \BB_I)$ and $\Perf^{[r, s]}_\et(\blank; \BB_I) \to \Perf^{[r, s]}_\tau(\blank; \BB_I)$ are fully faithful for any $r, s\in \ZZ\cup \{-\infty, \infty\}$ with $r \leq s$ and $\tau\in \{\proet, v\}$;
    \item\label{rmk:fully-faithful-finite-tor-amplitude-BBI-3}  We have an equivalence $\Perf^{[r,s]}\big(\BB_I(\blank)\big) \simeq \Perf^{[r, s]}_\et(\blank; \BB_I)$ for any $r, s \in \Z \cup \{-\infty, \infty\}$
    with $r \leq s$.
\end{enumerate}
\end{rmk}

Fix an affinoid perfectoid space $S$ over $\Spa(\QQ_p, \ZZ_p)$. Then we recall that a perfect complex $\cal{P} \in \Perf^{[r,s]}_\tau(S; \BB_I)$ defines a morphism of sheaves 
\[
[\cal{P}] \colon * \to \Perf_\tau^{[r, s]}(\blank; \BB_I) \colon \Perfd^{\aff, \op}_{/S} \to \Cat_\infty.
\]

Using this observation, we define the sheaf of lattices as follows: 

\begin{defn} The {\it $\tau$-sheaf of lattices} $\Latt_{\tau, \cal{P}}^{[r, s]}(\blank) \colon \Perfd^{\aff, \op}_{/S} \to \Cat_\infty$ is the fiber product
\[
\Latt^{[r,s]}_{\tau, \cal{P}}(\blank) \coloneqq * \times_{\Perf^{[r, s]}_\tau(-; \BB_I)} \Perf^{[r, s]}_\tau(-; \AA_I).
\]
\end{defn}

More informally, for every affinoid perfectoid $S$-space $X$, we have
\[
\Latt^{[r, s]}_{\tau, \cal{P}}(X)= \Big\{ \cal{P}^+ \in \Perf^{[r, s]}_\tau(X; \AA_I) \text{ and } \varphi \colon \cal{P}^+ \otimes^{\rL}_{\AA_I} \BB_I \xr{\sim} \cal{P} \Big\}.
\]

Below, we will need the following preliminary lemma.

\begin{lemma}\label{lemma:completion-colimit-AAI} Let $S$ be an affinoid perfectoid space over $\Spa(\QQ_p, \ZZ_p)$, let $p^\flat\in \O^{\flat, +}_{S}(S)$ be an element from \cref{notation:p-flat}, and let $\alpha, \beta\in \O^{\flat, +}_S(S)$ be element well-adapted to $p^\flat$ and $I$ (in the sense of \cref{defn:well-adapted}). Let $X=\lim_J X_j$ be a cofiltered limit in $\Perfd^{\aff}_{/S}$. Then $(\colim_{J^\op} \AA_I(X_j)\bigr)^{\wedge}_{(p)} \to \AA_I(X)$ is an isomorphism. 
\end{lemma}
\begin{proof}
        \cref{cor:almost-mathematics-AI} implies that $(\colim_{J^\op} \AA_I(X_j)\bigr)^{\wedge}_p \to \AA_I(X)$ is a homomorphism between $p$-adically complete $p$-torsionfree rings. Therefore, it suffices to show that this morphism is an isomorphism modulo $p$. Then \cref{cor:AI-mod-p-rings} and the fact that colimits commute with direct sums imply that it suffices to show that the natural map $\colim_{J^\op} \O_{X_j}^{\flat, +}(X_j)/\alpha \to \O_{X}^{\flat, +}(X)/\alpha$ is an isomorphism (and the similar statement for modulo $\beta$). Now a simple inductive argument using \cref{rmk:well-adapted-pseudo-uniformizer} and the fact that colimits commute with cokernels ensures that it suffices to show that the natural map $\colim_{J^\op} \O_{X_j}^{\flat, +}(X_j)/p^\flat \to \O_{X}^{\flat, +}(X)/p^\flat$ is an isomorphism. This is equivalent to showing that the natural morphism
        \[
        \colim_{J^\op} \O_{X_j}^{+}(X_j)/p \to \O_{X}^{+}(X)/p
        \]
        is an isomorphism. For this, we write $X=\Spa(R, R^+)$ and $X_j=\Spa(R_j, R_j^+)$. Then \cite[Prop.~6.4]{diamonds} implies that $R^+ \simeq (\colim_{J^\op} R_j^+)^{\wedge}_{(p)}$. Therefore, the result follows from the fact that the natural morphism $(\colim_{J^\op} R_j^+)/p \to (\colim_{J^\op} R_j^+)^{\wedge}_{(p)}/p$ is an isomorphism (see \cite[\href{https://stacks.math.columbia.edu/tag/05GG}{Tag 05GG}]{stacks-project}).
\end{proof}

\begin{lemma}\label{lemma:lattices-are-finitary-BBI} Let $S=\Spa(A, A^+)$ be an affinoid perfectoid over $\Spa(\QQ_p, \ZZ_p)$ and let $\cal{P}$ be an element of $\Perf^{[r, s]}_\tau(S; \BB_I)$ for some $r,s\in \Z \cup\{-\infty, \infty\}$
with $r \leq s$ and $\tau\in \{\proet, v\}$. Then the $\tau$-sheaf $\Latt^{[r, s]}_{\tau, \cal{P}}(\blank) \colon \Perfd^{\aff, \op}_{/S} \to \Cat_\infty$ is finitary.
\end{lemma}
\begin{proof}
    Arguing as in the proof of \cref{lemma:lattices-are-finitary}, we first reduce to the case when $r$ and $s$ are (finite) integers. Furthermore, \cref{lemma:finitariness-is-v-local} implies that we can check that $\Latt^{[r, s]}_{\tau, \cal{P}}(\blank)$ is finitary $\tau$-locally on $S$. Therefore, we can assume that $S$ is strictly totally disconnected and that $\cal{P} \simeq \widetilde{P}$ for some $P\in \Perf^{[r,s]}\big(\BB_I(S)\big)$. For brevity, we denote by $[P]\colon * \to \Perf^{[r, s]}\big(\BB_I(\blank)\big)$ and $[\cal{P}]\colon * \to \Perf_\tau^{[r, s]}\big(\blank; \BB_I\big)$ the corresponding morphisms in $\PShv(\Perfd^\aff_{/S}, \Cat_{(s-r+1, 1)})$. 
    
    Since the morphism $[\cal{P}]$ factors through the fully faithful morphism $\Perf^{[r, s]}\big(\BB_I(S)\big) \simeq \Perf^{[r,s]}_\et(S; \BB_I) \hookrightarrow \Perf^{[r,s]}_\tau(S; \BB_I)$ (see \cref{rmk:fully-faithful-finite-tor-amplitude-BBI}~\cref{rmk:fully-faithful-finite-tor-amplitude-BBI-2}),  \cref{cor:v-descent-A_I} implies that 
    \[
    \Latt^{[r, s]}_{\tau, \cal{P}}(\blank) \simeq * \times_{\Perf^{[r, s]}_\tau(\blank; \BB_I)} \Perf^{[r, s]}_\tau(\blank; \AA_I) \simeq * \times_{\Perf^{[r, s]}_\et(\blank; \BB_I)} \Perf^{[r, s]}_\et(\blank; \AA_I).
    \]
    Therefore, \cref{cor:finite-limits-commute-with-sheafifications} and \cref{thm:sheafification-is-finitary} imply that it suffices to show that the presheaf
    \[
    * \times_{\Perf^{[r, s]}\big(\BB_I(\blank) \big)} \Perf^{[r, s]}\big(\AA_I(\blank) \big) = \rm{fib}_{[P]}\Big(\Perf^{[r, s]}\big(\AA_I(\blank)\big) \to  \Perf^{[r, s]}\big(\BB_I(\blank)\big) \Big)
    \]
    is finitary. For this, we fix a cofiltered limit $X=\lim_J X_j$ in $\Perfd^\aff_{/S}$. We need to show that 
    \[
    \colim_{J^\op} \rm{fib}_{[P\otimes^\rL_{\BB_I(S)} \BB_I(X_j)]} \Bigl(\Perf^{[r, s]}\bigl(\AA_I(X_j)\bigr) \to  \Perf^{[r, s]}\bigl(\BB_I(X_j)\bigr) \Bigr) \to \rm{fib}_{[P\otimes^\rL_{\BB_I(S)} \BB_I(X)]} \Bigl(\Perf^{[r, s]}\bigl(\AA_I(X)\bigr) \to  \Perf^{[r, s]}\bigl(\BB_I(X)\bigr) \Bigr)
    \]
    is an equivalence. To see this, we recall that the Beauville--Laszlo gluing (see \cite[Prop.~5.6(2)]{Bhatt-Tannaka} or \cite[Prop.~7.4.2.1]{SAG}) together with \cref{lemma:completion-colimit-AAI} and \cref{rmk:often-well-adapted-pairs} imply that 
    \begin{align*}
    \rm{fib}_{[P\otimes^\rL_{\BB_I(S)} \colim_{J^\op} \BB_I(X_j)]} \Big(\Perf^{[r, s]}\big(\colim_{J^\op} \AA_I(X_j)\big)& \to  \Perf^{[r, s]}\big(\colim_{J^\op} \BB_I(X_j)\big) \Big) \simeq  \\
    & \rm{fib}_{[P\otimes^\rL_{\BB_I(S)} \BB_I(X)]} \Big(\Perf^{[r, s]}\big(\AA_I(X)\big) \to  \Perf^{[r, s]}\big(\BB_I(X)\big) \Big).
    \end{align*}
    Therefore, the desired result follows from \cref{lemma:perfect-complexes-are-finitary} and \cref{lemma:finite-limits-vs-filtered-colimits}. 
\end{proof}

Now we are ready to show that any perfect complex of $\BB_{I,v}$-modules admits a lattice \'etale locally.

\begin{lemma}\label{lemma:lattice-etale-locally-BBI} Let $S$ be an affinoid perfectoid over $\Spa(\QQ_p, \ZZ_p)$, let $\tau\in \{\proet, v\}$, and let $\cal{P}$ be an element of $\Perf^{[r, s]}_\tau(S; \BB_I)$ for some $r, s\in \Z \cup\{-\infty, \infty\}$. Then there is an \'etale covering $S' \to S$ and an $\cal{P}^+\in \Perf^{[r,s]}_\tau(S'; \AA_I)$ such that $\cal{P}^+\otimes^\rL_{\AA_I} \BB_I \simeq \cal{P}|_{S'}$.
\end{lemma}

\begin{proof}
    Arguing as in the proof of \cref{lemma:lattice-etale-locally} (and using \cref{lemma:lattices-are-finitary-BBI} in place of \cref{lemma:lattices-are-finitary}), we can reduce to the case when $S=\Spa(C, C^+)$ for an algebraically closed nonarchimedean field $C$ and an open bounded valuation subring $C^+\subset C$ and $\cP\simeq \widetilde{P}$ for some $P\in \Perf^{[r,s]}(\BB_I(S))$.  

    In this case, \cref{cor:BI-PID} ensures that $B_I\coloneqq \BB_I(S)$ is a principal ideal domain. In particular, any finite projective $B_I$-module is free. Therefore, we conclude that $P$ has a representative $K^\bullet$ such that $K$ is a finite complex, $K^i = 0$ if $i\notin [r, s]$, and each $K^i$ is a finite free $B_I$-module. Since $A_I\coloneqq \AA_I(S)$ is $p$-torsionfree and $B_I\simeq A_I[\frac{1}{p}]$, we note that after multiplying each differential by a large power of $p$, we can assume that there is a finite complex $K^{+, \bullet}$ of finite free $A_I$-modules such that $K^{+, \bullet} \otimes^\rL_{A_I} B_I \simeq K^\bullet$ and such that $K^{+, i}=0$ for $i\notin [r, s]$. Then $\widetilde{K^{+, \bullet}}$ is the desired lattice of $\cal{P}$. 
\end{proof}

\begin{thm}\label{thm:v-perfect-complexes-vs-etale-perfect-complexes-rationally-BBI} Let $X$ be a perfectoid space over $\Spa(\QQ_p, \ZZ_p)$, let $\tau\in \{\proet, v\}$, and let $\lambda_\tau \colon (X_\tau, \BB_I) \to (X_\et,  \BB_I)$ be the natural morphism of ringed sites. Then the following hold:
\begin{enumerate}[leftmargin=*,label={\upshape{(\roman*)}}]
    \item\label{thm:v-perfect-complexes-vs-etale-perfect-complexes-rationally-BBI-1} the functor $\rR\lambda_{\tau, *} \colon \cal{D}(X_\tau; \BB_I) \to \cal{D}(X_\et; \BB_I)$ restricts to a functor $\rR\lambda_{\tau, *}\colon \Perf_\tau(X; \BB_I) \to \Perf_\et(X; \BB_I)$;
    \item\label{thm:v-perfect-complexes-vs-etale-perfect-complexes-rationally-BBI-2} the functor $\lambda_\tau^* \colon \Perf_\et(X; \BB_I) \to \Perf_\tau(X; \BB_I)$ is an equivalence with the inverse functor $\rR\lambda_{\tau, *}$.
\end{enumerate}
\end{thm}
\begin{proof}
    The proof of \cref{thm:v-perfect-complexes-vs-etale-perfect-complexes-rationally}
    applies verbatim. The only difference is that one needs to use \cref{lemma:lattice-etale-locally-BBI}, \cref{prop:fully-faithful-rationally-BBI}, and \cref{thm:v-perfect-complexes-vs-etale-perfect-complexes-integrally-A_I} in place of \cref{lemma:lattice-etale-locally}, \cref{prop:fully-faithful-rationally}, and \cref{thm:v-perfect-complexes-vs-etale-perfect-complexes-integrally} respectively. 
\end{proof}

\begin{lemma}\label{lemma:surjective-BI} Let $X \to Y$ be a surjective morphism of affinoid perfectoid spaces over $\Spa(\Q_p, \Z_p)$. Then the natural morphism $\Spec \BB_I(X) \to \Spec \BB_I(Y)$ is surjective on closed points. 
\end{lemma}
\begin{proof}
    This follows from \cref{lemma:permanence-properties-Y} and \cref{lemma:surjectivity-on-closed-points-rationally}.    
\end{proof}

\begin{thm}\label{thm:v-descent-rationally-BBI} Let $r, s\in \Z \cup \{-\infty, \infty\}$. Then the presheaf $\Perf^{[r,s]}\big(\BB_I(\blank)\big) \colon \Perfd^{\aff, \op}_{/\QQ_p} \to \Cat_\infty$ is a hypercomplete $v$-sheaf of $\infty$-categories. In particular, the presheaf $\Perf\big(\BB_I(\blank)\big)$ is a hypercomplete $v$-sheaf. 
\end{thm}
\begin{proof}
    The proof of \cref{thm:v-descent-rationally} applies verbatim. The only difference is that one needs to use \cref{thm:v-perfect-complexes-vs-etale-perfect-complexes-rationally-BBI} and \cref{lemma:surjective-BI} in place of \cref{thm:v-perfect-complexes-vs-etale-perfect-complexes-rationally} and \cref{lemma:surjectivity-on-closed-points-rationally} respectively. 
\end{proof}

\begin{rmk} For $r=s=0$, \cref{thm:v-descent-rationally-BBI} recovers the fact that $\Vect\big(\BB_I(\blank)\big)$ is a $v$-stack of (ordinary) categories (see \cite[Prop.\,II.2.1]{Fargues-Scholze}). For $r=-\infty$ and $s=\infty$, this recovers \cite[Prop.\,2.3]{ALB} for standard opens $Y_{S, I}\subset Y_S$.
\end{rmk}

\begin{cor}
\label{cor:perfect-O-modules-on-affinoid-perfectoid-BBI}
Let $X$ be an affinoid perfectoid space over $\Spa(\QQ_p, \ZZ_p)$ and let $\tau\in \{\an, \et, \proet, v\}$. Then the functor $\widetilde{(\blank)}\colon \Perf\bigl(\BB_I(X)\bigr) \to \Perf_\tau(X; \BB_I)$ is an equivalence with inverse given by $\rR\Gamma_\tau(X, \blank)\colon \Perf_\tau(X, \BB_I) \to \Perf\bigl(\BB_I(X)\bigr)$.
\end{cor}
\begin{proof}
    The proof of \cref{cor:perfect-O-modules-on-affinoid-perfectoid} applies almost verbatim. The first difference is that one needs to use \cref{thm:v-descent-rationally-BBI} in place of \cref{thm:v-descent-rationally}. Second, one needs to justify that the natural morphism $\BB_I(X) \to \rR\Gamma_{\tau}(X, \BB_I)$ is an isomorphism. This follows from the proof of \cite[Prop.\,II.2.1]{Fargues-Scholze}. 
\end{proof}

\subsection{$v$-descent for perfect $\ud{\ZZ}_\ell$ and $\ud{\QQ}_\ell$-complexes}

Throughout this section, we fix a prime number $\ell$. We will primarily be interested in the case $\ell=p$, however we do not make this assumption in this subsection. 

Recall that the sheaves $\ud{\ZZ}_\ell, \ud{\QQ}_\ell\colon \Perfd^{\aff, \op} \to \Ab$ are defined via the rule $\ud{\ZZ}_\ell(X) = \Map_\cont(\abs{X}, \ZZ_\ell)$ and $\ud{\QQ}_\ell(X) = \Map_\cont(\abs{X}, \QQ_\ell)$, where $\ZZ_\ell$ and $\Q_\ell$ are endowed with their standard topologies. Alternatively, we have $\ud{\ZZ}_\ell \simeq \lim_n \ud{\ZZ/\ell^n}$ and $\ud{\Q}_\ell \simeq \ud{\ZZ}_{\ell}[\frac{1}{\ell}]$, where each $\ud{\ZZ/\ell^n}$ is the standard constant sheaf. 

The main goal of this section is to show that perfect $\ud{\ZZ}_\ell$ and $\ud{\QQ}_\ell$-complexes coincide in the pro-\'etale and $v$-topologies. Since many of these arguments are essentially identical to the arguments used in previous subsections, we sometimes only sketch the proofs. 

We begin the subsection by collecting certain facts about the sheaves $\ud{\ZZ}_\ell$ and $\ud{\QQ}_\ell$. 

\begin{lemma}\label{lemma:continuous-functions-complete} Let $X$ be a topological space. Then the $\ZZ_\ell$-module $\Map_\cont(X, \ZZ_\ell)$ is $\ell$-adically complete and $\Map_\cont(X, \Z_\ell)/\ell^n\simeq \Map_\cont(X, \ZZ/\ell^n)$. 
\end{lemma}
\begin{proof}
    Clearly, we have $\Map_\cont(X, \Z_\ell) \simeq \lim_n \Map_\cont(X, \ZZ/\ell^n)$. So we only need to show that the natural map $\alpha_n\colon \Map_\cont(X, \Z_\ell)/\ell^n \to \Map_\cont(X, \ZZ/\ell^n)$ is an isomorphism. Using the (strict) short exact sequence $0 \to \ZZ_\ell \xr{\cdot \ell^n} \ZZ_\ell \to \ZZ/\ell^n \to 0$ of topological groups, we see that $\alpha_n$ is injective. It is also surjective because $\Z_\ell \to \ZZ/\ell^n$ admits a continuous set-theoretic section. 
\end{proof}

\begin{cor}\label{cor:z-ell-derived-complete} Let $X$ be a perfectoid space, let $\tau\in \{\proet, v\}$. Then the natural morphisms $\ud{\ZZ}_\ell \to \lim_n \ud{\ZZ/\ell^n}$ and $\ud{\ZZ}_\ell \to \rR\lim_n \ud{\ZZ/\ell^n}$ are isomorphisms in $\cD(X_\tau; \Z)$. Furthermore, the $\ZZ_\ell$-module $\ud{\ZZ}_\ell(X)$ is $\ell$-adically complete and $\ud{\ZZ}_\ell(X)/\ell^n\simeq \ud{\ZZ/\ell^n}(X)$. 
\end{cor}
\begin{proof}
    It follows from the definition that $\ud{\ZZ}_\ell \to \lim_n \ud{\ZZ/\ell^n}$ is an isomorphism. Thus, in order to see that $\ud{\ZZ}_\ell \to \rR\lim_n \ud{\ZZ/\ell^n}$ is an isomorphism, it suffices to show that $\lim_n \ud{\ZZ/\ell^n} \to \rR\lim_n \ud{\ZZ/\ell^n}$ is an isomorphism. This follows from \cref{lemma:replete} and \cite[Prop.\,3.1.10]{proetale}. The facts that $\ud{\ZZ}_\ell(X)$ is $\ell$-adically complete and that $\ud{\ZZ}_\ell(X)/\ell^n\simeq \ud{\ZZ/\ell^n}(X)$ follows immediately from \cref{lemma:continuous-functions-complete}. 
\end{proof}

\begin{lemma}\label{lemma:completion-colimit-Zell} Let $X=\lim_I X_i$ be a cofiltered limit of spectral spaces with spectral transition maps. Then the natural map $(\colim_{I^\op} \Map_\cont(X_i, \ZZ_\ell)\bigr)^{\wedge}_{(\ell)} \to \Map_\cont(X, \Z_\ell)$ is an isomorphism. 
\end{lemma}
\begin{proof}
        \cref{lemma:continuous-functions-complete} implies that $(\colim_{I^\op} \Map_\cont(X_i, \ZZ_\ell)\bigr)^{\wedge}_{(\ell)} \to \Map_\cont(X, \Z_\ell)$ is a homomorphism between $\ell$-adically complete $\ell$-torsionfree rings. Therefore, it suffices to show that this map is an isomorphism modulo $\ell$. Furthermore, \cref{lemma:continuous-functions-complete} also implies that the mod-$\ell$ reduction of this map is equal to $\colim_{I^\op} \Map_\cont(X_i, \FF_\ell) \to \Map_\cont(X, \FF_\ell)$. The fact that this map is a bijection follows from \cite[\href{https://stacks.math.columbia.edu/tag/0A30}{Tag 0A30}]{stacks-project}. 
\end{proof}

\begin{cor}\label{cor:completion-geometric-colimit-Zell} Let $X=\lim_I X_i$ be a cofiltered limit in $\Perfd^{\aff}$. Then $(\colim_{I^\op} \ud{\ZZ}_\ell(X_i)\bigr)^{\wedge}_{(\ell)} \to \ud{\ZZ}_\ell(X)$ is an isomorphism. 
\end{cor}
\begin{proof}
    This follows immediately from \cref{lemma:completion-colimit-Zell} and the observation that $\abs{X} \to \lim_I \abs{X_i}$ is a homeomorphism (see \cite[Cor.~C.3.12]{zav-almost}).
\end{proof}

\begin{lemma}\label{lemma:Q-ell-functions} Let $X$ be a quasi-compact topological space. Then the natural morphism $\Map_\cont(X, \ZZ_\ell)[\frac{1}{\ell}] \to \Map_\cont(X, \QQ_\ell)$ is an isomorphism. 
\end{lemma}
\begin{proof}
    First, we note that $\Map_\cont(X, \ZZ_\ell)[\frac{1}{\ell}] \simeq \colim_\NN \Map_\cont(X, \frac{1}{\ell^n}\ZZ_\ell)$. In order to check that the natural morphism $\colim_\NN \Map_\cont(X, \frac{1}{\ell^n}\ZZ_\ell) \to \Map_\cont(X, \QQ_\ell)$ is an isomorphism, it suffices to show that any continuous morphism $X \to \QQ_\ell$ factors through $\frac{1}{\ell^n}\ZZ_\ell$ for some integer $n$. This follows from the observations that $X$ is quasi-compact and each $\frac{1}{\ell^n}\ZZ_\ell \hookrightarrow \QQ_\ell$ are open immersions. 
\end{proof}

\begin{cor}\label{cor:spectral-vs-profinite} Let $X$ be a spectral topological space and let $n\geq 1$ be an integer. Then the natural morphisms $\Map_\cont(\pi_0(X), \ZZ/\ell^n) \to \Map_\cont(X, \ZZ/\ell^n) $, $\Map_\cont(\pi_0(X), \ZZ_\ell) \to \Map_\cont(X, \ZZ_\ell)$, and $\Map_\cont(\pi_0(X), \QQ_\ell) \to \Map_\cont(X, \QQ_\ell)$ are isomorphisms. 
\end{cor}
\begin{proof}
    Since $\ZZ/\ell^n$ is a finite discrete set, we see that the natural morphism $\Map_\cont(\pi_0(X), \ZZ/\ell^n) \to \Map_\cont(X, \ZZ/\ell^n)$ is an isomorphism by the universal property of $\pi_0$. Thus, \cref{lemma:continuous-functions-complete} implies that the natural $\Map_\cont(\pi_0(X), \ZZ_\ell) \to \Map_\cont(X, \ZZ_\ell)$ is an isomorphism as well. Then \cref{lemma:Q-ell-functions} and \cite[\href{https://stacks.math.columbia.edu/tag/0906}{Tag 0906}]{stacks-project} imply that $\Map_\cont(\pi_0(X), \QQ_\ell) \to \Map_\cont(X, \QQ_\ell)$ is an isomorphism as well. 
\end{proof}

For a quasi-compact topological space $X$, we topologize $\Map_\cont(X, \ZZ_\ell)$ via the $\ell$-adic topology and $\Map_\cont(X, \Q_\ell)\simeq \Map_\cont(X, \ZZ_\ell)[\frac{1}{\ell}]$ via the unique topology that makes $\Map_\cont(X, \ZZ_\ell)\subset \Map_\cont(X, \Q_\ell)$ into an open bounded subset. Since $\ZZ_\ell\subset \QQ_\ell$ is integrally closed, it follows that $\bigl(\Map_\cont(X,\QQ_\ell), \Map_\cont(X, \ZZ_\ell)\bigr)$ is a complete uniform Tate--Huber pair. 

\begin{lemma}\label{lemma:surjective-for-profinite} Let $T\to S$ be a surjective map of profinite sets. Then the morphism 
\[
\Spa\bigl(\Map_\cont(T,\QQ_\ell), \Map_\cont(T, \ZZ_\ell)\bigr) \to \Spa\bigl(\Map_\cont(S,\QQ_\ell), \Map_\cont(S, \ZZ_\ell)\bigr)
\]
is surjective. 
\end{lemma}
\begin{proof}
    By \cite[Lecture~III, proof of Th.~3.2, p.~20]{Condensed}, the surjection
\(T\twoheadrightarrow S\) can be written as a cofiltered inverse
limit of surjections \(T_j\twoheadrightarrow S_j\) of finite
discrete sets; see also
\cite[Chap.~I, \S1, Ex.~1, p.~11]{Neukirch-Schmidt-Wingberg}.
    Therefore, \cref{lemma:completion-colimit-Zell} and \cite[Cor.~C.3.12]{zav-almost} imply that it suffices to prove the result under the additional assumption that $T$ and $S$ are finite. In this case, we have $\Map_\cont(T, \ZZ_\ell)\simeq \prod_{t\in T} \ZZ_\ell$, $\Map_\cont(T, \QQ_\ell)\simeq \prod_{t\in T} \QQ_\ell$, and the same for $S$. Therefore, the map $\Spa\bigl(\Map_\cont(T,\QQ_\ell), \Map_\cont(T, \ZZ_\ell)\bigr) \to \Spa\bigl(\Map_\cont(S,\QQ_\ell), \Map_\cont(S, \ZZ_\ell)\bigr)$ is simply equal to the map $T \to S$, so it is surjective by the very assumption. 
\end{proof}

\begin{cor}\label{cor:surjective-ell-adically} Let $f\colon X \to Y$ be a surjective morphism of qcqs adic spaces. Then the following statements hold:
\begin{enumerate}[leftmargin=*,label={\upshape{(\roman*)}}]
    \item\label{cor:surjective-ell-adically-1} the morphism $\Spa\bigl(\Map_\cont(\abs{X}, \QQ_\ell), \Map_\cont(\abs{X}, \ZZ_\ell)\bigr) \to \Spa\bigl(\Map_\cont(\abs{Y}, \QQ_\ell), \Map_\cont(\abs{Y}, \ZZ_\ell)\bigr)$ is surjective;
    \item\label{cor:surjective-ell-adically-2} the morphisms
    \[
    \Spec \bigl(\Map_\cont(\abs{X}, \QQ_\ell)\bigr) \to \Spec \bigl(\Map_\cont(\abs{Y}, \QQ_\ell)\bigr) \text{ and } \Spec \bigl(\Map_\cont(\abs{X}, \ZZ_\ell)\bigr) \to \Spec \bigl(\Map_\cont(\abs{Y}, \ZZ_\ell)\bigr)
    \]
    are surjective on closed points. 
\end{enumerate}
\end{cor}
\begin{proof}
    The first part follows immediately from \cref{cor:spectral-vs-profinite}, \cref{lemma:surjective-for-profinite}, and \cite[\href{https://stacks.math.columbia.edu/tag/0906}{Tag 0906}]{stacks-project}. The second part follows from the first part,
    \cref{lemma:surjectivity-on-closed-points-rationally}, and \cref{lemma:surjective-integrally}. 
\end{proof}

Now we begin the proof of $v$-descent for $\Perf_\proet(X; \ud{\ZZ}_\ell)$. 

\begin{cor}\label{cor:perfect-ell-complete} Let $X$ be a perfectoid space, let $\tau\in \{\proet, v\}$, and let $\cP\in \Perf_\tau(X; \ud{\ZZ}_\ell)$. Then $\cP$ is derived $\ell$-adically complete.
\end{cor}
\begin{proof}
    The question is $\tau$-local on $X$, we can assume that $\cP$ is a strictly perfect complex. Then a simple inductive argument reduces the question to the case when $\cP$ is a direct summand $\ud{\ZZ}_\ell^{\oplus d}$. This case can be further reduced to $\cP = \ud{\ZZ}_\ell$. In this case, the result follows from \cref{cor:z-ell-derived-complete}. 
\end{proof}

\begin{cor}\label{cor:check-perfect-mod-ell} Let $X$ be a perfectoid space, let $\tau\in \{\proet, v\}$, and let $\cP\in \cD(X_\tau; \ud{\ZZ}_\ell)$ be derived $\ell$-adic complete object such that $\cP/^\rL \ell \in \Perf_\tau(X; \ud{\FF}_\ell)$. 
\begin{enumerate}[leftmargin=*,label={\upshape{(\roman*)}}]
    \item\label{cor:check-perfect-mod-ell-1} Then $\cP$ lies in $\Perf_\tau(X; \ud{\ZZ}_\ell)$.
    \item\label{cor:check-perfect-mod-ell-2} If $X$ is strictly totally disconnected, then $\rR\Gamma_\tau(X, \cP)\in \Perf(\ud{\ZZ}_\ell(X))$ and the natural map $\widetilde{\rR\Gamma_\tau(X, \cP)} \to \cP$ is an isomorphism.
\end{enumerate}
\end{cor}
\begin{proof}
    Part\cref{cor:check-perfect-mod-ell-1} is $\tau$-local on $X$, so we can assume that $X$ is a strictly totally disconnected perfectoid space. In this case, it suffices to prove Part\cref{cor:check-perfect-mod-ell-2}.  Now \cite[\href{https://stacks.math.columbia.edu/tag/0A0G}{Tag 0A0G}]{stacks-project} implies that $\rR\Gamma_\tau(X, \cP)$ is derived $\ell$-adically complete, while \cref{cor:pullback-etale-finitary} and \cref{lemma:example-z-mod-n} guarantee that $\rR\Gamma_\tau(X, \cP)/^\rL\ell \simeq \rR\Gamma_\tau(X, \cP/^\rL\ell) \in \Perf(\ud{\FF}_\ell(X))$. Furthermore, \cref{cor:z-ell-derived-complete} ensures that $\ud{\ZZ}_\ell(X)$ is $\ell$-adically complete. Therefore, \cite[\href{https://stacks.math.columbia.edu/tag/09AW}{Tag 09AW}]{stacks-project} implies that $\rR\Gamma_\tau(X, \cP) \in \Perf\big(\ud{\ZZ}_\ell(X)\big)$. To finish the proof, it suffices to show that the natural morphism
    \[
    \widetilde{\rR\Gamma_\tau(X, \cP)} \to \cP
    \]
    is an equivalence. Now \cref{cor:perfect-ell-complete} implies that $\widetilde{\rR\Gamma_\tau(X, \cP)}$ is derived $\ell$-adically complete, so an easy inductive argument guarantees that it suffices to show that 
    \[
    \widetilde{\rR\Gamma_\tau(X, \cP)}/^\rL \ell \simeq \widetilde{\rR\Gamma_\tau(X, \cP/^\rL \ell)} \to \cP/^\rL\ell
    \]
    is an isomorphism. This follows from \cref{cor:etale-perfect-complex-on-strictly-totally-disconnected} and \cref{lemma:example-z-mod-n}. 
\end{proof}

Finally, we arrive at the first major result of this subsection. 

\begin{cor}\label{cor:pullback-Z-ell} Let $X$ be a perfectoid space and let $\mu \colon (X_v, \ud{\ZZ}_\ell) \to (X_\proet,  \ud{\ZZ}_\ell)$ be the natural morphism of ringed sites. Then the following statements hold: 
\begin{enumerate}[leftmargin=*,label={\upshape{(\roman*)}}]
    \item\label{cor:pullback-Z-ell-1} the functor $\rR\mu_{*}\colon \cal{D}(X_v; \ud{\ZZ}_\ell) \to \cal{D}(X_\proet; \ud{\ZZ}_\ell)$ restricts to a functor $\rR\mu_{*} \colon \Perf_v(X; \ud{\ZZ}_\ell) \to \Perf_\proet(X; \ud{\ZZ}_\ell)$;
    \item\label{cor:pullback-Z-ell-2} the functor $\mu^{*}\colon \Perf_\proet(X; \ud{\ZZ}_\ell) \to \Perf_v(X; \ud{\ZZ}_\ell)$ is an equivalence of $\infty$-categories with the inverse functor $\rR\mu_*$.
\end{enumerate}
\end{cor}
\begin{proof}
    First, we show that $\rR\mu_*$ preserves perfect complexes. Observe that \cite[\href{https://stacks.math.columbia.edu/tag/0A0G}{Tag 0A0G}]{stacks-project} and \cref{cor:check-perfect-mod-ell} imply that it suffices to show that $(\rR\mu_*\cP)/^\rL \ell \simeq \rR\mu_*(\cP/^\rL \ell) \in \Perf_\proet(X; \ud{\FF}_\ell)$ for any $\cP\in \Perf_v(X; \ud{\ZZ}_\ell)$. This follows immediately from \cref{cor:pullback-etale-finitary} and \cref{lemma:example-z-mod-n}. Now it suffices to show that, for any $\cP\in \Perf_v(X; \ud{\ZZ}_\ell)$ and $\cal{Q}\in \Perf_\proet(X; \ud{\ZZ}_\ell)$, the natural morphisms 
    \[
    \mu^* \rR \mu_* \cP \to \cP \text{ and } \cal{Q} \to \rR \mu_* \mu^*\cal{Q} 
    \]
    are isomorphisms. Since all complexes involved are perfect, they are derived $\ell$-adically complete due to \cref{cor:perfect-ell-complete}. Therefore, it suffices that both morphisms are isomorphisms modulo $\ell^n$ for all $n\geq 1$. This follows from \cref{cor:pullback-etale-finitary} and \cref{lemma:example-z-mod-n}. 
\end{proof}

Unlike for $\O^+$-perfect complexes (see \cref{section:descend-O-plus}), it is not true that $\Perf_\et(X; \ud{\ZZ}_\ell) \simeq \Perf_\proet(X; \ud{\ZZ}_\ell)$. 

\begin{cor}\label{cor:ell-adic-integral-perfect-complex-on-strictly-totally-disconnected} Let $X$ be a strictly totally disconnected perfectoid space and let $\tau\in\{\proet, v\}$. Then the functor $\widetilde{(\blank)} \colon \Perf(\ud{\ZZ}_\ell(X)) \to \Perf_\tau(X; \ud{\ZZ}_\ell)$ is an equivalence. Its inverse is given by 
\[
\rR\Gamma_\tau(X, \blank) \colon \Perf_\tau(X; \ud{\ZZ}_\ell) \to \Perf\bigl(\ud{\ZZ}_\ell(X)\bigr).
\]
\end{cor}
\begin{proof}
    \cref{cor:check-perfect-mod-ell}~\cref{cor:check-perfect-mod-ell-2} implies that $\rR\Gamma_\tau(X, \blank)$ indeed restricts to a functor $\rR\Gamma_\tau(X, \blank) \colon \Perf_\tau(X; \ud{\ZZ}_\ell) \to \Perf\bigl(\ud{\ZZ}_\ell(X)\bigr)$. Thus, it suffices to show, for $\cP\in \Perf_\tau(X; \ud{\ZZ}_\ell)$ and $Q\in \Perf\bigl(\ud{\ZZ}_\ell(X)\bigr)$, the natural morphisms 
    \[
    \widetilde{\rR\Gamma_\tau(X, \cP)} \to \cP \text{ and } Q \to \rR\Gamma_\tau(X, \widetilde{Q})
    \]
    are isomorphisms. As all involved (complexes of) sheaves and modules are derived $\ell$-adically complete, we can check that these morphisms are isomorphisms modulo $\ell$. In which case, the results follow immediately from \cref{cor:etale-perfect-complex-on-strictly-totally-disconnected}
    and \cref{lemma:example-z-mod-n}. 
\end{proof}

\begin{cor}\label{cor:v-descent-ell-integrally} The functor $\Perf_\proet(\blank; \ud{\ZZ}_\ell) \colon \Perfd^{\aff,\op} \to \Cat_\infty$ is a hypercomplete $v$-sheaf of $\infty$-categories.
\end{cor}
\begin{proof}
    It follows immediately from \cref{cor:pullback-Z-ell}~\cref{cor:pullback-Z-ell-2} and \cref{lemma:perfect-complexes-are-hypercomplete}. 
\end{proof}

For future reference, we also need to establish a version of \cref{cor:v-descent-ell-integrally} for perfect complexes of bounded tor-amplitude. For the rest of this section, we work with the $\tau$-sheaves
\[
\Perf^{[r,s]}_\tau(\blank; \ud{\ZZ}_\ell) \colon \Perfd^{\aff, \op} \to \Cat_\infty
\]
for some prime number $\ell$, elements $r,s\in \ZZ\cup \{\pm \infty\}$, and $\tau \in \{\proet, v\}$ (see \cref{notation:sheaf-perfect-complexes-finite-amplitude}).

\begin{cor}\label{cor:ell-adic-integral-perfect-complex-on-strictly-totally-disconnected-tor-amplitude} Let $X$ be a strictly totally disconnected perfectoid space, let $\tau\in\{\proet, v\}$, and let $r,s\in \ZZ\cup \{-\infty, \infty\}$. Then the equivalence $\widetilde{(\blank)} \colon \Perf(\ud{\ZZ}_\ell(X)) \xr{\sim} \Perf_\tau(X; \ud{\ZZ}_\ell)$ induces an equivalence $\widetilde{(\blank)} \colon \Perf^{[r,s]}(\ud{\ZZ}_\ell(X)) \xr{\sim} \Perf^{[r,s]}_\tau(X; \ud{\ZZ}_\ell)$
\end{cor}
\begin{proof}
    Let $P\in \Perf\bigl(\ud{\ZZ}_\ell(X)\bigr)$. It suffices to show that $P$ has tor-amplitude in $[r,s]$ if and only if $\widetilde{P}$ lies in $\Perf^{[r,s]}_\tau(X; \uZ_\ell)$. The forward direction follows immediately from \cite[\href{https://stacks.math.columbia.edu/tag/0658}{Tag 0658}]{stacks-project}, so we only need to deal with the reverse direction. 

    Assume that $\widetilde{P}$ lies in $\Perf^{[r,s]}_\tau(X; \uZ_\ell)$. Then \cref{lemma:locally-essentially-surjective} implies that there is a $\tau$-covering $f\colon X'\to X$ such that $f_\tau^*\widetilde{P}\simeq \widetilde{Q}$, where $Q\in \Perf^{[r,s]}\bigl(\uZ_\ell(X')\bigr)$. Without loss of generality, we can assume that $X'$ is strictly totally disconnected as well. Then \cref{cor:ell-adic-integral-perfect-complex-on-strictly-totally-disconnected} ensures that $P\otimes^\rL_{\uZ_\ell(X)} \uZ_\ell(X') \simeq Q$. Therefore, \cref{cor:surjective-ell-adically} and \cite[\href{https://stacks.math.columbia.edu/tag/068V}{Tag 068V}]{stacks-project} ensure that $P\in \Perf^{[r,s]}(\uZ_\ell(X))$ finishing the proof. 
\end{proof}

\begin{lemma}\label{lemma:test-tor-amplitude-v-locally-ell-adically} Let $f\colon X \to Y$ be a $v$-covering in $\Perfd^\aff$, let $r, s\in \Z \cup \{-\infty, \infty\}$, and let $\cal{P}\in \Perf^{[r, s]}_v(Y; \ud{\ZZ}_\ell)$ such that $f_v^*\cal{P} \in \Perf^{[r, s]}_\proet(X; \ud{\ZZ}_\ell)$. Then $\cal{P}$ lies in $\Perf^{[r, s]}_\proet(Y; \ud{\ZZ}_\ell)$ (see \cref{rmk:fully-faithful-finite-tor-amplitude-integral}).
\end{lemma}
\begin{proof}
    The proof is essentially identical to that of \cref{cor:test-tor-amplitude-v-locally}. We only need to use \cref{cor:pullback-Z-ell}~\cref{cor:pullback-Z-ell-2} in place of \cref{thm:v-perfect-complexes-vs-etale-perfect-complexes-integrally}\cref{thm:v-perfect-complexes-vs-etale-perfect-complexes-integrally-2} and \cref{cor:surjective-ell-adically} in place of \cref{lemma:surjective-integrally}. 
\end{proof}

\begin{cor}\label{cor:integral-descend-finite-tor-amplitude-ell-adically} The natural morphism $\Perf_\proet^{[r, s]}(\blank; \ud{\ZZ}_\ell) \to \Perf_v^{[r, s]}(\blank; \ud{\ZZ}_\ell)$ is an equivalence for any $r, s\in \Z \cup \{-\infty, \infty\}$. In particular, $\Perf_\proet^{[r, s]}(\blank; \ud{\ZZ}_\ell)$ is a hypercomplete $v$-sheaf. 
\end{cor} 
\begin{proof}
    This is a formal consequence of \cref{cor:v-descent-ell-integrally} and \cref{lemma:test-tor-amplitude-v-locally-ell-adically}; see \cref{cor:integral-descend-finite-tor-amplitude} for an identical argument. 
\end{proof}

Fix an affinoid perfectoid space $S$. Then a perfect complex $\cal{P} \in \Perf^{[r,s]}_\tau(S; \ud{\QQ}_\ell)$ defines a morphism of sheaves 
\[
[\cal{P}] \colon * \to \Perf_\tau^{[r, s]}(\blank; \ud{\QQ}_\ell) \colon \Perfd^{\aff, \op}_{/S} \to \Cat_\infty.
\]

Using this, we define the sheaf of lattices as follows: 

\begin{defn}\label{defn:lattices-qell-perfectoid} The {\it $\tau$-sheaf of lattices} $\Latt_{\tau, \cal{P}}^{[r, s]}(\blank) \colon \Perfd^{\aff, \op}_{/S} \to \Cat_\infty$ is the fiber product
\[
\Latt^{[r,s]}_{\tau, \cal{P}}(\blank) \coloneqq * \times_{\Perf^{[r, s]}_\tau(-; \ud{\QQ}_\ell)} \Perf^{[r, s]}_\tau(-; \ud{\ZZ}_\ell).
\]
\end{defn}

Similar to subsections~\cref{section:descent-O}~and~\cref{section:descent-BI}, we get the following results:

\begin{lemma}\label{lemma:lattices-are-finitary-Qell} Let $S$ be an affinoid perfectoid and let $\cal{P}$ be an element of $\Perf^{[r, s]}_\tau(S; \ud{\QQ}_\ell)$ for some $r,s\in \Z \cup\{-\infty, \infty\}$
with $r \leq s$ and $\tau\in \{\proet, v\}$. Then the $\tau$-sheaf $\Latt^{[r, s]}_{\tau, \cal{P}}(\blank) \colon \Perfd^{\aff, \op}_{/S} \to \Cat_\infty$ is finitary.
\end{lemma}
\begin{proof}
    The proof is completely analogous to the proofs of \cref{lemma:lattices-are-finitary} and \cref{lemma:lattices-are-finitary-BBI}. One only needs to use \cref{lemma:completion-colimit-Zell} in place of \cite[Prop.~6.4]{diamonds} or \cref{lemma:completion-colimit-AAI}. 
\end{proof}

\begin{lemma}\label{lemma:lattice-etale-locally-Qell} Let $S$ be an affinoid perfectoid, let $\tau\in \{\proet, v\}$, and let $\cal{P}$ be an element of $\Perf^{[r, s]}_\tau(S; \ud{\QQ}_\ell)$ for some $r, s\in \Z \cup\{-\infty, \infty\}$. Then there is an \'etale covering $S' \to S$ and an $\cal{P}^+\in \Perf^{[r,s]}_\tau(S'; \ud{\ZZ}_\ell)$ such that $\cal{P}^+\otimes^\rL_{\ud{\ZZ}_\ell} \ud{\QQ}_\ell \simeq \cal{P}|_{S'}$.
\end{lemma}
\begin{proof}
    The proof is identical to the proofs of \cref{lemma:lattice-etale-locally} and of \cref{lemma:lattice-etale-locally-BBI}. We sketch the main steps. First, arguing as in the proof of \cref{lemma:lattice-etale-locally} (and using \cref{lemma:lattices-are-finitary-Qell} in place of \cref{lemma:lattices-are-finitary}), we can reduce to the case when $S=\Spa(C, C^+)$ for an algebraically closed nonarchimedean field $C$ and an open bounded valuation subring $C^+\subset C$ and $\cP\simeq \widetilde{P}$ for some $P\in \Perf^{[r,s]}(\ud{\QQ}_\ell(S))$. In this case, the result follows
    because $\ud{\ZZ_\ell}(S)=\ZZ_\ell$ and $\ud{\QQ}_\ell(S)=\QQ_\ell$, so the natural functor $\Perf^{[r,s]}(\ZZ_\ell) \to \Perf^{[r,s]}(\QQ_\ell)$ is essentially surjective.  
\end{proof}

\begin{thm}\label{thm:v-perfect-complexes-vs-etale-perfect-complexes-rationally-Qell} 
Let $X$ be a perfectoid space and let $\mu \colon (X_v, \ud{\QQ}_\ell) \to (X_\proet,  \ud{\QQ}_\ell)$ be the natural morphism of ringed sites. Then the following hold:
\begin{enumerate}[leftmargin=*,label={\upshape{(\roman*)}}]
    \item\label{thm:v-perfect-complexes-vs-etale-perfect-complexes-rationally-Qell-1} the functor $\rR\mu_{*} \colon \cal{D}(X_v; \ud{\QQ}_\ell) \to \cal{D}(X_\proet; \ud{\QQ}_\ell)$ restricts to a functor $\rR\mu_*\colon \Perf_v(X; \ud{\QQ}_\ell) \to \Perf_\proet(X; \ud{\QQ}_\ell)$;
    \item\label{thm:v-perfect-complexes-vs-etale-perfect-complexes-rationally-Qell-2} the functor $\mu^* \colon \Perf_\proet(X; \ud{\QQ}_\ell) \to \Perf_v(X; \ud{\QQ}_\ell)$ is an equivalence with inverse functor $\rR\mu_*$;
    \item\label{thm:v-perfect-complexes-vs-etale-perfect-complexes-rationally-Qell-3} if $X$ is strictly totally disconnected and $\tau\in \{\proet, v\}$, the functor $\widetilde{(\blank)}\colon \Perf\bigl(\ud{\QQ}_\ell(X)\bigr) \to \Perf_\tau(X; \ud{\QQ}_\ell)$ is an equivalence with inverse given by $\rR\Gamma_\tau(X, \blank)\colon \Perf_\tau(X; \ud{\QQ}_\ell) \to \Perf\bigl(\ud{\QQ}_\ell(X)\bigr)$.
\end{enumerate}
\end{thm}
\begin{proof}
    The argument is basically identical to those of \cref{thm:v-perfect-complexes-vs-etale-perfect-complexes-rationally} and \cref{cor:perfect-O-modules-on-affinoid-perfectoid}. We explain the main steps. 

    Part~\cref{thm:v-perfect-complexes-vs-etale-perfect-complexes-rationally-Qell-1} is pro-\'etale local on $X$. Therefore, \cref{lemma:lattice-etale-locally-Qell} ensures that we can assume that $\cP$ admits a $\ud{\ZZ}_\ell$-lattice. In this case, the result follows from \cref{cor:pullback-Z-ell}~\cref{cor:pullback-Z-ell-1} and the observation that $\mu$ is a coherent morphism of coherent topoi (so $\rR\mu_*$ commutes with filtered colimits of universally bounded below complexes). 

    Now we address Part~\cref{thm:v-perfect-complexes-vs-etale-perfect-complexes-rationally-Qell-2}. We first show that $\mu^*$ is fully faithful. This is equivalent to showing that the natural map $\cP \to \rR\mu_*\mu^* \cP$ is an isomorphism for any $\cP\in \Perf_\proet(X; \ud{\QQ}_\ell)$. This question is pro-\'etale local on $X$, so we can reduce to the case when $\cP$ is strictly perfect and then to the case when $\cP=\uQ_\ell$. In this case, the result follows from the fact that $\ud{\ZZ}_\ell \to \rR\mu_*\mu^* \ud{\ZZ}_\ell$ is an isomorphism (see \cref{cor:pullback-Z-ell}~\cref{cor:pullback-Z-ell-2}). To show essential surjectivity, Part~\cref{thm:v-perfect-complexes-vs-etale-perfect-complexes-rationally-Qell-1} implies that it suffices to show that the natural morphism $\mu^*\rR\mu_*\cP \to \cP$ is an isomorphism for any $\cP\in \Perf_v(X; \ud{\QQ}_\ell)$. The claim is pro-\'etale local on $X$, so \cref{lemma:lattice-etale-locally-Qell} ensures that we can assume that $\cP$ admits a $\ud{\ZZ}_\ell$-lattice. In this case, the result follows from \cref{cor:pullback-Z-ell}~\cref{cor:pullback-Z-ell-2}. 

    For Part~\cref{thm:v-perfect-complexes-vs-etale-perfect-complexes-rationally-Qell-3}, we argue as in the proof of \cref{cor:perfect-O-modules-on-affinoid-perfectoid} to reduce the general claim to showing that $P \to \rR\Gamma_\proet(X, \widetilde{P})$ is an isomorphism for any $P\in \Perf(\ud{\QQ}_\ell(X))$. Then we reduce this to the case when $P=\ud{\QQ}_\ell(X)$. In this case, the question boils down to showing that $\ud{\QQ}_\ell(X) \to \rR\Gamma_\proet(X, \ud{\QQ}_\ell)$ is an isomorphism. Since $X_\proet$ is a coherent site, we can further reduce this question to showing that the map $\ud{\ZZ}_\ell(X) \to \rR\Gamma_\proet(X, \ud{\ZZ}_\ell)$ is an isomorphism. This follows from \cref{cor:ell-adic-integral-perfect-complex-on-strictly-totally-disconnected}. 
\end{proof}

\begin{cor}\label{cor:v-perfect-complexes-vs-etale-perfect-complexes-rationally-Qell-tor-amplitude} Let $X$ be a strictly totally disconnected perfectoid space, let $\tau\in\{\proet, v\}$, and let $r,s\in \ZZ\cup \{-\infty, \infty\}$. Then the equivalence $\widetilde{(\blank)} \colon \Perf(\uQ_\ell(X)) \xr{\sim} \Perf_\tau(X; \uQ_\ell)$ induces an equivalence $\widetilde{(\blank)} \colon \Perf^{[r,s]}(\uQ_\ell(X)) \xr{\sim} \Perf^{[r,s]}_\tau(X; \uQ_\ell)$
\end{cor}
\begin{proof}
    The proof is completely analogous to that of \cref{cor:ell-adic-integral-perfect-complex-on-strictly-totally-disconnected-tor-amplitude}.
\end{proof}

\begin{thm}\label{thm:v-descent-rationally-Qell} Let $r, s\in \Z \cup \{-\infty, \infty\}$. Then the presheaf $\Perf^{[r, s]}_\proet\big(\blank; \ud{\QQ}_\ell\big) \colon \Perfd^{\aff, \op} \to \Cat_\infty$ is a hypercomplete $v$-sheaf of $\infty$-categories. In particular, the presheaf $\Perf_\proet\big(\blank; \ud{\QQ}_\ell\big)$ is a hypercomplete $v$-sheaf. 
\end{thm}
\begin{proof}
    First, \cref{thm:v-perfect-complexes-vs-etale-perfect-complexes-rationally-Qell} and \cref{lemma:perfect-complexes-are-hypercomplete} ensure that $\Perf_\proet(\blank; \ud{\QQ}_\ell) \simeq \Perf_v(\blank; \ud{\QQ}_\ell)$ is a hypercomplete $v$-sheaf. In order to deduce the result for any $r, s\in \Z\cup\{-\infty, \infty\}$, one argues as in the proof of \cref{cor:integral-descend-finite-tor-amplitude} using \cref{cor:surjective-ell-adically}~\cref{cor:surjective-ell-adically-2} in place of \cref{lemma:surjective-integrally} (which is used in the proof of \cref{cor:test-tor-amplitude-v-locally}). 
\end{proof}

\end{appendices}

\providecommand{\bysame}{\leavevmode\hbox to3em{\hrulefill}\thinspace}
\providecommand{\MR}{\relax\ifhmode\unskip\space\fi MR }
\providecommand{\MRhref}[2]{%
  \href{http://www.ams.org/mathscinet-getitem?mr=#1}{#2}
}
\providecommand{\href}[2]{#2}

\end{document}